\documentclass[11pt,twoside,a4paper]{book}

\usepackage[utf8]{inputenc}
\usepackage[T1]{fontenc}
\usepackage[brazilian,english]{babel}
\usepackage[pdftex]{graphicx}           % usamos arquivos pdf/png como figuras
\usepackage{setspace}                   % espaçamento flexível
\usepackage{indentfirst}                % indentação do primeiro parágrafo
\usepackage{makeidx}                    % índice remissivo
\usepackage[nottoc]{tocbibind}          % acrescentamos a bibliografia/indice/conteudo no Table of Contents
\usepackage{courier}                    % usa o Adobe Courier no lugar de Computer Modern Typewriter
\usepackage{type1cm}                    % fontes realmente escaláveis
\usepackage{listings}                   % para formatar código-fonte (ex. em Java)
\usepackage{titletoc}
\usepackage[fixlanguage]{babelbib}
\usepackage[font=small,format=plain,labelfont=bf,up,textfont=it,up]{caption}
\usepackage[usenames,svgnames,dvipsnames]{xcolor}
\usepackage[a4paper,top=2.54cm,bottom=2.0cm,left=2.0cm,right=2.54cm]{geometry} % margens
\usepackage[pdftex,plainpages=false,pdfpagelabels,pagebackref,colorlinks=true,citecolor=DarkGreen,linkcolor=NavyBlue,urlcolor=DarkRed,filecolor=green,bookmarksopen=true]{hyperref} 
\usepackage[all]{hypcap}                 % soluciona o problema com o hyperref e capitulos

\fontsize{60}{62}\usefont{OT1}{cmr}{m}{n}{\selectfont}

\usepackage{amsmath,amssymb,latexsym,soul,cite,mathrsfs,amsthm}

\usepackage{color,enumitem,graphicx}

\usepackage[colorinlistoftodos]{todonotes}
\makeatletter
\providecommand\@dotsep{5}
\def\listtodoname{List of Todos}
\def\listoftodos{\@starttoc{tdo}\listtodoname}
\makeatother

\usepackage{verbatim}  % to using the command "\comment"

\usepackage{graphicx}            % Too draw ...
\usepackage[all,2cell]{xy}       % ... arrows
\xyoption{2cell}
\usepackage[makeroom]{cancel}    % I added this for cancelling

\newcommand{\LA}{\ensuremath{\mathcal{LA}}} % Shorthand for LA-groupoids
\newcommand{\F}{\ensuremath{\mathcal F}}    % Foliation, Pseudofunctor
\newcommand{\G}{\mathcal G}                 % Lie groupoid and the one integrating a Lie algebroid as a suffix
\newcommand{\E}{\ensuremath{\mathcal{E}}}   % A representation up to homotopy
\newcommand{\XX}{\mathfrak{X}}              % Stack and space of vector fields
\newcommand{\CC}[1]{\mathfrak{C}({#1})}     % For the core of #1
\newcommand{\hh}{\mathfrak{h}}              % Lie subalgebra
\renewcommand{\gg}{\mathfrak{g}}            % Lie algebra
\newcommand{\ggl}{\mathfrak{gl}}            % linear Lie algebra
\newcommand{\Lie}{\mathcal{L}}              % Lie derivative
\DeclareMathOperator{\ad}{ad}               % adjoint

\newcommand{\A}{\mathcal{A}}                % An associative algebra
\newcommand{\rk}{\text{\rm rk}\,}         % The rank of a vector bundle
\DeclareMathOperator{\Img}{Im}              % Image: \Im is \mathrac{I} the imaginary part!!
\DeclareMathOperator{\coker}{coker}         % Cokernel
\DeclareMathOperator{\Bis}{Bis}             % Group of bisections
\newcommand{\Nn}{\mathbb N}                 % The natural numbers
\newcommand{\Zz}{\mathbb Z}					% The integer numbers
\newcommand{\Rr}{\mathbb R}					% The real numbers		
\newcommand{\rest}[1]{\Big{\vert}_{#1}}		% Restriction or evaluation of a map
\newcommand{\vJoin}{\mathbin{\rotatebox[origin=c]{90}{$\Join$}}} %Vertical operation in a double groupoid
\newcommand{\Tp}[1]{\overline{\underline{#1}}}        %Structural maps of the top groupoid in a double groupoid
\newcommand{\Lf}[1]{\lvert{#1}\rvert}                 %Structural maps of the top groupoid in a double groupoid
\newcommand{\Rg}[1]{{#1}_V}                           %Structural maps of the top groupoid in a double groupoid
\newcommand{\inc}[1]{#1^{\blacktriangleleft}}         %The inclusion of a Poisson Lie group in its double group
\newcommand{\dinc}[1]{#1^{\mathbin{\rotatebox[origin=c]{90}{$\blacktriangleright$}}}} %The inclusion of the dual in the double

\usepackage{fancyhdr}
\renewcommand{\chaptermark}[1]{\markboth{\MakeUppercase{#1}}{}}
\renewcommand{\sectionmark}[1]{\markright{\MakeUppercase{#1}}{}}

\numberwithin{equation}{section}
\newtheorem{theo}{Theorem}[section]

\newtheorem{prop}[theo]{Proposition}
\newtheorem{cor}[theo]{Corollary}
\newtheorem{rmk}[theo]{Remark}
\theoremstyle{definition} \newtheorem{ex:}[theo]{Example}
\theoremstyle{definition} \newtheorem{Def}[theo]{Definition}
\newtheorem{theorem}{Theorem}[section]
\newtheorem{lemma}[theorem]{Lemma}
\newtheorem{proposition}[theorem]{Proposition}

\newtheorem{remark}[theorem]{Remark}
\newcommand\restr[2]{{% we make the whole thing an ordinary symbol
    \left.\kern-\nulldelimiterspace % automatically resize the bar with \right
    #1 % the function
    \vphantom{\big|} % pretend it's a little taller at normal size
    \right|_{#2} % this is the delimiter
}}
\providecommand{\abs}[1]{\lvert#1\rvert}
\begin{document}
\frontmatter 
% cabeçalho para as páginas das seções anteriores ao capítulo 1 (frontmatter)
\fancyhead[RO]{{\footnotesize\rightmark}\hspace{2em}\thepage}
\setcounter{tocdepth}{2}
\fancyhead[LE]{\thepage\hspace{2em}\footnotesize{\leftmark}}
\fancyhead[RE,LO]{}
\fancyhead[RO]{{\footnotesize\rightmark}\hspace{2em}\thepage}

\onehalfspacing  % espaçamento

% ---------------------------------------------------------------------------- %
% CAPA EDWIN
\thispagestyle{empty}
\vspace*{1.96cm}%4.96  3.96

\noindent
\hspace*{3.5cm}
%\mbox{\vrule width 0.02cm \hspace{0.5cm}
\begin{minipage}[c][6cm]{9.7cm}%10.7
    \centering
    \textbf{A cohomology theory for Lie $2$-algebras and Lie $2$-groups}
    
    \vspace{1cm}%2cm
    
    Camilo Andr\'es Angulo Santacruz
    
\end{minipage}
%}
 
\vspace{1cm}

%\noindent
%\hspace*{3.5cm}
%\framebox[107mm][c]{Hello world There is a lot of fun}

\begin{center}
\textsc{
    Tese apresentada\\[-0.1cm] 
    ao\\[-0.1cm]
    Instituto de Matem\'atica e Estat\'stica\\[-0.1cm]
    da\\[-0.1cm]
    Universidade de S\~ao Paulo\\[-0.1cm]
    para\\[-0.1cm]
    obten\c{c}\~ao do t\'itulo\\[-0.1cm]
    de\\[-0.1cm]
    Doutor em Ci\'encias}
    
    \vskip 1.5cm
    Programa: Matem\'atica\\
    Orientador: Prof. Dr. Cristi\'an Andr\'es Ortiz Gonz\'alez\\

   	\vskip 1cm
    \normalsize{Durante o desenvolvimento deste trabalho o autor recebeu aux\'ilio
    financeiro da CAPES e do CNPq}
    
    \vskip 0.5cm
    \normalsize{S\~ao Paulo, Janeiro de 2018}
\end{center}
\newpage

\newpage
\thispagestyle{empty}
    \begin{center}
        \vspace*{2.3 cm}
        \textbf{\Large{A cohomology theory for Lie $2$-algebras and Lie $2$-groups}}\\
        \vspace*{2 cm}
    \end{center}

    \vskip 2cm

    \begin{flushright}
	Esta vers\~ao da tese cont\'em as corre\c{c}\~oes e altera\c{c}\~oes sugeridas\\
	pela Comiss\~ao Julgadora durante a defesa da vers\~ao original do trabalho,\\
	realizada em 01/05/2018. Uma c\'opia da vers\~ao original est\'a dispon\'ivel no\\
	Instituto de Matem\'atica e Estat\'istica da Universidade de S\~ao Paulo.
	
    \vskip 2cm
    
    \end{flushright}
    
    \begin{quote}
    \noindent Comiss\~ao Julgadora:
    
    \begin{itemize}
		\item Prof. Dr. Cristi\'an Andr\'es Ortiz Gonz\'alez (orientador) - IME-USP %[sem ponto final]
		\item Prof. Dr. Henrique Bursztyn - IMPA %[sem ponto final]
		\item Prof. Dr. Olivier Brahic - UFPR %[sem ponto final]
		\item Prof. Dr. Alejandro Cabrera - UFRJ %[sem ponto final]
		\item Prof. Dr. Ivan Struchiner - USP %[sem ponto final]
    \end{itemize}
      
    \end{quote}
\pagebreak
%\end{comment}

\pagenumbering{roman}     % começamos a numerar 

% ---------------------------------------------------------------------------- %
% Agradecimentos:
% Se o candidato não quer fazer agradecimentos, deve simplesmente eliminar esta página 
%\chapter*{Agradecimentos}
%Texto texto texto texto texto texto texto texto texto texto texto texto texto
%texto texto texto texto texto texto texto texto texto texto texto texto texto
%texto texto texto texto texto texto texto texto texto texto texto texto texto
%texto texto texto texto. Texto opcional.

% ---------------------------------------------------------------------------- %
% Resumo
\chapter*{Resumo}

\noindent ANGULO, C. \textbf{A cohomology theory for Lie $2$-algebras and Lie $2$-groups}. 
2018. 68 f.
Tese (Doutorado) - Instituto de Matem\'atica e Estat\'istica,
Universidade de S\~ao Paulo, S\~ao Paulo, 2018.
\\

Nesta tese, n\'os introduzimos uma nova teoria de cohomologia associada \`as $2$-\'algebras de Lie e uma nova teoria de cohomologia associada aos $2$-grupos de Lie. Prova-se que estas teorias de cohomologia estendem as teorias de cohomologia cl\'assicas de \'algebras de Lie e grupos de Lie em que os seus segundos grupos classificam extens\~{o}es. Finalmente, usaremos estos fatos junto com um morfismo de van Est adaptado para encontrar uma nova prova da integrabilidade das $2$-\'algebras de Lie.
\\

\noindent \textbf{Palavras-chave:} Teoria de Lie, geometria de ordem superior, cohomologia.

% ---------------------------------------------------------------------------- %
% The previous text in actual portuguese: Nesta tese nós introduzimos uma nova teoria de cohomologia associada às 2-álgebras de Lie e uma nova teoria de cohomologia associada aos 2-grupos de Lie. Nosso trabalho consiste em quatro capítulos. O Capítulo 1 contém alguns fatos básicos sobre a teoria de cohomologia de grupoides e algebroides de Lie, assim como as definições de 2-álgebras e 2-grupos de Lie e algumas ferramentas algébricas a serem usadas. No Capítulo 2 nós definimos o complexo triplo de cocadeias de uma 2-álgebra de Lie. A cohomologia desse complexo triplo generaliza a cohomologia de álgebras e de algebroides de Lie em tanto que o segundo grupo classifica extensões. No Capítulo 3, introduzimos o complexo triplo de um 2-grupo de Lie cujo segundo grupo de cohomologia classifica extensões.  Finalmente, no Capítulo 4, nós relacionamos essas novas cohomologias usando uma aplicação de van Est e provamos uma serie de isomorfismos análogos ao teorema de van Est clássico. Como uma aplicação, nós mostramos que, usando essas teoria e esses teoremas de van Est generalizados, pode-se dar uma nova prova para a integrabilidade das 2-álgebras de Lie de dimensão finita  para 2-grupos de Lie.
% ---------------------------------------------------------------------------- %
% Abstract
\chapter*{Abstract}
\noindent ANGULO, C. \textbf{A cohomology theory for Lie $2$-algebras and Lie $2$-groups}. 
2018. 68 p.
Thesis (PhD) - Instituto de Matem\'atica e Estat\'istica,
Universidade de S\~ao Paulo, S\~ao Paulo, 2018.
\\

In this thesis, we introduce a new cohomology theory associated to a Lie $2$-algebras and a new cohomology theory associated to a Lie $2$-group. These cohomology theories are shown to extend the classical cohomology theories of Lie algebras and Lie groups in that their second groups classify extensions. We use this fact together with an adapted van Est map to prove the integrability of Lie $2$-algebras anew.
\\

\noindent \textbf{Keywords:} Lie theory, higher geometry, cohomology.

% ---------------------------------------------------------------------------- %
% Sumário
\tableofcontents    % imprime o sumário
% ---------------------------------------------------------------------------- %
% Capítulos do trabalho
\chapter*{Introduction}
\addcontentsline{toc}{chapter}{Introduction}
\chaptermark{Introduction}

Lie theory, understood as the relationship between global geometric structures and their corresponding infinitesimal counter-parts, has been a particularly useful and rich line of study. Arguably, it underwent a series of major achievements ever since the beginning of the $21$st century, when Crainic and Fernandes \cite{CF2} gave necessary and sufficient conditions for a Lie algebroid to be integrable. The study of several incarnations of a categorified version of Lie theory started to emerge thereafter. However, among these, double structures had been lying around for a while. For instance, already back in the late $1980$s, Jiang-Hua Lu and Alan Weinstein \cite{LuWein} had come across double structures, when looking for the symplectic groupoids of Poisson-Lie groups. These were also noticed to play a r\^{o}le for Poisson-Lie groupoids \cite{MckXu}. Another example of a fairly known double structure is that of a classic representation of a Lie groupoid, which corresponds to a certain type of vector bundle over the same Lie groupoid. The full category of vector bundles over Lie groupoids, in order, has been proven \cite{VB&Reps,CristMat} to be equivalent to a relevant subclass of the more modern representations up to homotopy \cite{AAC,AAC2}. Among the abovementioned incarnations, we highlight four classes that outline how different in flavor they all are. First, there are enriched categories, where the spaces of arrows are given themselves the structure of a category \cite{Cristian&Olivier}. There are also the so-called (Lie) $n$-groupoids, which are simplicial manifolds that verify certain filling conditions (see e.g. \cite{Lie2Algbds,nGpds}). To the best of our knowledge, the infinitesimal counter-parts of these first categorifications are still poorly understood or even missing from the literature. On the other hand, what we have been calling double structures are groupoid objects internal to geometric categories such as the category of Lie algebroids or Lie groupoids \cite{2ndOrdGeom}. Finally, there are categorifications in which the associativity of the multiplication and the Jacobi identity hold only up to isomorphism \cite{Lie2Gps,Lie2Alg} and have recently deserved enough attention themselves, as they have appeared in connection with various physical theories such as higher gauge theory \cite{LoopGps} and also as integrations for infinite dimensional Lie algebras \cite{IntInfDim}. The higher structures that we will be chiefly concerned with lie in the intersection of the last two classes, these are the so-called (strict) Lie $2$-groups and Lie $2$-algebras.  \\

Several steps have been taken in the direction of understanding the \textit{double} Lie theory, but a full understanding is still far off. In particular, the theory of integration has been particularly elusive. Successful strategies for Lie $III$ theorems have so far been \textit{ad hoc} and refuse generalizations. Moreover, the elegant theory of integration developed in \cite{CF2} does not account for certain topological peculiarities that appear in the categorified theory. Among the examples of integration, there is the integration of the cotangent groupoid of a Poisson-Lie group to a double Lie group \cite{LuWein}, which was later generalized to a class of double products of groupoids in \cite{2ndOrdGeom}. Also, the integration of VB-algebroids to VB-groupoids was settled in \cite{BCD}. Incidentally, in sight of the results of \cite{VB&Reps,DoubleVB}, this latter is a result on the integration of $2$-term representations up to homotopy of Lie algebroids to $2$-term representations up to homotopy of Lie groupoids. Yet another integration result is that of \cite{LucaPhD}, in which certain lifting properties are imposed to ensure that the theory of \cite{CF2} gives the integrating object right away. Finally, Lie $2$-algebras have been integrated to Lie $2$-groups in \cite{ZhuInt2Alg}. \\

This thesis sprung out of attempting to understand the integration of these higher structures via their cohomologies. The classic integration result of van Est \cite{VanEst} is built on relating the differential cohomology of Lie groups and the cohomology of Lie algebras and this relationship was extended to the arena of groupoids and algebroids in \cite{VanEstC}. Needless to say, the cohomology theories of Lie algebras and Lie algebroids have been studied abundantly in the literature and have found numerous applications. The main theme of this thesis is to introduce cohomology theories that extend those of Lie algebroids and Lie groupoids for the simplest of the double objects: Lie $2$-algebras and Lie $2$-groups. In particular, keeping integration as a goal in sight, we prescribe the second cohomolgy group to classify appropriate extensions. One advantage that this approach has got is that it still works in infinite dimensions. Indeed, just as much as the space of sections of a Lie algebroid is a Lie algebra, any \LA -groupoid has an infinite dimensional Lie $2$-algebra of multiplicative sections \cite{CristJames} for which the theory of this thesis can be applied. The cohomology theories that we will present are novel though there does exist a cohomology theory for Lie $2$-groups in the literature \cite{2Coh}. This existing cohomology theory is rather a deRham cohomology than a group cohomology for the $2$-group seen as a double groupoid. In particular, the second group of the referred cohomology theory does not correspond to any kind of extension. \\

We outline the contents. In Chapter \ref{preliminarieschapter}, we collect some generalities and results that we will be using in the rest of the work. In Chapters \ref{fractionalchapter} and \ref{fractionalchapter2}, we introduce the respective complexes of cochains associated to Lie $2$-algebras and Lie $2$-groups respectively and provide an interpretation for the lower dimensions of its cohomologies. Finally, in Chapter \ref{estimateschapter} we relate the above referred complexes using an adapted van Est map and we prove a van Est type theorem. The main result of this chapter is a new proof for the integrability of Lie $2$-algebras. We proceed to give a more detailed account of the contents of each chapter. \\

In Chapter \ref{preliminarieschapter}, we start including some background material on the cohomology theories of Lie groupoids and Lie algebroids. Right after, we introduce Lie $2$-groups and Lie $2$-algebras, introduce the notation and give a wealth of examples. We continue with an exposition of the simplicial structures naturally associated to these double objects. To do so, we introduce the more general double Lie groupoids and \LA -groupoids, and outline the constructions of their respective double complexes. The next section deals with a series of results on homological algebra. In it, we recast van Est theorem in a convenient way for our future purposes. \\

In Chapter \ref{fractionalchapter}, we study Lie $2$-algebras from a cohomological point of view. Starting with the natural double complex associated to a Lie $2$-algebra, we prove that its total cohomology classifies extensions by the unit $1$-dimensional Lie $2$-algebra. We use this as a justification for studying higher representations. We introduce the linear Lie $2$-algebra. Then, we proceed to define $2$-representations and verify that in a general extension of a Lie $2$-algebra by an abelian Lie $2$-algebra, the choice of a splitting induces such a structure. Finally, we re-write the extensions as semi-direct sums whose structure gets twisted by a series of maps. Conversely, we deduce the equations that such maps need to satisfy in order to define an extension of Lie $2$-algebras. \\
After having recognized the equations defining these $2$-cocycles, we move on to classifying them up to isomorphism of extensions. In doing so, we recognize that the isomorphism is reduced to a couple of maps verifying a series of equations that relate them to the respective $2$-cocycles of each of the extensions. In the process of understanding these equations as cocycle conditions, we will see the emergence of a ``triple complex'' whose total complex has got a differential which encapsulates the equations defining both $2$-cocycles and their equivalences. \\

Next up, in Chapter \ref{fractionalchapter2}, we go over the analogous theory for Lie $2$-groups, culminating in a corresponding complex whose second cohomology classifies extensions of Lie $2$-groups by $2$-vector spaces. \\

Lastly, in Chapter \ref{estimateschapter}, we start off by reviewing both the integrability of Lie $2$-algebras and the strategy that we will refer to as the van Est strategy for the integration of Lie algebras. We move on to study the map that differentiates an extension of Lie $2$-groups to an extension of Lie $2$-algebras. We are bound to call such a map a van Est map. We prove that under certain connectedness assumptions the van Est map induces isomorphisms in cohomology. This is done in two parts. First, we prove it for the cohomology of Lie $2$-groups and Lie $2$-algebras with values in an honest vector space, and then we prove the general result. Interestingly, as a corollary of the first part we recognize the algebraic compatibility condition for the adjoint $2$-representation to have values on an honest vector space.\\
Finally, as an application, we prove the integrability of Lie $2$-algebras using this machinery. \\

In the appendices, we collect some results and definitions that we considered relevant, but out of scope for the main body of the text. In appendix A, we recall the path strategy for the  integration of Lie algebroids. In appendix B, we recall the definitions of the cores of \LA -groupoid and double Lie groupoids. In appendix C, we go over the relation between VB-groupoids and $2$-term representations up to homotopy. In the final appendix, we recall that the space of multiplicative sections of an \LA -groupoid is naturally a category, and can further be endowed with the structure of a Lie $2$-algebra. In analogy, we prove that there is $2$-group structure on an appropriate category of bisections of a double Lie groupoid. 
\mainmatter
% cabeçalho para as páginas de todos os capítulos
\fancyhead[RE,LO]{\thesection}
\singlespacing              % espaçamento simples
%\onehalfspacing            % espaçamento um e meio
%% ------------------------------------------------------------------------- %%
\chapter{Preliminaries}\label{preliminarieschapter}
%% ------------------------------------------------------------------------- %%
In this chapter, we recall the definitions of the main characters of this dissertation and settle notations and some of their properties. We are to follow the following conventions: by $1$-connectedness, we mean both connectedness and simply connectedness. Accordingly, by $s1$-connectedness, we mean that source fibres are $1$-connected. We will write $G^{(n)}$ for the space of $n$-tuples of composable arrows of a Lie groupoid $\xymatrix{G \ar@<0.5ex>[r]\ar@<-0.5ex>[r] & M}$. We write $G_x$ for the isotropy Lie group at $x\in M$ and $\mathcal{O}^G_x$ for the orbit through it. Associated to this groupoid, we will write $A_G$ for its Lie algebroid, that is $A_G:=u^*\ker ds$. We write $\phi ':=Lie(\phi )$ for the derivative of a Lie groupoid morphism restricted to its Lie algebroid, this map is the induced map between Lie algebroids. For an abstract Lie algebroid $A$ over $M$, we write $\rho_A$ for its anchor and for its isotropy at $x\in M$, $\gg_x(A)=\ker(\rho_A)_x$. \\
For a map $\xymatrix{f:M \ar[r] & N}$ and a vector bundle $E$ over $N$, we invariantly write the pull-back $f^*E:=M\times_N E=\lbrace (x,e)\in M\times E:e\in E_{f(x)}\rbrace$. Analogously, for a Lie algebroid $A$ over $N$, the Lie algebroid pull-back will always be $f^!A:=TM\times_{TN}A=\lbrace (v,a)\in TM\times A:df(v)=\rho_A(a)\rbrace$. 

% Go through Lie theorems for groups, groupoids and double groupoids. Remark how Lie 2 is related to the topological issue above. In particular, comment on why this is related to the mentioned fact about quotients. Since this fails to be the case for double Lie groupoids, illustrate with an example.

\section{Cohomology of Lie groupoids and Lie algebroids}
In this section we go about the theory of Lie groupoid and Lie algebroid cohomology. The notations and theorems of this section will be used abundantly throughout the text. Ultimately, they will help us state an integrability theorem as well, the strategy of which will be mimicked during the last chapter. Recall that for a Lie groupoid $G$, there is a simplicial structure on the nerve whose maps are given by 
\begin{eqnarray*}
\partial _k (g_0,...,g_p)=
  \begin{cases}
    (g_1,...,g_p)                       & \quad \text{if } k=0  \\
    (g_0,...,g_{k-1}g_k,...,g_p)  & \quad \text{if } 0<k\leq p\\
    (g_0,...,g_{p-1})                   & \quad \text{if } k=p+1, 
  \end{cases}
\end{eqnarray*}
for a given element $(g_0,...,g_p)\in G^{(p+1)}$. With these, one builds the complex 
\begin{eqnarray*}
C^p(G):=\lbrace\varphi\in C^\infty(G^{(p)}):\phi(g_1,...,g_p)=0,\textnormal{ if }g_k\in u(M)\textnormal{ for some }k\rbrace
\end{eqnarray*} 
of (normalized) Lie groupoid cochains with differential \begin{eqnarray*}
\xymatrix{
\partial :C ^\bullet (G) \ar[r] & C ^{\bullet +1} (G)
}
\end{eqnarray*}
defined by the usual formula
\begin{eqnarray*}
(\partial\varphi)(g_0,...,g_p) = \sum_{k=0}^{p+1} (-1)^k \partial_k^*\varphi(g_0,...,g_p),
\end{eqnarray*}
for $\varphi\in C^p(G)$. \\
Thus defined, $(C^\bullet (G),\partial)$ is referred to as the \textit{groupoid complex} of $G$, and its cohomology is the so-called \textit{differentiable cohomology} of $G$. \\
Additionally, there is a graded product
\begin{eqnarray*}
\xymatrix{
C^{p}(G)\times C^{q}(G) \ar[r] & C^{p+q}(G):(\varphi_1,\varphi_2) \ar@{|->}[r] & \varphi_1\star\varphi_2
}
\end{eqnarray*}
given by 
\begin{align*}
    \varphi_1\star\varphi_2(g_1,...,g_p;g_1',...,g_q') & :=\varphi_1(g_1,...,g_p)\varphi_2(g_1',...,g_q'), \\
    \varphi_1\star\varphi_2(g_1,...,g_p) & :=\varphi_1(g_1,...,g_p)\varphi_2(s(g_p)) & \textnormal{ for }q=0, \\
    \varphi_1\star\varphi_2(g_1,...,g_q) & :=\varphi_1(t(g_1))\varphi_2(g_1,...,g_q) & \textnormal{ for }p=0, \\
    \varphi_1\star\varphi_2              & :=\varphi_1\varphi_2 & \textnormal{ for }p=q=0.\\
\end{align*}
Such a product is compatible with the differential via a graded Leibniz rule
\begin{eqnarray*}
\partial(\varphi_1\star\varphi_2)=(\partial\varphi_1)\star\varphi_2 +(-1)^{p}\varphi_1\star(\partial\varphi_2)
\end{eqnarray*}
for $\varphi_1\in C^{p}(G)$. On the other hand, for a Lie algebroid $\xymatrix{A \ar[r] & M}$ with anchor $\rho$ and bracket $[\cdot ,\cdot ]$, there is a complex 
\begin{align*}
\Omega^q(A)=\Gamma\Big{(}\bigwedge^q A^*\Big{)}
\end{align*}
with the associated exterior differential of $A$,
\begin{eqnarray*}
\xymatrix{
d_A :\Omega ^\bullet (A) \ar[r] & \Omega ^{\bullet +1} (A)
}
\end{eqnarray*}
defined by the usual formula
\begin{align*}
(d_A\omega)(\alpha)& =\sum_{j=0}^q(-1)^j\rho(a_j)\omega(\alpha(j))+\sum_{m<n}(-1)^{m+n}\omega([a_m,a_n],\alpha(m,n)),
\end{align*}
for $\omega\in\Omega^q(A)$ and $\alpha=(a_0,...,a_q)\in\Gamma(A)^{q+1}$. Here, we use the convention that
\begin{eqnarray*}
\alpha(j)=(a_0,...,a_{j-1},a_{j+1},...,a_{q}) & \textnormal{and} & \alpha(m,n)=(a_0,...,a_{m-1},a_{m+1},...,a_{n-1},a_{n+1},...,a_{q}),
\end{eqnarray*}
as opposed to the usual $\hat{\cdot}$ notation. We will stick to this convention in the entirety of this text. \\
Thus defined, $(\Omega^\bullet (A), d_A)$ is usually referred to as the \textit{Chevalley-Eilenberg complex} of $A$, and its cohomology is the so-called \textit{Lie algebroid cohomology} of $A$. This time around there is also an additional graded product structure 
\begin{eqnarray*}
\xymatrix{
\Omega^{p}(A)\times\Omega^{q}(A) \ar[r] & \Omega^{p+q}(A):(\omega_1,\omega_2) \ar@{|->}[r] & \omega_1\wedge\omega_2
}
\end{eqnarray*}
that turns out to be compatible with the differential. It is given by the formula
\begin{align*}
   \omega_1\wedge\omega_2(a_1,...,a_p;a_{p+1},...,a_{p+q}) & :=\sum_{\sigma\in S(p,q)}\abs{\sigma}\omega_1(a_{\sigma(1)},...,a_{\sigma(p)})\omega_2(a_{\sigma(p+1)},...,a_{\sigma(p+q)}), 
\end{align*}
where $S(p,q)$ is the set of the so-called $(p,q)$-shuffles, that is
\begin{eqnarray*}
S(p,q):=\lbrace\sigma\in S_{p+q}:\sigma(1)<...<\sigma(p),\quad\sigma(p+1)<...<\sigma(p+q)\rbrace
\end{eqnarray*}
and $\abs{\cdot}$ stands for the sign of a permutation. The graded Leibniz rule is the familiar formula
\begin{eqnarray*}
d_A(\omega_1\wedge\omega_2)=(d_A\omega_1)\wedge\omega_2 +(-1)^{p}\omega_1\wedge(d_A\omega_2)
\end{eqnarray*}
for $\omega_1\in\Omega^{p}(A)$. \\
These complexes can take values on representations. Recall then, that a \textit{representation} of $\xymatrix{G \ar@<-0.5ex>[r]\ar@<0.5ex>[r] & M}$ is a vector bundle $E$ over $M$, together with an action
\begin{eqnarray*}
\xymatrix{
G_{s}\times_M E \ar[r] & E:(g,e) \ar@{|->}[r] & \Delta_g e
}
\end{eqnarray*}
along the projection of the vector bundle. Infinitesimally, these correspond to flat $A$-connections on $E$. An $A$-\textit{connection} is a map
\begin{eqnarray*}
\xymatrix{
\Gamma(A)\otimes_{\Rr}\Gamma(E) \ar[r] & \Gamma(E):(a,\epsilon) \ar@{|->}[r] & \nabla_a\epsilon
}
\end{eqnarray*}
that is $C^{\infty}(M)$-linear in $\Gamma(A)$, and that verifies a Leibniz rule in $\Gamma(E)$:
\begin{eqnarray*}
\nabla_a(f\epsilon)=f\nabla_a\epsilon+(\Lie_{\rho(a)}f)\epsilon.
\end{eqnarray*}
Such an $A$-connection is called \textit{flat} (and a representation of $A$), if it is compatible with the bracket in that
\begin{eqnarray*}
\nabla_{[a_1,a_2]}=[\nabla_{a_1},\nabla_{a_2}]
\end{eqnarray*}
for all $a_1,a_2\in\Gamma(A)$. Before laying down the complexes with values in a representation, let us specify that representations can be pulled-back along homomorphisms. In fact, if 
\begin{eqnarray*}
\xymatrix{
H \ar@<-0.5ex>[d]\ar@<0.5ex>[d]\ar[r]^\phi & G \ar@<-0.5ex>[d]\ar@<0.5ex>[d] \\
N \ar[r]_f                                 & M
}
\end{eqnarray*}
is a Lie groupoid homomorphism, there is an action
\begin{eqnarray*}
\xymatrix{
H_{s}\times_N f^*E \ar[r] & f^*E:(h;y,e) \ar@{|->}[r] & \Delta_{\phi(h)} e
}.
\end{eqnarray*}
Not as straightforward, if 
\begin{eqnarray*}
\xymatrix{
B \ar[d]\ar[r]^F & A \ar[d] \\
N \ar[r]_f           & M
}
\end{eqnarray*}
is a Lie algebroid homomorphism, there is a flat $B$-connection
\begin{eqnarray*}
\xymatrix{
\nabla':\Gamma(B)\otimes_{\Rr}\Gamma(f^*E) \ar[r] & \Gamma(f^*E)
}
\end{eqnarray*}
defined as follows: Since $\Gamma(f^*E)\cong C^\infty(N)\otimes_{C^\infty(M)}\Gamma(E)$, an element $\zeta\in\Gamma(f^*E)$ is given by
\begin{eqnarray*}
\zeta =\zeta^n\otimes\epsilon_n,
\end{eqnarray*}
where $\zeta^a\in C^\infty(N)$, $\epsilon_a\in\Gamma(E)$ and we used the Einstein summation convention. Analogously, given a section $b\in\Gamma(B)$, $F(b)$ induces a section of $f^*A\cong C^\infty(N)\otimes_{C^\infty(M)}\Gamma(A)$; therefore, 
\begin{eqnarray*}
F(b) =\beta^m\otimes a_m,
\end{eqnarray*}
for some $\beta^m\in C^\infty(N)$ and $a_m\in\Gamma(A)$. Then,
\begin{eqnarray}
\nabla'_b\zeta:=\beta^m\zeta^n\otimes\nabla_{a_m}\epsilon_n+\Lie_{\rho(b)}\zeta^n\otimes\epsilon_n.
\end{eqnarray}
The fact that this yields a flat $B$-connection follows from the fact that $F$ is a Lie algebroid homomorphism. \\
We move on to define the complexes of Lie groupoid and Lie algebroid cochains with values on a representation $E$. For the Lie groupoid $G$,
\begin{eqnarray*}
C^p(G;E):=\Gamma(t_p^*E)
\end{eqnarray*} 
where the map $\xymatrix{t_p:G^{(p)}\ar[r] & E:(g_1,...g_p) \ar@{|->}[r] & t(g_1)}$ is the map that returns the final target of a $p$-tuple of composable arrows. The differential
\begin{eqnarray*}
\xymatrix{
\partial :C ^\bullet (G,E) \ar[r] & C ^{\bullet +1} (G,E)
}
\end{eqnarray*}
is defined by the same formula as before, though modified in the first term so that the sum can be performed
\begin{eqnarray*}
(\partial\varphi)(g_0,...,g_p) = \Delta_{g_0}\partial_0^*\varphi(g_0,...,g_p)+\sum_{k=1}^{p+1}(-1)^k\partial_k^*\varphi(g_0,...,g_p).
\end{eqnarray*}
For the Lie algebroid $A$,
\begin{eqnarray*}
\Omega^q(A,E):=\Gamma(\bigwedge^q A^*\otimes E)
\end{eqnarray*}
with differential
\begin{eqnarray*}
\xymatrix{
d_\nabla :\Omega ^\bullet (A,E) \ar[r] & \Omega ^{\bullet +1} (A,E)
}
\end{eqnarray*}
defined by an analogous formula, though using the representation instead of the anchor,
\begin{align*}
d_{\nabla}\omega(\alpha)& =\sum_{j=0}^q(-1)^j\nabla_{a_j}\omega(\alpha(j))+\sum_{m<n}(-1)^{m+n}\omega([a_m,a_n],\alpha(m,n)).
\end{align*}
The complexes $C(G,E)$ and $\Omega(A,E)$ can be endowed with the structure of right (graded) $C(G)$ and $\Omega(A)$-modules respectively. These are given by the formulas above, though having a slightly different meaning, as the sums and multiplications by scalars are taking place in a fibre of the vector bundle instead of taking place in $\Rr$. \\
Now, roughly speaking, Lie algebroid cochains can be seen as the infinitesimal version of groupoid cocycles. This statement is made precise by means of a map
\begin{eqnarray*}
\xymatrix{
\Phi: C^\bullet(G,E) \ar[r] & \Omega^\bullet(A,E)
}
\end{eqnarray*}
whose value at $\omega\in C^p(G,E)$ for $a_1,...,a_p\in\Gamma(A)$ is given by the formula
\begin{eqnarray*}
\Phi\omega(a_1,...,a_p):=\sum_{\sigma\in S_p}\abs{\sigma}R_{a_{\sigma(1)}}...R_{a_{\sigma(p)}}\omega .
\end{eqnarray*}
In turn, the maps
\begin{eqnarray*}
\xymatrix{
R_a: C^{p+1}(G,E) \ar[r] & C^{p}(G,E)
}
\end{eqnarray*}
for $a\in\Gamma(A)$ are defined by 
\begin{eqnarray*}
R_a\omega(g_1,...,g_p):=\overrightarrow{a}\omega_{(g_1,...,g_p)}(u\circ t(g_1)),
\end{eqnarray*}
where $\overrightarrow{a}$ is the right-invariant vector field generated by $a$; hence, $\overrightarrow{a}\in\Gamma(\ker ds)$. On the other hand, the map
$\xymatrix{\omega_{(g_1,...,g_p)}:s^{-1}(t(g_1)) \ar[r] & E_{t(g_1)}:g \ar@{|->}[r] & \Delta_{g^{-1}}\omega(g,g_1,...,g_p)}$. We call this map the \textit{van Est map} for Lie groupoids. According to the author of \cite{VanEstC}, it was first constructed for Lie groups in \cite{VanEstMap} and later extended to the realm of Lie groupoids in \cite{WeinXu}. Thus defined, this map turns out to be a map of complexes compatible with the module structure, and there is the following theorem.
\begin{theorem}\cite{VanEstC}\label{CrainicVanEst}
If the source fibres of $G$ are $k$-connected, the van Est map induces isomorphisms
\begin{eqnarray*}
\xymatrix{
\Phi^n:H^n(G,E) \ar[r] & H^n(A,E)
}
\end{eqnarray*}
for $n\leq k$, and it is injective for $n=k+1$.
\end{theorem}
We refer to this theorem as the Crainic-van Est theorem in the sequel. 

We close this section by recalling an integrability result for Lie algebroids that, though not forgotten, can be overseen in sight of the definitive \cite{CF2}. 

The application of this theorem in which we will be interested, and the reason because of which we listed it here is the following.
\begin{theorem}\cite{VanEstC}\label{MariusInt}
Let 
\begin{eqnarray*}
\xymatrix{
(0) \ar[r] & E \ar[r] & \Omega \ar[r] & A \ar[r] & (0)
}
\end{eqnarray*}
be an exact sequence of Lie algebroids with $E$ abelian. If $A$ admits a Hausdorff integration whose $s$-fibres are simply connected and have vanishing second cohomology groups, then $\Omega$ is integrable. 
\end{theorem}
We will outline the proof of this theorem and its underlying strategy in the beginning of chapter \ref{estimateschapter} in what we will refer to as the van Est strategy.

%-------------------------------------------
%-------------------------------------------
% SUCCESSFUL STRATEGIES ... not anymore.
%-------------------------------------------
%-------------------------------------------

\section{On Lie 2-groups and Lie 2-algebras}

Lie algebra-like structures have been paid much attention in the literature. In particular, categorifications of these have been considered in \cite{Lie2Gps,Lie2Alg,ZhuInt2Alg}. Categorifications of group structures are less ubiquitous than their infinitesimal counter-parts. Nevertheless, they have appear in the literature \cite{Lie2Gps,2Coh,IntInfDim,StackyInt}, mainly related to applications of higher category theory to physics, but interestingly, also related to integration problems. \\
The objects we are going to study are known as strict Lie $2$-algebras and strict Lie $2$-groups in some parts of the literature.

We start by recalling the definition of a Lie $2$-algebra.
\begin{Def}
A \textit{Lie $2$-algebra} is a groupoid object internal to the category of Lie algebras. 
\end{Def}
We will write a generic Lie $2$-algebra as
\begin{eqnarray*}
\xymatrix{
\gg_1\times_\hh \gg_1 \ar[r]^{\quad \hat{m}} & \gg_1 \ar@<0.5ex>[r]^{\hat{s}} \ar@<-0.5ex>[r]_{\hat{t}} \ar@(l,d)[]_{\hat{\iota}} & \hh \ar[r]^{\hat{u}} & \gg_1 .
}
\end{eqnarray*}
Correspondingly, 
\begin{Def}
A \textit{Lie $2$-group} is a groupoid object internal to the category of Lie groups.
\end{Def}
We will write a generic Lie $2$-group as
\begin{eqnarray*}
\xymatrix{
\G\times_H\G \ar[r]^{\quad m} & \G \ar@<0.5ex>[r]^{s} \ar@<-0.5ex>[r]_{t} \ar@(l,d)[]_{\iota} & H \ar[r]^{u} & \G .
}
\end{eqnarray*}
In order to make clear the difference between the group operation and the groupoid operation in $\G$, we assume the following convention: 
\begin{eqnarray*}
g_1\vJoin g_2 & & g_3\Join g_4,
\end{eqnarray*}
stand respectively for the group multiplication and the groupoid multiplication whenever $(g_1,g_2)\in\G^2$ and $(g_3,g_4)\in\G\times_H\G$. This intends to reflect the fact that the group multiplication is ``vertical'', whereas the groupoid multiplication is ``horizontal''.

It has been noted in \cite{Lie2Gps,Lie2Alg,ZhuInt2Alg} and elsewhere that the categories of Lie $2$-algebras and Lie $2$-groups are respectively equivalent to the categories of crossed modules of Lie algebras and of Lie groups.
\begin{Def}
A \textit{crossed module of Lie algebras} is a Lie algebra morphism $\xymatrix{\gg \ar[r]^\mu & \hh}$ together with a Lie algebra action by derivations $\xymatrix{\Lie :\hh \ar[r] & \ggl (\gg )}$ satisfying
\begin{eqnarray*}
\mu(\Lie _y x) & = & [y, \mu (x)] ,\\
\Lie _{\mu(x_0)}x_1 & = & [x_0 ,x_1 ].
\end{eqnarray*} 
These equations are referred to as equivariance and infinitesimal Peiffer respectively.
\end{Def} 
The equivalence between the category of Lie $2$-algebras and crossed modules of Lie algebras is given, at the level of objects, by the following. Associated to a crossed module $\xymatrix{\gg \ar[r]^\mu & \hh}$, the space of arrows of the Lie $2$-algebra is defined to be the semi-direct sum $\gg\oplus_\Lie\hh$, where the bracket is given by the usual formula
\begin{eqnarray*}
[(x_0 ,y_0 ),(x_1 ,y_1 )]_\Lie:=([x_0 ,x_1 ]+\Lie _{y_0}x_1 -\Lie _{y_1}x_0 ,[y_0 ,y_1 ]).
\end{eqnarray*}
The structural maps are given by 
\begin{align*}
\hat{s}(x,y)=y & \qquad\hat{t}(x,y)=y+\mu (x) & \hat{\iota}(x,y)=(-x,y+\mu (x)) & \qquad\hat{u}(y)=(0,y)            
\end{align*}
\begin{eqnarray*}
(x',y+\mu (x))\Join (x,y):=\hat{m}(x',y+\mu (x); x,y):=(x+x' ,y).
\end{eqnarray*}
Conversely, given a Lie $2$-algebra $\xymatrix{\gg_1 \ar@<0.5ex>[r] \ar@<-0.5ex>[r] & \hh}$, the associated crossed module is given by $\xymatrix{\ker\hat{s} \ar[r]^{\quad\hat{t}\vert_{\ker\hat{s}}} & \hh}$. The action for $\xi\in\ker \hat{s}$ and $y\in\hh$ is given by
\begin{eqnarray*}
\Lie _y\xi :=[\hat{u}(y),\xi]_1,
\end{eqnarray*} 
where $[,]_1$ is the bracket of $\gg_1$.\\
From here on out, we make no distinction between a Lie $2$-algebra and its associated crossed module.\\

\begin{Def}
A \textit{crossed module of Lie groups} is a Lie group homomorphism $\xymatrix{G \ar[r]^i & H}$ together with a right action of $H$ on $G$ by Lie group automorphisms satisfying
\begin{eqnarray*}
i(g^h)       & = & h^{-1}i(g)h,\\
g_1^{i(g_2)} & = & g_2^{-1}g_1g_2,
\end{eqnarray*} 
where we write $g^h$ for $h$ acting on $g$. We refer to the second relation as the Peiffer equation. 
\end{Def}

For future reference, we outline the correspondence. Given a crossed module $\xymatrix{G\ar[r]^i & H}$, the space of arrows of the Lie $2$-group associated to a given crossed module is the semi-direct product $G\rtimes H$ with respect to the $H$-action, that is
\begin{eqnarray*}
(g_1,h_1)\vJoin (g_2,h_2)=(g_1^{h_2}g_2 ,h_1 h_2),
\end{eqnarray*}
and the structural maps are the analogous
\begin{align*}
s(g,h)=h & \qquad t(g,h)=hi(g) & \iota(g,h)=(g^{-1},hi(g)) & \qquad u(h)=(e,h) 
\end{align*}
\begin{eqnarray*}
(g',hi(g))\Join (g,h):=(gg' ,h).
\end{eqnarray*}
Conversely, starting out with a Lie $2$-group $\xymatrix{\G \ar@<0.5ex>[r] \ar@<-0.5ex>[r] & H}$, the associated crossed module is given by $\xymatrix{\ker s \ar[r]^{\quad t\vert_{\ker s}} & H}$. The right action for $g\in\ker s$ by an element $h\in H$ is given by
\begin{eqnarray*}
g^h :=u(h)^{-1}\vJoin g\vJoin u(h),
\end{eqnarray*} 
where $\vJoin$ is the product in $\G$ and the $-1$ power stands for the inverse of this product. In other words, the action is given by conjugation in $\G$ by units. \\
From here on out, we make no distinction between a Lie $2$-group and its associated crossed module either. It should be clear that Lie $2$-groups and Lie $2$-algebras are related to one another by the Lie functor.\\
These objects abound in mathematics. Lie algebras and Lie groups are examples whose associated crossed modules are \begin{eqnarray*}
\xymatrix{
(0) \ar[r] & \gg
} & \textnormal{and} & \xymatrix{
1 \ar[r] & G .
}
\end{eqnarray*}
Ideals $\hh\trianglelefteq\gg$ and normal subgroups $N\trianglelefteq G$ are examples whose associated crossed modules are given by the inclusions 
\begin{eqnarray*}
\xymatrix{
\hh \ar[r] & \gg
} & \textnormal{and} &  \xymatrix{
N \ar[r] & G 
}
\end{eqnarray*}
and the actions by the adjoint and by conjugation respectively. Also, among these, one can find classical representations of Lie algebras and Lie groups on a vector space $V$; the associated crossed modules being
\begin{eqnarray*}
\xymatrix{
V \ar[r]^0 & \gg
} & \textnormal{and} &  \xymatrix{
V \ar[r]^1 & G 
}
\end{eqnarray*}
where the actions are the ones given by the representations themselves. We will see in the sequel that indeed these are morally all the examples, as any other will be a suitable combination. \\
Out of the crossed modules, one is able to read much of the groupoid theoretical data of both Lie $2$-groups and Lie $2$-algebras. Indeed, let $\gg_1$ be a Lie $2$-algebra with associated crossed module $\xymatrix{\gg \ar[r]^\mu & \hh}$ with action $\Lie$, then we have got the following collection of trivial observations.
\begin{lemma}
The orbit through $0$ is $\mu(\gg)$ and an ideal in $\hh$.
\end{lemma}
\begin{proof}
By definition, the orbit through $0$ is $\hat{t}(\hat{s}^{-1}(0))$. Now, $\hat{s}(x,y)=y$; thus, $\hat{s}^{-1}(y)=\gg\times\lbrace y\rbrace$. On the other hand, $\hat{t}(x,y)=y+\mu(x)$; therefore, $\hat{t}(\hat{s}^{-1}(0))=\hat{t}(\gg\times (0))=\mu(\gg)$ as claimed. As for the second part of the statement, the equivariance of $\mu$ yields, $[y,\mu(x)]=\mu(\Lie_yx)\in\mu(\gg)$.

\end{proof}
\begin{lemma}
The orbit through any point $y\in\hh$ is the coset of $\mu(\gg)$ in $\hh$ with respect to $y$; thus, every Lie $2$-algebra is regular and its leaf space is the Lie algebra $\hh/\mu(\gg)$.
\end{lemma}
One could be interested in Lie $2$-algebras precisely because they serve as models to study quotients. This will be made much clearer in the case of their global counter-parts, Lie $2$-groups, when the topology comes into play.
\begin{lemma}
The isotropy group of $\gg_1$ at $0$ is the central Lie subalgebra $\ker\mu$. All other isotropy groups for $y\in\hh$ are isomorphic.
\end{lemma}
\begin{proof}
First, notice that if $x_1,x_2\in\ker\mu$, by the infinitesimal Peiffer equation, $[x_1,x_2]=\Lie_{\mu(x_1)}x_2=0$; hence, $\ker\mu$ is indeed an abelian Lie algebra. By the same equation, taking $x_2$ an arbitrary element of $\gg$, it is central. Now, by definition the isotropy at $y$ is the intersection of its source and target fibres. Since, $\hat{t}^{-1}(y)=\lbrace(x,y')\in\gg\oplus\hh :y'+\mu(x)=y\rbrace$, $(\gg_1)_y=\lbrace (x,y)\in\gg\oplus\hh :y-\mu(x)=y\rbrace=\ker\mu\times\lbrace y\rbrace$ and the result follows.

\end{proof}
In fact, if one regards $\gg_1$ as an \LA -groupoid, the vector space $\ker\mu$ corresponds to a type of isotropy associated to the unique point $\ast$ in the double base.
\begin{lemma}
The associated action of $\hh$ on the second isotropy is the honest representation defined by the restriction of $\Lie$. Moreover, there is an honest representation $\overline{\Lie}$ of the orbit space on $\ker\mu$ of which the previous representation is a pull-back.
\end{lemma}
\begin{proof}
We just need to prove that the representation is well-defined. Let $x\in\ker\mu$ and $y\in\hh$, then by equivariance, $\mu(\Lie_yx)=[y,\mu(x)]=0$. As for the second statement, define $\overline{\Lie}_{[y]}v=\Lie_yv$. It is well-defined too, as for any other representative $y'=y+\mu(x)$, $\Lie_{y'}v=\Lie_yv +\Lie_\mu(x)v$ and the infinitesimal Peiffer equation implies $\Lie_\mu(x)v=[x,v]=-[v,x]=-\Lie_{\mu(v)}x=0$. The rest follows trivially.

\end{proof}
These observations make clear that the data of a crossed module of Lie algebras can be reduced to a $4$-tuple $(\hh,I,V,\rho)$ of a Lie algebra $\hh$, an ideal $I\trianglelefteq\hh$, a vector space $V$ and a representation of $\hh/I$ on $V$. The crossed module of such a $4$-tuple is
\begin{eqnarray}
\xymatrix{
V\oplus I \ar[r] & \hh:(v,x) \ar@{|->}[r] & x,
}
\end{eqnarray}
where the Lie algebra structure on $V\oplus I$ is the direct product, and the action $\Lie_y(v,x)=(\rho_{[y]}v,[y,x])$.
% This should reduce Chenchang's proof to integrating representations and integrating 'special' subobjects.
As in the case of Lie $2$-algebras, one is able to read much of the groupoid theoretical data by looking at the associated crossed module of a given Lie $2$-group. Indeed, let $\G$ be a Lie $2$-group with associated crossed module $\xymatrix{G \ar[r]^i & H}$, where we write the right action of $h\in H$ on $g\in G$ by $g^h$, then, we can collect the following easy observations.
\begin{lemma}
The orbit through the identity $1$ is $i(G)$ and a normal subgroup of $H$.
\end{lemma}
\begin{proof}
By definition, the orbit through $1$ is $t(s^{-1}(1))$. Now, $s(g,h)=h$; thus, $s^{-1}(h)=G\times\lbrace h\rbrace$. On the other hand, $t(g,h)=hi(g)$; therefore, $t(s^{-1}(1))=t(G\times\lbrace 1\rbrace)=i(G)$ as claimed. As for the second part of the statement, the equivariance of $i$ yields, $i(g^h)=h^{-1}i(g)h\in i(G)$.

\end{proof}
\begin{lemma}
The orbit through any element $h\in H$ is the coset of $i(G)$ in $H$ with respect to $h$; thus, every Lie $2$-group is regular and its leaf space inherits the structure of the group $H/i(G)$.
\end{lemma}
As we remarked above, in the case of Lie $2$-groups it is made patent that the study of Lie $2$-groups might be helpful to study quotients. Simple examples as the inclusion of an irrational slope subgroup in the torus, whose quotient is simultaneously a group and a badly behaved topological space, can be studied using tools of differential geometry when looking at its Lie $2$-group.
\begin{lemma}
The isotropy group of $\G$ at $1$ is the central Lie subgroup $\ker i$. All other isotropy groups for $h\in H$ are isomorphic.
\end{lemma}
\begin{proof}
Using the Peiffer equation for elements $g_1,g_2\in\ker i$, $g_2^{-1}g_1g_2=(g_1)^{i(g_2)}=g_1$. Hence, $\ker i$ is indeed an abelian Lie group. By the same argument, taking an arbitrary $g_1$, as opposed to one in the kernel of $i$, one sees that it is central. Now, since  $t^{-1}(h)=\lbrace(g,h')\in G\times H :h'i(g)=h\rbrace$, the isotropy group at $h$ is $\G_h=\lbrace (g,h)\in G\times H:hi(g)^{-1}=h\rbrace=\ker i\times\lbrace h\rbrace$ and the result follows.

\end{proof}
In fact, the vector space $\ker i$ can also be regarded as a type of isotropy of $\G$ at the unique point $\ast$, when regarded as a double Lie groupoid.
\begin{lemma}
The associated action of $H$ on $\ker i$ is the honest representation defined by the restriction of the right action to $\ker i$. Moreover, this representation is the pull-back of a representation of the orbit space. 
\end{lemma}
\begin{proof}
We prove that the representation can be restricted. Let $g\in\ker i$ and $h\in H$, then by equivariance, $i(g^h)=h^{-1}i(g)h=1$. As for the second statement, define $g^{[h]}=g^h$. This is well-defined, because for any other representative $h'=hi(\gamma)$, $g^{h'}=g^{hi(\gamma)}$. Thus, the Peiffer equation implies $g^{hi(\gamma)}=\gamma^{-1}g^h\gamma$, but we showed that $\ker i$ is central. The rest follows trivially.

\end{proof}
We can sum up these observations in an analogous conclusion to that which we drew with Lie $2$-algebras, that given the data of a $4$-tuple $(H,N,A,\rho)$ of a Lie group $H$, a normal subgroup $N\trianglelefteq H$, an abelian Lie group $A$ and a right action of $H/N$ on $A$, there is a crossed module
\begin{eqnarray*}
\xymatrix{
A\times N \ar[r] & H:(a,n) \ar@{|->}[r] & n,
}
\end{eqnarray*}
where the Lie group structure on $A\times N$ is the direct product, and the action $(a,n)^h=(a^{[h]},h^{-1}nh)$.

%-----------------------------------------
% THE UNDERLYING STRUCTURE
%-----------------------------------------
\subsection{The simplicial objects associated to doubles}
In this section, we point out that the cochain complexes of Lie $2$-groups and Lie $2$-algebras regarded as Lie groupoids are compatible with the respective complexes that each space of the nerve has got. We will do so by regarding Lie $2$-groups and Lie $2$-algebras as more general objects; namely, as double groupoids and \LA -groupoids respectively. 
\subsubsection{Double structures}
Double Lie groupoids are groupoid objects in the category of Lie groupoids. We write the definition explicitly.
\begin{Def}
A \textit{double Lie groupoid} consists of a square
\begin{eqnarray}\label{DoubleGpd}
\xymatrix{
 D \ar@<0.5ex>[r] \ar@<-0.5ex>[r] \ar@<0.5ex>[d] \ar@<-0.5ex>[d] & V \ar@<0.5ex>[d] \ar@<-0.5ex>[d]\\
 H \ar@<0.5ex>[r] \ar@<-0.5ex>[r] & M,
}
\end{eqnarray} 
where each side is a Lie groupoid, the structural maps are smooth functors, and such that the double source map
\begin{eqnarray*}
\xymatrix{
\mathbb{S}:=(\Lf{s},\Tp{s}): D \ar[r] & H_s\times _{\Rg{s}} V
}
\end{eqnarray*}
is a submersion. 
\end{Def}
The elements in $D$ can be interpreted as being squares whose vertical edges are arrows in $V$ and whose horizontal edges are arrows in $H$. Needless to say, the vertices will correspond to points in $M$. In order to recognize the structural maps then, we adopt the following mnemonic device. We write $\Tp{s}$, $\Tp{t}$, etc. for the structural maps of the top groupoid; $\Lf{s}$, $\Lf{t}$, etc. for the left vertical groupoid; $\Rg{s}$, $\Rg{t}$, etc. for the right vertical groupoid; and finally, the usual $s$, $t$, etc. for the bottom groupoid. For a given element $d\in D$, the square has the following edges
\begin{eqnarray*}
\xymatrix{
\bullet                    & \bullet \ar[l]_{\Lf{t}(d)} \\
\bullet \ar[u]^{\Tp{t}(d)} & \bullet \ar[l]^{\Lf{s}(d)}\ar[u]_{\Tp{s}(d)}
}
\end{eqnarray*} 
Additionally, as we did for Lie $2$-groups, we use the shorthand 
\begin{eqnarray*}
d_1\Join d_2:=\Tp{m}(d_1,d_2) & & d_3\vJoin d_4:=\Lf{m}(d_3,d_4),
\end{eqnarray*}
whenever $(d_1,d_2)\in D\times_V D$ and $(d_3,d_4)\in D\times_H D$ to reflect the fact that $d_1\Join d_2$ and $d_3\vJoin d_4$ are respectively
\begin{eqnarray*}
\xymatrix{
\bullet                      & \bullet \ar[l]_{\Lf{t}(d_1)}                                   & \bullet \ar[l]_{\Lf{t}(d_2)} \\
\bullet \ar[u]^{\Tp{t}(d_1)} & \bullet \ar[l]^{\Lf{s}(d_1)}\ar[u]^{\Tp{s}(d_1)}_{\Tp{t}(d_2)} & \bullet \ar[l]^{\Lf{s}(d_2)}\ar[u]_{\Tp{s}(d_2)}
} & \textnormal{and} & \xymatrix{
\bullet                      & \bullet \ar[l]_{\Lf{t}(d_3)} \\
\bullet \ar[u]^{\Tp{t}(d_3)} & \bullet \ar[l]_{\Lf{s}(d_3)}^{\Lf{t}(d_4)}\ar[u]_{\Tp{s}(d_3)} \\
\bullet \ar[u]^{\Tp{t}(d_4)} & \bullet \ar[l]^{\Lf{s}(d_4)}\ar[u]_{\Tp{s}(d_4)}
}
\end{eqnarray*}
With this notation, the fact that the multiplication in either groupoid is a groupoid homomorphism yields the formula
\begin{eqnarray*}
(d_1\Join d_2)\vJoin(d_3\Join d_4)=(d_1\vJoin d_3)\Join(d_2\vJoin d_4),
\end{eqnarray*} 
whenever it makes sense. We call this formula the interchange law.\\
In order, \LA -groupoids can be seen as first infinitesimal approximations of double Lie groupoids. That is, groupoid objects in the category of Lie algebroids.
\begin{Def}
An \LA -groupoid consists of a square
\begin{eqnarray}\label{LAGpd}
\xymatrix{
 \Omega \ar@<0.5ex>[r] \ar@<-0.5ex>[r] \ar[d]_{\pi}  & A \ar[d]^{p} \\
 H \ar@<0.5ex>[r] \ar@<-0.5ex>[r] & M ,
}
\end{eqnarray}
where $\Omega$ and $A$ are Lie algebroids over $H$ and $M$, the top and bottom sides are Lie groupoids, and the top structural maps are Lie algebroid morphisms covering the bottom ones. Further, the double source map
\begin{eqnarray*}
\xymatrix{
(\pi ,\hat{s}): \Omega \ar[r] & H_s\times _p A
}
\end{eqnarray*}
is surjective. 
\end{Def}
We will write the top structural maps by $\hat{s}$, $\hat{t}$, and so on.\\
The following proposition justifies us calling an \LA -groupoid an infinitesimal approximation of a double Lie groupoid.
\begin{prop}\cite{2ndOrdGeom}
Given a double Lie groupoid as \ref{DoubleGpd}, applying the Lie functor to the vertical groupoids yields a naturally associated \LA -groupoid
\begin{eqnarray*}
\xymatrix{
 A_D \ar@<0.5ex>[r] \ar@<-0.5ex>[r] \ar[d]  & A_V \ar[d] \\
 H \ar@<0.5ex>[r] \ar@<-0.5ex>[r] & M .
}
\end{eqnarray*} 
\end{prop}
We consider various examples.
\begin{ex:}\label{1}
Lie groupoids are zero \LA -groupoids. Naturally, their global counter-parts are unit double Lie groupoids.
\end{ex:}
\begin{ex:}\label{2}
Lie algebroids are unit \LA -groupoids. Again, and only if it is the case that the given Lie algebroid is integrable, their global counter-parts are unit double Lie groupoids.
\end{ex:}
\begin{ex:}\label{3}
The tangent prolongation of a Lie groupoid is an \LA -groupoid. It is constructed in the obvious way by applying the tangent functor to the diagram defining the given Lie groupoid. The fact that the tangent functor respects products together with the fact that the groupoid axioms can be given strictly in terms of diagrams imply that the tangent prolongation $\xymatrix{TG \ar@<0.5ex>[r] \ar@<-0.5ex>[r] & TM}$ is indeed a Lie groupoid. On the other hand, the tangent bundle is canonically a Lie algebroid, as the algebra of vector fields comes endowed with the commutator bracket. Moreover, the differential of a map is always a Lie algebroid morphism; hence, the structural maps of the tangent prolongation verify the conditions of the definition of \LA -groupoid. We will discuss their global counter-parts in detail in a section below; however, notice that the pair groupoid $\xymatrix{G\times G \ar@<0.5ex>[r] \ar@<-0.5ex>[r] & M\times M}$ can be given the structure of a double Lie groupoid which clearly differentiates to the tangent prolongation.
\end{ex:}
\begin{ex:}\label{4}
Multiplicative foliations are anchor-injective \LA -groupoids. These are involutive subbundles of the tangent prolongation which are also subgroupoids. %As we will later point out, each diagram as \ref{LAGpd} comes with an anchor map to the tangent prolongation. This anchor map naturally is a map of \LA -groupoids, and its image is therefore a subgroupoid of $TH$, i.e. a multiplicative foliation. \\
As a subexample, consider a right principal $G$-bundle $P$ with surjective submersive moment map. Assume further, that the moment map has connected fibres. Let $\F_\mu$ be the simple foliation induced on $P$ by the moment map. One then sees that $\F_\mu$ inherits the action of $G$. More precisely, the derivative of the right multiplication induces a map $\xymatrix{T_p\F_\mu \ar[r] & T_{pg}\F_\mu}$. The action groupoid for this action is a multiplicative foliation over the action groupoid $P\rtimes G$. 
\end{ex:}
\begin{ex:}\label{5}
VB-algebroids are \LA -groupoids for which both the top and bottom groupoid structures are flat-abelian, i.e. vector bundles. The global counter-parts of these are the so-called VB-groupoids, which can be interperted as being double Lie groupoids for which the horizontal groupoids are vector bundles. As it was mentioned in the introduction and will later be spelled out, VB-groupoids and VB-algebroids are also related to the study of representations. Indeed, VB-groupoids and VB-algebroids are in correspondence with certain representations up to homotopy. 
\end{ex:}
\begin{ex:}\label{6}
Importantly for us, Lie $2$-algebras are \LA -groupoids over the groupoid with one object and one arrow. The global counter-parts of these are of course Lie $2$-groups, which are double Lie groupoids for which the bottom groupoid is $\xymatrix{\ast \ar@<0.5ex>[r] \ar@<-0.5ex>[r] & \ast}$ too. 
\end{ex:}
\subsubsection{The double complex associated to a double Lie groupoid}
Let $D$ be a double Lie groupoid. There are two simplicial structures for $D$ given by each of its vertical and its horizontal groupoid structures. It turns out that the commutativity of all of the square diagram representing a double Lie groupoid is but a shadow of the general interaction between these two simplicial structures. \\
For this piece, let $\vec{\gamma}=(\gamma_0,...,\gamma_q)\in (D^{(p+1)})^{(q+1)}$ be represented by the matrix $\begin{pmatrix}
    \gamma_{00} & \gamma_{01} & ... & \gamma_{0q} \\
    \gamma_{10} & \gamma_{11} & ... & \gamma_{1q} \\
    \vdots      & \vdots      &     & \vdots \\
    \gamma_{p0} & \gamma_{p1} & ... & \gamma_{pq} \\
\end{pmatrix}$, where each $\gamma_{mn}\in D$. Also, in order to recognize between the simplicial maps of the vertical and horizontal groupoids, we will write the usual $\partial_k$ for the simplicial maps of the horizontal groupoid and $\delta_j$ for the simplicial maps of the vertical groupoid.
\begin{lemma}\label{(j,k)=(0,0)}
$\delta_0\partial_0=\partial_0\delta_0$.
\end{lemma}
\begin{proof}
A simple computation yields
\begin{align*}
\delta_0\partial_0\vec{\gamma} & =\delta_0(\partial_0\gamma_0,...,\partial_0\gamma_q) \\ 
							   & =\delta_0\begin{pmatrix}
    \gamma_{10} & \gamma_{11} & ... & \gamma_{1q} \\
    \gamma_{20} & \gamma_{21} & ... & \gamma_{2q} \\
    \vdots      & \vdots      &     & \vdots      \\
    \gamma_{p0} & \gamma_{p1} & ... & \gamma_{pq} \\
\end{pmatrix} \\
							   & =\begin{pmatrix}
    \gamma_{11} & \gamma_{12} & ... & \gamma_{1q} \\
    \gamma_{21} & \gamma_{22} & ... & \gamma_{2q} \\
    \vdots      & \vdots      &     & \vdots      \\
    \gamma_{p1} & \gamma_{p1} & ... & \gamma_{pq} \\
\end{pmatrix} \\
							   & =\partial_0\begin{pmatrix}
    \gamma_{01} & \gamma_{02} & ... & \gamma_{0q} \\
    \gamma_{11} & \gamma_{12} & ... & \gamma_{1q} \\
    \vdots      & \vdots      &     & \vdots      \\
    \gamma_{p1} & \gamma_{p1} & ... & \gamma_{pq} \\
\end{pmatrix} = \partial_0\delta_0\vec{\gamma}.
\end{align*}
\end{proof}
In `minor' notation, $\delta_0\partial_0\vec{\gamma}=\partial_0\delta_0\vec{\gamma}=\vec{\gamma}_{0,0}$.
\begin{lemma}\label{j=0,k}
For $0<k\leq p$, $\delta_0\partial_k=\partial_k\delta_0$.
\end{lemma}
\begin{proof}
Computing,
\begin{align*}
\delta_0\partial_k\vec{\gamma} & =\delta_0(\partial_k\gamma_0,...,\partial_k\gamma_q) \\
							   & =\delta_0\begin{pmatrix}
    \gamma_{00}                   & \gamma_{01}                   & ... & \gamma_{0q} \\
    \vdots                        & \vdots                        &     & \vdots      \\
    \gamma_{k-20}                 & \gamma_{k-21}                 & ... & \gamma_{k-2q} \\
    \gamma_{k-10}\Join\gamma_{k0} & \gamma_{k-11}\Join\gamma_{k1} & ... & \gamma_{k-1q}\Join\gamma_{kq} \\
    \gamma_{k+10}                 & \gamma_{k+11}                 & ... & \gamma_{k+1q} \\
    \vdots                        & \vdots                        &     & \vdots      \\
    \gamma_{p0}                   & \gamma_{p1}                   & ... & \gamma_{pq} \\
\end{pmatrix} \\ 
							   & =\begin{pmatrix}
    \gamma_{10}                   & \gamma_{11}                   & ... & \gamma_{1q} \\
    \vdots                        & \vdots                        &     & \vdots      \\
    \gamma_{k-20}                 & \gamma_{k-21}                 & ... & \gamma_{k-2q} \\
    \gamma_{k-10}\Join\gamma_{k0} & \gamma_{k-11}\Join\gamma_{k1} & ... & \gamma_{k-1q}\Join\gamma_{kq} \\
    \gamma_{k+10}                 & \gamma_{k+11}                 & ... & \gamma_{k+1q} \\
    \vdots                        & \vdots                        &     & \vdots      \\
    \gamma_{p0}                   & \gamma_{p1}                   & ... & \gamma_{pq} \\
\end{pmatrix} \\  
							   & =\partial_k\begin{pmatrix}
    \gamma_{10}   & \gamma_{11}   & ... & \gamma_{1q} \\
    \vdots        & \vdots        &     & \vdots      \\
    \gamma_{k0}   & \gamma_{k-11} & ... & \gamma_{k-1q} \\
    \gamma_{k0}   & \gamma_{k1}   & ... & \gamma_{kq} \\
    \gamma_{k+10} & \gamma_{k+11} & ... & \gamma_{k+1q} \\
    \vdots        & \vdots        &     & \vdots      \\
    \gamma_{p0}   & \gamma_{p1}   & ... & \gamma_{pq} \\
\end{pmatrix} = \partial_k\delta_0\vec{\gamma}. 
\end{align*}
\end{proof}
Since there is a symmetry between the horizontal and the vertical groupoids, from the proof of this lemma follows that $\delta_j\partial_0=\partial_0\delta_j$ as well.
\begin{lemma}\label{(j,k)=(0,p+1)}
$\delta_0\partial_{p+1}=\partial_{p+1}\delta_0$.
\end{lemma}
\begin{proof}
Again, a simple computation yields
\begin{align*}
\delta_0\partial_{p+1}\vec{\gamma} & =\delta_0(\partial_p\gamma_0,...,\partial_p\gamma_q) \\ 
							       & =\delta_0\begin{pmatrix}
    \gamma_{00}   & \gamma_{01}   & ... & \gamma_{0q}   \\
    \vdots        & \vdots        &     & \vdots        \\
    \gamma_{p-20} & \gamma_{p-21} & ... & \gamma_{p-2q} \\
    \gamma_{p-10} & \gamma_{p-11} & ... & \gamma_{p-1q} \\
\end{pmatrix} \\
							       & =\begin{pmatrix}
    \gamma_{01}   & \gamma_{02}   & ... & \gamma_{0q}   \\
    \vdots        & \vdots        &     & \vdots        \\
    \gamma_{p-21} & \gamma_{p-22} & ... & \gamma_{p-2q} \\
    \gamma_{p-11} & \gamma_{p-12} & ... & \gamma_{p-1q} \\
\end{pmatrix} \\
							       & =\partial_{p+1}\begin{pmatrix}
    \gamma_{01}   & \gamma_{02}   & ... & \gamma_{0q}   \\
    \vdots        & \vdots        &     & \vdots        \\
    \gamma_{p-11} & \gamma_{p-12} & ... & \gamma_{p-1q} \\
    \gamma_{p1}   & \gamma_{p2}   & ... & \gamma_{pq}   \\
\end{pmatrix} = \partial_p\delta_0\vec{\gamma}.
\end{align*}
\end{proof}
Exploiting the aforementioned symmetry, we also have $\delta_{q+1}\partial_0=\partial_0\delta_{q+1}$.
\begin{lemma}\label{non-zero jk}
For $0<j\leq q$ and $0<k\leq p$, $\delta_j\partial_k=\partial_k\delta_j$.
\end{lemma}
\begin{proof}
Recall that
\begin{eqnarray*}
\partial_k\gamma=\partial_k(\gamma_0,...,\gamma_p)=(\gamma_0,...,\gamma_{k-2},\gamma_{k-1}\Join\gamma_k,\gamma_{k+1},...,\gamma_p),
\end{eqnarray*}
and that the multiplication in $D^{(p+1)}$ is the one inherited from $D^{p+1}$. Consequently,
\begin{align*}
\delta_j\partial_k\vec{\gamma} & =\delta_j(\partial_k\gamma_0,...,\partial_k\gamma_q) =\delta_j\begin{pmatrix}
    \gamma_{00}                   & \gamma_{01}                   & ... & \gamma_{0q} \\
    \vdots                        & \vdots                        &     & \vdots      \\
    \gamma_{k-20}                 & \gamma_{k-21}                 & ... & \gamma_{k-2q} \\
    \gamma_{k-10}\Join\gamma_{k0} & \gamma_{k-11}\Join\gamma_{k1} & ... & \gamma_{k-1q}\Join\gamma_{kq} \\
    \gamma_{k+10}                 & \gamma_{k+11}                 & ... & \gamma_{k+1q} \\
    \vdots                        & \vdots                        &     & \vdots      \\
    \gamma_{p0}                   & \gamma_{p1}                   & ... & \gamma_{pq} \\
\end{pmatrix}, 
\end{align*}
and
\begin{align*}
\partial_k\delta_j\vec{\gamma} & =\partial_k\begin{pmatrix}
    \gamma_{00} & ... & \gamma_{0j-2} & \gamma_{0j-1}\vJoin\gamma_{0j} & \gamma_{0j+1} & ... & \gamma_{0q} \\
    \gamma_{10} & ... & \gamma_{1j-2} & \gamma_{1j-1}\vJoin\gamma_{1j} & \gamma_{1j+1} & ... & \gamma_{1q} \\
    \vdots      &     & \vdots        & \vdots                         & \vdots        & ... & \vdots      \\
    \gamma_{p0} & ... & \gamma_{pj-2} & \gamma_{pj-1}\vJoin\gamma_{pj} & \gamma_{pj+1} & ... & \gamma_{pq} \\
\end{pmatrix}. 
\end{align*}
Both these yield the matrix
\begin{eqnarray*}
\begin{pmatrix}
    \gamma_{00}   & ... & \gamma_{0j-2}   & \gamma_{0j-1}\vJoin\gamma_{0j}     & \gamma_{0j+1}   & ... & \gamma_{0q}   \\
    \vdots        &     & \vdots          & \vdots                             & \vdots          &     & \vdots        \\
    \gamma_{k-20} & ... & \gamma_{k-2j-2} & \gamma_{k-2j-1}\vJoin\gamma_{k-2j} & \gamma_{k-2j+1} & ... & \gamma_{k-2q} \\
\gamma_{k-10}\Join\gamma_{k0} & ... & \gamma_{k-1j-2}\Join\gamma_{kj-2} & \diamondsuit & \gamma_{k-1j+1}\Join\gamma_{kj+1} & ... & \gamma_{k-1q}\Join\gamma_{kq} \\
    \gamma_{k+10} & ... & \gamma_{k+1j-2} & \gamma_{k+1j-1}\vJoin\gamma_{k+1j} & \gamma_{k+1j+1} & ... & \gamma_{k+1q} \\
    \vdots        &     & \vdots          & \vdots                             & \vdots          &     & \vdots        \\
    \gamma_{p0}   & ... & \gamma_{pj-2}   & \gamma_{pj-1}\vJoin\gamma_{pj}     & \gamma_{pj+1}   & ... & \gamma_{pq} \\
\end{pmatrix},
\end{eqnarray*}						    
where $\diamondsuit$ is interchangeably $(\gamma_{k-1j-1}\Join\gamma_{kj-1})\vJoin(\gamma_{k-1j}\Join\gamma_{kj})$ for the first one, and $(\gamma_{k-1j-1}\vJoin\gamma_{k-1j})\Join(\gamma_{kj-1}\vJoin\gamma_{kj})$ for the second.

\end{proof}
\begin{lemma}\label{j k=p+1}
For $0<j\leq q$, $\delta_j\partial_{p+1}=\partial_{p+1}\delta_j$.
\end{lemma}
\begin{proof}
Computing,
\begin{align*}
\delta_j\partial_{p+1}\vec{\gamma} & =\delta_j(\partial_{p+1}\gamma_0,...,\partial_{p+1}\gamma_q) \\
							       & =\delta_j\begin{pmatrix}
    \gamma_{00}   & \gamma_{01}   & ... & \gamma_{0q} \\
    \vdots        & \vdots        &     & \vdots      \\
    \gamma_{p-20} & \gamma_{p-21} & ... & \gamma_{p-2q} \\
    \gamma_{p-10} & \gamma_{p-11} & ... & \gamma_{p-1q} \\
\end{pmatrix} \\ 
							   & =\begin{pmatrix}
    \gamma_{00}   & ... & \gamma_{0j-2}   & \gamma_{0j-1}\vJoin\gamma_{0j}     & \gamma_{0j+1}   & ... & \gamma_{0q}   \\
    \vdots        &     & \vdots          & \vdots                             & \vdots          & ... & \vdots        \\
    \gamma_{p-20} & ... & \gamma_{p-2j-2} & \gamma_{p-2j-1}\vJoin\gamma_{p-2j} & \gamma_{p-2j+1} & ... & \gamma_{p-2q} \\
    \gamma_{p-10} & ... & \gamma_{p-1j-2} & \gamma_{p-1j-1}\vJoin\gamma_{p-1j} & \gamma_{p-1j+1} & ... & \gamma_{p-1q} \\
\end{pmatrix} \\  
							   & =\partial_{p+1}\begin{pmatrix}
    \gamma_{00}   & ... & \gamma_{0j-2}   & \gamma_{0j-1}\vJoin\gamma_{0j}     & \gamma_{0j+1}   & ... & \gamma_{0q}   \\
    \vdots        &     & \vdots          & \vdots                             & \vdots          & ... & \vdots        \\
    \gamma_{p-10} & ... & \gamma_{p-1j-2} & \gamma_{p-1j-1}\vJoin\gamma_{p-1j} & \gamma_{p-1j+1} & ... & \gamma_{p-1q} \\
    \gamma_{p0}   & ... & \gamma_{pj-2}   & \gamma_{pj-1}\vJoin\gamma_{pj}     & \gamma_{pj+1}   & ... & \gamma_{pq}   \\
\end{pmatrix} \\   = \partial_{p+1}\delta_j\vec{\gamma}. 
\end{align*}
\end{proof}
Yet again, the proof lemma also implies $\delta_{q+1}\partial_k=\partial_k\delta_{q+1}$. Finally,
\begin{lemma}\label{(j,k)=(q+1,p+1)}
$\delta_{q+1}\partial_{p+1}=\partial_{q+1}\delta_{p+1}$.
\end{lemma}
\begin{proof}
A simple computation yields
\begin{align*}
\delta_{q+1}\partial_{p+1}\vec{\gamma} & =\delta_{q+1}(\partial_{p+1}\gamma_0,...,\partial_{p+1}\gamma_q) \\ 
							           & =\delta_0\begin{pmatrix}
    \gamma_{00}   & \gamma_{01}   & ... & \gamma_{0q}   \\
    \vdots        & \vdots        &     & \vdots        \\
    \gamma_{p-20} & \gamma_{p-21} & ... & \gamma_{p-2q} \\
    \gamma_{p-10} & \gamma_{p-11} & ... & \gamma_{p-1q} \\
\end{pmatrix} \\ 
							   & =\begin{pmatrix}
    \gamma_{00}   & ... & \gamma_{0q-2}   & \gamma_{0q-1}   \\
    \vdots        &     & \vdots          & \vdots          \\
    \gamma_{p-20} & ... & \gamma_{p-2q-2} & \gamma_{p-2q-1} \\
    \gamma_{p-10} & ... & \gamma_{p-1q-2} & \gamma_{p-1q-1} \\
\end{pmatrix} \\ 
							   & =\partial_{p+1}\begin{pmatrix}
    \gamma_{00}   & ... & \gamma_{0q-2}   & \gamma_{0q-1}   \\
    \vdots        &     & \vdots          & \vdots          \\
    \gamma_{p-10} & ... & \gamma_{p-1q-2} & \gamma_{p-1q-1} \\
    \gamma_{p0}   & ... & \gamma_{pq-2}   & \gamma_{pq-1}   \\
\end{pmatrix} = \partial_{p+1}\delta_{q+1}\vec{\gamma}.
\end{align*}
\end{proof}
In `minor' notation, $\delta_{q+1}\partial_{p+1}\vec{\gamma}=\partial_{q+1}\delta_{p+1}\vec{\gamma}=\vec{\gamma}_{q,p}$. Notice that the proof of this latter lemma is superfluous as groupoids have yet another symmetry, namely the inverse yields an automorphism; thus, some of this lemmas follow from each other. In particular lemma \ref{(j,k)=(q+1,p+1)} follows from \ref{(j,k)=(0,0)}. \\
An immediate consequence of these lemmas is the fact that the complexes of groupoid cochains for the horizontal and vertical groupoids fit into a double complex. Although this is an expected relation, we could not find a reference in the literature.
\begin{prop}\label{DoubleGpdDoubleCx}
Given a double Lie groupoid,
\begin{eqnarray*}
\xymatrix{
\vdots                              & \vdots                                 & \vdots                             &       \\ 
C(V^{(2)}) \ar[r]^{\partial}\ar[u]  & C(D\times_H D) \ar[r]^{\partial}\ar[u] & C(D_2^2)\ar[r]\ar[u]               & \dots \\
C(V) \ar[r]^{\partial}\ar[u]^\delta & C(D) \ar[r]^{\partial}\ar[u]^\delta    & C(D\times_V D) \ar[r]\ar[u]^\delta & \dots \\
C(M) \ar[r]^{\partial}\ar[u]^\delta & C(H) \ar[r]^{\partial}\ar[u]^\delta    & C(H^{(2)}) \ar[r]\ar[u]^\delta     & \dots 
}
\end{eqnarray*} 
is a double complex, where 
\begin{align*}
D_p^q:=\lbrace\begin{pmatrix}
d_{11} & d_{12} & ... & d_{1p} \\
d_{21} & d_{22} & ... & d_{2p} \\
\vdots & \vdots & ... & \vdots \\
d_{q1} & d_{q2} & ... & d_{qp} \\
\end{pmatrix} & \in M_{q\times p}(D): \\
               & \Tp{s}(d_{mn})=\Tp{t}(d_{mn+1}),\quad\Lf{s}(d_{mn})=\Lf{t}(d_{m+1n}), \\
               & \Tp{t}(d_{mn})=\Tp{s}(d_{mn-1}),\quad\Lf{t}(d_{mn})=\Lf{s}(d_{m-1n})\rbrace .
\end{align*}
\end{prop}
\begin{proof}
For $\omega\in C(D_p^q)$ and $\vec{\gamma}\in D_{p+1}^{q+1}$,
\begin{align*}
\partial\delta\omega(\vec{\gamma}) & =\sum_{k=0}^{p+1}(-1)^k\delta\omega(\partial_k\vec{\gamma}) \\
								   & =\sum_{k=0}^{p+1}(-1)^k\sum_{j=0}^{q+1}(-1)^j\omega(\delta_j\partial_k\vec{\gamma}) \\
								   & =\sum_{j=0}^{q+1}(-1)^j\sum_{k=0}^{p+1}(-1)^k\omega(\partial_k\delta_j\vec{\gamma}) \\
								   & =\sum_{j=0}^{q+1}(-1)^j\partial\omega(\delta_j\vec{\gamma})=\delta\partial\omega(\vec{\gamma}).
\end{align*}
\end{proof}
In the sequel, we will refer to this object as the \textit{double complex associated to} $D$, and the cohomology of its total complex $(C _{tot}^\bullet (D),d)$,
\begin{eqnarray*}
C _{tot}^k (D)=\bigoplus _{p+q=k}C(D_p^q)\qquad d=\partial + (-1)^q\delta ,
\end{eqnarray*}
the \textit{double groupoid cohomology of} $D$. We will be interested in the double complex associated to Lie $2$-groups in the upcoming chapters.

\subsubsection{The double complex associated to an \LA -groupoid}
Let now $A$ be a Lie algebroid, and recall that the $\Omega^\bullet(A)$ came with a exterior differential associated to $A$ and a graded product that verified a graded Leibniz type rule. Such a structure put together $(\Omega^\bullet(A),\wedge,d_A)$ is usually referred to as a differential graded algebra or DG-algebra for short.\\

Now, suppose we are given a vector bundle $A$ over $M$ and a differential $\delta$ making $(\Omega^\bullet (A),\wedge ,\delta)$ into a DG-algebra. Then, we claim that there is a Lie algebroid structure on $A$ induced by $\delta$.  Indeed, define the anchor to be a map of modules $\xymatrix{\Gamma(A) \ar[r] & \XX(M)}$, given by
\begin{eqnarray*}
\rho_\delta (a)f:=\delta f(a),
\end{eqnarray*}
where $a\in\Gamma(A)$ and $f\in C^\infty(M)=\Omega^0(A)$. We remark that this is well defined precisely because $\delta$ is a derivation. To define the bracket, we use the isomorphism $A\cong (A^*)^*$. Let $\omega\in\Omega^1(A)$, and $a_1,a_2\in\Gamma(A)$. Set
\begin{eqnarray*}
\omega([a_1,a_2]_\delta):=\rho_\delta (a_1)\omega(a_2)-\rho_\delta (a_2)\omega(a_1)-\delta\omega(a_1,a_2).
\end{eqnarray*}
It should be clear that the definitions for the anchor and the bracket came from their appearances in the formulae for the differential of the Chevalley-Eilenberg complex, and therefore, if we call $A_\delta$ the Lie algebroid thus defined, $d_{A_\delta}=\delta$ and conversely $A_{d_A}=A$. We summarize this in the following proposition.
\begin{prop}(see e.g.\cite{Vaintrob})
Given a vector bundle $A$ over $M$, there is a one-to-one correspondence between Lie algebroid structures on $A$ and DG-algebra  structures on $(\Omega^\bullet(A),\wedge)$. 
\end{prop}
This correspondence extends to an isomorphism of categories. According to the author of \cite{Ayala}, first appeared in \cite{Vaintrob}.
\begin{prop}
Let 
\begin{eqnarray*}
\xymatrix{
A_1 \ar[d]\ar[r]^F & A_2 \ar[d] \\
M_1 \ar[r]^f       & M_2
}
\end{eqnarray*} 
be a vector bundle map between the Lie algebroids $A_1$ and $A_2$. Then, $F$ is a Lie algebroid morphism if, and only if $F^*$ is a morphism between the complexes $(\Omega^\bullet (A_2 ),d_{A_2 })$ and $(\Omega^\bullet (A_1 ),d_{A_1 })$.
\end{prop}
In fact, notice that due to the compatibility with the wedge product, in order to proof that $F^*$ is a map of complexes, it is enough to see that the pull-back commutes with the exterior derivatives in degrees $0$ and $1$. Recall that the pull-back is defined by the formula
\begin{eqnarray*}
(F^*\omega)_x(a^1,...,a^k)=\omega _{f(x)}(F(a^1_x),...,F(a^k_x)),
\end{eqnarray*}
where $\omega\in\Omega^k(A_2)$, $x\in M_1$ and $a^j\in\Gamma(A_1)$ for each $j$. It is not straightforward to see why this definition makes sense. Indeed, it is not obvious that $x\mapsto (F^*\omega)_x$ is smooth, but let us set aside this issue for now. Next, consider an \LA -groupoid $\A$ schematized by
\begin{eqnarray*}
\xymatrix{
 A_1 \ar@<0.5ex>[r] \ar@<-0.5ex>[r] \ar[d]  & A_0 \ar[d] \\
 G_1 \ar@<0.5ex>[r] \ar@<-0.5ex>[r] & G_0 .
}
\end{eqnarray*}
We will omit the bottom structural maps, as they are induced from the top ones together with the zero section. We write down the rest of the structural maps in the following fashion
\begin{eqnarray*}
\xymatrix{
A_2 \ar[r]^{\hat{m}} & A_1 \ar[r]^{\hat{\iota}} & A_1 \ar@<0.5ex>[r]^{\hat{s}} \ar@<-0.5ex>[r]_{\hat{t}} & A_0 \ar[r]^{\hat{u}} & A_1,
}
\end{eqnarray*}
where $A_2:=A_1 \times_{A_0} A_1$, and this fibred product is the usual space of composable arrows. Using the above correspondence we get a \textit{co-groupoid} object in the category of DG-algebras associated to Lie algebroids; namely,
\begin{eqnarray*}
\xymatrix{
\vdots & \vdots  & \vdots & \vdots  & \vdots \\
\Omega ^2(A_2) \ar[u]^{d_2} & \Omega ^2(A_1)\ar[l]_{\hat{m}^*}\ar[u]^{d_1}  & \Omega ^2(A_1) \ar[l]_{\hat{\iota}^*}\ar[u]^{d_1} & \Omega ^2(A_0) \ar@<-0.5ex>[l]_{\hat{s}^*} \ar@<0.5ex>[l]^{\hat{t}^*}\ar[u]^{d_0}  & \Omega ^2(A_1) \ar[l]_{\hat{u}^*}\ar[u]^{d_1} \\
\Omega ^1(A_2) \ar[u]^{d_2} & \Omega ^1(A_1)\ar[l]_{\hat{m}^*}\ar[u]^{d_1}  & \Omega ^1(A_1) \ar[l]_{\hat{\iota}^*}\ar[u]^{d_1} & \Omega ^1(A_0) \ar@<-0.5ex>[l]_{\hat{s}^*} \ar@<0.5ex>[l]^{\hat{t}^*}\ar[u]^{d_0}  & \Omega ^1(A_1) \ar[l]_{\hat{u}^*}\ar[u]^{d_1} \\
C^\infty(G_1^{(2)}) \ar[u]^{d_2}  & C^\infty(G_1) \ar[l]_{\hat{m}^*}\ar[u]^{d_1}  & C^\infty(G_1) \ar[l]_{\hat{\iota}^*}\ar[u]^{d_1} & C^\infty(G_0) \ar@<-0.5ex>[l]_{\hat{s}^*} \ar@<0.5ex>[l]^{\hat{t}^*}\ar[u]^{d_0}  & C^\infty(G_1) \ar[l]_{\hat{u}^*}\ar[u]^{d_1}, 
}
\end{eqnarray*}
where $d_i$ is a shorthand for $d_{A_i}$. Perhaps more importantly, the nerve of the groupoid is a simplicial manifold internal to the category of Lie algebroids, as the multiplication and the projections are Lie algebroid morphisms. Therefore, out of the nerve
\begin{eqnarray*}
\xymatrix{
... \ar@<1.5ex>[r]\ar@<0.5ex>[r]\ar@<-0.5ex>[r]\ar@<-1.5ex>[r] & A_2 \ar[r]\ar@<1.0ex>[r]\ar@<-1.0ex>[r] &  A_1 \ar@<0.5ex>[r]\ar@<-0.5ex>[r] & A_0,
}
\end{eqnarray*}
we get 
\begin{eqnarray*}
\xymatrix{
... & \ar@<1.5ex>[l]\ar@<0.5ex>[l]\ar@<-0.5ex>[l]\ar@<-1.5ex>[l] \Omega^\bullet(A_2) & \ar[l]\ar@<1.0ex>[l]\ar@<-1.0ex>[l]\Omega^\bullet (A_1) & \ar@<0.5ex>[l]\ar@<-0.5ex>[l] \Omega^\bullet (A_0),
}
\end{eqnarray*}
by means of the above equivalence. Using the usual combinatorics, one then builds the double complex
\begin{eqnarray*}
\xymatrix{
\vdots                                     & \vdots                                     & \vdots            \\
\Omega ^2(A_0) \ar[u]^{d_0}\ar[r]^\partial & \Omega ^2(A_1) \ar[u]^{d_1}\ar[r]^\partial & \Omega ^2(A_2) \ar[u]^{d_2}\ar[r]^\partial & ... \\
\Omega ^1(A_0) \ar[u]^{d_0}\ar[r]^\partial & \Omega ^1(A_1) \ar[u]^{d_1}\ar[r]^\partial & \Omega ^1(A_2) \ar[u]^{d_2}\ar[r]^\partial & ... \\
C^\infty(G_0) \ar[u]^{d_0}\ar[r]^\partial  & C^\infty(G_1) \ar[u]^{d_1}\ar[r]^\partial  & C^\infty(G_1^{(2)}) \ar[u]^{d_2}\ar[r]^\partial & ... 
}
\end{eqnarray*}
We call this complex the \LA -\textit{double complex of} $\A$, and the cohomology of its total complex $(\Omega _{tot}^\bullet (\A),d)$,
\begin{eqnarray*}
\Omega _{tot}^k (\A)=\bigoplus_{p+q=k}\Omega^q(A_p)\qquad d=d_p + (-1)^q\partial ,
\end{eqnarray*}
the \LA -\textit{cohomology of} $\A =\xymatrix{A_1 \ar@<0.5ex>[r] \ar@<-0.5ex>[r] & A_0}$. As we pointed out before, Lie $2$-algebras coincide with \LA -groupoids over the unit groupoid of a point; hence, there is an \LA -double complex for these. We will be using it in the next chapter. \\
We close this subsection by pointing out two examples of this construction that appear in the literature.
\begin{ex:} 
When $\A$ is the tangent prolongation of a Lie groupoid, the double complex is the famous Bott-Shulman complex.
\end{ex:}
\begin{ex:}
The main character in the proof of the main theorem of \cite{VanEstC}, Crainic's extension of van Est theorem to Lie groupoids, is the \LA -complex of the \LA -groupoid of the subexample \ref{4}, with $P$ taken to be the trivial bundle, i.e. $G$ itself. Notice that as part of the construction and as part of the theorem, $s$-fibres are assumed to be connected.
\end{ex:}

\section{A short excursus on homological algebra}
In the sequel, it will be useful to have van Est theorem written in terms of the mapping cone. Though we could not find the result relating the cohomology of the mapping cone to the induced isomorphisms in cohomology literally (see Proposition \ref{ConeCoh} below), this follows from standard techniques that can be found in \cite{Weibel}. We adapt these to the case of double complexes. \\
Given $\xymatrix{\Phi:A^\bullet \ar[r] & B^\bullet}$, a map of complexes, the mapping cone complex of $\Phi$ is defined by 
\begin{eqnarray*}
C(\Phi):=(A[1]\oplus B,d_\Phi ),
\end{eqnarray*}
where $d_\Phi=\begin{pmatrix}
-d_A & 0 \\
\Phi & d_B 
\end{pmatrix}$. This complex fits in the obvious exact sequence
\begin{eqnarray*}
\xymatrix{
0 \ar[r] & B \ar[r]^j & C(\Phi) \ar[r]^\pi & A[1] \ar[r] & 0,
}
\end{eqnarray*}
where $j(b):=(0,b)$ and $\pi(a,b):=a$. We check that these are indeed maps of complexes
\begin{align*}
d_\Phi(j(b)) & =d_\Phi(0,b)           & \pi(d_\Phi(a,b)) & =\pi(-d_A(a),\Phi(a)+d_B(b)) \\
			 & =(0,d_B(b))=j(d_B(b)), &                  & =-d_A(a)=d_{A[1]}(\pi(a,b)).
\end{align*}
We will write $H(\Phi)$ for the cohomology of the mapping cone of $\Phi$.
\begin{prop}\label{ConeCoh}
If $\xymatrix{\Phi:A^\bullet \ar[r] & B^\bullet}$ is a map of complexes, the following are equivalent:
\begin{itemize}
\item[i)] $H^n(\Phi)=(0)$ for $n\leq k$.
\item[ii)] The induced map in cohomology
\begin{eqnarray*}
\xymatrix{
\Phi^n :H^n(A) \ar[r] & H^n(B),
}
\end{eqnarray*}
is an isomorphism for $q\leq k$, and it is injective for $q=k+1$.
\end{itemize}
\end{prop}
\begin{proof}
As remarked before the proof, $C(\Phi)$ fits into an exact sequence which induces a long exact sequence in cohomology
\begin{eqnarray*}
\xymatrix{
0 \ar[r] & \cancelto{0}{H^{-1}(B)} \ar[r]^{j^{-1}} & H^{-1}(\Phi) \ar[r]^{\pi^{-1}} & H^{-1}(A[1]) \ar[r] & H^0(B) \ar[r]^{j^0} & H^0(\Phi) \ar[r]^{\pi^0} & ...  
}\\
\qquad\qquad\xymatrix{
... \ar[r]^{\pi^0\quad} & H^0(A[1]) \ar[r] & H^1(B) \ar[r]^{j^1} & H^1(\Phi) \ar[r]^{\pi^1} & H^1(A[1]) \ar[r] & ...
}
\end{eqnarray*}
Naturally, $H^n(A[1])=H^{n+1}(A)$. We will show that the connecting homomorphism is $\Phi^*$. Recall that for $a\in Z^n(A)$, the connecting homomorphism is defined as follows. First, consider $(a,b)\in C^{n-1}(\Phi)$ which goes to $a$ under $\pi$. Since
\begin{align*}
\pi^{n-1}(d_\Phi(a,b))=\pi^{n-1}(0,\Phi(a)+d_B(b))=0,
\end{align*}
$d_\Phi(a,b)\in\ker\pi^{n-1}$. Then, there exists a unique element $\beta\in B^q$ such that $j(\beta)=d_\Phi(a,b)$. Of course, 
\begin{eqnarray*}
\beta =\Phi(a)+d_B(b).
\end{eqnarray*}
Now, the image of the class of $a$ under the connecting homomorphism is the class of $\beta$; therefore,
\begin{eqnarray*}
\xymatrix{
[a] \ar@{|->}[r] & [\Phi(a)]=\Phi^{n}([a]),
}
\end{eqnarray*}
as claimed. With this at hand the proof is rather straightforward. \\
($i)\Longrightarrow ii)$) If $H^n(\Phi)=(0)$ for $n\leq k$, the long exact sequence looks like 
\begin{eqnarray*}
\xymatrix{
0 \ar[r] & H^0(A) \ar[r]^{\Phi^0} & H^0(B) \ar[r] & 0 \ar[r] & H^1(A) \ar[r]^{\Phi^1} & H^1(B) \ar[r] & 0 \ar[r] & ...  
}\\
\qquad\xymatrix{
... \ar[r] & 0 \ar[r] & H^{k+1}(A) \ar[r]^{\Phi^{k+1}} & H^{k+1}(B) \ar[r]^{j^{k+1}} & H^{k+1}(\Phi) \ar[r]^{\pi^{k+1}} & H^{k+2}(A) \ar[r] & ...
}
\end{eqnarray*}
from which the second statement clearly follows .\\
($ii)\Longrightarrow i)$) Conversely, let's assume that $\Phi$ induces an isomorphism in cohomology for all degrees less than or equal to $k$ and an injective map for $k+1$. First, from 
\begin{eqnarray*}
\xymatrix{
0 \ar[r] & H^{-1}(\Phi) \ar[r]^{\pi^{-1}} & H^0(A) \ar[r]^{\Phi^0} & H^0(B) ,
}
\end{eqnarray*}
one sees that $\Img(\pi^0)=(0)$, but $\pi^0$ is injective; thus, $ H^{-1}(\Phi)=(0)$. Inductively, $\ker j^n=H^n(B)$, then $j^n\equiv 0$; therefore, $\ker\pi^n=(0)$ and since $\Img(\pi^n)=\ker\Phi^{n+1}=(0)$, $\ker\pi^n=H^n(\Phi)$ as well. Notice that the last step for which this argument runs is precisely $n=k$, as we have no information on $\ker\Phi^{k+2}$.

\end{proof}
We use this proposition to re-phrase van Est theorem.
\begin{theorem}\label{ClassicVanEst}
Let $H$ be a Lie group with Lie algebra $\hh$. If $H$ is $k$-connected and $\Phi$ is the van Est map, then
\begin{eqnarray*}
H^n(\Phi)=(0),\quad\textnormal{for all degrees } n\leq k.
\end{eqnarray*}
\end{theorem}

It will be relevant to us that the constructions above admit a generalization for double complexes. First, given a map $\xymatrix{\Phi:A^{\bullet,\bullet} \ar[r] & B^{\bullet,\bullet}}$ of double complexes, there is a mapping cone double complex constructed by putting the mapping cone complex of $\Phi\rest{A^{p,\bullet}}$ in the $p$th column and using the obvious differential in the rows.
\begin{lemma}
If $\xymatrix{\Phi:A^{\bullet,\bullet} \ar[r] & B^{\bullet,\bullet}}$ is a map of double complexes, then
\begin{eqnarray*}
\xymatrix{
\vdots                                                         & \vdots                                                         & \vdots                                         &       \\ 
A^{02}\oplus B^{01} \ar[r]^{\partial_\Phi}\ar[u]               & A^{12}\oplus B^{11} \ar[r]^{\partial_\Phi}\ar[u]               & A^{22}\oplus B^{21}\ar[r]\ar[u]                & \dots \\
A^{01}\oplus B^{00} \ar[r]^{\partial_\Phi}\ar[u]^{\delta_\Phi} & A^{11}\oplus B^{10} \ar[r]^{\partial_\Phi}\ar[u]^{\delta_\Phi} & A^{21}\oplus B^{20} \ar[r]\ar[u]^{\delta_\Phi} & \dots \\
A^{00} \ar[r]^{\partial_\Phi}\ar[u]^{\delta_\Phi}              & A^{10} \ar[r]^{\partial_\Phi}\ar[u]^{\delta_\Phi}              & A^{20} \ar[r]\ar[u]^{\delta_\Phi}              & \dots 
}
\end{eqnarray*}
with
\begin{eqnarray*}
\delta_\Phi=\begin{pmatrix}
-\delta_A & 0        \\
\Phi      & \delta_B 
\end{pmatrix}\textnormal{ and} & \partial_\Phi=\begin{pmatrix}
\partial_A & 0          \\
 0         & \partial_B 
\end{pmatrix}
\end{eqnarray*}
is a double complex itself.
\end{lemma}
Here and elsewhere, we assume the convention that $\delta$ represents the vertical differential and $\partial$ represents the horizontal differential in a double complex, if it is not stated otherwise .
\begin{proof}
We prove that the differentials commute. Let $(a,b)\in A^{p,q+1}\oplus B^{p,q}$, then
\begin{align*}
\delta_\Phi\partial_\Phi(a,b) & =\delta_\Phi(\partial_A(a),\partial_B(b)) \\
							  & =(-\delta_A\partial_A(a),\Phi(\partial_A(a))+\delta_B\partial_B(b)) \\
							  & =(-\partial_A\delta_A(a),\partial_B(\Phi(a))+\partial_B\delta_B(b)) \\
							  & =\partial_\Phi(-\delta_A(a),\Phi(a)+\delta_B(b))=\partial_\Phi\delta_\Phi(a,b).
\end{align*}
\end{proof}
We will denote the double complex of this lemma $C^{\bullet,\bullet}(\Phi)$, will refer to it also as the mapping cone of $\Phi$ and will write $H_{tot}(\Phi)$ for its total cohomology. This double complex also fits into an exact sequence
\begin{eqnarray*}
\xymatrix{
0 \ar[r] & B \ar[r]^j & C(\Phi) \ar[r]^\pi & A[0,1] \ar[r] & 0,
}
\end{eqnarray*}
where $A[0,1]^{p,q}:=A^{p,q+1}$ with $\partial_{A[0,1]}=\partial_A$ and $\delta_{A[0,1]}=-\delta_A$. This time round, though the maps are defined by the same formulas, thus making the sequence obviously exact again, we need to check that they are indeed compatible with the horizontal differentials, so that they constitute maps of double complexes. Indeed,
\begin{align*}
\partial_\Phi(j(b)) & =\partial_\Phi(0,b)                  & \pi(\partial_\Phi(a,b)) & =\pi(\partial_A(a),\partial_B(b)) \\
			        & =(0,\partial_B(b))=j(\partial_B(b)), &                         & =\partial_A(a)=\partial_{A[0,1]}(\pi(a,b)).
\end{align*}
Now, a map $\xymatrix{\Phi:A^{\bullet,\bullet} \ar[r] & B^{\bullet,\bullet}}$ of double complexes induces a map $\xymatrix{\Phi_{tot}:A_{tot} \ar[r] & B_{tot}}$ of complexes in the obvious way; therefore, there is a long exact sequence in the total cohomologies associated to a short exact sequence of double complexes as well. Thus, copying the proof of proposition \ref{ConeCoh} word for word, we have got the following.
\begin{prop}\label{mapConeCoh}
If $\xymatrix{\Phi:A^{\bullet,\bullet} \ar[r] & B^{\bullet,\bullet}}$ is a map of double complexes, the following are equivalent:
\begin{itemize}
\item[i)] $H_{tot}^n(\Phi)=(0)$ for $n\leq k$.
\item[ii)] The induced map in cohomology
\begin{eqnarray*}
\xymatrix{
\Phi^n :H_{tot}^n(A) \ar[r] & H_{tot}^n(B)
}
\end{eqnarray*}
is an isomorphism for $n\leq k$, and it is injective for $n=k+1$.
\end{itemize}
\end{prop}

%*** Where should the material for Spectral Sequences be inserted? Should it be inserted at all? Is Evans and Graham's The relative Hochschild-Serre spectral sequence a good reference for the main theorem...? *** \\

We will take advantage of this detour to prove this rather easy, but technical lemma that will be much used. Its contents can be seen as a simple application of the larger machinery of spectral sequences (see also \cite{Weibel}). 
\begin{lemma}\label{BelowDiag}
Let $(E^{p,q}_r,d_r)$ be the spectral sequence of a double complex $C^{\bullet,\bullet}$ filtrated by columns. If there is a page for which $E^{p,q}_r$ is zero for all $(p,q)$ satisfying $p+q\leq k$, then 
\begin{eqnarray*}
H_{tot}^n(C)=(0)\quad\textnormal{for }n\leq k.
\end{eqnarray*}
\end{lemma}
\begin{proof}
Of course, since the successive pages of the spectral sequence are defined by the cohomology of the $d_r$'s, once one reaches the page in which $E^{p,q}_r$ is zero \textit{below} $k$, it is safe to conclude that 
\begin{eqnarray*}
E^{p,q}_\infty =(0),\quad\forall(p,q)\textnormal{ such that }p+q\leq k.
\end{eqnarray*}
By definition, since $E^{p,q}_r\Rightarrow H_{tot}^{p+q}(C)$, 
\begin{eqnarray*}
E^{p,q}_\infty\cong gr_p(H_{tot}^{p+q}(C)):=F_p(H_{tot}^{p+q}(C))/F_{p+1}(H_{tot}^{p+q}(C)).
\end{eqnarray*}
Recall that the filtration by columns is given by
\begin{eqnarray*}
F_{p}(C_{tot}^{n}):=\bigoplus_{r\geq p}C^{r,n-r}
\end{eqnarray*} 
In general, this filtration starts to vanish at $n+1$, $F_{n+1}(C_{tot}^{n})=(0)$; thus implying, 
\begin{align*}
F_{n+1}(Z_{tot}^{n}(C)) & :=F_{n+1}C_{tot}^{n}\cap Z_{tot}^{n}(C)=(0), \\
F_{n+1}(H_{tot}^{n}(C)) & :=F_{n+1}Z_{tot}^{n}(C)/F_{n+1}B_{tot}^{n}(C)=(0)
\end{align*}
and ultimately $gr_n(H_{tot}^{n}(C))=F_n(H_{tot}^{n}(C))$. Hence, if $n\leq k$,
\begin{eqnarray*}
F_n(H_{tot}^{n}(C))\cong E^{n,0}_\infty =(0).
\end{eqnarray*}
By induction on $r$, for $p=n-r$ and $q=r$,
\begin{eqnarray*}
gr_{p}(H_{tot}^{n}(C))=F_{p}(H_{tot}^{n}(C))/\cancelto{I.H.}{F_{p+1}(H_{tot}^{n}(C))}=F_{p}(H_{tot}^{n}(C))\cong E^{p,q}_\infty =(0).
\end{eqnarray*}
Eventually, at $r=n$, this boils down to 
\begin{eqnarray*}
gr_0(H_{tot}^{n}(C))=H_{tot}^{n}(C)/\cancel{F_{1}(H_{tot}^{n}(C))}=H_{tot}^{n}(C)\cong E^{0,n}_\infty =(0).
\end{eqnarray*}
\end{proof}     
%% ------------------------------------------------------------------------- %%
\chapter{Cohomology: the Lie 2-algebra theory}\label{fractionalchapter}
\chaptermark{Cohomology: the Lie $2$-algebra theory}
\section{Introduction}
%---------------------------------------------------------------------------------------------------------------------------------------------------
In this chapter we introduce the cohomology of Lie $2$-algebras. We start by looking at the total cohomology of the \LA -double complex of the Lie $2$-algebra to discover that it classifies certain type of extensions. With the goal of classifying extensions in sight, we define the notion of a representation of a Lie $2$-algebra. We move on to study extensions by splitting them and looking for conditions for naturally defined maps to build an extension back up. We use these ``cocycle equations'' to read a complex associated to the Lie $2$-algebras with values in a representation. Finally, we give an interpretation for the lower dimensions of the cohomology of this complex and verify that its second cohomology indeed classifies extensions as prescribed.
 
\section{2-cohomology with trivial coefficients}\label{AlgDcx}

Let $\xymatrix{\gg _1 \ar@<0.5ex>[r] \ar@<-0.5ex>[r] & \hh}$ be a Lie $2$-algebra with associated crossed module $\xymatrix{\gg \ar[r]^\mu & \hh}$. Notice that $\gg_p :=\gg _1^{(p)}$ is a Lie subalgebra of $\gg _1^p$ for each $p$, and consider the \LA -double complex of $\gg_1$
\begin{eqnarray*}
\xymatrix{
\vdots                                             & \vdots                                            & 
\vdots                                                &      \\ 
\bigwedge ^3\hh ^*\ar[r]^{\partial}\ar[u]          & \bigwedge ^3\gg _1^* \ar[r]^{\partial}\ar[u]        & 
\bigwedge ^3\gg _2^*\ar[r]\ar[u]                      & \dots\\
\hh ^*\wedge \hh ^* \ar[r]^{\partial}\ar[u]^\delta & \gg _1^*\wedge\gg _1^* \ar[r]^{\partial}\ar[u]^\delta & 
\gg _2^*\wedge \gg _2^* \ar[r]\ar[u]^\delta           & \dots\\
\hh ^* \ar[r]^{\partial}\ar[u]^\delta              & \gg _1^* \ar[r]^{\partial}\ar[u]^\delta             & 
\gg _2^* \ar[r]\ar[u]^\delta                          & \dots 
}
\end{eqnarray*} 
Here, the $p$th column is the Chevalley-Eilenberg complex of the Lie algebra $\gg_p$ whose differential simplifies from that of the Lie algebroid to
\begin{eqnarray*}
\delta\omega(\Xi)=\sum _{m<n}(-1)^{m+n}\omega([\xi_m ,\xi_n],\Xi(m,n));
\end{eqnarray*}
here, $\omega\in\bigwedge ^q\gg_p ^*$ and $\Xi=(\xi_0,...,\xi_q)\in\gg_p^q$. In order, the $q$th row is the subcomplex of alternating multilinear groupoid cochains of the groupoid $\xymatrix{\gg _1^q \ar@<0.5ex>[r] \ar@<-0.5ex>[r] & \hh^q}$. \\
We will refer to the \LA -cohomology of the $2$-algebra simply as $2$-cohomology. The total complex has the simplified form 
\begin{eqnarray*}
\Omega ^k_{tot}(\gg _1)=\bigoplus _{p+q=k}\bigwedge ^q\gg _p^*,
\end{eqnarray*}
with differential $d=\delta +(-1)^{q}\partial$.

We give an interpretation of $H^2_{tot}(\gg_1 )$. A $2$-cocycle consists of a pair of functions $(\omega,\varphi)\in (\hh^*\wedge\hh^*)\oplus \gg_1^*$ such that:\\
1) $\delta \omega = 0$, i.e. $-\omega([y_0,y_1] ,y_2)+\omega([y_0 ,y_2], y_1)-\omega([y_1, y_2],y_0)=0$ for all triples $y_0,y_1,y_2\in \hh$ \\
2) $\partial \varphi =0$, i.e. $\varphi(x_2 ,y)-\varphi(x _1 +x _2 ,y)+\varphi(x_1 ,y+\mu (x_2))=0$ for all $x_1,x_2\in\gg$ and $y\in\hh$. Due to linearity, this boils down to $\varphi(0,y+\mu(x_2))=0$ for all $y\in\hh$ and $x_2\in\gg$.\\
3) $\partial \omega+\delta \varphi=0$, i.e. $\omega(y_0,y_1)-\omega(y_0+\mu (x_0),y_1+\mu (x_1))= \varphi([(x_0,y_0),(x_1,y_1)]_\Lie)$ for all pairs $(x_0,y_0),(x_1,y_1)\in\gg\oplus\hh$. Again, due to bilinearity of $\omega$ and linearity of $\varphi$, this equation can be rewritten as
\begin{eqnarray*}
-\omega(y_0,\mu(x_1))-\omega(\mu(x_0),\mu(x_1))-\omega(\mu(x_0),y_1)=\varphi(\Lie_{y_0}x_1-\Lie_{y_1}x_0,0)+\varphi([x_0,x_1],[y_0,y_1]).
\end{eqnarray*}
Notice that this equation is equivalent to the simpler $\varphi(\Lie_{y}x,0)=-\omega(y,\mu(x))$. \\
It is rather known that if $\delta \omega=0$, $\omega$ induces an (central) extension of $\hh$, 
\begin{eqnarray*}
\xymatrix{
0 \ar[r] & \Rr \ar[r]^{\bar{0}\times I\quad} & \hh\oplus ^\omega \Rr  \ar[r]^{\quad pr_1} & \hh \ar[r] & 0,
}
\end{eqnarray*}
where $\hh\oplus ^\omega\Rr$ is the semi-direct sum induced by $\omega$, that is the space $\hh\oplus\Rr$ endowed with the twisted bracket
\begin{eqnarray*}
[(y_0,\lambda_0),(y_1,\lambda_1)]_\omega=([y_0,y_1],-\omega(y_0,y_1)).
\end{eqnarray*}
\begin{lemma}
If $d(\omega,\varphi)=0$, then 
\begin{eqnarray*}
\xymatrix{
\mu _\varphi: \gg \ar[r] & \hh\oplus ^\omega\Rr : x \ar@{|->}[r] & (\mu(x),\varphi(x,0))
}
\end{eqnarray*}
defines a crossed module for the action $\Lie_{(y,\lambda)}x:=\Lie_y x$.
\end{lemma}
\begin{proof}
First, let us verify that $\mu _\varphi$ is indeed a Lie algebra homomorphism. It is clearly linear and
\begin{align*}
\mu _\varphi([x_0,x_1]) & = (\mu([x_0,x_1]),\varphi([x_0,x_1],0))        \\
                        & = ([\mu(x_0),\mu(x_1)],-\omega(\mu(x_0),\mu(x_1))) \\
                        & = [(\mu(x_0),\varphi(x_0,0)),(\mu(x_1),\varphi(x_1,0))]_\omega ,
\end{align*}
here $\partial \omega((0,x_0),(0,x_1))=\delta \varphi((0,x_0),(0,x_1))$ justifies the second equality. Now, the action is still a Lie algebra homomorphism and is still by derivations because $\lambda$ does not change the action
\begin{eqnarray*}
\Lie_{[(y_0,\lambda_0),(y_1,\lambda_1)]}x=\Lie_{[y_0,y_1]}x=[\Lie_{y_0},\Lie_{y_1}]x =[\Lie_{(y_0,\lambda_0)},\Lie_{(y_1,\lambda_1)}]x, 
\end{eqnarray*}
\begin{eqnarray*}
\Lie_{(y,\lambda)}[x_1,x_2]=\Lie_y[x_1,x_2]=[\Lie_y x_1,x_2]+[x_1,\Lie_y x_2]=[\Lie_{(y,\lambda)} x_1,x_2]+[x_1,\Lie_{(y,\lambda)} x_2].
\end{eqnarray*}
As for the equivariance of $\mu_\varphi$, on the one hand we have got
\begin{eqnarray*}
\mu _\varphi (\Lie_{(y,\lambda )}x)=(\mu(\Lie_y x),\varphi(\Lie_y x,0)),
\end{eqnarray*}
while on the other,
\begin{align*}
[(y,\lambda ),\mu_\varphi(x)]_\omega & = [(y,\lambda ),(\mu(x),\varphi(x,0))]_\omega \\
                                     & = ([y,\mu(x)],-\omega(y,\mu(x))).                                     
\end{align*}
The first entry coincides because $\mu$ is a crossed module morphism, whereas, as we pointed out, $\partial\omega+\delta\varphi=0$ is equivalent to
\begin{eqnarray*}
\varphi(\Lie_{y}x,0)=-\omega(y,\mu(x)).
\end{eqnarray*}
Finally, one sees that the infinitesimal Peiffer equation holds as,
\begin{align*}
\Lie_{\mu_\varphi(x_2)}x_1 & =\Lie_{(\mu(x_2),\varphi(x_2,0))}x_1 \\
                           & =\Lie_{\mu(x_2)}x_1 =[x_2,x_1];
\end{align*} 
thus proving the lemma.

\end{proof}
As a consequence of this lemma, we see that there is an induced short exact sequence of Lie $2$-algebras that we write using their associated crossed modules
\begin{eqnarray*}
\xymatrix{
0 \ar[r] & 0 \ar[d]\ar[r] & \gg \ar[d]^{\mu_\varphi}\ar[r]^{Id} & \gg \ar[d]^\mu \ar[r] & 0  \\
0 \ar[r] & \Rr \ar[r]     & \hh\oplus^\omega\Rr \ar[r]^{\quad pr_1}   & \hh \ar[r]            & 0.  
}
\end{eqnarray*}

\begin{prop}
Let $(\omega ,\varphi ),(\omega ',\varphi ')\in\Omega^2_{tot}(\gg_1 )$. If there exists a $\phi\in\Omega^1_{tot}(\gg _1)=\hh^*$, such that $(\omega,\varphi)-(\omega ',\varphi ')=d\phi$, then the induced extensions are isomorphic.
\end{prop}
\begin{proof}
First, recall that the object extensions are isomorphic via
\begin{eqnarray*}
\begin{pmatrix}
    I     & 0 \\
    -\phi & 1 
\end{pmatrix}\xymatrix{ 
:\hh\oplus^\omega\Rr \ar[r] & \hh\oplus^{\omega '}\Rr:(y,\lambda ) \ar@{|->}[r] & (y,\lambda-\phi (y)).
}
\end{eqnarray*}
We claim that this together with the identity of $\gg$ induce the isomorphism between the extensions. Indeed, using the notation from the previous lemma
\begin{align*}
\begin{pmatrix}
    I     & 0 \\
    -\phi & 1 
\end{pmatrix} \mu _\varphi(x) & =\begin{pmatrix}
    I     & 0 \\
    -\phi & 1 
\end{pmatrix} (\mu(x),\varphi(x,0)) \\
                                & =(\mu(x),\varphi(x,0)-\phi(\mu(x))) \\
                                & =(\mu(x),\varphi '(x,0)) = \mu _{\varphi '}(x),
\end{align*}
also, trivially $Id(\Lie_{(y,\lambda )}x)=\Lie_{(y,\lambda-\phi(y))}Id(x)$, thus finishing the proof.

\end{proof}
The fact that the second $2$-cohomology group classifies a certain type of extensions is a happy accident, because it suggests an extension of the classical theory of Lie algebra cohomology. Since in the classical theory extensions are classified by a cohomology with coefficients in a representation, we will devote the following section to develop an appropriate notion of representation for a Lie $2$-algebra.
%points in the direction of one of the items of the list in \ref{vanEstIngr}... the main character of Cartan's proof of Lie III using the van Est strategy

\section{Representations of Lie 2-algebras}
In order to define a representation, we are going to proceed along the lines of the general philosophy: ``A representation of an object in a category $\mathcal{C}$ is a morphism to the space of endomorphisms of a flat abelian object in $\mathcal{C}$''.
Although this approach should be reasonable with no further justification, let us write a couple of examples.
\begin{itemize}
\item Vector spaces are flat abelian Lie groups, and representations of Lie groups are indeed maps to $GL(V)$.
\item Vector spaces are also abelian Lie algebras, and representations of Lie algebras are maps to $\ggl (V)$.
\item Vector bundles are abelian Lie groupoids whose $s$-fibres are flat, and a representation of a Lie groupoid is a map to $\xymatrix{GL(E) \ar@<0.5ex>[r] \ar@<-0.5ex>[r] & M}$.
\end{itemize}
Notice that such spaces of endomorphisms might as well not be objects in the given category; therefore, to make sense of a representation when adopting this philosophy, one needs to ensure that the spaces of endomorphisms are indeed objects of $\mathcal{C}$.
\subsection{The linear Lie 2-algebra associated to a 2-vector space}
The structure of the linear Lie $2$-algebra has already appeared in \cite{IntSubLin2} and \cite{Cristian&Olivier}. In the former reference, the linear Lie $2$-algebra is a special case of the DGLA of endomorphisms of a complex of vector spaces, and in the latter, it coincides with the gauge 2-groupoid of a complex of vector spaces over the point manifold. We introduce this structure from a different perspective. \\
A $2$-vector space is the same as an abelian Lie $2$-algebra, that is a groupoid object in the category of vector spaces. Using the equivalence between Lie $2$-algebras and crossed modules, one realizes that a $2$-vector space is equivalent to a $2$-term complex of vector spaces
\begin{eqnarray*}
\xymatrix{
W \ar[r]^\phi & V.
}
\end{eqnarray*}
By studying the algebraic structure on the category of linear self functors of $\xymatrix{W\oplus V \ar@<0.5ex>[r] \ar@<-0.5ex>[r] & V}$ and linear natural transformations among them, one gets the linear Lie $2$-algebra, which we call $\ggl (\phi)$.\\
First, linear functors correspond to pairs of linear maps $(F,f)$ commuting with $\phi$, i.e.
\begin{eqnarray*}
\xymatrix{
W \ar[r]^F \ar[d]_{\phi} & W \ar[d]^{\phi} \\
V \ar[r]_f             & V.
}
\end{eqnarray*}
A natural transformation $\alpha$ between two functors $(F_1,f_1)$ and $(F_2,f_2)$ is a linear map 
\begin{eqnarray*}
\xymatrix{
\alpha : V \ar[r] & W\oplus V : v \ar@{|->}[r] & (\alpha_1(v),\alpha_2(v)),
}
\end{eqnarray*}
such that the $s(\alpha_1(v),\alpha_2(v))=\alpha_2(v)=f_1(v)$ and $t(\alpha_1(v),\alpha_2(v))=\alpha_2(v)+\phi(\alpha_1(v))=f_2(v)$ and such that for every $(w,v)\in W\oplus V$ the diagram
\begin{eqnarray*}
\xymatrix{
f_1(v) \ar[rrr]^{(F_1(w),f_1(v)) \qquad} \ar[d]_{(\alpha_1(v),\alpha_2(v))} & & & f_1(v)+\phi(F_1(w)) \ar[d]^{(\alpha_1(v+\phi(w)),\alpha_2(v+\phi(w)))} \\
f_2(v) \ar[rrr]_{(F_2(w),f_2(v)) \qquad}                                    & & & f_2(v)+\phi(F_2(w)).
}
\end{eqnarray*}
commutes inside the category $W\oplus V$. That is
\begin{eqnarray*}
\alpha_1(v)+F_2(w)=\alpha_1(v+\phi(w))+F_1(w);
\end{eqnarray*}
which, by linearity of $\alpha$, is
\begin{eqnarray*}
\alpha_1(\phi(w))=F_2(w)-F_1(w).
\end{eqnarray*}
Summing up, the space of objects of the category of endomorphisms of $\xymatrix{W \ar[r]^\phi & V}$, $\ggl (\phi )$, is
\begin{eqnarray*}
\ggl (\phi )_0 =\lbrace (F,f) \in End(W)\oplus End(V):\phi\circ F=f\circ\phi \rbrace 
\end{eqnarray*} 
which is a vector space. Arrows between $(F_1,f_1)$ and $(F_2,f_2)$ are given by  
\begin{eqnarray*}
\ggl (\phi )((F_1,f_1);(F_2,f_2))=\lbrace A \in Hom(V,W):\phi\circ A=f_2 -f_1 ,\quad A\circ\phi =F_2 -F_1 \rbrace ;
\end{eqnarray*} 
therefore, 
\begin{eqnarray*}
\ggl (\phi ) = Hom(V,W)\oplus \ggl (\phi )_0 .
\end{eqnarray*} 
The structural maps are:
\begin{align*}
s(A;F,f) & =(F,f)   & t(A;F,f)     & =(F+A\phi ,f+\phi A)\\
u(F,f)   & =(0;F,f) & \iota(A;F,f) & =(-A;F+A\phi ,f+\phi A)
\end{align*}
\begin{eqnarray*}
(A';F+A\phi ,f+\phi A)\ast_v (A;F,f):=(A+A';F,f).
\end{eqnarray*}
Naturally, the multiplication is given by the so-called vertical composition of natural transformations. Notice that, incidentally, we found out that the category of endomorphisms of a $2$-vector space is not only a category, but a $2$-vector space itself. In the sequel, we will be referring to its associated $2$-term complex, 
\begin{eqnarray*}
\xymatrix{
Hom(V,W) \ar[r]^{\quad\Delta} & \ggl (\phi)_0,
}
\end{eqnarray*} 
where $\Delta A=(A\phi ,\phi A)$.\\
As we pointed out, there is an additional algebraic structure on $\ggl (\phi)$ that will endow it with a Lie $2$-algebra structure. Indeed, there is a bracket on the space of objects inherited from $\ggl (W)\oplus\ggl (V)$. For $(F_1,f_1),(F_2,f_2)\in \ggl (\phi )_0$, a simple computation yields $([F_1,F_2],[f_1,f_2])\in \ggl (\phi )_0$. Indeed,
\begin{align*}
\phi\circ [F_1,F_2] & = \phi(F_1F_2 -F_2F_1)     \\
                    & = \phi F_1F_2 -\phi F_2F_1 \\ 
                    & = f_1\phi F_2 -f_2\phi F_1 \\
                    & = f_1f_2\phi -f_2f_1\phi = [f_1,f_2]\circ\phi ,
\end{align*}
thus proving $\ggl (\phi)_0$ is a Lie subalgebra.\\
On the other hand, looking at the horizontal composition of natural transformations, one gets 
\begin{align*}
((A;F,f)\ast_h(B;G,g))_v & = (B(fv+\phi Av),g(fv+\phi Av))\Join (G(Av),g(fv)) \\
                         & = ((Bf+B\phi A+GA)v,gf(v)),
\end{align*}
where $(A;F,f),(B;G,g)\in \ggl (\phi )$. This product is bi-linear and associative; therefore, it can be used to define a Lie bracket on $\ggl (\phi )$, 
\begin{eqnarray*}
[(A_1;F_1,f_1),(A_2;F_2,f_2)]:=(A_2;F_2,f_2)\ast_h(A_1;F_1,f_1)-(A_1;F_1,f_1)\ast_h(A_2;F_2,f_2).
\end{eqnarray*}
We claim that along with this structure, $\ggl (\phi )$ turns into a Lie $2$-algebra. We will prove this fact by looking at the associated $2$-term complex, which will thus inherit the structure of a crossed module of Lie algebras, where the bracket on $Hom(V,W)$ is defined by
\begin{eqnarray*}
[A_1,A_2]_\phi = A_1\phi A_2-A_2\phi A_1.
\end{eqnarray*}
When endowed with this additional structure, we will write $\ggl(\phi)_1$ instead of $Hom(V,W)$ to emphasize the presence of the bracket. Finally, the action of $(F,f)\in \ggl (\phi )_0$ on $A\in\ggl(\phi)_1$ is given by
\begin{eqnarray*}
\Lie^{\phi}_{(F,f)}A = FA-Af .
\end{eqnarray*}
\begin{prop}
Along with this bracket and this action,
\begin{eqnarray*}
\xymatrix{
\ggl(\phi)_1 \ar[r]^{\Delta} & \ggl(\phi)_0
}
\end{eqnarray*} 
is a crossed module of Lie algebras.
\end{prop}
\begin{proof}
This amounts to a routine check.
\begin{itemize}
\item $\Delta$ is a Lie algebra homomorphism: It is clearly linear and
\begin{align*}
\Delta [A_1,A_2]_\phi & = ([A_1,A_2]_\phi\phi ,\phi [A_1,A_2]_\phi) \\
                 & = ((A_1\phi A_2-A_2\phi A_1)\phi ,\phi (A_1\phi A_2-A_2\phi A_1)) \\
                 & = (A_1\phi A_2\phi -A_2\phi A_1\phi ,\phi A_1\phi A_2-\phi A_2\phi A_1) \\
                 & = ([A_1\phi, A_2\phi ],[\phi A_1\phi A_2])=[\Delta A_1,\Delta A_2].
\end{align*}
\item $\ggl(\phi)_0$ acts by derivations:
\begin{align*}
\Lie^\phi_{(F,f)}[A_1,A_2]_\phi & = F[A_1,A_2]_\phi-[A_1,A_2]_\phi f \\
                      & = F(A_1\phi A_2-A_2\phi A_1)-(A_1\phi A_2-A_2\phi A_1)f \\
                      & = FA_1\phi A_2-FA_2\phi A_1-A_1\phi A_2f+A_2\phi A_1f ,
\end{align*}
while on the other hand
\begin{align*}
[\Lie^\phi_{(F,f)}A_1,A_2]_\phi & = (\Lie^\phi_{(F,f)}A_1)\phi A_2 -A_2\phi\Lie^\phi_{(F,f)}A_1\\
                      & = FA_1\phi A_2-A_1f\phi A_2-A_2\phi FA_1+A_2\phi A_1f , 
\end{align*}
and
\begin{align*}
[A_1,\Lie^\phi_{(F,f)}A_2] & = A_1\phi\Lie^\phi_{(F,f)} A_2 -(\Lie^\phi_{(F,f)}A_2)\phi A_1\\
                      & = A_1\phi FA_2-A_1\phi A_2f-FA_2\phi A_1+A_2f\phi A_1 . 
\end{align*}
Hence, the desired equality follows by $\phi F=f\phi$.
\item $\xymatrix{\Lie^\phi :\ggl(\phi)_0 \ar[r] & \ggl (Hom(V,W))}$ is a Lie algebra homomorphism:
\begin{align*}
\Lie^\phi_{([F_1,F_2],[f_1,f_2])}A & = [F_1,F_2]A-A[f_1,f_2] \\
                              & = F_1F_2A-F_2F_1A-Af_1f_2+Af_2f_1 ,
\end{align*}
while on the other hand
\begin{align*}
\Lie^\phi_{(F_1,f_1)}\Lie^\phi_{(F_2,f_2)}A & = \Lie^\phi_{(F_1,f_1)}(F_2A-Af_2) \\
                                  & = F_1F_2A-F_1Af_2 -F_2Af_1+Af_2f_1 , 
\end{align*}
and
\begin{align*}
\Lie^\phi_{(F_2,f_2)}\Lie^\phi_{(F_1,f_1)}A & = \Lie^\phi_{(F_2,f_2)}(F_1A-Af_1) \\
                                  & = F_2F_1A-F_2Af_1 -F_1Af_2+Af_1f_2 . 
\end{align*}
Therefore, subtracting gives the desired identity.
\item $\Delta$ is equivariant:
\begin{align*}
\Delta(\Lie^\phi_{(F,f)}A) & = ((\Lie^\phi_{(F,f)}A)\phi ,\phi\Lie^\phi_{(F,f)}A) \\
                       & = (FA\phi -Af\phi ,\phi FA-\phi Af) \\
                       & = (FA\phi -A\phi F,f\phi A-\phi Af) \\
                       & = ([F,A\phi ],[f,\phi A])=[(F,f),\Delta A]
\end{align*}
\item Infinitesimal Peiffer: 
\begin{align*}
\Lie^\phi_{\Delta A_0}A_1 & = \Lie^\phi_{(A_0\phi ,\phi A_0)}A_1 \\
                     & = A_0\phi A_1-A_1\phi A_0 =[A_0,A_1]_\phi.
\end{align*}
\end{itemize}
\end{proof}

\subsection{2-representations}
\begin{Def}
A \textit{representation} of a Lie $2$-algebra $\gg_1 =\xymatrix{\gg \ar[r]^\mu & \hh}$ on $\mathbb{V} =\xymatrix{W \ar[r]^\phi & V}$ is a morphism of Lie $2$-algebras
\begin{eqnarray*}
\xymatrix{
\rho :\gg_1 \ar[r] & \ggl (\phi ).
}
\end{eqnarray*}
\end{Def}
We will refer to a map as the one above as a $2$-representation. We would like to remark that, when we say a morphism of Lie $2$-algebras, we mean a na\"ive linear functor respecting the Lie algebra structures, as opposed to more general types of morphisms (e.g. bibundles of Lie groupoids, or equivalently maps of $2$-term $L_\infty$-algebras). \\
By definition, a representation of a Lie $2$-algebra has two ``commuting'' representations of $\hh$ on $W$ and on $V$ coming from the map of Lie algebras at the level of objects,  
\begin{eqnarray*}
\xymatrix{
\rho_0 :\hh \ar[r] & \ggl (\phi )_0\leq\ggl (W)\oplus\ggl (V). 
}
\end{eqnarray*}
Its components $\xymatrix{\rho_0^1 :\hh \ar[r] & \ggl (W)}$, $\xymatrix{\rho_0^0 :\hh \ar[r] & \ggl (V)}$ are representations fitting in the diagram
\begin{eqnarray*}
\xymatrix{
W \ar[d]_\phi \ar[r]^{\rho_0^1(y)} & W \ar[d]^\phi \\
V \ar[r]_{\rho_0^0(y)}             & V,
}
\end{eqnarray*} 
for each $y\in\hh$. On the other hand, at the level of arrows, one has got the map 
\begin{eqnarray*}
\xymatrix{
\rho_1 :\gg \ar[r] & \ggl(\phi)_1.
}
\end{eqnarray*}
This map is a Lie algebra homomorphism; hence,
\begin{eqnarray*}
\rho_1([x_0,x_1])=\rho_1(x_0)\phi\rho_1(x_1)-\rho_1(x_1)\phi\rho_1(x_0),
\end{eqnarray*}
and it also verifies the following equations for all $x\in\gg$: 
\begin{eqnarray*}
\rho^0_0(\mu(x))=\phi\rho_1(x), & \rho^1_0(\mu(x))=\rho_1(x)\phi ,
\end{eqnarray*}
and for all $y\in\hh$
\begin{eqnarray*}
\rho_1(\Lie_y x)=\rho_0^1(y)\rho_1(x)-\rho_1(x)\rho_0^0(y). 
\end{eqnarray*}
We now build an honest representation of $\gg _1 =\gg\oplus_\Lie\hh$ on $W\oplus V$.
\begin{prop}\label{honestAlgRep}
Given a representation $2$-representation $\xymatrix{\rho :\gg_1 \ar[r] & \ggl (\phi )}$, there is an honest representation 
\begin{eqnarray*}
\xymatrix{
\bar{\rho}:\gg\oplus_\Lie\hh \ar[r] & \ggl (W\oplus V):(x,y) \ar@{|->}[r] & {}
}
\begin{pmatrix}
    \rho_0^1(y+\mu(x)) & \rho_1(x) \\
    0                  & \rho_0^0(y) 
\end{pmatrix}
\end{eqnarray*}
\end{prop}
\begin{proof}
$\bar{\rho}$ is clearly linear and
\begin{align*}
\bar{\rho}([(x_0,y_0),(x_1,y_1)]_\Lie) & = \bar{\rho}([x_0,x_1]+\Lie _{y_0}x_1-\Lie_{y_1}x_0,[y_0,y_1]) \\
                                       & = \begin{pmatrix}
             \rho_0^1([y_0,y_1]+\mu([x_0,x_1]+\Lie _{y_0}x_1-\Lie_{y_1}x_0)) & \rho_1([x_0,x_1]+\Lie _{y_0}x_1-\Lie_{y_1}x_0) \\
                                             0                               & \rho_0^0([y_0,y_1]) 
                                           \end{pmatrix}.
\end{align*}
The maps appearing in the diagonal of this matrix will agree right away with those in the diagonal of $[\bar{\rho}(x_0,y_0),\bar{\rho}(x_1,y_1)]$. The map in the top corner to the right will coincide as well, as is shown by the following computation:
\begin{align*}
\rho_1([x_0,x_1]+\Lie _{y_0}x_1-\Lie_{y_1}x_0) & = \rho_1([x_0,x_1])+\rho_1(\Lie _{y_0}x_1)-\rho_1(\Lie_{y_1}x_0) \\
                                               & = [\rho_1(x_0),\rho_1(x_1)]+\Lie _{\rho_0(y_0)}\rho_1(x_1)-\Lie_{\rho_0(y_1)}\rho_1(x_0), 
\end{align*}
but the bracket is
\begin{align*}
[\rho_1(x_0),\rho_1(x_1)] & = \rho_1(x_0)\phi\rho_1(x_1)-\rho_1(x_1)\phi\rho_1(x_0) \\
                          & = \rho_0^1(\mu(x_0))\rho_1(x_1)-\rho_0^1(\mu(x_1))\rho_1(x_0);
\end{align*}
therefore, since $\Delta\circ\rho_1=\rho_0\circ\mu$ and the action verifies
\begin{align*}
\Lie _{\rho_0(y_0)}\rho_1(x_1) & = \rho_0^1(y_0)\rho_1(x_1)-\rho_1(x_1)\rho_0^0(y_0),
\end{align*}
we have got
\begin{align*}
\rho_1([x_0,x_1]+\Lie _{y_0}x_1-\Lie_{y_1}x_0) & = \rho_0^1(y_0+\mu(x_0))\rho_1(x_1)-\rho_0^1(y_1+\mu(x_1))\rho_1(x_0)+ \\
                                               & \qquad\rho_1(x_1)\rho_0^0(y_0)-\rho_1(x_0)\rho_0^0(y_1),                                                  
\end{align*}
which is the entry in the top corner of $[\bar{\rho}(x_0,y_0),\bar{\rho}(x_1,y_1)]$.

\end{proof}
The semi-direct sum of Lie algebras with respect to this representation,  $\gg _1 {}_{\bar{\rho}}\ltimes\mathbb{V}$, can be endowed with a structure that we call semi-direct product of Lie $2$-algebras and whose associated crossed module is
\begin{eqnarray*}
\xymatrix{
\gg {}_{\rho_0^1\circ\mu}\oplus W \ar[r]^{\mu\times\phi} & \hh {}_{\rho_0^0}\oplus V,
}
\end{eqnarray*}
with action
\begin{eqnarray*}
\Lie^\rho_{(y,v)}(x,w)=(\Lie_y x,\rho_0^1(y)w-\rho_1(x)v).
\end{eqnarray*}
We now provide some examples.
\begin{ex:}
Trivial representations. Of course, all of the defining relations get trivially satisfied if $\rho\equiv 0$.
\end{ex:}
\begin{ex:}\label{unitRep}
Usual Lie algebra representations can also be seen as examples of $2$-representations. Setting $W=\lbrace 0\rbrace$, a $2$-representation is equivalent to a single representation $\rho$ of $\hh$ on $V$, such that $\rho_{\mu(x)}\equiv 0$ for every $x\in\gg$. Ultimately, then, a $2$-representation on the unit $2$-vector space is simply a representation of $\hh/\mu(\gg)$. In particular, $2$-representations of a unit Lie $2$-algebra on a unit $2$-vector space are usual Lie algebra representations. \\
Curiously enough, assuming $V=\lbrace 0\rbrace$, one gets the same prescription, i.e. $2$-representation on a Lie group internal to vector spaces consist also of a single representation of the orbit space.
\end{ex:}
\begin{ex:}
Representations up to homotopy. The data defining a $2$-representation of the unit Lie $2$-algebra $\xymatrix{\hh \ar@<0.5ex>[r] \ar@<-0.5ex>[r] & \hh}$ coincides with a representation up to homotopy with zero curvature. In general, notice that the semi-direct product of Lie $2$-algebras $\gg_1{}_{\bar{\rho}}\ltimes\mathbb{V}$ is also VB-groupoid over $\gg_1$; hence, there is an associated representation up to homotopy (cf. appendix \ref{appchapter}). However, this representation is oblivious of the fact that it comes from a $2$-representation. In other words, for any given $2$-representation, the associated representation up to homotopy is the same.
\end{ex:}
\begin{ex:}
The adjoint representation (!). We describe the maps defining the adjoint representation of a Lie $2$-algebra on itself.
\begin{eqnarray*}
\xymatrix{ad_1:\gg \ar[r] & \ggl(\mu)_1} & \ad_1(x)(u):=-\Lie_u x ,  \\
\xymatrix{ad_0^1:\hh \ar[r] & \ggl(\gg)}   & \ad_0^1(y)(v):=\Lie_y v , \\
\xymatrix{ad_0^0:\hh \ar[r] & \ggl(\hh)}   & \ad_0^0(y)(u):=[y,u].    
\end{eqnarray*}
\end{ex:}

We close this section by remarking that in contrast with the general groupoid case, the latter example above shows that Lie $2$-algebras admit an adjoint representation extending the classic one. We look for a cohomology theory of Lie $2$-algebras with values in these $2$-representations. Since we prescribed that the second cohomology group classify extensions as a requisite, we are prompted to begin the next section.

\section{Extensions of Lie 2-algebras or towards a 2-cohomology with coefficients}
We study abstract extensions to try and make an educated guess as to what are fine candidates for the spaces of $1$-cochains and $2$-cochains and the corresponding subspaces of $2$-cocycles and $2$-coboundaries. We look to express extensions in terms of simpler data and decide the equations these need to satisfy. \\ 
We start this section by proving that the notion of $2$-representations defined in the previous sections is the right one in the sense that it is the kind of action induced in an abstract extension.
\begin{Def}
An \textit{extension} of the Lie $2$-algebra $\xymatrix{\gg \ar[r]^{\mu} & \hh}$ by the $2$-vector space $\xymatrix{W \ar[r]^{\phi} & V}$ is a Lie $2$-algebra $\xymatrix{\mathfrak{e}_1 \ar[r]^{\epsilon} & \mathfrak{e}_0}$ that fits in 
\begin{eqnarray*}
\xymatrix{
0 \ar[r] & W \ar[d]_\phi\ar[r]^{j_1} & \mathfrak{e}_1 \ar[d]_\epsilon\ar[r]^{\pi_1} & \gg \ar[d]^\mu\ar[r] & 0  \\
0 \ar[r] & V \ar[r]_{j_0}            & \mathfrak{e}_0 \ar[r]_{\pi_0}                & \hh \ar[r]           & 0, 
}
\end{eqnarray*}
where the top and bottom rows are short exact sequences and the squares are maps of Lie $2$-algebras.
\end{Def}
\begin{prop}\label{Ind2Rep} %Ind stands for induced
Given a Lie $2$-algebra extension of $\xymatrix{\gg \ar[r]^{\mu} & \hh}$ by a $2$-vector space $\xymatrix{W \ar[r]^{\phi} & V}$,
\begin{eqnarray*}
\xymatrix{
0 \ar[r] & W \ar[d]_\phi\ar[r]^{j_1} & \mathfrak{e}_1 \ar[d]_\epsilon\ar[r]_{\pi_1} & \gg \ar[d]^\mu\ar[r]\ar@/_/[l]_{\sigma_1} & 0  \\
0 \ar[r] & V \ar[r]^{j_0}            & \mathfrak{e}_0 \ar[r]_{\pi_0}                & \hh \ar[r]\ar@/_/[l]_{\sigma_1}      & 0, 
}
\end{eqnarray*}
and a linear splitting $\sigma$, there is an induced $2$-representation $\xymatrix{\rho^{\epsilon}_{\sigma}:\gg_1 \ar[r] & \ggl (\phi)}$ given by
\begin{align*}
\rho^0_0(y)v     & =[\sigma_0(y),v]_{\mathfrak{e}_0}              & \rho^1_0(y)w & =\Lie^{\epsilon}_{\sigma_0(y)} w \\
\rho_1(x)v       & =-\Lie^{\epsilon}_v\sigma_1 (x),               &              & 
\end{align*}
for $y\in\hh$, $v\in V$, $w\in E$ and $x\in\gg$.
\end{prop}
The proof will consist of a series of computations that we will postpone to introduce several pieces of notation that will make easier its writing. \\
Given a $2$-extension as the one in the statement of proposition \ref{Ind2Rep} and adopting the convention that the injective maps are inclusions, we use the linear splitting to get the usual isomorphisms $\hh\oplus V\cong\mathfrak{e}_0$ and $\gg\oplus W\cong\mathfrak{e}_1$, given by
\begin{eqnarray*}
\xymatrix{
(z,a) \ar@{|->}[r] & a + \sigma_k (z)
}
\end{eqnarray*}
for $k=0$ and $k=1$ respectively. Their inverses are
\begin{eqnarray*}
\xymatrix{
e \ar@{|->}[r] & (\pi_k (e),e-\sigma_k (\pi_k (e))).
}
\end{eqnarray*}
In order to upgrade these isomorphisms of vector spaces to Lie algebra isomorphisms, we consider 
\begin{align*}
[v_0+\sigma_0(y_0),v_1+\sigma_0(y_1)]_{\mathfrak{e}_0} & =[v_0,v_1]+[\sigma_0(y_0),v_1]+[v_0,\sigma_0(y_1)]+[\sigma_0(y_0),\sigma_0(y_1)] \\
                                                       & = \rho^0_0(y_0)v_1 - \rho^0_0(y_1)v_0 +[\sigma_0(y_0),\sigma_0(y_1)].
\end{align*}
We use the inverse isomorphism to define the bracket on $\hh\oplus V$. Since $\rho^0_0(y_0)v_1 - \rho^0_0(y_1)v_0\in V$,
\begin{align*}
\pi_0(\rho^0_0(y_0)v_1 - \rho^0_0(y_1)v_0 +[\sigma_0(y_0),\sigma_0(y_1)]) & =\pi_0([\sigma_0(y_0),\sigma_0(y_1)]_{\mathfrak{e}_0}) \\
                                                                          & =[\pi_0\sigma_0(y_0),\pi_0\sigma_0(y_1)]_\hh =[y_0,y_1];
\end{align*}
thus, the bracket is
\begin{eqnarray*}
[(y_0,v_0),(y_1,v_1)]_0:=([y_0,y_1],\rho^0_0(y_0)v_1 - \rho^0_0(y_1)v_0 -\omega_0(y_0,y_1)),
\end{eqnarray*}
where $\omega_0(y_0,y_1)$ is shorthand for $\sigma_0([y_0,y_1])-[\sigma_0(y_0),\sigma_0(y_1)]$. This is nothing but the usual twisted semi-direct product $\hh {}_{\rho^0_0}\oplus^{\omega_0} V$ from the theory of Lie algebra extensions. Hence, if one supposes conversely, that the bracket was defined using an abstract $\omega_0\in (\hh^*\wedge\hh^*)\otimes V$, one is going to rediscover that in order for such bracket to satisfy the Jacobi identity, $\omega_0$ needs to be a $2$-cocycle in the Lie algebra cohomology of $\hh$ with values in $\rho^0_0$. Recall the general formula for the differential of the complex of Lie algebra cochains with values in a representation $\rho$ is the same as the one given for Lie algebroids in chapter \ref{preliminarieschapter}. Following the same reasoning, one finds out that $\mathfrak{e}_1\cong\gg {}_{\rho^1_0\circ\mu}\oplus^{\omega_1} W$ as Lie algebras, with $\omega_1(x_0,x_1)=\sigma_1([x_0,x_1])-[\sigma_1(x_0),\sigma_1(x_1)]$; however, we are going to be able to waive the necessity of $\omega_1$ using the rest of the crossed module structure. \\
We now turn to the homomorphism $\epsilon$. Consider
\begin{align*}
\epsilon(w+\sigma_1(x)) & =\phi(w)+\epsilon(\sigma_1(x)),
\end{align*}
and use the inverse isomorphism to define the crossed module map. Since $\pi$ is a crossed module map,
\begin{align*}
\pi_0(\phi(w)+\epsilon(\sigma_1(x))) & =\pi_0(\phi(w))+\pi_0(\epsilon(\sigma_1(x)))) \\
                                     & =\mu(\pi_1\sigma_1(x)))=\mu(x);
\end{align*}
thus, the crossed module map is
\begin{eqnarray*}
\xymatrix{
(x,w) \ar@{|->}[r] & (\mu(x),\phi(w)+\varphi(x)),
}
\end{eqnarray*}
where $\varphi(x):=\epsilon(\sigma_1(x))-\sigma_0(\mu(x))$. We repeat this strategy one last time to get the action,
\begin{align*}
\Lie^\epsilon_{v+\sigma_0(y)}(w+\sigma_1(x)) & =\Lie_v w +\Lie_{\sigma_0(y)}w+\Lie_v\sigma_1(x)+\Lie_{\sigma_0(y)}\sigma_1(x) \\
                                             & =\rho^1_0(y)w -\rho_1(x)v +\Lie_{\sigma_0(y)}\sigma_1(x) ,
\end{align*}
Since $\pi$ is a crossed module map,
\begin{align*}
\pi_1(\rho^1_0(y)w -\rho_1(x)v +\Lie_{\sigma_0(y)}\sigma_1(x)) & =\pi_1(\Lie^\epsilon_{\sigma_0(y)}\sigma_1(x)) \\
                                                               & =\Lie^\mu_{\pi_0(\sigma_0(y))}\pi_1(\sigma_1(x))=\Lie_y x
\end{align*}
thus, the action is given by the equation
\begin{eqnarray*}
\Lie_{(y,v)}(x,w):=(\Lie _y x,\rho^1_0(y)w -\rho_1(x)v -\alpha(y;x))
\end{eqnarray*}
where $\alpha(y;x):=\sigma_1(\Lie _y x)-\Lie^\epsilon_{\sigma_0(y)}\sigma_1(x)$.\\
Using this data and the infinitesimal Peiffer equation for $\mathfrak{e}$, one readily sees that 
\begin{eqnarray*}
\omega_1(x_0,x_1)=\rho_1(x_1)(\varphi(x_0))+\alpha(\mu(x_0);x_1).
\end{eqnarray*}
\begin{proof}(of Proposition \ref{Ind2Rep})
We make the computations necessary to prove that $\rho^{\epsilon}_{\sigma}$ is a $2$-representation.
\begin{itemize}
\item Well-defined: We use the exactness of the sequences to see that the maps land where they are supposed to.
\begin{eqnarray*}
\pi_0([\sigma_0(y),v]_{\mathfrak{e}_0})=[\pi_0(\sigma_0(y)),\pi_0(v)]_\hh =[y,0]_\hh =0 & \Longrightarrow & \rho_0^0(y)v\in V, \\
\pi_1(\Lie^{\epsilon}_{\sigma_0(y)} w)=\Lie_{\pi_0(\sigma_0(y))}\pi_1(w) =\Lie_y 0   =0 & \Longrightarrow & \rho_0^1(y)w\in W, \\
\pi_1(-\Lie^{\epsilon}_v\sigma_1 (x))=-\Lie_{\pi_0(v)}\pi_1(\sigma_1 (x)) =-\Lie_0 x =0 & \Longrightarrow & \rho_1(x)v\in W. 
\end{eqnarray*}
Further, thus defined, $\rho_0^0(y)\circ\phi=\phi\circ\rho_0^1(y)$ for each $y\in\hh$. Indeed,
\begin{align*}
\rho_0^0(y)(\phi(w)) & =[\sigma_0(y),\phi(w)]_{\mathfrak{e}_0} =[\sigma_0(y),\epsilon(w)]_{\mathfrak{e}_0} \\
                     & =\epsilon(\Lie^{\epsilon}_{\sigma_0(y)} w) \\
                     & =\phi(\Lie^{\epsilon}_{\sigma_0(y)} w)=\phi(\rho_0^1(y)w);                   
\end{align*}
in these equations, we used that $\epsilon\circ j_1=j_0\circ\phi$, that $\Lie^{\epsilon}_{\sigma_0(y)}w\in W$ and the equivariance for the crossed module $\epsilon$. This proves that for all $y\in \hh$, $\rho_0(y):=(\rho_0^0(y),\rho_0^1(y))\in\ggl (\phi)_0$, as desired.
\item $\rho^0_0$ Lie algebra homomorphism: We write each side of the equation, and then justify why their difference is zero.
\begin{align*}
\rho^0_0([y_0,y_1]_\hh)v & =[\sigma_0([y_0,y_1]_\hh),v]_{\mathfrak{e}_0},
\end{align*}
and
\begin{align*}
[\rho^0_0(y_0),\rho^0_0(y_1)]v & =\rho^0_0(y_0)\rho^0_0(y_1)v-\rho_0^0(y_1)\rho_0^0(y_0)v \\
                               & =\rho^0_0(y_0)([\sigma_0(y_1),v]_{\mathfrak{e}_0})-\rho_0^0(y_1)([\sigma_0(y_0),v]_{\mathfrak{e}_0}) \\
                               & =[\sigma_0(y_0),[\sigma_0(y_1),v]_{\mathfrak{e}_0}]_{\mathfrak{e}_0}-[\sigma_0(y_1),[\sigma_0(y_0),v]_{\mathfrak{e}_0}]_{\mathfrak{e}_0} \\
                               & =[[\sigma_0(y_0),\sigma_0(y_1)]_{\mathfrak{e}_0},v]_{\mathfrak{e}_0},
\end{align*}                               
where the last equality is due to the Jacobi identity. Then, considering the difference, we get
\begin{align*}
(\rho_0^0([y_0,y_1]_\hh)-[\rho^0_0(y_0),\rho^0_0(y_1)])v & =[\omega_0(y_0,y_1),v]_{\mathfrak{e}_0},
\end{align*}
which is zero given that $\omega_0(y_0,y_1)\in V$.
\item $\rho_0^1$ Lie algebra homomorphism: We proceed with the same strategy.
\begin{align*}
\rho_0^1([y_0,y_1]_\hh)w & =\Lie^{\epsilon}_{\sigma_0([y_0,y_1]_\hh)} w,
\end{align*}
and
\begin{align*}
[\rho_0^1(y_0),\rho_0^1(y_1)]w & =\rho_0^1(y_0)\rho_0^1(y_1)w-\rho_0^1(y_1)\rho_0^1(y_0)w \\
                               & =\Lie^{\epsilon}_{\sigma_0(y_0)}\Lie^{\epsilon}_{\sigma_0(y_1)}w-\Lie^{\epsilon}_{\sigma_0(y_1)}\Lie^{\epsilon}_{\sigma_0(y_0)}w \\
                               & =[\Lie^{\epsilon}_{\sigma_0(y_0)},\Lie^{\epsilon}_{\sigma_0(y_1)}]w =\Lie^{\epsilon}_{[\sigma_0(y_0),\sigma_0(y_1)]_{\mathfrak{e}_0}}w,
\end{align*}                               
where the last line is because $\Lie^\epsilon$ is a Lie algebra action. Then, considering the difference, we get
\begin{align*}
(\rho_0^1([y_0,y_1]_\hh)-[\rho_0^1(y_0),\rho_0^1(y_1)])w & =\Lie^{\epsilon}_{\omega_0(y_0,y_1)} w,
\end{align*}
which is zero given that $\omega_0(y_0,y_1)\in V$.
\item $\rho_1$ Lie algebra homomorphism: We proceed with the same strategy.
\begin{align*}
\rho_1([x_0,x_1]_\gg)v & =-\Lie^{\epsilon}_v\sigma_1([x_0,x_1]_\gg) ,
\end{align*}
and
\begin{align*}
[\rho_1(x_0),\rho_1(x_1)]_\phi v & =\rho_1(x_0)\phi\rho_1(x_1)v-\rho_1(x_1)\phi\rho_1(x_0)v \\
                                 & =-\rho_1(x_0)\epsilon(\Lie^{\epsilon}_v\sigma_1(x_1))+\rho_1(x_1)\epsilon(\Lie^{\epsilon}_v\sigma_1(x_0)) \\
                                 & =-\rho_1(x_0)[v,\epsilon(\sigma_1(x_1))]_{\mathfrak{e}_0}+\rho_1(x_1)[v,\epsilon(\sigma_1(x_0))]_{\mathfrak{e}_0} \\
                                 & =\Lie^{\epsilon}_{[v,\epsilon(\sigma_1(x_1))]_{\mathfrak{e}_0}}\sigma_1(x_0)-\Lie^{\epsilon}_{[v,\epsilon(\sigma_1(x_0))]_{\mathfrak{e}_0}}\sigma_1(x_1). 
\end{align*}                               
Now, 
\begin{align*}
\Lie^{\epsilon}_{[v,\epsilon(\sigma_1(x_1))]_{\mathfrak{e}_0}}\sigma_1(x_0) & 
                                    =\Lie^{\epsilon}_v\Lie^{\epsilon}_{\epsilon(\sigma_1(x_1))}\sigma_1(x_0)-\Lie^{\epsilon}_{\epsilon(\sigma_1(x_1))}\Lie^{\epsilon}_v\sigma_1(x_0) \\
                                  & =\Lie^{\epsilon}_v[\sigma_1(x_1),\sigma_1(x_0)]_{\mathfrak{e}_1}-[\sigma_1(x_1),\Lie^{\epsilon}_v\sigma_1(x_0)]_{\mathfrak{e}_1},
\end{align*} 
thanks to the infinitesimal Peiffer equation, and consequently,
\begin{align*}
[\rho_1(x_0),\rho_1(x_1)]_\phi v & =\Lie^{\epsilon}_v[\sigma_1(x_1),\sigma_1(x_0)]_{\mathfrak{e}_1}-[\sigma_1(x_1),\Lie^{\epsilon}_v\sigma_1(x_0)]_{\mathfrak{e}_1}\\
                                 & \qquad -\Lie^{\epsilon}_v[\sigma_1(x_0),\sigma_1(x_1)]_{\mathfrak{e}_1}+[\sigma_1(x_0),\Lie^{\epsilon}_v\sigma_1(x_1)]_{\mathfrak{e}_1} \\
                                 & =-2\Lie^{\epsilon}_v[\sigma_1(x_0),\sigma_1(x_1)]_{\mathfrak{e}_1}+\Lie^{\epsilon}_v[\sigma_1(x_0),\sigma_1(x_1)]_{\mathfrak{e}_1} \\
                                 & =-\Lie^{\epsilon}_v[\sigma_1(x_0),\sigma_1(x_1)]_{\mathfrak{e}_1},
\end{align*} 
where the second equality follows since $\Lie^\epsilon$ is an action by derivations. Then, considering the difference, we get
\begin{align*}
(\rho_1([x_0,x_1]_\gg)-[\rho_1(x_0),\rho_1(x_1)])v & =-\Lie^{\epsilon}_v\omega_1(x_0,x_1),
\end{align*}
which is zero given that $\omega_1(x_0,x_1)\in W$.
\item $\rho_0\circ\mu =\Delta\circ\rho_1$: This equation breaks into two components, one in $\ggl (W)$ and one in $\ggl (V)$; namely, 
\begin{eqnarray*}
\rho_0^1(\mu(x))=\rho_1(x)\circ\phi , & \rho_0^0(\mu(x))=\phi\circ\rho_1(x) ,
\end{eqnarray*} 
for each $x\in\gg$. These relations follow as, using the strategy above,
\begin{align*}
\rho_0^1(\mu(x))w = \Lie^{\epsilon}_{\sigma_0(\mu(x))}w,\quad & \quad\rho_0^0(\mu(x))v = [\sigma_0(\mu(x)),v]_{\mathfrak{e}_0} ,
\end{align*}
and
\begin{align*}
\rho_1(x)\phi(w) & =-\Lie^{\epsilon}_{\phi(w)}\sigma_1(x)\quad & \quad\phi(\rho_1(x)v) & =\phi(-\Lie^{\epsilon}_v\sigma_1(x))    \\
                 & =-\Lie^{\epsilon}_{\epsilon(w)}\sigma_1(x) &                   & =-\epsilon(\Lie^{\epsilon}_v\sigma_1(x))     \\
                 & =-[w,\sigma_1(x)]_{\mathfrak{e}_1}&                            & =-[v,\epsilon(\sigma_1(x))]_{\mathfrak{e}_0} \\
                 & =[\sigma_1(x),w]_{\mathfrak{e}_1}&                             & =[\epsilon(\sigma_1(x)),v]_{\mathfrak{e}_0}. \\
                 & =\Lie^{\epsilon}_{\epsilon(\sigma_1(x))}w,
\end{align*}
Thus, considering the respective differences, we get
\begin{align*}
(\rho_1(x)\phi -\rho_0^1(\mu(x)))w & =\Lie^{\epsilon}_{\varphi(x)}w,\quad & \quad(\phi\rho_1(x)-\rho_0^0(\mu(x)))v & = [\varphi(x),v]_{\mathfrak{e}_0}
\end{align*}
which are both zero, since $\varphi(x)\in V$ as desired.
\item $\rho_1$ respects the actions: One last time.
\begin{align*}
\rho_1(\Lie_y x)v & =-\Lie^{\epsilon}_v\sigma_1(\Lie_y x),
\end{align*}
and
\begin{align*}
(\Lie^{\phi}_{\rho_0(y)}\rho_1(x))v & =\rho_0^1(y)\rho_1(x)v-\rho_1(x)\rho_0^0(y)v \\
                                    & =\rho_0^1(y)(-\Lie^{\epsilon}_v\sigma_1(x))-\rho_1(x)[\sigma_0(y),v]_{\mathfrak{e}_0} \\
                                    & =-\Lie^{\epsilon}_{\sigma_0(y)}\Lie^{\epsilon}_v\sigma_1(x)+\Lie^{\epsilon}_{[\sigma_0(y),v]_{\mathfrak{e}_0}}\sigma_1(x) \\
                                    & =-\Lie^{\epsilon}_v\Lie^{\epsilon}_{\sigma_0(y)}\sigma_1(x)
\end{align*}                               
where the last equality is again since $\Lie^\epsilon$ is a Lie algebra action. Then, considering the difference, we get
\begin{align*}
(\rho_1(\Lie_y x)-\Lie^{\phi}_{\rho_0(y)}\rho_1(x))v & =-\Lie^{\epsilon}_v\alpha(y;x),
\end{align*}
which is zero given that $\alpha(y;x)\in W$.
\end{itemize}
\end{proof} 
%*** Check that the signs correspond to the ones we used before and the ones that appear later. *** 
Notice that in the case where the splitting can be taken to be a crossed module morphism, and accordingly $\omega$, $\alpha$ and $\varphi$ vanish, the given formulas for the crossed module structure on the extension coincide with those of the semi-direct sum defined in the previous section.
\begin{prop}\label{2-cocycles}
Let $\rho$ be a $2$-representation of $\xymatrix{\gg \ar[r]^\mu & \hh}$ on $\xymatrix{W \ar[r]^\phi & V}$. Given a triple $(\omega,\alpha,\varphi)\in((\hh^*\wedge\hh^*)\otimes V)\oplus(\hh^*\otimes\gg^*\otimes W)\oplus(\gg^*\otimes V)$, it defines a $2$-extension  \begin{eqnarray*}
\xymatrix{
0 \ar[r] & W \ar[d]_\phi\ar@{^{(}->}[r] & \gg {}_{\rho_0^1\circ\mu}\oplus^{\omega_1}W \ar[d]_\epsilon\ar[r]^{\qquad pr_1} & \gg \ar[d]^\mu\ar[r] & 0  \\
0 \ar[r] & V \ar@{^{(}->}[r]            & \hh {}_{\rho_0^0}\oplus^{\omega_1}V \ar[r]_{\qquad pr_1}                        & \hh \ar[r]           & 0, 
}
\end{eqnarray*}
with 
\begin{eqnarray*}
\omega_1(x_0,x_1)           & := & \rho_1(x_1)\varphi(x_0)+\alpha(\mu(x_0);x_1),  \\
\epsilon(x,w)                & = & (\mu(x),\phi(w)+\varphi(x)),                   \\
\Lie^{\epsilon}_{(y,v)}(x,w) & = & (\Lie _y x,\rho_0^1(y)w-\rho_1(x)v-\alpha(y;x))
\end{eqnarray*} 
if, and only if the following equations are satisfied 
\begin{itemize}
\item[i)] $\delta\omega =0$. Explicitely, for all triples $y_0,y_1,y_2\in\hh$, 
\begin{eqnarray*}
\rho^0_0(y_0)\omega(y_1,y_2)-\omega([y_0,y_1],y_2)+\circlearrowleft =0.
\end{eqnarray*}  
Here, $\circlearrowleft$ stands for cyclic permutations. 
\item[ii)] $\rho_1(x_0)\varphi(x_1)+\alpha(\mu(x_1);x_0)+\circlearrowleft =0$
\item[iii)] For all triples $x_0,x_1,x_2\in\gg$, 
\begin{eqnarray*}
\rho_0^1(\mu(x_0))(\rho_1(x_2)\varphi(x_1)+\alpha(\mu(x_2);x_1))-\rho_1(x_2)\varphi([x_0,x_1])-\alpha(\mu([x_0,x_1]);x_2)+\circlearrowleft =0.
\end{eqnarray*}
\item[iv)] $\omega(y,\mu(x))=\phi\circ\alpha(y;x)+\rho^0_0(y)\varphi(x)-\varphi(\Lie_y x)$
\item[v)] $\alpha([y_0,y_1];x)-\rho_1(x)\omega(y_0,y_1)=\rho^1_0(y_0)\alpha(y_1;x)-\alpha(y_1;\Lie_{y_0}x)+\circlearrowleft$
\item[vi)] For the contraction $a_y:=\alpha(y;-)\in\gg^*\otimes W$ seen as a $1$-cocycle with values in $\rho^1_0\circ\mu$,
\begin{eqnarray*}
\delta a_y(x_0,x_1)=\rho^1_0(y)\omega_1(x_0,x_1)-(\omega_1(\Lie _y x_0,x_1)-\omega_1(\Lie _y x_1,x_0))
\end{eqnarray*}
\end{itemize}
\end{prop}
There is nothing to these equations, in fact in most of the examples we have computed explicitly, they are redundant. In the proof, the reader might find the meaning of each of these equations.
\begin{proof}
We make the computations necessary to prove that a triple $(\omega ,\alpha ,\varphi)$ subject to the equations in the statement defines a $2$-extension. \\
First, the usual theory of Lie algebra extensions, tells us that item $i)$ says that $\omega$ is a $2$-cocycle with values in the representation $\rho_0^0$, and, as such, it defines a Lie bracket on $\hh\oplus V$. Analogously, $\omega_1$ defines an extension if, and only if $\omega_1$ is a $2$-cocycle with values in the representation $\rho_0^1\circ\mu$. With these, the equations in the statement have the following meaning
\begin{itemize}
\item[ii)] says $\omega_1$ is skew-symmetric. 
\item[iii)] says $\omega_1$ is a $2$-cocycle.
\item[iv)] evaluated at $y=\mu(x')$ says that $\epsilon$ is a Lie algebra homomorphism; indeed, 
\begin{align*}
\epsilon([(x',w'),(x,w)]) & =\epsilon([x',x],\rho_0^1(\mu(x'))w-\rho_0^1(\mu(x))w'-\omega_1(x',x)) \\
                          & =(\mu([x',x]),\phi(\rho_0^1(\mu(x'))w-\rho_0^1(\mu(x))w'-\omega_1(x',x))+\varphi([x',x])),
\end{align*}
while on the other hand,
\begin{align*}
[\epsilon(x',w'),\epsilon(x,w)] & =[(\mu(x'),\phi(w')+\varphi(x')),(\mu(x),\phi(w)+\varphi(x))] \\
                                & =([\mu(x'),\mu(x)],-\omega(\mu(x'),\mu(x))+ \\
                                & \qquad\qquad\rho_0^0(\mu(x'))(\phi(w)+\varphi(x))-\rho_0^0(\mu(x))(\phi(w')+\varphi(x'))).
\end{align*}
Clearly, the first entries coincide; moreover, since $\rho_0^0(\mu(x))=\phi\rho_1(x)$, $\epsilon$ is a Lie algebra homomorphism if, and only if
\begin{eqnarray*}
\varphi([x',x])-\phi\circ\omega_1(x',x)=\rho_0^0(\mu(x'))\varphi(x)-\rho_0^0(\mu(x))\varphi(x')-\omega(\mu(x'),\mu(x)).
\end{eqnarray*}
Replacing, the definition of $\omega_1$ yields
\begin{eqnarray*}
\varphi([x',x])-\phi\circ\alpha(\mu(x');x)=\rho_0^0(\mu(x'))\varphi(x)-\omega(\mu(x'),\mu(x)),
\end{eqnarray*}
which is precisely the equation in item $iv)$. Additionally, this equation also implies that $\epsilon$ is equivariant, as the following computation shows. 
\begin{align*}
\epsilon(\Lie^{\epsilon}_{(y,v)}(x,w)) & =\epsilon(\Lie _y x,\rho_0^1(y)w-\rho_1(x)v-\alpha(y;x)) \\
                                       & =(\mu(\Lie _y x),\phi(\rho_0^1(y)w-\rho_1(x)v-\alpha(y;x))+\varphi(\Lie _y x)),
\end{align*}
while on the other hand,
\begin{align*}
[(y,v),\epsilon(x,w)] & =[(y,v),(\mu(x),\phi(w)+\varphi(x))] \\
                      & =([y,\mu(x)],\rho_0^0(y)(\phi(w)+\varphi(x))-\rho_0^0(\mu(x))v-\omega(y,\mu(x))).
\end{align*}
Again, it is clear that the first entries coincide, and using again the relation $\rho_0^0(\mu(x))=\phi\rho_1(x)$ together with $\rho_0(y)\in\ggl (\phi)_0$, $\epsilon$ is equivariant if, and only if
\begin{eqnarray*}
\varphi(\Lie _y x)-\phi\circ\alpha(y;x)=\rho_0^0(y)\varphi(x)-\omega(y,\mu(x)).
\end{eqnarray*}
\item[v)] says $\Lie$ is a Lie algebra homomorphism and thus an action: 
\begin{align*}
\Lie^{\epsilon}_{[(y_0,v_0),(y_1,v_1)]}(x,w) & =\Lie^{\epsilon}_{([y_0,y_1],\rho_0^0(y_0)v_1-\rho_0^0(y_1)v_0-\omega(y_0,y_1))}(x,w) \\
                                             & =(\Lie _{[y_0,y_1]} x,\rho_0^1([y_0,y_1])w-\alpha([y_0,y_1];x) \\
                                             & \qquad\qquad-\rho_1(x)(\rho_0^0(y_0)v_1-\rho_0^0(y_1)v_0-\omega(y_0,y_1))). 
\end{align*}
On the other hand,
\begin{align*}
\Lie^{\epsilon}_{(y_0,v_0)}\Lie^{\epsilon}_{(y_1,v_1)}(x,w) 
                        & =\Lie^{\epsilon}_{(y_0,v_0)}(\Lie _{y_1} x,\rho_0^1(y_1)w-\rho_1(x)v_1-\alpha(y_1;x)) \\
                        & =(\Lie _{y_0} \Lie _{y_1} x,\rho_0^1(y_0)(\rho_0^1(y_1)w-\rho_1(x)v_1-\alpha(y_1;x))  \\
                        & \qquad\qquad-\rho_1(\Lie _{y_1} x)v_0-\alpha(y_0;\Lie _{y_1} x)).
\end{align*}
Then, considering the cyclic difference, one realizes that the first entries coincide. Using the fact that $\rho_0^1$ is a Lie algebra representation and the fact that $\rho_1$ respects is compatible with the actions, one is left with the relation
\begin{eqnarray*}
\rho_1(x)\omega(y_0,y_1)-\alpha([y_0,y_1];x)=\rho_0^1(y_1)\alpha(y_0;x)+\alpha(y_1;\Lie _{y_0} x))-\circlearrowleft.
\end{eqnarray*}
\item[vi)] says that $\Lie$ acts by derivations:
\begin{align*}
\Lie^{\epsilon}_{(y,v)}[(x_0,w_0),(x_1,w_1)] & =\Lie^{\epsilon}_{(y,v)}([x_0,x_1],\rho_0^1(\mu(x_0))w_1-\rho_0^1(\mu(x_1))w_0-\omega_1(x_0,x_1)) \\
                                             & =(\Lie _{y} [x_0,x_1],\rho_0^1(y)(\rho_0^1(\mu(x_0))w_1-\rho_0^1(\mu(x_1))w_0-\omega_1(x_0,x_1)) \\
                                             & \qquad\qquad -\rho_1([x_0,x_1])v-\alpha(y;[x_0,x_1])) 
\end{align*}
On the other hand,
\begin{align*}
[\Lie^{\epsilon}_{(y,v)}(x_0,w_0),(x_1,w_1)] 
                        & =[(\Lie _{y} x_0,\rho_0^1(y)w_0-\rho_1(x_0)v-\alpha(y;x_0)),(x_1,w_1)] \\
                        & =([\Lie _{y} x_0,x_1],\rho_0^1(\mu(\Lie _{y} x_0))w_1-\omega_1(\Lie _{y} x_0,x_1) \\
                        & \qquad\qquad -\rho_0^1(\mu(x_1))(\rho_0^1(y)w_0-\rho_1(x_0)v-\alpha(y;x_0))).
\end{align*}
Considering again the cyclic difference, one realizes that the first entries coincide. Since $\mu$ is equivariant and $\rho_0^1$ is a Lie algebra representation, we are left with the relation
\begin{eqnarray*}
-\rho_0^1(y)\omega_1(x_0,x_1)-\alpha(y;[x_0,x_1])=-\omega_1(\Lie _{y} x_0,x_1)+\rho_0^1(\mu(x_1))\alpha(y;x_0))-\circlearrowleft.
\end{eqnarray*}
\end{itemize}
Finally, the infinitesimal Peiffer equation holds by the very definition of $\omega_1$; indeed,
\begin{align*}
[(x_0,w_0),(x_1,w_1)] & =([x_0,x_1],\rho_0^1(\mu(x_0))w_1-\rho_0^1(\mu(x_1))w_0-\omega_1(x_0,x_1)),
\end{align*}
whereas
\begin{align*}
\Lie^{\epsilon}_{\epsilon(x_0,w_0)}(x_1,w_1) & =\Lie^{\epsilon}_{(\mu(x_0),\phi(w_0)+\varphi(x_0))}(x_1,w_1) \\
                                             & =(\Lie _{\mu(x_0)}x_1,\rho_0^1(\mu(x_0))w_1-\rho_1(x_1)(\phi(w_0)+\varphi(x_0))-\alpha(\mu(x_0),x_1)).
\end{align*}
Due to the relation $\rho_0^1(\mu(x))=\rho_1(x)\phi$, we are left precisely with the defining equation for $\omega_1$.

\end{proof}

Next, we analyze what happens when we have got equivalent extensions as in
\begin{eqnarray*}
\xymatrix{
         &                     & \mathfrak{e}_1 \ar[dd]\ar@{.>}[ddr]^{\psi_1}\ar[drr] &          &                    &    \\
0 \ar[r] & W \ar[dd]\ar[ur]\ar[drr] &                                 &                          & \gg \ar[r] \ar[dd] & 0  \\
         &                     & \mathfrak{e}_0 \ar[drr]\ar@{.>}[ddr]_{\psi_0} & \mathfrak{f}_1 \ar[dd]\ar[ur] &      &    \\
0 \ar[r] & V \ar[ur]\ar[drr]        &                                 &                          & \hh \ar[r]         & 0. \\
         &                          &                                          & \mathfrak{f}_0 \ar[ur]        &      &   
}
\end{eqnarray*}
%*** Should this be rephrased in terms of sections? Otherwise, how to address that the induced representations are the same? ***
Picking a splitting of either extension and composing it with the isomorphism, one gets a splitting for the other extension. In picking the splittings compatibly so, the induced $2$-representations are identical. We use these compatible splittings to identify both $\mathfrak{e}$ and $\mathfrak{f}$ with their respective semi-direct sums, and we write $\psi$ in these coordinates. Since both components of $\psi$ are linear, and respect both inclusions and projections
\begin{eqnarray}\label{anIsoOfExts}
\psi_k(z,a)=(z,a+\lambda_k(z))
\end{eqnarray}
for some linear maps $\xymatrix{\lambda_0:\hh \ar[r] & V}$ and $\xymatrix{\lambda_1:\gg \ar[r] & W}$.
\begin{prop}\label{2-coboundaries}
Let $\rho$ be a $2$-representation of $\xymatrix{\gg \ar[r]^\mu & \hh}$ on $\xymatrix{W \ar[r]^\phi & V}$. Given two $2$-cocycles  $(\omega_k,\alpha_k,\varphi_k)$ as in the previous proposition, the induced extensions are equivalent if, and only if there are linear maps $\xymatrix{\lambda_0:\hh \ar[r] & V}$ and $\xymatrix{\lambda_1:\gg \ar[r] & W}$ verifying  
\begin{itemize}
\item $\omega_2-\omega_1 =\delta\lambda_0$. Explicitly, for all triples $y_0,y_1\in\hh$, 
\begin{eqnarray*}
\omega_2(y_0,y_1)-\omega_1(y_0,y_1) =\rho^0_0(y_0)\lambda_0(y_1)-\rho^0_0(y_1)\lambda_0(y_0)-\lambda_0([y_0,y_1])
\end{eqnarray*}  
\item $\alpha_2(y;x)-\alpha_1(y;x)=\rho^1_0(y)(\lambda_1(x))-\lambda_1(\Lie _y x)-\rho_1(x)(\lambda_0(y))$
\item $\varphi_2(x)-\varphi_1(x)=\lambda_0(\mu(x))-\phi(\lambda_1(x))$
\end{itemize}
\end{prop}
\begin{proof}
The first equation is the classic identification of isomorphic extensions with cohomologous cocycles. \\
The second equation says that if the isomorphim of extensions $\psi$ is defined by the formula \ref{anIsoOfExts}, it respects the actions. Indeed, 
\begin{align*}
    \psi_1(\Lie_{(y,v)}(x,w)) & =\psi_1(\Lie _y x,\rho^1_0(y)w -\rho_1(x)v -\alpha_1(y;x)) \\
                              & =(\Lie_yx,\rho^1_0(y)w-\rho_1(x)v-\alpha_1(y;x)+\lambda_1(\Lie_yx)),
\end{align*}
whereas,
\begin{align*}
    \Lie_{\psi_0(y,v)}\psi_1(x,w) & =\Lie_{(y,v+\lambda_0(y))}(x,w+\lambda_1(x)) \\
                                  & =(\Lie_yx,\rho^1_0(y)(w+\lambda_1(x))-\rho_1(x)(v+\lambda_0(y))-\alpha_2(y;x)) .
\end{align*}
Thus, these expressions coincide if, and only if the second equation from the statement holds.
The third equation says that $\psi$ commutes with the structural maps of the crossed modules. Indeed,
\begin{align*}
 \psi_0(\epsilon_1(x,w)) & =\psi_0(\mu(x),\phi(w)+\varphi_1(x)) \\
                         & =(\mu(x),\phi(w)+\varphi_1(x)+\lambda_0(\mu(x))),
\end{align*}
while on the other hand,
\begin{align*}
 \epsilon_2(\psi_1(x,w)) & =\epsilon_2(x,w+\lambda_1(x)) \\
                         & =(\mu(x),\phi(w+\lambda_1(x))+\varphi_2(x)).
\end{align*}
Thus, these expressions coincide if, and only if the last equation of the statement holds. To conclude the proof, notice that if the equations in the statement are verified, the cocycles
\begin{eqnarray*}
\omega_k'(x_0,x_1)=\rho_1(x_1)\varphi_k(x_0)+\alpha_k(\mu(x_0);x_1),
\end{eqnarray*}
for $k\in\lbrace 1,2\rbrace$ defining the top extensions are cohomologous too. Indeed,
\begin{align*}
    (\omega_2'-\omega_1')(x_0,x_1) & =\rho_1(x_1)(\varphi_2-\varphi_1)(x_0)+(\alpha_2-\alpha_1)(\mu(x_0);x_1) \\
                                   & =\rho_1(x_1)(\lambda_0(\mu(x_0))-\phi(\lambda_1(x_0)))+\rho^1_0(\mu(x_0))(\lambda_1(x_1))+ \\
                                   & \qquad\qquad\qquad\qquad\qquad\qquad -\lambda_1(\Lie_{\mu(x_0)}x_1)-\rho_1(x_1)(\lambda_0(\mu(x_0)) \\
                                   & =-\rho_1(x_1)\phi(\lambda_1(x_0))+\rho^1_0(\mu(x_0))(\lambda_1(x_1))-\lambda_1([x_0,x_1])=\delta\lambda_1(x_0,x_1),
\end{align*}
where the last equality holds, because $\rho$ is a $2$-representation, and therefore, $\rho_1(x)\phi=\rho_0^1(\mu(x))$.

\end{proof}

\section{The complex of Lie 2-algebra cochains with values in a 2-representation}
Inspired by the case of trivial coefficients, one would like to see Lie $2$-algebra cohomology as the cohomology of the total complex of a double complex. Consider the case $W=\lbrace 0\rbrace$. Then one gets a double complex analogous to that in Section \ref{AlgDcx}
\begin{eqnarray*}
\xymatrix{
\vdots                                                        & \vdots                                                            & \vdots                                                 &      \\ 
\bigwedge ^3\hh ^*\otimes V\ar[r]^{\partial}\ar[u]            & \bigwedge ^3\gg _1^*\otimes V \ar[r]^{\partial}\ar[u]            & \bigwedge ^3\gg _2^*\otimes V \ar[r]\ar[u]             & \dots\\
(\hh ^*\wedge \hh ^*)\otimes V \ar[r]^{\partial}\ar[u]^\delta & (\gg _1^*\wedge\gg _1^*)\otimes V \ar[r]^{\partial}\ar[u]^\delta & (\gg _2^*\wedge \gg _2^*)\otimes V \ar[r]\ar[u]^\delta & \dots\\
\hh ^*\otimes V \ar[r]^{\partial}\ar[u]^\delta                & \gg _1^*\otimes V \ar[r]^{\partial}\ar[u]^\delta                 & \gg _2^*\otimes V \ar[r]\ar[u]^\delta                  & \dots 
}
\end{eqnarray*} 
Recall from Example \ref{unitRep} that a $2$-representation on $\xymatrix{(0)\ar[r] & V}$ amounts to a usual representation $\rho$ of $\hh/\mu(\gg)$ on $V$. In order to get a meaningful double complex, there need to be induced representations of $\gg_p$ on $V$. The representation $\rho_0$ for the first column is clear then, it is the pull-back of $\rho$ along the projection onto $\hh/\mu(\gg)$. However, for the columns there are various possible choices. We are going to fix $\rho_p$, the induced representation on $\gg_p$, to be the pull-back of $\rho_0$ along the ``final target'' map, $\xymatrix{\hat{t}_p:\gg_p \ar[r] & \hh :(x_1,...,x_p,y) \ar@{|->}[r] & y+\sum_{j=1}^p\mu(x_j)}$. In so, the columns are Chevalley-Eilenberg complexes with respect to $\rho_p$. It is of crucial importance that all of the representations vanish on the ideal $\mu(\gg)$; otherwise, the above double structure does not commute, and fails thus to be a double complex.  \\
\begin{remark}
Notice that, precisely because the representation vanishes on $\mu(\gg)$, we could have taken the pull-back along the ``initial source'' map, which is the reversal $\hat{s}_p=\hat{t}_p\circ\hat{\iota}$, and the complex would still commute.
\end{remark}
Consider the case $V=\lbrace 0\rbrace$. As it was pointed out, a $2$-representation on such a $2$-vector space also corresponds to a representation of $\hh/\mu(\gg)$, this time round, on $W$. When describing the algebraic data defining an extension with values in such a representation, there appears a second complex
\begin{eqnarray*}
\xymatrix{
\vdots                                                             & \vdots                                                            & \vdots                                     &      \\ 
\bigwedge ^2\hh ^*\otimes\gg ^*\otimes W \ar[r]^{\partial }\ar[u] & \bigwedge ^2 \gg _1^*\otimes\gg ^*\otimes W \ar[r]^{\partial }\ar[u]                  & \bigwedge ^2 \gg _2^*\otimes\gg ^*\otimes W \ar[r]\ar[u]             & \dots\\
\hh ^*\otimes \gg ^*\otimes W \ar[r]^{\partial }\ar[u]^{\delta '} & \gg _1^*\otimes \gg ^*\otimes W \ar[r]^{\partial }\ar[u]^{\delta '}                            & \gg _2^*\otimes\gg ^*\otimes W \ar[r]\ar[u]^{\delta '} & \dots\\
\gg ^*\otimes W \ar[r]^{0}\ar[u]^{\delta '}                & \gg ^*\otimes W \ar[r]^{id_{\gg ^*\otimes W}}\ar[u]^{\delta '}                     & \gg ^*\otimes W \ar[r]\ar[u]^{\delta '}    & \dots 
}
\end{eqnarray*} 
which is analogous to the one above. It formally is the double complex of the Lie $2$-algebra taking values in the representation 
\begin{eqnarray*}
\xymatrix{
\rho ':\hh \ar[r] & \ggl (\gg ^*\otimes W):y\ar@{|->}[r] & \rho '_y=\rho_0^1(y)-\Lie^*_y .
}
\end{eqnarray*}

In the general case, when the $2$-representation does not vanish on the ideal $\mu(\gg)$, both complexes above appear, though they cease to be double complexes in that they fail to commute. Moreover, they both appear intertwined with each other in the following sense: looking at the difference of the two ways to go around a given square in the first complex, one realizes that, when evaluated, $\delta\partial\omega$ and $\partial\delta\omega$ yield isomorphic elements in the $2$-vector space. Furthermore, the isomorphism is given by an element that comes from the second complex as we will make clear after appropriate exemplifications. For instance, consider $v\in V$, then
\begin{align*}
(\delta\partial -\partial\delta)v(x,y) & =\delta v(x,y)-\delta v(y+\mu(x))+\delta v(y) \\
                                       & =\rho_0^0(y+\mu(x))v-\rho_0^0(y)v \\
                                       & =\rho_0^0(\mu(x))v =\phi(\rho_1(x)v).
\end{align*}
Notice that there is a natural map $\xymatrix{\delta_{(1)}:V \ar[r] & \gg^*\otimes W}$ induced by $\rho_1$ and defined by $\delta_{(1)}v(x)=\rho_1(x)v$; hence, we can write the difference above as being $\delta\partial -\partial\delta=\phi\circ\delta_{(1)}$.
Also, for $\lambda\in\hh^*\otimes V$,
\begin{align*}
\delta\partial\lambda(x_0,y_0;x_1,y_1) & =\rho_0^0(y_0+\mu(x_0))\partial\lambda(x_1,y_1)-\rho_0^0(y_1+\mu(x_1))\partial\lambda(x_0,y_0)-\partial\lambda ([(x_0,y_0);(x_1,y_1)]) \\
                                       & =\rho_0^0(y_0+\mu(x_0))\lambda(\mu(x_1))-\rho_0^0(y_1+\mu(x_1))\lambda(\mu(x_0))-\lambda (\mu([x_0,x_1]+\Lie_{y_0}x_1-\Lie_{y_1}x_0), 
\end{align*}
and 
\begin{align*}
\partial\delta\lambda(x_0,y_0;x_1,y_1) & =\delta\lambda(y_0+\mu(x_0),y_1+\mu(x_1))-\delta\lambda(y_0,y_1) \\
                                       & =\delta\lambda(y_0,\mu(x_1))-\delta\lambda(y_1,\mu(x_0))+\delta\lambda(\mu(x_0),\mu(x_1)).                                        
\end{align*}
Each of these terms has the form 
\begin{align*}
\delta\lambda(y,\mu(x)) & =\rho_0^0(y)\lambda(\mu(x))-\rho_0^0(\mu(x))\lambda(y)-\lambda([y,\mu(x)]) \\
                        & =\rho_0^0(y)\lambda(\mu(x))-\rho_0^0(\mu(x))\lambda(y)-\lambda(\mu(\Lie_y x))
\end{align*}
and thus, the difference is
\begin{align*}
(\delta\partial -\partial\delta)\lambda(x_0,y_0;x_1,y_1) & =\rho_0^0(\mu(x_0))\lambda(y_1)-\rho_0^0(\mu(x_1))\lambda(y_0) \\
                                                         & =\phi(\rho_1(x_0)\lambda(y_1)-\rho_1(x_1)\lambda(y_0)).
\end{align*}
The map of the previous example can be generalized to $\xymatrix{\delta_{(1)}:\hh^*\otimes V \ar[r] & \hh^*\otimes\gg^*\otimes W}$, given by $\delta_{(1)}\lambda(y,x)=\rho_1(x)\lambda(y)$; therefore, this time around, we can write the difference as 
\begin{align*}
(\delta\partial -\partial\delta)\lambda(x_0,y_0;x_1,y_1) & =\phi(\delta_{(1)}\lambda(y_1,x_0)-\delta_{(1)}\lambda(y_0,x_1)).
\end{align*}
In general, since $\gg_{p+1}\cong\gg^{p+1}\oplus\hh$, one may identify $\xi\in\gg_{p+1}$ with a tuple $(x^0,...,x^p;y)$ and we can write the following relation.
\begin{prop}\label{AlgUpToHomotopy r=0}
Let $\omega\in\bigwedge^q\gg_p^*\otimes V$ and $\Xi=(\xi_0,...,\xi_q)\in\gg_{p+1}^{q+1}$. Then
\begin{eqnarray*}
\delta\partial\omega(\Xi)=\partial\delta\omega(\Xi)+\phi\Big{(}\sum_{j=0}^q(-1)^j\rho_1(x_j^0)\omega(\partial_0\Xi(j))\Big{)}
\end{eqnarray*}
\end{prop}
\begin{proof}
First of all, recall that by definition
\begin{eqnarray*}
\partial_k\Xi=(\partial_k\xi_0,...,\partial_k\xi_q);
\end{eqnarray*}
hence, trivially, $\partial_k(\Xi(j))=(\partial_k\Xi)(j)$, so we can drop the parenthesis. This shows that there is no ambiguity in the formula of the statement. Moreover, under the identification $\xi_j$ with $(x_j^0,...,x_j^p;y_j)$, we have got further identifications
\begin{eqnarray*}
\partial _k (\xi_j)\sim
  \begin{cases}
    (x_j^1,...,x_j^p;y_j)                     & \quad \text{if } k=0  \\
    (x_j^0,...,x_j^{k-1}+x_j^k,...,x_j^p;y_j) & \quad \text{if } 0<k\leq p\\
    (x_j^0,...,x_j^{p-1};y_j+\mu(x_j^p))      & \quad \text{if } k=p+1. 
  \end{cases}
\end{eqnarray*}
Now, computing
\begin{align*}
\delta\partial\omega(\Xi) & =\sum_{j=0}^{q}(-1)^j\rho_0^0(\hat{t}_{p+1}(\xi_j))\partial\omega(\Xi(j))+\sum_{m<n}(-1)^{m+n}\partial\omega([\xi_m,\xi_n],\Xi(m,n)) \\
						  & =\sum_{j=0}^{q}(-1)^j\rho_0^0(\hat{t}_{p+1}(\xi_j))\sum_{k=0}^{p+1}(-1)^k\omega(\partial_k\Xi(j))+\sum_{m<n}(-1)^{m+n}\sum_{k=0}^{p+1}(-1)^k\omega(\partial_k[\xi_m,\xi_n],\partial_k(\Xi(m,n)));
\end{align*}
while on the other hand,
\begin{align*}
\partial\delta\omega(\Xi) & =\sum_{k=0}^{p+1}(-1)^k\delta\omega(\partial_k\Xi) \\
						  & =\sum_{k=0}^{p+1}(-1)^k\Big{(}\sum_{j=0}^{q}(-1)^j\rho_0^0(\hat{t}_p(\partial_k\xi_j))\omega(\partial_k\Xi(j))+\sum_{m<n}(-1)^{m+n}\omega([\partial_k\xi_m,\partial_k\xi_n],(\partial_k\Xi)(m,n))\Big{)}.
\end{align*}
Notice that, again $\partial_k(\Xi(m,n))=(\partial_k\Xi)(m,n)$ trivially; therefore, since $\partial_k$ is a Lie algebra homomorphism for each $k$, the second term in the expressions above coincide. As for the first term, from the identifications above, one sees that 
\begin{eqnarray*}
\hat{t}_p(\partial_k\xi_j)=
  \begin{cases}
    y_j+\sum_{r=1}^p\mu(x_j^r)=\hat{t}_{p+1}(\xi_j)-\mu(x_j^0) & \quad \text{if } k=0  \\
    y_j+\sum_{r=0}^p\mu(x_j^r)=\hat{t}_{p+1}(\xi_j)            & \quad \text{otherwise}. 
  \end{cases}
\end{eqnarray*}
This, together with the relation $\rho_0^0(\mu(x))=\phi\rho_1(x)$ that comes from the $2$-representation, imply that the difference is 
\begin{align*}
    (\delta\partial-\partial\delta)\omega(\Xi) & =\sum_{j=0}^q(-1)^j\rho_0^0(\mu(x_j^0))\omega(\partial_0\Xi(j)) \\
                                               & =\phi\Big{(}\sum_{j=0}^q(-1)^j\rho_1(x_j^0)\omega(\partial_0\Xi(j))\Big{)},
\end{align*}
as desired. 

\end{proof}
To conclude that indeed the isomorphims come from the second complex consider the obvious generalization
\begin{eqnarray*}
\xymatrix{
\delta_{(1)}:\bigwedge^q\gg_p^*\otimes V \ar[r] & \bigwedge^q\gg_p^*\otimes\gg^*\otimes W
} \\
\delta_{(1)}\omega(\Xi;x)=\rho_1(x)\omega(\Xi)\qquad\qquad
\end{eqnarray*}
for $\Xi\in\gg_p^q$ and $x\in\gg$, and define the difference map
\begin{eqnarray*}
\xymatrix{
\Delta:\bigwedge^q\gg_p^*\otimes\gg^*\otimes W \ar[r] & \bigwedge^{q+1}\gg_{p+1}^*\otimes V 
} \\
\Delta\alpha(\Xi)=\phi\Big{(}\sum_{j=0}^q(-1)^j\alpha(\partial_0\Xi(j);x_j^0)\Big{)}
\end{eqnarray*}
for $\Xi=(\xi_0,...,\xi_q)\in\gg_{p+1}^{q+1}$ and under the identifications $\xi_j\sim(x_j^0,...,x_j^p;y_j)$. Then, the relation in proposition \ref{AlgUpToHomotopy r=0} is seen to coincide with
\begin{eqnarray*}
\delta\partial-\partial\delta=\Delta\circ\delta_{(1)}.
\end{eqnarray*}
As a consequence of this theorem, if the $2$-representation takes values on a Lie group bundle internal to the category of vector spaces, i.e. on a $2$-vector space $\xymatrix{W \ar[r]^0 & V}$, the first diagram above does commute. \\
We turn now to the difference $\partial\delta '-\delta '\partial$. As it turns out, this difference is not zero either. In fact, the same computation of the proof of the latter proposition, together with the infinitesimal Peiffer equation, show that for $\alpha\in\bigwedge^q\gg_p^*\otimes\gg^*\otimes W$, $x\in\gg$ and $\Xi$ as in the statement of the proposition,
\begin{align*}
    (\delta'\partial-\partial\delta')\alpha(\Xi;x) & =\sum_{j=0}^q(-1)^j\rho'(\mu(x_j^0))\alpha(\partial_0\Xi(j);x) \\
                                                   & =\sum_{j=0}^q(-1)^j\Big{(}\rho_0^1(\mu(x_j^0))\alpha(\partial_0\Xi(j);x)-\alpha(\partial_0\Xi(j);[x_j^0,x])\Big{)}.
\end{align*}
This expression is short of a few terms to be the (sum of the) Chevalley-Eilenberg differential(s) of $\alpha(\partial_0\Xi(j);-)\in\gg^*\otimes W$ with respect to $\rho_0^1\circ\mu$. What is more, the missing elements can be thought of as coming from the first complex in the sense that
\begin{align*}
    \delta_{(1)}\Delta\alpha(\Xi;x) & =\rho_1(x)\Delta(\Xi) \\
                                    & =\rho_1(x)\phi\Big{(}\sum_{j=0}^q(-1)^j\alpha(\partial_0\Xi(j);x_j^0)\Big{)};
\end{align*}
hence, using the fact that $\rho_1(x)\phi=\rho_0^0(\mu(x))$,
\begin{align*}
    (\delta'\partial-\partial\delta' & -\delta_{(1)}\Delta)\alpha(\Xi;x) \\
    & =\sum_{j=0}^q(-1)^j\Big{(}\rho_0^1(\mu(x_j^0))\omega(\partial_0\Xi(j);x)-\rho_0^0(\mu(x))\alpha(\partial_0\Xi(j);x_j^0)-\omega(\partial_0\Xi(j);[x_j^0,x])\Big{)} \\
    & =\sum_{j=0}^q(-1)^j\delta_{(1)}\alpha(\partial_0\Xi(j);x_j^0,x),
\end{align*}
where we write $\delta_{(1)}$ for the Chevalley-Eilenberg differential of $\gg$ with values in $\rho_0^1\circ\mu$. Define a new difference map
\begin{eqnarray*}
\xymatrix{
\Delta:\bigwedge^q\gg_p^*\otimes\bigwedge^2\gg^*\otimes W \ar[r] & \bigwedge^{q+1}\gg_{p+1}^*\otimes\gg^*\otimes W 
} \\
\Delta\omega(\Xi;x)=\sum_{j=0}^q(-1)^j\omega(\partial_0\Xi(j);x_j^0,x),\qquad\quad
\end{eqnarray*}
for $x\in\gg$, $\Xi=(\xi_0,...,\xi_q)\in\gg_{p+1}^{q+1}$ and under the identifications $\xi_j\sim(x_j^0,...,x_j^p;y_j)$. Then, the commutativity of the second diagram is revealed to be written as 
\begin{eqnarray*}
\delta'\partial-\partial\delta'=\delta_{(1)}\circ\Delta+\Delta\circ\delta_{(1)}.
\end{eqnarray*}
We can schematize this relation by the following diagram:
\begin{eqnarray*}
\xymatrix{
\circ \ar@{.}[rr]\ar@{.}[dr] & & \bigwedge^{q+1}\gg_p^*\otimes\gg^*\otimes W \ar@{.}[rr]\ar[dr]^\partial & & \circ\ar@{.}[dr] & \\
 & \bigwedge^{q+1}\gg_{p+1}^*\otimes V \ar[rr]^{\delta_{(1)}\qquad\quad} & & \bigwedge^{q+1}\gg_{p+1}^*\otimes\gg^*\otimes W\ar@{.}[rr] & & \circ \\
\circ \ar@{.}[rr]\ar@{.}[dr]\ar@{.}[uu] & & \bigwedge^q\gg_p^*\otimes\gg^*\otimes W \ar'[r]^{\quad\delta_{(1)}}[rr]\ar[dr]^\partial\ar[ul]^\Delta\ar'[u]_{\delta'}[uu] & & \bigwedge^q\gg_p^*\otimes\bigwedge^2\gg^*\otimes W\ar@{.}[dr]\ar@{.}[uu]\ar[ul]_\Delta & \\
 & \circ \ar@{.}[rr]\ar@{.}[uu] & & \bigwedge^q\gg_{p+1}^*\otimes\gg^*\otimes W\ar@{.}[rr]\ar[uu]_(.35){\delta'} & & \circ \ar@{.}[uu]
}
\end{eqnarray*}
Thus, we can regard the commutativity of the second diagram above as being controlled by elements in a third lattice of vector spaces 
\begin{eqnarray*}
\xymatrix{
\vdots                                               & \vdots                                                  & \vdots                                                                 &       \\ 
\bigwedge ^2\hh ^*\otimes\bigwedge ^2\gg ^*\otimes W & \bigwedge ^2 \gg _1^*\otimes\bigwedge ^2\gg ^*\otimes W & \bigwedge ^2 \gg _2^*\otimes\bigwedge ^2\gg ^*\otimes W & \dots \\
\hh ^*\otimes\bigwedge ^2\gg ^*\otimes W             & \gg _1^*\otimes\bigwedge ^2\gg ^*\otimes W              & \gg _2^*\otimes\bigwedge ^2\gg ^*\otimes W                            & \dots \\
\bigwedge ^2\gg ^*\otimes W                          & \bigwedge ^2\gg ^*\otimes W                             & \bigwedge ^2\gg ^*\otimes W                                            & \dots 
}
\end{eqnarray*} 
Of course, this can be endowed with the structure of a complex. Again, this will fail to be a double complex, and we will be able to encode its commutativity in further tensors of the double complex of the Lie $2$-algebra. This process will end up leaving a three dimensional lattice of vector spaces
\begin{eqnarray*}
C^{p,q}_r(\gg _1,\phi):=\bigwedge ^q\gg _p^*\otimes\bigwedge ^r\gg ^*\otimes W,
\end{eqnarray*}
where, we use the convention $\gg_0 =\hh$ and for $r=0$, the coefficients take values on $V$. For constant $r$, these vector spaces come with the structure of a complex, which is essentially the \LA -double complex of the Lie $2$-algebra tensored with $\bigwedge^r\gg^*\otimes W$. That is, it comes with groupoid differentials in the rows, and with Chevalley-Eilenberg differentials in the columns with respect to the pull-back along $\hat{t}_p$ of  
\begin{eqnarray*}
\xymatrix{
\rho ^{(r)}:\hh \ar[r] & \ggl (\bigwedge ^r\gg ^*\otimes W)
}
\end{eqnarray*}
given for $\omega\in\bigwedge ^r\gg ^*\otimes W$, $x_1,...,x_r$ and $y\in\hh$ by
\begin{eqnarray*}
\rho ^{(r)}(y)\omega(x_1,...,x_r):=\rho_0^1(y)\omega(x_1,...,x_r)-\sum_{k=1}^r\omega(x_1,...,\Lie _y x_k,...,x_r).
\end{eqnarray*}
\begin{remark}
Interestingly, these representations have already appeared in the literature \cite{HochSer}. The form in which they appear is restricted to the example of the inclusion of an ideal, regarded as a Lie $2$-algebra. % What are they used for? (Later: Just to define the Chevalley Eilenberg differential..., but it seems to be useful though, it appears in Evans work On the relative Hochschild-Serre). Is there a connection to this theory?
\end{remark}
As we remarked right after the definition of the page $r=0$ complex, we could have pulled back the representation along $\hat{s}_p$, there is no economy in working with either representation: one will always pay a computational price somewhere. Aside from these maps, for constant $q$, there is a second complex whose rows have been described already and whose columns are also Chevalley-Eilenberg differentials. As we saw, these will extend the differential of $\gg$ with values in $\rho_0^1\circ\mu$. That is, for constant $q$ and $p$, the complex is the Chevalley-Eilenberg complex of $\gg$ with values in
\begin{eqnarray*}
\xymatrix{
\rho _{(1)}:\gg \ar[r] & \ggl(\bigwedge^q\gg _p^*\otimes W)
}
\end{eqnarray*}
is given for $\alpha\in\bigwedge^q\gg _p^*\otimes W$, $\Xi\in\gg_p^q$ and $x\in\gg$ by
\begin{eqnarray*}
\rho _{(1)}(x)\alpha(\Xi):=\rho_0^1(\mu(x))\alpha(\Xi).
\end{eqnarray*}
%The heavy reliance on $\gg$ might rise some concerns in the reader, as this construction may seem ad hoc for Lie $2$-algebras; however, recall that $\gg$ is the core of the Lie $2$-algebra regarded as an \LA -groupoid. \\
We would like to point out that there is a heavy reliance on the space $\gg$ which, in principle, is not part of the Lie groupoid data. Nonetheless, if we regard the Lie $2$-algebra as an \LA -groupoid, $\gg$ corresponds to its core (cf. appendix \ref{appchapter}). \\
It is clear that, by the very definition, for constant $(q,r)$ and $(p,r)$, the referred maps form complexes. For constant $(p,q)$, this is the case as well, although the natural complex has $0$th degree given by $\bigwedge^q\gg_p^*\otimes W$, instead of the $\bigwedge^q\gg_p^*\otimes V$ that is appearing in our complex. This is solved by the following lemma.
\begin{lemma}\label{Alg r-cx}
The map 
\begin{eqnarray*}
\xymatrix{
\delta_{(1)}:\bigwedge^q\gg_p^*\otimes V \ar[r] & \bigwedge^q\gg_p^*\otimes\gg^*\otimes W
}
\end{eqnarray*}
defined above, fits in the complex
\begin{eqnarray*}
\xymatrix{
\bigwedge^q\gg_p^*\otimes V \ar[r]^{\delta_{(1)}\quad} & \bigwedge^q\gg_p^*\otimes\gg^*\otimes W \ar[r] & \bigwedge^q\gg_p^*\otimes\bigwedge^2\gg^*\otimes W \ar[r] & ...
}
\end{eqnarray*}
of $\gg$ with values in $\rho_{(1)}$
\end{lemma}
\begin{proof}
We prove that, for $\omega\bigwedge^q\gg_p^*\otimes V$, $\delta_{(1)}^2=0$. Let $\Xi\in\gg_p^q$ and $x_0,x_1\in\gg$, then
\begin{align*}
\delta_{(1)}\delta_{(1)}\omega(\Xi;x_0,x_1) & =\rho_{(1)}(x_0)\delta_{(1)}\omega(\Xi;x_1)-\rho_{(1)}(x_1)\delta_{(1)}\omega(\Xi;x_0)-\delta_{(1)}\omega(\Xi;[x_0,x_1]) \\
                                            & =\rho_0^1(\mu(x_0))\rho_1(x_1)\omega(\Xi)-\rho_0^1(\mu(x_1))\rho_1(x_0)\omega(\Xi)-\rho_1([x_0,x_1])\omega(\Xi).
\end{align*} 
The result follows from $\rho_0^1(\mu(x_k))=\rho_1(x_k)\circ\phi$, since $\rho_1$ is a Lie algebra homomorphism landing in $\ggl(\phi)_1$.

\end{proof}
The next proposition is the formalization of the heuristic observation that the successive $r$-pages of the triple lattice commute up to the $(r+1)$-page. 
\begin{prop}\label{starTop}
For $r>1$,
\begin{eqnarray*}
\delta^{(r)}\partial-\partial\delta^{(r)}=\delta_{(1)}\circ\Delta+\Delta\circ\delta_{(1)}.
\end{eqnarray*}
where 
\begin{eqnarray*}
\xymatrix{
\Delta:\bigwedge^q\gg_p^*\otimes\bigwedge^{r+1}\gg^*\otimes W \ar[r] & \bigwedge^{q+1}\gg_{p+1}^*\otimes\bigwedge^r\gg^*\otimes W 
} \\
\Delta\omega(\Xi;Z)=\sum_{j=0}^q(-1)^j\omega(\partial_0\Xi(j);x_j^0,Z),\qquad\qquad
\end{eqnarray*}
for $Z\in\gg^r$, $\Xi=(\xi_0,...,\xi_q)\in\gg_{p+1}^{q+1}$ and $\xi_j$ is identified with $(x_j^0,...,x_j^p;y_j)$.
\end{prop}
We already proved the special case $r=1$, which needed a slightly different formula for $\Delta$, one composed with $\phi$.
\begin{proof}
Just as in the proof of proposition \ref{AlgUpToHomotopy r=0} and the argument right after, we compute
\begin{align*}
(\delta^{(r)}\partial-\partial\delta^{(r)})\omega(\Xi ;Z) & =\sum_{j=0}^q(-1)^j\rho^{(r)}(\mu(x_j^0))\omega(\partial_0\Xi(j);Z) \\
                                          & =\sum_{j=0}^q(-1)^j\Big{(}\rho_0^1(\mu(x_j^0))\omega(\partial_0\Xi(j);Z)-\sum_{k=1}^r\omega(\partial_0\Xi(j);z_1,...,[x_j^0,z_k],...,z_r)\Big{)}.
\end{align*}
On the other hand, 
\begin{align*}
\delta_{(1)}\Delta\omega(\Xi;Z) & =\sum_{k=1}^r(-1)^k\rho_0^1(\mu(z_k))\Delta\omega(\Xi;Z(k))+\sum_{m<n}(-1)^{m+n}\Delta\omega(\Xi;[z_m,z_n],Z(m,n)) \\
                                & =\sum_{k=1}^r(-1)^k\rho_0^1(\mu(z_k))\sum_{j=0}^q(-1)^j\omega(\partial_0\Xi(j);x_j^0,Z(k))+ \\
                                & \qquad\qquad +\sum_{m<n}(-1)^{m+n}\sum_{j=0}^q(-1)^j\omega(\partial_0\Xi(j);x_j^0,[z_m,z_n],Z(m,n));
\end{align*}
hence, by introducing the vectors
\begin{eqnarray*}
Z^j=(x_j^0,z_1,...,z_r)\in\gg^{r+1},
\end{eqnarray*}
and rearranging, we have got
\begin{align*}
(\delta^{(r)} & \partial-\partial\delta^{(r)}-\delta_{(1)}\Delta)\omega(\Xi;Z) \\
    & =\sum_{j=0}^q(-1)^j\Big{(}\rho_0^1(\mu(x_j^0))\omega(\partial_0\Xi(j);Z)+\sum_{k=1}^r(-1)^{k+1}\rho_0^1(\mu(z_k))\omega(\partial_0\Xi(j);x_j^0,Z(k))+ \\
    &\qquad\qquad -\sum_{k=1}^r(-1)^{k-1}\omega(\partial_0\Xi(j);[x_j^0,z_k],Z(k))+\sum_{m<n}(-1)^{m+n+2}\omega(\partial_0\Xi(j);[z_m,z_n],x_j^0,Z(m,n))\Big{)} \\
    & =\sum_{j=0}^q(-1)^j\Big{(}\sum_{k=0}^r(-1)^{k}\rho_0^1(\mu(z^j_k))\omega(\partial_0\Xi(j);Z^j(k))+\sum_{m<n}(-1)^{m+n}\omega(\partial_0\Xi(j);[z^j_m,z^j_n],Z^j(m,n))\Big{)} \\
    & =\sum_{j=0}^q(-1)^j\delta_{(1)}\omega(\partial_0\Xi(j);Z^j)=\Delta\delta_{(1)}\omega(\Xi;Z),
\end{align*}
as desired. 

\end{proof}
One would hope that, put together, this data would define a triple complex of sorts, but it does not. If one introduces the obvious indexing by ``counter-diagonal planes'', i.e.
\begin{eqnarray}
C^k_{tot}(\gg_1,\phi):=\bigoplus_{p+q+r=k}C^{p,q}_r(\gg_1,\phi),
\end{eqnarray}
one can try and define a differential $d$ by taking the sum (or the alternated sum) of all the differentials we have described so far. Unfortunately, this does not build a complex, as the fundamental equation $d^2=0$ fails to hold. Let us show this more explicitely. Let $\omega\in C^{p,q}_r(\gg_1,\phi)$, without caring much about signs for the time being, and assuming $d=\partial+\delta^{(r)}+\delta_{(1)}+\Delta$, we have got the following diagram roughly representing the equations that are to vanish:
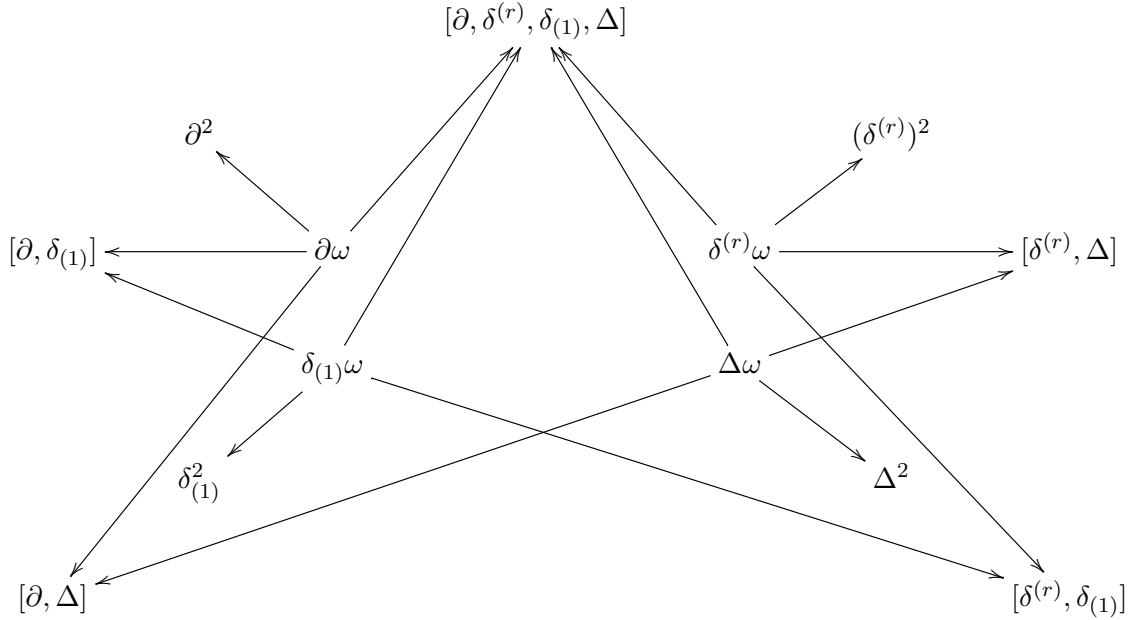
\begin{eqnarray*}
\xymatrix{
  &   &                                     & [\partial,\delta^{(r)},\delta_{(1)},\Delta] &                              &   &   \\
  & \partial^2 &                                               &   &                                                     & (\delta^{(r)})^2 & \\
[\partial,\delta_{(1)}] & & \partial\omega \ar[dddll]\ar[ll]\ar[ul]\ar[uur] & & \delta^{(r)}\omega \ar[uul]\ar[ur]\ar[rr]\ar[dddrr] & & [\delta^{(r)},\Delta] \\
  &   & \delta_{(1)}\omega \ar[dl]\ar[ull]\ar[ruuu]\ar[ddrrrr] &   & \Delta\omega \ar[ddllll]\ar[luuu]\ar[urr]\ar[dr]    &   &   \\
  & \delta_{(1)}^2 &                                           &   &                                                     & \Delta^2 &   \\
[\partial,\Delta] & &                                          &   &                                              & & [\delta^{(r)},\delta_{(1)}]
}
\end{eqnarray*}
We have shown that four of these equations are indeed zero; namely,
\begin{itemize}
    \item $\partial^2\omega=0\in C^{p+2,q}_r(\gg_1,\phi)$, because $\partial$ is the differential of groupoid cochains,
    \item $(\delta^{(r)})^2\omega=0\in C^{p,q+2}_r(\gg_1,\phi)$, because $\delta^{(r)}$ is a Chevalley-Eilenberg differential,
    \item $\delta_{(1)}^2\omega=0\in C^{p,q}_{r+2}(\gg_1,\phi)$, again, because $\delta_{(1)}$ is a Chevalley-Eilenberg differential and
    \item $(\delta^{(r)}\partial-\partial\delta^{(r)}-\delta_{(1)}\circ\Delta-\Delta\circ\delta_{(1)})\omega=0\in C^{p+1.q+1}_r(\gg_1,\phi)$ is the contents of proposition \ref{starTop}.
\end{itemize}
We will devote the rest of this subsection to study the remaining five equations.
\begin{prop}
$(\partial\delta_{(1)}-\delta_{(1)}\partial)\omega=0\in C^{p+1,q}_{r+1}(\gg_1,\phi)$; thus, for constant $q$, there is an actual double complex that we call the $q$-page.
\end{prop}
\begin{proof}
Let $\Xi=(\xi_1,...,\xi^q)\in\gg_{p+1}^q$ and $Z=(z_0,...,z^r)\in\gg^{r+1}$, then
\begin{align*}
\partial\delta_{(1)}\omega(\Xi;Z) & =\sum_{k=0}^{p+1}(-1)^k\delta_{(1)}\omega(\partial_k\Xi;Z) \\
                                  & =\sum_{k=0}^{p+1}(-1)^k\Big{(}\sum_{j=0}^r(-1)^j\rho_0^1(\mu(z_j))\omega(\partial_k\Xi;Z(j))+\sum_{m<n}(-1)^{m+n}\omega(\partial_k\Xi;[z_m,z_n],Z(m,n))\Big{)} \\
                                  & =\sum_{j=0}^r(-1)^j\rho_0^1(\mu(z_j))\sum_{k=0}^{p+1}(-1)^k\omega(\partial_k\Xi;Z(j))+ \\
                                  & \qquad\quad +\sum_{m<n}(-1)^{m+n}\sum_{k=0}^{p+1}(-1)^k\omega(\partial_k\Xi;[z_m,z_n],Z(m,n)) \\
                                  & =\sum_{j=0}^r(-1)^j\rho_0^1(\mu(z_j))\partial\omega(\Xi;Z(j))+\sum_{m<n}(-1)^{m+n}\partial\omega(\Xi;[z_m,z_n],Z(m,n)) =\delta_{(1)}\partial\omega(\Xi;Z)
\end{align*}
\end{proof}
\begin{prop}
$(\delta^{(r)}\delta_{(1)}-\delta^{(r+1)}\delta_{(1)})\omega=0\in C^{p,q+1}_{r+1}(\gg_1,\phi)$; thus, for constant $p$, there is an actual double complex that we call the $p$-page.
\end{prop}
\begin{proof}
Let $\Xi=(\xi_0,...,\xi^q)\in\gg_{p}^{q+1}$ and $Z=(z_0,...,z_r)\in\gg^{r+1}$, then
\begin{align*}
\delta & ^{(r+1)}\delta_{(1)}\omega(\Xi;Z)=\sum_{j=0}^q(-1)^j\rho^{(r+1)}(\hat{t}_p(\xi_j))\delta_{(1)}\omega(\Xi(j);Z)+\sum_{m<n}(-1)^{m+n}\delta_{(1)}\omega([\xi_m,\xi_n],\Xi(m,n);Z) \\
             & =\sum_{j=0}^q(-1)^j\rho^{(r+1)}(\hat{t}_p(\xi_j))\Big{(}\sum_{k=0}^r(-1)^k\rho_0^1(\mu(z_k))\omega(\Xi(j);Z(k))+\sum_{a<b}(-1)^{a+b}\omega(\Xi(j);[z_a,z_b],Z(a,b))\Big{)}+ \\
             & \qquad +\sum_{m<n}(-1)^{m+n}\Big{(}\sum_{k=0}^r(-1)^k\rho_0^1(\mu(z_k))\omega([\xi_m,\xi_n],\Xi(m,n);Z(k))+ \\
             & \qquad\qquad +\sum_{a<b}(-1)^{a+b}\omega([\xi_m,\xi_n],\Xi(m,n);[z_a,z_b],Z(a,b))\Big{)} ;
\end{align*}
whereas, on the other hand,
\begin{align*}
\delta & _{(1)}\delta^{(r)}\omega(\Xi;Z)=\sum_{k=0}^r(-1)^k\rho_0^1(\mu(z_k))\delta^{(r)}\omega(\Xi;Z(k))+\sum_{a<b}(-1)^{a+b}\delta^{(r)}\omega(\Xi;[z_a,z_b],Z(a,b)) \\
             & =\sum_{k=0}^r(-1)^k\rho_0^1(\mu(z_k))\Big{(}\sum_{j=0}^q(-1)^j\rho^{(r)}(\hat{t}_p(\xi_j))\omega(\Xi(j);Z(k))+\sum_{m<n}(-1)^{m+n}\omega([\xi_m,\xi_n],\Xi(m,n);Z(k))\Big{)}+ \\
             & \qquad +\sum_{a<b}(-1)^{a+b}\Big{(}\sum_{j=0}^q(-1)^j\rho^{(r)}(\hat{t}_p(\xi_j))\omega(\Xi(j);[z_a,z_b],Z(a,b))+ \\
             & \qquad\qquad +\sum_{m<n}(-1)^{m+n}\omega([\xi_m,\xi_n],\Xi(m,n);[z_a,z_b],Z(a,b))\Big{)} .
\end{align*}
As a consequence, these expressions will coincide if, and only if 
\begin{eqnarray*}
\rho^{(r+1)}(\hat{t}_p(\xi_j))\rho_0^1(\mu(z_k))\omega(\Xi(j);Z(k))=\rho_0^1(\mu(z_k))\rho^{(r)}(\hat{t}_p(\xi_j))\omega(\Xi(j);Z(k)).
\end{eqnarray*}
Indeed, for all $y\in\hh$, $z\in\gg$ and $X=(x_1,...,x_r)\in\gg^r$, using the fact that $\rho_0^1$ is a Lie algebra homomorphism, we have got
\begin{align*}
 \rho^{(r+1)}(y)\rho_0^1(\mu(z))\omega(\Xi(j);X) & = \rho_0^1(y)\rho_0^1(\mu(z))\omega(\Xi(j);X)-\rho_0^1(\mu(\Lie_yz))\omega(\Xi(j);X) + \\
                                               & \qquad\qquad -\sum_{k=1}^r\rho_0^1(\mu(z))\omega(\Xi(j);x_1,...,\Lie_yx_k,...) \\
                                               & =\rho_0^1(y)\rho_0^1(\mu(z))\omega(\Xi(j);X)-\rho_0^1([y,\mu(z)])\omega(\Xi(j);X) + \\
                                               & \qquad\qquad -\sum_{k=1}^r\rho_0^1(\mu(z))\omega(\Xi(j);x_1,...,\Lie_yx_k,...) \\
                                               & =\rho_0^1(\mu(z))\rho_0^1(y)\omega(\Xi(j);X)-\rho_0^1(\mu(z))\sum_{k=1}^r\omega(\Xi(j);x_1,...,\Lie_yx_k,...) \\
                                               & =\rho_0^1(\mu(z))\rho^{(r)}(y)\omega(\Xi(j);X),
\end{align*}
so the result follows.

\end{proof}
\begin{prop}\label{partialDelta}
$(\partial\Delta+\Delta\partial)\omega=0\in C^{p+2,q+1}_{r-1}(\gg_1,\phi)$.
\end{prop}
\begin{proof}
Let $Z\in\gg^{r+1}$, $\Xi=(\xi_0,...,\xi^q)\in\gg_{p+2}^{q+1}$ and let $\xi_j$ be identified with $(x_j^0,...,x_j^{p+1};y_j)$. Computing,
\begin{align*}
(\partial\Delta+\Delta\partial)\omega(\Xi;Z) & =\sum_{j=0}^q(-1)^j\partial\omega(\partial_0\Xi(j);x_j^0,Z)+\sum_{k=0}^{p+2}(-1)^k\Delta\omega(\partial_k\Xi;Z) \\
             & =\sum_{j=0}^q(-1)^j\sum_{k=0}^{p+1}(-1)^k\omega(\partial_k(\partial_0\Xi(j));x_j^0,Z)+\sum_{j=0}^q(-1)^j\omega(\partial_0(\partial_0\Xi)(j);x_j^1,Z) \\
             & \quad -\sum_{j=0}^q(-1)^j\omega(\partial_0(\partial_1\Xi)(j);x_j^0+x_j^1,Z)+\sum_{k=2}^{p+2}(-1)^k\sum_{j=0}^q(-1)^j\omega(\partial_0(\partial_k\Xi)(j);x_j^0,Z) \\
             & =\sum_{j=0}^q(-1)^j\Big{(}\sum_{k=0}^{p+1}(-1)^k\omega(\partial_k\partial_0\Xi(j);x_j^0,Z)-\omega(\partial_0(\partial_0\Xi)(j);x_j^0,Z) \\
             & \qquad\qquad +\sum_{k=2}^{p+2}(-1)^k\omega(\partial_0(\partial_k\Xi)(j);x_j^0,Z)\Big{)}.
\end{align*}
Now, either computing, or from the simplicial identities, we know that
\begin{eqnarray*}
\partial_k\partial_0(\xi_j)\sim
  \begin{cases}
    (x_j^2,...,x_j^{p+1};y_j)                       & \quad \text{if } k=0  \\
    (x_j^1,...,x_j^{k}+x_j^{k+1},...,x_j^{p+1};y_j) & \quad \text{if } 0<k\leq p\\
    (x_j^1,...,x_j^{p};y_j+\mu(x_j^{p+1}))          & \quad \text{if } k=p+1, 
  \end{cases}
\end{eqnarray*}
and
\begin{eqnarray*}
\partial_0\partial_k(\xi_j)\sim
  \begin{cases}
    (x_j^2,...,x_j^{p+1};y_j)                     & \quad \text{if } k\in\lbrace 0,1\rbrace  \\
    (x_j^1,...,x_j^k+x_j^{k+1},...,x_j^{p+1};y_j) & \quad \text{if } 1<k\leq p+1\\
    (x_j^1,...,x_j^{p};y_j+\mu(x_j^{p+1}))        & \quad \text{if } k=p+2; 
  \end{cases}
\end{eqnarray*}
therefore, the sums cancel one another and the result follows.

\end{proof}
Up until now, we were able to prove the vanishing of the equations in the star diagram above, but this is as good as it gets. The remaining two equations do not vanish in general, but there is a way out. 
\subsubsection{Higher difference maps}
We start by studying how the relation labeled $[\delta^{(r)},\Delta]$ behaves in the special case $r=1$. \\
Let $\alpha\in\bigwedge^q\gg_p^*\otimes\gg^*\otimes W$, $\Xi=(\xi_0,...,\xi_{q+1})\in\gg_{p+1}^{q+2}$ with the by now usual identification $\xi_j\sim(x_j^0,...,x_j^p;y_j)$, and let us compute separately $\delta\Delta\alpha$ and $\Delta\delta'\alpha$. On the one hand, we have got
\begin{align*}
\delta\Delta\alpha(\Xi) & =\sum_{k=0}^{q+1}(-1)^k\rho_0^0(\hat{t}_{p+1}(\xi_k))\Delta\alpha(\Xi(k))+\sum_{m<n}(-1)^{m+n}\Delta\alpha([\xi_m,\xi_n],\Xi(m,n)) .
\end{align*}
The first term in this expression is computed to be
\begin{align*}
 & \sum_{k=0}^{q+1}(-1)^k\rho_0^0(\hat{t}_{p+1}(\xi_k))\phi\Big{[}\sum_{j=0}^{k-1}(-1)^j\alpha(\partial_0\Xi(j,k);x_j^0)+\sum_{j=k+1}^{q+1}(-1)^{j+1}\alpha(\partial_0\Xi(k,j);x_j^0)\Big{]}.
\end{align*}
The second term is trickier though. Carefully computing it, one realizes it coincides with
\begin{align*}
 & \sum_{m<n}(-1)^{m+n}\phi\Big{(}\alpha(\partial_0\Xi(m,n);[x_m^0,x_n^0]+\Lie_{\hat{t}_p(\partial_0\xi_m)}x_n^0-\Lie_{\hat{t}_p(\partial_0\xi_n)}x_m^0)+ \\
 & \qquad +\sum_{j=0}^{m-1}(-1)^{j+1}\alpha(\partial_0([\xi_m,\xi_n]),\partial_0\Xi(j,m,n);x_j^0)+\sum_{j=m+1}^{n-1}(-1)^j\alpha(\partial_0([\xi_m,\xi_n]),\partial_0\Xi(m,j,n);x_j^0)+ \\
 & \qquad\qquad +\sum_{j=n+1}^q(-1)^{j+1}\alpha(\partial_0([\xi_m,\xi_n]),\partial_0\Xi(m,n,j);x_j^0)\Big{)} 
\end{align*}
The shape of the first term comes from the definition of the Lie algebra structure on $\gg_{p+1}$. Recall that the bracket is inherited from the product $\gg^{p+1}$; therefore, the first entry of $[\xi_m,\xi_n]$ is
\begin{align*}
\Big{[}\big{(}x_m^0,y_m & +\sum_{k=1}^p\mu(x_m^k)\big{)},\big{(}x_n^0,y_n+\sum_{k=1}^p\mu(x_n^k)\big{)}\Big{]}_1 \\
 & =\Big{(}[x_m^0,x_n^0]+\Lie_{y_m+\sum_{k=1}^p\mu(x_m^k)}x_n^0-\Lie_{y_n+\sum_{k=1}^p\mu(x_n^k)}x_m^0,\big{[}y_m+\sum_{k=1}^p\mu(x_m^k),y_n+\sum_{k=1}^p\mu(x_n^k)\big{]}\Big{)},
\end{align*}
and we already saw that $\hat{t}_p(\partial_0\xi_j)=y_j+\sum_{k=1}^p\mu(x_j^k)$. The second term is written in such a way that it has got the appropriate signs when written with our conventions. On the other hand, 
\begin{align*}
\Delta\delta'\alpha(\Xi) & =\phi\Big{(}\sum_{j=0}^q(-1)^j\delta'\alpha(\partial_0\Xi(j);x_j^0)\Big{)}.
\end{align*}
 After evaluating the differential $\delta'$, there will be two terms. The first will be
\begin{align*}         
\sum_{j=0}^{q+1}(-1)^j &  \phi\Big{(}\sum_{k=0}^{j-1}(-1)^k\rho'(\hat{t}_p(\partial_0\xi_k))\alpha(\partial_0\Xi(k,j);x_j^0)+\sum_{k=j+1}^q(-1)^{k+1}\rho'(\hat{t}_p(\partial_0\xi_k))\alpha(\partial_0\Xi(j,k);x_j^0)\Big{)}
\end{align*}
while the second,
\begin{align*}
\sum_{j=0}^{q+1}(-1)^j & \phi\Big{(}\sum_{m=0}^{j-1}\Big{[}\sum_{n=m+1}^{j-1}(-1)^{n}\alpha([\partial_0\xi_m,\partial_0\xi_n]),\partial_0\Xi(m,n,j);x_j^0)+ \\
 & \qquad +\sum_{n=j+1}^{q+1}(-1)^{n+1}\alpha([\partial_0\xi_m,\partial_0(\xi_n]),\partial_0\Xi(m,j,n);x_j^0)\Big{]}+ \\
 & \qquad\qquad +\sum_{m=j+1}^{q+1}\sum_{n=m+1}^{q+1}(-1)^{n}\alpha([\partial_0\xi_m,\partial_0\xi_n],\partial_0\Xi(m,n,j);x_j^0)\Big{)}. 
\end{align*}
When considering the sum of the whole of the two expressions, this second term will cancel out with the second part of the second term above simply due to the fact that $\partial_0$ is a Lie algebra homomorphism. Furthermore, since $\rho_0^0(y)\phi=\phi\rho_0^1(y)$ for all $y\in\hh$, we will have again the difference 
\begin{eqnarray*}
\rho_0^1(\hat{t}_{p}(\xi_k))-\rho_0^1(\hat{t}_p(\partial_0\xi_k))=\rho_0^1(\mu(x_k^0)),
\end{eqnarray*} 
and in the meantime, the second terms of the representation $\rho'$,
\begin{eqnarray*}
\alpha(\partial_0\Xi(m,n);\Lie_{\hat{t}_p(\partial_0\xi_m)}x_n^0)
\end{eqnarray*}
will cancel the ones coming from the Lie bracket in $\gg_{p+1}$. Thus, we are left behind with a series of terms of the form
\begin{eqnarray*}
\rho_0^1(\mu(x_k))\alpha(\partial_0\Xi(j,k);x_j^0)
\end{eqnarray*}
and another series of terms of the form
\begin{eqnarray*}
\alpha(\partial_0\Xi(m,n);[x_m^0,x_n^0]).
\end{eqnarray*}
What is remarkable is that these come with the right signs to build up differentials in the image of $\delta_{(1)}$. In fact, the preceding discussion can be summarized in the following lemma.
\begin{lemma}
Let 
\begin{eqnarray*}
\xymatrix{
\Delta_2:\bigwedge^q\gg_p^*\otimes\bigwedge^2\gg^*\otimes W \ar[r] & \bigwedge^{q+2}\gg_{p+1}^*\otimes V
}
\end{eqnarray*}
be defined by
\begin{eqnarray*}
\Delta_2\alpha(\Xi):=\phi\Big{(}\sum_{m<n}(-1)^{m+n}\alpha(\partial_0\Xi(m,n);x_m^0,x_n^0)\Big{)}.
\end{eqnarray*}
Then
\begin{eqnarray*}
\delta\Delta+\Delta\delta'=\Delta_2\circ\delta_{(1)}.
\end{eqnarray*}
\end{lemma}
These second difference maps will help us simultaneously decide the behaviour of the missing equations, as shown by the following proposition.
\begin{prop}
For $r>2$, let
\begin{eqnarray*}
\xymatrix{
\Delta_2:\bigwedge^q\gg_p^*\otimes\bigwedge^r\gg^*\otimes W \ar[r] & \bigwedge^{q+2}\gg_{p+1}^*\otimes\bigwedge^{r-2}\gg^*\otimes W 
}
\end{eqnarray*}
be defined by
\begin{eqnarray*}
\Delta_2\omega(\Xi;Z):=\sum_{m<n}(-1)^{m+n}\omega(\partial_0\Xi(m,n);x_m^0,x_n^0,Z).
\end{eqnarray*}
Then
\begin{eqnarray*}
\delta^{(r-1)}\circ\Delta+\Delta\circ\delta^{(r)}=\delta_{(1)}\circ\Delta_2-\Delta_2\circ\delta_{(1)};
\end{eqnarray*}
furthermore, 
\begin{eqnarray*}
\Delta^2=-(\Delta_2\circ\partial+\partial\circ\Delta_2).
\end{eqnarray*}
\end{prop}
Before diving into the proof, we would like to remark that although this proposition settles the behaviour of the missing points of the star shaped diagram above, by introducing the new differential $\Delta_2$, there appear a number of new equations that do not vanish either. This failure will lead us to define the higher difference maps. With them, we will ultimately get a differential that will turn the three-dimensional lattice into a complex.
\begin{proof}
Let $Z=(z_1,...,z_{r-1})\in\gg^{r-1}$, $\Xi=(\xi_0,...,\xi_{q+1})\in\gg_{p+1}^{q+2}$ and identify $\xi_j$ with $(x_j^0,...,x_j^p;y_j)$. We write the expressions for the terms of the left hand side of the first equation
\begin{align*}
\delta^{(r-1)}\Delta\omega(\Xi;Z) & =\sum_{k=0}^{q+1}(-1)^k\rho^{(r-1)}(\hat{t}_{p+1}(\xi_k))\Delta(\Xi(k);Z)+\sum_{m<n}(-1)^{m+n}\Delta([\xi_m,\xi_n],\Xi(k);Z), \\
\Delta\delta^{(r)}\omega(\Xi;Z)   & =\sum_{j=0}^{q+1}(-1)^j\delta^{(r)}\omega(\partial_0\Xi(j);x_j^0,Z).
\end{align*}
For each pair $(j,k)$, if $j>k$, there are terms 
\begin{eqnarray*}
(-1)^k\rho^{(r-1)}(\hat{t}_{p+1}(\xi_k))\big{[}(-1)^{j+1}\omega(\partial_0\Xi(k,j);x_j^0,Z)\big{]} & \textnormal{and} &
(-1)^{j}(-1)^k\rho^{(r)}(\hat{t}_p(\partial_0\xi_k))\omega(\partial_0\Xi(k,j);x_j^0,Z)
\end{eqnarray*}
in each of the equations above respectively. If, on the other hand, $j<k$, one has got the following terms: 
\begin{eqnarray*}
(-1)^k\rho^{(r-1)}(\hat{t}_{p+1}(\xi_k))\big{[}(-1)^{j}\omega(\partial_0\Xi(j,k);x_j^0,Z)\big{]} & \textnormal{and} &
(-1)^{j}(-1)^{k+1}\rho^{(r)}(\hat{t}_p(\partial_0\xi_k))\omega(\partial_0\Xi(j,k);x_j^0,Z).
\end{eqnarray*}
We write the sums of these pairs of terms, 
\begin{align*}
\sum_{\tau\in S_2} & (-1)^{(j+k)}  \abs{\tau}\big{[}\rho^{(r-1)}(\hat{t}_{p+1}(\xi_k))-\rho^{(r)}(\hat{t}_p(\partial_0\xi_k))\big{]}\omega(\partial_0\Xi(\tau(k,j));x_j^0,Z) \\
             & =\sum_{\tau\in S_2}(-1)^{(j+k+1)}\abs{\tau}\big{[}\omega(\partial_0\Xi(\tau(k,j));\Lie_{\hat{t}_p(\partial_0\xi_k)}x_j^0,Z)+ \\
             & \qquad +\rho_0^1(\mu(x_k^0))\omega(\partial_0\Xi(\tau(k,j));x_j^0,Z)-\sum_{i=1}^{r-1}\omega(\partial_0\Xi(\tau(k,j));x_j^0,z_1,...,[x_k^0,x_i],...,x_{r-1})\big{]}.
\end{align*}
The second sum $\delta^{(r-1)}\Delta\omega(\Xi;Z)$, when expanding $\Delta$ will have first term 
\begin{eqnarray*}
\sum_{m<n}(-1)^{m+n}\omega(\Xi(m,n);[x_m^0,x_n^0]+\Lie_{\hat{t}_p(\partial_0\xi_m)}x_n^0-\Lie_{\hat{t}_p(\partial_0\xi_n)}x_m^0,Z).
\end{eqnarray*}
All the remaining terms will cancel out in a similar fashion to the discussion preceding the statement of the theorem. In so, 
\begin{align*}
(\delta^{(r-1)}\Delta & +\Delta\delta^{(r)})\omega(\Xi; Z) \\
  & =-\sum_{m<n}(-1)^{m+n}\Big{[}\rho_0^1(\mu(x_m^0))\omega(\partial_0\Xi(m,n);x_n^0,Z)-\rho_0^1(\mu(x_m^0))\omega(\partial_0\Xi(m,n);x_n^0,Z)+ \\
  & \qquad -\omega(\partial_0\Xi(m,n);[x_m^0,x_n^0],Z)+ \\
  & \qquad\quad+\sum_{k=1}^{r-1}(-1)^{k+1}\big{(}\omega(\partial_0\Xi(m,n);[x_m^0,z_k],x_n^0,Z(k))-\omega(\partial_0\Xi(m,n);[x_n^0,z_k],x_n^0,Z(k))\big{)}\Big{]}
\end{align*}
Now, 
\begin{align*}
\delta_{(1)}\Delta_2\omega(\Xi;Z) & =\sum_{k=1}^{r-1}(-1)^{k+1}\rho_0^1(\mu(z_k))\Delta_2\omega(\Xi;Z(k))+\sum_{a<b}(-1)^{a+b}\Delta_2\omega(\Xi;[z_a,z_b],Z(a,b)) \\
                                  & =\sum_{k=1}^{r-1}(-1)^{k+1}\sum_{m<n}(-1)^{m+n}\rho_0^1(\mu(z_k))\omega(\partial_0\Xi(m,n);x_m^0,x_n^0,Z(k))+ \\
                                  & \qquad\quad +\sum_{a<b}(-1)^{a+b}\sum_{m<n}(-1)^{m+n}\omega(\partial_0\Xi(m,n);x_m^0,x_n^0,[z_a,z_b],Z(a,b)) ;
\end{align*}
thus, one sees that 
\begin{eqnarray*}
(\delta^{(r-1)}\Delta+\Delta\delta^{(r)}-\delta_{(1)}\Delta_2)\omega(\Xi; Z)=-\Delta_2\delta_{(1)}\omega(\Xi; Z),
\end{eqnarray*}
and the first equation of the statement holds. \\
As for the second equation, let $\Xi=(\xi_0,...,\xi_{q+1})\in\gg_{p+2}^{q+2}$ with $\xi_j\sim(x_j^0,...,x_j^{p+1};y_j)$ and $Z$ as before. Computing,
\begin{align*}
\Delta_2\partial\omega(\Xi;Z) & =\sum_{m<n}(-1)^{m+n}\partial\omega(\partial_0\Xi(m,n);x_m^0,x_n^0,Z) \\
                              & =\sum_{m<n}(-1)^{m+n}\sum_{k=0}^{p+1}(-1)^k\omega(\partial_k\partial_0\Xi(m,n);x_m^0,x_n^0,Z);
\end{align*}
whereas, on the other hand, 
\begin{align*}
\partial\Delta_2\omega(\Xi;Z) & =\sum_{k=0}^{p+2}(-1)^k\Delta_2\omega(\partial_k\Xi;Z) \\
                              & = \sum_{m<n}(-1)^{m+n}\Big{(}\omega(\partial_0\partial_0\Xi(m,n);x_m^1,x_n^1,Z)-\omega(\partial_0\partial_1\Xi(m,n);x_m^0+x_m^1,x_n^0+x_n^1,Z)\Big{)}+ \\
                              & \qquad +\sum_{k=2}^{p+2}(-1)^k\sum_{m<n}(-1)^{m+n}\omega(\partial_0\partial_k\Xi(m,n);x_m^0,x_n^0,Z).
\end{align*}
Due to the simplicial identities that we outlined in the proof of proposition \ref{partialDelta}, we know that $\partial_0\partial_0=\partial_0\partial_1$ so we can define
\begin{eqnarray*}
\varpi_{mn}(w_1,w_2):=\omega(\partial_0\partial_0\Xi(m,n);w_1,w_2,Z)
\end{eqnarray*}
and write the sum
\begin{align*}
(\Delta_2\partial+\partial\Delta_2)\omega(\Xi;Z) & =\sum_{m<n}(-1)^{m+n}\Big{(}\varpi(x_m^1,x_n^1)+\varpi(x_m^1,x_n^1)-\varpi(x_m^0+x_m^1,x_n^0+x_n^1)\Big{)} \\
                                                 & =\sum_{m<n}(-1)^{m+n}\Big{(}\varpi(x_n^1,x_m^0)-\varpi(x_m^1,x_n^0)\Big{)}
\end{align*}
which is precisely $-\Delta^2(\Xi;Z)$. %If doubtful, see 11.03.2018. The blue print.

\end{proof}
Trying and defining the differential for the complex $C(\gg_1,\phi)$ as
\begin{eqnarray*}
d=\delta^{(r)}+\delta_{(1)}+\partial+\Delta+\Delta_{2}
\end{eqnarray*}
up to signs, will still fail to verify $d^2=0$, as it was stated before the proof of the latter proposition. Inductively, studying how the $(k-1)$th difference map $\Delta_{k-1}$ fails to commute with $\delta^{(r)}$, one finds a $k$th difference map  
\begin{eqnarray*}
\xymatrix{
\Delta_{k}:C^{p,q}_r(\gg _1,\phi) \ar[r] & C^{p+1,q+k}_{r-k}(\gg _1,\phi).
}
\end{eqnarray*}
These maps will be defined in general as
\begin{eqnarray*}
\Delta_k\omega(\Xi;Z):=\sum_{a_1<...<a_k}(-1)^{a_1+...+a_k}\omega(\partial_0\Xi(a_1,...,a_k);x_{a_1}^0,...,x_{a_k}^0,Z)
\end{eqnarray*}
for $Z\in\gg^{r-k}$ and $\Xi$ . In the special case $r=k$, the map is essentially defined by the same formula, but composed with $\phi$, so that it takes values in the right vector space. With these maps we get the main theorems of this chapter.
\begin{theorem}\label{algDiffs}
Thus defined, the maps $\Delta_k$ verify the following set of equations for elements in $C^{p,q}_r(\gg_1,\phi)$
\begin{itemize}
    \item[i)] $\delta^{(r-k)}\Delta_k+\Delta_k\delta^{(r)}=\delta_{(1)}\Delta_{k+1}+\Delta_{k+1}\delta_{(1)}$ for all $1\leq k<r$,
    \item[ii)]  $\delta\Delta_r+\Delta_r\delta^{(r)}=\Delta_{r+1}\delta_{(1)}$ and
    \item[iii)] assuming the convention that $\Delta_0=\partial$,
    \begin{eqnarray*}
    \sum_{i=0}^k\Delta_{k-i}\Delta_i=0
    \end{eqnarray*}
    for all $1\leq k\leq r$.
\end{itemize}
\end{theorem}
Notice that with respect to the diagonal grading, all difference maps are homogeneous of degree $+1$. Then, we can use them as differentials and indeed, we have got the following.
\begin{theorem}\label{The2AlgCx}
$C(\gg_1,\phi)$ graded by
\begin{eqnarray*}
C^{n}_{tot}(\gg_1,\phi)=\bigoplus_{p+q+r=n}C^{p,q}_r(\gg_1,\phi)
\end{eqnarray*}
together with the differential
\begin{eqnarray*}
\nabla=\delta^{(r)}+(-1)^q\delta_{(1)}+(-1)^{q+r}\partial+(-1)^r\sum_k\Delta_k
\end{eqnarray*}
is a complex.
\end{theorem}
\begin{remark}
Recall that, at the beginning of the section, we set out to look for a double complex and the non-commutativity of the squares forced us to consider the successive $r$-pages. As it turns out, in sight of theorem \ref{algDiffs} and of their grading, we can think of the difference maps as being maps of double complexes between the $p$-pages. Doing so, one builds an honest double complex, whose columns are the total complexes of the $p$-pages and whose horizontal differentials are the maps of complexes induced by the difference maps. The total complex of this double complex coincides with that of theorem \ref{The2AlgCx}. 
\end{remark}

%------------------------------

\section[2-cohomology of Lie 2-algebras]{2-cohomology of Lie 2-algebras}\label{2AlgCoh}
\sectionmark{Cohomology}
%---------------------------------------------------------------------------------------------------------------------------------------------------

In this section, we study the cohomology of the complex from theorem \ref{The2AlgCx}. \\
Start out with an element $v\in V=C^0(\gg_1,\phi)=C^{0,0}_0(\gg_1,\phi)$, then its differential is given by
\begin{eqnarray*}
\nabla v=(\delta v,\delta_{(1)}v,\cancelto{0}{\partial v})\in C^1(\gg_1,\phi)
\end{eqnarray*}
If $v$ is a $0$-cocycle, then for all $y\in\hh$ and for all $x\in\gg$,
\begin{eqnarray*}
\rho_0^0(y)v=0 & \textnormal{and} & \rho_1(x)v=0;
\end{eqnarray*}
therefore,
\begin{eqnarray*}
H^0_\nabla(\gg_1,\phi)=V^{\gg_1}:=\lbrace v\in V:\bar{\rho}_{(x,y)}(0,v)=0,\quad\forall (x,y)\in\gg\oplus_\Lie\hh\rbrace ,
\end{eqnarray*}
where $\bar{\rho}$ is the honest representation of proposition \ref{honestAlgRep}. \\
A $1$-cochain $\lambda$ is a triple $(\lambda_0,\lambda_1,v)\in(\hh^*\otimes V)\oplus(\gg^*\otimes W)\oplus V$, and its differential has six entries that we write using their coordinates, i.e. $\nabla\lambda^{p,q}_r\in C^{p,q}_r(\gg _1,\phi)$
\begin{align*}
\nabla\lambda^{0,2}_0 & =\delta\lambda_0                             &                       & & &                               \\
\nabla\lambda^{1,1}_0 & =-\partial\lambda_0+\delta v-\Delta\lambda_1 & \nabla\lambda^{0,1}_1 & =-\delta_{(1)}\lambda_0+\delta'\lambda_1 & & \\
\nabla\lambda^{2,0}_0 & =\partial v=v                                & \nabla\lambda^{1,0}_1 & =\delta_{(1)}v-\partial\lambda _1=\delta_{(1)}v & \nabla\lambda^{0,0}_2 & =\delta_{(1)}\lambda_1 
\end{align*}
Schematically,
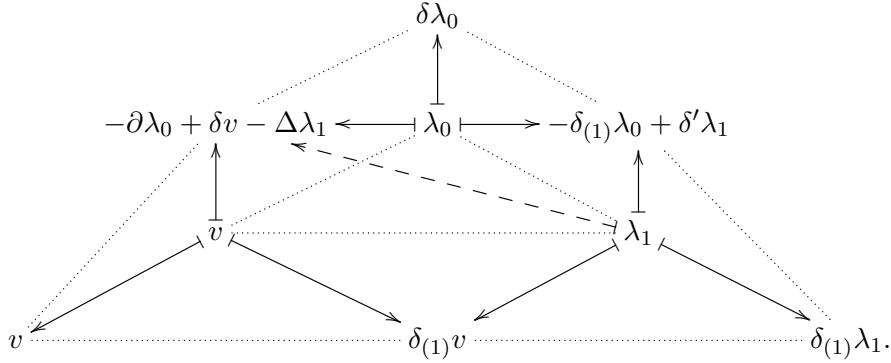
\begin{eqnarray*}
\xymatrix{
  & & \delta\lambda_0 \ar@{.}[dl]\ar@{.}[dr] & &  \\
  & -\partial\lambda_0+\delta v-\Delta\lambda_1 & \lambda_0 \ar@{|->}[u]\ar@{|->}[r]\ar@{|->}[l]\ar@{.}[dr] & -\delta_{(1)}\lambda_0+\delta'\lambda_1 & \\
  & v \ar@{|->}[dl]\ar@{|->}[u]\ar@{|->}[dr]\ar@{.}[ur]\ar@{.}[rr] &  & \lambda_1 \ar@{|->}[dr]\ar@{|->}[u]\ar@{|->}[dl]\ar@{|-->}[ull] & \\
v \ar@{.}[uur] \ar@{.}[rr] & & \delta_{(1)}v \ar@{.}[rr] & & \delta_{(1)}\lambda_1 .\ar@{.}[uul]
}
\end{eqnarray*}
Here the solid arrows represent the differentials, the dashed arrow represents the difference map, and the pointed polygons represent elements of the same degree. If $\lambda$ is a $1$-cocylce, then $v=0$, $(\lambda_0,\lambda_1)\in Der(\hh,V)\oplus Der(\gg,W)$ and the following relations hold for all pairs $(y,x)\in\hh\times\gg$:
\begin{eqnarray*}
\phi(\lambda_1(x))=\lambda_0(\mu(x)), \\
\lambda_1(\Lie_y x)=\rho_1(x)\lambda_0(y)-\rho_0^1(y)\lambda_1(x).
\end{eqnarray*}
The first of these relations says that 
\begin{eqnarray*}
\xymatrix{
\bar{\lambda}:\gg\oplus\hh \ar[r] & W\oplus V:(x,y) \ar@{|->}[r] & (\lambda_1(x),\lambda_0(y))
}
\end{eqnarray*}
is a map of $2$-vector spaces. The second, says that $\bar{\lambda}\in Der(\gg\oplus_\Lie\hh,W\oplus V)$ with respect to $\bar{\rho}$. Notice that since $\nabla v$ can also be seen as $\bar{\rho}_{(x,y)}(0,v)$; by analogy, we can define:
\begin{Def}
The \textit{space of derivations} of a Lie $2$-algebra $\gg_1$ with respect to a $2$-representation $\rho$ on the $2$-vector $\mathbb{V}=\xymatrix{W \ar[r]^\phi & V}$ is defined to be
\begin{eqnarray*}
Der(\gg_1,\phi):=\lbrace (\lambda_1,\lambda_0)\in Hom_{2-Vect}(\gg_1,\mathbb{V}):(\lambda_1,\lambda_0)\in Der(\gg\oplus_\Lie\hh,W\oplus V)\rbrace .
\end{eqnarray*}
The space of \textit{inner derivations}, on the other hand, is defined to be
\begin{eqnarray*}
Inn(\gg_1,\phi):=\lbrace (\lambda_1,\lambda_0)\in Der(\gg_1,\phi):(\lambda_1(x),\lambda_0(y))=\bar{\rho}_{(x,y)}(0,v)\textnormal{ for some }v\in V\rbrace .
\end{eqnarray*}
\end{Def}
With these definitions, we can write
\begin{eqnarray*}
H^1_\nabla(\gg_1,\phi)=Out(\gg_1,\phi):=Der(\gg_1,\phi)/Inn(\gg_1,\phi).
\end{eqnarray*}
A $2$-cochain $\vec{\omega}$ is a $6$-tuple $(\omega_0,\alpha,\varphi,\omega_1,\lambda,v)$ where
\begin{align*}
\omega_0 & \in\bigwedge^2\hh^*\otimes V &          &                               &          &   \\
\varphi  & \in\gg_1^*\otimes V          & \alpha   & \in\hh^*\otimes\gg^*\otimes W &          &   \\
v        & \in V                        & \lambda  & \in\gg^*\otimes W             & \omega_1 & \in\bigwedge^2\gg^*\otimes W . 
\end{align*} 
Let us compute the differential using coordinates as before. First, $\nabla\vec{\omega}^{3,0}_0=\partial v=0$; for the other values,
\begin{align*}
\nabla\vec{\omega}^{0,3}_0 & =\delta\omega_0 & & & &  \\
\nabla\vec{\omega}^{1,2}_0 & =\partial\omega_0+\delta\varphi-\Delta\alpha+\Delta_2\omega_1 & \nabla\vec{\omega}^{0,2}_1 & =\delta_{(1)}\omega_0+\delta'\alpha & & \\
\nabla\vec{\omega}^{2,1}_0 & =-\partial\varphi+\delta v-\Delta\lambda & \nabla\vec{\omega}^{1,1}_1 & =-\delta_{(1)}\varphi+\partial\alpha+\delta'\lambda+\Delta\omega_1 & \nabla\vec{\omega}^{0,1}_2 & =-\delta_{(1)}\alpha+\delta^{(2)}\omega_1 \\
\nabla\vec{\omega}^{2,0}_1 & =\delta_{(1)}v-\cancelto{\lambda}{\partial\lambda} & \nabla\vec{\omega}^{1,0}_2 & =\cancelto{0}{\partial\omega_1}+\delta_{(1)}\lambda & \nabla\vec{\omega}^{0,0}_3 & =\delta_{(1)}\omega_1
\end{align*}
If we assume  $v=0$, $\lambda$ vanishes too and these equations correspond to those in the statement of proposition \ref{2-cocycles}. First, $\partial\varphi(x_0,x_1,y)=\varphi(0,y+\mu(x_1))=0$ says that $\varphi$ does not depend on $\hh$, as was supposed in the statement. Then, the equation $\nabla\vec{\omega}^{1,1}_1=0$ coincides with the definition of $\omega_1$ in the statement. Finally, for the numbering of the equations in the statement of proposition \ref{2-cocycles}:
\begin{itemize}
\item Equation $i)$ is literally $\delta\omega_0 =0$.  
\item Equation $ii)$ just said that $\omega_1$ was skew-symmetric, so we've got it already, as our $\omega_1$ is alternating by definition.
\item Equation $iii)$ said that $\omega_1$ was a cocycle, so it follows from $\nabla\vec{\omega}^{0,0}_3=0$.
\item Equation $iv)$ is equivalent to $\nabla\vec{\omega}^{1,2}_0=0$ evaluated at $\begin{pmatrix}
0 & x \\
y & 0
\end{pmatrix}\in\gg_1^2$. % We're cheating a little bit. There is a sign missing: On \Delta\alpha. (not anymore)
\item Equation $v)$ is exactly $\nabla\vec{\omega}^{0,2}_1=0$.
\item Equation $vi)$ is exactly $\nabla\vec{\omega}^{0,1}_2=0$.
\end{itemize}
Moreover, for a pair of $2$-cocycles 
\begin{eqnarray*}
(\vec{\omega})^k=(\omega_0^k,\alpha^k,\varphi^k,\omega_1^k,\lambda^k,v^k), & k\in\lbrace 1,2\rbrace
\end{eqnarray*}
with $v^k=0$, if they are cohomologous we recover the equations in proposition \ref{2-coboundaries}. Indeed, if
\begin{eqnarray*}
(\vec{\omega})^2-(\vec{\omega})^1=\nabla(\lambda_0,\lambda_1,v),
\end{eqnarray*}
coordinate-wise we've got
\begin{align*}
\omega_0^2-\omega_0^1 & =\delta\lambda_0                             &                   & & &                               \\
\varphi^2-\varphi^1   & =-\partial\lambda_0+\delta v-\Delta\lambda_1 & \alpha^2-\alpha^1 & =-\delta_{(1)}\lambda_0+\delta'\lambda_1 & & \\
0                     & =v                                           & 0                 & =\delta_{(1)}v                           & \omega_1^2-\omega_1^1 &  =\delta_{(1)}\lambda_1 ;
\end{align*}
thus, explicitely 
\begin{align*}
    (\varphi^2-\varphi^1)(x) & =\lambda_0(\mu(x))-\phi(\lambda_1(x)) \\
    (\alpha^2-\alpha^1)(y;x) & =-\rho_1(x)\lambda_0(y)+\rho_0^1(y)\lambda_1(x)-\lambda_1(\Lie_y x)
\end{align*}
which coincide with the referred equations, as we claimed. \\
In the end, this is the theorem this theory is dedicated to. 
\begin{theorem}\label{H2Alg}
$H^2_\nabla(\gg_1,\phi)$ classifies $2$-extensions.
\end{theorem}

% We hope to use this algebraic description to prove a general van Est theorem for Lie $2$-algebras. Next step is to see how far this theory can be stretched to include \LA -groups. The issue of defining a representation seems already complicated, as there is no ``flat abelian'' candidate in the quotient category of \LA -groups over a prescribed Lie group.
%***Of the things that this section lacks, but would be nice for it to have.3) Particular cases. The second cohomology of the ideal with values in a trivial representation, in a representation with values in the pair groupoid.4) All too hopeful, does Hochschild-Serre shows up? Look for it! 

%% ------------------------------------------------------------------------- %%
\chapter{Comomology: the Lie 2-group theory}\label{fractionalchapter2}
\chaptermark{Comomology: the Lie $2$-group theory}
\section{Introduction}
%---------------------------------------------------------------------------------------------------------------------------------------------------

In this chapter, we will introduce the cohomology of Lie $2$-groups. We start by looking at the total cohomology of the double complex of the Lie $2$-group to discover that it classifies certain type of extensions. With the goal of classifying extensions in sight, we define the notion of a representation of a Lie $2$-group. We move on to study split extensions by picking a splitting and looking for conditions for naturally defined maps to build an extension back up. We use these ``cocycle equations'' to read a complex associated to the Lie $2$-groups with values in a representation. Finally, we give an interpretation for the lower dimensions of the cohomology of this complex and verify that its second cohomology indeed classifies extensions as prescribed.

% ------------------------------------------
% ------------------------------------------
% LIE 2-GROUP THEORY
% ------------------------------------------
%-------------------------------------------

%1- Introduction - Outline of the main problem and motivations.
Along this section, we repeat the run-the-mill proofs that we had for Lie $2$-algebras in the case of their global counter-parts.

\section{2-cohomology with trivial coefficients}\label{GpDcx}
Let $\xymatrix{\G \ar@<0.5ex>[r] \ar@<-0.5ex>[r] & H}$ be a Lie $2$-group with associated crossed module $\xymatrix{G \ar[r]^i & H}$. We write $\G_p :=\G^{(p)}$ and remark that this is a Lie subgroup of $\G ^p$ for each $p$. Consider the associated double complex 
\begin{eqnarray*}
\xymatrix{
\vdots                                & \vdots                                  & \vdots                         &       \\ 
C(H^3) \ar[r]^{\partial}\ar[u]        & C(\G ^3) \ar[r]^{\partial}\ar[u]        & C(\G _2^3)\ar[r]\ar[u]         & \dots \\
C(H^2) \ar[r]^{\partial}\ar[u]^\delta & C(\G ^2) \ar[r]^{\partial}\ar[u]^\delta & C(\G _2^2) \ar[r]\ar[u]^\delta & \dots \\
C(H) \ar[r]^{\partial}\ar[u]^\delta   & C(\G) \ar[r]^{\partial}\ar[u]^\delta    & C(\G _2) \ar[r]\ar[u]^\delta   & \dots 
}
\end{eqnarray*} 
In it, columns are usual complexes of Lie group cochains, and rows are complexes of groupoid cochains for powers of the Lie $2$-group. For future reference, we write the total complex 
\begin{eqnarray*}
\Omega ^k_{tot}(\G)=\bigoplus _{p+q=k}C(\G _p^q),
\end{eqnarray*}
with differential $d=\delta +(-1)^{q}\partial$. We will refer to the double groupoid cohomology of $\G$ simply as $2$-cohomology.

We give an interpretation of $H^2_{tot}(\G )$. A $2$-cocycle consists of a pair of functions $(F,f)\in C(H^2)\oplus C(\G )$ such that:\\
1) $\delta F = 0$, i.e. $F(h_1 ,h_2)+F(h_0 ,h_1 h_2)=F(h_0 h_1 ,h_2)+F(h_0 ,h_1 )$ for all triples $h_0,h_1,h_2\in H$ \\
2) $\partial f =0$, i.e.  $f(\gamma _1 \Join\gamma _2 )=f(\gamma _1 )+f(\gamma _2 )$ for all $(\gamma _1 ,\gamma _2 )\in\G _2$. Using the decomposition $\G \cong G\rtimes H$, this is equivalent to $f(g_2g_1,h)=f(g_2,h)+f(g_1,hi(g_2))$ for all $h\in H$ and $g_1,g_2\in G$.\\
3) $\partial F+\delta f=0$, i.e. $f(\gamma _1)-f(\gamma _0\vJoin\gamma _1)+f(\gamma _0)=F(t(\gamma _0),t(\gamma _1))-F(s(\gamma _0),s(\gamma _1))$ for all pairs $\gamma _0,\gamma _1\in\G$. Again, using, this equation can be rewritten as
\begin{eqnarray*}
f(g_1,h_1)-f(g_0^{h_1}g_1,h_0h_1)+f(g_0,h_0)=F(h_0i(g_0),h_1i(g_1))-F(h_0,h_1).
\end{eqnarray*}
Notice that, by making $h_0=h_1=1$, this equation yields 
\begin{eqnarray*}
f(g_0g_1,1)=f(g_0,1)+f(g_1,1)-F(i(g_0),i(g_1)).
\end{eqnarray*}
Also as a consequence of the second item above, replacing $g_2$ by the identity element, we have got that $f(1,h)=0$ for every $h\in H$. \\
It is rather known that if $\delta F=0$, $F$ induces a (central) extension of $H$, 
\begin{eqnarray*}
\xymatrix{
1 \ar[r] & \Rr \ar[r]^{\bar{1}\times I\quad} & H\ltimes^F\Rr  \ar[r]^{\quad pr_1} & H \ar[r] & 1,
}
\end{eqnarray*}
where $H\ltimes^F \Rr$ is the semi-direct product induced by $F$, that is the space $H\times\Rr$ endowed with the twisted product
\begin{eqnarray*}
(h_0,\lambda_0)\odot_F(h_1,\lambda_1):=(h_0h_1,\lambda_0+\lambda_1+F(h_0,h_1)).
\end{eqnarray*}
\begin{lemma}\label{easyExt}
If $d(F,f)=0$, then 
\begin{eqnarray*}
\xymatrix{
\psi_f: G \ar[r] & H\ltimes^F\Rr : g \ar@{|->}[r] & (i(g),-f(g,1))
}
\end{eqnarray*}
defines a crossed module for the action $g^{(h,\lambda)}:=g^h$.
\end{lemma}
\begin{proof}
First, let us verify that $\psi_f$ is indeed a Lie group homomorphism. 
\begin{align*}
\psi_f(g_0)\odot_F\psi_f(g_1) & = (i(g_0),-f(g_0,1))\odot_F(i(g_1),-f(g_1,1))         \\
                              & = (i(g_0)i(g_1),-f(g_0,1)-f(g_1,1)+F(i(g_0),i(g_1))) \\
                              & = (i(g_0g_1),-f(g_0g_1,1))=\psi_f(g_0g_1)
\end{align*}
here $\partial F=\delta f$ justifies the third equality as remarked. Now, the action thus defined is still a right action by automorphisms 
\begin{eqnarray*}
g^{(h_0,\lambda_0)\odot_F(h_1,\lambda_1)}=g^{(h_0h_1,\lambda_0+\lambda_1+F(h_0,h_1))}=g^{h_0h_1}=(g^{h_0})^{h_1}, 
\end{eqnarray*}
\begin{eqnarray*}
(g_1g_2)^{(h,\lambda)}=(g_1g_2)^h=g_1^hg_2^h=g_1^{(h,\lambda)}g_2^{(h,\lambda)}.
\end{eqnarray*}
Peiffer identity holds as a consequence of the independence of the variable in $\Rr$ as well
\begin{eqnarray*}
g_1^{\psi_f(g_2)}=g_1^{(i(g_2),-f(g_2,1))}=g_1^{i(g_2)}.
\end{eqnarray*}
As for the equivariance of $\psi_f$, on the one hand we have got
\begin{eqnarray*}
\psi_f(g^{(h,\lambda)})=(i(g^h),-f(g^h,1)),
\end{eqnarray*}
while on the other,
\begin{align*}
(h,\lambda)^{-1}\odot_F\psi_f(g)\odot_F(h,\lambda) &= (h^{-1},-\lambda-F(h^{-1},h))\odot_F(i(g),-f(g,1))\odot_F(h,\lambda) \\
                                                  &= (h^{-1}i(g),-\lambda-F(h^{-1},h)-f(g,1)+F(h^{-1},i(g)))\odot_F(h,\lambda) \\
                                                  &= (h^{-1}i(g)h,-\cancel{\lambda}-F(h^{-1},h)-f(g,1)+F(h^{-1},i(g))+\cancel{\lambda}+F(h^{-1}i(g),h))                    
\end{align*}
The first entry coincides because $i$ is the structural morphism of a crossed module, whereas using $\partial F+\delta f=0$ evaluated on $((g,h^{-1}),(1,h))$, we get
\begin{align*}
f(1,h)-f(g^h,1)+f(g,h^{-1}) & =F(h^{-1}i(g),h)-F(h^{-1},h) \\
-f(g^h,1)                   & =-F(h^{-1},h)-f(g,h^{-1})+F(h^{-1}i(g),h).
\end{align*}
Here, the term $f(1,h)$ vanished as a consequence of $\partial f=0$, as we saw right after spelling out the cocycle equations. Using $\partial F+\delta f=0$ once again, though evaluating at $((1,h^{-1}),(g,1))$ this time around, we get
\begin{align*}
f(g,1)-f(g,h^{-1})+f(1,h^{-1}) & =F(h^{-1},i(g))-F(h^{-1},1)  \\
    -f(g,h^{-1})               & =-f(g,1)+F(h^{-1},i(g));
\end{align*} 
thus proving the lemma.

\end{proof}
As a consequence of this lemma, we see that there is an induced short exact sequence of Lie $2$-groups that we write using their associated crossed modules
\begin{eqnarray*}
\xymatrix{
1 \ar[r] & 1 \ar[d]\ar[r] & G \ar[d]_{\psi_f}\ar[r]^{Id}      & G \ar[d]^i \ar[r]   & 1  \\
1 \ar[r] & \Rr \ar[r]     & H\ltimes^F\Rr \ar[r]_{\quad pr_1} & H \ar[r]            & 1.  
}
\end{eqnarray*}

\begin{lemma}\label{isoEasyExt}
Let $(F,f),(F',f')\in\Omega^2_{tot}(\G)$. If there exists a $\phi\in\Omega^1_{tot}(\G)=C(H)$, such that $(F,f)-(F',f')=d\phi$, then the induced extensions are isomorphic.
\end{lemma}
\begin{proof}
First, recall that the object extensions are isomorphic via
\begin{eqnarray*}
\alpha\xymatrix{ 
:H\ltimes^F\Rr \ar[r] & H\ltimes^{F'}\Rr:(h,\lambda) \ar@{|->}[r] & (h,\lambda+\phi(h)).
}
\end{eqnarray*}
We claim that this together with the identity of $G$ induce the isomorphism between the extensions. Indeed, using the notation from the previous lemma
\begin{align*}
\alpha(\psi_f(g)) &= \alpha(i(g),-f(g,1)) \\
                  &= (i(g),-f(g,1)+\phi(i(g))) \\
                  &= (i(g),-f'(g,1)) = \psi_{f'}(g),
\end{align*}
also, trivially $id_G(g^{(h,\lambda)})=id_G(g)^{\alpha(h,\lambda+\phi(h))}$, thus finishing the proof.

\end{proof}

\section{Representations of Lie 2-groups}
We turn now to defining what a representation of a Lie $2$-group is. In order to do so, we are going to suitably integrate the linear Lie $2$-algebra $\ggl (\phi)$. At this point we would like to single out the fact that even when considering this process for the subcategories of Lie groups and Lie algebras, one doesn't use the connected, nor the simply connected integration of the $\ggl (n)$'s. In fact, the spaces one takes are the $GL(n)$'s, which always have two connected components for $n>0$ and which are simply connected if, and only if $n=1$ (\textit{In case you are wondering, $\pi_1(GL(2)^0)=\Zz$ and the fundamental group of the identity component is $\Zz_2$ for any other $n$}). Motivated by this we consider an integration of $\ggl (\phi)$ whose space of objects is a subgroup of $GL(W)\times GL(V)$.\\

\subsection{The General Linear Lie 2-group of a 2-vector space}
The General Linear Lie $2$-group has already appeared in \cite{IntSubLin2}, we go over its construction from a slightly different point of view and add it for completeness. We construct the crossed module whose associated Lie $2$-group will play the r\^ole of space of automorphisms of our given $2$-vector space $\xymatrix{W \ar[r]^\phi & V}$:
\begin{eqnarray*}
\xymatrix{
GL(\phi)_1 \ar[r]^\Delta & GL(\phi)_0.
}
\end{eqnarray*}
The space of objects is the space of invertible self functors
\begin{eqnarray*}
GL(\phi)_0 =\lbrace (F,f)\in GL(W)\times GL(V) : \phi\circ F=f\circ\phi\rbrace .
\end{eqnarray*}
Notice that this is a Lie subgroup of $GL(W)\times GL(V)$ as claimed, since it contains the identity $(I,I)$ and given $(F_1,f_1),(F_2,f_2)\in GL(\phi)_0$,
\begin{align*}
\phi\circ (F_1F_2) & =(\phi\circ F_1)F_2 \\
                   & =f_1\circ\phi F_2 =(f_1f_2)\circ\phi .
\end{align*}
Now, the group of arrows is given by the set
\begin{eqnarray*}
GL(\phi)_1 =\lbrace A\in Hom(V,W):(I+A\phi ,I+\phi A)\in GL(W)\times GL(V)\rbrace ,
\end{eqnarray*}
endowed with the operation
\begin{eqnarray*}
A_1\odot A_2 := A_1+A_2+A_1\phi A_2,
\end{eqnarray*}
for which the identity element is the $0$ map, and inverses are given by either
\begin{eqnarray*}
A^{\dagger} = -A(I+\phi A)^{-1}=-(I+A\phi )^{-1}A.
\end{eqnarray*}
We write $\dagger$ instead of $-1$ to avoid any possible overlap of notation with the actual inverse of a matrix.
Instead of verifying that this defines a group operation, we just point out that it is drawn out of the diagonal embedding in $GL(W\oplus V)$, thus realizing $GL(\phi)_1$ as a Lie subgroup.\\
The map $\Delta$ is the projection onto the diagonal components, that is
\begin{eqnarray*}
\Delta A=(I+A\phi ,I+\phi A).
\end{eqnarray*}
This is  well defined since by definition it lands in $GL(W)\times GL(V)$, and 
\begin{align*}
\phi(I+A\phi) & = \phi +\phi A\phi \\
              & = (I+\phi A)\phi .
\end{align*}
Finally, the right action of $GL(\phi)_0$ on $GL(\phi)_1$ is given by
\begin{eqnarray*}
A^{(F,f)}=F^{-1}Af.
\end{eqnarray*}
\begin{prop}\label{GLPhi}
Along with this group structure and this action,
\begin{eqnarray*}
\xymatrix{
GL(\phi)_1 \ar[r]^{\Delta} & GL(\phi)_0}
\end{eqnarray*} 
is a crossed module of Lie groups.
\end{prop}
\begin{proof}
Again, this amounts to a routine check.
\begin{itemize}
\item $\Delta$ is a Lie group homomorphism: 
\begin{align*}
\Delta (A_1\odot A_2) & = (I+(A_1\odot A_2)\phi ,I +\phi (A_1\odot A_2)) \\
                       & = (I + (A_1+ A_2+A_1\phi A_2)\phi ,I+\phi (A_1+A_2+A_1\phi A_2)) \\
                       & = (I +A_1\phi + A_2\phi +A_1\phi A_2\phi ,I+\phi A_1+\phi A_2+\phi A_1\phi A_2) \\
                       & = ((I+A_1\phi )(I+A_2\phi ),(I+\phi A_1)(I+\phi A_2))=\Delta A_1\Delta A_2.
\end{align*}
\item $GL(\phi)_0$ acts by automorphisms:
\begin{align*}
(A_1\odot A_2)^{(F,f)} & = F^{-1}(A_1\odot A_2)f \\
                        & = F^{-1}(A_1+ A_2+A_1\phi A_2)f \\
                        & = F^{-1}A_1f+F^{-1}A_2f+F^{-1}A_1\phi (FF^{-1})A_2f \\
                        & = F^{-1}A_1f+F^{-1}A_2f+F^{-1}A_1f\phi F^{-1}A_2f =A_1^{(F,f)}\odot A_2^{(F,f)}  ,
\end{align*}
where we pass to the last line using $\phi F=f\phi$.
\item The map $\xymatrix{GL(\phi)_0 \ar[r] & Aut(GL(\phi)_1)}$ is a indeed a right action: $A^{(I,I)}=A$, and
\begin{align*}
A^{(F_1,f_1)(F_2,f_2)} & = (F_1F_2)^{-1}A(f_1f_2) \\
                       & = F_2^{-1}(F_1^{-1}Af_1)f_2 =(A^{(F_1,f_1)})^{(F_2,f_2)}.
\end{align*}
\item $\Delta$ is equivariant:
\begin{align*}
\Delta (A^{(F,f)}) & = (I+A^{(F,f)}\phi ,I+\phi A^{(F,f)}) \\
                   & = (I+F^{-1}Af\phi ,I+\phi F^{-1}Af) \\
                   & = (F^{-1}F+F^{-1}A\phi F,f^{-1}f+f^{-1}\phi Af) \\
                   & = (F^{-1}(I+A\phi )F,f^{-1}(I+\phi A)f)=(F,f)^{-1}\Delta A(F,f) ,
\end{align*}
where we used the obvious $\phi F^{-1}=f^{-1}\phi$ in the third line.
\item Peiffer identity: 
\begin{align*}
A_1^{\Delta A_2} & = A_1^{(I+A_2\phi ,I+\phi A_2)}\\
                 & = (I+ A_2\phi )^{-1}A_1(I+\phi A_2) ,
\end{align*}
while on the other hand 
\begin{align*}
A_2^{\dagger}\odot A_1\odot A_2 & = A_2^{\dagger}\odot (A_1+A_2+A_1\phi A_2) \\
                                & = A_2^{\dagger} + A_1+A_2+A_1\phi A_2 + A_2^{\dagger}\phi(A_1+A_2+A_1\phi A_2) \\
                           & = A_1(I+\phi A_2) + A_2^{\dagger}\phi A_1(I+\phi A_2)+A_2 +A_2^{\dagger}(I+\phi A_2) .
\end{align*}
Now, $A_2^{\dagger}(I+\phi A_2)=-A_2$; thus, 
\begin{align*}
A_2^{\dagger}\odot A_1\odot A_2 & = (I+A_2^{\dagger}\phi )A_1(I+\phi A_2) .
\end{align*}
We finish the proof by showing that $I+A_2^{\dagger}\phi =(I+ A_2\phi )^{-1}$; indeed,
\begin{align*}
(I+A_2^{\dagger}\phi )(I+ A_2\phi ) & = I+A_2^{\dagger}\phi +A_2\phi +A_2^{\dagger}\phi A_2\phi \\
                               & = I+(A_2^{\dagger} +A_2 +A_2^{\dagger}\phi A_2)\phi \\
                               & = I+(A_2^{\dagger}(I+\phi A_2)+A_2)\phi \\
                               & = I+(-A_2+A_2)\phi =I .
\end{align*}
\end{itemize}
\end{proof}
In the sequel, we will adopt the notation $GL(\phi )$ for the crossed module of Proposition \ref{GLPhi} and its associated Lie $2$-group as well. We prove that this object serve as the global counterpart of $\ggl (\phi )$.
\begin{prop}\label{Int'ingFullLinear}
The Lie $2$-algebra of $GL(\phi )$ is $\ggl (\phi )$ .
\end{prop}
\begin{proof}
The Lie algebra of $GL(\phi)_0$ clearly is $\ggl(\phi)_0$, as the condition defining them is linear and the same. For the Lie algebra of $GL(\phi)_1$, we look at it as diagonally embedded in $GL(W\oplus V)$. Then, the derivative of this embedding will give us the inclusion of the Lie algebra of $GL(\phi)_1$ in $\ggl (W\oplus V)$. Consider a curve $\xymatrix{T:(-\varepsilon ,\varepsilon) \ar[r] & GL(\phi)_1}$ with  $T(0)=0$. Then, 
\begin{eqnarray*}
\frac{d}{d\lambda}\bigg\rvert_{\lambda =0}\begin{pmatrix}
    I+T(\lambda)\phi & 0 \\
    0                & I+\phi T(\lambda) 
\end{pmatrix}
=\begin{pmatrix}
    (\frac{d}{d\lambda}\big\rvert_{\lambda =0}T(\lambda))\phi & 0 \\
    0                                                          & \phi\frac{d}{d\lambda}\big\rvert_{\lambda =0}T(\lambda) 
\end{pmatrix}.
\end{eqnarray*}
To get the desired conclusion, we just need to point out that we can regard $\ggl(\phi )_1$ as a Lie subalgebra of $\ggl (W\oplus V)$, exactly by means of 
\begin{eqnarray*}
\xymatrix{
\ggl(\phi )_1 \ar[r] & \ggl (W\oplus V):A \ar@{|->}[r] & {}
}\begin{pmatrix}
    A\phi & 0 \\
    0     & \phi A 
\end{pmatrix};
\end{eqnarray*} %It is not true that this map yields an inclusion, nor that the "diagonal embedding" is an embedding. Both these maps have got a kernel. They are given by the same subspace of Hom(V,W) of maps that send elements in the cokernel of \phi, to the kernel of \phi. This is simply a subspace, i.e. it is abelian when regarded as an algebra and as a group; therefore, the spaces are (direct) products of these so-called embeddings with the kernel. So the proof is still an argument for the Lie-correspondence.
thus, we conclude not only that $\ggl(\phi )_1$ is the Lie algebra of $GL(\phi)_1$, but incidentally saw that the derivative at the identities of the crossed module map is the right one. We are left to verify that the (second) derivative of the action is $\Lie$. This actually follows from the embedding as well, since the right action of $(F,f)$ on $T$ is nothing but the conjugation
\begin{align*}
\begin{pmatrix}
    F & 0 \\
    0 & f 
\end{pmatrix}^{-1}
\begin{pmatrix}
    I+T\phi & 0 \\
    0       & I+\phi T 
\end{pmatrix}
\begin{pmatrix}
    F & 0 \\
    0 & f 
\end{pmatrix}
& =\begin{pmatrix}
    F^{-1}(I+T\phi )F & 0 \\
    0                 & f^{-1}(I+\phi T)f  
\end{pmatrix} \\
& =\begin{pmatrix}
    I+(F^{-1}Tf)\phi & 0 \\
    0              & I+\phi (F^{-1}Tf)
\end{pmatrix};
\end{align*}
whereas, $\Lie_{(F,f)}A$ is the bracket
\begin{align*}
\bigg{[}\begin{pmatrix}
    F & 0 \\
    0 & f 
\end{pmatrix},
\begin{pmatrix}
    A\phi & 0 \\
    0       & \phi A 
\end{pmatrix}\bigg{]}
& =\begin{pmatrix}
    FA\phi & 0 \\
    0      & f\phi A  
\end{pmatrix}-
\begin{pmatrix}
    A\phi F & 0 \\
    0       & \phi Af  
\end{pmatrix} \\
& =\begin{pmatrix}
    (FA-Af)\phi & 0 \\
    0           & \phi (FA-Af)
\end{pmatrix};
\end{align*}
thus, the result follows.

\end{proof}
We use the strategy of the proof of the previous proposition to come up with a formula for the the exponential of $\ggl(\phi)_1$ to $GL(\phi)_1$.
\begin{prop}\label{TheExpOfGL(phi)}
The exponential $\xymatrix{\exp_{GL(\phi)_1} :\ggl(\phi)_1 \ar[r] & GL(\phi)_1}$ is given by the formula
\begin{eqnarray*}
\exp_{GL(\phi)_1}(A)=A\sum_{n=0}^\infty\frac{(\phi A)^n}{(n+1)!}=\sum_{n=0}^\infty\frac{(A\phi)^n}{(n+1)!}A.
\end{eqnarray*}
\end{prop}
\begin{proof}
We start by noticing that the two defining equations coincide. Indeed,
\begin{eqnarray*}
A\sum_{n=0}^\infty\frac{(\phi A)^n}{(n+1)!}=\sum_{n=0}^\infty\frac{A(\phi A)^n}{(n+1)!}.
\end{eqnarray*}
We prove inductively that $A(\phi A)^n=(A\phi )^nA$. This is obviously true for $n=0$ and $n=1$, so suppose the equation holds for $k$, then
\begin{eqnarray*}
A(\phi A)^{k+1}=A(\phi A)^{k}(\phi A)=(A\phi )^k(A\phi)A=(A\phi)^{k+1}A
\end{eqnarray*}
as desired. Then, of course,
\begin{eqnarray*}
\sum_{n=0}^\infty\frac{A(\phi A)^n}{(n+1)!}=\sum_{n=0}^\infty\frac{(A\phi)^n}{(n+1)!}A.
\end{eqnarray*}
Now, recall that $\ggl(\phi)_1$ and $GL(\phi)_1$ could be respectively regarded as embedded in $\ggl(W\oplus V)$ and $GL(W\oplus V)$ by
\begin{eqnarray*}
\xymatrix{
A \ar@{|->}[r] & }\begin{pmatrix}
A\phi & 0 \\
0 & \phi A
\end{pmatrix}
 & \textnormal{and} & \xymatrix{
A \ar@{|->}[r] & }\begin{pmatrix}
I+A\phi & 0 \\
0 & I+\phi A
\end{pmatrix};
\end{eqnarray*}
thus, for $A\in\ggl(\phi)_1$,
\begin{align*}
\exp\begin{pmatrix}
A\phi & 0 \\
0 & \phi A
\end{pmatrix} & =\begin{pmatrix}
\exp(A\phi) & 0 \\
0 & \exp(\phi A)
\end{pmatrix} \\
              & =\begin{pmatrix}
I+A\phi +\frac{(A\phi)^2}{2}+\frac{(A\phi)^3}{3!}+\frac{(A\phi)^4}{4!}+... & 0 \\
0 & I+\phi A+\frac{(\phi A)^2}{2}+\frac{(\phi A)^3}{3!}+\frac{(\phi A)^4}{4!}+...
\end{pmatrix} \\
			  & =\begin{pmatrix}
I+(A+\frac{A\phi A}{2}+\frac{A(\phi A)^2}{3!}+\frac{A(\phi A)^3}{4!}+...)\phi & 0 \\
0 & I+\phi (A+\frac{A\phi A}{2}+\frac{(A\phi)^2A}{3!}+\frac{(A\phi)^3A}{4!}+...)
\end{pmatrix} \\
			  & =\begin{pmatrix}
I+A\sum_{n=0}^\infty\frac{(\phi A)^n}{(n+1)!}\phi & 0 \\
0 & I+\phi\sum_{n=0}^\infty\frac{(A\phi)^n}{(n+1)!}A
\end{pmatrix} 
\end{align*}
and the result follows. 

\end{proof}

\subsection{2-representations}

Now that we have got a Lie $2$-group of endomorphisms of a $2$-vector space $\xymatrix{W \ar[r]^\phi & V}$, we define a representation of a Lie $2$-group.
\begin{Def}
A \textit{representation} of a Lie $2$-group $\G =\xymatrix{G \ar[r]^i & H}$ on $\xymatrix{W \ar[r]^\phi & V}$ is a morphism of Lie $2$-groups
\begin{eqnarray*}
\xymatrix{
\rho :\G \ar[r] & GL(\phi ).
}
\end{eqnarray*}
\end{Def}
We will refer again to such a map as a $2$-representation. The same remarks after the definition of a Lie $2$-algebra representation apply to the case of Lie $2$-groups. In particular, by morphism of Lie $2$-groups, we mean a na\"ive smooth functor respecting the Lie group structures. 
From the very construction and as a consequence of proposition \ref{Int'ingFullLinear}, we have got the following obvious  proposition.
\begin{prop}\label{RepGoToRep}
If $\xymatrix{\rho :\G \ar[r] & GL(\phi )}$ is Lie $2$-group representation, then its derivative 
\begin{eqnarray*}
\xymatrix{
\gg \ar[d]_\mu \ar[r]^{d_1\rho_1\quad} & \ggl(\phi)_1 \ar[d]^\Delta \\
\hh \ar[r]_{d_1\rho_0\quad}            & \ggl(\phi)_0 
}
\end{eqnarray*}
is a Lie $2$-algebra representation.
\end{prop}
Given a $2$-representation, one has got again two ``commuting'' honest representations on $W$ and on $V$, though this time around, they are Lie group representations of $H$. That is, the map at objects
\begin{eqnarray*}
\xymatrix{
\rho_0 :H \ar[r] & GL(\phi)_0\leq GL(W)\times GL(V), 
}
\end{eqnarray*}
consists of two representations $\xymatrix{\rho_0^1 :H \ar[r] & GL(W)}$, $\xymatrix{\rho_0^0 :H \ar[r] & GL(V)}$ fitting in
\begin{eqnarray*}
\xymatrix{
W \ar[d]_\phi \ar[r]^{(\rho_0^1)_h} & W \ar[d]^\phi \\
V \ar[r]_{(\rho_0^0)_h}             & V,
}
\end{eqnarray*} 
for each $h\in H$.

Additionally, one has got the map 
\begin{eqnarray*}
\xymatrix{
\rho_1 :G \ar[r] & GL(\phi)_1\subset Hom(V,W), 
}
\end{eqnarray*}
at the level of arrows. This map is a Lie group homomorphism; hence,
\begin{eqnarray*}
\rho_1(g_0g_1)=\rho_1(g_0)+\rho_1(g_1)+\rho_1(g_0)\phi\rho_1(g_1),
\end{eqnarray*}
and it also verifies the following equations for all $g\in G$:
\begin{eqnarray*}
\rho^0_0(i(g))=I+\phi\rho_1(g), & & \rho^1_0(i(g))=I+\rho_1(g)\phi ,
\end{eqnarray*}
and for all $h\in H$
\begin{eqnarray*}
\rho_1(g^h)=\rho_0^1(h)^{-1}\rho_1(g)\rho_0^0(h). 
\end{eqnarray*}
We now build an honest representation of $\G\cong G\rtimes H$ on $W\oplus V$.

\begin{prop}\label{honestGpRep}
Given a representation $2$-representation $\xymatrix{\rho :\G \ar[r] & GL(\phi )}$, there is an honest representation 
\begin{eqnarray*}
\xymatrix{
\bar{\rho}:G\rtimes H \ar[r] & GL(W\oplus V):(g,h) \ar@{|->}[r] & {}
}
\begin{pmatrix}
    \rho_0^1(hi(g)) & \rho_0^1(h)\rho_1(g) \\
    0               & \rho_0^0(h) 
\end{pmatrix}
\end{eqnarray*}
\end{prop}
\begin{proof}
Consider the product
\begin{align*}
\bar{\rho}(g_0,h_0)\bar{\rho}(g_1,h_1) & = \begin{pmatrix}
                                               \rho_0^1(h_0i(g_0)) & \rho_0^1(h_0)\rho_1(g_0) \\
                                               0                   & \rho_0^0(h_0) 
                                           \end{pmatrix}
                                           \begin{pmatrix}
                                               \rho_0^1(h_1i(g_1)) & \rho_0^1(h_1)\rho_1(g_1) \\
                                               0                   & \rho_0^0(h_1) 
                                           \end{pmatrix} \\
                                       & = \begin{pmatrix}
     \rho_0^1(h_0i(g_0))\rho_0^1(h_1i(g_1)) & \rho_0^1(h_0i(g_0))\rho_0^1(h_1)\rho_1(g_1) +\rho_0^1(h_0)\rho_1(g_0)\rho_0^0(h_1) \\
                               0            & \rho_0^0(h_0)\rho_0^0(h_1) 
                                           \end{pmatrix}   . 
\end{align*}
The bottom row will agree with the bottom row of $\bar{\rho}((g_0,h_0)\vJoin (g_1,h_1))=\bar{\rho}(g_0^{h_1}g_1,h_0h_1)$, because $\rho_0^0$ is a group homomorphism. The first entries coincide too, put simply, because the target is a group homomorphism as well:
\begin{align*}
\rho_0^1(h_0h_1i(g_0^{h_1}g_1)) & = \rho_0^1(h_0h_1h_1^{-1}i(g_0)h_1i(g_1)) \\
                                & = \rho_0^1(h_0i(g_0))\rho_0^1(h_1i(g_1)) .
\end{align*}
Finally, for the remaining entries, we use the fact that $\rho_1$ is a homomorphism, and that $\rho$ respects the actions:
\begin{align*}
\rho_0^1(h_0h_1)\rho_1(g_0^{h_1}g_1) & = \rho_0^1(h_0h_1)(\rho_1(g_0^{h_1})(I+\phi\rho_1(g_1))+\rho_1(g_1)) \\
                                     & = \rho_0^1(h_0h_1)\rho_0^1(h_1)^{-1}\rho_1(g_0)\rho_0^0(h_1)(I+\phi\rho_1(g_1))+\rho_0^1(h_0h_1)\rho_1(g_1) \\
                                     & = \rho_0^1(h_0)\rho_1(g_0)\rho_0^0(h_1)+\rho_0^1(h_0)\rho_1(g_0)\phi\rho_0^1(h_1)\rho_1(g_1)+\rho_0^1(h_0)\rho_0^1(h_1)\rho_1(g_1) \\
                                     & = \rho_0^1(h_0)\rho_1(g_0)\rho_0^0(h_1)+\rho_0^1(h_0)(\rho_1(g_0)\phi +I)\rho_0^1(h_1)\rho_1(g_1) ,
\end{align*}
but $\Delta\circ\rho_1=\rho_0\circ i$, so 
\begin{align*}
I+\rho_1(g_0)\phi & = \rho_0^1(i(g_0)) ,
\end{align*}
and the result follows.

\end{proof}
Using this representation, we now build the semi-direct product of Lie $2$-groups, 
\begin{eqnarray*}
\xymatrix{
(G\rtimes H){}_{\bar{\rho}}\ltimes (W\oplus V) \ar@<0.5ex>[r] \ar@<-0.5ex>[r] & H{}_{\rho_0^0}\ltimes V.
}
\end{eqnarray*}
We specify the structure by describing the associated crossed module
\begin{eqnarray*}
\xymatrix{
G{}_{\rho_0^1\circ i}\ltimes W \ar[r]^{i\times\phi} & H{}_{\rho_0^0}\ltimes V,
}
\end{eqnarray*}
where the right action is given by
\begin{eqnarray*}
(g,w)^{(h,w)}=(g^h,\rho_0^1(h)^{-1}(w+\rho_1(g)v)).
\end{eqnarray*}
Forgetting for the time being the Lie group structure,
\begin{eqnarray*}
\xymatrix{
\G{}_{\bar{\rho}}\ltimes\mathbb{V} \ar@<-0.5ex>[r]\ar@<0.5ex>[r]\ar[d] & H{}_{\rho_0^0}\ltimes V \ar[d] \\
\G \ar@<-0.5ex>[r]\ar@<0.5ex>[r]                                       & H
}
\end{eqnarray*}
has got the structure of a VB-groupoid. Indeed, the structural maps are defined to cover those of the Lie $2$-group we started with and are defined using the representations and the very structure of $\xymatrix{W \ar[r]^\phi & V}$. Thus, there is an associated representation up to homotopy (cf. appendix \ref{appchapter}). In fact, since all vector bundles are trivial, there is an obvious splitting of the core sequence 
\begin{eqnarray*}
\xymatrix{
(0) \ar[r] & \G\times W \ar[r] & \G{}_{\bar{\rho}}\ltimes\mathbb{V} \ar[r] & \G\times V \ar[r] & (0).
}
\end{eqnarray*}
given by
\begin{eqnarray*}
\xymatrix{
\mathrm{h}:\G\times V \ar[r] & \G{}_{\bar{\rho}}\ltimes\mathbb{V}:(g,h;v) \ar@{|->}[r] & \mathrm{h}_{(g,h)}(h,v):=(g,h;0,v),
} 
\end{eqnarray*} 
which verifies $\mathrm{h}_{u(h)}(h,v)=\mathrm{h}_{(1,h)}(h,v)=(1,h;0,v)=\hat{u}(h,v)$. This splitting will provide a canonical representation up to homotopy associated to the $2$-representation. Such a representation up to homotopy is defined by
\begin{eqnarray*}
\xymatrix{
\eth :H\times W \ar[r] & H\times V:(h,w) \ar@{|->}[r] & \hat{t}(1,h;w,0)=(h,0)(i(1),\phi(w))=(h,\rho_0^0(h)\phi(w)),
}
\end{eqnarray*}
together with the quasi-actions
\begin{align*}
\Delta^V_{(g,h)}(h,v) & = \hat{t}(\mathrm{h}_{(g,h)}(h,v))=\hat{t}(g,h;0,v)=(h,v)(i(g),0)=(hi(g),v) \\
\Delta^W_{(g,h)}(h,w) & = \mathrm{h}_{(g,h)}(\partial(h,w))\hat{\Join}(1,h;w,0)\hat{\Join}(g^{-1},hi(g);0,0) \\
                      & = (g,h;0,\rho_0^0(h)\phi(w))\hat{\Join}(g^{-1},hi(g);\rho_0^1(i(g))^{-1}w,0)=(1,hi(g);\rho_0^1(i(g))^{-1}w,0)
\end{align*}
and finally,
\begin{align*}
\Omega_{(g_1,g_2,h)}(h,v) & =(\mathrm{h}_{(g_2g_1,h)}(h,v)-\mathrm{h}_{(g_1,hi(g_2))}(\Delta^V_{(g_2,h)}(h,v))\hat{\Join}\mathrm{h}_{(g_2,h)}(h,v))\hat{\Join} 0_{(g_2g_1,h)^{-1}} \\                          
                          & =((g_2g_1,h;0,v)-\mathrm{h}_{(g_1,hi(g_2))}(hi(g_2),v)\hat{\Join}(g_2,h;0,v))\hat{\Join} 0_{((g_2g_1)^{-1},hi(g_2g_1))} \\
                          & =((g_2g_1,h;0,v)-(g_1,hi(g_2);0,v)\hat{\Join}(g_2,h;0,v))\hat{\Join}((g_2g_1)^{-1},hi(g_2g_1);0,0) \\
                          & =((g_2g_1,h;0,v)-(g_2g_1,hi(g_2);0,v))\hat{\Join}((g_2g_1)^{-1},hi(g_2g_1);0,0) \\
                          & =(g_2g_1,h;0,0)\hat{\Join}((g_2g_1)^{-1},hi(g_2g_1);0,0)=(1,hi(g_2g_1);0,0). 
\end{align*}
Here, we used $\eth$ instead of the usual $\partial$, since this latter symbol will have a different meaning in this section. Notice that since the curvature form $\Omega$ is zero, the quasi-actions define actual representations of the Lie $2$-group $\G$ on the corresponding vector bundles. \\
We end this section with some examples.
\begin{ex:}
Trivial representations. Of course, all of the defining equations for a $2$-representation get satisfied trivially if $(\rho_1,\rho_0)\equiv(0,I)$.
\end{ex:}
\begin{ex:}\label{unitGpRep}
Just as in the case of Lie $2$-algebras, usual Lie group representations can be regarded as $2$-representations of the unit Lie $2$-group on the unit $2$-vector space. Indeed, if $W=(0)$, then a $2$-representation is equivalent to a single representation $\rho$ of $H$ on $V$, such that $\rho_{i(g)}\equiv I$ for every $g\in G$. One has got an analogous situation if $V=(0)$, in which a $2$-representation reduces to a single representation of $H$ on $W$ with the same vanishing property on $i(G)$. Ultimately, any $2$-representation on either of these $2$-vector spaces turns out to be simply a representation of $H/i(G)$. 
\end{ex:}
\begin{ex:}
The adjoint representation: Let $\xymatrix{\gg \ar[r]^\mu & \hh}$ be the Lie $2$-algebra of the Lie $2$-group $\xymatrix{G \ar[r]^i & H}$, then we have 
\begin{eqnarray*}
\xymatrix{
G \ar[d] \ar[r]^{Ad_1\quad} & GL(\mu)_1\ar[d]\\
H \ar[r]_{Ad_0\quad}        & GL(\mu)_0 ,
}
\end{eqnarray*}
where $(Ad_0)_h=((-)^{h^{-1}},Ad_{h^{-1}})$, and $(Ad_1)_g=d_e({}_g\wedge)$ with
\begin{eqnarray*}
\xymatrix{
{}_g\wedge : H \ar[r] & G:h \ar@{|->}[r] & g(g^{-1})^{h^{-1}}. 
}
\end{eqnarray*}
\end{ex:}

\section{Extensions of Lie 2-groups and 2-cohomology with coefficients}
We look for a cohomology theory of Lie $2$-groups. As a requisite, we want the second cohomology to classify extensions. We will write down the equations defining the differentials, and the spaces of $1$-cochains and $2$-cochains. Eventually, we make sense of them understanding them as part of a larger complex as it was the case in the Lie $2$-algebra case. \\
We start this section by proving that the notion of $2$-representations defined in the previous section is the correct one in the sense that it is the kind of action induced in an abstract (split) extension.
\begin{Def}
An \textit{extension} of the Lie $2$-group $\xymatrix{G \ar[r]^{i} & H}$ by the $2$-vector space $\xymatrix{W \ar[r]^{\phi} & V}$ is a Lie $2$-group $\xymatrix{E_1 \ar[r]^{\epsilon} & E_0}$ that fits in 
\begin{eqnarray*}
\xymatrix{
1 \ar[r] & W \ar[d]_\phi\ar[r]^{j_1} & E_1 \ar[d]_\epsilon\ar[r]^{\pi_1} & G \ar[d]^i\ar[r] & 1  \\
1 \ar[r] & V \ar[r]_{j_0}            & E_0 \ar[r]_{\pi_0}                & H \ar[r]         & 1, 
}
\end{eqnarray*}
where the top and bottom rows are short exact sequences and the squares are maps of Lie $2$-groups.
\end{Def}
\begin{prop}\label{Ind2GpRep}
Given a split Lie $2$-group extension of $\xymatrix{G \ar[r]^{i} & H}$ by a $2$-vector space $\xymatrix{W \ar[r]^{\phi} & V}$,
\begin{eqnarray*}
\xymatrix{
1 \ar[r] & W \ar[d]_\phi\ar[r]^{j_1} & E_1 \ar[d]_\epsilon\ar[r]_{\pi_1} & G \ar[d]^i\ar[r]\ar@/_/[l]_{\sigma_1} & 1  \\
1 \ar[r] & V \ar[r]^{j_0}            & E_0 \ar[r]_{\pi_0}                & H \ar[r]\ar@/_/[l]_{\sigma_1}         & 1, 
}
\end{eqnarray*}
and a smooth splitting $\sigma$, there is an induced $2$-representation $\xymatrix{\rho^{\epsilon}_{\sigma}:\G \ar[r] & GL(\phi)}$ given by
\begin{align*}
\rho^0_0(h)v     & =\sigma_0(h)v\sigma_0(h)^{-1}   & \rho^1_0(h)w & =w^{\sigma_0(h)^{-1}} \\
\rho_1(g)v       & =\sigma_1(g)^v\sigma_1(g)^{-1}, &              & 
\end{align*}
for $h\in H$, $v\in V$, $w\in W$ and $g\in G$.
\end{prop}
%Since the proof of Proposition \ref{Ind2GpRep} consists of straighforward computations, consider postponing it to the Appendix.\\ 
In order to write the proof in an easier manner, we introduce the following notation. What follows will also serve us to get at the right space of $2$-cocycles. Consider a $2$-extension as the one from the proposition and use the splitting to establish diffeomorphisms $H\times V\cong E_0$ and $G\times W \cong E_1$, given by
\begin{eqnarray*}
\xymatrix{
(z,a) \ar@{|->}[r] & a\sigma_k (z)
}
\end{eqnarray*}
for $k=0$ and $k=1$ respectively. Their inverses are
\begin{eqnarray*}
\xymatrix{
e \ar@{|->}[r] & (\pi_k (e),e\sigma_k (\pi_k (e))^{-1}).
}
\end{eqnarray*}
Upgrading these diffeomorphisms to Lie group isomorphisms, we get 
\begin{align*}
(h_0,v_0)(h_1,v_1) & \sim v_0\sigma_0(h_0)v_1\sigma_0(h_1) \\
				   & \sim (\pi_0(v_0\sigma_0(h_0)v_1\sigma_0(h_1)),v_0\sigma_0(h_0)v_1\sigma_0(h_1)\sigma_0(\pi_0 (v_0\sigma_0(h_0)v_1\sigma_0(h_1)))^{-1}).
\end{align*}
Here, $\sim$ stands for the corresponding images under the above isomorphisms. Evaluating,
\begin{align*}
\pi_0(v_0\sigma_0(h_0)v_1\sigma_0(h_1)) & =\cancel{\pi_0(v_0)}\pi_0(\sigma_0(h_0))\cancel{\pi_0(v_1)}\pi_0(\sigma_0(h_1)) \\
                                        & =h_0h_1,
\end{align*}
and
\begin{align*}
v_0\sigma_0(h_0)v_1\sigma_0(h_1)\sigma_0(h_0h_1)^{-1} & = v_0\sigma_0(h_0)v_1\sigma_0(h_0)^{-1}\sigma_0(h_0)\sigma_0(h_1)\sigma_0(h_0h_1)^{-1};
\end{align*}
thus, the multiplication gets defined by
\begin{eqnarray*}
(h_0,v_0)(h_1,v_1):=(h_0h_1,v_0+\rho_0^0(h_0)v_1+\omega(h_0,h_1)).
\end{eqnarray*}
Here, $\omega_0(h_0,h_1)$ is shorthand for $\sigma_0(h_0)\sigma_0(h_1)\sigma_0(h_0h_1)^{-1}$, $\rho_0^0$ is defined as in the previous proposition, and we replace the product notation by sum notation since each factor lies in $V$. This is nothing but the usual twisted semi-direct product $H{}_{\rho^0_0}\ltimes^{\omega_0} V$ from the theory of Lie group extensions. Hence, if one supposes conversely that the product was defined using an abstract $\omega_0\in C(H^2,V)$, one is bound to rediscover that in order for such product to be associative, $\omega_0$ needs to be a $2$-cocycle for the Lie group cohomology of $H$ with values in $\rho^0_0$. The general formula for the differential of the complex of Lie group cochains with values in a representation is an instance of the differential of groupoid cochains with values in a representation that we laid down in chapter \ref{preliminarieschapter}. \\
Analogously, one finds out that there is an isomorphism of Lie groups $E_1\cong G{}_{\rho^1_0\circ i}\ltimes^{\omega_1}W$, with $\omega_1(g_0,g_1)=\sigma_1(g_0)\sigma_1(g_1)\sigma_1(g_0g_1)^{-1}$. In contrast with the case of Lie $2$-algebras, this time around we will not be able to waive the necessity for $\omega_1$ using the rest of the structure of crossed module. \\
We now turn to the homomorphism $\epsilon$. Consider
\begin{align*}
\epsilon(w\sigma_1(g)) & =\phi(w)\epsilon(\sigma_1(g)),
\end{align*}
we use the inverse isomorphism to define the crossed module map. Since $\pi$ is a crossed module map,
\begin{align*}
\pi_0(\phi(w)\epsilon(\sigma_1(g))) & =\cancel{\pi_0(\phi(w))}\pi_0(\epsilon(\sigma_1(g)))) \\
                                    & =i(\pi_1(\sigma_1(g)))=i(g);
\end{align*}
thus, the crossed module map is
\begin{eqnarray*}
\xymatrix{
(g,w) \ar@{|->}[r] & (i(x),\phi(w)+\varphi(g)),
}
\end{eqnarray*}
where $\varphi(g):=\epsilon(\sigma_1(g))\sigma_0(i(g))^{-1}$. We repeat this argument one last time to get the action,
\begin{align*}
(w\sigma_1(g))^{v\sigma_0(h)} & =w^{v\sigma_0(h)}\sigma_1(g)^{v\sigma_0(h)} \\
                              & =w^{\sigma_0(h)}\sigma_1(g)^{v\sigma_0(h)}.
\end{align*}
Since $\pi$ is a crossed module map,
\begin{align*}
\pi_1\big{(}w^{\sigma_0(h)}\sigma_1(g)^{v\sigma_0(h)}\big{)} & =\pi_1\big{(}w^{\sigma_0(h)}\big{)}\pi_1\big{(}\sigma_1(g)^{v\sigma_0(h)}\big{)} \\
                                                 & =\cancel{\pi_1(w)}^{\pi_0(\sigma_0(h))}\pi_1(\sigma_1(g))^{\cancel{\pi_0(v)}\pi_0(\sigma_0(h))}=g^h, 
\end{align*}
and
\begin{align*}
w^{\sigma_0(h)}\sigma_1(g)^{v\sigma_0(h)}\sigma_1(g^h)^{-1} & =(w\sigma_1(g)^v)^{\sigma_0(h)}\sigma_1(g^h)^{-1} \\
															& =\big{(}w\sigma_1(g)^v\sigma_1(g)^{-1}\big{)}^{\sigma_0(h)}\sigma_1(g)^{\sigma_0(h)}\sigma_1\big{(}g^h\big{)}^{-1};
\end{align*}
thus, the action is given by the equation
\begin{eqnarray*}
(g,w)^{(h,v)}:=(g^h,\rho^1_0(h)^{-1}(w+\rho_1(g)v)+\alpha(h;g))
\end{eqnarray*}
where $\alpha(h;g):=\sigma_1(g)^{\sigma_0(h)}\sigma_1\big{(}g^h\big{)}^{-1}$.
\begin{proof}(of Proposition \ref{Ind2GpRep})
We make the computations necessary to prove that $\rho^{\epsilon}_{\sigma}$ is a $2$-representation.
\begin{itemize}
\item Well-defined: We use the exactness of the sequences to see that the maps land where they are supposed to.
\begin{eqnarray*}
\pi_0\big{(}\sigma_0(h)v\sigma_0(h)^{-1}\big{)}=\pi_0(\sigma_0(h))\pi_0(v)\pi_0\big{(}\sigma_0(h)^{-1}\big{)}=h1h^{-1}=1 & \Longrightarrow & \rho_0^0(h)v\in V, \\
\pi_1\big{(}w^{\sigma_0(h)^{-1}}\big{)}=\pi_1(w)^{\pi_0\big{(}\sigma_0(h)^{-1}\big{)}}=1^{h^{-1}}=1 & \Longrightarrow & \rho_0^1(h)w\in W, \\
\pi_1\big{(}\sigma_1(g)^v\sigma_1(g)^{-1}\big{)}=\pi_1(\sigma_1(g))^{\pi_1(v)}\pi_1\big{(}\sigma_1(g)^{-1}\big{)}=g^1g^{-1}=1 & \Longrightarrow & \rho_1(g)v\in W. 
\end{eqnarray*}
Further, thus defined, $\rho_0^0(h)\circ\phi=\phi\circ\rho_0^1(h)$ for each $h\in H$. Indeed,
\begin{align*}
\rho_0^0(h)(\phi(w)) & =\sigma_0(h)\phi(w)\sigma_0(h)^{-1} =\sigma_0(h)\epsilon(w)\sigma_0(h)^{-1} \\
                     & =\epsilon\big{(}w^{\sigma_0(h)^{-1}}\big{)} =\phi(\rho_0^1(h)w);                   
\end{align*}
in these equations, we used that $\epsilon\circ j_1=j_0\circ\phi$, that $w^{\sigma_0(y)}\in W$ and the equivariance for the crossed module map $\epsilon$. This proves that for all $h\in H$, $\rho_0(h):=(\rho_0^0(y),\rho_0^1(y))\in GL(\phi)_0$, as desired.\\
Finally, the components of $\rho$, when evaluated, are all smooth and
\begin{align*}
\rho_0^0(h)(v_1+v_2) & =\sigma_0(h)(v_1v_2)\sigma_0(h)^{-1} =\sigma_0(h)v_1\sigma_0(h)^{-1}+\sigma_0(h) v_2\sigma_0(h)^{-1}  \\
\rho_0^1(h)(w_1+w_2) & =(w_1w_2)^{\sigma_0(h)^{-1}} =w_1^{\sigma_0(h)^{-1}}+w_2^{\sigma_0(h)^{-1}} \\
\rho_1(g)(v_1+v_2)   & =\sigma_1(g)^{v_1v_2}\sigma_1(g)^{-1}=(\sigma_1(g)^{v_1})^{v_2}\sigma_1(g)^{-1} \\
 					 & =(\sigma_1(g)^{v_1}\sigma_1(g)^{-1})^{v_2}\sigma_1(g)^{v_2}\sigma_1(g)^{-1} =\rho_1(g)(v_1)+\rho_1(v_2);
\end{align*}
thus proving they are linear. It is probably helpful to point out that in these equations we regard $V$ and $W$ as being subgroups of $E_0$ and $E_1$ respectively; therefore, we can interchange the additive notation with the multiplicative notation. As a second consequence, the crossed module action when restricted to $W$ is the action of the $2$-vector space regarded as a Lie $2$-group, which is always trivial.
\item $\rho^0_0$ is a Lie group homomorphism: Let $h_1,h_2\in H$ and $v\in V$, then
\begin{align*}
\rho^0_0(h_1)\rho^0_0(h_2)v & =\sigma_0(h_1)(\rho^0_0(h_2)v)\sigma_0(h_1)^{-1} \\
							& =\sigma_0(h_1)\sigma_0(h_2)v\sigma_0(h_2)^{-1}\sigma_0(h_1)^{-1} \\
							& =\sigma_0(h_1)\sigma_0(h_2)\sigma_0(h_1h_2)^{-1}\sigma_0(h_1h_2)v\sigma_0(h_1h_2)^{-1}\sigma_0(h_1h_2)\sigma_0(h_2)^{-1}\sigma_0(h_1)^{-1} \\
							& =\omega_0(h_1,h_2)(\rho^0_0(h_1h_2)v)\omega_0(h_1,h_2)^{-1}.
\end{align*}
Since both $\omega_0(h_1,h_2),\rho^0_0(h_1h_2)v\in V$ and conjugation in a vector space is trivial, we get the desired result.
\item $\rho_0^1$ is a Lie group homomorphism:  Let $h_1,h_2\in H$ and $w\in W$, then
\begin{align*}
\rho^1_0(h_1)\rho^1_0(h_2)w & =(\rho^1_0(h_2)w)^{\sigma_0(h_1)^{-1}} \\
							& =w^{\sigma_0(h_2)^{-1}\sigma_0(h_1)^{-1}} \\
							& =\big{(}w^{\sigma_0(h_1h_2)^{-1}}\big{)}^{\sigma_0(h_1h_2)\sigma_0(h_2)^{-1}\sigma_0(h_1)^{-1}} \\
							& =(\rho^1_0(h_1h_2)w)^{\omega_0(h_1,h_2)^{-1}}.
\end{align*}
This time around, $\omega_0(h_1,h_2)\in V$ and $\rho^1_0(h_1h_2)w\in W$; hence, as pointed out before, the action is trivial and we get the desired equation.
\item $\rho_1$ is a Lie group homomorphism: Let $g_1,g_2\in G$ and $v\in V$, then
\begin{align*}
\rho_1(g_1)\odot\rho_1(g_2)v & =(\rho_1(g_1)+\rho_1(g_2)+\rho_1(g_1)\phi\rho_1(g_2))v \\
							 & =\sigma_1(g_1)^{v}\sigma_1(g_1)^{-1}\sigma_1(g_2)^{v}\sigma_1(g_2)^{-1}\rho_1(g_1)\epsilon(\sigma_1(g_2)^{v}\sigma_1(g_2)^{-1}) \\
							 & =\sigma_1(g_1)^{v}\sigma_1(g_1)^{-1}\sigma_1(g_2)^{v}\sigma_1(g_2)^{-1}\sigma_1(g_1)^{\epsilon\big{(}\sigma_1(g_2)^{v}\sigma_1(g_2)^{-1}\big{)}}\sigma_1(g_1)^{-1}.
\end{align*}               
Using Peiffer equation, we get 
\begin{eqnarray*}
\sigma_1(g_1)^{\epsilon\big{(}\sigma_1(g_2)^{v}\sigma_1(g_2)^{-1}\big{)}}=\big{(}\sigma_1(g_2)^{v}\sigma_1(g_2)^{-1}\big{)}^{-1}\sigma_1(g_1)\big{(}\sigma_1(g_2)^{v}\sigma_1(g_2)^{-1}\big{)};
\end{eqnarray*} 
thus yielding,
\begin{align*}
\rho_1(g_1)\odot\rho_1(g_2)v & =\sigma_1(g_1)^{v}\sigma_1(g_1)^{-1}\sigma_1(g_1)\sigma_1(g_2)^{v}\sigma_1(g_2)^{-1}\sigma_1(g_1)^{-1} \\
							 & =\sigma_1(g_1)^{v}\sigma_1(g_2)^{v}\sigma_1(g_2)^{-1}\sigma_1(g_1)^{-1} \\
							 & =\sigma_1(g_1)^{v}\sigma_1(g_2)^{v}\big{(}\sigma_1(g_1g_2)^{v}\big{)}^{-1}\sigma_1(g_1g_2)^{v}\sigma_1(g_1g_2)\sigma_1(g_1g_2)^{-1}\sigma_1(g_2)^{-1}\sigma_1(g_1)^{-1} \\
							 & =(\omega_1(g_1,g_2))^v(\rho_1(g_1g_2)v)\omega_1(g_1,g_2).
\end{align*}    
Applying the remarks above, given that each factor lies in $W$, both the action and the conjugation are trivial, implying the desired equality.       
\item $\rho_0\circ i=\Delta\circ\rho_1$: This equation breaks into two components, one in $GL(V)$ and one in $GL(W)$; namely, 
\begin{eqnarray*}
\rho_0^0(i(g))=I+\phi\circ\rho_1(g) & \textnormal{and} & \rho_0^1(i(g))=I+\rho_1(g)\circ\phi 
\end{eqnarray*} 
for each $g\in G$. These relations hold as,
\begin{align*}
(I+\phi\rho_1(g))v & =v+\phi\big{(}\sigma_1(g)^v\sigma_1(g)^{-1}\big{)} \\              
                   & =v\epsilon\big{(}\sigma_1(g)^v\big{)}\epsilon\big{(}\sigma_1(g)^{-1}\big{)} \\
                   & =vv^{-1}\epsilon(\sigma_1(g))v\epsilon(\sigma_1(g))^{-1} \\
                   & =\epsilon(\sigma_1(g))v\epsilon(\sigma_1(g))^{-1} \\
                   & =\epsilon(\sigma_1(g))\sigma_0(i(g))^{-1}\sigma_0(i(g))v\sigma_0(i(g))^{-1}\sigma_0(i(g))\epsilon(\sigma_1(g))^{-1} \\      
                   & =\cancel{\varphi(g)}(\rho_0^0(i(g))v)\cancel{\varphi(g)^{-1}}=\rho_0^0(i(g))v,
\end{align*}
and
\begin{align*}
(I+\rho_1(g)\phi)w & =w+\sigma_1(g)^{\phi(w)}\sigma_1(g)^{-1} \\
				   & =w\sigma_1(g)^{\epsilon(w)}\sigma_1(g)^{-1} \\
                   & =ww^{-1}\sigma_1(g)w\sigma_1(g)^{-1} \\
                   & =w^{\epsilon\big{(}\sigma_1(g)^{-1}\big{)}} \\
                   & =\big{(}w^{\sigma_0(i(g))^{-1}}\big{)}^{\sigma_0(i(g))\epsilon(\sigma_1(g))^{-1}} \\
                   & =(\rho_0^1(i(g))w)^{\varphi(g)^{-1}}=\rho_0^1(i(g))w,
\end{align*}
where the last line of each sequence of equalities follows from $\varphi(g)\in V$.
\item $\rho_1$ respects the actions: Let $g\in G$, $h\in H$ and $v\in V$, then
\begin{align*}
\rho_1(g)^{\rho_0(h)}v & =\rho_0^1(h)^{-1}\rho_1(g)\rho_0^0(h)v \\
                       & =\rho_0^1(h)^{-1}\rho_1(g)\big{(}\sigma_0(h)v\sigma_0(h)^{-1}\big{)} \\
                       & =\rho_0^1(h)^{-1}\big{(}\sigma_1(g)^{\sigma_0(h)v\sigma_0(h)^{-1}}\sigma_1(g)^{-1}\big{)} \\
                       & =\big{(}\sigma_1(g)^{\sigma_0(h)v\sigma_0(h)^{-1}}\sigma_1(g)^{-1}\big{)}^{\sigma_0(h)} \\
                       & =\sigma_1(g)^{\sigma_0(h)v}\big{(}\sigma_1(g)^{-1}\big{)}^{\sigma_0(h)} \\
                       & =\big{(}\sigma_1(g)^{\sigma_0(h)}\sigma_1(g^h)^{-1}\big{)}^v\sigma_1\big{(}g^h\big{)}^v\big{(}\sigma_1(g)^{\sigma_0(h)}\big{)}^{-1} \\
                       & =\alpha(h;g)^v\sigma_1\big{(}g^h\big{)}^v\sigma_1(g^h)^{-1}\sigma_1\big{(}g^h\big{)}\big{(}\sigma_1(g)^{\sigma_0(h)}\big{)}^{-1} \\
                       & =\alpha(h;g)^v\big{(}\rho_1(g^h)v\big{)}\alpha(h;g)^{-1}=\rho_1\big{(}g^h\big{)}v,
\end{align*}                               
where the last equality is again due the fact that each factor lies in $W$.
\end{itemize}
\end{proof}

Starting out with a split extension of a Lie $2$-group by a $2$-vector space and a given splitting, we spelled out the whole structure in terms of the induced $2$-representation of Proposition \ref{Ind2GpRep}, together with the maps $\omega_0$, $\omega_1$, $\varphi$ and $\alpha$. Again, in case such splitting can be taken to be a crossed module map, the maps $\omega_0$, $\omega_1$, $\varphi$ and $\alpha$ will all vanish and we recover the semi-direct product structure defined above. The next proposition explores the converse to this discussion; namely, what equations should these maps satisfy in order for the structure defined by these formulae to be a $2$-extension. 
\begin{prop}\label{Gp2-cocycles}
Let $\rho$ be a $2$-representation of $\xymatrix{G \ar[r]^i & H}$ on $\xymatrix{W \ar[r]^\phi & V}$. Given a $4$-tuple $(\omega_0,\omega_1,\alpha,\varphi)\in C(H^2,V)\oplus C(G^2,W)\oplus C(H\times G,W)\oplus C(G,V)$, it defines a $2$-extension  
\begin{eqnarray*}
\xymatrix{
1 \ar[r] & W \ar[d]_\phi\ar@{^{(}->}[r] & G{}_{\rho_0^1\circ i}\ltimes^{\omega_1}W \ar[d]_\epsilon\ar[r]^{\qquad pr_1} & G \ar[d]^i \ar[r] & 1  \\
1 \ar[r] & V \ar@{^{(}->}[r]            & H{}_{\rho_0^0}\ltimes^{\omega_0}V \ar[r]_{\qquad pr_1}                       & H  \ar[r]          & 1, 
}
\end{eqnarray*}
with 
\begin{eqnarray*}
\epsilon(g,w) & = & (i(g),\phi(w)+\varphi(g)),                     \\
(g,w)^{(h,v)} & = & (g^h,\rho^1_0(h)^{-1}(w+\rho_1(g)v)+\alpha(h;g))
\end{eqnarray*} 
if, and only if the following equations are satisfied 
\begin{itemize}
\item[i)] $\delta\omega_0 =0$. Explicitly, for all triples $h_0,h_1,h_2\in H$, 
\begin{eqnarray*}
\rho^0_0(h_0)\omega_0(h_1,h_2)-\omega(h_0h_1,h_2)+\omega(h_0,h_1h_2)-\omega(h_0,h_1)=0.
\end{eqnarray*} 
\item[ii)] $\delta\omega_1 =0$. Explicitly, for all triples $g_0,g_1,g_2\in G$, 
\begin{eqnarray*}
\rho^1_0(i(g_0))\omega_1(g_1,g_2)-\omega(g_0g_1,g_2)+\omega(g_0,g_1g_2)-\omega(g_0,g_1)=0.
\end{eqnarray*}  
\item[iii)] $\phi(\omega_1(g_1,g_2))-\omega_0(i(g_1),i(g_2))=\rho^0_0(i(g_1))\varphi(g_2)-\varphi(g_1g_2)+\varphi(g_1)$
\item[iv)] $\rho^1_0(h_1h_2)^{-1}\rho_1(g)\omega_0(h_1,h_2)=\rho^1_0(h_2)^{-1}\alpha(h_1;g)-\alpha(h_1h_2;g)+\alpha(h_2;g^{h_1})$
\item[v)] $\varphi(g^h)-\rho^0_0(h^{-1})\varphi(g)+\phi(\alpha(h;g))=\rho^0_0(h^{-1})\omega_0(i(g),h)+\omega_0(h^{-1},i(g)h)-\omega_0(h^{-1},h)$
\item[vi)] $\rho^1_0(i(g_2))^{-1}\rho_1(g_1)\varphi(g_2)+\alpha(i(g_2);g_1)=\rho^1_0(i(g_2))^{-1}\omega_1(g_1,g_2)+\omega_1(g_2^{-1},g_1g_2)-\omega_1(g_2^{-1},g_2)$
\item[vii)] $\rho^1_0(h)^{-1}\omega_1(g_1,g_2)-\omega_1(g_1^h,g_2^h)=\rho^1_0(i(g_1^h))\alpha(h;g_2)-\alpha(h;g_1g_2)+\alpha(h;g_1)$ 
\end{itemize}
\end{prop}
There is nothing to these equations. We will make the computations necessary to prove that a $4$-tuple $(\omega_0,\omega_1,\alpha,\varphi)$ subject to the equations in the statement defines a $2$-extension. %For this proposition and its proof it would be useful to have some tool for breaking equations. Two options are the package breaqn in the ambient begin{dmath}, or \emph{multline} breaking them manually.
\begin{proof}
First, the usual theory of Lie group extensions, tells us that items $i)$ and $ii)$ say that $\omega_0$ and $\omega_1$ are $2$-cocycles with values in the representation $\rho_0^0$ and $\rho^1_0\circ i$ respectively. As such, these first two equations say that $H{}_{\rho^0_0}\ltimes^{\omega_0}V$ and $G{}_{\rho^1_0\circ i}\ltimes^{\omega_1}W$ are Lie groups, specifically that each of their product rules is associative. The other equations in the statement have the following meaning
\begin{itemize}
\item[iii)] says $\epsilon$ is a Lie group homomorphism:
\begin{align*}
\epsilon((g_1,w_1)\vJoin_{\omega_1}(g_2,w_2)) & =\epsilon(g_1g_2,w_1+\rho^1_0(i(g_1))w_2+\omega_1(g_1,g_2)) \\
											  & =(i(g_1g_2),\phi(w_1+\rho^1_0(i(g_1))w_2+\omega_1(g_1,g_2))+\varphi(g_1g_2)),	  
\end{align*} 
and
\begin{align*}
\epsilon(g_1,w_1)\epsilon(g_2,w_2) & =(i(g_1),\phi(w_1)+\varphi(g_1))(i(g_2),\phi(w_2)+\varphi(g_2)) \\
								   & =(i(g_1)i(g_2),\phi(w_1)+\varphi(g_1)+\rho^0_0(i(g_1))(\phi(w_2)+\varphi(g_2))+\omega_0(i(g_1),i(g_2))).	  
\end{align*} 
Since $i$ is a Lie group homomorphism, $\phi$ is linear and $\phi\circ\rho^0_0=\rho^1_0\circ\phi$, the expressions above are equal if, and only if the equation in item $iii)$ is satisfied.
\item[iv)] says that the formula for $(g,w)^{(h,v)}$ defines a right action: 
\begin{align*}
(g,w)^{(h_1,v_1)(h_2,v_2)} & =(g,w)^{(h_1h_2,v_1+\rho^0_0(h_1)v_2+\omega_0(h_1,h_2))} \\
                           & =(g^{h_1h_2},\rho_0^1(h_1h_2)^{-1}(w+\rho_1(g)(v_1+\rho^0_0(h_1)v_2+\omega_0(h_1,h_2)))+\alpha(h_1h_2;g)),
\end{align*}
and
\begin{align*}
\big{(}(g,w)^{(h_1,v_1)}\big{)}^{(h_2,v_2)} & =\big{(}g^{h_1},\rho^1_0(h_1)^{-1}(w+\rho_1(v_1))+\alpha(h_1;g)\big{)}^{(h_2,v_2)} \\
                                & =\big{(}\big{(}g^{h_1}\big{)}^{h_2},\rho^1_0(h_2)^{-1}\big{(}\rho^1_0(h_1)^{-1}(w+\rho_1(g)v_1)+\alpha(h_1;g)+\rho_1\big{(}g^{h_1}\big{)}v_2\big{)}+\alpha(h_2;g^{h_1})\big{)}.
\end{align*}
Since $g^h$ is a right action, $\rho^1_0$ is a representation, and $\rho$ is a map of crossed modules, the expressions above are equal if, and only if the equation in item $iv)$ is satisfied. 
\item[v)] says $\epsilon$ is equivariant. 
\begin{align*}
\epsilon\big{(}(g,w)^{(h,v)}\big{)} & =\epsilon\big{(}g^h,\rho_0^1(h)^{-1}(w+\rho_1(g)v)+\alpha(h;g)\big{)} \\
						& =\big{(}i(g^h),\phi(\rho_0^1(h)^{-1}(w+\rho_1(g)v)+\alpha(h;g))+\varphi(g^h)\big{)},
\end{align*}
and
\begin{align*}
(h,v) & ^{-1}\epsilon(g,w)(h,v) \\
           & =(h^{-1},-\rho^0_0(h^{-1})v-\omega_0(h^{-1},h))(i(g),\phi(w)+\varphi(g))(h,v) \\
		   & =(h^{-1},-\omega_0(h^{-1},h)-\rho^0_0(h^{-1})v)(i(g)h,\phi(w)+\varphi(g)+\rho^0_0(i(g))v+\omega_0(i(g),h)) \\
           & =\big{(}h^{-1}i(g)h,\rho^0_0(h^{-1})\big{(}\phi(w)-v+\varphi(g)+\rho^0_0(i(g))v+\omega_0(i(g),h)\big{)}+ \omega_0(h^{-1},h)+\omega_0(h^{-1},i(g)h)\big{)}.
\end{align*}
Since $i$ is equivariant, $\phi\circ\rho^0_0=\rho^1_0\circ\phi$, and $\rho^0_0(i(g))=I+\phi\rho_1(g)$, the expressions above are equal if, and only if the equation in item $v)$ is satisfied. 
\item[vi)] says that the Peiffer equation is satisfied.
\begin{align*}
(g_1,w_1)^{\epsilon(g_2,w_2)} & =(g_1,w_1)^{(i(g_2),\phi(w_2)+\varphi(g_2))} \\
						      & =\big{(}g_1^{i(g_2)},\rho_0^1(i(g_2))^{-1}\big{(}w_1+\rho_1(g_1)(\phi(w_2)+\varphi(g_2))\big{)}+\alpha(i(g_2);g_1)\big{)},
\end{align*}
and
\begin{align*}
(g_2,w_2) & ^{-1}\vJoin_{\omega_1}(g_1,w_1)\vJoin_{\omega_1}(g_2,w_2) \\
		  & =\big{(}g_2^{-1},-\rho^1_0(i(g_2^{-1}))w_2-\omega_1(g_2^{-1},g_2)\big{)}\vJoin_{\omega_1}\big{(}g_1g_2,w_1+\rho^1_0(i(g_1))w_2+\omega_1(g_1,g_2)\big{)} \\
		  & =\big{(}g_2^{-1}g_1g_2,\rho^1_0(i(g_2^{-1}))\big{(}w_1-w_2+\rho^1_0(i(g_1))w_2+\omega_1(g_1,g_2)\big{)}-\omega_1(g_2^{-1},g_2)+\omega_1(g_2^{-1},g_1g_2)\big{)}.
\end{align*}
Because of the Peiffer equation and of the relation $\rho^1_0(i(g))=I+\rho_1(g)\phi$, the expressions above are equal if, and only if the equation in item $vi)$ is satisfied. 
\item[vii)] says that the formula for $(g,w)^{(h,v)}$ defines an action by automorphisms.
\begin{align*}
((g_1,w_1) & \vJoin_{\omega_1}(g_2,w_2))^{(h,v)} \\
           & =\big{(}g_1g_2,w_1+\rho^1_0(i(g_1))w_2+\omega_1(g_1,g_2)\big{)}^{(h,v)} \\
		   & =\big{(}(g_1g_2)^h,\rho_0^1(h)^{-1}\big{(}w_1+\rho^1_0(i(g_1))w_2+\omega_1(g_1,g_2)+\rho_1(g_1g_2)v\big{)}+\alpha(h;g_1g_2)\big{)},
\end{align*}
and
\begin{align*}
(g_1,w_1) & ^{(h,v)}\vJoin_{\omega_1}(g_2,w_2)^{(h,v)} \\
          & =\big{(}g_1^h,\rho^1_0(h)^{-1}(w_1+\rho_1(g_1)v)+\alpha(h;g_1)\big{)}\vJoin_{\omega_1} \big{(}g_2^h,\rho^1_0(h)^{-1}(w_2+\rho_1(g_2)v)+\alpha(h;g_2)\big{)} \\
          & =\big{(}g_1^hg_2^h,\rho^1_0(h)^{-1}(w_1+\rho_1(g_1)v)+\alpha(h;g_1)+ \\
          & \qquad\qquad+\rho^1_0(i(g_1^h))\big{(}\rho^1_0(h)^{-1}(w_2+\rho_1(g_2)v)+\alpha(h;g_2)\big{)}+\omega_1(g_1^h,g_2^h)\big{)}.
\end{align*}
Since the action of $H$ on $G$ is by Lie group automorphisms, $i$ is equivariant, $\rho_1$ is a Lie group homomorphism and $\rho^1_0(i(g))=I+\rho_1(g)\phi$, the expressions above are equal if, and only if the equation in item $vii)$ is satisfied.
\end{itemize}
\end{proof}
Next, we analyze what happens when we have got equivalent extensions. 
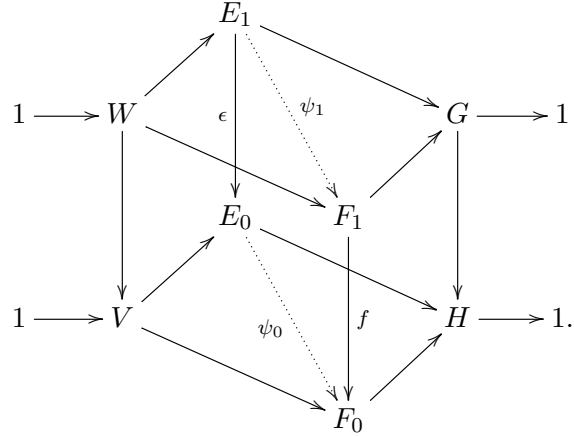
\begin{eqnarray*}
\xymatrix{
         &                          & E_1 \ar[dd]_\epsilon\ar@{.>}[ddr]^{\psi_1}\ar[drr] &             &                  &    \\
1 \ar[r] & W \ar[dd]\ar[ur]\ar[drr] &                                           &                      & G \ar[r] \ar[dd] & 1  \\
         &                          & E_0 \ar[drr]\ar@{.>}[ddr]_{\psi_0}        & F_1 \ar[dd]^f\ar[ur] &                  &    \\
1 \ar[r] & V \ar[ur]\ar[drr]        &                                           &                      & H \ar[r]         & 1. \\
         &                          &                                           & F_0 \ar[ur]          &                  &   
}
\end{eqnarray*}
We need to suppose that these are split extensions to apply the decomposition that we outlined above. As in the case of Lie $2$-algebras, it suffices to pick a splitting of either extension to induce a splitting of the other one. In this manner, the induced $2$-representations are the same. We use these compatible splittings to identify both $\mathcal{E}$ and $\mathcal{F}$ with their respective semi-direct products, and we write $\psi$ in these coordinates 
\begin{eqnarray*}
\psi_k(z,a)=(\psi_k^{Gp}(z,a),\psi_k^{Vec}(z,a)).
\end{eqnarray*}
Since both $\psi_0$ and $\psi_1$ respect inclusions and projections, $\psi_k(1,a)=(1,a)$ and $\psi_k^{Gp}(z,a)=z$. Therefore, using the factorization $(z,a)=(1,a)\vJoin(z,0)$ and the fact that both $\psi_k$'s are group homomorphisms,
\begin{align*}
\psi_k(z,a) & =\psi_k((1,a)\vJoin(z,0)) \\
			& =\psi_k(1,a)\vJoin\psi_k(z,0)\\
			& =(1,a)\vJoin(z,\psi_k^{Vec}(z,0)) \\
			& =(z,a+\cancel{\rho_k(1)}(\psi_k^{Vec}(z,0))+\cancel{\omega_k(1,z)}).
\end{align*}
As a consequence,
\begin{eqnarray*}
\psi_k(z,a)=(z,a+\lambda_k(z)),
\end{eqnarray*}
where the maps $\xymatrix{\lambda_0:H \ar[r] & V}$ and $\xymatrix{\lambda_1:G \ar[r] & W}$ are defined by $\psi_k^{Vec}(z,0)$ for $k=0$ and $k=1$ respectively.
\begin{prop}\label{Gp2-coboundaries}
Let $\rho$ be a $2$-representation of $\xymatrix{G \ar[r]^i & H}$ on $\xymatrix{W \ar[r]^\phi & V}$. Given two $2$-cocycles  $(\omega_0,\omega_1,\alpha,\varphi)$, $(\omega_0',\omega_1',\alpha',\varphi')$ as in the previous proposition, the induced extensions are equivalent if, and only if there are maps $\lambda_0\in C(H,V)$ and $\lambda_1\in C(G,W)$ verifying  
\begin{itemize}
\item[i)] $\omega_0-\omega_0'=\delta\lambda_0$. Explicitly, for all pairs $h_0,h_1\in H$, 
\begin{eqnarray*}
\omega_0(h_0,h_1)-\omega_0'(h_0,h_1)=\rho^0_0(h_0)\lambda_0(h_1)-\lambda_0(h_0h_1)+\lambda_0(h_0)
\end{eqnarray*}  
\item[ii)] $\omega_1-\omega_1'=\delta\lambda_1$. Explicitly, for all pairs $g_0,g_1\in G$, 
\begin{eqnarray*}
\omega_1(g_0,g_1)-\omega_1'(g_0,g_1)=\rho^1_0(i(g_0))\lambda_1(g_1)-\lambda_1(g_0g_1)+\lambda_1(g_0)
\end{eqnarray*}
\item[iii)] $\varphi(g)-\varphi'(g)=\phi(\lambda_1(g))-\lambda_0(i(g))$
\item[iv)] $\alpha(h;g)-\alpha'(h;g)=\rho^1_0(h)^{-1}(\lambda_1(g)+\rho_1(g)\lambda_0(h))-\lambda_1(g^h)$
\end{itemize}
\end{prop}
\begin{proof}
Items $i)$ and $ii)$ are the usual relations for two group cocycles to be cohomologous. As for the other two, we define $\psi_k(z,a):=(z,a+\lambda_k(z))$, and we ask for $\psi$ to be a map of crossed modules. Thus, it has got to commute with the structural maps and respect the actions. In symbols,
\begin{eqnarray*}
\xymatrix{
G{}_{\rho^1_0\circ i}\ltimes^{\omega_1}W \ar[r]^{\psi_1}\ar[d] & G{}_{\rho^1_0\circ i}\ltimes^{\omega_1'}W \ar[d] \\
H{}_{\rho^0_0}\ltimes^{\omega_0}V \ar[r]_{\psi_0}              & H{}_{\rho^0_0}\ltimes^{\omega_0'}V,
}
\end{eqnarray*}
and $\psi_1\big{(}(g,w)^{(h,v)}\big{)}=\psi_1(g,w)^{\psi_0(h,v)}$. The commutativity of the square reads
\begin{eqnarray*}
\big{(}i(g),\phi(w+\lambda_1(g))+\varphi'(g)\big{)}=\big{(}i(g),\phi(w)+\varphi(g)+\lambda_0(i(g))\big{)},
\end{eqnarray*}
which is clearly equivalent to the equation in item $iii)$. On the other hand, 
\begin{align*}
\psi_1\big{(}(g,w)^{(h,v)}\big{)} & =\psi_1\big{(}g^h,\rho^1_0(h)^{-1}(w+\rho_1(g)v)+\alpha(h;g)\big{)} \\
					  & =\big{(}g^h,\rho^1_0(h)^{-1}(w+\rho_1(g)v)+\alpha(h;g)+\lambda_1(g^h)\big{)},
\end{align*}
and
\begin{align*}
\psi_1(g,w)^{\psi_0(h,v)} & =\big{(}g,w+\lambda_1(g)\big{)}^{(h,v+\lambda_0(h))} \\
						  & =\big{(}g^h,\rho^1_0(h)^{-1}\big{(}w+\lambda_1(g)+\rho_1(g)(v+\lambda_0(h))\big{)}+\alpha'(h;g)\big{)}
\end{align*}
Clearly, these expressions are equal if, and only if the equation in item $iv)$ holds, thus completing the proof.

\end{proof}

\section{The complex of Lie 2-group cochains with values in a 2-representation}
We want to interpret the equation of propositions \ref{Gp2-cocycles} and \ref{Gp2-coboundaries} as the cocycle equations for some complex. At this point, we would like to remark that that proposition \ref{Gp2-cocycles} applied to the representation on $\xymatrix{(0) \ar[r] & \Rr}$ for which the only surviving component is the identity is equivalent to lemma \ref{easyExt}. Indeed, for this particular, case the only non-trivial equations in proposition \ref{Gp2-cocycles} are 
\begin{itemize}
    \item[i)] $\delta\omega_0=0$,
    \item[ii)] $-\omega_0(i(g_1),i(g_2))=\varphi(g_2)-\varphi(g_1g_2)+\varphi(g_1)$ and
    \item[iii)] $\varphi(g^h)-\varphi(g)=\omega_0(i(g),h)+\omega_0(h^{-1},i(g)h)-\omega_0(h^{-1},h)$.
\end{itemize}
Assuming the hypothesis $d(F,f)=0$ and under the identifications $\omega_0=F$ and $\varphi(g)=-f(g,1)$, we have got the first equation right away; the second equation is the one we singled out, evaluating $\partial F+\delta f=0$ at $(\gamma_1,\gamma_2)$, with $\gamma_k=(g_k,1)$. Finally, the third equations is implied from this latter equation, together with $\delta F=0$ in the following way. Evaluating $\partial F+\delta f=0$ at
\begin{align*}
    \begin{pmatrix}
    g & 1 \\
    h^{-1} & h
    \end{pmatrix}: & & f(g,h^{-1})-f(g^h,1)+\cancel{f(1,h)} & =F(h^{-1}i(g),h)-F(h^{-1},h), \\
    \begin{pmatrix}
    1 & g \\
    h^{-1} & 1
    \end{pmatrix}: & & \cancel{f(1,h^{-1})}-f(g,h^{-1})+f(g,1) & =F(h^{-1},i(g))-\cancel{F(h^{-1},1)} \\
\end{align*}
and summing these expressions yields,
\begin{eqnarray*}
f(g,1)-f(g^h,1)=F(h^{-1}i(g),h)-F(h^{-1},h)+F(h^{-1},i(g)),
\end{eqnarray*}
but $\delta F=0$, then
\begin{eqnarray}
F(h^{-1}i(g),h)+F(h^{-1},i(g))=F(i(g),h)+F(h^{-1},i(g)h).
\end{eqnarray}
Thus, replacing the identifications for $\omega_0$ and $\varphi$, there appears the third equation in the list above. Conversely, define $F=\omega_0$ and $f(g,h):=-\omega_0(h,i(g))-\varphi(g)$ and assume the equations in proposition \ref{Gp2-cocycles} hold. Naturally, the cocycle equation for $F$ holds right away. As for the equation $\partial f=0$, 
\begin{align*}
f(g_2g_1) & =-\omega_0(h,i(g_2g_1))-\varphi(g_2g_1) \\
          & =-\omega_0(h,i(g_2g_1))-\omega_0(i(g_2),i(g_1))-\varphi(g_2)-\varphi(g_1) \\
          & =-\omega_0(hi(g_2),i(g_1))-\phi(g_1)-\omega_0(h,i(g_2))-\phi(g_2)=f(g_1,hi(g_2))+f(g_2,h),
\end{align*}
where we used equation $iii)$ in the list above to pass to the second line and the cocycle equation $\delta\omega_0=0$ to pass to the third line. Finally, for $\partial F+\delta f=0$, we compute
\begin{align*}
    -f(g_1^{h_2}g_2,h_1h_2) & =\omega_0(h_1h_2,i(g_1^{h_2}g_2))+\varphi(g_1^{h_2}g_2) \\
                            & =\omega_0(h_1h_2,i(g_1^{h_2}g_2))+\omega_0(i(g_1^{h_2}),i(g_2))+\varphi(g_1^{h_2})+\varphi(g_2) \\
                            & =\omega_0(h_1h_2i(g_1^{h_2}),i(g_2))+\omega_0(h_1h_2,i(g_1^{h_2}))+ \\
                            & \qquad +\omega_0(i(g_1),h_2)+\omega_0(h_2^{-1},i(g_1)h_2)-\omega_0(h_2^{-1},h_2)+\varphi(g_1)+\varphi(g_2).
\end{align*}
 Here, we used again equation $ii)$ in the list above to pass to the second line and equation $iii)$ together with the cocycle equation $\delta\omega_0=0$ to pass to the third line. Summing this expression with
 \begin{eqnarray}
 f(g_2,h_2)+f(g_1,h_1)=-\omega_0(h_2,i(g_2))-\varphi(g_2)-\omega_0(h_1,i(g_1))-\varphi(g_1), 
 \end{eqnarray}
 we get
 \begin{align*}
     \delta f\begin{pmatrix}
     g_1 & g_2 \\
     h_1 & h_2
     \end{pmatrix} & =\omega_0(h_1i(g_1)h_2,i(g_2))+\omega_0(h_1h_2,i(g_1^{h_2}))+ \\
                   & \qquad +\omega_0(i(g_1),h_2)+\omega_0(h_2^{-1},i(g_1)h_2)-\omega_0(h_2^{-1},h_2)-\omega_0(h_2,i(g_2))-\omega_0(h_1,i(g_1)),
 \end{align*}
where, in writing the first term, we used the equivariance of $i$. Now, we are going to use the cocycle equation $\delta\omega_0=0$ successively to get the desired result. Indeed, 
\begin{align*}
     \delta\omega_0(h_1i(g_1),h_2,i(g_2)) & =\omega_0(h_2,i(g_2))-\omega_0(h_1i(g_1)h_2,i(g_2))+\omega_0(h_1i(g_1),h_2i(g_2))-\omega_0(h_1i(g_1),h_2)=0, \\
     \delta\omega_0(h_1,i(g_1),h_2)       & =\omega_0(i(g_1),h_2)-\omega_0(h_1i(g_1),h_2)+\omega_0(h_1,i(g_1)h_2)-\omega_0(h_1,i(g_1))=0, \\
\delta\omega_0(h_1h_2,h_2^{-1},i(g_2)h_2) & =\omega_0(h_2^{-1},i(g_2)h_2)-\omega_0(h_1,i(g_2)h_2)+\omega_0(h_1h_2,i(g_1^{h_2})i(g_2))-\omega_0(h_1h_2,h_2^{-1})=0, \\
     \delta\omega_0(h_1,h_2,h_2^{-1})     & =\omega_0(h_2,h_2^{-1})-\omega_0(h_1h_2,h_2^{-1})-\omega_0(h_1,h_2)=0, \\
    \delta\omega_0(h_2^{-1},h_2,h_2^{-1}) & =\omega_0(h_2,h_2^{-1})-\omega_0(h_2^{-1},h_2)=0; \\
\end{align*}
therefore, 
\begin{align*}
\delta f\begin{pmatrix}
 g_1 & g_2 \\
 h_1 & h_2
\end{pmatrix} & =\omega_0(h_1i(g_1),h_2i(g_2))+\big{(}-\omega_0(h_1i(g_1),h_2)+\omega_0(i(g_1),h_2)-\omega_0(h_1,i(g_1))\big{)}+ \\
              & \qquad +\big{(}\omega_0(h_2^{-1},i(g_1)h_2)-\omega_0(h_1h_2,i(g_1^{h_2}))\big{)}-\omega_0(h_2^{-1},h_2) \\
              & =\omega_0(h_1i(g_1),h_2i(g_2))-\omega_0(h_1,i(g_1)h_2)+\omega_0(h_1,i(g_1)h_2)+\omega_0(h_1h_2,h_2^{-1})-\omega_0(h_2^{-1},h_2) \\
              & =\omega_0(h_1i(g_1),h_2i(g_2))+\omega_0(h_2,h_2^{-1})-\omega_0(h_1,h_2)-\omega_0(h_2^{-1},h_2)=-\partial\omega_0\begin{pmatrix}
 g_1 & g_2 \\
 h_1 & h_2
\end{pmatrix}.
 \end{align*}
Of course, it is also true that restricted to this case, proposition \ref{Gp2-coboundaries} relating isomorphic extensions in terms of their cocycles is equivalent to lemma \ref{isoEasyExt}. In sight of this correspondence, one would expect the equations defining a cocycle to come from the total complex of a double complex. However, although eventually we will be using a type of double complex, the underlying object that we will come across first is a kind of triple complex, as it was the case for Lie $2$-algebras. Again, not all the faces of this triple complex will commute, but they will do so up to homotopy, and the homotopies relating them will be encoded in a series of additional maps.  \\
%in fact the columns of such complex are total complexes of double complexes themselves and there is an induced double complex... This seems to be false...
Suppose now that $W=(0)$ or equivalently that the $2$-representation takes values on an honest vector space, then we get a double complex analogous to that of Section \ref{GpDcx}
\begin{eqnarray}\label{page r=0}
\xymatrix{
\vdots                                  & \vdots                                    & \vdots                           &       \\ 
C(H^3,V) \ar[r]^{\partial}\ar[u]        & C(\G ^3,V) \ar[r]^{\partial}\ar[u]        & C(\G _2^3,V) \ar[r]\ar[u]        & \dots \\
C(H^2,V) \ar[r]^{\partial}\ar[u]^\delta & C(\G ^2,V) \ar[r]^{\partial}\ar[u]^\delta & C(\G _2^2,V) \ar[r]\ar[u]^\delta & \dots \\
C(H,V) \ar[r]^{\partial}\ar[u]^\delta   & C(\G ,V) \ar[r]^{\partial}\ar[u]^\delta   & C(\G _2,V) \ar[r]\ar[u]^\delta   & \dots \\
V \ar[r]^{0}\ar[u]^\delta               & V \ar[r]^{I_V}\ar[u]^\delta               & V	\ar[r]^{0}\ar[u]^\delta        & \dots  
}
\end{eqnarray} 
From example \ref{unitGpRep}, we know that a $2$-representation on $\xymatrix{(0)\ar[r] & V}$ amounts to a usual representation $\rho$ of $H/i(G)$ on $V$. In order to get a meaningful double complex, there needs to be induced representations of $\G _p$ on $V$. The first column has a single candidate, namely, the pull-back representation of $\rho$ along the projection, $\rho_0$. However, we face multiple possible choices for the other columns. Again, we are going to fix the representation of the $p$th column, $\rho_p$, to be the pull-back of $\rho_0$ along $\xymatrix{t_p:\G _p \ar[r] & H :(g_1,...,g_p,h) \ar@{|->}[r] & hi(g_p...g_1)}$. In so, each column is a complex of  Lie group cochains with values in $\rho_p$. It is of crucial importance that all of these representations vanish on the ideal $i(G)$; otherwise, the squares in the diagram above would not commute, thus failing to be a double complex.
\begin{remark}
Notice that, precisely because the representation vanishes on $i(G)$, taking the pull-back along the initial source map, $s_p=t_p\circ\iota$, would not affect the commutativity of the complex.
\end{remark}
In the general case, when the $2$-representation takes values on a general $2$-vector space, the squares in the diagram above cease to commute. Interestingly though, when applied to a $q$-tuple, the elements of $V$ that result of going through each of the sides of a square are isomorphic in the $2$-vector space. In symbols, given a cochain $\varphi\in C(\G _{p-1}^{q-1},V)$ and a $q$-tuple $(\gamma_1,...,\gamma_q)\in\G _{p}^{q}$,
\begin{eqnarray*}
\delta\partial\varphi(\gamma_1,...,\gamma_q)\cong\partial\delta\varphi(\gamma_1,...,\gamma_q) & \textnormal{in }V.
\end{eqnarray*} 
In fact, by looking at the difference of the two ways to go around, one sees that not only are $\delta\partial\varphi$ and $\partial\delta\varphi$ isomorphic, but they are so in a coherent way. That is, the isomorphism joining them goes through a second complex. For instance, for $v\in V$ and $(g,h)\in G\rtimes H\cong\G$,
\begin{align*}
(\partial\delta-\delta\cancel{\partial})v(g,h) & =\delta v(h)-\delta v(hi(g)) \\
                                               & =\rho_0^0(h)v-v-(\rho_0^0(hi(g))v-v) \\
                                               & =\rho_0^0(h)(v-\rho_0^0(i(g))v) \\
                                               & =\rho_0^0(h)(v-v-\phi\rho_1(g)v)=-\rho_0^0(h)\phi\rho_1(g)v.
\end{align*}
Also, for $(g_1,g_2,h)\in G^2\times H\cong\G _2$,
\begin{align*}
(\partial\delta-\delta\partial)v(g_1,g_2,h) & =\delta v(g_2,h)-\delta v(g_2g_1,h)+\delta v(g_1,hi(g_2))-(\rho_0^0(hi(g_2g_1))v-v) \\
                                            & =\rho_0^0(hi(g_2))v-v-(\rho_0^0(hi(g_2g_1))v-v)+ \\
                                            & \qquad\qquad +\rho_0^0(hi(g_2)i(g_1))v-v-\rho_0^0(hi(g_2g_1))v+v \\
                                            & =\rho_0^0(hi(g_2))(v-\rho_0^0(i(g_1))v) \\
                                            & =\rho_0^0(hi(g_2))(v-v-\phi\rho_1(g_1)v)=-\rho_0^0(hi(g_2))\phi\rho_1(g_1)v.
\end{align*}
These computations prompt the following lemma.
\begin{lemma}\label{upToHomotopy0}
For $v\in V$ and $(g_0,...,g_p,h)\in G^{p+1}\times H\cong\G_{p+1}$,
\begin{eqnarray*}
\delta\partial v(g_0,...,g_p,h)=\partial\delta v(g_0,...,g_p,h)+\rho_0^0(t_p(g_1,...,g_p,h))\phi\rho_1(g_0)v
\end{eqnarray*}
\end{lemma}
\begin{proof}
First, let us compute the first term on the right hand side of the equation:
\begin{align*}
\partial\delta v(g_0,...,g_p,h) & =\delta v(g_1,...,g_p,h)+ \\
							    & \qquad +\sum_{j=0}^{p-1}(-1)^{j+1}\delta v(g_0,...,g_{j+1}g_j,...,g_p,h)+(-1)^{p+1}\delta v(g_0,...,g_{p-1},hi(g_p)) \\
								& =\rho_0^0(hi(g_p...g_1))v-v+ \\
								& \qquad +\sum_{j=0}^{p-1}(-1)^{j+1}(\rho_0^0(hi(g_p...g_0))v-v)+(-1)^{p+1}(\rho_0^0(hi(g_p)i(g_{p-1}...g_0))v-v).
\end{align*}
This divides then in two cases. If $p$ is even,
\begin{align*}
\partial\delta v(g_0,...,g_p,h) & =\rho_0^0(hi(g_p...g_1))(v-\rho_0^0(i(g_0))v) \\
								& =\rho_0^0(hi(g_p...g_1))(v-v-\phi\rho_1(g_0)v) =-\rho_0^0(hi(g_p...g_1))\phi\rho_1(g_0)v;
\end{align*}
otherwise,
\begin{align*}
\partial\delta v(g_0,...,g_p,h) & =\rho_0^0(hi(g_p...g_1))v-v.
\end{align*}
Now, the left hand side also divides in two cases. If $p$ is even, $\partial v=0$ and
\begin{align*}
\delta\partial v(g_0,...,g_p,h) & =\rho_0^0(hi(g_p...g_0))\cancel{\partial v}-\cancel{\partial v}=0;
\end{align*}
otherwise, $\partial v=v$ and
\begin{align*}
\delta\partial v(g_0,...,g_p,h) & =\rho_0^0(hi(g_p...g_0))v-v.
\end{align*}
\end{proof}
Before moving on, a couple of comments are in order. First, since $\rho$ is a $2$-representation, $\rho_0^0\circ\phi=\phi\circ\rho_0^1$; hence, lemma \ref{upToHomotopy0} indeed says that there is an isomorphism in the $2$-vector space as claimed. Secondly, we can schematize the relation in lemma \ref{upToHomotopy0} with the following diagram
\begin{eqnarray*}
\xymatrix{
C(\G_p,V) \ar[rr]^\partial                             &                       & C(\G_{p+1},V)    \\
                                                       & C(G,W) \ar[ur]^\Delta &                  \\
V \ar[uu]^\delta\ar[ur]^{\delta_{(1)}}\ar[rr]_\partial &                       & V \ar[uu]_\delta ,  
}
\end{eqnarray*}
where $\delta_{(1)} v(g):=\rho_1(g)v$ and $\Delta\lambda(g_0,...,g_p,h):=\rho_0^0(hi(g_p...g_1))\phi\lambda(g_0)$. Notice that when $p=0$, the formula for $\Delta$ can be read out to be the differential of the representation up to homotopy induced by the $2$-representation,
\begin{eqnarray*}
\xymatrix{
\eth :H\times W \ar[r] & H\times V:(h,w) \ar@{|->}[r] & (h,\rho_0^0(h)\phi(w)).
}
\end{eqnarray*}
In fact, $\lambda\in C(G,W)$ defines a map of bundles from the Lie group bundle $H\times G$ to the vector bundle of $H\times W$,
\begin{eqnarray*}
\xymatrix{
H\times G \ar[r]^{\bar{\lambda}}\ar[d]_{pr_1} & H\times W\ar[d] \\
H \ar@{=}[r]                                  & H
}
\end{eqnarray*}
given by $\bar{\lambda}(h;g):=(h,\lambda(g))$, which composed with $\eth$ yields $\Delta$. Furthermore, for any given value of $p$, $\Delta$ coincides with the structural map of the pull-back of this $2$-term complex along $t_p$ composed with the map of bundles induced by $\lambda$.\\
Finally, we justify the seeming overlap in notation for the diagonal map $\delta_{(1)}$ with the differential of the group cochain complex of $G$ with values in $\rho_0^1\circ i$.
\begin{lemma}\label{G-coh[V]}
Defined as above, $\delta_{(1)}$ makes 
\begin{eqnarray*}
\xymatrix{
V \ar[r]^{\delta_{(1)}\qquad} & C(G,W) \ar[r] & C(G^2,W) \ar[r] & ...
}
\end{eqnarray*}
into a complex.
\end{lemma}
\begin{proof}
The only thing one needs to prove is that, for $v\in V$, $\delta_{(1)}\delta_{(1)} v=0$.
\begin{align*}
\delta_{(1)}\delta_{(1)} v(g_1,g_2) & =\rho_0^1(i(g_1))\delta_{(1)} v(g_2)-\delta_{(1)} v(g_1g_2)+\delta_{(1)} v(g_1) \\
									& =\rho_0^1(i(g_1))\rho_1(g_2)v-\rho_1(g_1g_2)v+\rho_1(g_1)v \\
									& =(I+\rho_1(g_1)\phi)\rho_1(g_2)v-\rho_1(g_1g_2)v+\rho_1(g_1)v,
\end{align*}
and $\rho_1$ is a Lie group homomorphism with codomain $GL(\phi)_1$.

\end{proof}
The previous discussion serves as a full description of what is needed to ensure the commutativity of the first row of squares in the ``double complex'' above. Let us see what happens in the next row, let $\lambda\in C(H,V)$, then
\begin{align*}
\partial\delta\lambda\begin{pmatrix}
g_0 & g_1 \\
h_0 & h_1
\end{pmatrix} & =\delta\lambda(h_0,h_1)-\delta\lambda(h_0i(g_0),h_1i(g_1)) \\
              & =\rho_0^0(h_0)\lambda(h_1)-\lambda(h_0h_1)+\lambda(h_0)+ \\
              & \qquad\qquad -\big{(}\rho_0^0(h_0i(g_0))\lambda(h_1i(g_1))-\lambda(h_0i(g_0)h_1i(g_1))+\lambda(h_0i(g_0))\big{)},                                        
\end{align*}
and 
\begin{align*}
\delta\partial\lambda\begin{pmatrix}
g_0 & g_1 \\
h_0 & h_1
\end{pmatrix} & =\rho_0^0(h_0i(g_0))\partial\lambda(g_1,h_1)-\partial\lambda(g_0^{h_1}g_1,h_0h_1)+\partial\lambda(g_0,h_0) \\
              & =\rho_0^0(h_0i(g_0))(\lambda(h_1)-\lambda(h_1i(g_1)))+ \\
              & \qquad\qquad -\big{(}\lambda(h_0h_1)-\lambda(h_0h_1i(g_0^{h_1}g_1))\big{)}+\lambda(h_0)-\lambda(h_0i(g_0)). 
\end{align*}
Considering their difference, we get
\begin{align*}
(\delta\partial -\partial\delta)\lambda\begin{pmatrix}
g_0 & g_1 \\
h_0 & h_1
\end{pmatrix} & =\rho_0^0(h_0i(g_0))\lambda(h_1)-\rho_0^0(h_0)\lambda(h_1) \\
              & =\rho_0^0(h_0)\big{(}\rho_0^0(i(g_0))\lambda(h_1)-\lambda(h_1)\big{)} \\
              & =\rho_0^0(h_0)\big{(}\lambda(h_1)+\phi\rho_1(g_0)\lambda(h_1)-\lambda(h_1)\big{)}=\rho_0^0(h_0)\phi\rho_1(g_0)\lambda(h_1) .
\end{align*}
This hints the statement of a lemma analogous to lemma \ref{upToHomotopy0}, which can be schematized with the following diagram
\begin{eqnarray*}
\xymatrix{
C(\G_p^2,V) \ar[rr]^\partial                       &                                  & C(\G_{p+1}^2,V)    \\
                                                   & C(\G_p\times G,W) \ar[ur]^\Delta &                    \\
C(\G_p,V) \ar[uu]^\delta\ar[ur]^{\partial'}\ar[rr] &                                  & C(\G_{p+1},V) \ar[uu]_\delta ,  
}
\end{eqnarray*}
where $\partial'\lambda(\gamma;g):=\rho_0^1(t_p(\gamma))^{-1}\rho_1(g)\lambda(\gamma)$ and 
\begin{eqnarray*}
\Delta\alpha(\gamma_0,\gamma_1)=\Delta\alpha \begin{pmatrix}
    g_{00} & g_{01} \\
    \vdots & \vdots \\
    g_{p0} & g_{p1} \\
    h_0    & h_1 
\end{pmatrix}:=\rho_0^0(t_p(\partial_0\gamma_0)t_p(\partial_0\gamma_1))\phi\big{(}\alpha(\partial_0\gamma_1;g_{00})\big{)}.
\end{eqnarray*}
In fact, with pretty much the same amount of effort, we can state a proposition whose schematization is 
\begin{eqnarray*}
\xymatrix{
C(\G_p^{q+1},V) \ar[rr]^\partial                              &                                    & C(\G_{p+1}^{q+1},V)    \\
                                                              & C(\G_p^q\times G,W) \ar[ur]^\Delta &                        \\
C(\G_p^q,V) \ar[uu]^\delta\ar[ur]^{\partial'}\ar[rr]_\partial &                                    & C(\G_{p+1}^q,V) \ar[uu]_\delta ,  
}
\end{eqnarray*}
where $\partial'\omega(\gamma_1,...,\gamma_q;g):=\rho_0^1(t_p(\gamma_1)...t_p(\gamma_q))^{-1}\rho_1(g)\omega(\gamma_1,...,\gamma_q)$ and 
\begin{align*}
    \Delta\alpha(\vec{\gamma}) & :=\rho_0^0(t_p(\partial_0\gamma_0)...t_p(\partial_0\gamma_q))\phi\big{(}\alpha(\partial_0\delta_0\vec{\gamma};g_{00})\big{)}.
\end{align*}
Here,
\begin{align*}
\vec{\gamma}=(\gamma_0,...,\gamma_q) & =\begin{pmatrix}
    g_{00} & g_{01} & ... & g_{0q} \\
    g_{10} & g_{11} & ... & g_{1q} \\
   \vdots  & \vdots &     & \vdots \\
    g_{p0} & g_{p1} & ... & g_{pq} \\
    h_0    & h_1    & ... & h_q
\end{pmatrix}\in\G_{p+1}^{q+1}.
\end{align*}									
%\begin{pmatrix}
%    g_{11} & g_{12} & ... & g_{1q} \\
%   \vdots  & \vdots &     & \vdots \\
%    g_{p1} & g_{p2} & ... & g_{pq} \\
%    h_1    & h_2    & ... & h_q
%\end{pmatrix}
%*** Before diving into the proof, consider improving notation. Alternatives are:\\
%\begin{itemize}
%\item Using the group/groupoid differentials, so that 
%\begin{eqnarray*}
%\partial_0\gamma=\partial_0 \begin{pmatrix}
%    g_0 \\
%    g_1 \\
%    \vdots \\
%    g_p \\
%    h  
%\end{pmatrix}:=\begin{pmatrix}
%    g_1 \\
   % g_2 \\
%    \vdots \\
%    g_p \\
%    h  
%\end{pmatrix}.
%\end{eqnarray*}
%Using this notation, 
%\begin{eqnarray*}
%\Delta\alpha(\gamma_0,...,\gamma_q):=\rho_0^0(t_p(\partial_0\gamma_0)...t_p(\partial_0\gamma_q))\phi(\alpha(\partial_0\gamma_1,...,\partial_0\gamma_q;g_{00}))
%\end{eqnarray*}
%\item 
Understanding the vector $\vec{\gamma}=(\gamma_0,...,\gamma_q)$ as this latter matrix, %$\begin{pmatrix}
%    g_{00} & g_{01} & ... & g_{0q} \\
%    g_{10} & g_{11} & ... & g_{1q} \\
%   \vdots  & \vdots &     & \vdots \\
%    g_{p0} & g_{p1} & ... & g_{pq} \\
%    h_0    & h_1    & ... & h_q
%\end{pmatrix}$
we can use ``minor'' notation to abbreviate
\begin{eqnarray*}
\partial_i\delta_j\vec{\gamma}=\vec{\gamma}_{i,j}:=\begin{pmatrix}
    g_{00}   & ... & g_{0j-1}   & g_{0j+1}   & ... & g_{0q} \\
    \vdots   &     & \vdots     & \vdots     &     & \vdots \\
    g_{i-10} & ... & g_{i-1j-1} & g_{i-1j+1} & ... & g_{i-1q} \\
    g_{i+10} & ... & g_{i+1j-1} & g_{i+1j+1} & ... & g_{i+1q} \\
    \vdots   &     & \vdots     & \vdots     &     & \vdots \\
    g_{p0}   & ... & g_{pj-1}   & g_{pj+1}   & ... & g_{pq} \\
    h_0      & ... & h_{j-1}    & h_{j+1}    & ... & \vdots \\
\end{pmatrix}.
\end{eqnarray*}
%With this, the argument of $\alpha$ in the definition of $\Delta$ is %simply $(\vec{\gamma}_{0,0};g_{00})$. ***
%\end{itemize}
The following proposition is redundant with the previous lemmas of this subsection. It specifies up to which isomorphisms in the $2$-vector space $\xymatrix{W \ar[r]^\phi & V}$ the diagram \ref{page r=0} commutes.
\begin{prop}\label{upToHomotopy r=0}
Let $\omega\in C(\G_p^q,V)$ and $\vec{\gamma}=(\gamma_0,...,\gamma_q)\in\G_{p+1}^{q+1}$. Then
\begin{eqnarray*}
\delta\partial\omega(\vec{\gamma})=\partial\delta\omega(\vec{\gamma})+\phi\Big{(}\rho_0^1(t_p(\partial_0\gamma_0))\rho_1(g_{00})\omega(\vec{\gamma}_{0,0})\Big{)}
\end{eqnarray*}
\end{prop}
\begin{proof}
Essentially, the formula follows from the commutativity of $\delta_j$ and $\partial_k$ for all values of $(j,k)$ by taking care of the representations that appear at $(0,0)$.
\begin{align*}
\delta\partial\omega(\vec{\gamma}) & = \rho_0^0(t_{p+1}(\gamma_0))\delta_0^*\partial\omega(\vec{\gamma})+\sum_{j=1}^{q+1}(-1)^j\delta_j^*\partial\omega(\vec{\gamma}) \\
								   & = \rho_0^0(t_{p+1}(\gamma_0))\partial\omega(\delta_0\vec{\gamma})+\sum_{j=1}^{q+1}(-1)^j\partial\omega(\delta_j\vec{\gamma}) \\
								   & = \rho_0^0(t_{p+1}(\gamma_0))\sum_{k=0}^{p+1}(-1)^k\partial_k^*\omega(\delta_0\vec{\gamma})+\sum_{j=1}^{q+1}(-1)^j\sum_{k=0}^{p+1}(-1)^k\partial_k^*\omega(\delta_j\vec{\gamma}) \\
								   & = \sum_{k=0}^{p+1}(-1)^k\rho_0^0(t_{p+1}(\gamma_0))\omega(\partial_k\delta_0\vec{\gamma})+\sum_{j=1}^{q+1}\sum_{k=0}^{p+1}(-1)^{j+k}\omega(\partial_k\delta_j\vec{\gamma}). 
\end{align*}
On the other hand,
\begin{align*}
\partial\delta\omega(\vec{\gamma}) & = \sum_{k=0}^{p+1}(-1)^k\partial_k^*\delta\omega(\vec{\gamma}) \\
								   & = \sum_{k=0}^{p+1}(-1)^k\delta\omega(\partial_k\vec{\gamma}) \\
								   & = \sum_{k=0}^{p+1}(-1)^k\Big{(}\rho_0^0(t_p((\partial_k\vec{\gamma})_0))\delta_0^*\omega(\partial_k\vec{\gamma})+\sum_{j=1}^{q+1}(-1)^j\delta_j^*\omega(\partial_k\vec{\gamma})\Big{)} \\
								   & = \sum_{k=0}^{p+1}(-1)^k\rho_0^0(t_p((\partial_k\vec{\gamma})_0))\omega(\delta_0\partial_k\vec{\gamma})+\sum_{k=0}^{p+1}\sum_{j=1}^{q+1}(-1)^{k+j}\omega(\delta_j\partial_k\vec{\gamma}). 
\end{align*}
As stated, the double sums in the above expressions will coincide thanks to lemma \ref{non-zero jk}. Furthermore, for $0<k\leq p$,
\begin{align*}
(\partial_k\vec{\gamma})_0 & := \partial_k\gamma_0 \\
					       & = (\gamma_{00},...,\gamma_{k-2},\gamma_{k-10}\Join\gamma_{k0},\gamma_{k+10},...,\gamma_{p0}) \\
					       & \sim (g_{00},...,g_{k-20},g_{k0}g_{k-10},g_{k+10},...,g_{p0},h_0)
\end{align*}
since $\gamma_{k-10}\Join\gamma_{k0}=(g_{k-10},h_0i(g_{p0}...g_{k0}))\Join(g_{k0},h_0i(g_{p0}...g_{k+10}))=(g_{k0}g_{k-10},h_0i(g_{p0}...g_{k+10}))$; therefore,
\begin{align*}
t_p((\partial_k\vec{\gamma})_0) & = t_p(h_0i(g_{p0}...g_{k+10}(g_{k0}g_{k-10})g_{k-20}...g_{00}))=t_p(\gamma_0).
\end{align*}
Analogously, 
\begin{align*}
(\partial_{p+1}\vec{\gamma})_0 & := \partial_{p+1}\gamma_0 \\
					           & = (\gamma_{00},...,\gamma_{p-10}) \\
					           & \sim (g_{00},...,g_{p-10},h_0i(p_{p0}))
\end{align*}
since $\gamma_{p-10}=(g_{p-10},h_0i(g_{p0}))$; thus implying,
\begin{align*}
t_p((\partial_{p+1}\vec{\gamma})_0) & = t_p((h_0i(g_{p0}))i(g_{p-10}...g_{00}))=t_p(\gamma_0).
\end{align*}
Hence, using lemma \ref{(j,k)=(0,0)} and $\rho_0^0(i(g))=I+\phi\rho_1(g)$,
\begin{align*}
(\delta\partial-\partial\delta)\omega(\vec{\gamma}) & =\rho_0^0(t_p(\gamma_0))\omega(\vec{\gamma}_{0,0})-\rho_0^0(t_p(\partial_0\gamma_0))\omega(\vec{\gamma}_{0,0}) \\
													& =\rho_0^0(h_0i(g_{p0}...g_{10}))\Big{[}\rho_0^0(i(g_{00}))-I\Big{]}\omega(\vec{\gamma}_{0,0}) \\
													& =\rho_0^0(h_0i(g_{p0}...g_{10}))\phi\rho_1(g_{00})\omega(\vec{\gamma}_{0,0})=\phi\Big{(}\rho_0^1(t_p(\partial_0\gamma_0))\rho_1(g_{00})\omega(\vec{\gamma}_{0,0})\Big{)}.
\end{align*}
\end{proof}
We close this discussion remarking that the equation in the statement of the latter lemma reads
\begin{eqnarray*}
\delta\partial-\partial\delta=\Delta\circ\partial'.
\end{eqnarray*}
\subsubsection{Pieces of the triple complex: the natural emergence}
In sight of proposition \ref{upToHomotopy r=0}, we have got a lattice of spaces
\begin{eqnarray*}
\xymatrix{
\vdots           & \vdots            & \vdots              &       \\ 
C(H^3\times G,W) & C(\G^3\times G,W) & C(\G_2^3\times G,W) & \dots \\
C(H^2\times G,W) & C(\G^2\times G,W) & C(\G_2^2\times G,W) & \dots \\
C(H\times G,W)   & C(\G\times G ,W)  & C(\G_2\times G,W)   & \dots \\
C(G,W)           & C(G,W)            & C(G,W)	           & \dots  
}
\end{eqnarray*} 
hanging around. In fact, in the specific sense of the aforementioned proposition, diagram \ref{page r=0} commutes up to this lattice. Nevertheless, it should come as no surprise that this object can be filled with maps turning it into a sort of double complex itself: 
\begin{eqnarray}\label{page r=1}
\xymatrix{
\vdots                                             & \vdots                                              & \vdots                           &       \\ 
C(H^3\times G,W) \ar[r]^{\partial}\ar[u]           & C(\G^3\times G,W) \ar[r]^{\partial}\ar[u]           & C(\G_2^3\times G,W) \ar[r]\ar[u]           & \dots \\
C(H^2\times G,W) \ar[r]^{\partial}\ar[u]^{\delta'} & C(\G^2\times G,W) \ar[r]^{\partial}\ar[u]^{\delta'} & C(\G_2^2\times G,W) \ar[r]\ar[u]^{\delta'} & \dots \\
C(H\times G,W) \ar[r]^{\partial}\ar[u]^{\delta'}   & C(\G\times G ,W) \ar[r]^{\partial}\ar[u]^{\delta'}  & C(\G_2\times G,W) \ar[r]\ar[u]^{\delta'} & \dots \\
C(G,W) \ar[r]^{0}\ar[u]^{\delta'}                  & C(G,W) \ar[r]^{Id}\ar[u]^{\delta'}                  & C(G,W)	\ar[r]^{0}\ar[u]^{\delta'} & \dots  
}
\end{eqnarray} 
Here, the vertical $\delta'$ maps are differentials for the groupoid cochain complex of the (right) transformation groupoid 
\begin{eqnarray*}
\xymatrix{\G_p\ltimes G \ar@<0.5ex>[r]\ar@<-0.5ex>[r] & G},
\end{eqnarray*}
with respect to the right representation
\begin{eqnarray*}
(g;w)\cdot(\gamma;g):=(g^{t_p(\gamma)};\rho_0^1(t_p(\gamma))^{-1}w)
\end{eqnarray*}
on the trivial bundle $\xymatrix{G\times W \ar[r]^{\quad pr_1} & G}$; whereas the horizontal maps are differentials for the groupoid cochain complex of the product groupoid 
\begin{eqnarray*}
\xymatrix{\G^q\ltimes G \ar@<0.5ex>[r]\ar@<-0.5ex>[r] & H^q\times G},
\end{eqnarray*}
with respect to the representation
\begin{eqnarray}\label{r=1 p-rep}
(\gamma_1,...,\gamma_q;g)\cdot(h_1,...,h_q;g,w):=(h_1i(g_1),...,h_qi(g_q);g,\rho_0^1(i(pr_G(\gamma_1\vJoin ...\vJoin\gamma_q)))^{-1}w), 
\end{eqnarray}
where $\gamma_k=(g_k,h_k)$, on the trivial bundle $\xymatrix{H^q\times G\times W \ar[r]^{\quad pr_1} & H^q\times G}$. Notice that for $q=1$, this representation coincides with the pull-back of $\Delta^W$ of the representation up to homotopy induced by the $2$-representation along the projection onto $\xymatrix{\G \ar@<0.5ex>[r]\ar@<-0.5ex>[r] & H}$. \\
We wrote that diagram \ref{page r=1} is a sort of double complex, because in fact it does not commute either. By studying how it fails to commute, one is bound to discover that there is yet a third ``double complex'' which controls this failure in a similar fashion to that of the discussion we just laid down. Inductively, this third lattice will fail to commute prompting a sequence of two-dimensional lattices that will commute only up to homotopy, and that, put together, form a three-dimensional lattice. By looking at the sum of all the differentials and carefully introducing difference maps, we will turn this into an appropriate complex whose second cohomology controls extensions of the Lie $2$-group by $2$-vector spaces. \\
Before moving on to formalize this latter paragraph, let us first set the horizontal representations in stone
\begin{lemma}\label{r=1 p-cx}
Equation \ref{r=1 p-rep} defines a representation.
\end{lemma}
In order to make cleaner computations, we introduce the following auxiliary straightforward lemma.
\begin{lemma}\label{multiprods}
Let $\gamma_1,...,\gamma_q\in\G$. If $\gamma_k=(g_k,h_k)$, then
\begin{eqnarray*}
\gamma_1\vJoin ...\vJoin\gamma_q =(g_1^{h_2...h_q}g_2^{h_3...h_q}...g_{q-2}^{h_{q-1}h_q}g_{q-1}^{h_q}g_q,h_1...h_q).
\end{eqnarray*}
\end{lemma}
\begin{proof}
By induction on $q$, for $q=2$ the formula is nothing but the definition of the product in $\G$. Now, suppose the equation holds for $q-1$ elements, then
\begin{align*}
\gamma_1\vJoin ...\vJoin\gamma_q & =(\gamma_1\vJoin ...\vJoin\gamma_{q-1})\vJoin\gamma_q \\
								 & =(g_1^{h_2...h_{q-1}}g_2^{h_3...h_{q-1}}...g_{q-3}^{h_{q-2}h_{q-1}}g_{q-2}^{h_{q-1}}g_{q-1},h_1...h_{q-1})\vJoin (g_q,h_q) \\
								 & =((g_1^{h_2...h_{q-1}}g_2^{h_3...h_{q-1}}...g_{q-3}^{h_{q-2}h_{q-1}}g_{q-2}^{h_{q-1}}g_{q-1})^{h_q}g_q,h_1...h_{q-1}h_q),
\end{align*} 
and the result follows since the action of $H$ is by automorphisms.

\end{proof}
\begin{proof}(\textit{of lemma} \ref{r=1 p-cx})
We defined the first coordinate of the action to coincide with the action, and hence, to respect the moment map. Because of this, we can focus our attention on the $W$ coordinate. \\ 
Thus defined, the action respects units: Let $\gamma_1,...,\gamma_q\in\G$ and $w\in W$. If $g_k=1$ for all $k$, the formula reads
\begin{align*}
\rho_0^1(i(pr_G(\gamma_1\vJoin ...\vJoin\gamma_q)))^{-1}w & =\rho_0^1(i(1^{h_2...h_q}1^{h_3...h_q}...1^{h_{q-1}h_q}1^{h_q}1))^{-1}w \\
														  & =\rho_0^1(i(1...1))^{-1}w =w.
\end{align*}
As for the groupoid multiplication, let $\gamma'_1,...,\gamma'_q\in\G$ be such that $(\gamma'_k,\gamma_k)\in\G^{(2)}$ for every $k$, then $\gamma'_k=(g'_k,h_ki(g_k))$ and 
\begin{eqnarray*}
(\gamma'_1,...,\gamma'_q;g)(\gamma_1,...,\gamma_q;g):=\Big{(}\begin{pmatrix}
g_1g'_1 & ... & g_qg'_q \\
h_1     & ... & h_q
\end{pmatrix};g\Big{)}
\end{eqnarray*}
as $\gamma'_k\Join\gamma_k=(g_kg'_k,h_k)$. The compatibility will follow then from the formula
\begin{align*}
(g_1g'_1)^{h_2...h_q}(g_2g'_2)^{h_3...h_q}...(g_{q-1} & g'_{q-1})^{h_q}g_qg'_q = \\
                                                      & g_1^{h_2...h_q}g_2^{h_3...h_q}...g_{q-1}^{h_q}g_q(g'_1)^{h_2i(g_2)...h_qi(g_q)}(g'_2)^{h_3i(g_3)...h_qi(g_q)}...(g'_{q-1})^{h_qi(g_q)}g'_q .
\end{align*} 
We proceed by induction on $q$. For $q=2$,
\begin{align*}
(g_1g'_1)^{h_2}g_2g'_2 & = g_1^{h_2}(g'_1)^{h_2}g_2g'_2 \\
                       & = g_1^{h_2}g_2g_2^{-1}(g'_1)^{h_2}g_2g'_2 \\
                       & = g_1^{h_2}g_2(g'_1)^{h_2i(g_2)}g'_2 .
\end{align*}
Suppose now that the equation holds for $q-1$, then
\begin{align*}
\prod_{k=1}^q (g_kg'_k)^{\prod_{j=k+1}^q h_j} & =\Bigg{(}\prod_{k=1}^{q-1} (g_kg'_k)^{\Big{(}\prod_{j=k+1}^{q-1} h_j\Big{)}h_q}\Bigg{)}g_qg'_q \\
											  & =\Bigg{(}\prod_{k=1}^{q-1} (g_kg'_k)^{\prod_{j=k+1}^{q-1} h_j}\Bigg{)}^{h_q}g_qg'_q \\
					(I.H.)   				  & =\Bigg{(}\prod_{k=1}^{q-1} g_k^{\prod_{j=k+1}^{q-1} h_j}\prod_{k=1}^{q-1}(g'_k)^{\prod_{j=k+1}^{q-1} h_ji(g_j)}\Bigg{)}^{h_q}g_qg'_q \\
											  & =\Bigg{(}\prod_{k=1}^{q-1} g_k^{\prod_{j=k+1}^{q-1} h_j}\Bigg{)}^{h_q}g_qg_q^{-1}\Bigg{(}\prod_{k=1}^{q-1}(g'_k)^{\prod_{j=k+1}^{q-1} h_ji(g_j)}\Bigg{)}^{h_q}g_qg'_q \\
											  & =\prod_{k=1}^{q} g_k^{\prod_{j=k+1}^{q} h_j}\Bigg{(}\prod_{k=1}^{q-1}(g'_k)^{\prod_{j=k+1}^{q-1} h_ji(g_j)}\Bigg{)}^{h_qi(g_q)}g'_q , 
\end{align*}
which is precisely what we wanted.

\end{proof}
In order to get acquainted with the type of relations that will control the non-commutativity of diagram \ref{page r=1}, let us consider $(\partial\delta'-\delta'\partial)\lambda$ for $\lambda\in C(G,W)$. Let $\gamma=(g,h)\in\G$ and $f\in G$, then
\begin{align*}
(\partial\delta'-\delta'\cancel{\partial})\lambda(\gamma;f) & =\rho_0^1(i(g))^{-1}\delta'\lambda(h;f)-\delta'\lambda(hi(g);f) \\
                                                            & =\rho_0^1(i(g))^{-1}\big{(}\lambda(f^h)-\rho_0^1(h)^{-1}\lambda(f)\big{)}-\big{(}\lambda(f^{hi(g)})-\rho_0^1(hi(g))^{-1}\lambda(f)\big{)} \\
                                                            & =\rho_0^1(i(g))^{-1}\lambda(f^h)-\lambda(f^{hi(g)}).
\end{align*}
Clearly, this expression need not be zero. There is a third way to go from $C(G,W)$ to $C(\G\times G,W)$, though. Indeed, since the differentials in the first page did not commute, we will need to add the maps we laid down in the previous section: the difference maps $\Delta$ and the differentials labeled $\partial'$. Now, $\partial'\Delta\lambda\in C(\G\times G,W)$ and
\begin{align*}
\partial'\Delta\lambda(\gamma;f) & =\rho_0^1(hi(g))^{-1}\rho_1(f)\Delta\lambda(\gamma) \\
                                 & =\rho_0^1(hi(g))^{-1}\rho_1(f)\rho_0^0(h)\phi\lambda(g) \\
                                 & =\rho_0^1(i(g))^{-1}\rho_1(f^h)\phi\lambda(g) \\
                                 & =\rho_0^1(i(g))^{-1}[\rho_0^1(i(f^h))-I]\lambda(g). 
\end{align*}
Unfortunately, summing this expression does not help cancelling the difference $(\partial\delta'-\delta'\partial)\lambda$; however, when added, they will yield a familiar expression.
\begin{align*}
(\partial\delta'-\delta' & \partial+\partial'\Delta)\lambda(\gamma;f) \\
          & =\rho_0^1(i(g))^{-1}\lambda(f^h)-\lambda(f^{hi(g)})+\rho_0^1(i(g))^{-1}\rho_0^1(i(f^h))\lambda(g)-\rho_0^1(i(g))^{-1}\lambda(g) \\
          & =\rho_0^1(i(g))^{-1}\big{(}\rho_0^1(i(f^h))\lambda(g)-\lambda(f^hg)+\lambda(f^h)\big{)}+ \\
          & \qquad +\big{(}\rho_0^1(i(g))^{-1}\lambda(f^hg)-\lambda(f^{hi(g)})+\lambda(g^{-1})\big{)}-\big{(}\rho_0^1(i(g))^{-1}\lambda(g)+\lambda(g^{-1})\big{)} \\
          & =\rho_0^1(i(g))^{-1}\delta_{(1)}\lambda(f^h,g)+\delta_{(1)}\lambda(g^{-1},f^{h}g)-\delta_{(1)}\lambda(g^{-1},g),
\end{align*}
where $\delta_{(1)}$ is the differential of group cochains for $G$ with values in $\rho_0^1\circ i$. Thus, defining
\begin{eqnarray*}
\xymatrix{
\Delta :C(G^2,W) \ar[r] & C(\G\times G,W)
}\qquad\qquad\qquad \\
\Delta\omega(\gamma;f):=\rho_0^1(i(g))^{-1}\omega(f^h,g)+\omega(g^{-1},f^{h}g)-\omega(g^{-1},g),
\end{eqnarray*}
we get the relation
\begin{eqnarray*}
\partial\circ\delta'-\delta'\circ\partial=\Delta\circ\delta_{(1)}-\partial'\circ\Delta.
\end{eqnarray*}
Analogously, for $\alpha\in C(H\times G,W)$, for $j\in\lbrace 0,1\rbrace$ let $\gamma_j=(g_j,h_j)\in\G$ and $f\in G$, then
\begin{align*}
\partial\delta'\alpha(\gamma_0,\gamma_1;f) & =\rho_0^1(i(g_0^{h_1}g_1))^{-1}\delta'\alpha(h_0,h_1;f)-\delta'\alpha(h_0i(g_0),h_1i(g_1);f) \\
                                           & =\rho_0^1(i(g_0^{h_1}g_1))^{-1}\big{(}\alpha(h_1;f^{h_0})-\alpha(h_0h_1;f)+\rho_0^1(h_1)^{-1}\alpha(h_0;f)\big{)}+ \\
                                           & \qquad -\big{(}\alpha(h_1i(g_1);f^{h_0i(g_0)})-\alpha(h_0i(g_0)h_1i(g_1);f)+\rho_0^1(h_1i(g_1))^{-1}\alpha(h_0i(g_0);f)\big{)}, \\
\end{align*}
and also
\begin{align*}
\delta'\partial\alpha(\gamma_0,\gamma_1;f) & =\partial\alpha(\gamma_1;f^{h_0i(g_0)})-\partial\alpha(\gamma_0\vJoin\gamma_1;f)+\rho_0^1(h_1i(g_1))^{-1}\partial\alpha(\gamma_0;f) \\
                                           & =\rho_0^1(i(g_1))^{-1}\alpha(h_1;f^{h_0i(g_0)})-\alpha(h_1i(g_1);f^{h_0i(g_0)})+ \\
                                           & \qquad -\big{(}\rho_0^1(i(g_0^{h_1}g_1))^{-1}\alpha(h_0h_1;f)-\alpha(h_0h_1i(g_0g_1);f)\big{)}+ \\
                                           & \qquad\quad +\rho_0^1(h_1i(g_1))^{-1}\big{(}\rho_0^1(i(g_0))^{-1}\alpha(h_0;f)-\alpha(h_0i(g_0);f)\big{)}; 
\end{align*}
thus taking their difference yields
\begin{eqnarray*}
(\partial\delta'-\delta'\partial)\alpha(\gamma_0,\gamma_1;f)=\rho_0^1(i(g_1))^{-1}\big{[}\rho_0^1(i(g_0^{h_1}))^{-1}\alpha(h_1;f^{h_0})-\alpha(h_1;f^{h_0i(g_0)})\big{]}.
\end{eqnarray*}
On the other hand,
\begin{align*}
\partial'\Delta\alpha(\gamma_0,\gamma_1;f) & =\rho_0^1(h_0i(g_0)h_1i(g_1))^{-1}\rho_1(f)\Delta\alpha(\gamma_0,\alpha_1) \\
                                           & =\rho_0^1(h_0i(g_0)h_1i(g_1))^{-1}\rho_1(f)\rho_0^0(h_0h_1)\phi\alpha(h_1;g_0) \\
                                           & =\rho_0^1(i(g_0^{h_1}g_1))^{-1}\rho_1(f^{h_0h_1})\phi\alpha(h_1;g_0) \\
                                           & =\rho_0^1(i(g_0^{h_1}g_1))^{-1}[\rho_0^1(i(f^{h_0h_1}))-I]\alpha(h_1;g_0). 
\end{align*}
Adding these expressions will make another difference map emerge:
\begin{align*}
(\partial & \delta'-\delta'\partial+\partial'\Delta)\alpha(\gamma_0,\gamma_1;f) \\
          & =\rho_0^1(i(g_1))^{-1}\big{[}\rho_0^1(i(g_0^{h_1}))^{-1}\alpha(h_1;f^{h_0})-\alpha(h_1;f^{h_0i(g_0)})\big{]}+\rho_0^1(i(g_0^{h_1}g_1))^{-1}[\rho_0^1(i(f^{h_0h_1}))-I]\alpha(h_1;g_0) \\
          & =\rho_0^1(i(g_1))^{-1}\Big{[}\rho_0^1(i(g_0^{h_1}))^{-1}\big{(}\rho_0^1(i(f^{h_0h_1}))\alpha(h_1;g_0)-\alpha(h_1;f^{h_0}g_0)+\alpha(h_1;f^{h_0})\big{)}+ \\
          & \qquad +\rho_0^1(i(g_0^{h_1}))^{-1}\alpha(h_1;f^{h_0}g_0)-\alpha(h_1;f^{h_0i(g_0)})+\alpha(h_1;g_0^{-1})+\rho_0^1(i(g_0^{h_1}))^{-1}\alpha(h_1;g_0)-\alpha(h_1;g_0^{-1})\Big{]} \\
          & =\rho_0^1(i(g_1))^{-1}\Big{[}\rho_0^1(i(g_0^{h_1}))^{-1}\delta_{(1)}\alpha(h_1;f^{h_0},g_0)+\delta_{(1)}\alpha(h_1;g_0^{-1},f^{h_0}g_0)-\delta_{(1)}\alpha(h_1;g_0^{-1},g_0)\Big{]},
\end{align*}
where this time around $\delta_{(1)}$ is the differential of groupoid cochains for the bundle of Lie groups $\xymatrix{H\times G \ar@<-0.5ex>[r]\ar@<0.5ex>[r] & H}$ with values on the trivial vector bundle $\xymatrix{H\times W \ar[r] & H}$ endowed with the left action
\begin{eqnarray*}
(h;g)\cdot(h;w):=(h;\rho_0^1(i(g^h))w).
\end{eqnarray*} 
Thus, if we define
\begin{eqnarray*}
\xymatrix{
\Delta :C(H\times G^2,W) \ar[r] & C(\G^2\times G,W)
}\qquad\qquad\qquad\qquad\qquad \\
\Delta\omega(\gamma_0,\gamma_1;f):=\rho_0^1(i(g_1))^{-1}\Big{[}\rho_0^1(i(g_0^{h_1}))^{-1}\omega(h_1;f^{h_0},g_0)+\omega(h_1;g_0^{-1},f^{h_0}g_0)-\omega(h_1;g_0^{-1},g_0)\Big{]},
\end{eqnarray*}
we get the relation
\begin{eqnarray*}
\partial\circ\delta'-\delta'\circ\partial=\Delta\circ\delta_{(1)}-\partial'\circ\Delta
\end{eqnarray*}
again. Vastly generalizing, we get the upcoming proposition, which formalizes the following schematization
\begin{eqnarray*}
\xymatrix{
\circ \ar@{.}[rr]\ar@{.}[dr] & & C(\G_p^{q+1}\times G,W) \ar@{.}[rr]\ar[dr]^\partial & & \circ\ar@{.}[dr] & \\
 & C(\G_{p+1}^{q+1},V) \ar[rr]^{\partial'\qquad\qquad} & & C(\G_{p+1}^{q+1}\times G,W) \ar@{.}[rr] & & \circ \\
\circ\ar@{.}[rr]\ar@{.}[dr]\ar@{.}[uu] && C(\G_p^q\times G,W)\ar'[r]^{\qquad\delta_{(1)}}[rr]\ar[dr]_\partial\ar[ul]^\Delta\ar'[u]_{\delta'}[uu] && C(\G_p^q\times G^2,W)\ar@{.}[dr]\ar@{.}[uu]\ar[ul]_\Delta & \\
 & \circ \ar@{.}[rr]\ar@{.}[uu] & & C(\G_{p+1}^q\times G,W) \ar@{.}[rr]\ar[uu]_(.35){\delta'} & & \circ \ar@{.}[uu]
}
\end{eqnarray*}
\begin{prop}
Let 
\begin{eqnarray*}
\xymatrix{
\Delta :C(\G_p^q\times G^2,W) \ar[r] & C(\G_{p+1}^{q+1}\times G,W)
}
\end{eqnarray*}
be defined by the formula
\begin{align*}
\Delta\omega(\vec{\gamma};f):=\rho_0^1(i(pr_G(\gamma_{01}\vJoin ...\vJoin\gamma_{0q})))^{-1}\Big{[}\rho_0^1( & i(g_{00}^{h_{01}...h_{0q}}))^{-1}\omega(\vec{\gamma}_{0,0};f^{h_{00}},g_{00})+ \\
                                                                    & +\omega(\vec{\gamma}_{0,0};g_{00}^{-1},f^{h_{00}}g_{00})-\omega(\vec{\gamma}_{0,0};g_{00}^{-1},g_{00})\Big{]},    
\end{align*}
for $f\in G$ and $\vec{\gamma}=(\gamma_0,...,\gamma_q)\in\G_{p+1}^{q+1}$, where 
\begin{eqnarray*}
\gamma_b=\begin{pmatrix}
\gamma_{0b} \\
\vdots    \\
\gamma_{pb}
\end{pmatrix}=\begin{pmatrix}
g_{0b} & h_{0b} \\
\vdots & \vdots   \\
g_{pb} & h_{pb}
\end{pmatrix}\sim\begin{pmatrix}
g_{0b} \\
\vdots \\
g_{pb} \\
h_b
\end{pmatrix}.
\end{eqnarray*}
Then 
\begin{eqnarray*}
\partial\circ\delta'-\delta'\circ\partial=\Delta\circ\delta_{(1)}-\partial'\circ\Delta ,
\end{eqnarray*}
where $\delta_{(1)}$ is the differential of groupoid cochains for the bundle of Lie groups $\xymatrix{\G_{p}^{q}\times G \ar@<-0.5ex>[r]\ar@<0.5ex>[r] & \G_{p}^{q}}$ with values on the trivial vector bundle $\xymatrix{\G_{p}^{q}\times W \ar[r] & \G_{p}^{q}}$ endowed with the left action
\begin{eqnarray*}
(\gamma_1,...,\gamma_q;f)\cdot(\gamma_1,...,\gamma_q;w):=(\gamma_1,...,\gamma_q;\rho_0^1(i(f^{t_p(\gamma_1)...t_p(\gamma_q)}))w).
\end{eqnarray*} 
\end{prop}
\begin{proof}
First of all, let us point out that under the identification $\G_{p+1}\cong G^{p+1}\times H$ referred in the statement, $\gamma_{ab}=(g_{ab},h_{ab})=\Big{(}g_{ab},h_bi\big{(}\prod_{k=0}^{p-a-1}g_{(p-k)b}\big{)}\Big{)}$. Now, let $\omega\in C(\G_p^q\times G,W)$ and $(\vec{\gamma};f)$ as in the statement. Let us introduce the shorthand
\begin{eqnarray*}
\rho^{q+1}(\vec{\gamma}_0):=\rho_0^1(i(pr_G(\gamma_{00}\vJoin...\vJoin\gamma_{0q})))^{-1}
\end{eqnarray*}
for the representation of the $q$th row. We compute the left hand side of the equation of the statement. On the one hand
\begin{align*}
\partial\delta'\omega(\vec{\gamma};f) & =\rho^{q+1}(\vec{\gamma}_0)\delta'\omega(\partial_0\vec{\gamma};f)+\sum_{j=1}^{p+1}(-1)^j\delta'\omega(\partial_j\vec{\gamma};f) \\
                                      & =\rho^{q+1}(\vec{\gamma}_0)\Big{(}\omega(\delta_0\partial_0\vec{\gamma};f^{t_p(\partial_0\gamma_0)})+\sum_{k=1}^q(-1)^k\omega(\delta_k\partial_0\vec{\gamma};f)+(-1)^{q+1}\rho_0^1(t_p(\partial_0\gamma_q))^{-1}\omega(\delta_{q+1}\partial_0\vec{\gamma};f)\Big{)}+ \\
                                      & \quad +\sum_{j=1}^{p+1}(-1)^j\Big{(}\omega(\delta_0\partial_j\vec{\gamma};f^{t_p(\partial_j\gamma_0)})+\sum_{k=1}^q(-1)^k\omega(\delta_k\partial_j\vec{\gamma};f)+(-1)^{q+1}\rho_0^1(t_p(\partial_j\gamma_q))^{-1}\omega(\delta_{q+1}\partial_j\vec{\gamma};f)\Big{)}, \\
\end{align*}
whereas on the other
\begin{align*}
\delta'\partial\omega(\vec{\gamma};f) & =\partial\omega(\delta_0\vec{\gamma};f^{t_{p+1}(\gamma_0)})+\sum_{k=1}^q(-1)^k\partial\omega(\delta_k\vec{\gamma};f)+(-1)^{q+1}\rho_0^1(t_{p+1}(\gamma_q))^{-1}\omega(\delta_{q+1}\vec{\gamma};f) \\
                                      & =\rho^{q}((\delta_0\vec{\gamma})_0)\omega(\partial_0\delta_0\vec{\gamma};f^{t_{p+1}(\gamma_0)})+\sum_{j=1}^{p+1}(-1)^j\omega(\partial_j\delta_0\vec{\gamma};f^{t_{p+1}(\gamma_0)})+ \\
                                      & \qquad +\sum_{k=1}^q(-1)^k\Big{(}\rho^{q}((\delta_k\vec{\gamma})_0)\omega(\partial_0\delta_k\vec{\gamma};f)+\sum_{j=1}^{p+1}(-1)^j\omega(\partial_j\delta_k\vec{\gamma};f)\Big{)}+ \\
                                      & \qquad\quad +(-1)^{q+1}\rho_0^1(t_{p+1}(\gamma_q))^{-1}\Big{(}\rho^{q}((\delta_{q+1}\vec{\gamma})_0)\omega(\partial_0\delta_{q+1}\vec{\gamma};f)+\sum_{j=1}^{p+1}(-1)^j\omega(\partial_j\delta_{q+1}\vec{\gamma};f)\Big{)}. \\
\end{align*}
We claim that their difference is
\begin{eqnarray*}
 (\partial\delta'-\delta'\partial)\omega(\vec{\gamma};f)=\rho^{q+1}(\vec{\gamma}_0)\omega(\delta_0\partial_0\vec{\gamma};f^{t_p(\partial_0\gamma_0)})-\rho^{q}((\delta_0\vec{\gamma})_0)\omega(\partial_0\delta_0\vec{\gamma};f^{t_{p+1}(\gamma_0)}).
\end{eqnarray*}
This will follow first and foremost from the commutativity of all simplicial maps $\delta_k\partial_j=\partial_j\delta_k$, and from the following identities:
\begin{itemize}
    \item $\rho^{q+1}(\vec{\gamma}_0)=\rho^{q}((\delta_k\vec{\gamma})_0)$ for $1\leq k\leq q$. Indeed, for the ranging values of $k$, 
    \begin{eqnarray*}
    (\delta_k\vec{\gamma})_0=(\gamma_{00},...,\gamma_{0(k-1)}\vJoin\gamma_{0k},...,\gamma_{0q});
    \end{eqnarray*} 
    therefore, 
    \begin{eqnarray*}
    \rho^{q}((\delta_k\vec{\gamma})_0)=\rho_0^1(i(pr_G(\gamma_{00}\vJoin...\vJoin(\gamma_{0(k-1)}\vJoin\gamma_{0k})\vJoin...\vJoin\gamma_{0q})))^{-1}=\rho^{q+1}(\vec{\gamma}_0).
    \end{eqnarray*}
    \item $\rho^{q+1}(\vec{\gamma}_0)\rho_0^1(t_p(\partial_0\gamma_q))^{-1}=\rho_0^1(t_{p+1}(\gamma_q))^{-1}\rho^{q}((\delta_{q+1}\vec{\gamma})_0)$. Indeed, from lemma \ref{multiprods}, we know that 
    \begin{eqnarray*}
    \rho^{q+1}(\vec{\gamma}_0)=\rho_0^1(i(g_{00}^{h_{01}...h_{0q}}g_{01}^{h_{02}...h_{0q}}...g_{0q-2}^{h_{0q-1}h_{0q}}g_{0q-1}^{h_{0q}}g_{0q}))^{-1}.
    \end{eqnarray*}
    Now, by the very definition $t_p(\partial_0\gamma_q)=t(\gamma_{1q})=s(\gamma_{0q})=h_{0q}$; thus,
    \begin{align*}
      \rho^{q+1}(\vec{\gamma}_0)\rho_0^1(t_p(\partial_0\gamma_q))^{-1} & =\rho_0^1(i(g_{00}^{h_{01}...h_{0q}}g_{01}^{h_{02}...h_{0q}}...g_{0q-2}^{h_{0q-1}h_{0q}}g_{0q-1}^{h_{0q}}g_{0q}))^{-1} \rho_0^1(h_{0q})^{-1} \\
                                                                       & =\rho_0^1(h_{0q}i\big{(}(g_{00}^{h_{01}...h_{0q-1}}g_{01}^{h_{02}...h_{0q-1}}...g_{0q-2}^{h_{0q-1}}g_{0q-1})^{h_{0q}}\big{)}i(g_{0q}))^{-1} \\
                                                                       & =\rho_0^1(i\big{(}g_{00}^{h_{01}...h_{0q-1}}g_{01}^{h_{02}...h_{0q-1}}...g_{0q-2}^{h_{0q-1}}g_{0q-1}\big{)}h_{0q}i(g_{0q}))^{-1} \\
                                                                       & =\rho_0^1(h_{0q}i(g_{0q}))^{-1}\rho_0^1(i(g_{00}^{h_{01}...h_{0q-1}}g_{01}^{h_{02}...h_{0q-1}}...g_{0q-2}^{h_{0q-1}}g_{0q-1}))^{-1} \\
                                                                       & =\rho_0^1(t(\gamma_{0q}))^{-1}\rho^q((\delta_{q+1}\vec{\gamma})_0)=\rho_0^1(t_{p+1}(\gamma_q))^{-1}\rho^q((\delta_{q+1}\vec{\gamma})_0).
    \end{align*}
    \item As we have already seen, $t_p(\partial_j\gamma_b)=t_{p+1}(\gamma_b)$ for $1\leq j\leq q+1$.
\end{itemize}
Using lemma \ref{multiprods} again, it is easy to see that 
\begin{align*}
    \rho^{q+1}(\vec{\gamma}_0) & =\rho_0^1(i(g_{00}^{h_{01}...h_{0q}}g_{01}^{h_{02}...h_{0q}}...g_{0q-2}^{h_{0q-1}h_{0q}}g_{0q-1}^{h_{0q}}g_{0q}))^{-1} \\
                               & =\rho_0^1(i(g_{01}^{h_{02}...h_{0q}}...g_{0q-2}^{h_{0q-1}h_{0q}}g_{0q-1}^{h_{0q}}g_{0q}))^{-1}\rho_0^1(i(g_{00}^{h_{01}...h_{0q}}))^{-1} =\rho^q((\delta_0\vec{\gamma})_0)\rho_0^1(i(g_{00}^{h_{01}...h_{0q}}))^{-1};
\end{align*}
hence, we can rewrite the difference as
\begin{eqnarray*}
 (\partial\delta'-\delta'\partial)\omega(\vec{\gamma};f)=\rho^q((\delta_0\vec{\gamma})_0)\big{[}\rho_0^1(i(g_{00}^{h_{01}...h_{0q}}))^{-1}\omega(\vec{\gamma}_{0,0};f^{h_{00}})-\omega(\vec{\gamma}_{0,0};f^{h_{00}i(g_{00})})\big{]}.
\end{eqnarray*}
We proceed to compute the second term on the right hand side of the equation in the statement,
\begin{align*}
\partial'\Delta\omega(\vec{\gamma};f) & =\rho_0^1(t_{p+1}(\gamma_0)...t_{p+1}(\gamma_q))^{-1}\rho_1(f)\Delta\omega(\vec{\gamma}) \\
                                      & =\rho_0^1(t_{p+1}(\gamma_0)...t_{p+1}(\gamma_q))^{-1}\rho_1(f)\rho_0^0(t_p(\partial_0\gamma_0)...t_p(\partial_0\gamma_q))\phi(\omega(\vec{\gamma}_{0,0};g_{00})) \\
                                      & =\rho_0^1(t(\gamma_{00})...t(\gamma_{0q}))^{-1}\rho_1(f)\rho_0^0(t(\gamma_{10})...t(\gamma_{1q}))\phi(\omega(\vec{\gamma}_{0,0};g_{00})) \\
                                      & =\rho_0^1(t(\gamma_{00}\vJoin...\vJoin\gamma_{0q}))^{-1}\rho_1(f)\rho_0^0(s(\gamma_{00})...s(\gamma_{0q}))\phi(\omega(\vec{\gamma}_{0,0};g_{00})) \\
                                      & =\rho_0^1(h_{00}...h_{0q}i(pr_G(\gamma_{00}\vJoin...\vJoin\gamma_{0q})))^{-1}\rho_1(f)\rho_0^0(h_{00}...h_{0q})\phi(\omega(\vec{\gamma}_{0,0};g_{00})) \\
                                      & =\rho_0^1(i(pr_G(\gamma_{00}\vJoin...\vJoin\gamma_{0q})))^{-1}\rho_1(f^{h_{00}...h_{0q}})\phi(\omega(\vec{\gamma}_{0,0};g_{00})) \\
                                      & =\rho^{q+1}(\vec{\gamma}_0)[\rho_0^1(i(f^{h_{00}...h_{0q}}))-I]\omega(\vec{\gamma}_{0,0};g_{00}). 
\end{align*}
Using the factorization for $\rho^{q+1}(\vec{\gamma}_0)$ and summing, we get
\begin{align*}
(\partial & \delta'-\delta'\partial+\partial'\Delta)\omega(\vec{\gamma};f) \\
          & =\rho^q((\delta_0\vec{\gamma})_0)\Big{[}\rho_0^1(i(g_{00}^{h_{01}...h_{0q}}))^{-1}\big{(}\omega(\vec{\gamma}_{0,0};f^{h_{00}})+[\rho_0^1(i(f^{h_{00}...h_{0q}}))-I]\omega(\vec{\gamma}_{0,0};g_{00})\big{)}-\omega(\vec{\gamma}_{0,0};f^{h_{00}i(g_{00})})\Big{]} \\
          & =\rho^q((\delta_0\vec{\gamma})_0)\Big{[}\rho_0^1(i(g_{00}^{h_{01}...h_{0q}}))^{-1}\big{(}\rho_0^1(i(f^{h_{00}...h_{0q}}))\omega(\vec{\gamma}_{0,0};g_{00})-\omega(\vec{\gamma}_{0,0};f^{h_{00}}g_{00})+\omega(\vec{\gamma}_{0,0};f^{h_{00}})\big{)}+ \\
          & \qquad +\rho_0^1(i(g_{00}^{h_{01}...h_{0q}}))^{-1}\omega(\vec{\gamma}_{0,0};f^{h_{00}}g_{00})-\omega(\vec{\gamma}_{0,0};f^{h_{00}i(g_{00})})+\omega(\vec{\gamma}_{0,0};g_{00}^{-1})+ \\
          & \qquad\qquad +\rho_0^1(i(g_{00}^{h_{01}...h_{0q}}))^{-1}\omega(\vec{\gamma}_{0,0};g_{00})-\omega(\vec{\gamma}_{0,0};g_{00}^{-1})\Big{]} \\
          & =\rho^q((\delta_0\vec{\gamma})_0)\Big{[}\rho_0^1(i(g_{00}^{h_{01}...h_{0q}}))^{-1}\delta_{(1)}\omega(\vec{\gamma}_{0,0};f^{h_{00}},g_{00})+\delta_{(1)}\omega(\vec{\gamma}_{0,0};g_{00}^{-1},f^{h_{00}}g_{00})-\delta_{(1)}\omega(\vec{\gamma}_{0,0};g_{00}^{-1},g_{00})\Big{]} \\
          & =\Delta\delta_{(1)}\omega(\vec{\gamma};f).
\end{align*}
\end{proof}
As announced, we regard this latter proposition as saying that the commutativity of diagram \ref{page r=1} is controlled by elements in a third lattice  
\begin{eqnarray*}
\xymatrix{
\vdots           & \vdots            & \vdots              &       \\ 
C(H^3\times G^2,W) & C(\G^3\times G^2,W) & C(\G_2^3\times G^2,W) & \dots \\
C(H^2\times G^2,W) & C(\G^2\times G^2,W) & C(\G_2^2\times G^2,W) & \dots \\
C(H\times G^2,W)   & C(\G\times G^2 ,W)  & C(\G_2\times G^2,W)   & \dots \\
C(G^2,W)           & C(G^2,W)            & C(G^2,W)	           & \dots  
}
\end{eqnarray*} 
that can too be filled with maps turning it into a non-commuting double complex. Inductively, we will be able to encode its commutativity in further pages filled with appropriate groupoid cochain complexes in rows and columns. This process will end up leaving a three dimensional lattice of spaces
\begin{eqnarray*}
C^{p,q}_r(\G,\phi):=C(\G_p^q\times G^r,W),
\end{eqnarray*}
where, we use the convention $\G_0 =H$ and for $r=0$, the coefficients take values on $V$. We will devote the rest of this section to spelling out how to endow this lattice with the structure of a complex whose second cohomology classifies the extensions of the Lie $2$-group by the $2$-vector space. 

\subsubsection{Pieces of the triple complex: $p$-pages}
In this subsubsection, we spell out one of the main pieces that form the referred triple complex, the constant $p$-pages. \\
We start by outlining some general facts about actions of Lie groups and representations. Let $H$ and $G$ be Lie groups, and let $H$ act on the right on $G$. Now, consider the following array
\begin{eqnarray*}
\xymatrix{
 H\ltimes G \ar@<0.5ex>[r] \ar@<-0.5ex>[r] \ar@<0.5ex>[d] \ar@<-0.5ex>[d] & H \ar@<0.5ex>[d] \ar@<-0.5ex>[d] \\
 G \ar@<0.5ex>[r] \ar@<-0.5ex>[r] & \ast ,
} 
\end{eqnarray*}
where the top groupoid is a bundle of Lie groups and the left groupoid is the (right) action groupoid for the action of $H$ on $G$ that we write
\begin{eqnarray*}
\Lf{s}(h;g)=g^{h}.
\end{eqnarray*}
We claim that this arrangement is actually a double Lie group. A word of warning though, given that the notation coincide, the reader might be inclined to believe that there is a Lie $2$-group somewhere in this diagram, but there is none. In this particular case, the action of $H$ on $G$ is not even asked to be by automorphisms, nor is any crossed module appearing. Now, since each of the sides of the square above is clearly a Lie groupoid, we will prove the claim by showing that the horizontal structural maps are indeed groupoid morphisms. First of all, the square commutes as all compositions land in $\ast$. We verify that $\Tp{s}=\Tp{t}$ respects units and compositions,
\begin{align*}
\Tp{s}(1;g) & =1\textnormal{ and} & \Tp{s}((h_1;g)\vJoin(h_2;g^{h_1})) & =\Tp{s}(h_1h_2;g) \\
            &	 & 													   & =h_1h_2=\Tp{s}(h_1;g)\Tp{s}(h_2;g^{h_1}).
\end{align*}
As for the unit map, 
\begin{align*}
\Lf{s}(\Tp{u}(h)) & =\Lf{s}(h;1) & \Lf{t}(\Tp{u}(h)) & =\Lf{t}(h;1) \\
                  & =1^{h}=1=u(\ast),  &                   & =1=u(\ast),
\end{align*}
and it is also compatible with units and multiplication:
\begin{align*}
\Tp{u}(1) & =(1;1)\textnormal{ and} & \Tp{u}(h_1h_2) & =(h_1h_2;1)\\
          &	                        &  				 & =(h_1;1)\vJoin(h_2;1^{h_1})=\Tp{u}(h_1)\vJoin\Tp{u}(h_2).
\end{align*}
Finally, $(H\ltimes G)_{\Tp{s}}\times_{\Tp{t}}(H\ltimes G)\cong H\times G^2$ and $\xymatrix{H\ltimes G^2\ar@<0.5ex>[r] \ar@<-0.5ex>[r] & G^2}$ is the action groupoid for the diagonal action. We verify that the multiplication is a Lie groupoid homomorphism:
\begin{align*}
\Lf{s}(\Tp{m}(h;g_1,g_2)) & =\Lf{s}(h;g_1g_2) \\
                          & =(g_1g_2)^h \\
                          & =g_1^hg_2^h=m(\Lf{s}^2(h;g_1,g_2)).
\end{align*}
and
\begin{align*}
\Lf{t}(\Tp{m}(h;g_1,g_2)) & =\Lf{t}(h;g_1g_2) \\
                          & =g_1g_2=m(\Lf{t}^2(h;g_1,g_2)).
\end{align*} 
Furthermore, the multiplication is also compatible with units and multiplication:
\begin{align*}
\Tp{m}(1;g_1,g_2) & =(1;g_1g_2),
\end{align*}
and
\begin{align*}
((h_1;g_1)\vJoin(h_2;g_1^{h_1})) & \Join ((h_1;g_2)\vJoin(h_2;g_2^{h_1})) \\
								 & =(h_1h_2;g_1)\Join (h_1h_2;g_2) \\
                                 & =\Tp{m}(h_1h_2;g_1,g_2) \\
                                 & =(h_1h_2;g_1g_2) \\
                                 & =(h_1;g_1g_2)\vJoin(h_2;(g_1g_2)^{h_1}) \\
                                 & =\Tp{m}(h_1;g_1,g_2)\vJoin\Tp{m}(h_2;g_1^{h_1},g_2^{h_1}) \\
                                 & =((h_1;g_1)\Join(h_1;g_2))\vJoin ((h_2;g_1^{h_1})\Join(h_2;g_2^{h_1})).
\end{align*}
At this point, we remark that for a vector space $W$, the previous verification obviously works too for the action by conjugation of $GL(W)$ on itself. To be consistent with the convention above, we would need to consider the right action of $GL(W)$ on $GL(W)$; either way, we conclude that 
\begin{eqnarray*}
\xymatrix{
 GL(W)\ltimes GL(W) \ar@<0.5ex>[r] \ar@<-0.5ex>[r] \ar@<0.5ex>[d] \ar@<-0.5ex>[d] & GL(W) \ar@<0.5ex>[d] \ar@<-0.5ex>[d] \\
 GL(W) \ar@<0.5ex>[r] \ar@<-0.5ex>[r] & \ast 
} 
\end{eqnarray*}
is a double Lie groupoid.
\begin{lemma}\label{p-pageModel}
Let $\xymatrix{H\ltimes G \ar[r] & GL(W)\ltimes GL(W):(h;g) \ar@{|->}[r] & (\rho_H(h),\rho_G(g))}$ be a map of double Lie groups. Then 
\begin{eqnarray*}
\xymatrix{
\vdots                                          & \vdots                                                & \vdots                       & \\ 
C(H^3,W) \ar[r]^{\delta'}\ar[u]          & C(H^3\times G,W) \ar[r]^{\delta'}\ar[u]          & C(H^3\times G^2,W) \ar[r]\ar[u]          & \dots \\
C(H^2,W) \ar[r]^{\delta'}\ar[u]^{\delta} & C(H^2\times G,W) \ar[r]^{\delta'}\ar[u]^{\delta} & C(H^2\times G^2,W) \ar[r]\ar[u]^{\delta} & \dots \\
C(H,W) \ar[r]^{\delta'}\ar[u]^{\delta}   & C(H\times G,W) \ar[r]^{\delta'}\ar[u]^{\delta}   & C(H\times G^2,W) \ar[r]\ar[u]^{\delta}   & \dots \\
W \ar[r]^{\delta'}\ar[u]^{\delta}        & C(G,W) \ar[r]^{\delta'}\ar[u]^{\delta}           & C(G^2,W)	\ar[r]\ar[u]^{\delta}          & \dots  
}
\end{eqnarray*}
where the rows and the columns are respectively the complexes of Lie groupoid cochains for the Lie group bundles
\begin{eqnarray*}
\xymatrix{
H^q\times G \ar@<0.5ex>[r] \ar@<-0.5ex>[r] & H^q 
} 
\end{eqnarray*}
on $\xymatrix{H^q\times W \ar[r] & H^q}$ and for the (right) transformation groupoids
\begin{eqnarray*}
\xymatrix{
 H\ltimes G^r \ar@<0.5ex>[r] \ar@<-0.5ex>[r] & G^r 
} 
\end{eqnarray*}
on $\xymatrix{G^r\times W \ar[r] & G^r}$ is a double complex.
\end{lemma}
%We remark that, though it is hidden in this notation, we are always considering the opposite group $H^{op}$.
\begin{proof}
We start by pointing out that if the given map is a map of double Lie groups, its restriction to the bottom and right groupoids give respectively representations $\rho_G$ and $\rho_H$ of $G$ and $H$ on $W$. Thus, the first row and first column are indeed complexes. In order to get the representations of the other rows and columns, we are going to pull-back these Lie group representations along the groupoid homomorphisms $\Lf{s}_q$ and $\Tp{t}_r$. Specifically, we define the representation $\rho_G^q$ of $\xymatrix{H^q\times G \ar@<0.5ex>[r] \ar@<-0.5ex>[r] & H^q}$ on $\xymatrix{H^q\times W \ar[r] & H^q}$ by 
\begin{eqnarray*}
\rho_G^q(h_1,...,h_q;g):=\Lf{s}_q^*\rho_G(h_1,...,h_q;g)=\rho_G(g^{h_1...h_q}) ,
\end{eqnarray*}
where indeed $(s_H)_q^*W=H^q\times W$. Analogously, we define the representation $\rho_H^r$ of $\xymatrix{H\ltimes G^r \ar@<0.5ex>[r] \ar@<-0.5ex>[r] & G^r}$ on $\xymatrix{G^r\times W \ar[r] & G^r}$ by 
\begin{eqnarray*}
\rho_H^r(h;g_1,...,g_r):=(\iota_H\circ\Tp{t}_r)^*\rho_H(h;g_1,...,g_r)=\rho_H(h)^{-1} ,
\end{eqnarray*}
and this time round, $t_r^*W=G^r\times W$. This, of course, is a right representation and, accordingly the complex appearing in each column is a complex of groupoid cochains with values in one of these. In the sequel, we write $(\vec{g})^{h}:=(g_1^{h},...,g_r^{h})$ for $h\in H$ and $\vec{g}\in G^r$. \\
Having clarified what the representations are, we show that the spaces of cochains are the right ones. On the one hand,
\begin{eqnarray*}
(H^q\times G)^{(r)}=\lbrace (\vec{h}_1,g_1;...;\vec{h}_r,g_r)\in(H^q\times G)^{r}:\Tp{s}(\vec{h}_j,g_j)=\vec{h}_j=\vec{h}_{j+1}=\Tp{t}(\vec{h}_{j+1},g_{j+1})\rbrace ;
\end{eqnarray*}
therefore, $(H^q\times G)^{(r)}\cong H^q\times G^r$ and the diffeomorphism will obviously be given by assigning $\xymatrix{(\vec{h},g_1;...;\vec{h},g_r) \ar@{|->}[r] & (\vec{h};g_1,...,g_r)}$. Since the representation is taken on a vector bundle which is trivial, its pull-back along the final target map will be trivial as well, and its sections will coincide with smooth functions to $W$, i.e. $C^r(H^q\times G,H^q\times W)=C(H^q\times G^r,W)$. On the other hand, 
\begin{eqnarray*}
(H\ltimes G^r)^{(q)}=\lbrace (h_1,\vec{g}_1;...;h_q,\vec{g}_r)\in(H\times G^r)^{q}:\Lf{s}(h_j,\vec{g}_j)=\vec{g}_j^{h_j}=\vec{g}_{j+1}=\Lf{t}(h_{j+1},\vec{g}_{j+1})\rbrace ;
\end{eqnarray*}
therefore, $(H\ltimes G^r)^{(q)}\cong H^q\times G^r$ and this time around the diffeomorphism will be given by the assignment $\xymatrix{(h_1,\vec{g};h_{2},(\vec{g})^{h_1};...;h_q,(\vec{g})^{h_1...h_{q-1}}) \ar@{|->}[r] & (h_1,...,h_q;\vec{g})}$. Since, again, the representation is taken on a vector bundle which is trivial, the pull-back along the initial source is trivial as well, and its sections will coincide with smooth functions to $W$ and, in so, with the space of $r$-cochains that we reviewed in the previous paragraph, i.e. $C^q(H\ltimes G^r,G^r\times W)=C(H^q\times G^r,W)$. \\
We will use the diffeomorphisms of the latter discussion to write formulas for the face maps of the simplicial structure: For $\vec{h}=(h_0,...,h_q)\in H^{q+1}$ and $\vec{g}=(g_0,...,g_r)\in G^{r+1}$, since
\begin{eqnarray*}
\delta _j (\vec{h})=
  \begin{cases}
    (h_1,...,h_q)                               & \quad \text{if } j=0  \\
    (h_0,...,h_{j-2},h_{j-1}h_j,h_{j+1}...,h_q) & \quad \text{if } 0<j\leq q\\
    (h_0,...,h_{q-1})                           & \quad \text{if } j=q+1 
  \end{cases}
\end{eqnarray*}
and
\begin{eqnarray*}
\delta'_k (\vec{g})=
  \begin{cases}
    (g_1,...,g_r)                                 & \quad \text{if } k=0  \\
    (g_0,...,g_{k-2},g_{k-1}g_k,g_{k+1},...,g_r)  & \quad \text{if } 0<k\leq r\\
    (g_0,...,g_{r-1})                             & \quad \text{if } k=r+1 ,
  \end{cases}
\end{eqnarray*}
\begin{eqnarray*}
\delta_j(\vec{h};\vec{g})=
  \begin{cases}
    (\delta_0\vec{h};(\vec{g})^{h_0}) & \quad \text{if } j=0 \\
    (\delta_{j}\vec{h};\vec{g})       & \quad \text{otherwise} 
  \end{cases}
\end{eqnarray*}
and
\begin{eqnarray*}
\delta'_k(\vec{h};\vec{g}) & =(\vec{h};\delta'_k\vec{g}).
\end{eqnarray*}
The only thing left to prove is the statement itself, that the generic square
\begin{eqnarray}\label{generic p-square}
\xymatrix{ 
C(H^{q+1}\times G^r,W) \ar[r]^{\delta'}            & C(H^{q+1}\times G^{r+1},W)          \\
C(H^q\times G^r,W) \ar[r]^{\delta'}\ar[u]^{\delta} & C(H^q\times G^{r+1},W) \ar[u]^{\delta} 
}
\end{eqnarray}
commutes. Indeed, let $\omega\in C(H^q\times G^r,W)$ and $\vec{h}$ and $\vec{g}$ as above. Then, 
\begin{align*}
\delta'\delta\omega(\vec{h};\vec{g}) & =\rho_G^{q+1}(\vec{h};g_0)\delta\omega(\vec{h};\delta'_0\vec{g})+\sum_{k=1}^{r+1}(-1)^{k}\delta\omega(\vec{h};\delta'_k\vec{g}),
\end{align*}
while
\begin{align*}
\delta\delta'\omega(\vec{h};\vec{g}) & =\delta'\omega(\delta_0\vec{h};(\vec{g})^{h_0})+\sum_{j=1}^{q}(-1)^{j}\delta'\omega(\delta_j\vec{h};\vec{g})+(-1)^{q+1}\rho_H^{r+1}(h_q;\vec{g})\delta'\omega(\delta_{q+1}\vec{h};\vec{g}).
\end{align*}
We expand further to make evident the common terms:
\begin{align*}
\delta'\delta\omega(\vec{h};\vec{g}) & =\rho_G^{q+1}(\vec{h};g_0)\Big{[}\omega(\delta_0\vec{h};(\delta'_0\vec{g})^{h_0})+\sum_{j=1}^{q}(-1)^{j}\omega(\delta_j\vec{h};\delta'_0\vec{g})+(-1)^{q+1}\rho_H^r(h_q;\delta'_0\vec{g})\omega(\delta_{q+1}\vec{h};\delta'_0\vec{g})\Big{]}+ \\
									 & \quad +\sum_{k=1}^{r+1}(-1)^{k}\Big{[}\omega(\delta_0\vec{h};(\delta'_k\vec{g})^{h_0})+\sum_{j=1}^{q}(-1)^{j}\omega(\delta_j\vec{h};\delta'_k\vec{g})+(-1)^{q+1}\rho_H^r(h_q;\delta'_k\vec{g})\omega(\delta_{q+1}\vec{h};\delta'_k\vec{g})\Big{]}											   
\end{align*}
and
\begin{align*}
\delta\delta'\omega(\vec{h};\vec{g}) & =\rho_G^q(\delta_0\vec{h};g_0^{h_0})\omega(\delta_0\vec{h};\delta'_0(\vec{g})^{h_0})+\sum_{k=1}^{r+1}(-1)^{k}\omega(\delta_0\vec{h};\delta'_k(\vec{g})^{h_0})+ \\
									 & \quad +\sum_{j=1}^{q}(-1)^{j}\Big{[}\rho_G^q(\delta_j\vec{h};g_0)\omega(\delta_j\vec{h};\delta'_0\vec{g})+\sum_{k=1}^{r+1}(-1)^{k}\omega(\delta_j\vec{h};\delta'_k\vec{g})\Big{]}+ \\
									 & \quad +(-1)^{q+1}\rho_H^{r+1}(h_q;\vec{g})\Big{[}\rho_G^q(\delta_{q+1}\vec{h};g_0)\omega(\delta_{q+1}\vec{h};\delta'_0\vec{g})+\sum_{k=1}^{r+1}(-1)^{k}\omega(\delta_{q+1}\vec{h};\delta'_k\vec{g})\Big{]}.
\end{align*}
The desired equality follows now by noticing the following identities. First, we obviously have that
\begin{eqnarray*}
(\delta'_k\vec{g})^{h_0}=\delta'_k(\vec{g})^{h_0}.
\end{eqnarray*}
Then, since for every $(h;\vec{g})\in H\times G^r$,
\begin{eqnarray*}
\rho_H^r(h_q;\delta'_k\vec{g})=\rho_H(h_q)^{-1}=\rho_H^{r+1}(h_q;\vec{g}).
\end{eqnarray*}
Further, since for every $(h_1,...,h_q;g)\in H^q\times G$,
\begin{eqnarray*}
\rho_G^{q+1}(\vec{h};g_0)=\rho_G(g_0^{h_0...h_q})=\rho_G((g_0^{h_0})^{h_{1}...h_q})=\rho_G^{q}(\delta_{0}\vec{h};g_0^{h_0}).
\end{eqnarray*}
Also, for all values $0<j\leq q$
\begin{eqnarray*}
\rho_G^{q+1}(\vec{h};g_0)=\rho_G(g_0^{h_0...h_q})=\rho_G(g_0^{h_0...(h_{j-1}h_j)...h_q})=\rho_G^{q}(\delta_{j}\vec{h};g_0).
\end{eqnarray*}
Finally, as $\rho_H\times\rho_G$ is a double Lie group map
\begin{align*}
\rho_G(\Lf{s}(h;g)) & =\Lf{s}(\rho_H(h),\rho_G(g)) \\
\rho_G(g^h)         & =\rho_H(h)^{-1}\rho_G(g)\rho_H(h);
\end{align*}
thereby implying,
\begin{align*}
\rho_G^{q+1}(\vec{h};g_0) & =\rho_G(g_0^{h_0...h_q}) \\
						 & =\rho_G((g_0^{h_0...h_{q-1}})^{h_q})=\rho_H(h_q)^{-1}\rho_G(g_0^{h_0...h_{q-1}})\rho_H(h_q),
\end{align*}
\begin{eqnarray*}
\rho_G^{q+1}(\vec{h};g_0)\rho_H^r(h_q;\delta'_0\vec{g})=\rho_H^{r+1}(h_q;\vec{g})\rho_G^q(\delta_{q+1}\vec{h};g_0)
\end{eqnarray*}
and the commutativity of the square \ref{generic p-square}.

\end{proof} 
With this latter lemma at hand, we go back to treating the $p$-pages of the triple complex of Lie $2$-group $\G=\xymatrix{G \ar[r]^i & H}$ with values in a $2$-representation $\rho$ on the $2$-vector space $\xymatrix{W \ar[r]^\phi & V}$. For any given $p$, $\G_p$ acts on the right on $G$, by pulling back the action of $H$ along the ``final target'' map $t_p$, that is 
\begin{eqnarray*}
g^{\gamma}:=g^{t_p(\gamma)},
\end{eqnarray*}
for $g\in G$ and $\gamma\in\G_p$. We prove that the representations 
\begin{align*}
\rho_G(g)           & :=\rho_0^1(i(g)) \\
\rho_{\G_p}(\gamma) & :=\rho_0^1(t_p(\gamma)) \\
\end{align*} 
verify the hypothesis of the lemma. First of all, of course these are indeed representations as they are the pull-back of the representation $\rho^1_0$ along the homomorphisms $i$ and $t_p$ respectively. Thus, after carefully looking at the compatibility with the whole structure, one realizes the only thing left to prove is that the vertical source maps are respected, or what is the same, that the equation
\begin{eqnarray*}
\rho_G(\Lf{s}(\gamma;g))=\Lf{s}(\rho_{\G_p}(\gamma),\rho_G(g))
\end{eqnarray*}
holds. However, this follows easily as a consequence of the equations defining the crossed module structure,
\begin{align*}
\rho_G(\Lf{s}(\gamma;g)) & =\rho_0^1(i(g^{t_p(\gamma)})) \\
						 & =\rho_0^1(t_p(\gamma)^{-1}i(g)t_p(\gamma)) \\
						 & =\rho_0^1(t_p(\gamma))^{-1}\rho_0^1(i(g))\rho_0^1(t_p(\gamma))=\rho_{\G_p}(\gamma)^{-1}\rho_G(g)\rho_{\G_p}(\gamma).
\end{align*}
%Maybe, it is of help to point out that in the definition of $\rho_{\G_p}$ we take the inverse in order to get a representation of $\G_p^{op}$ as in the discussion above. \\
As a matter of fact, the $p$-pages of the triple complex will coincide with the ones that we get out of this lemma, with the caveat that the first column takes values in $V$ instead of W.
\begin{lemma}\label{r-cx}
Let $\gamma_1,...,\gamma_q\in\G_p$ and $g\in G$ and $\omega\in C(\G_p^q,V)$. Set
\begin{eqnarray*}
\partial'\omega(\gamma_1,...,\gamma_q;g)=\rho_0^1(t_p(\gamma_1)...t_p(\gamma_q))^{-1}\rho_1(g)\omega(\gamma_1,...,\gamma_q),
\end{eqnarray*}
then the first term of the $q$th row of the $p$-page of the triple complex can be replaced by 
\begin{eqnarray*}
\xymatrix{
C(\G_p^q,V) \ar[r]^{\partial'\quad} & C(\G_p^q\times G,W) \ar[r]^{\delta'} & C(\G_p^q\times G^2,W) \ar[r] & ...
}
\end{eqnarray*}
still yielding a complex.
\end{lemma}
\begin{proof}
Clearly, the only thing to prove is that, for $\omega\in C(\G_p^q,V)$, $\delta'\partial'\omega=0$.
\begin{align*}
\delta'\partial' & \omega(\gamma_1,...,\gamma_q;g_0,g_1) \\
                 & =\rho_G^{q}(\gamma_1,...,\gamma_q;g_0)\partial'\omega(\gamma_1,...,\gamma_q;g_1)-\partial'\omega(\gamma_1,...,\gamma_q;g_0g_1)+\partial'\omega(\gamma_1,...,\gamma_q;g_0) \\
                 & =\rho_0^1(i(g_0^{t_p(\gamma_1)...t_p(\gamma_q)}))\rho_0^1(t_p(\gamma_1)...t_p(\gamma_q))^{-1}\rho_1(g_1)\omega(\gamma_1,...,\gamma_q)+ \\
				 & \qquad-\rho_0^1(t_p(\gamma_1)...t_p(\gamma_q))^{-1}\rho_1(g_0g_1)\omega(\gamma_1,...,\gamma_q)+\rho_0^1(t_p(\gamma_1)...t_p(\gamma_q))^{-1}\rho_1(g_0)\omega(\gamma_1,...,\gamma_q) \\
				 & =\rho_0^1(t_p(\gamma_1)...t_p(\gamma_q))^{-1}\Big{(}\rho_0^1(i(g_0))\rho_1(g_1)\omega(\gamma_1,...,\gamma_q)-\rho_1(g_0g_1)\omega(\gamma_1,...,\gamma_q)+\rho_1(g_0)\omega(\gamma_1,...,\gamma_q)\Big{)}, 
\end{align*} 
and the result follows from the relation $\rho_0^1(i(g_0))=I+\rho_1(g_0)\circ\phi$ and the fact that $\rho_1$ is a Lie group homomorphism landing in $GL(\phi)_1$.

\end{proof}
We now combine these latter two lemmas to get the sought after $p$-pages of the triple complex.
\begin{prop}\label{(q,r)-doubleCx}
For each $p$, 
\begin{eqnarray*}
\xymatrix{
\vdots                                             & \vdots                                                   & \vdots                           &       \\ 
C(\G_p^3,V) \ar[r]^{\partial'\quad}\ar[u]          & C(\G_p^3\times G,W) \ar[r]^{\delta'}\ar[u]          & C(\G_p^3\times G^2,W) \ar[r]\ar[u]          & \dots \\
C(\G_p^2,V) \ar[r]^{\partial'\quad}\ar[u]^{\delta} & C(\G_p^2\times G,W) \ar[r]^{\delta'}\ar[u]^{\delta} & C(\G_p^2\times G^2,W) \ar[r]\ar[u]^{\delta} & \dots \\
C(\G_p,V) \ar[r]^{\partial'\quad}\ar[u]^{\delta}   & C(\G_p\times G,W) \ar[r]^{\delta'}\ar[u]^{\delta}   & C(\G_p\times G^2,W) \ar[r]\ar[u]^{\delta} & \dots \\
V \ar[r]^{\delta_{(1)}\quad}\ar[u]^{\delta}        & C(G,W) \ar[r]^{\delta'}\ar[u]^{\delta}              & C(G^2,W)	\ar[r]\ar[u]^{\delta}        & \dots  
}
\end{eqnarray*}
is a double complex.
\end{prop}
\begin{proof}
Due to lemmas \ref{G-coh[V]} and \ref{r-cx}, each row is a complex, and clearly so is each column. Disregarding the first column of squares, lemma \ref{p-pageModel} says that we have got a double complex. Now, in order to finish the proof, one needs to check that the generic square 
\begin{eqnarray*}
\xymatrix{ 
C(\G_p^{q+1},V) \ar[r]^{\partial'\quad}            & C(\G_p^{q+1}\times G,W)          \\
C(\G_p^q,V) \ar[r]^{\partial'\quad}\ar[u]^{\delta} & C(\G_p^q\times G,W) \ar[u]^{\delta} 
}
\end{eqnarray*}
in the first column commutes. First, for $q=0$, let $\gamma\in\G_p$, $g\in G$ and $v\in V$,
\begin{align*}
\partial'\delta v(\gamma;g) & =\rho_0^1(t_p(\gamma))^{-1}\rho_1(g)\delta v(\gamma) \\
							& =\rho_0^1(t_p(\gamma))^{-1}\rho_1(g)(\rho_0^0(t_p(\gamma))v-v),
\end{align*}
whereas
\begin{align*}
\delta\delta_{(1)}v(\gamma;g) & =\delta_{(1)}v(g^{t_p(\gamma)})-\rho_{\G_p}^1(\gamma;g)\delta_{(1)}v(g) \\
							  & =\rho_1(g^{t_p(\gamma)})v-\rho_0^1(t_p(\gamma))^{-1}\rho_1(g)v \\
							  & =\rho_0^1(t_p(\gamma))^{-1}\rho_1(g)(\rho_0^0(t_p(\gamma))v-v).
\end{align*}
As for the other values of $q$, let $\vec{\gamma}=(\gamma_0,...,\gamma_q)\in\G_p^{q+1}$, then
\begin{align*}
\partial'\delta\omega(\vec{\gamma};g) & =\rho_0^1(t_p(\gamma_0)...t_p(\gamma_q))^{-1}\rho_1(g)\delta\omega(\vec{\gamma}) \\
							          & =\rho_0^1(t_p(\gamma_0)...t_p(\gamma_q))^{-1}\rho_1(g)\Big{(}\rho_0^0(t_p(\gamma_0))\omega(\delta_{0}\vec{\gamma})+\sum_{j=1}^{q+1}(-1)^{j}\omega(\delta_j\vec{\gamma})\Big{)},
\end{align*}
and
\begin{align*}
\delta\partial'\omega(\vec{\gamma};g) & =\partial'\omega(\delta_0\vec{\gamma};g^{t_p(\gamma_0)})+\sum_{j=1}^{q}(-1)^{j}\partial'\omega(\delta_j\vec{\gamma};g)+(-1)^{q+1}\rho_{\G_p}^1(\gamma_q;g)\partial'\omega(\delta_{q+1}\vec{\gamma};g) \\
							          & =\rho_0^1(t_p(\gamma_1)...t_p(\gamma_q))^{-1}\rho_1(g^{t_p(\gamma_0)})\omega(\delta_0\vec{\gamma})+ \\
							          & \qquad +\sum_{j=1}^{q}(-1)^{j}\rho_0^1(t_p(\gamma_q)...t_p(\gamma_{j-1}\gamma_j)...t_p(\gamma_q))^{-1}\rho_1(g)\omega(\delta_j\vec{\gamma})+ \\
							          & \qquad\qquad +(-1)^{q+1}\rho_0^1(t_p(\gamma_q))^{-1}\rho_0^1(t_p(\gamma_0)...t_p(\gamma_{q-1}))^{-1}\rho_1(g)\omega(\delta_{q+1}\vec{\gamma}) \\
							          & =\rho_0^1(t_p(\gamma_0)...t_p(\gamma_q))^{-1}\rho_1(g)\Big{(}\rho_0^0(t_p(\gamma_0))\omega(\delta_{0}\vec{\gamma})+\sum_{j=1}^{q+1}(-1)^{j}\omega(\delta_j\vec{\gamma})\Big{)}.
\end{align*}
\end{proof}
In the course of this section, we stuck to the convention that $\delta$ and $\delta'$ were respectively the vertical and the horizontal differential. In the sequel we will part from this, and assume a convention closer to the one in the triple complex of Lie $2$-algebra cochains. Specifically, the differentials in the $r$-direction, that played the part of horizontal differentials in this section, will be written $\delta_{(1)}$ instead.  

\subsubsection{Pieces of the triple complex: $p$-complexes and difference maps}
We start by spelling out the differentials in the $p$-direction. The generic complex
\begin{eqnarray}\label{Gp p-cx}
\xymatrix{
C(H^q\times G^r,W) \ar[r]^{\partial} & C(\G^q\times G^r,W) \ar[r]^{\partial} & C(\G_2^q\times G^r,W) \ar[r] & \dots 
}
\end{eqnarray} 
is the groupoid cochain complex of the product groupoid 
\begin{eqnarray*}
\xymatrix{\G^q\ltimes G^r \ar@<0.5ex>[r]\ar@<-0.5ex>[r] & H^q\times G^r},
\end{eqnarray*}
with respect to the representation
\begin{eqnarray}\label{p-rep}
(\gamma_1,...,\gamma_q;\vec{f})\cdot(h_1,...,h_q;\vec{f},w):=(h_1i(g_1),...,h_qi(g_q);\vec{f},\rho_0^1(i(pr_G(\gamma_1\vJoin ...\vJoin\gamma_q)))^{-1}w), 
\end{eqnarray}
where $\gamma_k=(g_k,h_k)\in\G$ and $\vec{f}\in G^r$, on the trivial bundle $\xymatrix{H^q\times G\times W \ar[r]^{\quad pr_1} & H^q\times G}$. This is essentially the same representation that we had on the rows of the page \ref{page r=1}; therefore, it is a pull-back of $\Delta^W$ of the representation up to homotopy induced by the $2$-representation along the projection onto $\xymatrix{\G^q \ar@<0.5ex>[r]\ar@<-0.5ex>[r] & H^q}$ composed with the multiplication. This differential fits conveniently as the $q$-pages turn into honest double complexes.
\begin{prop}\label{Gp q-pageDoubleCx}
With the differentials described and for constant $q$, $C(\G_p^q\times G^r,W)$ is a double complex.
\end{prop}
\begin{proof}
Let rows and columns be respectively defined by making $r$ and $p$ constant. Then rows are defined to be complexes of Lie groupoid cochains, and, by lemma \ref{r-cx}, the columns are complexes as well. It is left to see that the differentials commute. We will use the identification 
\begin{eqnarray}
\gamma_b=\begin{pmatrix}
\gamma_{0b} \\
\vdots \\
\gamma_{pb}
\end{pmatrix}=\begin{pmatrix}
g_{0b} & h_{0b} \\
\vdots & \vdots \\
g_{pb} & h_{pb}
\end{pmatrix} 
\end{eqnarray}
for $\vec{\gamma}=(\gamma_1,...,\gamma_q)\in\G_{p+1}^q$. Now, let $\omega\in C(\G_p^q,V)$ and $f\in G$, then
\begin{align*}
\partial\partial'\omega(\vec{\gamma};f) & =\rho_0^1(i(pr_G(\gamma_{01}\vJoin ...\vJoin\gamma_{0q})))^{-1}\partial'\omega(\partial_0\vec{\gamma};f)+\sum_{j=1}^{p+1}(-1)^j\partial'\omega(\partial_j\vec{\gamma};f) \\
                                        & =\rho_0^1(i(pr_G(\gamma_{01}\vJoin ...\vJoin\gamma_{0q})))^{-1}\rho_0^1(t_{p}(\partial_0\gamma_1)...t_{p}(\partial_0\gamma_q))^{-1}\rho_1(f)\omega(\partial_0\vec{\gamma})+ \\
                                        & \qquad +\sum_{j=1}^{p+1}(-1)^j\rho_0^1(t_{p}(\partial_j\gamma_1)...t_{p}(\partial_j\gamma_q))^{-1}\rho_1(f)\omega(\partial_j\vec{\gamma}) \\
                                        & =\rho_0^1(i(pr_G(\gamma_{01}\vJoin ...\vJoin\gamma_{0q})))^{-1}\rho_0^1(s(\gamma_{01})...s(\gamma_{0q}))^{-1}\rho_1(f)\omega(\partial_0\vec{\gamma})+ \\
                                        & \qquad +\sum_{j=1}^{p+1}(-1)^j\rho_0^1(t(\gamma_{01})...t(\gamma_{0q}))^{-1}\rho_1(f)\omega(\partial_j\vec{\gamma}) \\
                                        & =\rho_0^1(s(\gamma_{01}\vJoin ...\vJoin\gamma_{0q})i(pr_G(\gamma_{01}\vJoin ...\vJoin\gamma_{0q})))^{-1}\rho_1(f)\omega(\partial_0\vec{\gamma})+ \\
                                        & \qquad +\sum_{j=1}^{p+1}(-1)^j\rho_0^1(t(\gamma_{01})...t(\gamma_{0q}))^{-1}\rho_1(f)\omega(\partial_j\vec{\gamma}) \\
                                        & =\rho_0^1(t(\gamma_{01})...t_p(\gamma_{0q}))^{-1}\rho_1(f)\sum_{j=0}^{p+1}(-1)^j\omega(\partial_j\vec{\gamma}) \\
                                        & =\rho_0^1(t_{p+1}(\gamma_1)...t_{p+1}(\gamma_q))^{-1}\rho_1(f)\partial\omega(\vec{\gamma})=\partial'\partial\omega(\vec{\gamma};f).
\end{align*}
On the other hand, for $\omega\in C(\G_p^q\times G^r,V)$ and $\vec{f}=(f_0,...,f_r)\in G^{r+1}$,
\begin{align*}
\delta_{(1)}\partial\omega(\vec{\gamma};\vec{f}) & =\rho_0^1(i(f_0^{t_{p+1}(\gamma_0)...t_{p+1}(\gamma_q)}))\partial\omega(\vec{\gamma};\delta_0\vec{f})+\sum_{k=1}^{r+1}(-1)^k\partial\omega(\vec{\gamma};\delta_k\vec{f}) \\
                                        & =\rho_0^1(i(f_0^{t_{p+1}(\gamma_0)...t_{p+1}(\gamma_q)}))\Big{(}\rho_0^1(i(pr_G(\gamma_{01}\vJoin ...\vJoin\gamma_{0q})))^{-1}\omega(\partial_0\vec{\gamma};\delta_0\vec{f})+\sum_{j=1}^{p+1}(-1)^j\omega(\partial_j\vec{\gamma};\delta_0\vec{f})\Big{)}+ \\
                                        & \qquad +\sum_{k=1}^{r+1}(-1)^k\Big{(}\rho_0^1(i(pr_G(\gamma_{01}\vJoin ...\vJoin\gamma_{0q})))^{-1}\omega(\partial_0\vec{\gamma};\delta_k\vec{f})+\sum_{j=1}^{p+1}(-1)^j\omega(\partial_j\vec{\gamma};\delta_k\vec{f})\Big{)} .
\end{align*}
Now, given that 
\begin{align*}
\rho_0^1(i(f_0^{t_{p+1}(\gamma_0)...t_{p+1}(\gamma_q)})) & =\rho_0^1(i(f_0^{t(\gamma_{00})...t(\gamma_{0q})})) \\
                                                 & =\rho_0^1(i(f_0^{t(\gamma_{00}\vJoin ...\vJoin\gamma_{0q})})) \\
                                                 & =\rho_0^1(i(f_0^{h_{00}...h_{0q}i(pr_G(\gamma_{00}\vJoin ...\vJoin\gamma_{0q}))})) \\
                                                 & =\rho_0^1(i(pr_G(\gamma_{00}\vJoin ...\vJoin\gamma_{0q})))^{-1}\rho_0^1(i(f_0^{h_{00}...h_{0q}}))\rho_0^1(i(pr_G(\gamma_{00}\vJoin ...\vJoin\gamma_{0q}))); 
\end{align*}
we can group the terms 
\begin{align*}
\delta_{(1)}\partial\omega(\vec{\gamma};\vec{f}) & =\rho_0^1(i(pr_G(\gamma_{01}\vJoin ...\vJoin\gamma_{0q})))^{-1}\Big{(}\rho_0^1(i(f_0^{h_{00}...h_{0q}}))\omega(\partial_0\vec{\gamma};\delta_0\vec{f})+\sum_{k=1}^{r+1}(-1)^k\omega(\partial_0\vec{\gamma};\delta_k\vec{f})\Big{)} \\
                                           & \qquad +\sum_{j=1}^{p+1}(-1)^j\Big{(}\rho_0^1(i(f_0^{t_{p+1}(\gamma_0)...t_{p+1}(\gamma_q)}))\omega(\partial_j\vec{\gamma};\delta_0\vec{f})+\sum_{k=1}^{r+1}(-1)^k\omega(\partial_j\vec{\gamma};\delta_k\vec{f})\Big{)} .
\end{align*}
Since $t_p(\partial_0\gamma_b)=t(\gamma_{1b})=s(\gamma_{0b})=h_{0b}$ and $t_p(\partial_j\gamma_b)=t(\gamma_{0b})=t_{p+1}(\gamma_b)$ for any other $j$, 
\begin{align*}
\delta_{(1)}\partial\omega(\vec{\gamma};\vec{f}) & =\rho_0^1(i(pr_G(\gamma_{01}\vJoin ...\vJoin\gamma_{0q})))^{-1}\delta_{(1)}\omega(\partial_0\vec{\gamma};\vec{f})+\sum_{j=1}^{p+1}(-1)^j\delta_{(1)}\omega(\partial_j\vec{\gamma};\vec{f})= \partial\delta_{(1)}\omega(\vec{\gamma};\vec{f})
\end{align*}
thus yielding the desired equality and the proposition.

\end{proof}
Now that we have got the differential in the $p$-direction, we have the grid that forms the triple complex. Recall that we claimed that the $(r+1)$th page appeared in attempting to understand the non-commutativity of the $r$th page. The upcoming propositions are the formalization of this comment. We will be using the following notation: For an element $\vec{f}=(f_1,...,f_r)\in G^r$ and integers $1\leq a<b\leq r$, define
\begin{eqnarray*}
\vec{f}_{[a,b]}:=(f_a,f_{a+1},...,f_{b-1},f_b) & \textnormal{and} & \vec{f}_{[a,b)}:=(f_a,f_{a+1},...,f_{b-2},f_{b-1}).
\end{eqnarray*}
With this shorthand, we define for $1\leq n\leq r$
\begin{eqnarray*}
\xymatrix{
\Delta^n:\G\times G^r \ar[r] & G^r
}
\end{eqnarray*}
at the level of spaces. If $\gamma=(g,h)^T\in\G$ and $\vec{f}=(f_1,...,f_r)\in G^r$,
\begin{align*}
\Delta^n(\gamma;\vec{f}):=\Big{(}\big{(}\vec{f}_{[1,n)}\big{)}^{hi(g)},g^{-1},\big{(}\vec{f}_{[n,r)}\big{)}^h\Big{)}
\end{align*}
and we are ready to state and prove the following result.
\begin{prop}\label{GpStarTop(p,q)=(0,0)}
For $r>1$,
\begin{eqnarray*}
(-1)^r(\delta\partial-\partial\delta)=\Delta\circ\delta_{(1)}-\delta_{(1)}\circ\Delta,
\end{eqnarray*}
where 
\begin{eqnarray*}
\xymatrix{
\Delta:C(G^{r+1},W) \ar[r] & C(\G\times G^r,W) 
} 
\end{eqnarray*}
is defined for $\vec{f}\in G^r$ and $\gamma=(g,h)^T\in\G$ by the formula
\begin{align*}
\Delta\omega(\gamma;\vec{f}) & =\rho_0^1(i(g))^{-1}\omega((\vec{f})^h,g)+\sum_{n=1}^r(-1)^{r-n}\big{[}\omega(\Delta^n(\gamma;\vec{f}),f_r^hg)-\omega(\Delta^n(\gamma;\vec{f}),g)\big{]} .   
\end{align*}
\end{prop}
\begin{proof}
We start by computing the left hand side of the equation. Let $\omega\in C(G^r,W)$ and recall that $\partial\omega$ is defined to be zero. On the other hand,
\begin{align*}
\partial\delta\omega(\gamma;\vec{f}) & =\rho_0^1(i(g))^{-1}\delta\omega(h;\vec{f})-\delta\omega(hi(g);\vec{f}) \\
                                     & =\rho_0^1(i(g))^{-1}\big{(}\omega((\vec{f})^h)-\rho_0^1(h)^{-1}\omega(\vec{f})\big{)}-\big{(}\omega((\vec{f})^{hi(g)})-\rho_0^1(hi(g))^{-1}\omega(\vec{f})\big{)} \\
                                     & =\rho_0^1(i(g))^{-1}\omega((\vec{f})^h)-\omega((\vec{f})^{hi(g)}). 
\end{align*}
Now, for the left hand side we consider
\begin{align*}
\Delta\delta_{(1)}\omega(\gamma;\vec{f}) & =\rho_0^1(i(g))^{-1}\delta_{(1)}\omega((\vec{f})^h,g)+\sum_{n=1}^r(-1)^{r-n}\big{[}\delta_{(1)}\omega(\Delta^n(\gamma;\vec{f}),f_r^hg)-\delta_{(1)}\omega(\Delta^n(\gamma;\vec{f}),g)\big{]}    
\end{align*}
and
\begin{align*}
\delta_{(1)}\Delta\omega(\gamma;\vec{f}) & =\rho_0^1(i(f_1^{hi(g)}))\Delta\omega(\gamma;\delta_0\vec{f})+\sum_{k=1}^r(-1)^k\Delta\omega(\gamma;\delta_k\vec{f}).
\end{align*}
Expanding further the first term of $\Delta\delta_{(1)}\omega$, 
\begin{align*}
\rho_0^1(i(g))^{-1}\delta_{(1)}\omega((\vec{f})^h,g)=\rho_0^1(i(g))^{-1}\Big{(}\rho_0^1(i(f_1^h)) & \omega(\delta_0(\vec{f})^h,g)+\sum_{k=1}^{r-1}(-1)^k\omega(\delta_k(\vec{f})^h,g)+ \\
& \qquad +(-1)^r\omega(\delta_r(\vec{f})^h,f_r^hg)+(-1)^{r+1}\omega((\vec{f})^h)\Big{)}
\end{align*}
one realizes that all but the last two terms get cancelled by the first terms of each $\Delta\omega(\gamma;\delta_k\vec{f})$; indeed, all of these are given by the expression
\begin{align*}
\rho_0^1(i(f_1^{hi(g)})\rho_0^1(i(g))^{-1}\omega((\delta_0\vec{f})^h,g)+\sum_{k=1}^r(-1)^k\rho_0^1(i(g))^{-1}\omega((\delta_k\vec{f})^h,g)
\end{align*}
and of course
\begin{eqnarray*}
(\delta_k\vec{f})^h=\delta_k(\vec{f})^h
\end{eqnarray*}
as well as $g^{-1}f_1^h=f_1^{hi(g)}g^{-1}$. Hence, out of these terms, the following remains in the difference
\begin{eqnarray}\label{leftOver1}
\rho_0^1(i(g))^{-1}\Big{(}(-1)^r\omega(\delta_r(\vec{f})^h,f_r^hg)+(-1)^{r+1}\omega((\vec{f})^h)\Big{)}-(-1)^r\rho_0^1(i(g))^{-1}\omega((\delta_r\vec{f})^h,g).
\end{eqnarray}
We proceed to consider the term $n=1$ in $\Delta\delta_{(1)}\omega$,
\begin{align*}
(- & 1)^{r-1}\Big{(}\delta_{(1)}\omega(\Delta^1(\gamma;\vec{f}),f_r^hg)-\delta_{(1)}\omega(\Delta^1(\gamma;\vec{f}),g)\Big{)}.
\end{align*}
Since $\Delta^1(\gamma;\vec{f})=(g^{-1},(\delta_r\vec{f})^h)$, this expression is
\begin{align*}
 & (-1)^{r-1}\Big{(}\rho_0^1(i(g))^{-1}\omega((\delta_r\vec{f})^h,f_r^hg)-\omega(g^{-1}f_1^h,\delta_0(\delta_r\vec{f})^h,f_r^hg)+\sum_{k=1}^{r-2}(-1)^{k+1}\omega(g^{-1},\delta_k(\delta_r\vec{f})^h,f_r^hg)+ \\
 & \qquad\qquad +(-1)^{r}\omega(g^{-1},\delta_{r-1}(\delta_r\vec{f})^h,(f_{r-1}f_r)^hg)+(-1)^{r+1}\omega(g^{-1},(\delta_r\vec{f})^h)+ \\
 & \qquad -\rho_0^1(i(g))^{-1}\omega((\delta_r\vec{f})^h,g)+\omega(g^{-1}f_1^h,\delta_0(\delta_r\vec{f})^h,g)+\sum_{k=1}^{r-2}(-1)^{k}\omega(g^{-1},\delta_k(\delta_r\vec{f})^h,g)+ \\
 & \qquad\qquad +(-1)^{r-1}\omega(g^{-1},\delta_{r-1}(\delta_r\vec{f})^h,f_{r-1}^hg)+(-1)^{r}\omega(g^{-1},(\delta_r\vec{f})^h)\Big{)}.
\end{align*}
The leading terms (the ones that are being acted on) cancel out with two of the terms in equation \ref{leftOver1} leaving behind $(-1)^{r+1}\rho_0^1(i(g))^{-1}\omega((\vec{f})^h)$ which is the first term of what we computed for the right hand side. The last terms of each piece will cancel with each other and the terms inside the sum will cancel one another with the second term of each of the $\Delta\omega(\gamma;\delta_k\vec{f})$ in the sum of $\Delta\delta_{(1)}\omega$. Indeed, these have got the form 
\begin{align*}
    (-1)^{k+1}(-1)^{r-2}(\omega(g^{-1},\delta_{r}(\delta_k\vec{f})^h,f_r^{h}g)-\omega(g^{-1},\delta_{r}(\delta_k\vec{f})^h,g)).
\end{align*}
Updating the difference, it reads
\begin{align*}
(-1)^{r+1}\rho_0^1(i(g))^{-1}\omega((\vec{f})^h) & +(-1)^r\Big{(}\omega(g^{-1}f_1^h,\delta_0(\delta_r\vec{f})^h,f_r^hg)-\omega(g^{-1}f_1^h,\delta_0(\delta_r\vec{f})^h,g) \Big{)} \Big{)}+ \\
 & \qquad -\omega(g^{-1},\delta_{r-1}(\delta_r\vec{f})^h,(f_{r-1}f_r)^hg)+\omega(g^{-1},\delta_{r-1}(\delta_r\vec{f})^h,f_{r-1}^hg).
\end{align*}
Inductively, the terms in the parenthesis are going to cancel out with the first terms of the next $n$; indeed, consider $n=2$ in $\Delta\delta_{(1)}\omega$:
\begin{align*}
   (- & 1)^{r-2}\Big{(}\delta_{(1)}\omega(\Delta^2(\gamma;\vec{f}),f_r^hg)-\delta_{(1)}\omega(\Delta^2(\gamma;\vec{f}),g)\Big{)} \\
    & =(-1)^{r-2}\Big{(}\rho_0^1(i(f_1^{hi(g)}))^{-1}\omega(g^{-1},(\vec{f}_{[2,r)})^h,f_r^hg)+ \\
    & \qquad -\omega(f_1^{hi(g)}g^{-1},(\vec{f}_{[2,r)})^h,f_r^hg)+\omega(f_1^{hi(g)},g^{-1}f_2^h,(\vec{f}_{[3,r)})^h,f_r^hg)+ \\
    & \qquad\quad +\sum_{k=1}^{r-3}(-1)^{k}\omega(f_1^{hi(g)},g^{-1},\delta_k(\vec{f}_{[2,r)})^h,f_r^hg)+ \\
    & \qquad\qquad +(-1)^{r-2}\omega(f_1^{hi(g)},g^{-1},(\vec{f}_{[2,r-1)})^h,(f_{r-1}f_r)^hg)+(-1)^{r-1}\omega(f_1^{hi(g)},g^{-1},(\vec{f}_{[2,r)})^h)+ \\
    & \quad -\rho_0^1(i(f_1^{hi(g)}))^{-1}\omega(g^{-1},(\vec{f}_{[2,r)})^h,g)+ \\
    & \qquad +\omega(f_1^{hi(g)}g^{-1},(\vec{f}_{[2,r)})^h,g)-\omega(f_1^{hi(g)},g^{-1}f_2^h,(\vec{f}_{[3,r)})^h,g)+ \\
    & \qquad\quad +\sum_{k=1}^{r-3}(-1)^{k+1}\omega(f_1^{hi(g)},g^{-1},\delta_k(\vec{f}_{[2,r)})^h,g)+ \\
    & \qquad\qquad +(-1)^{r-1}\omega(f_1^{hi(g)},g^{-1},(\vec{f}_{[2,r-1)})^h,f_{r-1}^hg)+(-1)^{r}\omega(f_1^{hi(g)},g^{-1},(\vec{f}_{[2,r)})^h)\Big{)};
\end{align*}
hence, using the identity $g^{-1}f^h=f^{hi(g)}g^{-1}$, we see that the leading terms cancel with the terms in the neglected $\rho_0^1(i(f_1^{hi(g)}))\Delta\omega(\gamma;\delta_0\vec{f})$. We also see that the second pair of terms cancel with the ones surviving from $n=1$ as announced, as well as the terms in the sum cancel one another with the third term of each of the $\Delta\omega(\gamma;\delta_k\vec{f})$ in the sum of $\Delta\delta_{(1)}\omega$. Inductively, we are left with the terms remaining from $n=r$ in $\Delta\delta_{(1)}\omega$, these are
\begin{eqnarray*}
(-1)^{r-r}\Big{(}(-1)^{r}\omega((\delta_r\vec{f})^{hi(g)},g^{-1}f_r^hg)-(-1)^{r+1}\omega((\delta_r\vec{f})^{hi(g)},g^{-1}g) \Big{)},
\end{eqnarray*}
and since $g^{-1}f_r^hg=f_r^{hi(g)}$ and $\omega((\delta_r\vec{f})^{hi(g)},1)\equiv 0$, the result follows.

% These computations can be found in "zokuji 29."
\end{proof}
\begin{prop}\label{GpStarTop q=0}
For $r>1$,
\begin{eqnarray*}
(-1)^r(\delta\partial-\partial\delta)=\Delta\circ\delta_{(1)}-\delta_{(1)}\circ\Delta,
\end{eqnarray*}
where 
\begin{eqnarray*}
\xymatrix{
\Delta:C(G^{r+1},W) \ar[r] & C(\G_{p+1}\times G^r,W) 
} 
\end{eqnarray*}
is defined for $\vec{f}\in G^r$ and $\gamma=(\gamma_0,...,\gamma_p)^T\in\G_{p+1}$ with each $\gamma_{a}$ identified with $(g_{a},h_{a})$ by the formula
\begin{align*}
\Delta\omega(\gamma;\vec{f}) & =\rho_0^1(i(g_1)))^{-1}\omega((\vec{f})^{h_0},g_0)+\sum_{n=1}^r(-1)^{r-n}\big{(}\omega(\Delta^n(\gamma_0;\vec{f}),f_r^{h_0}g_0)-\omega(\Delta^n(\gamma_0;\vec{f}),g_0)\big{)} .   
\end{align*}
\end{prop}
\begin{proof}
The proof of this proposition gets established as 
\begin{align*}
\partial\delta\omega(\gamma;\vec{f}) & =\rho_0^1(i(g_0))^{-1}\delta\omega(\partial_0\gamma;\vec{f})+\sum_{j=1}^{p+1}(-1)^j\delta\omega(\partial_j\gamma;\vec{f}) \\
                                     & =\rho_0^1(i(g_0))^{-1}\big{(}\omega((\vec{f})^{h_0})-\rho_0^1(h_0)^{-1}\omega(\vec{f})\big{)}+\sum_{j=1}^{p+1}(-1)^j\big{(}\omega((\vec{f})^{h_0i(g_0)})-\rho_0^1(h_0i(g_0))^{-1}\omega(\vec{f})\big{)};
\end{align*}
thus, if we split in cases for $p$ even and odd
\begin{eqnarray*}
\delta\partial\omega(\gamma;\vec{f}) & =\begin{cases}
    \omega((\vec{f})^{h_0i(g_0)})-\rho_0^1(h_0i(g_0))^{-1}\omega(\vec{f})                & \quad p\text{ odd}   \\
    0 & \quad \text{otherwise}, 
  \end{cases}
\end{eqnarray*}
we get that 
\begin{align*}
 (\delta\partial-\partial\delta)\omega(\gamma;\vec{f}) & =\omega((\vec{f})^{h_0i(g_0)})-\rho_0^1(i(g_0))^{-1}\omega((\vec{f})^{h_0}),
\end{align*}                         
and the rest of the previous proof applies to this case as well.

\end{proof}
\begin{prop}\label{GpStarTop}
For $r>1$,
\begin{eqnarray*}
(-1)^r(\delta\partial-\partial\delta)=\Delta\circ\delta_{(1)}-\delta_{(1)}\circ\Delta,
\end{eqnarray*}
where 
\begin{eqnarray*}
\xymatrix{
\Delta:C(\G_p^q\times G^{r+1},W) \ar[r] & C(\G_{p+1}^{q+1}\times G^r,W) 
}
\end{eqnarray*}
is defined for $\vec{f}\in G^r$, $\vec{\gamma}=(\gamma_0,...,\gamma_q)\in\G_{p+1}^{q+1}$ with each $\gamma_b=(\gamma_{0b},...,\gamma_{pb})^T$ and each $\gamma_{ab}=(g_{ab},h_{ab})$ by the formula
\begin{align*}
\Delta(\vec{\gamma};\vec{f})  & =\rho_0^1(i(pr_G(\gamma_{01}\vJoin ...\vJoin\gamma_{0q})))^{-1}\Big{[}\rho_0^1(i(g_{00}^{h_{01}...h_{0q}}))^{-1}\omega(\vec{\gamma}_{0,0};(\vec{f})^{h_{00}},g_{00})+ \\
& \qquad +\sum_{n=1}^r(-1)^{r-n}\big{(}\omega(\vec{\gamma}_{0,0};\Delta^n(\gamma_{00};\vec{f}),f_r^{h_{00}}g_{00})-\omega(\vec{\gamma}_{0,0};\Delta^n(\gamma_{00};\vec{f}),g_{00})\big{)}\Big{]}.
\end{align*}
\end{prop}
We already proved the special case $r=1$, which needed a slightly different formula for $\Delta$; also, the previous propositions settle the cases $r>1$ and $q=0$.
\begin{proof}
This statement follows from a similar argument as the ones in the previous proofs after noticing the following. First, computing the right hand side yields
\begin{align*}
\delta\partial\omega(\vec{\gamma};\vec{f}) & =\partial\omega(\delta_0\vec{\gamma};(\vec{f})^{h_{00}i(g_{00})})+\sum_{k=1}^{q}(-1)^k\partial\omega(\delta_k\vec{\gamma};\vec{f})+(-1)^{q+1}\rho_0^1(h_{0q}i(g_{0q}))^{-1}\partial\omega(\delta_{q+1}\vec{\gamma};\vec{f}) \\
                                           & =\rho_0^1(i(pr_G(\gamma_{01}\vJoin ...\vJoin\gamma_{0q})))^{-1}\omega(\vec{\gamma}_{0,0};(\vec{f})^{h_{00}i(g_{00})})+\sum_{j=1}^{p+1}(-1)^j\omega(\partial_j\delta_0\vec{\gamma};(\vec{f})^{h_{00}i(g_{00})})+ \\
                                           & \quad +\sum_{k=1}^{q}(-1)^k\Big{[}\rho_0^1(i(pr_G(\gamma_{00}\vJoin ...\vJoin\gamma_{0q})))^{-1}\omega(\partial_0\delta_k\vec{\gamma};\vec{f})+\sum_{j=1}^{p+1}(-1)^j\omega(\partial_j\delta_k\vec{\gamma};\vec{f})\Big{]}+ \\
                                           & \quad +(-1)^{q+1}\rho_0^1(h_{0q}i(g_{0q}))^{-1}\Big{[}\rho_0^1(i(pr_G(\gamma_{00}\vJoin ...\vJoin\gamma_{0q-1})))^{-1}\omega(\vec{\gamma}_{0,q};\vec{f})+\sum_{j=1}^{p+1}(-1)^j\omega(\partial_j\delta_{q+1}\vec{\gamma};\vec{f})\Big{]},
\end{align*}
and  
\begin{align*}
\partial\delta\omega(\vec{\gamma};\vec{f}) & =\rho_0^1(i(pr_G(\gamma_{00}\vJoin ...\vJoin\gamma_{0q})))^{-1}\delta\omega(\partial_0\vec{\gamma};\vec{f})+\sum_{j=1}^{p+1}(-1)^j\delta\omega(\partial_j\vec{\gamma};\vec{f}) \\
                                           & =\rho_0^1(i(pr_G(\gamma_{00}\vJoin ...\vJoin\gamma_{0q})))^{-1}\Big{[}\omega(\vec{\gamma}_{0,0};(\vec{f})^{h_{00}})+\sum_{k=1}^{q}(-1)^k\omega(\delta_k\partial_0\vec{\gamma};\vec{f})+(-1)^{q+1}\rho_0^1(h_{0q})^{-1}\omega(\vec{\gamma}_{0,q};\vec{f})\Big{]}+ \\
                                          % & \quad -\Big{(}\omega(\delta_0\partial_1\vec{\gamma};(\vec{f})^{h_{00}})+\sum_{k=1}^{q}(-1)^k\omega(\delta_k\partial_1\vec{\gamma};\vec{f})+(-1)^{q+1}\rho_0^1(h_{0q})^{-1}\omega(\vec{\gamma}_{0,q};\vec{f})\Big{)}+ \\
                                           & \quad +\sum_{j=1}^{p+1}(-1)^j\Big{[}\omega(\delta_0\partial_j\vec{\gamma};(\vec{f})^{h_{00}i(g_{00})})+\sum_{k=1}^{q}(-1)^k\omega(\delta_k\partial_j\vec{\gamma};\vec{f})+(-1)^{q+1}\rho_0^1(h_{0q}i(g_{0q}))^{-1}\omega(\delta_{q+1}\partial_j\vec{\gamma};\vec{f})\Big{]}.
\end{align*}
Clearly, the latter line will cancel with the terms that appear to the right in the last lines of $\delta\partial\omega$. Also, since 
\begin{eqnarray*}
\rho_0^1(hi(g))^{-1}\rho_0^1(i(f))^{-1}=\rho_0^1(i(f)hi(g))^{-1}=\rho_0^1(hi(f^h)i(g)))^{-1}=\rho_0^1(i(f^hg))^{-1}\rho_0^1(h)^{-1},
\end{eqnarray*}
replacing $f$ by $pr_G(\gamma_{00}\vJoin ...\vJoin\gamma_{0q-1})=g_{00}^{h_{01}...h_{0q-1}}g_{01}^{h_{02}...h_{0q-1}}...g_{0q-2}^{h_{0q-1}}g_{0q-1}$, the only surviving terms when subtracting are
\begin{align*}
(\delta\partial-\partial\delta)\omega(\vec{\gamma};\vec{f}) & =\rho_0^1(i(pr_G(\gamma_{01}\vJoin ...\vJoin\gamma_{0q})))^{-1}\Big{(}\omega(\vec{\gamma}_{0,0};(\vec{f})^{h_{00}i(g_{00})})-\rho_0^1(i(g_{00}^{h_{01}...h_{0q-1}}))^{-1}\omega(\vec{\gamma}_{0,0};(\vec{f})^{h_{00}})\Big{)}.
\end{align*}
The rest of the proof follows using the same argument of the proofs of the previous two propositions.

\end{proof}
Attempting to define a complex structure for the indexing by ``counter-diagonal planes''
\begin{eqnarray}
C^k_{tot}(\G,\phi):=\bigoplus_{p+q+r=k}C^{p,q}_r(\G,\phi)
\end{eqnarray}
using the differentials we have got so far, one runs into the same issue as in the case of Lie $2$-algebras, i.e. $d^2=0$ fails to hold. Specifically, if one lets $\omega\in C^{p,q}_r(\G,\phi)$ and assumes $d=\partial+\delta^{(r)}+\delta_{(1)}+\Delta$ (up to signs), we roughly represent the equations that are to vanish in the following diagram:
\begin{eqnarray*}
\xymatrix{
  & & \Delta\circ\Delta & &  \\
  & [\partial,\Delta] & \Delta\omega \ar[l]\ar[u]\ar[r]\ar[d] & [\Delta,\delta] &  \\
\partial^2 & \partial\omega \ar[l]\ar[u]\ar[r]\ar[d] & [\partial,\Delta,\delta_{(1)},\delta]  & \delta\omega \ar[l]\ar[u]\ar[r]\ar[d] & \delta^2 \\
  & [\partial,\delta_{(1)}]  & \delta_{(1)}\omega \ar[l]\ar[u]\ar[r]\ar[d] & [\delta_{(1)},\delta] &  \\
  & & \delta_{(1)}^2 & &  
}
\end{eqnarray*}
We list the equations that we have already encountered and proved that are indeed zero:
\begin{itemize}
    \item $\partial^2\omega=0\in C^{p+2,q}_r(\G,\phi)$, because $\partial$ is the differential of groupoid cochains,
    \item $\delta^2\omega=0\in C^{p,q+2}_r(\G,\phi)$, again, because $\delta$ is a differential of groupoid cochains,
    \item $\delta_{(1)}^2\omega=0\in C^{p,q}_{r+2}(\G,\phi)$, one last time, because $\delta_{(1)}$ is a differential of groupoid cochains, 
    \item $(\delta\delta_{(1)}-\delta_{(1)}\delta)\omega=0\in C^{p,q+1}_{r+1}(\G,\phi)$, because the $p$-pages are honest double complexes (cf. Proposition \ref{(q,r)-doubleCx}), 
    \item $(\partial\delta_{(1)}-\delta_{(1)}\partial)\omega=0\in C^{p+1,q}_{r+1}(\G,\phi)$, because the $q$-pages are honest double complexes (cf. Proposition \ref{Gp q-pageDoubleCx}) and
    \item $((-1)^r(\delta\partial-\partial\delta)-\Delta\circ\delta_{(1)}+\delta_{(1)}\circ\Delta)\omega=0\in C^{p+1,q+1}_{r-1}(\G,\phi)$ is the contents of proposition \ref{GpStarTop}.
\end{itemize}
% ---------
In contrast with the case of Lie $2$-algebras, we cannot ensure the vanishing of any other relation. Instead, we are going to have two separate classes of higher difference maps.

% ---------

%$\xymatrix{\Delta:C^{\bullet,\bullet}_1(\G,\phi) \ar[r] & C^{\bullet,\bullet}_0(\G,\phi)}$,

\subsubsection{Higher difference maps 1: landing on the front page corners}
We claim that the relation labeled $[\partial,\Delta]$ does not vanish. In fact, for $\omega\in C(G,W)$ and $\gamma=(g_0,...,g_{p+1},h)\in\G_{p+2}$, we have got
\begin{align*}
\partial\Delta\omega(\gamma) & =\sum_{j=0}^{p+2}(-1)^j\Delta\omega(\partial_j\gamma) \\
                             & =\rho_0^0(t_{p}(\partial_0\partial_0\gamma))\phi(\omega(g_1))-\rho_0^0(t_{p}(\partial_0\partial_1\gamma))\phi(\omega(g_1g_0))+\sum_{j=2}^{p+2}(-1)^j\rho_0^0(t_{p}(\partial_0\partial_j)\gamma)\phi(\omega(g_0)).
\end{align*}
Now,
\begin{eqnarray*}
\partial_0\partial_j\gamma=
  \begin{cases}
    (g_2,...,g_{p+1};h)                & \quad \text{if } j\in\lbrace 0,1\rbrace  \\
    (g_1,...,g_jg_{j-1},...,g_{p+1};h) & \quad \text{if } 1<j\leq p+1\\
    (g_1,...,g_{p};hi(g_{p+1}))      & \quad \text{if } j=p+2; 
  \end{cases}
\end{eqnarray*}
therefore, 
\begin{align*}
\partial\Delta\omega(\gamma) & =\rho_0^0(hi(g_{p+1}...g_2))\phi\big{(}\omega(g_1)-\omega(g_1g_0)\big{)}+\sum_{j=2}^{p+2}(-1)^j\rho_0^0(hi(g_{p+1}...g_1))\phi(\omega(g_0)) \\
                             & =\begin{cases}
    \rho_0^0(hi(g_{p+1}...g_2))\phi\big{(}\omega(g_1)-\omega(g_1g_0)+\rho_0^1(i(g_1))\omega(g_0)\big{)} & \quad \text{if } p\text{ even},  \\
    \rho_0^0(hi(g_{p+1}...g_2))\phi\big{(}\omega(g_1)-\omega(g_1g_0)\big{)}                              & \quad \text{otherwise}.
  \end{cases}
\end{align*}
On the other hand, 
\begin{align*}
\Delta\partial\omega(\gamma) & =\rho_0^0(hi(g_{p+1}...g_1))\phi(\partial\omega(g_0)) \\
                             & =\begin{cases}
    \rho_0^0(hi(g_{p+1}...g_1))\phi(\omega(g_0)) & \quad \text{if } p\text{ odd},  \\
    0                                        & \quad \text{otherwise};
  \end{cases}
\end{align*}
thus, if we define
\begin{eqnarray*}
\xymatrix{
\Delta_2^q:C(G^2,W) \ar[r] & C(\G_{p+2},V)
} \\
\Delta\omega(\gamma):=\rho_0^0(t_{p}(\partial_0\partial_0\gamma))\phi(\omega(g_1,g_0))
\end{eqnarray*}
for $\gamma=(g_0,...,g_{p+1},h)\in\G_{p+2}$, we have that
\begin{eqnarray*}
\partial\circ\Delta+\Delta\circ\partial=\Delta_2^q\circ\delta_{(1)}.
\end{eqnarray*}
In fact, such a relation holds for all difference maps $\xymatrix{\Delta:C^{\bullet,\bullet}_1(\G,\phi) \ar[r] & C^{\bullet,\bullet}_0(\G,\phi)}$.
\begin{lemma}\label{IV atSch}
Let 
\begin{eqnarray*}
\xymatrix{
\Delta_2^q:C(\G_p^q\times G^2,W) \ar[r] & C(\G_{p+2}^{q+1},V)
} 
\end{eqnarray*}
be defined by
\begin{eqnarray*}
\Delta_2^q\omega(\vec{\gamma}):=\rho_0^0(t_{p}(\partial_0\partial_0\gamma_0)...t_{p}(\partial_0\partial_0\gamma_q))\phi(\omega(\partial_0\vec{\gamma}_{0,0};g_{10},g_{00}))
\end{eqnarray*}
for $\vec{\gamma}=(\gamma_0,...,\gamma_q)\in\G_{p+2}^{q+1}$, where $\gamma_b=(\gamma_{0b},...,\gamma_{p+1b})^T$ and $\gamma_{ab}=(g_{ab},h_{ab})$. Then
\begin{eqnarray*}
\partial\circ\Delta+\Delta\circ\partial=\Delta_2^q\circ\delta_{(1)}.
\end{eqnarray*}
\end{lemma}
\begin{proof}
We already proved that the equation holds for $q=0$, so assume $q\geq 1$. Let $\omega\in C(\G_p^q\times G,W)$ and $\gamma$ as in the statement, then
\begin{align*}
\partial\Delta\omega(\vec{\gamma}) & =\sum_{j=0}^{p+2}(-1)^j\Delta\omega(\partial_j\vec{\gamma}) \\
                             & =\rho_0^0(t_{p}(\partial_0\partial_0\gamma_0)...t_{p}(\partial_0\partial_0\gamma_q))\phi(\omega((\partial_0\vec{\gamma})_{0,0};g_{10}))+ \\
                             & \qquad -\rho_0^0(t_{p}(\partial_0\partial_1\gamma_0)...t_{p}(\partial_0\partial_1\gamma_q))\phi(\omega((\partial_1\vec{\gamma})_{0,0};g_{10}g_{00}))+ \\
                             & \qquad\quad +\sum_{j=2}^{p+2}(-1)^j\rho_0^0(t_{p}(\partial_0\partial_j\gamma_0)...t_{p}(\partial_0\partial_j\gamma_q))\phi(\omega((\partial_j\vec{\gamma})_{0,0};g_{00})).
\end{align*}
Now, the simplicial relations that we wrote explicitly before the statement of the lemma hold for each $\gamma_b$, and therefore, 
\begin{align*}
\partial\Delta\omega(\vec{\gamma}) & =\rho_0^0(s(\gamma_{10})...s(\gamma_{1q}))\phi\big{(}\omega((\partial_0\vec{\gamma})_{0,0};g_{10})-\omega((\partial_0\vec{\gamma})_{0,0};g_{10}g_{00})\big{)}+ \\
                             & \qquad\quad +\sum_{j=2}^{p+2}(-1)^j\rho_0^0(t(\gamma_{10})...t(\gamma_{1q}))\phi(\omega((\partial_j\vec{\gamma})_{0,0};g_{00})).
\end{align*}
In contrast with the case $q=0$, this time around, there is an explicit dependence on $\partial_j\vec{\gamma}$, and, as a consequence, these terms do not cancel with one another. Nonetheless, on the other hand we have got 
\begin{align*}
\Delta\partial\omega(\vec{\gamma}) & =\rho_0^0(t_{p+1}(\partial_0\gamma_0)...t_{p+1}(\partial_0\gamma_q))\phi(\partial\omega(\vec{\gamma}_{0,0};g_{00})) \\
                             & =\rho_0^0(t(\gamma_{10})...t(\gamma_{1q}))\phi\big{(}\rho^q((\vec{\gamma}_{0,0})_0)\omega(\partial_0\vec{\gamma}_{0,0};g_{00})+\sum_{j=1}^{p+1}(-1)^j\omega(\partial_j\vec{\gamma}_{0,0};g_{00})\big{)},
\end{align*}
where we used the shorthand introduced above
\begin{eqnarray*}
\rho^{q}((\vec{\gamma}_{0,0})_0)=\rho_0^1(i(pr_G(\gamma_{11}\vJoin ...\vJoin\gamma_{1q})))^{-1}.
\end{eqnarray*}
Thus, we discover the sum to be
\begin{align*}
(\partial\Delta+\Delta\partial)\omega(\vec{\gamma}) & =\rho_0^0(s(\gamma_{10})...s(\gamma_{1q}))\phi\big{(}\omega((\partial_0\vec{\gamma})_{0,0};g_{10})-\omega((\partial_0\vec{\gamma})_{0,0};g_{10}g_{00})\big{)}+ \\
                                              & \qquad +\rho_0^0(t(\gamma_{10})...t(\gamma_{1q}))\phi\big{(}\rho^q((\vec{\gamma}_{0,0})_0)\omega(\partial_0\vec{\gamma}_{0,0};g_{00})\big{)}.
\end{align*}
and since $t(\gamma)=s(\gamma)i(pr_G(\gamma))$ and $t$ is a Lie group homomorphism
\begin{align*}
(\partial\Delta+\Delta\partial)\omega(\vec{\gamma}) & =\rho_0^0(s(\gamma_{10})...s(\gamma_{1q}))\phi\Big{(}\omega((\partial_0\vec{\gamma})_{0,0};g_{10})-\omega((\partial_0\vec{\gamma})_{0,0};g_{10}g_{00})\big{)}+ \\
                                              & \qquad +\rho_0^1\big{(}i(g_{10}^{h_{11}...h_{1q}})\big{)}\omega(\partial_0\vec{\gamma}_{0,0};g_{00})\Big{)}.
\end{align*}
Indeed, $h_{1b}=s(\gamma_{1b})=t(\gamma_{2b})=t_{p}(\partial_0\partial_0\gamma_b)=t_{p}(\partial_0(\vec{\gamma}_{0,0})_b)$; thus, the result follows.

\end{proof}
We proceed to study how the relation labeled $[\delta,\Delta]$ behaves in the special case $r=1$. \\
Let $\omega\in C(\G_p^q\times G,W)$ and consider $\vec{\gamma}=(\gamma_0,...,\gamma_{q+1})\in\G_{p+1}^{q+2}$ under the identification $\gamma_b=(\gamma_{0b},...,\gamma_{pb})^T$ with $\gamma_{ab}\sim(g_{ab},h_{ab})$. Let us compute separately $\delta\Delta\omega$ and $\Delta\delta\omega$. On the one hand, we've got
\begin{align*}
\delta\Delta\omega(\vec{\gamma}) & =\rho_0^0(t_{p+1}(\gamma_0))\Delta\omega(\delta_0\vec{\gamma})+\sum_{j=1}^{q+2}(-1)^j\Delta\omega(\delta_j\vec{\gamma}) \\
                                 & =\rho_0^0(h_{00}i(g_{00}))\big{(}\rho_0^0(t_{p}(\partial_0\gamma_1)...t_{p}(\partial_0\gamma_{q+1}))\phi\omega((\delta_0\vec{\gamma})_{0,0};g_{01})\big{)}+ \\
                                 & \qquad -\rho_0^0(t_{p}(\partial_0(\gamma_0\vJoin\gamma_1))...t_{p}(\partial_0\gamma_{q+1}))\phi\omega((\delta_1\vec{\gamma})_{0,0};g_{00}^{h_{01}}g_{01})+ \\
                                 & \qquad\qquad +\sum_{j=2}^{q+1}(-1)^j\rho_0^0(t_{p}(\partial_0\gamma_1)...t_{p}(\partial_0(\gamma_{j-1}\vJoin\gamma_j))...t_{p}(\partial_0\gamma_{q+1}))\phi\omega((\delta_j\vec{\gamma})_{0,0};g_{00})+ \\
                                 & \qquad\qquad\qquad +(-1)^{q+2}\rho_0^0(t_{p}(\partial_0\gamma_1)...t_{p}(\partial_0\gamma_{q}))\phi\omega((\delta_{q+2}\vec{\gamma})_{0,0};g_{00}) \\
                                 & =\rho_0^0(h_{00}...h_{0q+1})\phi\Big{(}\rho_0^1(i(g_{00}^{h_{01}...h_{0q+1}}))\omega(\delta_0\vec{\gamma}_{0,0};g_{01})-\omega(\delta_0\vec{\gamma}_{0,0};g_{00}^{h_{01}}g_{01})+ \\
                                 & \qquad\qquad +\sum_{j=2}^{q+1}(-1)^j\omega((\delta_j\vec{\gamma})_{0,0};g_{00})+(-1)^{q+2}\rho_0^1(h_{0q+1})^{-1}\omega((\delta_{q+2}\vec{\gamma})_{0,0};g_{00})\Big{)}.
\end{align*}
On the other hand,
\begin{align*}
\Delta\delta\omega(\vec{\gamma}) & =\rho_0^0(t_{p}(\partial_0\gamma_0)...t_{p}(\partial_0\gamma_{q+1}))\phi\big{(}\delta\omega(\vec{\gamma}_{0,0};g_{00})\big{)} \\
                                 & =\rho_0^0(h_{00}...h_{0q+1})\phi\Big{(}\omega(\delta_0\vec{\gamma}_{0,0};g_{00}^{t_p(\partial_0\gamma_1)})+ \\
                                 & \qquad\qquad +\sum_{j=1}^{q+1}(-1)^j\omega(\delta_j\vec{\gamma}_{0,0};g_{00})+(-1)^{q+1}\rho_0^1(t_p(\partial_0\gamma_{q+1}))^{-1}\omega(\delta_{q+1}\vec{\gamma}_{0,0};g_{00})\Big{)}.
\end{align*}
Now, since $(\delta_j\vec{\gamma})_{0,0}=\delta_{j-1}\vec{\gamma}_{0,0}$ for $j\geq 2$ and $t_p(\partial_0\gamma_1)=h_{01}$, 
\begin{align*}
(\delta\Delta+\Delta\delta)\omega(\vec{\gamma}) & =\rho_0^0(h_{00}...h_{0q+1})\phi\Big{(}\rho_0^1(i(g_{00}^{h_{01}...h_{0q+1}}))\omega(\delta_0\vec{\gamma}_{0,0};g_{01})-\omega(\delta_0\vec{\gamma}_{0,0};g_{00}^{h_{01}}g_{01})+\omega(\delta_0\vec{\gamma}_{0,0};g_{00}^{h_{01}})\Big{)} \\
                                                & =\rho_0^0(h_{00}...h_{0q+1})\phi\big{(}\delta_{(1)}\omega(\delta_0\vec{\gamma}_{0,0};g_{00}^{h_{01}},g_{01})\big{)}.
\end{align*}
We summarize this discussion in the following lemma.
\begin{lemma}\label{V atSch}
Let 
\begin{eqnarray*}
\xymatrix{
\Delta_2^p:C(\G_p^q\times G^2,W) \ar[r] & C(\G_{p+1}^{q+2},V)
} 
\end{eqnarray*}
be defined by
\begin{eqnarray*}
\Delta_2^p\omega(\vec{\gamma}):=\rho_0^0(t_{p}(\partial_0\gamma_0)...t_{p}(\partial_0\gamma_{q+1}))\phi(\omega(\delta_0\vec{\gamma}_{0,0};g_{00}^{h_{01}},g_{01}))
\end{eqnarray*}
for $\vec{\gamma}=(\gamma_0,...,\gamma_{q+1})\in\G_{p+1}^{q+2}$, where $\gamma_b=(\gamma_{0b},...,\gamma_{pb})^T$ and $\gamma_{ab}=(g_{ab},h_{ab})$. Then
\begin{eqnarray*}
\delta\circ\Delta+\Delta\circ\delta=\Delta_2^p\circ\delta_{(1)}.
\end{eqnarray*}
\end{lemma}
By adding these second difference maps to our candidate for differential, we partially solved the issue of it not squaring to zero. Nevertheless, simultaneously, we are adding more relations that need to be satisfied. For instance, we would now need to make sure that
\begin{eqnarray*}
\partial\Delta_2^q-\Delta_2^q\partial=0 & \textnormal{and} & \delta\Delta_2^q-\Delta_2^q\delta=0
\end{eqnarray*}
for elements in $C(\G_p^q\times G^2,W)$. This is not the case though. Instead, these differences will be controlled by a sequence of higher difference maps that we proceed to define. \\
Let 
\begin{eqnarray*}
\xymatrix{
\Delta_r^q:C(\G_p^q\times G^r,W) \ar[r] & C(\G_{p+r}^{q+1},V)
} 
\end{eqnarray*}
be defined by
\begin{eqnarray*}
\Delta_r^q\omega(\vec{\gamma}):=\rho_0^0(t_{p}(\partial_0^r\gamma_0)...t_{p}(\partial_0^r\gamma_q))\phi\big{(}\omega(\partial_0^{r-1}\vec{\gamma}_{0,0};g_{r-10},...,g_{00})\big{)}
\end{eqnarray*}
for $\vec{\gamma}=(\gamma_0,...,\gamma_q)\in\G_{p+r}^{q+1}$, where $\gamma_b=(\gamma_{0b},...,\gamma_{p+r-1b})^T$ and $\gamma_{ab}=(g_{ab},h_{ab})$, and let
\begin{eqnarray*}
\xymatrix{
\Delta_r^p:C(\G_p^q\times G^r,W) \ar[r] & C(\G_{p+1}^{q+r},V)
} 
\end{eqnarray*}
be defined by
\begin{eqnarray*}
\Delta_r^p\omega(\vec{\gamma}):=\rho_0^0(t_{p}(\partial_0\gamma_0)...t_{p}(\partial_0\gamma_{q+r-1}))\phi\big{(}\omega(\delta^{r-1}_0\vec{\gamma}_{0,0};g_{00}^{h_{01}...h_{0r-1}},g_{01}^{h_{02}...h_{0r-1}},...,g_{0r-2}^{h_{0r-1}},g_{0r-1})\big{)}
\end{eqnarray*}
for $\vec{\gamma}=(\gamma_0,...,\gamma_{q+r-1})\in\G_{p+1}^{q+r}$, where $\gamma_b=(\gamma_{0b},...,\gamma_{pb})^T$ and $\gamma_{ab}=(g_{ab},h_{ab})$. Thus defined we have got the following relations.
\begin{proposition}\label{FrontPageCorners}
If $\omega\in C(\G_p^q\times G^r,W)$, then
\begin{eqnarray*}
(-1)^{r+1}\partial\circ\Delta_r^q+\Delta_r^q\circ\partial=\Delta_{r+1}^q\circ\delta_{(1)} & \textnormal{and} & 
\delta\circ\Delta_r^p+(-1)^{r+1}\Delta_r^p\circ\delta=\Delta_{r+1}^p\circ\delta_{(1)}.
\end{eqnarray*}
\end{proposition}
\begin{proof}
For the first relation, we need to consider two separate cases $q=0$ and $q>0$. Let $\gamma=(\gamma_0,...,\gamma_{p+r})^T\in\G_{p+r+1}$ be identified with $\gamma\sim (g_0,...,g_{p+r},h)$, then
\begin{align*}
\partial\Delta_r^q\omega(\gamma) & =\sum_{j=0}^{p+r+1}(-1)^j\Delta_r^q\omega(\partial_j\gamma) \\
                             & =\rho_0^0(hi(g_{r+p}...g_{r+1}))\phi\Big{(}\omega(g_r,...,g_1)-\omega(g_r,...,g_2,g_1g_0)+...+(-1)^{r-1}\omega(g_r,g_{r-2}g_{r-1},...,g_1,g_0)+ \\
                             & \qquad +(-1)^{r}\omega(g_rg_{r-1},g_{r-2}...,g_1,g_0)\Big{)}+\sum_{j=r+1}^{r+p}(-1)^j\rho_0^0(hi(g_{r+p}...g_{r}))\phi(\omega(g_{r-1},...,g_0))+ \\
                             & \qquad\qquad +(-1)^{r+p+1}\rho_0^0(hi(g_{r+p})i(g_{r+p-1}...g_{r}))\phi(\omega(g_{r-1},...,g_0)).
\end{align*}
The values that the representation takes in are given by the following identifications:
\begin{eqnarray*}
\partial_0^r\partial_j\gamma\sim
  \begin{cases}
    (g_{r+1},...,g_{p+r};h)            & \quad \text{if } 0\leq j\leq r  \\
    (g_r,...,g_jg_{j-1},...,g_{p+r};h) & \quad \text{if } r<k\leq r+p\\
    (g_r,...,g_{p+r-1};hi(g_{p+r}))      & \quad \text{if } k=r+p+1. 
  \end{cases}
\end{eqnarray*}
Now, notice that the elements in the sum have got no dependence on $j$ and therefore are going to cancel with one another yielding 
\begin{align*}
\partial\Delta_r^q\omega(\gamma) & =\begin{cases}
    (-1)^{r+1}\rho_0^0(hi(g_{p+r}...g_{r+1}))\phi\big{(}\delta_{(1)}\omega(g_{r},...,g_0)-\rho_0^1(i(g_r))\omega(g_{r-1},...,g_0)\big{)} & \quad \text{if } p\text{ odd},  \\
    (-1)^{r+1}\rho_0^0(hi(g_{p+r}...g_{r+1}))\phi\big{(}\delta_{(1)}\omega(g_{r},...,g_0)\big{)}                                         & \quad \text{otherwise}.
  \end{cases}
\end{align*}
On the other hand, 
\begin{align*}
\Delta_r^q\partial\omega(\gamma) & =\rho_0^0(hi(g_{p+r}...g_r))\phi(\partial\omega(g_{r-1},...,g_0)) \\
                             & =\begin{cases}
    \rho_0^0(hi(g_{p+r}...g_r))\phi(\omega(g_{r-1},...,g_0)) & \quad \text{if } p\text{ odd},  \\
    0                                        & \quad \text{otherwise};
  \end{cases}
\end{align*}
thus, 
\begin{eqnarray*}
\partial\circ\Delta_r^q+(-1)^{r+1}\Delta_r^q\circ\partial=(-1)^{r+1}\Delta_{r+1}^q\circ\delta_{(1)}.
\end{eqnarray*}
Assume now that $q\geq 1$ and let $\vec{\gamma}=(\gamma_0,...,\gamma_q)\in\G_{p+r+1}^{q+1}$. Then
\begin{align*}
\partial\Delta_r^q\omega(\vec{\gamma}) & =\sum_{j=0}^{p+r+1}(-1)^j\Delta_r^q\omega(\partial_j\vec{\gamma}) \\
                             & =\rho_0^0(t_{p}(\partial_0^{r+1}\gamma_0)...t_{p}(\partial_0^{r+1}\gamma_q))\phi(\omega(\partial_0^{r-1}(\partial_0\vec{\gamma})_{0,0};g_{r0},...,g_{10}))+ \\
                             & \qquad -\rho_0^0(t_{p}(\partial_0^r\partial_1\gamma_0)...t_{p}(\partial_0^r\partial_1\gamma_q))\phi(\omega(\partial_0^{r-1}(\partial_1\vec{\gamma})_{0,0};g_{r0},...,g_{20},g_{10}g_{00}))+...+ \\
                             & \qquad\quad +(-1)^{r-1}\rho_0^0(t_{p}(\partial_0^{r}\partial_{r-1}\gamma_0)...t_{p}(\partial_0^{r}\partial_{r-1}\gamma_q))\phi(\omega(\partial_0^{r-1}(\partial_{r-1}\vec{\gamma})_{0,0};g_{r0},g_{r-20}g_{r-10},...,g_{10},g_{00}))+ \\
                             & \qquad\qquad +(-1)^{r}\rho_0^0(t_{p}(\partial_0^{r}\partial_{r}\gamma_0)...t_{p}(\partial_0^{r}\partial_{r}\gamma_q))\phi(\omega(\partial_0^{r-1}(\partial_{r}\vec{\gamma})_{0,0};g_{r0}g_{r-10},g_{r-20}...,g_{10},g_{00}))+ \\
                             & \qquad\qquad\quad +\sum_{j=r+1}^{p+r+1}(-1)^j\rho_0^0(t_{p}(\partial_0^{r}\partial_j\gamma_0)...t_{p}(\partial_0^{r}\partial_j\gamma_q))\phi(\omega(\partial_0^{r-1}(\partial_j\vec{\gamma})_{0,0};g_{r-10},...,g_{00})).
\end{align*}
Using the simplicial relations spelled out at the beginning of the proof coordinate-wise, together with the relation $h_{ab}=h_{a+1b}i(g_{a+1b})=h_{0b}i(g_{p+rb}...g_{a+1b})$, we get 
\begin{align*}
\partial\Delta_r^q\omega(\vec{\gamma}) & =(-1)^{r+1}\rho_0^0(h_{r0}...h_{rq}))\phi\Big{(}\delta_{(1)}\omega(\partial_0^{r}\vec{\gamma}_{0,0};g_{r0},...,g_{00})-\rho_0^1(i(g_{r0}^{h_{r1}...h_{rq}}))\omega(\partial_0^{r}\vec{\gamma}_{0,0};g_{r-10},...,g_{00})\Big{)}+ \\
                             & \qquad\quad +\sum_{j=r+1}^{p+r+1}(-1)^j\rho_0^0(h_{r0}i(g_{r0})...h_{rq}i(g_{rq}))\phi(\omega(\partial_0^{r-1}(\partial_j\vec{\gamma})_{0,0};g_{00})).
\end{align*}
As opposed to the case $q=0$, there is an explicit dependence on $j$ for the terms in the sum this time around. As a consequence, these terms won't cancel with one another. Computing the second term of the relation, we get 
\begin{align*}
\Delta_r^q\partial\omega(\vec{\gamma}) & =\rho_0^0(t_{p+1}(\partial_0^r\gamma_0)...t_{p+1}(\partial_0^r\gamma_q))\phi\big{(}\partial\omega(\partial_0^{r-1}\vec{\gamma}_{0,0};g_{r-10},...,g_{00})\big{)} \\
                             & =\rho_0^0(h_{r0}i(g_{r0})...h_{rq}i(g_{rq}))\phi\Big{(}\rho_0^1(i(g_{r1}^{h_{r2}...h_{rq}}g_{r2}^{h_{r3}...h_{rq}}...g_{rq-1}^{h_{rq}}g_{rq}))^{-1}\omega(\partial_0^r\vec{\gamma}_{0,0};g_{r-10},...,g_{00})+ \\
                             & \qquad +\sum_{j=1}^{p+1}(-1)^j\omega(\partial_j\partial_0^{r-1}\vec{\gamma}_{0,0};g_{r-10},...,g_{00})\Big{)}.
\end{align*}
Of course, 
\begin{eqnarray*}
h_{r0}i(g_{r0})...h_{rq}i(g_{rq})=h_{r0}...h_{rq}i(g_{r0}^{h_{r1}...h_{rq}}g_{r1}^{h_{r2}...h_{rq}}...g_{rq-1}^{h_{rq}}g_{rq});
\end{eqnarray*}
thus, by introducing a factor of $(-1)^{r+1}$ and carefully looking at the indices, we get
\begin{align*}
(\partial\Delta_r^q+(-1)^{r+1}\Delta_r^q\partial)\omega(\vec{\gamma}) & =(-1)^{r+1}\rho_0^0(h_{r0}...h_{rq})\phi\big{(}\delta_{(1)}\omega(\partial_0^{r}\vec{\gamma}_{0,0};g_{r0},...,g_{00})\big{)}=(-1)^{r+1}\Delta_{r+1}^q\delta_{(1)}\omega(\vec{\gamma})
\end{align*}
as desired. \\
As for the second relation, let $\vec{\gamma}=(\gamma_0,...,\gamma_{q+r})\in\G_{p+1}^{q+r+1}$ with the usual identifications. Let us compute separately $\delta\Delta_r^p\omega$ and $\Delta_r^p\delta\omega$. On the one hand, we have got
\begin{align*}
\delta\Delta_r^p & \omega(\vec{\gamma})=\rho_0^0(t_{p+1}(\gamma_0))\Delta_r^p\omega(\delta_0\vec{\gamma})+\sum_{j=1}^{q+r+1}(-1)^j\Delta_r^p\omega(\delta_j\vec{\gamma}) \\
                 & =\rho_0^0(h_{00}i(g_{00}))\Big{(}\rho_0^0(t_{p}(\partial_0\gamma_1)...t_{p}(\partial_0\gamma_{q+r}))\phi\big{(}\omega(\delta_0^{r-1}(\delta_0\vec{\gamma})_{0,0};g_{01}^{h_{02}...h_{0r}},...,g_{0r-1}^{h_{0r}},g_{0r})\big{)}\Big{)}+ \\
                 & \quad -\rho_0^0(t_{p}(\partial_0(\gamma_0\vJoin\gamma_1))...t_{p}(\partial_0\gamma_{q+r}))\phi\big{(}\omega(\delta_0^{r-1}(\delta_1\vec{\gamma})_{0,0};(g_{00}^{h_{01}}g_{01})^{h_{02}...h_{0r}},...,g_{0r-1}^{h_{0r}},g_{0r})\big{)}+ \\
                 & \quad +\sum_{j=2}^{r}(-1)^j\rho_0^0(t_{p}(\partial_0\gamma_0)...t_{p}(\partial_0(\gamma_{j-1}\vJoin\gamma_j))...t_{p}(\partial_0\gamma_{q+r}))\phi\big{(}\omega(\delta_0^{r-1}(\delta_j\vec{\gamma})_{0,0};...,(g_{0j-1}^{h_{0j}}g_{0j})^{h_{0j+1}...h_{0r}},...)\big{)}+ \\
                 & \quad +\sum_{j=r+1}^{r+q}(-1)^j\rho_0^0(t_{p}(\partial_0\gamma_0)...t_{p}(\partial_0(\gamma_{j-1}\vJoin\gamma_j))...t_{p}(\partial_0\gamma_{q+r}))\phi\big{(}\omega(\delta_0^{r-1}(\delta_j\vec{\gamma})_{0,0};g_{00}^{h_{01}...h_{0r-1}},...,g_{0r-2}^{h_{r-1}},g_{0r-1})\big{)}+ \\
                 & \quad\qquad +(-1)^{r+q+1}\rho_0^0(t_{p}(\partial_0\gamma_0)...t_{p}(\partial_0\gamma_{q+r-1}))\phi\big{(}\omega(\delta_0^{r-1}(\delta_{r+q+1}\vec{\gamma})_{0,0};g_{00}^{h_{01}...h_{0r-1}},...,g_{0r-2}^{h_{r-1}},g_{0r-1})\big{)} \\
                 & =\rho_0^0(h_{00}...h_{0q+r})\phi\Big{(}\rho_0^1(i(g_{00}^{h_{01}...h_{0q+r}}))\omega(\delta_0^r\vec{\gamma}_{0,0};g_{01}^{h_{02}...h_{0r}},...,g_{0r-1}^{h_{0r}},g_{0r}) \\
                 & \qquad +\sum_{k=1}^r(-1)^k\omega(\delta_0^r\vec{\gamma}_{0,0};\delta_k(g_{00}^{h_{01}...h_{0r}},...,g_{0r-1}^{h_{0r}},g_{0r}))+ \\
                 & \qquad\quad +\sum_{j=r+1}^{r+q}(-1)^j\omega(\delta_0^{r-1}(\delta_j\vec{\gamma})_{0,0};g_{00}^{h_{01}...h_{0r-1}},...,g_{0r-2}^{h_{0r-1}},g_{0r-1})) \\
                 & \qquad\qquad +(-1)^{r+q+1}\rho_0^1(h_{0q+r})^{-1}\omega(\delta_0^{r-1}(\delta_{r+q+1}\vec{\gamma})_{0,0};g_{00}^{h_{01}...h_{0r-1}},...,g_{0r-2}^{h_{r-1}},g_{0r-1})\Big{)}
\end{align*}
On the other hand,
\begin{align*}
\Delta_r^p\delta\omega(\vec{\gamma}) & =\rho_0^0(t_{p}(\partial_0\gamma_0)...t_{p}(\partial_0\gamma_{q+r}))\phi\big{(}\delta\omega(\delta_0^{r-1}\vec{\gamma}_{0,0};g_{00}^{h_{01}...h_{0r-1}},...,g_{0r-2}^{h_{0r-1}},g_{0r-1})\big{)} \\
                                 & =\rho_0^0(h_{00}...h_{0q+r})\phi\Big{(}\omega(\delta_0^r\vec{\gamma}_{0,0};(g_{00}^{h_{01}...h_{0r-1}})^{h_{0r}},...,(g_{0r-2}^{h_{0r-1}})^{h_{0r}},g_{0r-1}^{h_{0r}}))+ \\
                                 & \qquad\quad +\sum_{j=1}^{q}(-1)^j\omega(\delta_j\delta_0^{r-1}\vec{\gamma}_{0,0};g_{00}^{h_{01}...h_{0r-1}},...,g_{0r-2}^{h_{r-1}},g_{0r-1})+ \\
                                 & \qquad\qquad +(-1)^{q+1}\rho_0^1(t_p(\partial_0\gamma_{q+r}))^{-1}\omega(\delta_{q+1}\delta_0^{r-1}\vec{\gamma}_{0,0};g_{00}^{h_{01}...h_{0r-1}},...,g_{0r-2}^{h_{r-1}},g_{0r-1})\Big{)}.
\end{align*}
Now, $\delta_0^{r-1}(\delta_j\vec{\gamma})_{0,0}=\delta_{j-r}\delta_0^{r-1}\vec{\gamma}_{0,0}$ for $j\geq 1$ and $t_p(\partial_0\gamma_{q+r})=h_{0q+r}$; hence, 
\begin{align*}
(\delta\Delta_r^p+(-1)^{r+1}\Delta_r^p\delta)\omega(\vec{\gamma}) & =\rho_0^0(h_{00}...h_{0q+r})\phi\big{(}\delta_{(1)}\omega(\delta_0^r\vec{\gamma}_{0,0};g_{00}^{h_{01}...h_{0r}},...,g_{0r-1}^{h_{0r}},g_{0r}))\big{)} =\Delta_{r+1}^p\delta_{(1)}\omega(\vec{\gamma})
\end{align*}
and the result follows.

\end{proof}
There is still a piece of data missing. The previous proposition hints that, modulo some appropriate signs, the sum of the differentials of the triple complex together and the difference maps squares to zero when it gets to the front page; however, there are still a number of relations that are introduced in the process. We will exemplify this briefly. Let $\omega\in C(G^2,W)$, then the following diagram represents its differential
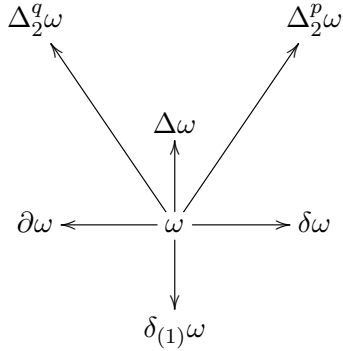
\begin{eqnarray*}
\xymatrix{
\Delta_2^q\omega & & \Delta_2^p\omega \\
         & \Delta\omega  &            \\
\partial\omega & \omega \ar[d]\ar[l]\ar[r]\ar[u]\ar[uul]\ar[uur] & \delta\omega \\   
       & \delta_{(1)}\omega &         
}
\end{eqnarray*}
and, using pointed arrows, we represent the second differential by
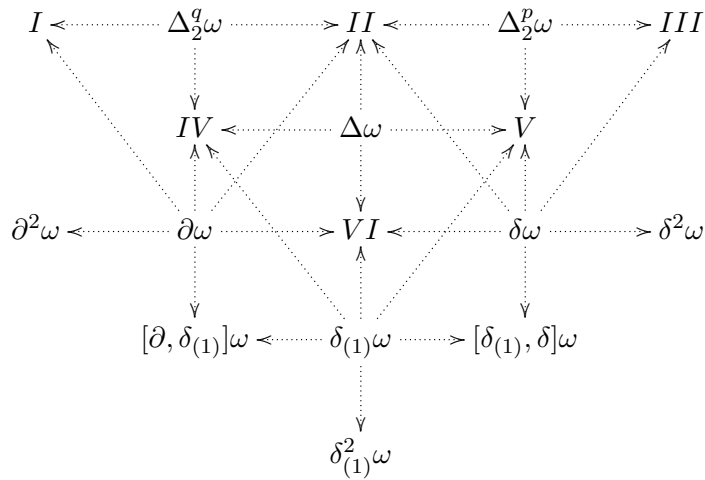
\begin{eqnarray*}
\xymatrix{
I & \Delta_2^q\omega \ar@{.>}[d]\ar@{.>}[l]\ar@{.>}[r] & II & \Delta_2^p\omega \ar@{.>}[d]\ar@{.>}[l]\ar@{.>}[r] & III \\
 &    IV   & \Delta\omega \ar@{.>}[d]\ar@{.>}[l]\ar@{.>}[r]\ar@{.>}[u] &    V     & \\
\partial^2\omega & \partial\omega \ar@{.>}[d]\ar@{.>}[l]\ar@{.>}[r]\ar@{.>}[u]\ar@{.>}[uul]\ar@{.>}[uur] & VI & \delta\omega \ar@{.>}[d]\ar@{.>}[l]\ar@{.>}[r]\ar@{.>}[u]\ar@{.>}[uul]\ar@{.>}[uur] & \delta^2\omega \\   
 & [\partial,\delta_{(1)}]\omega & \delta_{(1)}\omega \ar@{.>}[d]\ar@{.>}[l]\ar@{.>}[r]\ar@{.>}[u]\ar@{.>}[uul]\ar@{.>}[uur] & [\delta_{(1)},\delta]\omega & \\
 &        &                               \delta_{(1)}^2\omega                                         &       &
}
\end{eqnarray*}
We tailored the second difference maps to fill in $IV$ (cf. lemma \ref{IV atSch}) and $V$ (cf. lemma \ref{V atSch}). Also, the full generality of proposition \ref{FrontPageCorners} ensures the vanishing of $I$ and $III$. From this, it is clear that the only relation left to vanish is $II$. This relation involves 
\begin{eqnarray}\label{II atSch}
\delta\Delta_2^q+\Delta_2^p\partial-\Delta\circ\Delta+\Delta_2^q\delta-\partial\Delta_2^p.
\end{eqnarray}
As it turns out, there is yet another difference map
\begin{eqnarray*}
\xymatrix{
\Delta_{2,2}:C(G^3,W) \ar[r] & C(\G^2_2,V)
}
\end{eqnarray*}
landing on the front page, defined by
\begin{align*}
\Delta_{2,2}\omega(\gamma_0,\gamma_1) & =\rho_0^0(h_0h_1)\phi\Big{(}\omega(g_{10}^{h_1},g_{00}^{h_1},g_{11})+\omega((g_{10}g_{00})^{h_1},g_{11},g_{01})-\omega(g_{10}^{h_1}g_{11},g_{00}^{h_1i(g_{11})},g_{01})+ \\
 & \qquad\qquad\qquad\qquad -\omega(g_{10}^{h_1}g_{11},g_{11}^{-1},g_{00}^{h_1}g_{11})+\omega(g_{10}^{h_1}g_{11},g_{11}^{-1},g_{11})\Big{)}    
\end{align*}
making the linear combination \ref{II atSch} equal to $\Delta_{2,2}\delta_{(1)}$. This phenomenon will occur ubiquitously and the fact that such an expression involves $\Delta\circ\Delta$ precludes the involvement of (higher) difference maps not landing on the front page; therefore, in order to compute the values of the missing difference maps, we need to see first how the difference maps look like when they land away from the front page. 

\subsubsection{Higher difference maps 2: corner maps away from the front page}
For this part, we will compute the formulas for $\Delta_2^q$ and $\Delta_2^p$ that land on the page $r=1$ by looking at the square of the differential applied to an element $\omega\in C(\G_p^q\times G^3,W)$. Immediately thereafter, we successively generalize and get formulas for $\Delta_k^q$ and $\Delta_k^p$ landing on the page $r=1$ and the formulas for arbitrary $r$. \\
First, let $f\in G$ and $\vec{\gamma}=(\gamma_0,...,\gamma_q)\in\G_{p+2}^{q+1}$ with the usual identifications $\gamma_b=(\gamma_{0b},...,\gamma_{p+1b})^T$ and $\gamma_{ab}\sim(g_{ab},h_{ab})$. We compute
\begin{align*}
\partial\Delta\omega(\vec{\gamma};f) & =\rho_0^1(i(g_{00}^{h_{01}...h_{0q}}...g_{0q}))^{-1}\Delta\omega(\partial_0\vec{\gamma};f)+\sum_{j=1}^{p+2}(-1)^j\Delta\omega(\partial_j\vec{\gamma};f) \\
                                     & =\rho_0^1(i(g_{00}^{h_{01}...h_{0q}}...g_{0q}))^{-1}\rho_0^1(i(g_{11}^{h_{12}...h_{1q}}...g_{1q}))^{-1}\Big{(}\rho_0^1(i(g_{10}^{h_{11}...h_{1q}}))^{-1}\omega(\partial_0\vec{\gamma}_{0,0};f^{h_{10}},g_{10})+  \\
                                     & \qquad\qquad\quad +\omega(\partial_0\vec{\gamma}_{0,0};g_{10}^{-1},f^{h_{10}}g_{10})-\omega(\partial_0\vec{\gamma}_{0,0};g_{10}^{-1},g_{10})\Big{)}+ \\
                                     & \qquad -\rho_0^1(i((g_{11}g_{01})^{h_{12}...h_{1q}}...(g_{1q}g_{0q})))^{-1}\Big{(}\rho_0^1(i((g_{10}g_{00})^{h_{11}...h_{1q}}))^{-1}\omega((\partial_1\vec{\gamma})_{0,0};f^{h_{10}},g_{10}g_{00})+  \\
                                     & \qquad\qquad\quad +\omega((\partial_1\vec{\gamma})_{0,0};(g_{10}g_{00})^{-1},f^{h_{10}}g_{10}g_{00})-\omega((\partial_1\vec{\gamma})_{0,0};(g_{10}g_{00})^{-1},g_{10}g_{00})\Big{)}+ \\
                                     & \qquad +\sum_{j=2}^{p+2}(-1)^j\rho_0^1(i(g_{01}^{h_{02}...h_{0q}}...g_{0q}))^{-1}\Big{(}\rho_0^1(i(g_{00}^{h_{01}...h_{0q}}))^{-1}\omega((\partial_j\vec{\gamma})_{0,0};f^{h_{00}},g_{00})+  \\
                                     & \qquad\qquad\quad +\omega((\partial_j\vec{\gamma})_{0,0};g_{00}^{-1},f^{h_{00}}g_{00})-\omega((\partial_j\vec{\gamma})_{0,0};g_{00}^{-1},g_{00})\Big{)}
\end{align*}
and 
\begin{align*}
\Delta\partial\omega(\vec{\gamma};f) & =\rho_0^1(i(g_{01}^{h_{02}...h_{0q}}...g_{0q}))^{-1}\Big{(}\rho_0^1(i(g_{00}^{h_{01}...h_{0q}}))^{-1}\partial\omega(\vec{\gamma}_{0,0};f^{h_{00}},g_{00})+ \\
                                     & \qquad\qquad\quad +\partial\omega(\vec{\gamma}_{0,0};g_{00}^{-1},f^{h_{00}}g_{00})-\partial\omega(\vec{\gamma}_{0,0};g_{00}^{-1},g_{00})\Big{)} \\
                                     & =\rho_0^1(i(g_{01}^{h_{02}...h_{0q}}...g_{0q}))^{-1}\rho_0^1(i(g_{00}^{h_{01}...h_{0q}}))^{-1}\Big{(}\rho_0^1(i(g_{11}^{h_{12}...h_{1q}}...g_{1q}))^{-1}\omega(\partial_0\vec{\gamma}_{0,0};f^{h_{00}},g_{00})+ \\
                                     & \qquad\qquad\quad +\sum_{j=1}^{p+1}(-1)^j\omega(\partial_j\vec{\gamma}_{0,0};f^{h_{00}},g_{00}) \\
                                     & \qquad +\rho_0^1(i(g_{01}^{h_{02}...h_{0q}}...g_{0q}))^{-1}\Big{(}\rho_0^1(i(g_{11}^{h_{12}...h_{1q}}...g_{1q}))^{-1}\omega(\partial_0\vec{\gamma}_{0,0};g_{00}^{-1},f^{h_{00}}g_{00})+ \\
                                     & \qquad\qquad +\sum_{j=1}^{p+1}(-1)^j\omega(\partial_j\vec{\gamma}_{0,0};g_{00}^{-1},f^{h_{00}}g_{00})-\rho_0^1(i(g_{11}^{h_{12}...h_{1q}}...g_{1q}))^{-1}\omega(\partial_0\vec{\gamma}_{0,0};g_{00}^{-1},g_{00})+ \\
                                     & \qquad\qquad\quad +\sum_{j=1}^{p+1}(-1)^j\omega(\partial_j\vec{\gamma}_{0,0};g_{00}^{-1},g_{00})\Big{)}.
\end{align*}
Using the simplicial identity $(\partial_j\vec{\gamma})_{0,0}=\partial_{j-1}\vec{\gamma}_{0,0}$, we see that by adding these terms together, the only terms surviving are the ones with first entry $\partial_0\gamma_{0,0}$ (of course, this includes $(\partial_1\vec{\gamma})_{0,0}$). There is a third path to get from $C(\G_p^q\times G^2,W)$ to $C(\G_{p+2}^{q+1}\times G,W)$; namely,
\begin{align*}
\delta_{(1)}\Delta_2^q\omega(\vec{\gamma};f) & =\rho_0^1(t_{p+2}\gamma_0...t_{p+2}\gamma_q)^{-1}\rho_1(f)\Delta_2^q\omega(\vec{\gamma}) \\
                                             & =\rho_0^1(h_{00}i(g_{00})...h_{0q}i(g_{0q}))^{-1}\rho_1(f)\rho_0^0(t_p(\partial_0^2\gamma_0)...t_p(\partial_0^2\gamma_q))\phi(\omega(\partial_0\vec{\gamma}_{0,0};g_{10},g_{00})).
\end{align*}
We notice that, on the one hand, $t_p(\partial_0^2\gamma_b)=t(\gamma_{2b})=h_{2b}i(g_{2b})=h_{1b}$; while on the other, $h_{0b}=h_{1b}i(g_{1b})$, implying 
\begin{align*}
\rho_0^1(h_{00}i(g_{00})...h_{0q}i(g_{0q}))^{-1} & =\rho_0^1(h_{00}...h_{0q}i(g_{00}^{h_{01}...h_{0q}}...g_{0q}))^{-1} \\
                                                 & =\rho_0^1(i(g_{00}^{h_{01}...h_{0q}}...g_{0q}))^{-1}\rho_0^1(h_{00}...h_{0q})^{-1} \\
                                                 & =\rho_0^1(i(g_{00}^{h_{01}...h_{0q}}...g_{0q}))^{-1}\rho_0^1(h_{10}i(g_{10})...h_{1q}i(g_{1q}))^{-1} \\
                                                 & =\rho_0^1(i(g_{00}^{h_{01}...h_{0q}}...g_{0q}))^{-1}\rho_0^1(i(g_{10}^{h_{11}...h_{1q}}...g_{1q}))^{-1}\rho_0^1(h_{10}...h_{1q})^{-1} 
\end{align*}
and
\begin{align*}
\delta_{(1)}\Delta_2^q\omega(\vec{\gamma};f) & =\rho_0^1(i(g_{00}^{h_{01}...h_{0q}}...g_{0q}))^{-1}\rho_0^1(i(g_{10}^{h_{11}...h_{1q}}...g_{1q}))^{-1}\rho_1(f^{h_{10}...h_{1q}})\phi(\omega(\partial_0\vec{\gamma}_{0,0};g_{10},g_{00})) \\
                                             & =\rho_0^1(i(g_{00}^{h_{01}...h_{0q}}...g_{0q}))^{-1}\rho_0^1(i(g_{10}^{h_{11}...h_{1q}}...g_{1q}))^{-1}\Big{[}\rho_0^1(i(f^{h_{10}...h_{1q}}))-I\Big{]}\omega(\partial_0\vec{\gamma}_{0,0};g_{10},g_{00}).
\end{align*}
If one considers the difference $\delta_{(1)}\Delta_2^q-(\partial\Delta+\Delta\partial)$, all terms come in the form
\begin{eqnarray*}
\rho_0^1(i(g_{01}^{h_{02}...h_{0q}}...g_{0q}))^{-1}\rho_0^1(i(g_{11}^{h_{12}...h_{1q}}...g_{1q}))^{-1}\rho_0^1(i(\bullet^{h_{11}...h_{0q}}))\omega(\partial_0\vec{\gamma}_{0,0};\bullet ,\bullet ).
\end{eqnarray*}
We arrange the surviving terms, writing down just what the bullets in this expression are filled with. We use the terms that have a first bullet different from $1$ as top seed:
\begin{eqnarray*}
\begin{pmatrix}
    \textnormal{Top Seed} &  &  & \\
  \rho((g_{10}g_{00})^{-1}f^{h_{10}})(g_{10},g_{00}) & & & \\
   -\rho(g_{10}g_{00})^{-1}(g_{10},g_{00}) & +(g_{00}^{-1},g_{00}) & -((g_{10}g_{00})^{-1},g_{10}g_{00}) & \\
   -\rho(g_{00})^{-1}(f^{h_{00}},g_{00}) & & -(g_{00}^{-1},f^{h_{00}}g_{00}) & \\
   +\rho(g_{10}g_{00})^{-1}(f^{h_{10}},g_{10}g_{00}) & & +((g_{10}g_{00})^{-1},f^{h_{10}}g_{10}g_{00}) & \\
   -\rho(g_{10}g_{00})^{-1}(f^{h_{10}},g_{10}) & & & \\
   -\rho(g_{00})^{-1}(g_{10}^{-1},f^{h_{10}}g_{10}) & & & \\
   +\rho(g_{00})^{-1}(g_{10}^{-1},g_{10}) & & & 
\end{pmatrix}
\end{eqnarray*}
The strategy now is to complete cocycles; to do so, we are going to add and subtract what the top seed are missing. The remarkable fact is that this is even possible; indeed,
\begin{eqnarray*}
\begin{pmatrix}
    \textnormal{Top Seed} & (ab,c) & (a,bc) & (a,b) \\
    \rho((g_{10}g_{00})^{-1}f^{h_{10}})(g_{10},g_{00}) & -\epsilon_1 & +\epsilon_2 & -\epsilon_3 \\
   -\rho(g_{10}g_{00})^{-1}(g_{10},g_{00}) & +(g_{00}^{-1},g_{00}) & -((g_{10}g_{00})^{-1},g_{10}g_{00}) & +\epsilon_4 \\
   -\rho(g_{00})^{-1}(f^{h_{00}},g_{00}) & +\epsilon_1 & -(g_{00}^{-1},f^{h_{00}}g_{00}) & +\epsilon_5 \\
   +\rho(g_{10}g_{00})^{-1}(f^{h_{10}},g_{10}g_{00}) & -\epsilon_2 & +((g_{10}g_{00})^{-1},f^{h_{10}}g_{10}g_{00}) & -\epsilon_6 \\
   -\rho(g_{10}g_{00})^{-1}(f^{h_{10}},g_{10}) & +\epsilon_3 & -\epsilon_7 & +\epsilon_6 \\
   -\rho(g_{00})^{-1}(g_{10}^{-1},f^{h_{10}}g_{10}) & +\epsilon_7 & -\epsilon_5 & +\epsilon_8 \\
   +\rho(g_{00})^{-1}(g_{10}^{-1},g_{10}) & -\epsilon_4 & \circ & -\epsilon_8 
\end{pmatrix}
\end{eqnarray*}
where in the slot of $\circ$ there is a zero as $\omega(\partial_0\vec{\gamma}_{0,0};g_{00}^{-1},1)\equiv 0$. Using the cocycle equation, we make the replacement 
\begin{align*}
\alpha((g_{10}g_{00})^{-1}f^{h_{10}},g_{10},g_{00})+ & \alpha((g_{10}g_{00})^{-1},f^{h_{10}},g_{10}g_{00})-\alpha((g_{10}g_{00})^{-1},f^{h_{10}},g_{10}) \\
 & =\rho(g_{10}g_{00})^{-1}\alpha(f^{h_{10}},g_{10},g_{00})+\alpha((g_{10}g_{00})^{-1},f^{h_{10}}g_{10},g_{00})
\end{align*}
where we used the shorthand $\alpha(a,b,c)=\delta_{(1)}\omega(\partial_0\vec{\gamma}_{0,0};a,b,c)$ and therefore $\rho(a)$ actually stands for $\rho_0^1(i(a^{h_{11}...h_{1q}}))$. We summarize this discussion in the following lemma.
\begin{lemma}
Let
\begin{eqnarray*}
\xymatrix{
\Delta_2^q:C(\G_p^q\times G^2,W) \ar[r] & C(\G_{p+2}^{q+1}\times G,W)
}
\end{eqnarray*}
be defined by
\begin{align*}
\Delta_2^q\omega(\vec{\gamma};f) & :=\rho_0^1(i(pr_G((\gamma_{01}\Join\gamma_{11})\vJoin ...\vJoin(\gamma_{0q}\Join\gamma_{1q}))))^{-1}\Big{(}\rho_0^1(i((g_{10}g_{00})^{h_{11}...h_{1q}}))^{-1}\omega(\partial_0\vec{\gamma}_{0,0};f^{h_{10}},g_{10},g_{00})+ \\
    & \qquad -\omega(\partial_0\vec{\gamma}_{0,0};g_{00}^{-1},f^{h_{00}},g_{00})+ \\
    & \qquad\quad +\omega(\partial_0\vec{\gamma}_{0,0};(g_{10}g_{00})^{-1},f^{h_{10}}g_{10},g_{00})-\omega(\partial_0\vec{\gamma}_{0,0};(g_{10}g_{00})^{-1},g_{10},g_{00})+ \\ 
    & \qquad\quad -\omega(\partial_0\vec{\gamma}_{0,0};g_{00}^{-1},g_{10}^{-1},f^{h_{10}}g_{10})+\omega(\partial_0\vec{\gamma}_{0,0};g_{00}^{-1},g_{10}^{-1},g_{10})\Big{)}   
\end{align*}
for $\vec{\gamma}=(\gamma_0,...,\gamma_q)\in\G_{p+2}^{q+1}$, $\gamma_b=(\gamma_{0b},...,\gamma_{p+1b})^T$ and $\gamma_{ab}\sim(g_{ab},h_{ab})$. Then
\begin{eqnarray*}
\partial\circ\Delta+\Delta\circ\partial=\delta_{(1)}\circ\Delta_2^q-\Delta_2^q\circ\delta_{(1)}.
\end{eqnarray*}
\end{lemma}
We now turn to the map $\Delta_2^p$. For this purpose, let $f\in G$ and $\vec{\gamma}=(\gamma_0,...,\gamma_{q+1})\in\G_{p+1}^{q+2}$ with the usual identifications $\gamma_b=(\gamma_{0b},...,\gamma_{pb})^T$ and $\gamma_{ab}\sim(g_{ab},h_{ab})$. We compute
\begin{align*}
\delta\Delta\omega(\vec{\gamma};f) & =\Delta\omega(\delta_0\vec{\gamma};f^{h_{00}i(g_{00})})+\sum_{j=1}^{q+1}(-1)^j\Delta\omega(\delta_j\vec{\gamma};f)+(-1)^{q+2}\rho_0^1(h_{0q+1}i(g_{0q+1}))^{-1}\Delta\omega(\delta_{q+2}\vec{\gamma};f) \\
                                   & =\rho_0^1(i(g_{02}^{h_{03}...h_{0q+1}}...g_{0q+1}))^{-1}\Big{(}\rho_0^1(i(g_{01}^{h_{02}...h_{0q+1}}))^{-1}\omega(\delta_0\vec{\gamma}_{0,0};f^{h_{00}i(g_{00})h_{01}},g_{01})+  \\
                                   & \qquad\qquad\quad +\omega(\delta_0\vec{\gamma}_{0,0};g_{01}^{-1},f^{h_{00}i(g_{00})h_{01}}g_{01})-\omega(\delta_0\vec{\gamma}_{0,0};g_{01}^{-1},g_{01})\Big{)}+ \\
                                   & \qquad -\rho_0^1(i(g_{02}^{h_{03}...h_{0q+1}}...g_{0q+1}))^{-1}\Big{(}\rho_0^1(i((g_{00}^{h_{01}}g_{01})^{h_{02}...h_{0q+1}}))^{-1}\omega((\delta_1\vec{\gamma})_{0,0};f^{h_{00}h_{01}},g_{00}^{h_{01}}g_{01})+  \\
                                   & \qquad\qquad\quad +\omega((\delta_1\vec{\gamma})_{0,0};(g_{00}^{h_{01}}g_{01})^{-1},f^{h_{00}h_{01}}g_{00}^{h_{01}}g_{01})-\omega((\delta_1\vec{\gamma})_{0,0};(g_{00}^{h_{01}}g_{01})^{-1},g_{00}^{h_{01}}g_{01})\Big{)}+ \\
                                   & \qquad +\sum_{j=2}^{q+1}(-1)^j\rho_0^1(i(g_{01}^{h_{02}...h_{0q+1}}...g_{0q+1}))^{-1}\Big{(}\rho_0^1(i(g_{00}^{h_{01}...h_{0q+1}}))^{-1}\omega((\delta_j\vec{\gamma})_{0,0};f^{h_{00}},g_{00})+  \\
                                   & \qquad\qquad\quad +\omega((\delta_j\vec{\gamma})_{0,0};(g_{00})^{-1},f^{h_{00}}g_{00})-\omega((\delta_j\vec{\gamma})_{0,0};(g_{00})^{-1},g_{00})\Big{)}+ \\
                                   & \qquad +(-1)^{q+2}\rho_0^1(h_{0q+1}i(g_{0q+1}))^{-1}\rho_0^1(i(g_{01}^{h_{02}...h_{0q}}...g_{0q}))^{-1}\Big{(}\rho_0^1(i(g_{00}^{h_{01}...h_{0q}}))^{-1}\omega((\delta_{q+2}\vec{\gamma})_{0,0};f^{h_{00}},g_{00})+  \\
                                   & \qquad\qquad\quad +\omega((\delta_{q+2}\vec{\gamma})_{0,0};(g_{00})^{-1},f^{h_{00}}g_{00})-\omega((\delta_{q+2}\vec{\gamma})_{0,0};(g_{00})^{-1},g_{00})\Big{)}+ \\
\end{align*}
and 
\begin{align*}
\Delta\delta\omega(\vec{\gamma};f) & =\rho_0^1(i(g_{01}^{h_{02}...h_{0q+1}}...g_{0q+1}))^{-1}\Big{(}\rho_0^1(i(g_{00}^{h_{01}...h_{0q+1}}))^{-1}\delta\omega(\vec{\gamma}_{0,0};f^{h_{00}},g_{00})+ \\
                                   & \qquad\qquad\quad +\delta\omega(\vec{\gamma}_{0,0};g_{00}^{-1},f^{h_{00}}g_{00})-\delta\omega(\vec{\gamma}_{0,0};g_{00}^{-1},g_{00})\Big{)} \\
                                   & =\rho_0^1(i(g_{00}^{h_{01}...h_{0q+1}}...g_{0q+1}))^{-1}\Big{(}\omega(\delta_0\vec{\gamma}_{0,0};f^{h_{00}h_{01}},g_{00}^{h_{01}})+ \\
                                   & \qquad\qquad\quad +\sum_{j=1}^{q}(-1)^j\omega(\delta_j\vec{\gamma}_{0,0};f^{h_{00}},g_{00})+(-1)^{q+1}\rho_0^1(h_{0q+1})^{-1}\omega(\delta_{q+1}\vec{\gamma}_{0,0};f^{h_{00}},g_{00})\Big{)} \\
                                   & \qquad +\rho_0^1(i(g_{01}^{h_{02}...h_{0q+1}}...g_{0q+1}))^{-1}\Big{(}\omega(\delta_0\vec{\gamma}_{0,0};(g_{00}^{h_{01}})^{-1},f^{h_{00}h_{01}}g_{00}^{h_{01}})+ \\
                                   & \qquad\qquad +\sum_{j=1}^{q}(-1)^j\omega(\delta_j\vec{\gamma}_{0,0};g_{00}^{-1},f^{h_{00}}g_{00})+(-1)^{q+1}\rho_0^1(h_{0q+1})^{-1}\omega(\delta_{q+1}\vec{\gamma}_{0,0};g_{00}^{-1},f^{h_{00}}g_{00}) + \\
                                   & \qquad -\omega(\delta_0\vec{\gamma}_{0,0};(g_{00}^{h_{01}})^{-1},g_{00}^{h_{01}})+ \\
                                   & \qquad\qquad +\sum_{j=1}^{q}(-1)^{j+1}\omega(\delta_j\vec{\gamma}_{0,0};g_{00}^{-1},g_{00})+(-1)^{q}\rho_0^1(h_{0q+1})^{-1}\omega(\delta_{q+1}\vec{\gamma}_{0,0};g_{00}^{-1},g_{00})\Big{)}.
\end{align*}
Using the simplicial identity $(\delta_j\vec{\gamma})_{0,0}=\delta_{j-1}\vec{\gamma}_{0,0}$, we see that by adding these terms together, the only terms surviving are the ones with first entry $\delta_0\gamma_{0,0}$ (of course, this includes $(\delta_1\vec{\gamma})_{0,0}$). There is a third path to get from $C(\G_p^q\times G^2,W)$ to $C(\G_{p+1}^{q+2}\times G,W)$; namely,
\begin{align*}
\delta_{(1)}\Delta_2^p\omega(\vec{\gamma};f) & =\rho_0^1(t_{p+1}\gamma_0...t_{p+1}\gamma_{q+1})^{-1}\rho_1(f)\Delta_2^p\omega(\vec{\gamma}) \\
                                             & =\rho_0^1(h_{00}i(g_{00})...h_{0q+1}i(g_{0q+1}))^{-1}\rho_1(f)\rho_0^0(t_p(\partial_0\gamma_0)...t_p(\partial_0\gamma_{q+1}))\phi(\omega(\delta_0\vec{\gamma}_{0,0};g_{00}^{01},g_{01})).
\end{align*}
Again, since on the one hand, $t_p(\partial_0\gamma_b)=t(\gamma_{1b})=h_{1b}i(g_{1b})=h_{0b}$ and on the other,  
\begin{align*}
\rho_0^1(h_{00}i(g_{00})...h_{0q+1}i(g_{0q+1}))^{-1} & =\rho_0^1(h_{00}...h_{0q+1}i(g_{00}^{h_{01}...h_{0q+1}}...g_{0q+1}))^{-1} \\
                                                 & =\rho_0^1(i(g_{00}^{h_{01}...h_{0q+1}}...g_{0q+1}))^{-1}\rho_0^1(h_{00}...h_{0q+1})^{-1},
\end{align*}
we have got
\begin{align*}
\delta_{(1)}\Delta_2^p\omega(\vec{\gamma};f) & =\rho_0^1(i(g_{00}^{h_{01}...h_{0q+1}}...g_{0q+1}))^{-1}\rho_1(f^{h_{00}...h_{0q+1}})\phi(\omega(\delta_0\vec{\gamma}_{0,0};g_{00}^{h_{01}},g_{01})) \\
                                             & =\rho_0^1(i(g_{00}^{h_{01}...h_{0q+1}}...g_{0q+1}))^{-1}\Big{[}\rho_0^1(i(f^{h_{00}...h_{0q+1}}))-I\Big{]}\omega(\delta_0\vec{\gamma}_{0,0};g_{00}^{h_{01}},g_{01}).
\end{align*}
If one considers the difference $\delta_{(1)}\Delta_2^p-(\delta\Delta+\Delta\delta)$, all terms come in the form
\begin{eqnarray*}
\rho_0^1(i(g_{02}^{h_{03}...h_{0q+1}}...g_{0q+1}))^{-1}\rho_0^1(i(\bullet^{h_{02}...h_{0q+1}}))\omega(\delta_0\vec{\gamma}_{0,0};\bullet ,\bullet ).
\end{eqnarray*}
We arrange the surviving terms as before, writing down just what the bullets in this expression are filled with. We use the terms that have a first bullet different from $1$ as top seed:
\begin{eqnarray*}
\begin{pmatrix}
    \textnormal{Top Seed} & & & \\
    \rho((g_{00}^{h_{01}}g_{01})^{-1}f^{h_{00}h_{10}})(g_{00}^{h_{01}},g_{00}) & & & \\
   -\rho(g_{00}^{h_{01}}g_{01})^{-1}(g_{00}^{h_{01}},g_{00}) & +(g_{01}^{-1},g_{01}) & -((g_{00}^{h_{01}}g_{01})^{-1},g_{00}^{h_{01}}g_{01}) & \\
   -\rho(g_{01})^{-1}(f^{h_{00}i(g_{00})h_{01}},g_{01}) & & -(g_{01}^{-1},f^{h_{00}i(g_{00})h_{01}}g_{01}) & \\
   +\rho(g_{00}^{h_{01}}g_{01})^{-1}(f^{h_{00}h_{01}},g_{00}^{h_{01}}g_{01}) & & +((g_{00}^{h_{01}}g_{01})^{-1},f^{h_{00}h_{01}}g_{00}^{h_{01}}g_{01}) & \\
   -\rho(g_{00}^{h_{01}}g_{01})^{-1}(f^{h_{00}h_{01}},g_{00}^{h_{01}}) & & & \\
   -\rho(g_{01})^{-1}((g_{00}^{h_{10}})^{-1},f^{h_{00}h_{01}}g_{00}^{h_{01}}) & & & \\
   +\rho(g_{01})^{-1}((g_{00}^{h_{10}})^{-1},g_{00}^{h_{01}}) & & & 
\end{pmatrix}
\end{eqnarray*}
Amazingly, we will be able to complete each cocycle again with nothing being left out:
\begin{eqnarray*}
\begin{pmatrix}
    \textnormal{Top Seed} & (ab,c) & (a,bc) & (a,b) \\
    \rho((g_{00}^{h_{01}}g_{01})^{-1}f^{h_{00}h_{10}})(g_{00}^{h_{01}},g_{01}) & -\epsilon_1 & +\epsilon_2 & -\epsilon_3 \\
   -\rho(g_{00}^{h_{01}}g_{01})^{-1}(g_{00}^{h_{01}},g_{00}) & +(g_{01}^{-1},g_{01}) & -((g_{00}^{h_{01}}g_{01})^{-1},g_{00}^{h_{01}}g_{01}) & +\epsilon_4 \\
   -\rho(g_{01})^{-1}(f^{h_{00}i(g_{00})h_{01}},g_{01}) & +\epsilon_1 & -(g_{01}^{-1},f^{h_{00}i(g_{00})h_{01}}g_{01}) & +\epsilon_5 \\
   +\rho(g_{00}^{h_{01}}g_{01})^{-1}(f^{h_{00}h_{01}},g_{00}^{h_{01}}g_{01}) & -\epsilon_2 & +((g_{00}^{h_{01}}g_{01})^{-1},f^{h_{00}h_{01}}g_{00}^{h_{01}}g_{01}) & -\epsilon_6 \\
   -\rho(g_{00}^{h_{01}}g_{01})^{-1}(f^{h_{00}h_{01}},g_{00}^{h_{01}}) & +\epsilon_3 & -\epsilon_7 & +\epsilon_6 \\
   -\rho(g_{01})^{-1}((g_{00}^{h_{10}})^{-1},f^{h_{00}h_{01}}g_{00}^{h_{01}}) & +\epsilon_7 & -\epsilon_5 & +\epsilon_8 \\
   +\rho(g_{01})^{-1}((g_{00}^{h_{10}})^{-1},g_{00}^{h_{01}}) & -\epsilon_4 & \circ & -\epsilon_8
\end{pmatrix}
\end{eqnarray*}
where there is a zero in the slot of $\circ$, as $\omega(\delta_0\vec{\gamma}_{0,0};(g_{01})^{-1},1)\equiv 0$. Using the cocycle equation, we make the replacement 
\begin{align*}
\alpha((g_{00}^{h_{01}}g_{01})^{-1}f^{h_{00}h_{01}},g_{00}^{h_{01}},g_{01})+ & \alpha((g_{00}^{h_{01}}g_{01})^{-1},f^{h_{00}h_{01}},g_{00}^{h_{01}}g_{01})-\alpha((g_{00}^{h_{01}}g_{01})^{-1},f^{h_{00}h_{01}},g_{00}^{h_{01}}) \\
 & =\rho(g_{00}^{h_{01}}g_{01})^{-1}\alpha(f^{h_{00}h_{01}},g_{00}^{h_{01}},g_{01})+\alpha((g_{00}^{h_{01}}g_{01})^{-1},f^{h_{00}h_{01}}g_{00}^{h_{01}},g_{01})
\end{align*}
where we used the shorthand $\alpha(a,b,c)=\delta_{(1)}\omega(\delta_0\vec{\gamma}_{0,0};a,b,c)$ and therefore $\rho(a)$ actually stands for $\rho_0^1(i(a^{h_{02}...h_{0q+1}}))$. We summarize this discussion in the following lemma.
\begin{lemma}
Let
\begin{eqnarray*}
\xymatrix{
\Delta_2^p:C(\G_p^q\times G^2,W) \ar[r] & C(\G_{p+1}^{q+2}\times G,W)
}
\end{eqnarray*}
be defined by
\begin{align*}
\Delta_2^p\omega(\vec{\gamma};f) & :=\rho_0^1(i(pr_G(\gamma_{02}\vJoin ...\vJoin\gamma_{0q+1})))^{-1}\Big{(}\rho_0^1(i((g_{00}^{h_{01}}g_{01})^{h_{02}...h_{0q+1}}))^{-1}\omega(\delta_0\vec{\gamma}_{0,0};f^{h_{00}h_{01}},g_{00}^{h_{01}},g_{01})+ \\
    & \qquad -\omega(\delta_0\vec{\gamma}_{0,0};g_{01}^{-1},f^{h_{00}i(g_{00})h_{01}},g_{01})+ \\
    & \qquad\quad +\omega(\delta_0\vec{\gamma}_{0,0};(g_{00}^{h_{01}}g_{01})^{-1},f^{h_{00}h_{01}}g_{00}^{h_{01}},g_{01})-\omega(\delta_0\vec{\gamma}_{0,0};(g_{00}^{h_{01}}g_{01})^{-1},g_{00}^{h_{01}},g_{01})+ \\ 
    & \qquad\quad -\omega(\delta_0\vec{\gamma}_{0,0};g_{01}^{-1},(g_{00}^{h_{10}})^{-1},f^{h_{00}h_{01}}g_{00}^{h_{01}})+\omega(\partial_0\vec{\gamma}_{0,0};g_{01}^{-1},(g_{00}^{h_{10}})^{-1},g_{00}^{h_{01}})\Big{)}   
\end{align*}
for $\vec{\gamma}=(\gamma_0,...,\gamma_{q+1})\in\G_{p+1}^{q+2}$, $\gamma_b=(\gamma_{0b},...,\gamma_{pb})^T$ and $\gamma_{ab}\sim(g_{ab},h_{ab})$. Then
\begin{eqnarray*}
\delta\circ\Delta+\Delta\circ\delta=\delta_{(1)}\circ\Delta_2^p-\Delta_2^p\circ\delta_{(1)}.
\end{eqnarray*}
\end{lemma}

As announced, we move now to define all difference maps landing on the page $r=1$. We introduce pieces of notation that will be of help in writing things down efficiently. First, for $\vec{\gamma}\in\G_{p+k}^{q+1}$ and $k>0$, we define
\begin{eqnarray*}
\rho_k^q(\vec{\gamma}):=\rho_0^1\Big{(}i\big{(}pr_G\big{(}(\gamma_{01}\Join ...\Join\gamma_{k-11})\vJoin ...\vJoin(\gamma_{0q}\Join ...\Join\gamma_{k-1q})\big{)}\big{)}\Big{)}^{-1}.
\end{eqnarray*}
Next, at the level of spaces, we define two sequences of maps. We will assume our usual convention $\gamma=(\gamma_0,...,\gamma_{k+p-1})\in\G_{p+k}$ and $\gamma_b=(g_b,h_b)$, so that, if $h=h_{k+p-1}$, $\gamma\sim(g_0,...,g_{k+p-1},h)$. For $1<m\leq k$, we define
\begin{eqnarray*}
\xymatrix{
Q^t_m:\G_{p+k}\times G \ar[r] & G^{k+1}}\qquad\qquad\qquad\qquad\qquad \\
Q^t_m(\gamma ;f)=(g_{0}^{-1},g_{1}^{-1},...,g_{k-m-1}^{-1},(g_{k-2}...g_{k-m})^{-1},f^{h_{k-2}},g_{k-2},...,g_{k-m}),
\end{eqnarray*}
and for $1<m\leq k$, we define
\begin{eqnarray*}
\xymatrix{
Q^s_m:\G_{p+k}\times G \ar[r] & G^{k+1}}\qquad\qquad\qquad\qquad\qquad\qquad \\
Q^s_m(\gamma ;f)=(g_{0}^{-1},g_{1}^{-1},...,g_{k-m-1}^{-1},(g_{k-1}...g_{k-m})^{-1},f^{h_{k-1}}g_{k-1},g_{k-2},...,g_{k-m}).
\end{eqnarray*}
With these, let
\begin{eqnarray*}
\xymatrix{
\Delta_k^q:C(\G_p^q\times G^{k+1},W) \ar[r] & C(\G_{p+k}^{q+1}\times G,W)
}
\end{eqnarray*}
be defined by
\begin{align*}
\Delta & _k^q(\vec{\gamma};f)=\rho_k^q(\vec{\gamma})\Big{(}\rho_0^1(i((g_{k-10}...g_{00})^{}))^{-1}\omega(\partial_0^{k-1}\vec{\gamma}_{0,0};f^{h_{k-10}},g_{k-10},...,g_{00})+ \\
       &  +\sum_{m=2}^{k}(-1)^{p+m-k}\omega(\partial_0^{k-1}\vec{\gamma}_{0,0};Q^t_m(\gamma_0;f))+\sum_{m=1}^{k}(-1)^{m-k}\big{(}\omega(\partial_0^{k-1}\vec{\gamma}_{0,0};Q^s_m(\gamma_0;f))-\omega(\partial_0^{k-1}\vec{\gamma}_{0,0};Q^s_m(\gamma_0;1))\big{)}\Big{)}.
\end{align*}
Analogously, for $\vec{\gamma}\in\G_{p+1}^{q+k}$ and $k>0$, we define
\begin{eqnarray*}
\rho_k^p(\vec{\gamma}):=\rho_0^1\Big{(}i\big{(}pr_G(\gamma_{0k}\vJoin ...\vJoin\gamma_{0q+r-1})\big{)}\Big{)}^{-1}.
\end{eqnarray*}
Also, the sequences of maps at the level of spaces are defined assuming $\vec{\gamma}=(\gamma_0,...,\gamma_{k+q-1})\in\G^{q+k}$ and $\gamma_b=(g_b,h_b)$. For $1\leq m<k$, we define
\begin{eqnarray*}
\xymatrix{
P^t_m:\G^{q+k}\times G \ar[r] & G^{k+1}
}
\end{eqnarray*}
\begin{align*}
P^t_m(\vec{\gamma};f)=\Big{(}(g_{k-1})^{-1},...,(g_{m+1}^{h_{m+2}...h_{k-1}})^{-1}, & (g_{1}^{h_{2}...h_{k-1}}...g_{m}^{h_{m+1}...h_{k-1}})^{-1}, \\  
 & \quad ,f^{h_{0}i(g_{0})h_{1}...h_{k-1}},g_{1}^{h_{2}...h_{k-1}},...,g_{m}^{h_{m+1}...h_{k-1}}\Big{)},
\end{align*}
and for $0\leq m<k$, we define
\begin{eqnarray*}
\xymatrix{
P^s_m:\G^{q+k}\times G \ar[r] & G^{k+1}
}
\end{eqnarray*}
\begin{align*}
P^s_m(\vec{\gamma};f)=\Big{(}(g_{k-1})^{-1},...,(g_{m+1}^{h_{m+2}...h_{k-1}})^{-1}, & (g_{0}^{h_{1}...h_{k-1}}...g_{m}^{h_{m+1}...h_{k-1}})^{-1}, \\  
 & \quad ,(f^{h_{0}}g_0)^{h_{1}...h_{k-1}},g_{1}^{h_{2}...h_{k-1}},...,g_{m}^{h_{m+1}...h_{0k-1}}\Big{)},
\end{align*}
With these, let
\begin{eqnarray*}
\xymatrix{
\Delta_k^p:C(\G_p^q\times G^{k+1},W) \ar[r] & C(\G_{p+1}^{q+k}\times G,W)
}
\end{eqnarray*}
be defined by
\begin{align*}
\Delta & _k^p(\vec{\gamma};f)=\rho_k^p(\vec{\gamma})\Big{(}\rho_0^1(i(g_{00}^{h_{01}...h_{0k-1}}...g_{0k-1}))^{-1}\omega(\delta_0^{k-1}\vec{\gamma}_{0,0};f^{h_{00}...h_{0k-1}},g_{00}^{h_{01}...h_{0k-1}},...,g_{0k-1}) + \\
       & +\sum_{m=1}^{k-1}(-1)^{q+m-k}\omega(\delta_0^{k-1}\vec{\gamma}_{0,0};P^t_m(\vec{\gamma}_0;f))+\sum_{m=0}^{k-1}(-1)^{m-k}\big{(}\omega(\delta_0^{k-1}\vec{\gamma}_{0,0};P^s_m(\vec{\gamma}_0;f))-\omega(\delta_0^{k-1}\vec{\gamma}_{0,0};P^s_m(\vec{\gamma}_0;1))\big{)}\Big{)},
\end{align*}
where $\vec{\gamma}_0$ is shorthand for $(\gamma_{00},...,\gamma_{0q+k-1})$. These maps verify the equations
\begin{eqnarray*}
\partial\Delta_k^q+\Delta_k^q\partial=\delta_{(1)}\Delta_{k+1}^q-\Delta_{k+1}^q\delta_{(1)} & \textnormal{and} & \partial\Delta_k^p+\Delta_k^p\partial=\delta_{(1)}\Delta_{k+1}^p-\Delta_{k+1}^p\delta_{(1)}. 
\end{eqnarray*}
They also are key to understanding how the maps landing on an arbitrary $r$-page are defined. Indeed, in the general formula, each of the terms in the sums defining $\Delta_k^q$ and $\Delta_k^p$ will be shuffled by the $\Delta^n$ operators that we introduced when we defined the first difference maps. That is, e.g. in the place of the term $\omega(\partial_0^{k-1}\vec{\gamma}_{0,0};Q^s_m(\vec{\gamma};f))$, when evaluating at a general term $\vec{f}=(f_1,...,f_{r-k})$ instead, there will be a sum
\begin{eqnarray*}
\sum_{n=0}^{r-k-1}(-1)^n\omega(\partial_0^{k-1}\vec{\gamma}_{0,0};Q^s_m(\vec{\gamma};f_{r-k})_{[0,k-m)},\Delta^n(\gamma_{k-m0}\Join ...\Join\gamma_{k-10};\vec{f}),Q^s_m(\vec{\gamma};f_{r-k})_{[k-m+1,k+1]}).
\end{eqnarray*}
The general formula for 
\begin{eqnarray*}
\xymatrix{
\Delta_k^q:C(\G_p^q\times G^r,W) \ar[r] & C(\G_{p+k}^{q+1}\times G^{r-k},W)
}
\end{eqnarray*}
is then given by the rather cumbersome
\begin{align*}
\Delta & _k^q(\vec{\gamma};\vec{f})=\rho_k^q(\vec{\gamma})\Big{(}\rho_0^1(i((g_{k-10}...g_{00})^{}))^{-1}\omega(\partial_0^{k-1}\vec{\gamma}_{0,0};\vec{f}^{h_{k-10}},g_{k-10},...,g_{00})+ \\
       &  +\sum_{m=2}^{k}\sum_{n=0}^{r-k-1}(-1)^{p+m-k+n}\omega(\partial_0^{k-1}\vec{\gamma}_{0,0};Q^t_m(\vec{\gamma};f_{r-k})_{[0,k-m)},\Delta^n(\gamma_{k-m0}\Join ...\Join\gamma_{k-10};\vec{f}),Q^t_m(\vec{\gamma};f_{r-k})_{[k-m+1,k+1]})+ \\
       & +\sum_{m=1}^{k}\sum_{n=0}^{r-k-1}(-1)^{m-k+n}\big{(}\omega(\partial_0^{k-1}\vec{\gamma}_{0,0};Q^s_m(\vec{\gamma};f_{r-k})_{[0,k-m)},\Delta^n(\gamma_{k-m0}\Join ...\Join\gamma_{k-10};\vec{f}),Q^s_m(\vec{\gamma};f_{r-k})_{[k-m+1,k+1]})+ \\
       & \qquad\qquad -\omega(\partial_0^{k-1}\vec{\gamma}_{0,0};Q^s_m(\vec{\gamma};f_{r-k})_{[0,k-m)},\Delta^n(\gamma_{k-m0}\Join ...\Join\gamma_{k-10};\vec{f}),Q^s_m(\vec{\gamma};1)_{[k-m+1,k+1]})\big{)}\Big{)}.
\end{align*}
and the formula for $\Delta_k^p$ is defined likewise.

\subsubsection{Higher difference maps 3: the general case}
We already saw that trying and defining the differential for the complex $C(\G,\phi)$ as
\begin{eqnarray*}
d=(-1)^r\delta+(-1)^{q+1}\partial+(-1)^{r}\delta_{(1)}+\sum_k(-1)^{p+q+k-1}(\Delta_k^p+\Delta_k^q)
\end{eqnarray*}
up to signs, will still fail to verify $d^2=0$, as there are still relations unaccounted for. Inductively, studying how these relations behave, one finds the missing difference maps 
\begin{eqnarray*}
\xymatrix{
\Delta_{a,b}:C^{p,q}_r(\G,\phi) \ar[r] & C^{p+a,q+b}_{r-k}(\G,\phi),
}
\end{eqnarray*}
where $1\leq k\leq r$ and $(a,b)$ is a pair of strictly positive natural numbers satisfying $a+b=k+1$. These maps will be defined in general as
\begin{eqnarray*}
\Delta_{a,b}\omega(\vec{\gamma};\vec{f}):=\sum_{n}(-1)^{\alpha}\omega(\partial_0^{a-1}\delta_0^{b-1}\vec{\gamma}_{0,0};[Q_n,P_n](\vec{\gamma};\vec{f}))
\end{eqnarray*}
for $\vec{\gamma}\in\G_{p+a}^{q+b}$ and $\vec{f}\in G^{r-k}$. These maps will include the difference maps we have already got. For instance, when $k=1$, $\Delta_{1,1}$ is the first difference map $\Delta$. In the special cases $(a,b)\in\lbrace(1,k),(k,1)\rbrace$, these maps are given by the formulas of $\Delta_k^p$ and $\Delta_k^q$ respectively. 
\begin{theorem}\label{GpDiffs}
Thus defined, the maps $\Delta_{a,b}$ verify the following set of equations for elements in $C^{p,q}_r(\G,\phi)$
\begin{itemize}
    \item[i)] $\delta\Delta_{k}^p+\Delta_{k}^p\delta=\delta_{(1)}\Delta_{k+1}^p+\Delta_{k+1}^p\delta_{(1)}$ for all $1\leq k<r$, and
    $\delta\Delta_{r}^p+\Delta_{r}^p\delta=\delta_{(1)}\Delta_{r+1}^p$, 
    \item[ii)]$\partial\Delta_{k}^q+\Delta_{k}^q\partial=\delta_{(1)}\Delta_{k+1}^q+\Delta_{k+1}^q\delta_{(1)}$ for all $1\leq k<r$, and
    $\partial\Delta_{r}^q+\Delta_{r}^q\partial=\delta_{(1)}\Delta_{r+1}^q$, 
    \item[iii)] assuming the convention that $\Delta_0^p=\partial$, for each $a$
    \begin{eqnarray*}
    \sum_{i=0}^k[\Delta_{a,b-i},\Delta_i^p]=[\delta,\Delta_{b,a}]+[\Delta_{a,b},\delta_{(1)}]
    \end{eqnarray*}
    for all $1\leq k\leq r$, and a symmetric set of equations for $q$.
\end{itemize}
\end{theorem}
Notice that with respect to the diagonal grading, all difference maps are homogeneous of degree $+1$. Then, we can use them as differentials and indeed, we've got the following.
\begin{theorem}\label{The2GpCx}
$C(\G,\phi)$ graded by
\begin{eqnarray*}
C^{n}_{tot}(\G,\phi)=\bigoplus_{p+q+r=n}C^{p,q}_r(\G,\phi)
\end{eqnarray*}
together with the differential
\begin{eqnarray*}
\nabla=(-1)^r\delta+(-1)^{q+1}\partial+(-1)^{r}\delta_{(1)}+\sum_k(-1)^{p+q+k-1}(\Delta_k^p+\Delta_k^q) +\sum_{a+b=k+1}\Delta_{a,b}
\end{eqnarray*}
is a complex.
\end{theorem}
\begin{remark}\label{DblCxForGpsRmk}
Recall that, at the beginning of the section, we set out to look for a double complex and the non-commutativity of the squares forced us to consider the successive $r$-pages. As it turns out, in sight of theorem \ref{GpDiffs} and of their grading, we see that there is no way of thinking of the difference maps as being maps of double complexes between the $p$-pages in contrast with the case of Lie $2$-algebras. What we got instead is that the total complexes of the $p$-pages form a $2$-dimensional lattice that we will provisionally call $C_{p-pag}^{ab}$, together with differentials
\begin{eqnarray*}
\xymatrix{
d_r^{ab}:C_{p-pag}^{ab} \ar[r] & C_{p-pag}^{a+r,b-r+1}
}
\end{eqnarray*}
whose sum yields an actual differential on the total complex of $C_{p-pag}$.
\end{remark}

\section{2-cohomology of Lie 2-groups}
In this section, we study the cohomology of the complex from theorem \ref{The2GpCx}. \\
Start out with an element $v\in V=C^0(\G,\phi)=C^{0,0}_0(\G,\phi)$, then its differential is given by
\begin{eqnarray*}
\nabla v=(\delta v,\delta_{(1)}v,-\cancelto{0}{\partial v})\in C^1(\G,\phi)
\end{eqnarray*}
If $v$ is a $0$-cocycle, then for all $h\in H$ and for all $g\in G$,
\begin{eqnarray*}
\rho_0^0(h)v=v & \textnormal{and} & \rho_1(g)v=0;
\end{eqnarray*}
therefore, since $\rho_0^1(h)$ is an automorphism of $W$,
\begin{eqnarray*}
H^0_\nabla(\G,\phi)=V^{\G}:=\lbrace v\in V:\bar{\rho}_{(g,h)}(0,v)=(0,v),\quad\forall (g,h)\in G\rtimes H\cong\G\rbrace ,
\end{eqnarray*}
where $\bar{\rho}$ is the honest representation of proposition \ref{honestGpRep}. \\
A $1$-cochain $\lambda$ is a triple $(\lambda_0,\lambda_1,v)\in C(H,V)\oplus C(G,W)\oplus V$, and its differential has six entries that we write using their coordinates, i.e. $\nabla\lambda^{p,q}_r\in C^{p,q}_r(\G,\phi)$
\begin{align*}
\nabla\lambda^{0,2}_0 & =\delta\lambda_0                             &                       & & &                               \\
\nabla\lambda^{1,1}_0 & =\partial\lambda_0+\delta v+\Delta\lambda_1  & \nabla\lambda^{0,1}_1 & =\delta_{(1)}\lambda_0-\delta\lambda_1 & & \\
\nabla\lambda^{2,0}_0 & =-\partial v=-v                              & \nabla\lambda^{1,0}_1 & =\delta_{(1)}v+\partial\lambda _1=\delta_{(1)}v & \nabla\lambda^{0,0}_2 & =\delta_{(1)}\lambda_1 
\end{align*}
Schematically,
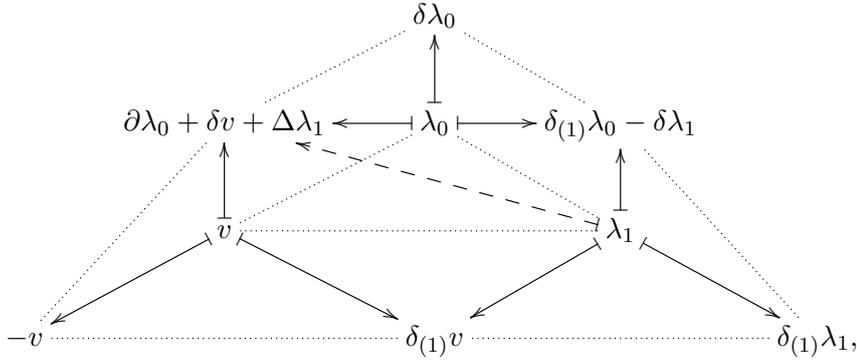
\begin{eqnarray*}
\xymatrix{
  & & \delta\lambda_0 \ar@{.}[dl]\ar@{.}[dr] & &  \\
  & \partial\lambda_0+\delta v+\Delta\lambda_1 & \lambda_0 \ar@{|->}[u]\ar@{|->}[r]\ar@{|->}[l]\ar@{.}[dr] & \delta_{(1)}\lambda_0-\delta\lambda_1 & \\
  & v \ar@{|->}[dl]\ar@{|->}[u]\ar@{|->}[dr]\ar@{.}[ur]\ar@{.}[rr] &  & \lambda_1 \ar@{|->}[dr]\ar@{|->}[u]\ar@{|->}[dl]\ar@{|-->}[ull] & \\
-v \ar@{.}[uur] \ar@{.}[rr] & & \delta_{(1)}v \ar@{.}[rr] & & \delta_{(1)}\lambda_1 , \ar@{.}[uul]
}
\end{eqnarray*}
where the solid arrows represent the differentials, the dashed arrow represents the difference map, and the pointed polygons represent elements of the same degree. If $\lambda$ is a $1$-cocylce, then $v=0$, $\lambda_0$ and $\lambda_1$ are respectively crossed homomorphisms of $H$ and $G$ into $V$ and $W$ with respect to $\rho_0^0$ and $\rho_0^1\circ i$. In symbols,
\begin{align*}
    \lambda_0(h_0h_1) & =\lambda_0(h_0)+\rho_0^0(h_0)\lambda_0(h_1),    & \forall h_0,h_1\in H \\
    \lambda_1(g_0g_1) & =\lambda_1(g_0)+\rho_0^1(i(g_0))\lambda_1(g_1), & \forall g_0,g_1\in G.
\end{align*}
Additionally, the following relations hold for each $\gamma=(g,h)\in\G$:
\begin{align*}
    (\partial\lambda_0+\Delta\lambda_1)(\gamma) & =0 \\
    \lambda_0(h)+\phi(\rho_0^1(h)\lambda_1(g))  & =\lambda_0(hi(g)),
\end{align*}
and
\begin{align*}
    (\delta_{(1)}\lambda_0-\delta\lambda_1)(\gamma)                    & =0 \\
    \rho_0^1(h)^{-1}\rho_1(g)\lambda_0(h)+\rho_0^1(h)^{-1}\lambda_1(g) & =\lambda_1(g^h).
\end{align*}
The first of these relations says that there is a map
\begin{eqnarray*}
\xymatrix{
\bar{\lambda}:\G \ar[r] & W\oplus V:(g,h) \ar@{|->}[r] & (\rho_0^1(h)\lambda_1(g),\lambda_0(h))
}
\end{eqnarray*}
that respects the source and the target. In fact, notice that when $h=1$, the relation reduces to the commutativity of 
\begin{eqnarray*}
\xymatrix{
G \ar[r]^{\lambda_1}\ar[d]_i & W \ar[d]^\phi \\
H \ar[r]_{\lambda_0}         & V.
}
\end{eqnarray*}
Moreover, together with the fact that $\lambda_1$ is a crossed homomorphism, it implies that $\bar{\lambda}$ is a functor. Indeed, if $(\gamma_0,\gamma_1)\in\G_2$
\begin{align*}
    \bar{\lambda}(\gamma_0\Join\gamma_1) & =\bar{\lambda}(g_1g_0,h) \\
                                         & =(\rho_0^1(h)\lambda_1(g_1g_0),\lambda_0(h)) \\
                                         & =(\rho_0^1(h)\big{(}\lambda_1(g_1)+\rho_0^1(i(g_1))\lambda_1(g_0)\big{)},\lambda_0(h)) \\
                                         & =\big{(}\rho_0^1(hi(g_1))\lambda_1(g_0),\lambda_0(hi(g_1))\big{)}\big{(}\rho_0^1(h)\lambda_1(g_1),\lambda_0(h)\big{)}=\bar{\lambda}(\gamma_0)\bar{\lambda}(\gamma_1),
\end{align*}
where the composition in the last line makes sense, precisely because $\lambda_0$ is a crossed homomorphism, so that 
\begin{eqnarray*}
    \lambda_0(h)+\phi\big{(}\rho_0^1(h)\lambda_1(g)\big{)}=\lambda_0(h)+\rho_0^0(h)\phi(\lambda_1(g))=\lambda_0(h)+\rho_0^0(h)\lambda_0(i(g))=\lambda_0(hi(g)).
\end{eqnarray*}
The second relation, says that $\bar{\lambda}$ is a crossed homomorphism into $W\oplus V$ with respect to $\bar{\rho}$. To see this, let $\gamma_k=(g_k,h_k)\in G\rtimes H$ for $k\in\lbrace 0,1\rbrace$, then
\begin{align*}
    \bar{\lambda}(\gamma_0\vJoin\gamma_1) & =\bar{\lambda}(g_0^{h_1}g_1,h_0h_1) \\
                                          & =(\rho_0^1(h_0h_1)\lambda_1(g_0^{h_1}g_1),\lambda_0(h_0h_1)) \\  
                                          & =(\rho_0^1(h_0h_1)\big{(}\lambda_1(g_0^{h_1})+\rho_0^1(i(g_0^{h_1}))\lambda_1(g_1)\big{)},\lambda_0(h_0)+\rho_0^0(h_0)\lambda_0(h_1)) \\ 
                                          & =(\rho_0^1(h_0)\big{(}\rho_1(g_0)\lambda_0(h_1)+\lambda_1(g_0) +\rho_0^1(i(g_0)h_1)\lambda_1(g_1)\big{)},\lambda_0(h_0)+\rho_0^0(h_0)\lambda_0(h_1)) \\
                                          & =\bar{\lambda}(\gamma_0)+\bar{\rho}_{(g_0,h_0)}\bar{\lambda}(\gamma_1)
\end{align*}

A coboundary $\nabla v$, will induce a crossed homomorphism-functor that can also be seen as $\bar{\rho}_{(g,h)}(0,v)$; by analogy, we call these \textit{principal crossed homomorphisms}. We introduce the following notation
\begin{Def}
The \textit{space of crossed homomorphisms} of a Lie $2$-group $\G$ with respect to a $2$-representation $\rho$ on the $2$-vector $\mathbb{V}=\xymatrix{W \ar[r]^\phi & V}$ is defined to be
\begin{align*}
CrHom(\G,\phi):=\lbrace \bar{\lambda}\in Hom_{Gpd}(\G,\mathbb{V}) & :\bar{\lambda}(g,h)=(\rho_0^1(h)\lambda_1(g),\lambda_0(h)) \\
 & \qquad\textnormal{ is a crossed homomorphism with respect to }\bar{\rho}\rbrace .
\end{align*}
The \textit{space of principal crossed homomorphisms} is defined to be
\begin{eqnarray*}
PrCrHom(\G,\phi):=\lbrace \bar{\lambda}\in CrHom(\G,\phi):\bar{\lambda}(g,h)=\bar{\rho}_{(g,h)}(0,v)\textnormal{ for some }v\in V\rbrace .
\end{eqnarray*}
\end{Def}
With these definitions, we can write
\begin{eqnarray*}
H^1_\nabla(\G,\phi)=CrHom(\G,\phi)/PrCrHom(\G,\phi).
\end{eqnarray*}
A $2$-cochain $\vec{\omega}$ is a $6$-tuple $(\omega_0,\alpha,\varphi,\omega_1,\lambda,v)$ where
\begin{align*}
\omega_0 & \in C(H^2,V) &          &                    &          &   \\
\varphi  & \in C(\G,V)  & \alpha   & \in C(H\times G,W) &          &   \\
v        & \in V        & \lambda  & \in C(G,W)         & \omega_1 & \in C(G^2,W) . 
\end{align*} 
Let's compute the differential using coordinates as before. First, $\nabla\vec{\omega}^{3,0}_0=-\partial v=0$; for the other values, 
\begin{align*}
\nabla\vec{\omega}^{0,3}_0 & =\delta\omega_0 & & & &  \\
\nabla\vec{\omega}^{1,2}_0 & =-\partial\omega_0+\delta\varphi-\Delta\alpha-\Delta_2^p\omega_1 & \nabla\vec{\omega}^{0,2}_1 & =\delta_{(1)}\omega_0-\delta\alpha & & \\
\nabla\vec{\omega}^{2,1}_0 & =\partial\varphi+\delta v-\Delta\lambda-\Delta_2^q\omega_1 & \nabla\vec{\omega}^{1,1}_1 & =-\delta_{(1)}\varphi+\partial\alpha+\delta\lambda+\Delta\omega_1 & \nabla\vec{\omega}^{0,1}_2 & =\delta_{(1)}\alpha+\delta\omega_1 \\
\nabla\vec{\omega}^{2,0}_1 & =\delta_{(1)}v-\cancelto{\lambda}{\partial\lambda} & \nabla\vec{\omega}^{1,0}_2 & =-\cancelto{0}{\partial\omega_1}-\delta_{(1)}\lambda & \nabla\vec{\omega}^{0,0}_3 & =\delta_{(1)}\omega_1
\end{align*}
If we assume  $v=0$, $\lambda$ vanishes too and these equations correspond to those in the statement of proposition \ref{Gp2-cocycles}. First, notice that under these assumptions, the equation at $C^{2,1}_0(\G,\phi)$ says 
\begin{eqnarray}
\partial\varphi=\Delta_2^q\omega_1.
\end{eqnarray} 
Evaluated at an arbitrary element $(g_0,g_1,h)\in\G_2$, we get
\begin{eqnarray}
\varphi\begin{pmatrix}
g_1g_0 \\
h
\end{pmatrix}=\varphi\begin{pmatrix}
g_1 \\
h
\end{pmatrix}+\varphi\begin{pmatrix}
g_0 \\
hi(g_1)
\end{pmatrix}-\phi\big{(}\rho_0^1(h)\omega_1(g_1,g_0)\big{)}.
\end{eqnarray} 
In particular, making $g_0=1$, we recover the fact that $\varphi$ vanishes on the space of units $H$. Since the equations in proposition \ref{Gp2-cocycles} involve a function that depends only on $G$ with values in $V$, in order to compare them to the cocycle equations above, we introduce \begin{eqnarray*}
\Check{\varphi}(g):=\varphi\begin{pmatrix}
g \\
1
\end{pmatrix}.
\end{eqnarray*}
Now we are ready. For the numbering of the equations in the statement of proposition \ref{Gp2-cocycles}, we've got:
\begin{itemize}
\item Equation $i)$ is literally $\delta\omega_0 =0$.  
\item Equation $ii)$ is literally $\delta_{(1)}\omega_1=0$.
\item Equation $iii)$ is $\nabla\vec{\omega}^{1,2}_0=0$ evaluated at $\begin{pmatrix}
g_1 & g_2 \\
1   & 1
\end{pmatrix}\in\G^2$. 
\item Equation $iv)$ is exactly $\nabla\vec{\omega}^{0,2}_1=0$. % We're cheating a little bit. There is a sign missing: On \Delta\alpha. (not anymore)
\item Equation $v)$ follows from $\nabla\vec{\omega}^{1,2}_0=0$ in the following manner. Evaluating at
\begin{align*}
    & \begin{pmatrix}
g      & 1 \\
h^{-1} & h
\end{pmatrix}: & -\big{(}\omega_0(h^{-1},h)-\omega_0(h^{-1}i(g),h)\big{)}-\Check{\varphi}(g^h)+\varphi\begin{pmatrix}
g \\
h^{-1}
\end{pmatrix}-\phi(\alpha(h;g)) & =0 \\
    & \begin{pmatrix}
1      & g \\
h^{-1} & 1
\end{pmatrix}: & \omega_0(h^{-1},i(g))+\rho_0^0(h)^{-1}\Check{\varphi}(g)-\varphi\begin{pmatrix}
g \\
h^{-1}
\end{pmatrix} & =0
\end{align*}
and summing yields
\begin{eqnarray*}
-\omega_0(h^{-1},h)+\omega_0(h^{-1}i(g),h)+\omega_0(h^{-1},i(g))=\Check{\varphi}(g^h)-\rho_0^0(h)^{-1}\Check{\varphi}(g)+\phi(\alpha(h;g)),
\end{eqnarray*}
but $\delta\omega_0=0$; therefore,
\begin{eqnarray*}
\rho_0^0(h)^{-1}\omega_0(i(g),h)+\omega_0(h^{-1}i(g),h)=\omega_0(h^{-1}i(g),h)+\omega_0(h^{-1},i(g)),
\end{eqnarray*}
so replacing one gets the desired equation.
\item Equation $vi)$ is $\nabla\vec{\omega}^{1,1}_1=0$ evaluated at $(\gamma;g_1)\in\G\times G$, where $\gamma=\begin{pmatrix}
g_2 \\
1  
\end{pmatrix}\in\G$.
\item Equation $vii)$ is exactly $\nabla\vec{\omega}^{0,1}_2=0$.
\end{itemize}
Conversely, let $(\omega_0,\omega_1,\alpha,\varphi)$ be a $4$-tuple as in the statement of proposition \ref{Gp2-cocycles} and verifying the equations therein. Define
\begin{eqnarray*}
\hat{\varphi}\begin{pmatrix}
g \\
h
\end{pmatrix}:=\omega_0(h,i(g))+\rho_0^0(h)\varphi(g).
\end{eqnarray*}
Then, if we define $\vec{\omega}=(\omega_0,\alpha,\hat{\varphi},\omega_1,0,0)$, from the discussion above, we just need to check that $\nabla\vec{\omega}^{2,1}_0$, $\nabla\vec{\omega}^{1,1}_1$ and $\nabla\vec{\omega}^{1,2}_0$ vanish. Indeed, consider
\begin{align*}
\hat{\varphi}\begin{pmatrix}
g_1g_0 \\
h
\end{pmatrix} & =\omega_0(h,i(g_1g_0))+\rho_0^0(h)\varphi(g_1g_0). 
\end{align*}
From the cocycle equation $\delta\omega_0=0$, we know that 
\begin{eqnarray*}
\omega_0(h,i(g_1g_0))=-\rho_0^0(h)\omega_0(i(g_1),i(g_0))+\omega_0(hi(g_1),g_0)+\omega_0(h,i(g_1));
\end{eqnarray*}
thus, replacing and using equation $iii)$ in the proposition, 
\begin{align*}
\hat{\varphi}\begin{pmatrix}
g_1g_0 \\
h
\end{pmatrix} & =\omega_0(hi(g_1),g_0)+\omega_0(h,i(g_1))+\rho_0^0(h)\big{[}\rho_0^1(i(g_1))\varphi(g_0)+\varphi(g_1)-\phi(\omega_1(g_1,g_0))\big{]} \\
              & =\hat{\varphi}\begin{pmatrix}
g_0 \\
hi(g_1)
\end{pmatrix}+\hat{\varphi}\begin{pmatrix}
g_1 \\
h
\end{pmatrix}-\phi(\rho_0^1(h)\omega_1(g_1,g_0))
\end{align*}
which we saw to coincide with $\nabla\vec{\omega}^{2,1}_0=0$. Next, let $\gamma=(g,h)\in\G$ and $f\in G$ and consider
\begin{align*}
\delta_{(1)}\hat{\varphi}(\gamma;f) & =\rho_0^1(hi(g))^{-1}\rho_1(f)\hat{\varphi}(\gamma) \\
                                    & =\rho_0^1(hi(g))^{-1}\rho_1(f)\big{(}\omega_0(h,i(g))+\rho_0^0(h)\varphi(g)\big{)} \\
                                    & =\rho_0^1(hi(g))^{-1}\rho_1(f)\omega_0(h,i(g))+\rho_0^1(i(g))^{-1}\rho_1(f^h)\varphi(g). 
\end{align*}
Using equations $iv)$ and $vi)$ in the proposition,
\begin{align*}
\delta_{(1)}\hat{\varphi}(\gamma;f) & =\rho_0^1(i(g))^{-1}\alpha(h;f)-\alpha(hi(g);f)+\alpha(i(g);f^h)+ \\
                                    & \qquad\quad -\alpha(i(g),f^h)+\rho_0^1(i(g))^{-1}\omega_1(f^h,g)+\omega_1(g^{-1},f^hg)-\omega_1(g^{-1},g) \\
                                    & =\partial\alpha(\gamma;f)+\Delta\omega_1(\gamma;f)
\end{align*}
which is precisely $\nabla\vec{\omega}^{1,1}_1=0$. As for $\nabla\vec{\omega}^{1,2}_0$, consider
\begin{align*}
\delta\hat{\varphi}\begin{pmatrix}
g_0 & g_1 \\
h_0 & h_1
\end{pmatrix} & =\rho_0^0(h_0i(g_0))\hat{\varphi}\begin{pmatrix}
g_1 \\
h_1
\end{pmatrix}-\hat{\varphi}\begin{pmatrix}
g_0^{h_1}g_1 \\
h_0h_1
\end{pmatrix}+\hat{\varphi}\begin{pmatrix}
g_0 \\
h_0
\end{pmatrix}\\
              & =\rho_0^0(h_0i(g_0))\big{[}\omega_0(h_1,i(g_1))+\rho_0^0(h_1)\varphi(g_1)\big{]}+ \\
              & \qquad -\big{(}\omega_0(h_0h_1,i(g_0^{h_1}g_1))+\rho_0^0(h_0h_1)\varphi(g_0^{h_1}g_1)\big{)}+\omega_0(h_0,i(g_0))+\rho_0^0(h_0)\varphi(g_0)
\end{align*}
Using successively equations $iii)$ and $v)$ in the proposition
\begin{align*}
\varphi(g_0^{h_1}g_1) & =\rho_0^0(i(g_0^{h_1}))\varphi(g_1)+\varphi(g_0^{h_1})+\omega_0(i(g_0^{h_1}),i(g_1))-\phi(\omega_1(g_0^{h_1},g_1)) \\
                      & =\rho_0^0(i(g_0^{h_1}))\varphi(g_1)+\rho_0^0(h_1)^{-1}\omega_0(i(g_0),h_1)+\omega_0(h_1^{-1},i(g_0)h_1)-\omega_0(h_1^{-1},h_1)+ \\
                      & \qquad +\rho_0^0(h_1)^{-1}\varphi(g_0)-\phi(\alpha(h_1;g_0))+\omega_0(i(g_0^{h_1}),i(g_1))-\phi(\omega_1(g_0^{h_1},g_1)). 
\end{align*}
Replacing,
\begin{align*}
\delta\hat{\varphi}\begin{pmatrix}
g_0 & g_1 \\
h_0 & h_1
\end{pmatrix} & =\rho_0^0(h_0h_1)\phi(\alpha(h_1;g_0)+\omega_1(g_0^{h_1},g_0))+R(\omega_0) \\
              & =\Delta\alpha\begin{pmatrix}
g_0 & g_1 \\
h_0 & h_1
\end{pmatrix}+\Delta_2^p\omega_1\begin{pmatrix}
g_0 & g_1 \\
h_0 & h_1
\end{pmatrix}+R(\omega_0),
\end{align*}
where 
\begin{align*}
R(\omega_0) & =\rho_0^0(h_0i(g_0))\omega_0(h_1,i(g_1))-\omega_0(h_0h_1,i(g_0^{h_1}g_1))+\omega_0(h_0,i(g_0)) \\
            & \qquad -\rho_0^0(h_0h_1)\big{(}\rho_0^0(h_1)^{-1}\omega_0(i(g_0),h_1)+\omega_0(h_1^{-1},i(g_0)h_1)-\omega_0(h_1^{-1},h_1)+\omega_0(i(g_0^{h_1}),i(g_1))\big{)}.
\end{align*}
We claim that $R(\omega_0)=\partial\omega_0\begin{pmatrix}
g_0 & g_1 \\
h_0 & h_1
\end{pmatrix}$, and consequently $\nabla\vec{\omega}^{1,2}_0=0$. Indeed, this will follow from repeatedly applying the cocycle equation $\delta\omega_0$ to
\begin{align*}
(h_0h_1,i(g_0^{h_1}),i(g_1)) & : \\ 
-\rho_0^0(h_0h_1)\omega_0 & (i(g_0^{h_1}),i(g_1))-\omega_0(h_0h_1,i(g_0^{h_1}g_1)) =-\omega_0(h_0i(g_0)h_1,i(g_1))-\omega_0(h_0h_1,i(g_0^{h_1})) \\
(h_1^{-1},i(g_0),h_1)        & : \\ 
\rho_0^0(h_1^{-1}) & \omega_0(i(g_0),h_1)+\omega_0(h_1^{-1},i(g_0)h_1)=\omega_0(h_1^{-1}i(g_0),h_1)+\omega_0(h_1^{-1},i(g_0)) \\ 
(h_0i(g_0),h_1,i(g_1))       & : \\ 
\rho_0^0(h_0i(g_0)) & \omega_0(h_1,i(g_1))-\omega_0(h_0i(g_0)h_1,i(g_1))=-\omega_0(h_0i(g_0),h_1i(g_1))+\omega_0(h_0i(g_0),h_1)  \\ 
(h_0h_1,h_1^{-1},i(g_0))     & : \\ 
-\rho_0^0(h_0h_1) & \omega_0(h_1^{-1},i(g_0))+\omega_0(h_0,i(g_0))=\omega_0(h_0h_1,h_1^{-1}i(g_0))-\omega_0(h_0h_1,h_1^{-1}) \\ 
(h_0h_1,h_1^{-1},h_1)        & : \\ 
& \rho_0^0(h_0h_1)\omega_0(h_1^{-1},h_1)=\omega_0(h_0,h_1)+\omega_0(h_0h_1,h_1^{-1}) \\ 
(h_0h_1,h_1^{-1}i(g_0),h_1)  & : \\
-\rho_0^0(h_0h_1)\omega_0(h_1^{-1} & i(g_0),h_1)+\omega_0(h_0i(g_0),h_1)-\omega_0(h_0h_1,h_1^{-1}i(g_0)h_1)+\omega_0(h_0h_1,h_1^{-1}i(g_0))=0 
\end{align*}
Then,
\begin{align*}
R(\omega_0) & =\rho_0^0(h_0i(g_0))\omega_0(h_1,i(g_1))-\omega_0(h_0i(g_0)h_1,i(g_1))-\omega_0(h_0h_1,h_1^{-1}i(g_0)h_1)+\omega_0(h_0,i(g_0))+ \\
            & \qquad -\rho_0^0(h_0h_1)\big{(}\omega_0(h_1^{-1}i(g_0),h_1)+\omega_0(h_1^{-1},i(g_0))-\omega_0(h_1^{-1},h_1)\big{)} \\
            & =-\omega_0(h_0i(g_0),h_1i(g_1))+\omega_0(h_0i(g_0),h_1)-\omega_0(h_0h_1,h_1^{-1}i(g_0)h_1)+ \\
            & \qquad +\omega_0(h_0h_1,h_1^{-1}i(g_0))-\omega_0(h_0h_1,h_1^{-1})-\rho_0^0(h_0h_1)\omega_0(h_1^{-1}i(g_0),h_1)+\omega_0(h_0,h_1)+\omega_0(h_0h_1,h_1^{-1}) \\
            & =\partial\omega_0\begin{pmatrix}
g_0 & g_1 \\
h_0 & h_1
\end{pmatrix}    
\end{align*}
as claimed.
Furthermore, for a pair of $2$-cocycles 
\begin{eqnarray*}
(\vec{\omega})^k=(\omega_0^k,\alpha^k,\varphi^k,\omega_1^k,\lambda^k,v^k), & k\in\lbrace 1,2\rbrace
\end{eqnarray*}
with $v^k=0$, if they are cohomologous we recover the equations in proposition \ref{Gp2-coboundaries}. Indeed, if
\begin{eqnarray*}
(\vec{\omega})^2-(\vec{\omega})^1=\nabla(\lambda_0,\lambda_1,v),
\end{eqnarray*}
coordinate-wise we have got
\begin{align*}
\omega_0^2-\omega_0^1 & =\delta\lambda_0                            &                   & & &                               \\
\varphi^2-\varphi^1   & =\partial\lambda_0+\delta v+\Delta\lambda_1 & \alpha^2-\alpha^1 & =\delta_{(1)}\lambda_0-\delta'\lambda_1 & & \\
0                     & =-v                                         & 0                 & =\delta_{(1)}v                           & \omega_1^2-\omega_1^1 &  =\delta_{(1)}\lambda_1 ;
\end{align*}
thus, explicitly for $\gamma=(g,h)\in\G$ and for $(h;f)\in H\times G$
\begin{align*}
    (\varphi^2-\varphi^1)(\gamma) & =\lambda_0(h)-\lambda_0(hi(g))+\phi(\rho_0^1(h)\lambda_1(g)) \\
    (\alpha^2-\alpha^1)(h;f) & =\rho_0^1(h)^{-1}\rho_1(f)\lambda_0(h)-\big{(}\lambda_1(f^h)-\rho_0^1(h)^{-1}\lambda_1(f)\big{)}.
\end{align*}
Evaluating the first of these equations at $\gamma=(g,1)$, we recover equation $iii)$ of the referred proposition. The second of these equations is exactly equation $iv)$ in the proposition. Conversely, given two $4$-tuples whose difference verifies the equations in proposition \ref{Gp2-coboundaries},
\begin{align*}
((\hat{\varphi})^2-(\hat{\varphi})^1)(\gamma) & =\omega_0^2(h,i(g))+\rho_0^0(h)\varphi^2(g)-\omega_0^1(h,i(g))-\rho_0^0(h)\varphi^1(g) \\    
                                      & =\delta\lambda_0(h,i(g))+\rho_0^0(h)(\phi(\lambda_1(g))-\lambda_0(i(g))) \\
                                      & =\rho_0^0(h)\lambda_0(i(g))-\lambda_0(hi(g))+\lambda_0(h)+\phi(\rho_0^1(h)\lambda_1(g))-\rho_0^0(h)\lambda_0(i(g)) \\
                                      & =\partial\lambda_0(\gamma)+\Delta\lambda_1(\gamma).
\end{align*}
In the end, this is the theorem this theory is dedicated to.
\begin{theorem}\label{H2Gp}
$H^2_\nabla(\G,\phi)$ classifies $2$-extensions.
\end{theorem}

% ------------------------------------

% We hope to use this algebraic description to prove a general van Est theorem for Lie $2$-algebras. Next step is to see how far this theory can be stretched to include \LA -groups. The issue of defining a representation seems already complicated, as there is no ``flat abelian'' candidate in the quotient category of \LA -groups over a prescribed Lie group.
%***Of the things that this section lacks, but would be nice for it to have.3) Particular cases. The second cohomology of the ideal with values in a trivial representation, in a representation with values in the pair groupoid.4) All too hopeful, does Hochschild-Serre shows up? Look for it! 
\chapter{Relating the 2-group and 2-algebra theories}\label{estimateschapter}
\chaptermark{Relating the $2$-group and $2$-algebra theories}
%---------------------------------------------------------------------------------------------------------------------------------------------------
\section{Introduction}
%---------------------------------------------------------------------------------------------------------------------------------------------------
During the last section of each of the previous two chapters, we saw that the cohomology theories we constructed have got the property that their second cohomology group classify extensions in the corresponding categories. In this chapter, we will approach the integration of Lie $2$-algebras making use of this fact. More precisely, we will introduce a map from the complex of Lie $2$-group cochains to the complex of Lie $2$-algebra cochains and prove a series of van Est type theorems. 
%that there are good reasons to study double Lie groups. Also, that the problem in most occasions is topological in nature in a way that is controled by certain cohomologies. Given that the only strategy for integration that isn't crippled from the start out is van Est's, we try and understand the integration of Lie $2$-algebras from a cohomological point of view. We restrict our attention to this case of study, as it has already proven to have positive solution. We do so with the aim in mind of extending this procedure to \LA -groups. \\
The nice thing about this approach is that it has shown the following historical advantages: \\
It is functorial at each step. \\ 
It also works in infinite dimensions. \\
It was used (in infinite-dimensions) to show that certain Lie algebras do not integrate \cite{NonIntInfDim}. \\ % According to this survey https://pdfs.semanticscholar.org/f421/bb5ae85531224bde23d08c29cca97c015fd1.pdf , van Est and Korthagen showed that \gg is integrable if and only if a certain subgroup Π(\gg) ≤ \mathfrak{z}(\gg) of the center of \gg (the period group) is discrete (!) Doesn't this look a lot like Marius and Rui's result?
It was used in the non-integrable case to show how non-integrable Lie algebras maybe still integrate to \'etale Lie 2-groups. \cite{IntInfDim}.

%\begin{remark}
%Inequality $f(t)/t^3 < 1/s$ in \eqref{f_{10}} is a compatible condition with \eqref{f_{12}}. 
%\textcolor{red}{que significa?}
%In fact \eqref{eq:ARSP}  implies that $f(t)\geq C t^{\theta -1}$ with $\theta >4$.
%  if the two conditions hold then
%\[
%C t^{\theta -4}\leq \frac{f(t)}{t^3} \leq \frac{1}{s}
%\]
%But this leads to a contradiction, because the left inequality implies $\lim_{t\rightarrow \infty} f(t)/t^3= \infty$, while the right inequality and \eqref{f_{10}} imply $\lim_{t\rightarrow \infty} f(t)/t^3= 1/s<\infty$.
%\end{remark}

%Our result is the following
%\begin{theorem}\label{theorem}
% Under the conditions \eqref{f_{7}}-\eqref{f_{12}} there exists a positive constant $K_{0}$ such that,
% for every $\lambda\in(0,K_0/s]$ problem \eqref{ourproblem} has a radial sign-changing ground state solution that changes sign exactly once in $\mathbb R^3$.
%\end{theorem}
%

%\begin{rem}
%    Note that \eqref{f_{7}}-\eqref{f_{12}} holds true for the function $sf(\cdot)$,with $s$ as in \eqref{f_{10}}, replacing $f(\cdot)$. For this reason, we assume $s=1$ in  and henceforth in our work.
%\end{rem}
%The van Est map. Van Est theorem for ``honest'' coefficients. Van Est theorem for arbitrary coefficients. Integration as a corollary. Applications(?)

%---------------------------------------------------------------------------------------------------------------------------------------------------

\section{The 2-van Est map}

If one starts out with an extension of Lie $2$-groups
\begin{eqnarray*}
\xymatrix{
1 \ar[r] & W \ar[d]_\phi\ar[r]^{j_1} & E_1 \ar[d]_\epsilon\ar[r]^{\pi_1} & G \ar[d]^i\ar[r] & 1  \\
1 \ar[r] & V \ar[r]_{j_0}            & E_0 \ar[r]_{\pi_0}                & H \ar[r]         & 1, 
}
\end{eqnarray*}
one can consider its associated extension of Lie $2$-algebras
\begin{eqnarray*}
\xymatrix{
0 \ar[r] & W \ar[d]_\phi\ar[r]^{d_0j_1} & \mathfrak{e}_1 \ar[d]_{d_1\epsilon}\ar[r]^{d_1\pi_1} & \gg \ar[d]^ \mu\ar[r] & 0  \\
0 \ar[r] & V \ar[r]_{d_0j_0}             & \mathfrak{e}_0 \ar[r]_{d_1\pi_0}                    & \hh \ar[r]         & 0 
}
\end{eqnarray*}
by linearizing, i.e. by differentiating each entry and each map. In sight of theorems \ref{H2Alg} and \ref{H2Gp}, this procedure defines a map 
\begin{eqnarray*}
\xymatrix{
H_\nabla^2(\G,\phi) \ar[r] & H_\nabla^2(\gg_1,\phi).
}
\end{eqnarray*}
In fact, when appropriately relating the $2$-representations by proposition \ref{RepGoToRep}, %*** Do I want it to say also something about the $2$-representations induced by a splitting? ***
this map will be defined at the level of cochains and due to its nature, we are bound to call it a van Est map. We proceed to define it in full generality, not only for $2$-cochains, but for any $n$-cochain. At the level of the triple complexes 
\begin{eqnarray*}
\xymatrix{
\Phi :C^{p,q}_r(\G ,\phi) \ar[r] & C^{p,q}_r(\gg_1,\phi)
}
\end{eqnarray*}
is defined for $\omega\in C^{p,q}_r(\G ,\phi)$, $\xi_1,...,\xi_q\in\gg_p$ and $x_1,...,x_r\in\gg$ by
\begin{eqnarray*}
\Phi\omega(\xi_1,...,\xi_q ;x_1,...,x_r)=\sum_{\sigma\in S_q}\sum_{\varrho\in S_r}\abs{\sigma}\abs{\varrho}R_{\xi_{\sigma(q)}}...R_{\xi_{\sigma(1)}}R_{x_{\varrho(r)}}...R_{x_{\varrho(1)}}\omega .
\end{eqnarray*}
We recall that $\abs{\cdot}$ stands for the sign of a permutation, and for $x\in\gg$ and $\xi\in\gg_p$
\begin{align*}
R_x   & :\xymatrix{C(\G_p^q\times G^r ,W) \ar[r] & C(\G_p^q\times G^{r-1},W)}, \\
R_\xi & :\xymatrix{C(\G_p^q ,W) \ar[r] & C(\G_p^{q-1},W)}
\end{align*}
are given by
\begin{align*}
R_x\omega(\vec{\gamma};\vec{g})     & :=\Lie_{\vec{x}}\omega_{(\vec{\gamma};\vec{g})}(1) \\
R_\xi\varphi(\gamma_2,...,\gamma_q) & :=\Lie_{\vec{\xi}}\varphi_{(\gamma_2,...,\gamma_q)}(1)
\end{align*}
respectively, for $\vec{\gamma}=(\gamma_1,...,\gamma_q)\in\G_p^q$ and $\vec{g}\in G^{r-1}$. Perhaps, it is adequate to spell out the other pieces of notation. First, $\vec{x}\in\XX(G)$ and $\vec{\xi}\in\XX(\G_p)$ stand for the right-invariant vector fields generated by $x$ and $\xi$ respectively. The symbol $\Lie$ stands for the actual Lie derivative this time round; therefore, when applied to the functions 
\begin{align*}
\omega_{(\vec{\gamma};\vec{g})}   & :\xymatrix{G \ar[r] & W :g \ar@{|->}[r] & \omega(\vec{\gamma};g,\vec{g})}, \\
\varphi_{(\gamma_2,...,\gamma_q)} & :\xymatrix{\G_p \ar[r] & W :\gamma \ar@{|->}[r] & \varphi(\gamma,\gamma_2,...,\gamma_q)},
\end{align*}
these derivations yield a second pair of functions, which are later evaluated at the respective units. Notice that these considerations only make sense for $r>0$; nonetheless, for $r=0$, the formulas will have the same shape with the sole difference that the derivative will be taken in $V$ instead of $W$. Indeed, the formulas for the page $r=0$ will coincide column-wise with the classic van Est map from Lie group cochains to Lie algebra cochains. \\
Let us comment on why this map is well-defined. First, notice that these latter $R$ maps are linear in the subscript, because they are defined in terms of vector fields. This will also imply that a succession of $R$'s is going to yield a multi-linear map. Finally, the alternated permutations will make $\Phi\omega$ totally skew-symmetric in the first $q$ variables $\xi_1,...,\xi_q $ and independently in the next $r$ variables $x_1,...,x_r$ as well; thus, $\Phi\omega\in\bigwedge^q\gg_p^*\otimes\bigwedge^r\gg^*\otimes W=C^{p,q}_r(\gg_1,\phi)$. \\
We turn now to show that this is the right Ansatz for a van Est map in that the image of the Lie $2$-group $2$-cochain representing an extension coincides with the Lie $2$-algebra $2$-cochain representing the associated linearized extension, as given at the beginning of this section. In order to do so, let us spell out the $R$ maps and introduce a piece of notation to handle efficiently the upcoming computations. \\
For $\omega\in C^{p,q}_r(\G ,\phi)$, $\xi_1,...,\xi_q\in\gg_p$ and $x_1,...,x_r\in\gg$, we will write 
\begin{eqnarray*}
\overrightarrow{R}_\xi\overrightarrow{R}_x\omega :=R_{\xi_q}...R_{\xi_1}R_{x_r}...R_{x_1}\omega.
\end{eqnarray*}
Now, let $x\in\gg$, $\vec{\gamma}=(\gamma_1,...,\gamma_q)\in\G_p^q$ and $\vec{g}\in G^{r-1}$, and consider 
\begin{align*}
R_x\omega(\vec{\gamma};\vec{g}) & =\Lie_{\vec{x}}\omega_{(\vec{\gamma};\vec{g})}(1) \\
                                & =d_1\omega_{(\vec{\gamma};\vec{g})}(\vec{x}_1) \\
                                & =\frac{d}{d\tau}\rest{\tau=0}\omega_{(\vec{\gamma};\vec{g})}(\exp_G(\tau x)) \\
                                & =\frac{d}{d\tau}\rest{\tau=0}\omega(\vec{\gamma};\exp_G(\tau x),\vec{g}).
\end{align*}
Inductively, one then gets that 
\begin{eqnarray*}
\overrightarrow{R}_x\omega(\vec{\gamma})=\frac{d}{d\tau_r}\rest{\tau_r=0}...\frac{d}{d\tau_1}\rest{\tau_1=0}\omega(\vec{\gamma};\exp_G(\tau_1 x_1),...,\exp_G(\tau_r x_r)) 
\end{eqnarray*}
We introduce the following notation for this kind of multi-derivation
\begin{eqnarray*}
\frac{d^J}{d\tau_J}\rest{\tau =0}:=\frac{d}{d\tau_r}\rest{\tau_r=0}...\frac{d}{d\tau_1}\rest{\tau_1=0},
\end{eqnarray*}
and also,
\begin{eqnarray*}
\exp(\tau\cdot x):=(\exp_G(\tau_1 x_1),...,\exp_G(\tau_r x_r)).
\end{eqnarray*}
Thus, we analogously deduce that
\begin{eqnarray*}
\overrightarrow{R}_\xi\overrightarrow{R}_x\omega =\frac{d}{d\lambda_q}\rest{\lambda_q=0}...\frac{d}{d\lambda_1}\rest{\lambda_1=0}\frac{d^J}{d\tau_J}\rest{\tau =0}\omega(\exp_{\G_p}(\lambda_1\xi_1),...,\exp_{\G_p}(\lambda_q\xi_q);\exp_G(\tau\cdot x)), 
\end{eqnarray*}
and we write
\begin{eqnarray*}
\overrightarrow{R}_\xi\overrightarrow{R}_x\omega =\frac{d^I}{d\lambda_I}\rest{\lambda =0}\frac{d^J}{d\tau_J}\rest{\tau =0}\omega(\exp(\lambda\cdot\xi );\exp(\tau\cdot x)). 
\end{eqnarray*}
Admittedly, this is quite an ambiguous notation, but we hope it will be clear from the context, where the exponentials are taken and what the multi-indices $I$ and $J$ are going to be. Moreover, every time we judge there is room for confusion, we will write the necessary pieces unabbreviated. \\
As announced, let $(\omega_0,\omega_1,\alpha,\hat{\varphi})$ be a Lie $2$-group $2$-cocycle. According to the previous chapter, this cocycle defines a $2$-extension of Lie $2$-groups. We will see that $(\Phi\omega_0,\Phi\omega_1,\Phi\alpha,\Phi\hat{\varphi})$ is the Lie $2$-algebra cocycle associated to the linearized $2$-extension. \\
First of all, recall that as we said, the $2$-van Est map will be a map of complexes given that the $2$-representations are appropriately related. We will write then, $\rho$ for the Lie $2$-group representation and $\dot{\rho}$ for its associated Lie $2$-algebra representation. As we also remarked above, for $r=0$, the $2$-van Est map coincides with the classic van Est map; and hence, $\Phi\omega_0$ is the Lie algebra $2$-cochain defining the Lie algebra of $H{}_{\rho^0_0}\ltimes^{\omega_0}V$ which is itself an extension of $\hh$ by $V$: $\hh_{\dot{\rho}^0_0}\oplus^{\Phi\omega_0}V$. Furthermore, for $q=p=0$ and $r>0$ the formula for the $2$-van Est map coincides with the classic van Est map for the Lie group $G$; hence, $\Phi\omega_1$ also defines the right Lie algebra $2$-cochain, pretty much in the same fashion as $\Phi\omega_0$. \\
Now that we have got the right spaces, we just need to analyze the crossed module structure. Recall from the proof of theorem \ref{H2Gp} that the structural map of the Lie $2$-group defined by the $2$-cocycle $(\omega_0,\omega_1,\alpha,\hat{\varphi})$ is
\begin{eqnarray*}
\xymatrix{G{}_{\rho^1_0\circ i}\ltimes^{\omega_1}W \ar[r] & H{}_{\rho^0_0}\ltimes^{\omega_0}V:(g,w) \ar@{|->}[r] & (i(g),\phi(w)+\hat{\varphi}{\begin{pmatrix}
g \\
1
\end{pmatrix}})};
\end{eqnarray*}
thus, its derivative will be
\begin{eqnarray*}
\xymatrix{\gg{}_{\dot{\rho}^1_0\circ\mu}\oplus^{\Phi\omega_1}W \ar[r] & \hh{}_{\dot{\rho}^0_0}\oplus^{\Phi\omega_0}V:(x,w) \ar@{|->}[r] & (\mu(x),\phi(w)+d_{(1,1)}\hat{\varphi}{\begin{pmatrix}
x \\
0
\end{pmatrix}})}.
\end{eqnarray*}
Of course, 
\begin{eqnarray*}
d_{(1,1)}\hat{\varphi}{\begin{pmatrix}
x \\
0
\end{pmatrix}}=\frac{d}{d\lambda}\rest{\lambda=0}\hat{\varphi}(\exp_{\G}\lambda{\begin{pmatrix}
x \\
0
\end{pmatrix}})=\Phi\hat{\varphi}{\begin{pmatrix}
x \\
0
\end{pmatrix}};
\end{eqnarray*}
thus, we see the structural map of the crossed module of Lie algebras coincide with the one in the proof of theorem \ref{H2Alg}. Finally, we look at the action induced by $\alpha$,
\begin{eqnarray*}
(g,w)^{(h,v)}=(g^h,\rho_0^1(h)^{-1}(w+\rho_1(g)v)+\alpha(h,g)).
\end{eqnarray*}
For each $(h,v)\in H{}_{\rho^0_0}\ltimes^{\omega_0}V$, this latter formula defines a Lie group automorphism 
\begin{eqnarray*}
\xymatrix{(-)^{(h,v)}:G{}_{\rho^1_0\circ i}\ltimes^{\omega_1}W \ar[r] & G{}_{\rho^1_0\circ i}\ltimes^{\omega_1}W}.
\end{eqnarray*}
Differentiating at the identity, one gets a Lie algebra automorphism, for which we will use the same notation
\begin{eqnarray*}
\xymatrix{(-)^{(h,v)}:\gg{}_{\dot{\rho}^1_0\circ\mu}\oplus^{\Phi\omega_1}W \ar[r] & \gg{}_{\dot{\rho}^1_0\circ\mu}\oplus^{\Phi\omega_1}W}.
\end{eqnarray*}
We write down the formula for this action explicitly. First, for $w\in W$,
\begin{align*}
(0,w)^{(h,v)} & =\frac{d}{d\tau}\rest{\tau=0}(1^h,\rho_0^1(h)^{-1}(\tau w+\rho_1(1)v)+\alpha(h,1)) \\
			  & =\frac{d}{d\tau}\rest{\tau=0}(1,\tau \rho_0^1(h)^{-1}(w))=(0,\rho_0^1(h)^{-1}(w)).
\end{align*}
Then, for $x\in\gg$,
\begin{align*}
(x,0)^{(h,v)} & =\frac{d}{d\tau}\rest{\tau=0}(\exp_G(\tau x)^h,\rho_0^1(h)^{-1}(\rho_1(\exp_G(\tau x))v)+\alpha(h,\exp_G(\tau x))) \\
              & =(x^h,\rho_0^1(h)^{-1}(\dot{\rho}_1(x)v)+\frac{d}{d\tau}\rest{\tau=0}\alpha(h,\exp(\tau x))),
\end{align*}
where we abused notation again and wrote $x^h$ for the induced action of $H$ on $\gg$. Put together, 
\begin{eqnarray*}
(x,w)^{(h,v)}=(x^h,\rho_0^1(h)^{-1}(w+\dot{\rho}_1(x)v)+\frac{d}{d\tau}\rest{\tau=0}\alpha(h,\exp(\tau x))).
\end{eqnarray*}
Now, consider the map
\begin{eqnarray*}
\xymatrix{H{}_{\rho^0_0}\ltimes^{\omega_0}V \ar[r] & Aut(\gg{}_{\dot{\rho}^1_0\circ\mu}\oplus^{\Phi\omega_1}W):(h,v) \ar@{|->}[r] & (-)^{(h,v)}},
\end{eqnarray*}
whose differential is by definition the structural action of the crossed module associated to the Lie $2$-algebra of the extension $\xymatrix{G{}_{\rho^1_0\circ i}\ltimes^{\omega_1}W \ar[r] & H{}_{\rho^0_0}\ltimes^{\omega_0}V}$. Then, for $v\in V$,
\begin{align*}
\Lie_{(0,v)}(x,w) & :=\frac{d}{d\lambda}\rest{\lambda=0}(x,w)^{(1,\lambda v)} \\
                  & =\frac{d}{d\lambda}\rest{\lambda=0}(x^1,\rho_0^1(1)^{-1}(w+\dot{\rho}_1(x)(\lambda v))+\frac{d}{d\tau}\rest{\tau=0}\alpha(1,\exp(\tau x)))=(0,\dot{\rho}_1(x)v);
\end{align*}
whereas, on the other hand, for $y\in\hh$,
\begin{align*}
\Lie_{(y,0)}(x,w) & :=\frac{d}{d\lambda}\rest{\lambda=0}(x,w)^{(\exp_H(\lambda y),0)} \\
                  & =\frac{d}{d\lambda}\rest{\lambda=0}(x^{\exp_H(\lambda y)},\rho_0^1(\exp_H(\lambda y))^{-1}(w)+\frac{d}{d\tau}\rest{\tau=0}\alpha(\exp_H(\lambda y),\exp_G(\tau x))) \\
                  & =(\Lie_y x,\dot{\rho}_0^1(-y)w+\frac{d}{d\lambda}\rest{\lambda=0}\frac{d}{d\tau}\rest{\tau=0}\alpha(\exp(\lambda y),\exp(\tau x)))).
\end{align*}
Ultimately, this yields
\begin{eqnarray*}
\Lie_{(y,v)}(x,w)=(\Lie_y x,\dot{\rho}_1(x)v-\dot{\rho}_0^1(y)w+\Phi\alpha(y,x))
\end{eqnarray*}
which is indeed the formula for the action by derivations induced by $\Phi\alpha$. \\

In the sequel, we will prove that $\Phi$ defines a map of triple complexes. It will be educative to consider two cases separately. First, we will deal with the case where the coefficients of the $2$-representation are taken on an honest vector space. This will have the effect of killing off all the $(p,q)$-pages but the front one; hence, to prove that $\Phi$ yields a map of triple complexes in this case, reduces to prove that it induces a map of the double complexes defined by $r=0$. Thereafter, we will deal with the general case.

\subsubsection{Vector space coefficients}

As we mention right before introducing both the triple complex of Lie $2$-algebra cochains and the triple complex of Lie $2$-group cochains, the page $r=0$ does not commute in general. Instead, the images of going around the squares in the two possible ways yield isomorphic elements in the $2$-vector space. Therefore, when taking coefficients on an honest vector space, the first pages commute on the nose, or in other words, are double compĺexes. In order, this implies that we just need to prove that the generic cube
\begin{eqnarray*}
\xymatrix@!0{
& & C(\G_{p}^{q+1},V) \ar[rrrrr]^\Phi\ar[ddll]_\partial & & & & & \bigwedge^{q+1}\gg_{p}^*\otimes V \ar[ddll]_\partial  \\
 \\
C(\G_{p+1}^{q+1},V) \ar[rrrrr]^\Phi & & & & & \bigwedge^{q+1}\gg _{p+1}^*\otimes V  \\
& & C(\G_{p}^{q},V) \ar'[rrr]^\Phi[rrrrr]\ar'[u][uuu]^\delta\ar[ddll]_\partial & & & & & \bigwedge^{q}\gg_{p}^*\otimes V \ar[uuu]_\delta\ar[ddll]^\partial \\
\\
C(\G_{p+1}^{q},V) \ar[rrrrr]^\Phi  \ar[uuu]^\delta & & & & & \bigwedge^{q}\gg_{p+1}^*\otimes V  \ar[uuu]^\delta
}
\end{eqnarray*}
commutes, to verify that $\Phi$ yields a cochain map. 
In fact, as the left and right squares lie in their respective double complexes, they commute. Also, the front and back squares commute, as they are a piece of the classic Van Est map for the Lie group cochains of $\G_{p+1}$ and $\G_p$ respectively. Thus, the only thing left to prove is the following lemma.
\begin{lemma}
In the latter cubic diagram, $\Phi\partial=\partial\Phi$.
\end{lemma} 
\begin{proof}
This is a consequence of the fact that the face maps $\partial_k$ of the simplicial structure of the nerve of the Lie $2$-algebra associated to any given Lie $2$-group are the derivatives of the face maps $\partial_k$ of the simplicial structure of its nerve. In symbols, for $\xi\in\gg_{p+1}$ and for each $k$,
\begin{eqnarray*}
\partial_k\xi =\frac{d}{d\lambda}\rest{\lambda=0}\partial_k(\exp_{\G_{p+1}}(\lambda\xi)),
\end{eqnarray*}
and consequently, 
\begin{eqnarray*}
\exp_{\G_{p}}(\partial_k\xi)=\partial_k(\exp_{\G_{p+1}}(\xi)).
\end{eqnarray*}
Now, let $\omega\in C(\G_p^q,V)$ and $\xi_1,...,\xi_q\in\gg_{p+1}$, then
\begin{align*}
R_{\xi_q}...R_{\xi_1}(\partial\omega) & =\frac{d^I}{d\lambda_I}\rest{\lambda=0}\partial\omega(\exp(\lambda\cdot\xi)) \\
									  & =\frac{d^I}{d\lambda_I}\rest{\lambda=0}\sum_{k=0}^{p+1}(-1)^k\partial_k^*\omega(\exp(\lambda\cdot\xi)) \\
									  & =\sum_{k=0}^{p+1}(-1)^k\frac{d^I}{d\lambda_I}\rest{\lambda=0}\omega\big{(}\partial_k(\exp_{\G_{p+1}}(\lambda_1\xi_1)),...,\partial_k(\exp_{\G_{p+1}}(\lambda_q\xi_q))\big{)} \\
									  & =\sum_{k=0}^{p+1}(-1)^k\frac{d^I}{d\lambda_I}\rest{\lambda=0}\omega(\exp_{\G_{p}}(\lambda_1\partial_k\xi_1),...,\exp_{\G_{p}}(\lambda_q\partial_k\xi_q)) \\
									  & =\sum_{k=0}^{p+1}(-1)^kR_{\partial_k\xi_q}...R_{\partial_k\xi_1}\omega .
\end{align*}
Therefore,
\begin{align*}
\Phi(\partial\omega)(\xi_1,...,\xi_q) & =\sum_{\sigma\in S_q}\abs{\sigma}R_{\xi_{\sigma(q)}}...R_{\xi_{\sigma(1)}}(\partial\omega) \\
									  & =\sum_{\sigma\in S_q}\abs{\sigma}\sum_{k=0}^{p+1}(-1)^kR_{\partial_k\xi_{\sigma(q)}}...R_{\partial_k\xi_{\sigma(1)}}\omega \\
									  & =\sum_{k=0}^{p+1}(-1)^k\sum_{\sigma\in S_q}\abs{\sigma}R_{\partial_k\xi_{\sigma(q)}}...R_{\partial_k\xi_{\sigma(1)}}\omega \\
									  & =\sum_{k=0}^{p+1}(-1)^k\Phi\omega(\partial_k\xi_{1},...,\partial_k\xi_{q})=\partial(\Phi\omega)(\xi_1,...,\xi_q).
\end{align*}
\end{proof}
As claimed, with respect to a $2$-representation with values on an honest vector space, the $2$-van Est map $\Phi$ turned out to be a map of double complexes and we've got the following van Est type theorem.

\begin{theorem}\label{2vE-vs} % vs stands for vector space
Let $\G=\xymatrix{G \ar[r] & H}$ represent a Lie $2$-group with Lie $2$-algebra represented by $\gg_1=\xymatrix{\gg \ar[r] & \hh}$, and $\rho$ be a $2$-representation of $\G$ on an honest vector space $V$. If $H$ is $k$-connected, $G$ is $(k-1)$-connected and $\Phi$ is the $2$-van Est map, then
\begin{eqnarray*}
H^n_{tot}(\Phi)=(0),\quad\textnormal{for all degrees } n\leq k.
\end{eqnarray*}
\end{theorem}
\begin{proof}
In order to compute the cohomology of the mapping cone of the $2$-van Est map, we use the spectral sequence of the double complex filtrated by columns, whose first page is
\begin{eqnarray*}
\xymatrix{
           & \vdots                                        & \vdots                                        & \vdots                           &                      \\ 
           & H^2_{\delta_\Phi}(C^{0,\bullet}(\Phi)) \ar[r] & H^2_{\delta_\Phi}(C^{1,\bullet}(\Phi)) \ar[r] & H^2_{\delta_\Phi}(C^{2,\bullet}(\Phi)) \ar[r] & \dots \\
E^{p,q}_1: & H^1_{\delta_\Phi}(C^{0,\bullet}(\Phi)) \ar[r] & H^1_{\delta_\Phi}(C^{1,\bullet}(\Phi)) \ar[r] & H^1_{\delta_\Phi}(C^{2,\bullet}(\Phi)) \ar[r] & \dots \\
           & H^0_{\delta_\Phi}(C^{0,\bullet}(\Phi)) \ar[r] & H^0_{\delta_\Phi}(C^{1,\bullet}(\Phi)) \ar[r] & H^0_{\delta_\Phi}(C^{2,\bullet}(\Phi)) \ar[r] & \dots 
}
\end{eqnarray*}
Now, the $2$-van Est map is defined column-wise by the classic van Est map; consequently, the columns of the mapping cone double complex coincide with the mapping cones of these. By means of the classic van Est theorem (cf. Theorem \ref{ClassicVanEst}), we know that the $p$th column of the latter diagram is zero \textit{below the connectedness of} $\G_p$. For instance, as $H$ is $k$-connected, the first column is zero below $k$. \\
We use the K\"unneth formula to figure out how connected $\G_p$ is. Since $\G_p\cong G^p\times H$,
\begin{eqnarray*}
H^q_{dR}(\G_p)=\bigoplus_{r+s=q}H^r_{dR}(G^p)\otimes H^s_{dR}(H),
\end{eqnarray*}
but $H$ is $k$-connected; thus yielding,
\begin{eqnarray*}
H^q_{dR}(\G_p)=H^q_{dR}(G^p),
\end{eqnarray*}
for all $q\leq k$. Now, inductively,
\begin{eqnarray*}
H^q_{dR}(G^p)=\bigoplus_{q_1+...+q_p=q}H^{q_1}_{dR}(G)\otimes ...\otimes H^{q_p}_{dR}(G);
\end{eqnarray*}
therefore, since $G$ is $(k-1)$-connected, so is $G^p$. Although it won't be relevant in the sequel, notice that 
\begin{eqnarray*}
H^k_{dR}(G^p)=(H^k_{dR}(G))^{\oplus p},
\end{eqnarray*}
so $(k-1)$-connectedness is as far as we can deduce. In sight of this discussion, $\G_p$ is always $(k-1)$-connected, and all columns in the first page of the spectral sequence above vanish below $k-1$. Given that
\begin{eqnarray*}
E^{p,q}_1\Rightarrow H^{p+q}_{tot}(\Phi),
\end{eqnarray*}
it follows from lemma \ref{BelowDiag} that $H_{tot}^n(\Phi)=(0)$ for $n\leq k$ as desired.

\end{proof}
\begin{cor}
Under the hypothesis of the previous theorem, the $2$-van Est map induces isomorphisms
\begin{eqnarray*}
\xymatrix{
\Phi^n :H_\nabla^n(\G ,V) \ar[r] & H_\nabla^n(\gg_1 ,V),
}
\end{eqnarray*}
for $n\leq k$, and it is injective for $n=k+1$.
\end{cor}
\begin{proof}
This is nothing but the re-phrasing of Proposition \ref{mapConeCoh}.

\end{proof}

\subsubsection{$2$-vector space coefficients}

We move now to the general case. We will break the proof that $\Phi$ defines an appropriate map of triple complexes in several parts. Due to proposition \ref{(q,r)-doubleCx}, we know that each constant $p$-page is a double complex; therefore, $\Phi$ needs to restrict to a map of double complexes when $p$ is left constant. In other words, the generic cube 
\begin{eqnarray*}
\xymatrix@!0{
& & & C(\G_{p}^{q+1}\times G^r,W) \ar[rrrrrr]^\Phi\ar[ddlll]_{\delta_{(1)}} & & & & & & \bigwedge^{q+1}\gg_{p}^*\otimes\bigwedge^{r}\gg^*\otimes W \ar[ddlll]_{\delta_{(1)}}  \\
 \\
C(\G_{p}^{q+1}\times G^{r+1},W) \ar[rrrrrr]^{\Phi\quad} & & & & & & \bigwedge^{q+1}\gg _{p}^*\otimes\bigwedge^{r+1}\gg^*\otimes W  \\
& & & C(\G_{p}^{q}\times G^r,W) \ar'[rrr]^{\qquad\quad\Phi}[rrrrrr]\ar'[u][uuu]_\delta\ar[ddlll]_{\delta_{(1)}} & & & & & & \bigwedge^{q}\gg_{p}^*\otimes\bigwedge^{r}\gg^*\otimes W \ar[uuu]_\delta\ar[ddlll]^{\delta_{(1)}} \\
\\
C(\G_{p}^{q}\times G^{r+1},W) \ar[rrrrrr]^\Phi  \ar[uuu]^\delta & & & & & & \bigwedge^{q}\gg_{p}^*\otimes\bigwedge^{r+1}\gg^*\otimes W  \ar[uuu]^\delta
}
\end{eqnarray*}
needs to commute. Since eventually we are to take alternating sums, we introduce the following partitions of the symmetric group $S_{q+1}$. First,
\begin{eqnarray*}
S_{q+1}=\bigcup_{n=0}^qS_q(m\vert n),
\end{eqnarray*} 			
where $S_q(m\vert n)$ is the set of permutations that fix the $m$th element to be $n$. In symbols,
\begin{eqnarray*}
S_q(m\vert n)=\lbrace\sigma\in S_{q+1}:\sigma(m)=n\rbrace .
\end{eqnarray*}				 
Among these partitions, we are going to be interested in the ones fixing the first element and the ones fixing the last element. Indeed, every element $\sigma$ of $S_q(0\vert j)$ can be factored as
\begin{eqnarray*}
\sigma=\sigma'\sigma_j^0,
\end{eqnarray*}			
where
\begin{eqnarray*}
\sigma_j^0:=\begin{pmatrix}
0 & 1 & 2 & ... & j-1 & j   & j+1 & ... & q-1 & q \\
j & 0 & 1 & ... & j-2 & j-1 & j+1 & ... & q-1 & q
\end{pmatrix},
\end{eqnarray*}	
and $\sigma'$ is the residual permutation that leaves the first element alone and shifts the remaining $q$ elements. In so, we can regard $\sigma'$ as belonging to $Sym(\lbrace 0,1,...,j-1,j+1,...,q\rbrace)\cong S_q$, incidentally justifying the notation. Now, $\sigma_j^0$ is clearly the composition of $j$-many transpositions; therefore, due to the factorization, 
\begin{eqnarray*}
\abs{\sigma}=\abs{\sigma'}\abs{\sigma_j^0}=(-1)^j\abs{\sigma'}.
\end{eqnarray*}	
Analogously, every element $\sigma$ of $S_q(q\vert j)$ can be factored as
\begin{eqnarray*}
\sigma=\sigma'\sigma_j^q,
\end{eqnarray*}			
where
\begin{eqnarray*}
\sigma_j^q:=\begin{pmatrix}
0 & 1 & ... & j-1 & j   & j+1 & ... & q-2 & q-1 & q \\
0 & 1 & ... & j-1 & j+1 & j+2 & ... & q-1 & q   & j
\end{pmatrix},
\end{eqnarray*}	
and $\sigma'$ is the residual permutation that leaves the last element alone and shifts the remaining $q$ elements. We regard $\sigma'$ as belonging to $Sym(\lbrace 0,1,...,j-1,j+1,...,q\rbrace)\cong S_q$ as well. Now, $\sigma_j^q$ is clearly the composition of $(q-j)$-many transpositions; therefore, due to the factorization, this time around we have got
\begin{eqnarray*}
\abs{\sigma}=\abs{\sigma'}\abs{\sigma_j^0}=(-1)^{q-j}\abs{\sigma'}.
\end{eqnarray*}	
The other class of partitions we will be considering is 
\begin{eqnarray*}
S_{q+1}=\bigcup_{m,n}S_{q-1}(j\vert mn),
\end{eqnarray*} 
where $S_{q-1}(j\vert mn)$ is the set of all permutations sending $j$ to $m$ and $j+1$ to $n$. In symbols,
\begin{eqnarray*}
S_{q-1}(j\vert mn)=\lbrace\sigma\in S_{q+1}:\sigma(j)=m,\quad\sigma(j+1)=n\rbrace .
\end{eqnarray*}	
For instance, for $m<n$, each element $\sigma$ of $S_{q-1}(0\vert mn)$ can be factored again as
\begin{eqnarray*}
\sigma=\sigma'\sigma_{mn}^{0},
\end{eqnarray*}			
where $\sigma_{mn}^{0}$ is the permutation
\begin{eqnarray*}
\begin{pmatrix}
0\quad 1\quad 2\quad 3\quad ...\quad m-1\quad m\quad\qquad m+1\quad m+2\quad ...\quad n-1\quad n \quad\qquad n+1\quad ...\quad q-1\quad q \\
m\quad n\quad 0\quad 1\quad ...\quad m-3\quad m-2\quad m-1\quad m+1\quad ...\quad n-2\quad n-1\quad n+1\quad ...\quad q-1\quad q
\end{pmatrix},
\end{eqnarray*}	
and $\sigma'$ is the residual permutation that leaves the first two elements fixed and shifts the remaining $q-1$ elements. We thus regard $\sigma'$ as belonging to the symmetric group on the set of $q+1$ elements without $m$ and $n$. That is, say
\begin{eqnarray*}
\Omega(q):=\lbrace j\in\Nn :0\leq j\leq q\rbrace ,
\end{eqnarray*}
then $\sigma'\in Sym(\Omega(q)\setminus\lbrace m,n\rbrace)\cong S_{q-1}$. To compute the sign of $\sigma\in S_{q-1}(0\vert mn)$, still for $m<n$, notice that $\sigma_{mn}^{0}$ is the product of $\sigma_m^0$ as above, with $(n-1)$-many transpositions; hence, we compute the sign to be
\begin{eqnarray*}
\abs{\sigma}=\abs{\sigma'}\abs{\sigma_{mn}^{0}}=(-1)^{m+n-1}\abs{\sigma'}.
\end{eqnarray*}	
In general, for $m<n$, an element $\sigma$ of $S_{q-1}(j\vert mn)$ will be factored as
\begin{eqnarray*}
\sigma=\sigma'\sigma_{mn}^{j},
\end{eqnarray*}			
where $\sigma_{mn}^{j}$ is the permutation that takes $j$ and $j+1$ to $m$ and $n$ respectively and leaves the rest of the elements ordered increasingly, and $\sigma'$ is the residual permutation that shuffles around the remaining $q-1$ elements. As before, we regard $\sigma'$ as belonging to $Sym(\Omega(q)\setminus\lbrace m,n\rbrace)\cong S_{q-1}$, and we compute the sign of $\sigma$ by noticing that $\sigma_{mn}^j$ is the product of the transpositions necessary to take $j$ to $m$ and the ones to take $j+1$ to $n$ 
\begin{eqnarray*}
\abs{\sigma}=\abs{\sigma'}\abs{\sigma_{mn}^{0}}=(-1)^{m-j}(-1)^{n-j-1}\abs{\sigma'}.
\end{eqnarray*}	
We will also need the following lemma.
\begin{lemma}\label{Multilin}
If $\xymatrix{R:V\times ...\times V \ar[r] & W}$ is an $r$-multilinear map and $H_\lambda$ is a differentiable path of automorphisms of $V$ with $H_0=Id_V$, then
\begin{eqnarray*}
\frac{d}{d\lambda}\rest{\lambda=0}R(H_\lambda(v_1),...,H_\lambda(v_r))=\sum_{k=1}^rR(v_1,...,v_{k-1},\frac{d}{d\lambda}\rest{\lambda=0}H_{\lambda}(v_k),v_{k+1},...,v_r)
\end{eqnarray*}
\end{lemma}
\begin{proof}
Let $\lbrace e_i\rbrace_{i=1}^n$ be a basis for $V$. During this proof, we will use Einstein's summation convention. In these coordinates, set
\begin{eqnarray*}
R_{a_1...a_r}:=R(e_{a_1},...,e_{a_r}) & \textnormal{and} & H_\lambda(e_a)=H_a^b(\lambda)e_b.
\end{eqnarray*}
Notice that since $H_0=Id_V$, $H_a^b(0)=\delta_a^b$. Clearly, it suffices to prove the equation for basic elements, then consider
\begin{align*}
\frac{d}{d\lambda}\rest{\lambda=0}R(H_\lambda(e_{a_1}),...,H_\lambda(e_{a_r})) & =\frac{d}{d\lambda}\rest{\lambda=0}R(H_{a_1}^{b_1}(\lambda)e_{b_1},...,H_{a_r}^{b_r}(\lambda)e_{b_r}) \\
 & =\frac{d}{d\lambda}\rest{\lambda=0}H_{a_1}^{b_1}(\lambda)...H_{a_r}^{b_r}(\lambda)R_{b_1...b_r} \\
 & =R_{b_1...b_r}\sum_{k=1}^rH_{a_1}^{b_1}(0)...H_{a_{k-1}}^{b_{k-1}}(0)\dot{H}_{a_k}^{b_k}(0)H_{a_{k+1}}^{b_{k+1}}(0)...H_{a_r}^{b_r}(0) \\
 & =\sum_{k=1}^rR_{b_1...b_r}\delta_{a_1}^{b_1}...\delta_{a_{k-1}}^{b_{k-1}}\dot{H}_{a_k}^{b_k}(0)\delta_{a_{k+1}}^{b_{k+1}}...\delta_{a_r}^{b_r} \\
 & =\sum_{k=1}^rR_{a_1...a_{k-1}b_ka_{k+1}...a_r}\dot{H}_{a_k}^{b_k}(0) \\
  & =\sum_{k=1}^rR(e_{a_1},...,e_{a_{k-1}},\frac{d}{d\lambda}\rest{\lambda=0}H_{a_k}^{b_k}(\lambda)e_{b_k},e_{a_{k+1}},...,e_{a_r}). 
\end{align*}
\end{proof}		
Since, yet again, both right and left faces of the generic cube commute right away, as they belong to their respective double complexes, the fact that $\Phi$ restricts to a map of double complexes for constant $p$ will follow from the sequence of lemmas ahead.		 
\begin{lemma}
The back and front faces commute: $\Phi\delta=\delta\Phi$.
\end{lemma}
\begin{proof}
In the sequel, let $\omega\in C(\G_p^q\times G^r,W)$, $\xi_0,...,\xi_{q}\in\gg_{p}$ and $x_1,...,x_{r}\in\gg$. Consider
\begin{eqnarray*}
R_{\xi_q}...R_{\xi_{0}}R_{x_r}...R_{x_{1}}(\delta\omega) & =\frac{d^I}{d\lambda_I}\rest{\lambda=0}\frac{d^J}{d\tau_J}\rest{\tau=0}\delta\omega(\exp(\lambda\cdot\xi);\exp(\tau\cdot x)). 
\end{eqnarray*}
We call the three types of term from $\delta\omega$ as follows: 
\begin{align*}
I   & :=\frac{d^I}{d\lambda_I}\rest{\lambda=0}\frac{d^J}{d\tau_J}\rest{\tau=0}\omega(\exp(\lambda_1\xi_1),...,\exp(\lambda_{q}\xi_{q});\exp(\tau\cdot x)^{t_p(\exp(\lambda_0\xi_0))}) \\
II  & :=\frac{d^I}{d\lambda_I}\rest{\lambda=0}\frac{d^J}{d\tau_J}\rest{\tau=0}\omega(\exp(\lambda_0\xi_0),...,\exp(\lambda_{j-1}\xi_{j-1})\exp(\lambda_j\xi_j),...,\exp(\lambda_q\xi_q);\exp(\tau\cdot x)) \\
III & :=\frac{d^I}{d\lambda_I}\rest{\lambda=0}\frac{d^J}{d\tau_J}\rest{\tau=0}\rho_0^1(t_p(\exp(\lambda_q\xi_q)))^{-1}\omega(\exp(\lambda_0\xi_0),...,\exp(\lambda_{q-1}\xi_{q-1});\exp(\tau\cdot x)).
\end{align*}
In the first type, $(g_1,...,g_r)^h:=(g_1^h,...,g_r^h)$ for $g_1,...,g_r\in G$ and $h\in H$. By the same remark that the face maps of the simplicial structure of a Lie $2$-group and those of the one of its Lie $2$-algebra are related by derivation, we conclude that for $\xi\in\gg_p$,
\begin{eqnarray*}
\exp_H(\hat{t}_p(\xi))=t_p(\exp_{\G_p}(\xi)).
\end{eqnarray*}
Hence, we get that
\begin{align*}
III & =\frac{d^{I}}{d\lambda_{I}}\rest{\lambda=0}\frac{d^J}{d\tau_J}\rest{\tau=0}\rho_0^1(\exp(\lambda_q\hat{t}_p(\xi_q))^{-1})\omega(\exp(\lambda_0\xi_0),...,\exp(\lambda_{q-1}\xi_{q-1});\exp(\tau\cdot x)) \\
  & =\frac{d}{d\lambda_q}\rest{\lambda_q=0}\rho_0^1(\exp(\lambda_q\hat{t}_p(-\xi_q)))R_{\xi_{q-1}}...R_{\xi_{0}}R_{x_r}...R_{x_{1}}\omega \\
  & =-\dot{\rho}_0^1(\hat{t}_p(\xi_q))R_{\xi_{q-1}}...R_{\xi_{0}}R_{x_r}...R_{x_{1}}\omega .
\end{align*}
Considering the alternating sum of terms of this type and using the partition of $S_{q+1}$ by $S_q(q\vert j)$'s, we get
\begin{align*}
\sum_{\sigma\in S_{q+1}} & \sum_{\varrho\in S_r} -\abs{\sigma}\abs{\varrho}\dot{\rho}_0^1(\hat{t}_p(\xi_{\sigma(q)}))R_{\xi_{\sigma(q-1)}}...R_{\xi_{\sigma(0)}}R_{x_{\varrho(r)}}...R_{x_{\varrho(1)}}\omega \\
						 & =\sum_{j=0}^q(-1)^{q-j+1}\dot{\rho}_0^1(\hat{t}_p(\xi_{j}))\sum_{\sigma'\in S_{q}(q\vert j)}\sum_{\varrho\in S_r}\abs{\sigma'}\abs{\varrho}R_{\xi_{\sigma'(q-1)}}...R_{\xi_{\sigma'(j+1)}}R_{\xi_{\sigma'(j-1)}}...R_{\xi_{\sigma'(0)}}R_{x_{\varrho(r)}}...R_{x_{\varrho(1)}}\omega \\
						 & =\sum_{j=0}^q(-1)^{q-j+1}\dot{\rho}_0^1(\hat{t}_p(\xi_{j}))\Phi\omega(\xi_0,...,\xi_{j-1},\xi_{j+1},...,\xi_q;x_1,...,x_r).
\end{align*}
Next, we consider the terms of type $I$,
\begin{align*}
I & =\frac{d^I}{d\lambda_I}\rest{\lambda=0}\frac{d^J}{d\tau_J}\rest{\tau=0}\omega(\exp(\lambda_1\xi_1),...,\exp(\lambda_{q}\xi_{q});\exp(\tau\cdot x)^{\exp(\lambda_0\hat{t}_p(\xi_0))}) \\
    & =\frac{d^I}{d\lambda_I}\rest{\lambda=0}R_{x_r^{\exp(\lambda_0\hat{t}_p(\xi_0))}}...R_{x_1^{\exp(\lambda_0\hat{t}_p(\xi_0))}}\omega(\exp(\lambda_1\xi_1),...,\exp(\lambda_{q}\xi_{q})) .
\end{align*}
Since, for fixed $\vec{\gamma}\in\G_p^q$,
\begin{align*}
R_{x_r}...R_{x_1}\omega(\vec{\gamma})
\end{align*}
is $r$-multilinear in the $x$ variables and for $y\in\hh$, $(-)^{\exp(\lambda y)}$ is a differentiable path of automorphisms through the identity, invoking lemma \ref{Multilin}, we have got that 
\begin{align*}
\frac{d}{d\lambda}\rest{\lambda=0}R_{x_r^{\exp(\lambda y)}}...R_{x_1^{\exp(\lambda y)}}\omega(\vec{\gamma}) & =\sum_{k=1}^r R_{x_r^{\exp(0)}}...R_{x_{k+1}^{\exp(0)}}R_{\frac{d}{d\lambda}\vert_{\lambda=0}x_k^{\exp(\lambda y)}}R_{x_{k-1}^{\exp(0)}}...R_{x_1^{\exp(0)}}\omega(\vec{\gamma}) \\
 & =\sum_{k=1}^r R_{x_r}...R_{x_{k+1}}R_{\Lie_y x_k}R_{x_{k-1}}...R_{x_1}\omega(\vec{\gamma}). 
\end{align*}
Consequently,
\begin{align*}
\frac{d}{d\lambda_0}\rest{\lambda_0=0} & R_{x_r^{\exp(\lambda_0\hat{t}_p(\xi_0))}}...R_{x_1^{\exp(\lambda_0\hat{t}_p(\xi_0))}}\omega(\exp(\lambda_1\xi_1),...,\exp(\lambda_{q}\xi_{q}))= \\
                     &\sum_{k=1}^r R_{x_r}...R_{x_{k+1}}R_{\Lie_{\hat{t}_p(\xi_0)} x_k}R_{x_{k-1}}...R_{x_1}\omega(\exp(\lambda_1\xi_1),...,\exp(\lambda_{q}\xi_{q})),
\end{align*}
and
\begin{eqnarray*}
I=\sum_{k=1}^r R_{\xi_q}...R_{\xi_{1}}R_{x_r}...R_{x_{k+1}}R_{\Lie_{\hat{t}_p(\xi_0)} x_k}R_{x_{k-1}}...R_{x_1}\omega .
\end{eqnarray*}
As for the alternating sum of terms of this type, we use the partition of $S_{q+1}$ by $S_q(0\vert j)$'s; thus yielding,
\begin{align*}
& \sum_{\sigma\in S_{q+1}}\sum_{\varrho\in S_r}\abs{\sigma}\abs{\varrho}\sum_{k=1}^r R_{\xi_{\sigma(q)}}...R_{\xi_{\sigma(1)}}R_{x_{\varrho(r)}}...R_{x_{\varrho(k+1)}}R_{\Lie_{\hat{t}_p(\xi_{\sigma(0)})} x_{\varrho(k)}}R_{x_{\varrho(k-1)}}...R_{x_{\varrho(1)}}\omega \\
& =\sum_{j=0}^q(-1)^{j}\sum_{k=1}^r\sum_{\sigma'\in S_{q}(0\vert j)}\sum_{\varrho\in S_r}\abs{\sigma'}\abs{\varrho}R_{\xi_{\sigma'(q)}}...R_{\xi_{\sigma'(0)}}R_{x_{\varrho(r)}}...R_{x_{\varrho(k+1)}}R_{\Lie_{\hat{t}_p(\xi_j)} x_{\varrho(k)}}R_{x_{\varrho(k-1)}}...R_{x_{\varrho(1)}}\omega \\
& =\sum_{j=0}^q(-1)^{j}\sum_{k=1}^r\Phi\omega(\xi_0,...,\xi_{j-1},\xi_{j+1},...,\xi_q;x_1,...,\Lie_{\hat{t}_p(\xi_j)} x_k,...,x_r).
\end{align*}
Since, the terms of type $I$ come together with a $(-1)^{q+1}$ coefficient in $\delta\omega$, putting together the two types of terms we have analyzed so far, we get
\begin{eqnarray*}
\sum_{j=0}^q(-1)^{j+1}\rho'_{\xi_{j}}(\Phi\omega)(\xi_0,...,\xi_{j-1},\xi_{j+1},...,\xi_q;x_1,...,x_r).
\end{eqnarray*}
In order to move forward with the proof, let us remark that by definition, for $\varphi\in C(\G_p,W)$, 
\begin{eqnarray*}
R_{[\xi_1,\xi_2]}\varphi :=\Lie_{\overrightarrow{[\xi_1,\xi_2]}}\varphi(1).
\end{eqnarray*}				
By construction, and since the space of right-invariant vector fields is a Lie subalgebra of $\XX(\G_p)$,
\begin{eqnarray*}
\Lie_{\overrightarrow{[\xi_1,\xi_2]}}\varphi(1)=\Lie_{[\vec{\xi}_1,\vec{\xi}_2]}\varphi(1)=[\Lie_{\vec{\xi}_1},\Lie_{\vec{\xi}_2}]\varphi(1).
\end{eqnarray*}
Computing,
\begin{align*}
\Lie_{\vec{\xi}_1}\Lie_{\vec{\xi}_2}\varphi(1) & =d_1(\Lie_{\vec{\xi}_2}\varphi)(\xi_1) \\
                                               & =\frac{d}{d\lambda_1}\rest{\lambda_1=0}(\Lie_{\vec{\xi}_2}\varphi)(\exp(\lambda_1\xi_1)) \\
                                               & =\frac{d}{d\lambda_1}\rest{\lambda_1=0}d_{\exp(\lambda_1\xi_1)}\varphi((\vec{\xi}_2)_{\exp(\lambda_1\xi_1)}). 
\end{align*}
By definition, 
\begin{eqnarray*}
(\vec{\xi}_2)_{\exp(\lambda_1\xi_1)}=dR_{\exp(\lambda_1\xi_1)}\xi_2 
\end{eqnarray*}
and as always,
\begin{eqnarray*}
\exp(dR_{\exp(\lambda_1\xi_1)}(\lambda_2\xi_2))=R_{\exp(\lambda_1\xi_1)}(\exp(\lambda_2\xi_2));
\end{eqnarray*}
thus,
\begin{eqnarray*}
\Lie_{\vec{\xi}_1}\Lie_{\vec{\xi}_2}\varphi(1) & =\frac{d}{d\lambda_1}\rest{\lambda_1=0}\frac{d}{d\lambda_2}\rest{\lambda_2=0}\varphi(\exp(\lambda_2\xi_2)\exp(\lambda_1\xi_1)),
\end{eqnarray*}
ultimately implying that 
\begin{eqnarray*}
R_{[\xi_1,\xi_2]}\varphi & =\frac{d}{d\lambda_1}\rest{\lambda_1=0}\frac{d}{d\lambda_2}\rest{\lambda_2=0}\varphi(\exp(\lambda_2\xi_2)\exp(\lambda_1\xi_1))-\frac{d}{d\lambda_2}\rest{\lambda_2=0}\frac{d}{d\lambda_1}\rest{\lambda_1=0}\varphi(\exp(\lambda_1\xi_1)\exp(\lambda_2\xi_2)).
\end{eqnarray*}
We now carefully consider the alternating sums of elements of type $II$. First, let $j=1$ and consider the transposition $(0\quad 1)$, which obviously comes with a minus sign, then, in sight of the previous discussion, the sum
\begin{align*}
\frac{d^I}{d\lambda_I} & \rest{\lambda=0}\frac{d^J}{d\tau_J}\rest{\tau=0}\omega(\exp(\lambda_{0}\xi_{0})\exp(\lambda_1\xi_1),\exp(\lambda_2\xi_2),...,\exp(\lambda_q\xi_q);\exp(\tau\cdot x)) \\
  & -\frac{d^I}{d\lambda_I}\rest{\lambda=0}\frac{d^J}{d\tau_J}\rest{\tau=0}\omega(\exp(\lambda_{1}\xi_{1})\exp(\lambda_0\xi_0),\exp(\lambda_2\xi_2),...,\exp(\lambda_q\xi_q);\exp(\tau\cdot x))
\end{align*}
equals
\begin{eqnarray*}
R_{\xi_q}...R_{\xi_2}R_{[\xi_1,\xi_0]}R_{x_r}...R_{x_1}\omega .
\end{eqnarray*}
Let $m<n$, summing elements of type $II$ with $j=1$ over the union of $S_{q-1}(0\vert mn)$ and $S_{q-1}(0\vert nm)$ and well-aware that the elements in the latter come with a sign from the transposition of $m$ and $n$, we get
\begin{eqnarray*}
\sum_{\sigma'\in S_{q-1}(0\vert mn)}\sum_{\varrho\in S_r}(-1)^{m+n}\abs{\sigma'}\abs{\varrho}\overrightarrow{R}_{\sigma'(\xi(m,n))}R_{[\xi_m,\xi_n]}\overrightarrow{R}_{\varrho(x)}\omega ,
\end{eqnarray*}
where we used the following short-hands
\begin{align*}
\sigma'(\xi(m,n)) & :=(\xi_{\sigma'(q)},\xi_{\sigma'(q-1)},...,\xi_{\sigma'(n+1)},\xi_{\sigma'(n-1)},...,\xi_{\sigma'(m+1)},\xi_{\sigma'(m-1)},...,\xi_{\sigma'(1)},\xi_{\sigma'(0)}), \\
\varrho(x)        & :=(x_{\varrho(r)},...,x_{\varrho(1)}).
\end{align*}
Following this line of thought, still with $m<n$, the sum over $S_{q-1}(j\vert mn)\cup S_{q-1}(j\vert nm)$ of elements of type $II$ for a fixed $j$ is
\begin{eqnarray*}
\sum_{\sigma'\in S_{q-1}(j\vert mn)}\sum_{\varrho\in S_r}(-1)^{m+n}\abs{\sigma'}\abs{\varrho}\overrightarrow{R}_{\sigma'(\xi(m,n))_1^j}R_{[\xi_m,\xi_n]}\overrightarrow{R}_{\sigma'(\xi(m,n))_0^j}\overrightarrow{R}_{\varrho(x)}\omega ,
\end{eqnarray*}
where this time around, $\sigma'(\xi(m,n))_0^j$ and $\sigma'(\xi(m,n))_1^j$ stand for the two pieces that come out when cutting $\sigma'(\xi(m,n))$ after the $j$th term, i.e. $\sigma'(\xi(m,n))_0^j$ is the string with the first $j$-many terms of $\sigma'(\xi(m,n))$ and $\sigma'(\xi(m,n))_1^j$ is the string with the remaining $q-1-j$ terms.
Thus, considering the alternating sum of all elements of type $II$ with their respective coefficients from $\delta\omega$, we get
\begin{align*}
\sum_{m<n} & (-1)^{m+n}\sum_{j=0}^{q-1}\sum_{\sigma'\in S_{q-1}}\sum_{\varrho\in S_r}(-1)^{j+1}\abs{\sigma'}\abs{\varrho}\overrightarrow{R}_{\sigma'(\xi(m,n))_1^j}R_{[\xi_m,\xi_n]}\overrightarrow{R}_{\sigma'(\xi(m,n))_0^j}\overrightarrow{R}_{\varrho(x)}\omega \\
 & =\sum_{m<n}(-1)^{m+n}\Phi\omega([\xi_m,\xi_n],\xi_{0},...,\xi_{m-1},\xi_{m+1},...,\xi_{n-1},\xi_{n+1},...,\xi_{q};x_1,...,x_r),
\end{align*}
and the result follows.
						 
\end{proof}						 
\begin{lemma}
The top and bottom faces commute: $\Phi\delta_{(1)}=\delta_{(1)}\Phi$.
\end{lemma}
\begin{proof}
In the sequel, let $\omega\in C(\G_p^q\times G^r,W)$, $\xi_1,...,\xi_{q}\in\gg_{p}$ and $x_0,...,x_{r}\in\gg$. Consider
\begin{eqnarray*}
R_{\xi_q}...R_{\xi_{1}}R_{x_r}...R_{x_{0}}(\delta_{(1)}\omega) & =\frac{d^I}{d\lambda_I}\rest{\lambda=0}\frac{d^J}{d\tau_J}\rest{\tau=0}\delta_{(1)}\omega(\exp(\lambda\cdot\xi);\exp(\tau\cdot x)). 
\end{eqnarray*}
We call the three types of term from $\delta_{(1)}\omega$ as follows: 
\begin{align*}
I   & :=\frac{d^I}{d\lambda_I}\rest{\lambda=0}\frac{d^J}{d\tau_J}\rest{\tau=0}\rho_0^1(i(\exp(\tau_0x_0)^{t_p(\prod\exp(\lambda\cdot\xi))}))\omega(\exp(\lambda\cdot\xi);\exp(\tau_1x_1),...,\exp(\tau_rx_r)) \\
II  & :=\frac{d^I}{d\lambda_I}\rest{\lambda=0}\frac{d^J}{d\tau_J}\rest{\tau=0}\omega(\exp(\lambda\cdot\xi);\exp(\tau_0x_0),...,\exp(\tau_{k}x_k)\exp(\tau_{k+1}x_{k+1}),...,\exp(\tau_rx_r)) \\
III & :=\frac{d^I}{d\lambda_I}\rest{\lambda=0}\frac{d^J}{d\tau_J}\rest{\tau=0}\omega(\exp(\lambda\cdot\xi);\exp(\tau_0x_0),...,\exp(\tau_{r-1}x_{r-1})).
\end{align*}
This time, we start by considering the terms of type $III$; indeed, since
\begin{align*}
\frac{d}{d\tau_r}\rest{\tau_r=0}\omega(\exp(\lambda\cdot\xi);\exp(\tau_0x_0),...,\exp(\tau_{r-1}x_{r-1})) =0,
\end{align*}
we can disregard them completely. \\
As for the terms of type $II$, the considerations of the proof of the previous lemma carry over; thus implying that the alternating sum of all of them, ranging over all values $0\leq k\leq r-1$, is
\begin{align*}
\sum_{m<n} & (-1)^{m+n}\sum_{k=0}^{r-1}\sum_{\sigma\in S_{q}}\sum_{\varrho'\in S_{r-1}}(-1)^{k+1}\abs{\sigma}\abs{\varrho'}\overrightarrow{R}_{\sigma(\xi)}\overrightarrow{R}_{\varrho'(x(m,n))_1^j}R_{[x_m,x_n]}\overrightarrow{R}_{\varrho'(x(m,n))_0^j}\omega \\
 & =\sum_{m<n}(-1)^{m+n}\Phi\omega(\xi_{1},...,\xi_{q};[x_m,x_n],x_0,...,x_{m-1},x_{m+1},...,x_{n-1},x_{n+1},...,x_r),
\end{align*}
where we used the notation from the previous proof as well.  \\%Notice that there is a sign missing. This is due to the fact that the products in the argument of the terms of type $II$ have reversed order compared to the ones in the previous lemma. 
To conclude the proof, we make the following observations. First, since $t_p$ is a composition of maps in the simplicial structure of a Lie $2$-group, it is a homomorphism; and therefore,
\begin{eqnarray*}
t_p\Big{(}\prod\exp(\lambda\cdot\xi)\Big{)}=\prod t_p(\exp(\lambda\cdot\xi)).
\end{eqnarray*}
Also, as before, for $\xi\in\gg_p$,
\begin{eqnarray*}
t_p(\exp_{\G_p}(\xi))=\exp_H(\hat{t}_p(\xi));
\end{eqnarray*}
thereby yielding,
\begin{eqnarray*}
t_p\Big{(}\prod\exp(\lambda\cdot\xi)\Big{)}=\prod\exp(\lambda\cdot\hat{t}_p(\xi)).
\end{eqnarray*}
On the other hand, for all $g\in G$ and $h\in H$, the equivariance of the structural map of the crossed module $i$ and the fact that $\rho_0^1$ is a representation makes
\begin{eqnarray*}
\rho_0^1(i(g^h))=\rho_0^1(h^{-1}i(g)h)=\rho_0^1(h^{-1})\rho_0^1(i(g))\rho_0^1(h).
\end{eqnarray*}
Differentiating this expression (twice) will yield the commutator of $\ggl(W)$. We use this to get to a formula for the terms of type $I$. For $y\in\hh$ and $x\in\gg$, consider
\begin{align*}
\frac{d}{d\lambda}\rest{\lambda=0}\frac{d}{d\tau}\rest{\tau=0}\rho_0^1(i(\exp(\tau x)^{\exp(\lambda y)})) & =\frac{d}{d\lambda}\rest{\lambda=0}\dot{\rho}_0^1(\mu(x^{\exp(\lambda y)})).
\end{align*}
Using the formula instead,
\begin{align*}
\frac{d}{d\lambda}\rest{\lambda=0}\frac{d}{d\tau}\rest{\tau=0}\rho_0^1(i(\exp(\tau x)^{\exp(\lambda y)})) & =\frac{d}{d\lambda}\rest{\lambda=0}\frac{d}{d\tau}\rest{\tau=0}\rho_0^1(\exp(\lambda y)^{-1})\rho_0^1(i(\exp(\tau x)))\rho_0^1(\exp(\lambda y)) \\
 & =\frac{d}{d\lambda}\rest{\lambda=0}\rho_0^1(\exp(-\lambda y))\dot{\rho}_0^1(\mu(x))\rho_0^1(\exp(\lambda y)) \\
 & =-\dot{\rho}_0^1(y)\dot{\rho}_0^1(\mu(x))\cancelto{I}{\rho_0^1(\exp(0))}+\cancelto{I}{\rho_0^1(\exp(0))}\dot{\rho}_0^1(\mu(x))\dot{\rho}_0^1(y) \\
 & =[\dot{\rho}_0^1(\mu(x)),\dot{\rho}_0^1(y)]=\dot{\rho}_0^1([\mu(x),y]);
\end{align*}
and thus, using the equivariance of $\mu$,
\begin{align*}
\frac{d}{d\lambda}\rest{\lambda=0}\dot{\rho}_0^1(\mu(x^{\exp(\lambda y)}))=-\dot{\rho}_0^1(\mu(\Lie_yx)).
\end{align*}
Hence, we get that, 
\begin{align*}
I & =\frac{d^I}{d\lambda_I}\rest{\lambda=0}\frac{d}{d\tau_0}\rest{\tau_0=0}\rho_0^1(i(\exp(\tau_0x_0)^{\prod\exp(\lambda\cdot\hat{t}_p(\xi))}))R_{x_r}...R_{x_1}\omega(\exp(\lambda\cdot\xi)) \\
  & =\frac{d^I}{d\lambda_I}\rest{\lambda=0}\dot{\rho}_0^1(\mu(x_0^{\prod\exp(\lambda\cdot\hat{t}_p(\xi))}))R_{x_r}...R_{x_1}\omega(\exp(\lambda\cdot\xi)).
\end{align*}
Setting $h_j\in H$ to be $h_j:=\exp(\lambda_{j}\hat{t}_p(\xi_{j}))...\exp(\lambda_q\hat{t}_p(\xi_q))$,
\begin{align*}
\frac{d}{d\lambda_1} & \rest{\lambda_1=0}\dot{\rho}_0^1(\mu(x_0^{h_1}))R_{x_r}...R_{x_1}\omega(\exp(\lambda_1\xi_1),...,\exp(\lambda_q\xi_q)) \\
  & =-\dot{\rho}_0^1(\mu((\Lie_{\hat{t}_p(\xi_1)}x_0)^{h_2}))R_{x_r}...R_{x_1}\omega(\exp(0),\exp(\lambda_{2}\xi_{2}),...,\exp(\lambda_q\xi_q))+ \\
  & \qquad\qquad\qquad +\dot{\rho}_0^1(\mu(x_0^{\exp(0)h_2}))R_{\xi_1}R_{x_r}...R_{x_1}\omega(\exp(\lambda_{2}\xi_{2}),...,\exp(\lambda_q\xi_q)) \\
  & =\dot{\rho}_0^1(\mu(x_0^{h_2}))R_{\xi_1}R_{x_r}...R_{x_1}\omega(\exp(\lambda_{2}\xi_{2}),...,\exp(\lambda_q\xi_q));
\end{align*}
here, the last equation follows from the fact that $R_{x_1}...R_{x_r}\omega\in C(\G_p^q,W)$. Inductively,
\begin{align*}
I & =\frac{d}{d\lambda_q}\rest{\lambda_q=0}...\frac{d}{d\lambda_j}\rest{\lambda_j=0}\dot{\rho}_0^1(\mu(x_0^{h_{j}}))R_{\xi_{j-1}}...R_{\xi_1}R_{x_r}...R_{x_1}\omega(\exp(\lambda_j\xi_j),...,\exp(\lambda_q\xi_q)) \\
  & = ... =\frac{d}{d\lambda_q}\rest{\lambda_q=0}\dot{\rho}_0^1(\mu(x_0^{\exp(\lambda_q\hat{t}_p(\xi_q))}))R_{\xi_{q-1}}...R_{\xi_1}R_{x_r}...R_{x_1}\omega(\exp(\lambda_q\xi_q)) \\
  & =\dot{\rho}_0^1(\mu(x_0))R_{\xi_{q}}...R_{\xi_1}R_{x_r}...R_{x_1}\omega
\end{align*}
Considering the alternating sum of terms of this type and using the partition of $S_{r+1}$ by $S_r(0\vert k)$'s, we get
\begin{align*}
\sum_{\sigma\in S_{q}} & \sum_{\varrho\in S_{r+1}}\abs{\sigma}\abs{\varrho}\dot{\rho}_0^1(\mu(x_{\varrho(0)}))R_{\xi_{\sigma(q)}}...R_{\xi_{\sigma(1)}}R_{x_{\varrho(r)}}...R_{x_{\varrho(1)}}\omega \\
						 & =\sum_{k=0}^r(-1)^{k}\dot{\rho}_0^1(\mu(x_{k}))\sum_{\sigma\in S_{q}}\sum_{\varrho'\in S_r(0\vert k)}\abs{\sigma}\abs{\varrho'}R_{\xi_{\sigma(q)}}...R_{\xi_{\sigma(1)}}R_{x_{\varrho'(r)}}...R_{x_{\varrho'(k+1)}}R_{x_{\varrho'(k-1)}}...R_{x_{\varrho'(1)}}\omega \\
						 & =\sum_{k=0}^r(-1)^{k}\dot{\rho}_0^1(\mu(x_{k}))\Phi\omega(\xi_1,...,\xi_q;x_0,...,x_{k-1},x_{k+1},...,x_r),
\end{align*}
and the result follows.
						 
\end{proof}

\begin{cor}
For $p$ constant, $\Phi$ restricts to a map of double complexes.
\end{cor}
\begin{proof}
Thanks to the previous lemmas, together with the fact that for $r=0$, $\Phi$ coincides with the classic van Est map, we just need to prove that 
\begin{eqnarray*}
\xymatrix{
C(\G_p^q\times G,W) \ar[r]^\Phi            & \bigwedge^q\gg_p^*\otimes\gg^*\otimes W \\
C(\G_p^q ,V) \ar[u]^{\partial'}\ar[r]_\Phi & \bigwedge^q\gg_p^*\otimes V\ar[u]_{\delta_{(1)}}
}
\end{eqnarray*}
commutes. For $q=0$, $v\in V$ and $x\in\gg$,
\begin{align*}
\Phi\delta_{(1)}v(x) & =\frac{d}{d\tau}\rest{\tau=0}\delta_{(1)}v(\exp(\tau x)) =\frac{d}{d\tau}\rest{\tau=0}\rho_1(\exp(\tau x))v=\dot{\rho}_1(x)v=\delta_{(1)}v(x),
\end{align*}
and $\Phi$ is defined to be the identity for $(p,q,r)=(0,0,0)$. For any other value of $q$, $\omega\in C(\G_p^q ,V)$ and $\xi_1,...,\xi_q\in\gg_p$,
\begin{align*}
R_{\xi_q}...R_{\xi_1}R_x\partial'\omega & =\frac{d^I}{d\lambda_I}\rest{\lambda=0}\frac{d}{d\tau}\rest{\tau=0}\partial'\omega(\exp(\lambda\cdot\xi);\exp(\tau x)) \\
										& =\frac{d^I}{d\lambda_I}\rest{\lambda=0}\frac{d}{d\tau}\rest{\tau=0}\rho_0^1(t_p(\prod\exp(\lambda\cdot\xi)))^{-1}\rho_1(\exp(\tau x))\omega(\exp(\lambda\cdot\xi)) \\
										& =\frac{d^I}{d\lambda_I}\rest{\lambda=0}\rho_0^1(\prod\exp(\lambda\cdot\hat{t}_p(\xi)))^{-1}\dot{\rho}_1(x)\omega(\exp(\lambda\cdot\xi)) .										
\end{align*}
Now, using the convention again that $h_j\in H$ is $h_j:=\exp(\lambda_{j}\hat{t}_p(\xi_{j}))...\exp(\lambda_q\hat{t}_p(\xi_q))$ and computing
\begin{align*}
\frac{d}{d\lambda_1} & \rest{\lambda_1=0}\rho_0^1(\prod\exp(\lambda\cdot\hat{t}_p(\xi)))^{-1}\dot{\rho}_1(x)\omega(\exp(\lambda\cdot\xi))\\
					 & =\frac{d}{d\lambda_1}\rest{\lambda_1=0}\rho_0^1(h_2)^{-1}\rho_0^1(\exp(-\lambda_1\hat{t}_p(\xi_1)))\dot{\rho}_1(x)\omega(\exp(\lambda_1\xi_1),...,\exp(\lambda_q\xi_q))	\\
                     & =\rho_0^1(h_2)^{-1}\Big{[}-\dot{\rho}_0^1(\hat{t}_p(\xi_1))\dot{\rho}_1(x)\cancelto{0}{\omega(1,\exp(\lambda_{2}\xi_{2}),...,\exp(\lambda_q\xi_q))}+ \\
                     & \qquad\qquad\qquad\qquad\qquad\qquad +\cancelto{I}{\rho_0^1(\exp(0))}\dot{\rho}_1(x)R_{\xi_1}\omega(\exp(\lambda_{2}\xi_{2}),...,\exp(\lambda_q\xi_q))\Big{]}	\\
                     & =\rho_0^1(\prod_{j=2}^{q}\exp(\lambda_j\hat{t}_p(\xi_j)))^{-1}\dot{\rho}_1(x)R_{\xi_1}\omega(\exp(\lambda_{2}\xi_{2}),...,\exp(\lambda_q\xi_q));
\end{align*}
inductively implying,
\begin{align*}
R_{\xi_q}...R_{\xi_1}R_x\partial'\omega & =\dot{\rho}_1(x)R_{\xi_q}...R_{\xi_1}\omega .										
\end{align*}
Finally, considering the alternating sum
\begin{align*}
\Phi\partial'\omega(\xi_1,...,\xi_q;x) & =\sum_{\sigma\in S_q}\abs{\sigma}R_{\xi_{\sigma(q)}}...R_{\xi_{\sigma(1)}}R_x\partial'\omega \\
                                       & =\sum_{\sigma\in S_q}\abs{\sigma}\dot{\rho}_1(x)R_{\xi_{\sigma(q)}}...R_{\xi_{\sigma(1)}}\omega \\
                                       & =\dot{\rho}_1(x)\Phi\omega(\xi_1,...,\xi_q)=\delta_{(1)}\Phi\omega(\xi_1,...,\xi_q;x),
\end{align*}
and the result follows.

\end{proof}
Next up, we prove the compatibility of $\Phi$ with all the other differentials. That is the contents of the upcoming sequence of lemmas. \\
\begin{lemma}
For $q$ constant, $\Phi$ restricts to a map of double complexes between the $q$-pages.
\end{lemma}
\begin{proof}
Due to the lemma that states that $\Phi$ commutes with $\delta_{(1)}$, we just need to see that $\Phi\partial=\partial\Phi$. We already saw in the case of coefficients on an honest vector space that the relation holds for $\omega\in C(\G_p^q,V)$. In fact, in the course of the proof, we found the relation
\begin{eqnarray*}
\partial_k(\exp(\xi))=\exp(\partial_k\xi)
\end{eqnarray*}
to hold for all $\xi\in\gg_{p+1}$. Then, let $\omega\in C(\G_p^q\times G^r,W)$, $\xi_1,...,\xi_q\in\gg_{p+1}$ and $z_1,...,z_r\in\gg^r$ and assume the identification $\xi_b\sim(x_{0b},...,x_{pb},h_b)$. Computing,
\begin{align*}
R_{\xi_q}...R_{\xi_1}R_{z_r}...R_{z_1}(\partial\omega) & =\frac{d^I}{d\lambda_I}\rest{\lambda=0}\frac{d^J}{d\tau_J}\rest{\tau=0}\partial\omega(\exp(\lambda\cdot\xi);\exp(\tau\cdot z)). 
\end{align*}
In order to compute, we introduce the following notation: 
\begin{align*}
\gamma:=\begin{pmatrix}
\gamma_0 \\
\vdots   \\
\gamma_p
\end{pmatrix}=\exp(\lambda_1\xi_1)\vJoin ...\vJoin\exp(\lambda_q\xi_q)\in\G_{p+1} .   
\end{align*}
Then,
\begin{align*}
R_{\xi_q}...R_{\xi_1}R_{z_r}...R_{z_1}(\partial\omega) & =\frac{d^I}{d\lambda_I}\rest{\lambda=0}\frac{d^J}{d\tau_J}\rest{\tau=0}\Big{[}\rho_0^1(i(pr_G(\gamma_0))^{-1}\omega(\partial_0\exp(\lambda\cdot\xi);\exp(\tau\cdot z)) \\
									  & \qquad\qquad +\sum_{j=1}^{p+1}(-1)^j\omega(\partial_j\exp(\lambda\cdot\xi);\exp(\tau\cdot x))\Big{]} \\
									  & =\frac{d^I}{d\lambda_I}\rest{\lambda=0}\rho_0^1(i(pr_G(\gamma_0))^{-1}R_{z_r}...R_{z_1}\omega\big{(}\partial_0(\exp_{\G_{p+1}}(\lambda_1\xi_1)),...,\partial_0(\exp_{\G_{p+1}}(\lambda_q\xi_q))\big{)} \\
									  & \qquad\qquad +\sum_{j=1}^{p+1}(-1)^j\frac{d^I}{d\lambda_I}\rest{\lambda=0}R_{z_r}...R_{z_1}\omega\big{(}\partial_k(\exp_{\G_{p+1}}(\lambda_1\xi_1)),...,\partial_k(\exp_{\G_{p+1}}(\lambda_q\xi_q))\big{)} \\
									  & =\frac{d^I}{d\lambda_I}\rest{\lambda=0}\rho_0^1(i(pr_G(\gamma_0))^{-1}R_{z_r}...R_{z_1}\omega\big{(}\exp_{\G_{p+1}}(\lambda_1\partial_0\xi_1),...,\exp_{\G_{p+1}}(\lambda_q\partial_0\xi_q)\big{)} \\
									  & \qquad\qquad +\sum_{j=0}^{p+1}(-1)^j\frac{d^I}{d\lambda_I}\rest{\lambda=0}R_{z_r}...R_{z_1}\omega(\exp_{\G_{p}}(\lambda_1\partial_k\xi_1),...,\exp_{\G_{p}}(\lambda_q\partial_k\xi_q)). 
\end{align*}
Now,
\begin{align*}
\frac{d}{d\lambda_1}\rest{\lambda_1=0} & \rho_0^1(i(pr_G(\gamma_0))^{-1}R_{z_r}...R_{z_1}\omega\big{(}\exp(\lambda_1\partial_0\xi_1),...,\exp(\lambda_q\partial_0\xi_q)\big{)} \\
                                       & =\Big{(}\frac{d}{d\lambda_1}\rest{\lambda_1=0}\rho_0^1(i(pr_G(\gamma_0))^{-1}\Big{)}R_{z_r}...R_{z_1}\omega\big{(}\exp(0),\exp(\lambda_2\partial_0\xi_2),...,\exp(\lambda_q\partial_0\xi_q)\big{)}+ \\
                                       & \qquad\qquad +\rho_0^1(i(pr_G(\gamma_0'))^{-1}R_{\partial_0\xi_1}R_{z_r}...R_{z_1}\omega\big{(}\exp(\lambda_2\partial_0\xi_2),...,\exp(\lambda_q\partial_0\xi_q)\big{)},
\end{align*}
where
\begin{align*}
\gamma':=\begin{pmatrix}
\gamma_0' \\
\vdots   \\
\gamma_p'
\end{pmatrix}=\exp(0)\vJoin\exp(\lambda_2\xi_2)\vJoin ...\vJoin\exp(\lambda_q\xi_q)\in\G_{p+1};  
\end{align*}
thus, inductively,
\begin{align*}
R_{\xi_q}...R_{\xi_1}R_{z_r}...R_{z_1}(\partial\omega)	& =\sum_{j=0}^{p+1}(-1)^jR_{\partial_k\xi_q}...R_{\partial_k\xi_1}R_{z_r}...R_{z_1}\omega 
\end{align*}
and as a consequence
\begin{align*}
\Phi(\partial\omega)(\xi_1,...,\xi_q;z_1,...,z_r) & =\sum_{\sigma\in S_q}\sum_{\varrho\in S_r}\abs{\sigma}\abs{\varrho}R_{\xi_{\sigma(q)}}...R_{\xi_{\sigma(1)}}R_{z_{\varrho(r)}}...R_{z_{\varrho(1)}}(\partial\omega) \\
									  & =\sum_{\sigma\in S_q}\sum_{\varrho\in S_r}\abs{\sigma}\abs{\varrho}\sum_{j=0}^{p+1}(-1)^jR_{\partial_j\xi_{\sigma(q)}}...R_{\partial_j\xi_{\sigma(1)}}R_{z_{\varrho(r)}}...R_{z_{\varrho(1)}}\omega \\
									  & =\sum_{j=0}^{p+1}(-1)^j\sum_{\sigma\in S_q}\sum_{\varrho\in S_r}\abs{\sigma}\abs{\varrho}R_{\partial_j\xi_{\sigma(q)}}...R_{\partial_j\xi_{\sigma(1)}}R_{z_{\varrho(r)}}...R_{z_{\varrho(1)}}\omega \\
									  & =\sum_{j=0}^{p+1}(-1)^j\Phi\omega(\partial_j\xi_{1},...,\partial_j\xi_{q};z_1,...,z_r)=\partial(\Phi\omega)(\xi_1,...,\xi_q;z_1,...,z_r)
\end{align*}
as desired.

\end{proof}
\begin{proposition}
$\Phi$ commutes with all difference maps.
\end{proposition}
This latter result implies that there is an induced map of triple complexes. In particular, there are honest maps of double complexes between the $p$-pages and we can consider the mapping cone associated to them. Recall that the total complexes of the $p$-pages form a double complex for Lie $2$-algebras, whereas, for Lie $2$-groups there was a bit more structure (cf. Remark \ref{DblCxForGpsRmk}). Due to the extra differentials in the ``double complex'' of $p$-pages of the triple complex of Lie $2$-group cochains, putting the mapping cones together does not quite build a double complex, but rather a structure that shadows that of $C_{p-pag}$ in the referred remark. Nevertheless, when filtrating by columns, the induced spectral sequence has the same first page as the spectral sequence of an honest double complex. We are going to use this idea to go about computing when there are induced isomorphisms in cohomology. This is the contents of what follows.
%*** Insert diagram? Should I have written this elsewhere: the fact that the triple complex isn't quite so, but it induces a double complex where the columns are given by the total complexes of the $p$-pages? *** \\

%Close remarking that all of the ingredients are in place, and that the only piece missing is a van Est theorem for the $2$-van Est map.

%---------------------------------------------
%---------------------------------------------
\section{A collection of van Est type theorems}
\sectionmark{The main theorem}
%---------------------------------------------
%---------------------------------------------
In this section, we will prove the main theorem of this thesis, which is a van Est type result:
\begin{theorem}\label{2-vanEstTheo}
Let $\G=\xymatrix{G \ar[r] & H}$ represent a Lie $2$-group with Lie $2$-algebra represented by $\gg_1=\xymatrix{\gg \ar[r] & \hh}$, and $\rho$ be a $2$-representation of $\G$ on a $2$-vector space $\xymatrix{W \ar[r]^\phi & V}$. If $H$ and $G$ are both $k$-connected, then the $2$-van Est map induces isomorphisms
\begin{eqnarray*}
\xymatrix{
\Phi^n :H_\nabla^n(\G ,\phi) \ar[r] & H_\nabla^n(\gg_1 ,\phi),
}
\end{eqnarray*}
for $n\leq k$, and it is injective for $n=k+1$.
\end{theorem}
The strategy of the proof will go along the lines of the proof of theorem \ref{2vE-vs}. First, we are to consider the map of honest double complexes induced by $\Phi$. In this setting, we look at the mapping cone double complex and using the spectral sequence of the filtration by columns, the result will follow from understanding how $\Phi$ behaves when restricted to these columns, i.e. when restricted to the total complex of the $p$-pages. The missing piece of this proof is, thus, that the cohomology of the columns to vanishes below the diagonal. We collect the necessary results to conclude this statement, which, in order, will let us argue that the total cohomology vanishes, thus implying the result. As it turns out, most of these results are themselves van Est type theorems, hence the name of the section. \\

Suppose we are back in the situation we had at the beginning of the section where we laid down the shape of the $p$-pages of the triple complex of Lie $2$-group cochains with values in a $2$-representation, i.e. that we have got two Lie groups $H$ and $G$, with $H$ acting on the right on $G$, a vector space $W$ and a map of double Lie groups 
\begin{eqnarray*}
\xymatrix{
 & H\ltimes G \ar@<0.5ex>[dl]\ar@<-0.5ex>[dl]\ar@<0.5ex>[dd]\ar@<-0.5ex>[dd]\ar[r] & GL(W)\ltimes GL(W) \ar@<0.5ex>[dr]\ar@<-0.5ex>[dr]\ar@<0.5ex>[dd]\ar@<-0.5ex>[dd] & \\
H \ar[rrr]^{\rho_H\qquad\quad}\ar@<0.5ex>[dd]\ar@<-0.5ex>[dd] & & & GL(W) \ar@<0.5ex>[dd]\ar@<-0.5ex>[dd]   \\
 &             G \ar@<0.5ex>[dl]\ar@<-0.5ex>[dl]\ar[r]^{\rho_G\quad} & GL(W) \ar@<0.5ex>[dr]\ar@<-0.5ex>[dr] & \\       
 \ast \ar[rrr]                                                     & & & \ast                   
}
\end{eqnarray*} 
As it turns out, associated to each of these double Lie groups, there are two \LA -groups,
\begin{eqnarray*}
\xymatrix{
\hh\ltimes G \ar@<0.5ex>[r]\ar@<-0.5ex>[r]\ar[d] & \hh \ar[d] \\       
G \ar@<0.5ex>[r]\ar@<-0.5ex>[r]                  & \ast   } & \xymatrix{
H\ltimes \gg \ar@<0.5ex>[r]\ar@<-0.5ex>[r]\ar[d] & \gg \ar[d] \\       
H \ar@<0.5ex>[r]\ar@<-0.5ex>[r]                  & \ast
}
\end{eqnarray*}
of $H\ltimes G$ and 
\begin{eqnarray*}
\xymatrix{
\ggl(W)\ltimes GL(W) \ar@<0.5ex>[r]\ar@<-0.5ex>[r]\ar[d] & \ggl(W) \ar[d] \\       
GL(W) \ar@<0.5ex>[r]\ar@<-0.5ex>[r]                  & \ast   } & \xymatrix{
GL(W)\ltimes \ggl(W) \ar@<0.5ex>[r]\ar@<-0.5ex>[r]\ar[d] & \ggl(W) \ar[d] \\       
GL(W) \ar@<0.5ex>[r]\ar@<-0.5ex>[r]                  & \ast
}
\end{eqnarray*}
of $GL(W)\ltimes GL(W)$. Bear in mind that the notation $\ltimes$ stands alternatively for the transformation groupoid associated to an action of a Lie group and the action Lie algebroid associated to the action of a Lie algebra. Correspondingly, there are two morphisms of \LA -groups which correspond to the differentiation of the map above in each direction. We will prove that on the one hand, differentiating in the vertical direction, we get an \LA -double complex values on $W$; whereas, adding the hypothesis that the action of $H$ is by Lie group automorphisms, differentiating in the horizontal direction yields a second \LA -double complex values on $W$. 
\begin{lemma}\label{Ver p-page}
Let $\xymatrix{\hh\ltimes G \ar[r] & \ggl(W)\ltimes GL(W):(y;g) \ar@{|->}[r] & (\rho_\hh(y),\rho_G(g))}$ be a map of \LA -groups. Then 
\begin{eqnarray*}
\xymatrix{
\vdots                                                    & \vdots                                                         & \vdots                        & \\ 
\bigwedge^3\hh^*\otimes W \ar[r]^{\delta'}\ar[u]          & C(G,\bigwedge^3\hh^*\otimes W) \ar[r]^{\delta'}\ar[u]          & C(G^2,\bigwedge^3\hh^*\otimes W) \ar[r]\ar[u]          & \dots \\
\bigwedge^2\hh^*\otimes W \ar[r]^{\delta'}\ar[u]^{\delta} & C(G,\bigwedge^2\hh^*\otimes W) \ar[r]^{\delta'}\ar[u]^{\delta} & C(G^2,\bigwedge^2\hh^*\otimes W) \ar[r]\ar[u]^{\delta} & \dots \\
\hh^*\otimes W  \ar[r]^{\delta'}\ar[u]^{\delta}           & C(G,\hh^*\otimes W) \ar[r]^{\delta'}\ar[u]^{\delta}            & C(G^2,\hh^*\otimes W) \ar[r]\ar[u]^{\delta}                 & \dots \\
W \ar[r]^{\delta'}\ar[u]^{\delta}                         & C(G,W) \ar[r]^{\delta'}\ar[u]^{\delta}                         & C(G^2,W)	\ar[r]\ar[u]^{\delta}                            & \dots  
}
\end{eqnarray*}
where the rows are subcomplexes of Lie groupoid cochains for the Lie group bundles
\begin{eqnarray*}
\xymatrix{
\hh^q\times G \ar@<0.5ex>[r] \ar@<-0.5ex>[r] & \hh^q 
} 
\end{eqnarray*}
on $\xymatrix{\hh^q\times W \ar[r] & \hh^q}$ and the columns are complexes of Lie algebroid cochains for the action Lie algebroids
\begin{eqnarray*}
\xymatrix{
 \hh\ltimes G^r \ar[r] & G^r 
} 
\end{eqnarray*}
on $\xymatrix{G^r\times W \ar[r] & G^r}$ is a double complex.
\end{lemma}
\begin{proof}
We start by pointing out that since the given map is a map of \LA -groups, its restriction to the zero section and to the unit Lie algebra give respectively representations $\rho_G$ and $\rho_\hh$ of $G$ and $\hh$ on $W$. Thus, the first row and first column are indeed complexes. In order to get the representations of the other rows and columns, we proceed to pull-back the Lie group representation along projection $\pi_q$ and the Lie algebra representation along the Lie algebroid map $\hat{t}_r$. Specifically, we define the representation $\rho_G^q$ of $\xymatrix{\hh^q\times G \ar@<0.5ex>[r] \ar@<-0.5ex>[r] & \hh^q}$ on $\xymatrix{\hh^q\times W \ar[r] & \hh^q}$ by 
\begin{eqnarray*}
\rho_G^q :=\pi_q^*\rho_G ,
\end{eqnarray*}
where indeed $p_q^*W=\hh^q\times W$. Analogously, we define the representation $\rho_\hh^r$ of $\xymatrix{\hh\ltimes G^r \ar[r] & G^r}$ on $\xymatrix{G^r\times W \ar[r] & G^r}$ by 
\begin{eqnarray*}
\rho_\hh^r :=\hat{t}_r^*\rho_\hh ,
\end{eqnarray*}
and this time round, $t_r^*W=G^r\times W$. \\
Having clarified what the representations are, we show that the spaces of cochains are the right ones, and exhibit the sub-complex that we will be using. On the one hand,
\begin{eqnarray*}
(\hh^q\times G)^{(r)}=\lbrace (Y_1,g_1;...;Y_r,g_r)\in(\hh^q\times G)^{r}:\hat{s}(Y_k,g_k)=\hat{t}(Y_{k+1},g_{k+1})\rbrace ;
\end{eqnarray*}
therefore, $(\hh^q\times G)^{(r)}\cong \hh^q\times G^r$ and the diffeomorphism is given again by the obvious assignment $\xymatrix{(Y,g_1;...;Y,g_r) \ar@{|->}[r] & (Y;g_1,...,g_r)}$. Since the representation is taken on a vector bundle which is trivial, its pull-back will be trivial as well, and its sections will coincide with smooth functions to $W$, i.e. $C^r(\hh^q\times G,\hh^q\times W)=C(\hh^q\times G^r,W)$. On the other hand, the columns are complexes of Lie algebroid cochains with values in a representation which is trivial, as the pull-back of a trivial vector bundle; therefore,
\begin{eqnarray*}
C^q(\hh\ltimes G^r,G^r\times W)=\Gamma(\bigwedge^q(\hh\ltimes G^r)^*\otimes(G^r\times W))=\Gamma(G^r\times(\bigwedge^q\hh^*\otimes W))
\end{eqnarray*}
In so, the space of $q$-cochains is $C^q(\hh\ltimes G^r,G^r\times W)=C(G^r,\bigwedge^q\hh^*\otimes W)$. In order for these to coincide with the spaces of $r$-cochains of the previous paragraph, we consider the sub-spaces of $r$-cochains which are alternating $q$-multilinear in the $\hh$-coordinates. In symbols,
\begin{eqnarray*}
C^r_{lin}(\hh^q\times G,W):=\lbrace\omega\in C(\hh^q\times G^r,W):\omega(-;\vec{g})\in\bigwedge^q\hh^*\otimes W,\forall\vec{g}\in G^r\rbrace .
\end{eqnarray*}
We will use the diffeomorphism of the latter discussion to write formulas for the face maps of the simplicial structure and then prove that $C^\bullet_{lin}(\hh^q\times G,W)$ is indeed a subcomplex: For $Y=(y_0,...,y_q)\in \hh^{q}$ and $\vec{g}=(g_0,...,g_r)\in G^{r+1}$, since
\begin{eqnarray*}
\delta'_k (\vec{g})=
  \begin{cases}
    (g_1,...,g_r)                                 & \quad \text{if } k=0  \\
    (g_0,...,g_{k-2},g_{k-1}g_k,g_{k+1},...,g_r)  & \quad \text{if } 0<k\leq r\\
    (g_0,...,g_{r-1})                             & \quad \text{if } k=r+1 ,
  \end{cases}
\end{eqnarray*}
\begin{eqnarray*}
\delta'_k(Y;\vec{g}) & =(Y;\delta'_k\vec{g}).
\end{eqnarray*}
Thus, for $\omega\in C^r_{lin}(\hh^q\times G,W)$,
\begin{eqnarray*}
\delta'\omega(Y;\vec{g})=\rho_G^q(Y;g_0)\omega(Y;\delta'_0\vec{g})+\sum_{k=1}^{r+1}(-1)^k\omega(Y;\delta'_k\vec{g}),
\end{eqnarray*}
but $\rho_G^q(Y;g_0)=\rho_G(\pi_q(Y;g_0))=\rho_G(g_0)$; hence, $\delta'\omega(-;\vec{g})$ is a linear combination of alternating $q$-multilinear maps and $\delta'\omega\in C^{r+1}_{lin}(\hh^q\times G,W)$. Now, by the very definition, 
\begin{eqnarray*}
C^r_{lin}(\hh^q\times G,W)=C(G^r,\bigwedge^q\hh^*\otimes W),
\end{eqnarray*}
so we have the right spaces of cochains. \\
The only thing left to prove is the statement itself, that the generic square
\begin{eqnarray}\label{generic Ver p-square}
\xymatrix{ 
C(G^r,\bigwedge^{q+1}\hh^*\otimes W) \ar[r]^{\delta'}            & C(G^{r+1},\bigwedge^{q+1}\hh^*\otimes W)          \\
C(G^r,\bigwedge^q\hh^*\otimes W) \ar[r]_{\delta'}\ar[u]^{\delta} & C(G^{r+1},\bigwedge^q\hh^*\otimes W) \ar[u]_{\delta} 
}
\end{eqnarray}
commutes. Notice that since the bracket of the Lie algebroids involved is completely determined by the bracket of $\hh$, it suffices to prove the commutativity of the latter diagram for constant sections. Indeed, let $\omega\in C(G^r,\bigwedge^q\hh^*\otimes W)$ and $Y=(y_0,...,y_q)\in\hh^q$ and $\vec{g}$ as above. Then, 
\begin{align*}
\delta'\delta\omega(Y;\vec{g}) & =\rho_G^{q+1}(Y;g_0)\delta\omega(Y;\delta'_0\vec{g})+\sum_{k=1}^{r+1}(-1)^{k}\delta\omega(Y;\delta'_k\vec{g}),
\end{align*}
while
\begin{align*}
\delta\delta'\omega(Y;\vec{g}) & =\sum_{m<n}(-1)^{m+n}\delta'\omega([y_m,y_n],Y(m,n);\vec{g})+\sum_{j=0}^{q}(-1)^{j+1}\rho_\hh^{r+1}(y_j)\delta'\omega(Y(j);\vec{g}).
\end{align*}
We expand further to make evident the common terms:
\begin{align*}
\delta'\delta & \omega(Y;\vec{g}) \\
              & =\rho_G^{q+1}(Y;g_0)\Big{[}\sum_{m<n}(-1)^{m+n}\omega([y_m,y_n],Y(m,n);\delta'_0\vec{g})+\sum_{j=0}^{q}(-1)^{j+1}\rho_\hh^{r}(y_j)\omega(Y(j);\delta'_0\vec{g})\Big{]}+ \\
			  & \qquad +\sum_{k=1}^{r+1}(-1)^{k}\Big{[}\sum_{m<n}(-1)^{m+n}\omega([y_m,y_n],Y(m,n);\delta'_k\vec{g})+\sum_{j=0}^{q}(-1)^{j+1}\rho_\hh^{r}(y_j)\omega(Y(j);\delta'_k\vec{g})\Big{]},										   
\end{align*}
and
\begin{align*}
\delta\delta' & \omega(Y;\vec{g}) \\
              & =\sum_{m<n}(-1)^{m+n}\Big{[}\rho_G^{q}([y_m,y_n],Y(m,n);g_0)\omega([y_m,y_n],Y(m,n);\delta'_0\vec{g})+\sum_{k=1}^{r+1}(-1)^{k}\omega([y_m,y_n],Y(m,n);\delta'_k\vec{g})\Big{]}+ \\
              & \qquad +\sum_{j=0}^{q}(-1)^{j+1}\rho_\hh^{r+1}(y_j)\Big{[}\rho_G^{q}(Y(j);g_0)\omega(Y(j);\delta'_0\vec{g})+ \sum_{k=1}^{r+1}(-1)^{k}\omega(Y(j);\delta'_k\vec{g})\Big{]}
\end{align*}
The desired equality follows now by noticing the following identities. First, we obviously have that
\begin{eqnarray*}
\rho_G^{q+1}(Y;g_0)=\rho_G^q([y_m,y_n],Y(m,n);g_0),
\end{eqnarray*}
as they are both equal to $\rho_G(g_0)$. Next, let $\lbrace e_a\rbrace$ be a basis for $W$, and 
\begin{align*}
\omega(Y(j);\delta'_0\vec{g})=\omega^a e_a ,
\end{align*}
where $\omega^a=\omega^a(Y(j);\delta'_0\vec{g})$ and we are using the Einstein summation convention. Since the anchor of $\hh\ltimes G^r$ is given by
\begin{align*}
d_{(1;\vec{g})}\Lf{t}(y,0) & =\frac{d}{d\lambda}\rest{\lambda=0}(\vec{g})^{\exp(\lambda y)}\in T_{\vec{g}}G^r  
\end{align*} 
and for $(y;\vec{g})\in\hh\ltimes G^r$,
\begin{align*}
\hat{t}_r(y;\vec{g})=y ,
\end{align*} 
then, by definition
\begin{eqnarray*}
\rho_\hh^r(y_j)\omega(Y(j);\delta'_0\vec{g})=\omega^a\rho_\hh(y_j)e_a+(\frac{d}{d\lambda}\rest{\lambda=0}\omega^a(Y(j);(\delta'_0\vec{g})^{\exp(\lambda y_j)}))e_a ,
\end{eqnarray*}
whereas
\begin{eqnarray*}
\rho_\hh^{r+1}(y_j)\omega(Y(j),\delta'_0\vec{g})=\omega^a\rho_\hh(y_j)e_a+(\frac{d}{d\lambda}\rest{\lambda=0}\omega^a(Y(j),\delta'_0(\vec{g})^{\exp(\lambda y_j)}))e_a .
\end{eqnarray*}
Hence, they clearly coincide as for every $h\in H$ and every ranging value of $k$,
\begin{eqnarray*}
(\delta'_k\vec{g})^{h}=\delta'_k(\vec{g})^{h}.
\end{eqnarray*}
Finally, as $\rho_\hh\times\rho_G$ is an \LA -group map, 
\begin{eqnarray*}
\xymatrix{
TG \ar[r]^{d\rho_G\qquad}                              & TGL(W) \\
\hh\ltimes G \ar[r]_{\rho_\hh\times\rho_G\qquad}\ar[u] & \ggl(W)\ltimes GL(W)\ar[u] ,
}
\end{eqnarray*}
where the vertical maps are the anchors, commutes. Now, right side anchor is 
\begin{align*}
d_{(I,A)}\Lf{t}(X,0) & =\frac{d}{d\lambda}\rest{\lambda=0}\exp(\lambda A)X\exp(\lambda A)^{-1}=AX-XA\in T_AGL(W);
\end{align*}
thereby yielding,
\begin{align*}
\frac{d}{d\lambda}\rest{\lambda=0}\rho_G(g^{\exp(\lambda y)}) & =\rho_G(g)\rho_\hh(y)-\rho_\hh(y)\rho_G(g).
\end{align*}
Computing,
\begin{eqnarray*}
\rho_G^{q+1}(Y;g_0)\rho_\hh^r(y_j)\omega(Y(j);\delta'_0\vec{g})=\rho_G(g_0)\Big{[}\omega^a\rho_\hh(y_j)e_a+\big{(}\frac{d}{d\lambda}\rest{\lambda=0}\omega^k(Y(j);(\delta'_0\vec{g})^{\exp(\lambda y_j)}\big{)}e_a\Big{]};
\end{eqnarray*}
while on the other hand, $\rho_G^{q}(Y;g_0)\omega(Y(j);\delta'_0\vec{g})=\omega^a\rho_G(g_0)e_a$ and
\begin{align*}
\rho_\hh^{r+1} & (y_j)\rho_G^{q}(Y;g_0)\omega(Y(j);\delta'_0\vec{g}) \\
  & =\omega^a\rho_\hh(y_j)\rho_G(g_0)e_a+\big{(}\frac{d}{d\lambda}\rest{\lambda=0}\omega^k(Y(j);(\delta'_0\vec{g})^{\exp(\lambda y_j)}\rho_G(g_0^{\exp(\lambda y_j)})\big{)}e_a \\
  & =\omega^a\rho_\hh(y_j)\rho_G(g_0)e_a+\big{(}\frac{d}{d\lambda}\rest{\lambda=0}\omega^k(Y(j);(\delta'_0\vec{g})^{\exp(\lambda y_j)}\big{)}\rho_G(g_0^{\exp(0)})e_a+ \\
  & \qquad\qquad +\omega^a(Y(j);(\delta'_0\vec{g})^{\exp(0)})\big{(}\frac{d}{d\lambda}\rest{\lambda=0}\rho_G(g_0^{\exp(\lambda y_j)})\big{)}e_a \\
  & =\omega^a\rho_\hh(y_j)\rho_G(g_0)e_a+\rho_G(g_0)\Big{[}\big{(}\frac{d}{d\lambda}\rest{\lambda=0}\omega^k(Y(j);(\delta'_0\vec{g})^{\exp(\lambda y_j)}\big{)}e_a+\omega^a\rho_\hh(y_j)e_a\Big{]}-\omega^a\rho_\hh(y_j)\rho_G(g_0)e_a,
\end{align*}
ultimately implying the commutativity of the square \ref{generic Ver p-square}. 

\end{proof} 
\begin{lemma}\label{Hor p-page}
Let $H$ act on the right on $\gg$ by Lie algebra automorphisms. Further suppose that $\xymatrix{H\ltimes \gg \ar[r] & GL(W)\ltimes\ggl(W):(h;x) \ar@{|->}[r] & (\rho_H(h),\rho_\gg(x))}$ is a map of \LA -groups. Then 
\begin{eqnarray*}
\xymatrix{
\vdots                                   & \vdots                                                & \vdots                                & \\ 
C(H^3,W) \ar[r]^{\delta'}\ar[u]          & C(H^3,\gg^*\otimes W) \ar[r]^{\delta'}\ar[u]          & C(H^3,\bigwedge^2\gg^*\otimes W) \ar[r]\ar[u]          & \dots \\
C(H^2,W) \ar[r]^{\delta'}\ar[u]^{\delta} & C(H^2,\gg^*\otimes W) \ar[r]^{\delta'}\ar[u]^{\delta} & C(H^2,\bigwedge^2\gg^*\otimes W) \ar[r]\ar[u]^{\delta} & \dots \\
C(H,W) \ar[r]^{\delta'}\ar[u]^{\delta}   & C(H,\gg^*\otimes W) \ar[r]^{\delta'}\ar[u]^{\delta}   & C(H,\bigwedge^2\gg^*\otimes W) \ar[r]\ar[u]^{\delta}       & \dots \\
W \ar[r]^{\delta'}\ar[u]^{\delta}        & \gg^*\otimes W \ar[r]^{\delta'}\ar[u]^{\delta}        & \bigwedge^2\gg^*\otimes W	\ar[r]\ar[u]^{\delta}       & \dots  
}
\end{eqnarray*}
where the rows are the complexes of Lie algebroid cochains for the Lie algebra bundles
\begin{eqnarray*}
\xymatrix{
H^q\times \gg \ar@<0.5ex>[r] \ar@<-0.5ex>[r] & H^q 
} 
\end{eqnarray*}
on $\xymatrix{H^q\times W \ar[r] & H^q}$ and and the columns are sub-complexes of Lie groupoid cochains for the transformation groupoids
\begin{eqnarray*}
\xymatrix{
 H^{op}\ltimes \gg^r \ar@<0.5ex>[r] \ar@<-0.5ex>[r] & \gg^r 
} 
\end{eqnarray*}
on $\xymatrix{\gg^r\times W \ar[r] & \gg^r}$ is a double complex.
\end{lemma}
\begin{proof}
We start again by pointing out that since the given map is a map of \LA -groups, its restriction to the zero section and to the unit Lie algebra give respectively representations $\rho_H$ and $\rho_\gg$ of $H$ and $\gg$ on $W$. Thus, the first row and first column are indeed complexes. In order to get the representations of the other rows and columns, we proceed to pull-back the Lie group representation along projection $\pi_r$ and the Lie algebra representation along the Lie algebroid map $\hat{t}_q$. Specifically, we define the representation $\rho_H^r$ of $\xymatrix{H\ltimes\gg^r \ar@<0.5ex>[r]\ar@<-0.5ex>[r] & \gg^r}$ on $\xymatrix{\gg^r\times W \ar[r] & \gg^r}$ by 
\begin{eqnarray*}
\rho_H^r :=\pi_r^*\rho_H ,
\end{eqnarray*}
where indeed $p_r^*W=\gg^r\times W$. Analogously, we define the representation $\rho_\gg^q$ of $\xymatrix{H^q\times\gg \ar[r] & H^q}$ on $\xymatrix{H^q\times W \ar[r] & H^q}$ by 
\begin{eqnarray*}
\rho_\gg^q :=\hat{t}_q^*\rho_\gg ,
\end{eqnarray*}
and this time round, $t_q^*W=H^q\times W$. \\
Having clarified what the representations are, we use the same strategy of the proof of the previous lemma and define the columns to be the sub-complexes of groupoid cochains that are alternating $r$-multilinear in the $\gg$-coordinates. In order to prove that these indeed yields a sub-complexes, consider
\begin{eqnarray*}
(H\ltimes\gg^r)^{(q)}=\lbrace (h_1,X_1;...;h_q,X_q)\in(H\times\gg^r)^{q}:\hat{s}(h_j,X_j)=\hat{t}(h_{j+1},X_{j+1})\rbrace ;
\end{eqnarray*}
therefore, $(H\times\gg^r)^{(q)}\cong H^q\times\gg^r$, where the diffeomorphism is 
\begin{eqnarray*}
\xymatrix{(h_1,X^{h_q...h_2};h_2,X^{h_q...h_3};...;h_{q-1},X^{h_q};h_q,X) \ar@{|->}[r] & (h_1,...,h_q;X).}
\end{eqnarray*} 
The considerations regarding the triviality of the vector bundle still hold; thus, $C^q(H\ltimes\gg^r,\gg^r\times W)=C(H^q\times\gg^r,W)$. As stated before, define
\begin{eqnarray*}
C^q_{lin}(H\ltimes\gg^r ,W):=\lbrace\omega\in C(H^q\times\gg^r,W):\omega(\vec{h};-)\in\bigwedge^r\gg^*\otimes W,\forall\vec{h}\in H^q\rbrace .
\end{eqnarray*}
We use the diffeomorphism to write formulas for the face maps of the simplicial structure: For $\vec{h}=(h_0,...,h_q)\in H^{q+1}$ and $X=(x_1,...,x_r)\in\gg^{r}$, since
\begin{eqnarray*}
\delta_j (\vec{h})=
  \begin{cases}
    (h_1,...,h_q)                                 & \quad \text{if } j=0  \\
    (h_0,...,h_{j-2},h_jh_{j-1},h_{j+1},...,h_q)  & \quad \text{if } 0<j\leq q\\
    (g_0,...,g_{r-1})                             & \quad \text{if } j=q+1 ,
  \end{cases}
\end{eqnarray*}
the formula for the face maps of $H\ltimes\gg^r$ are
\begin{eqnarray*}
\delta_j(\vec{h};X)=
  \begin{cases}
    (\delta_j\vec{h};X)             & \quad \text{if } 0\leq j\leq q\\
    (\delta_{q+1}\vec{h};X^{h_q})   & \quad \text{otherwise} 
  \end{cases}
\end{eqnarray*}
Thus, for $\omega\in C^q_{lin}(H\ltimes\gg^r ,W)$,
\begin{eqnarray*}
\delta\omega(\vec{h};X)=\rho_H^r(h_0;X^{h_q...h_1})\omega(\delta_0\vec{h};X)+\sum_{j=1}^{q}(-1)^j\omega(\delta_j\vec{h};X)+(-1)^{q+1}\omega(\delta_{q+1}\vec{h};X^{h_q}),
\end{eqnarray*}
but $\rho_H^r(h_0;X^{h_q...h_1})=\rho_H(\pi_r(h_0;X^{h_q...h_1}))=\rho_H(h_0)$, $X^{h_q}:=(x_1^{h_q},...,x_r^{h_q})$ and $(-)^{h_q}$ is linear; hence, $\delta\omega(\vec{h};-)$ is a linear combination of alternating $r$-multilinear maps and $\delta\omega\in C^{q+1}_{lin}(H\ltimes\gg^r,W)$. Again, by the very definition, 
\begin{eqnarray*}
C^q_{lin}(H\ltimes\gg^r,W)=C(H^q,\bigwedge^r\gg^*\otimes W),
\end{eqnarray*}
so these coincide with the spaces of $r$-cochains for the Lie algebroid complexes of the rows. Indeed, because the pull-back vector bundle is trivial,
\begin{eqnarray*}
C^r(H^q\times\gg ,H^q\times W)=\Gamma(\bigwedge^r(H^q\times\gg)^*\otimes(H^q\times W))=\Gamma(H^q\times(\bigwedge^r\gg^*\otimes W))
\end{eqnarray*}
To prove the statement itself, we see that the generic square
\begin{eqnarray}\label{generic Hor p-square}
\xymatrix{ 
C(H^{q+1},\bigwedge^r\gg^*\otimes W) \ar[r]^{\delta'}            & C(H^{q+1},\bigwedge^{r+1}\gg^*\otimes W)          \\
C(H^q,\bigwedge^r\gg^*\otimes W) \ar[r]^{\delta'}\ar[u]^{\delta} & C(H^q,\bigwedge^{r+1}\gg^*\otimes W) \ar[u]^{\delta} 
}
\end{eqnarray}
commutes. Now, the Lie algebroids involved are trivial bundles of Lie algebras, so the brackets are completely determined by the bracket of $\gg$; thus, it suffices to prove the commutativity of the latter diagram for constant sections. Indeed, let $\omega\in C(H^q,\bigwedge^r\gg^*\otimes W)$ and $X=(x_0,...,x_r)\in\gg^r$ and $\vec{h}$ as above. Then, 
\begin{align*}
\delta'\delta\omega(\vec{h};X) & =\sum_{m<n}(-1)^{m+n}\delta\omega(\vec{h};[x_m,x_n],X(m,n))+\sum_{k=0}^{r}(-1)^{k+1}\rho_\gg^{q+1}(x_k)\delta\omega(\vec{h};X(k)),
\end{align*}
while
\begin{align*}
\delta\delta'\omega(Y;\vec{g}) & =\rho_H^{r+1}(h_0;X^{h_q...h_1})\delta'\omega(\delta_0\vec{h};X)+\sum_{j=1}^{q}(-1)^{j}\delta'\omega(\delta_j\vec{h};X)+(-1)^{q+1}\delta'\omega(\delta_{q+1}\vec{h};X^{h_q}).
\end{align*}
We expand further to make evident the common terms:
\begin{align*}
\delta'\delta\omega(\vec{h};X) & =\sum_{m<n}(-1)^{m+n}\Big{[}\rho_H^{r}(h_0;[x_m,x_n]^{h_q...h_1},X(m,n)^{h_q...h_1})\omega(\delta_0\vec{h};[x_m,x_n],X(m,n))+ \\
							   & \qquad\qquad +\sum_{j=1}^{q}(-1)^{j}\omega(\delta_j\vec{h};[x_m,x_n],X(m,n))+(-1)^{q+1}\omega(\delta_{q+1}\vec{h};[x_m,x_n]^{h_q},X(m,n)^{h_q})\Big{]}+ \\
						       & \quad +\sum_{k=0}^{r}(-1)^{k+1}\rho_\gg^{q+1}(x_k)\Big{[}\rho_H^{r}(h_0;X(k)^{h_q...h_1})\omega(\delta_0\vec{h};X(k))+ \\
						       & \qquad\qquad\qquad\qquad +\sum_{j=1}^{q}(-1)^{j}\omega(\delta_j\vec{h};X(k))+(-1)^{q+1}\omega(\delta_{q+1}\vec{h};X(k)^{h_q})\Big{]}								   
\end{align*}
and
\begin{align*}
\delta\delta' & \omega(\vec{h};X) \\
              & =\rho_H^{r+1}(h_0;X^{h_q...h_1})\Big{[}\sum_{m<n}(-1)^{m+n}\omega(\delta_0\vec{h};[x_m,x_n],X(m,n))+\sum_{k=0}^{r}(-1)^{k+1}\rho_\gg^{q}(x_k)\omega(\delta_0\vec{h};X(k))\Big{]}+ \\
              & \qquad +\sum_{j=1}^{q}(-1)^{j}\Big{[}\sum_{m<n}(-1)^{m+n}\omega(\delta_j\vec{h};[x_m,x_n],X(m,n))+\sum_{k=0}^{r}(-1)^{k+1}\rho_\gg^{q}(x_k)\omega(\delta_j\vec{h};X(k))\Big{]}+ \\
              & \qquad\quad +(-1)^{q+1}\Big{[}\sum_{m<n}(-1)^{m+n}\omega(\delta_{q+1}\vec{h};[x_m^{h_q},x_n^{h_q}],X^{h_q}(m,n))+\sum_{k=0}^{r}(-1)^{k+1}\rho_\gg^{q}(x_k^{h_q})\omega(\delta_{q+1}\vec{h};X^{h_q}(k))\Big{]} ,	
\end{align*}
The desired equality eventually follows after noticing the following identities. First, we obviously have that
\begin{eqnarray*}
\rho_H^{r}(h_0;[x_m,x_n]^{h_q...h_1},X(m,n)^{h_q...h_1})=\rho_H^{r+1}(h_0;X^{h_q...h_1}),
\end{eqnarray*}
as they are both equal to $\rho_H(h_0)$. Next, thanks to the fact that the action of $H$ is by Lie algebra automorphisms
\begin{eqnarray*}
[x_m,x_n]^{h_q}=[x_m^{h_q},x_n^{h_q}],
\end{eqnarray*}
and also, we obviously have
\begin{eqnarray*}
X(m,n)^{h_q}=X^{h_q}(m,n) & X(k)^{h_q}=X^{h_q}(k).
\end{eqnarray*}
Since $H^q\times\gg$ is a bundle of Lie algebras, its anchor is zero and for $(h_1,...,h_q;x)\in H^q\times\gg$,
\begin{align*}
\hat{t}_q(h_1,....,h_q;x)=x^{h_q...h_1} ,
\end{align*} 
then, by definition
\begin{eqnarray*}
\rho_\gg^{q+1}(x_k)\delta\omega(\vec{h};X(k))=\rho_\gg(x_k^{h_q...h_0})\delta\omega(\vec{h};X(k)).
\end{eqnarray*}
On the other hand, 
\begin{eqnarray*}
\rho_\gg^q(x_k)\omega(\delta_j\vec{h};X(k))=\rho_\gg(x_k^{h_q...(h_jh_{j-1})...h_0})\omega(\delta_j\vec{h};X(k)) ;
\end{eqnarray*}
thus, for all ranging values $0<j\leq q$
\begin{eqnarray*}
\rho_\gg^{q+1}(x_k)\omega(\delta_j\vec{h};X(k))=\rho_\gg^q(x_k)\omega(\delta_j\vec{h};X(k)),
\end{eqnarray*}
whereas
\begin{align*}
\rho_\gg^{q+1}(x_k)\omega(\delta_{q+1}\vec{h};X(k)^{h_q}) & =\rho_\gg((x_k^{h_q})^{h_{q-1}...h_0})\omega(\delta_{q+1}\vec{h};X^{h_q}(k)) \\
                                                          & =\rho_\gg^q(x_k^{h_q})\omega(\delta_{q+1}\vec{h};X^{h_q}(k)) .
\end{align*}
Finally, as $\rho_H\times\rho_\gg$ is an \LA -group map, 
\begin{align*}
\rho_\gg(\Lf{t}(h;x)) & =\Lf{t}(\rho_H(h),\rho_\gg(x)) \\
\rho_\gg(x^h)         & =Ad_{\rho_H(h)}\rho_\gg(x)=\rho_H(h)\rho_\gg(x)\rho_H(h)^{-1};
\end{align*}
thereby yielding,
\begin{align*}
\rho_\gg(x_k^{h_q...h_0}) & =\rho_\gg((x_k^{h_q...h_1})^{h_0}) \\
						  & =\rho_H(h_0)\rho_\gg(x_k^{h_q...h_1})\rho_H(h_0)^{-1},
\end{align*}
which in order, together with $\rho_\gg^q(x_k)\omega(\delta_0\vec{h};X(k))=\rho_\gg(x_k^{h_q...h_1})\omega(\delta_0\vec{h};X(k))$, implies
\begin{eqnarray*}
\rho_\gg^{q+1}(x_k)\rho_H^r(h_0;X(k)^{h_q...h_1})\omega(\delta_0\vec{h};X(k))=\rho_H^{r+1}(h_0;X^{h_q...h_1})\rho_\gg^q(x_k)\omega(\delta_0\vec{h};X(k))
\end{eqnarray*}
ultimately proving the commutativity of the square \ref{generic Hor p-square}. 

\end{proof} 
We argue that starting out with the the double complex of lemma \ref{p-pageModel}, there are maps going to the double complexes of the latter lemmas. Indeed, the $r$-th column of the double complex of lemma \ref{p-pageModel} is the complex of Lie groupoid cochains of the transformation groupoid
\begin{eqnarray*}
\xymatrix{
H\ltimes G^r \ar@<0.5ex>[r] \ar@<-0.5ex>[r] & G^r , 
} 
\end{eqnarray*} 
whose Lie algebroid is precisely the action Lie algebroid
\begin{eqnarray*}
\xymatrix{
\hh\ltimes G^r \ar[r] & G^r . 
} 
\end{eqnarray*} 
As stated, the $r$-th column of the double complex in lemma \ref{Ver p-page} is the Chevalley-Eilenberg complex of this Lie algebroid; thus, there is a van Est map taking us from the former to the latter. Assembling these column-wise van Est maps we get a map
\begin{eqnarray*}
\xymatrix{
\Phi_{col}:C(H^q\times G^r,W) \ar[r] & C(G^r,\bigwedge^q\hh^*\otimes W) . 
} 
\end{eqnarray*}
Accordingly, the $q$-th row of the double complex of lemma \ref{p-pageModel} is the complex of Lie groupoid cochains of the Lie group bundle 
\begin{eqnarray*}
\xymatrix{
H^q\times G \ar@<0.5ex>[r] \ar@<-0.5ex>[r] & H^q  
} 
\end{eqnarray*} 
with Lie algebroid the Lie algebra bundle
\begin{eqnarray*}
\xymatrix{
H^q\times\gg \ar[r] & H^q  
} 
\end{eqnarray*} 
whose complex lies in the $q$-th row of the double complex of lemma \ref{Hor p-page}. Assembling the row-wise van Est maps relating there, we get
\begin{eqnarray*}
\xymatrix{
\Phi_{row}:C(H^q\times G^r,W) \ar[r] & C(H^q,\bigwedge^r\gg^*\otimes W) . 
} 
\end{eqnarray*}
We claim that both $\Phi_{col}$ and $\Phi_{row}$ are maps of double complexes. To prove so, we write explicit formulas for them. \\
In the sequel, let $\omega\in C(H^q\times G^r,W)$. Consider $y_1,...,y_q\in\Gamma(\hh\ltimes G^r)$ and let $\lbrace u_a\rbrace$ be a basis for $\hh$, then for $\vec{g}\in G^r$
\begin{eqnarray*}
y_j(\vec{g})=y_j^a(\vec{g})u_a ,
\end{eqnarray*}  
where we are using the Einstein summation convention. Now, the $s$-fibres of the transformation groupoid are $s^{-1}(\vec{g})=H\times\lbrace\vec{g}\rbrace$; therefore, right multiplication by an element $(h;\vec{g})$ is defined as
\begin{eqnarray*}
\xymatrix{
R_{(h;\vec{g})}:H\times\lbrace(\vec{g})^h\rbrace \ar[r] & H\times\lbrace\vec{g}\rbrace :(h';(\vec{g})^h) \ar@{|->}[r] & (hh';\vec{g}).
} 
\end{eqnarray*} 
Consequently, for each $j$,
\begin{align*}
(\vec{y}_j)_{(h;\vec{g})} & :=d_{(1;(\vec{g})^h)}R_{(h;\vec{g})}(y_j^a((\vec{g})^h)u_a) \\
					      & =y_j^a((\vec{g})^h)\frac{d}{d\lambda_a}\rest{\lambda_a=0}R_{(h;\vec{g})}(\exp(\lambda_a u_a);(\vec{g})^h) \\
					      & =y_j^a((\vec{g})^h)\frac{d}{d\lambda_a}\rest{\lambda_a=0}(h\exp(\lambda_a u_a);\vec{g})=y_j^a((\vec{g})^h)dL_h(u_a)\in T_hH\leq T_hH\oplus T_{\vec{g}}G^r. 
\end{align*}
Let $h_2,...,h_q\in H$, then
\begin{align*}
R_{y_q}\omega(h_2,...,h_q;\vec{g}) & =\vec{y}_q\omega_{(h_2,...,h_q;\vec{g})}(1;(\vec{g})^{h_q...h_2}) \\
  								   & =d_{(1;(\vec{g})^{h_q...h_2})}\omega_{(h_2,...,h_q;\vec{g})}((\vec{y}_q)_{(1;(\vec{g})^{h_q...h_2})}) \\
  								   & =y_q^a((\vec{g})^{h_q...h_2})d_{(1;(\vec{g})^{h_q...h_2})}\omega_{(h_2,...,h_q;\vec{g})}(u_a) \\
  								   & =y_q^a((\vec{g})^{h_q...h_2})\frac{d}{d\lambda_a}\rest{\lambda_a=0}\omega(\exp(\lambda_au_a),h_2,...,h_q;\vec{g}).
\end{align*}
Moving on,
\begin{align*}
R_{y_{q-1}}R_{y_q}\omega(h_3,...,h_q;\vec{g}) & =\vec{y}_{q-1}(R_{y_q}\omega)_{(h_3,...,h_q;\vec{g})}(1;(\vec{g})^{h_q...h_3}) \\
  								   & =d_{(1;(\vec{g})^{h_q...h_3})}(R_{y_q}\omega)_{(h_3,...,h_q;\vec{g})}((\vec{y}_{q-1})_{(1;(\vec{g})^{h_q...h_3})}) \\
  								   & =y_{q-1}^b((\vec{g})^{h_q...h_3})d_{(1;(\vec{g})^{h_q...h_3})}(R_{y_q}\omega)_{(h_3,...,h_q;\vec{g})}(u_b) \\
  								   & =y_{q-1}^b((\vec{g})^{h_q...h_3})\frac{d}{d\lambda_{b}}\rest{\lambda_{b}=0}R_{y_q}\omega(\exp(\lambda_{b}u_b),h_3,...,h_q;\vec{g}),
\end{align*}
but
\begin{align*}
\frac{d}{d\lambda_{b}} & \rest{\lambda_{b}=0}R_{y_q}\omega(\exp(\lambda_{b}u_b),h_3,...,h_q;\vec{g}) \\
 & =\frac{d}{d\lambda_{b}}\rest{\lambda_{b}=0}\Big{(}y_q^a((\vec{g})^{h_q...h_3\exp(\lambda_{b}u_b)})\frac{d}{d\lambda_a}\rest{\lambda_a=0}\omega(\exp(\lambda_a u_a),\exp(\lambda_{b}u_b),h_3,...,h_q;\vec{g})\Big{)} \\
 & =\Big{(}\frac{d}{d\lambda_{b}}\rest{\lambda_{b}=0}y_q^a((\vec{g})^{h_q...h_3\exp(\lambda_{b}u_b)})\Big{)}\frac{d}{d\lambda_a}\rest{\lambda_a=0}\omega(\exp(\lambda_a u_a),\exp(0),h_3,...,h_q;\vec{g})+ \\
 & \quad +y_q^a((\vec{g})^{h_q...h_3\exp(0)})\frac{d}{d\lambda_{b}}\rest{\lambda_{b}=0}\frac{d}{d\lambda_a}\rest{\lambda_a=0}\omega(\exp(\lambda_a u_a),\exp(\lambda_{b}u_b),h_3,...,h_q;\vec{g}) ;
\end{align*}
therefore,
\begin{align*}
R_{y_{q-1}}R_{y_q}\omega(h_3,..., & h_q;\vec{g})= \\
     							  & (y_{q-1}^by_q^a)((\vec{g})^{h_q...h_3})\frac{d}{d\lambda_{b}}\rest{\lambda_{b}=0}\frac{d}{d\lambda_a}\rest{\lambda_a=0}\omega(\exp(\lambda_a u_a),\exp(\lambda_{b}u_b),h_3,...,h_q;\vec{g}).
\end{align*}
Inductively,
\begin{align*}
R_{y_{1}}...R_{y_q}\omega(\vec{g}) & =(y_{1}^{a_1}...y_q^{a_q})(\vec{g})\frac{d^I}{d\lambda_{I}}\rest{\lambda=0}\omega(\exp(\lambda_qu_{a_q}),...,\exp(\lambda_{1}u_{a_1});\vec{g}).
\end{align*}
Of course, if we let $Y=(y_1,...,y_q)$ and we recur to the notation that we have been using, this is equal to $\frac{d^I}{d\lambda_{I}}\rest{\lambda=0}\omega(\exp(\lambda\cdot Y(\vec{g}));\vec{g})$ and we conclude
\begin{align*}
\Phi_{col}\omega(y_1,...,y_q;\vec{g}) & =\sum_{\sigma\in S_q}\abs{\sigma}\frac{d^I}{d\lambda_{I}}\rest{\lambda=0}\omega(\exp(\lambda_{\sigma(q)}y_{\sigma(q)}(\vec{g})),...,\exp(\lambda_{\sigma(1)}y_{\sigma(1)}(\vec{g}));\vec{g}).
\end{align*}
As for $\Phi_{row}$, consider $x_1,...,x_r\in\Gamma(H^q\times\gg)$ and $\vec{h}\in H^q$. The $s$-fibres in the Lie group bundle are $s^{-1}(\vec{h})=\lbrace\vec{h}\rbrace\times G$; therefore, right multiplication by an element $(\vec{h};g)$ is defined as
\begin{eqnarray*}
\xymatrix{
R_{(\vec{h};g)}:\lbrace\vec{h}\rbrace\times G \ar[r] & \lbrace\vec{h}\rbrace\times G :(\vec{h};g') \ar@{|->}[r] & (\vec{h};gg').
} 
\end{eqnarray*} 
Consequently, for each $k$,
\begin{align*}
(\vec{x}_k)_{(\vec{h};g)} & :=d_{(\vec{h};1)}R_{(h;\vec{g})}(x_k(\vec{h})) \\
					      & =\frac{d}{d\tau_k}\rest{\tau_k=0}R_{(h;\vec{g})}(\vec{h};\exp(\tau_k x_k(\vec{h}))) \\
					      & =\frac{d}{d\tau_k}\rest{\tau_k=0}(\vec{h};\exp(\tau_k x_k(\vec{h}))g)=dR_g(x_k(\vec{h}))\in T_gG\leq T_{\vec{h}}H^q\oplus T_{g}G. 
\end{align*}
From this, it is immediate that the desired formula is
\begin{align*}
\Phi_{row}\omega(\vec{h};x_1,...,x_r) & =\sum_{\varrho\in S_r}\abs{\varrho}\frac{d^J}{d\tau_{J}}\rest{\tau=0}\omega(\vec{h};\exp(\tau_{\varrho(r)}y_{\varrho(r)}(\vec{h})),...,\exp(\tau_{\varrho(1)}y_{\varrho(1)}(\vec{h}))).
\end{align*}
\begin{lemma}\label{partialLAVanEsts}
$\Phi_{col}$ and $\Phi_{row}$ are maps of double complexes.
\end{lemma}
\begin{proof}
Since, respectively, $\Phi_{col}$ and $\Phi_{row}$ are assembled by column-wise and row-wise complex maps, the lemma follows from proving the compatibility with the other differentials, i.e. we are to prove that
\begin{eqnarray*}
\Phi_{col}\delta'=\delta'\Phi_{col} & \textnormal{and} & \Phi_{row}\delta=\delta\Phi_{row} .
\end{eqnarray*}
Let $Y=(y_1,...,y_q)\in\hh^q$ and $\vec{g}=(g_0,...,g_r)\in G^r$, and consider
\begin{align*}
\overrightarrow{R}_Y\delta'\omega(\vec{g}) & =\frac{d^I}{d\lambda_{I}}\rest{\lambda=0}\delta'\omega(\exp(\lambda\cdot Y);\vec{g}) \\
 								           & =\frac{d^I}{d\lambda_{I}}\rest{\lambda=0}\rho_G^{q}(\exp(\lambda\cdot Y);g_0)\omega(\exp(\lambda\cdot Y);\delta'_0\vec{g})+\sum_{k=1}^{r+1}(-1)^{k}\omega(\exp(\lambda\cdot Y);\delta'_k\vec{g}) \\
 							               & =\Big{(}\frac{d^I}{d\lambda_{I}}\rest{\lambda=0}\rho_G^{q}(\exp(\lambda\cdot Y);g_0)\omega(\exp(\lambda\cdot Y);\delta'_0\vec{g})\Big{)}+\sum_{k=1}^{r+1}(-1)^{k}\overrightarrow{R}_Y\omega(\delta'_k\vec{g}) .
\end{align*}
However,
\begin{align*}
\frac{d}{d\lambda_{q}} & \rest{\lambda_q=0}\rho_G^{q}(\exp(\lambda\cdot Y);g_0)\omega(\exp(\lambda\cdot Y);\delta'_0\vec{g}) \\
 & =\Big{(}\frac{d}{d\lambda_{q}}\rest{\lambda_q=0}\rho_G^{q}(\exp(\lambda\cdot Y);g_0)\Big{)}\omega(\exp(0),\exp(\lambda_{q-1}y_{q-1}),...,\exp(\lambda_{1}y_{1});\delta'_0\vec{g}) + \\
 & \qquad + \rho_G^q(\exp(0),\exp(\lambda_{q-1}y_{q-1}),...,\exp(\lambda_{1}y_{1});g_0)\frac{d}{d\lambda_{q}}\rest{\lambda_q=0}\omega(\exp(\lambda_{q}y_{q}),...,\exp(\lambda_{1}y_{1});\delta'_0\vec{g})\\
 & =\rho_G(g_0^{\exp(\lambda_{1}y_{1})...\exp(\lambda_{q-1}y_{q-1})})\frac{d}{d\lambda_{q}}\rest{\lambda_q=0}\omega(\exp(\lambda_{q}y_{q}),...,\exp(\lambda_{1}y_{1});\delta'_0\vec{g});
\end{align*}
therefore, inductively,
\begin{align*}
\frac{d^I}{d\lambda_{I}}\rest{\lambda=0}\rho_G^{q}(\exp(\lambda\cdot Y);g_0)\omega(\exp(\lambda\cdot Y);\delta'_0\vec{g})=\rho_G(g_0)\overrightarrow{R}_Y\omega(\delta'_0\vec{g}).
\end{align*}
Taking the alternating sum over $S_q$ and recalling $\rho_G^q(Y;g_0)=\rho_G(g_0)$,
\begin{align*}
\Phi_{col}\delta'\omega(Y;\vec{g}) & =\rho_G^q(Y;g_0)\Phi_{col}\omega(Y;\delta'_0\vec{g})+\sum_{k=1}^{r+1}(-1)^{k}\Phi_{col}\omega(Y;\delta'_k\vec{g})=\delta'\Phi_{col}\omega(Y;\vec{g}).
\end{align*}
As for the other equation, let $\vec{h}=(h_0,...,h_q)\in H^q$ and $X=(x_1,...,x_r)\in\gg^r$, and consider
\begin{align*}
\overrightarrow{R}_X\delta\omega(\vec{h}) & =\frac{d^J}{d\tau_{J}}\rest{\tau=0}\delta\omega(\vec{h};\exp(\tau\cdot X)) \\
 								          & =\frac{d^J}{d\tau_{J}}\rest{\tau=0}\rho_H^{r}(h_0;\exp(\tau\cdot X)^{h_q...h_1})\omega(\delta_0\vec{h};\exp(\tau\cdot X))+ \\
 								          & \qquad +\sum_{j=1}^{q}(-1)^{j}\omega(\delta_j\vec{h};\exp(\tau\cdot X))+(-1)^{q+1}\omega(\delta_{q+1}\vec{h};\exp(\tau\cdot X)^{h_q}) \\
 							              & =\rho_H(h_0)\overrightarrow{R}_X\omega(\delta_0\vec{h})+\sum_{j=1}^{q}(-1)^{j}\overrightarrow{R}_X\omega(\delta_j\vec{h})+(-1)^{q+1}\overrightarrow{R}_{X^{h_q}}\omega(\delta_{q+1}\vec{h}) .
\end{align*}
In the last line, we used that by definition $\rho_H^{r}(h_0;\exp(\tau\cdot X)^{h_q...h_1})=\rho_H(h_0)$. Since $\rho_H^{r}(h_0;X^{h_q...h_1})$ is also $\rho_H(h_0)$, taking the alternating sum over $S_r$, we get
\begin{align*}
\Phi_{row}\delta\omega(\vec{h};X) & =\rho_H^r(h_0;X^{h_q...h_1})\Phi_{row}\omega(\delta_0\vec{h};X)+ \\
								  & \qquad +\sum_{j=1}^{q}(-1)^{j}\Phi_{row}\omega(\delta_j\vec{h};X)+(-1)^{q+1}\Phi_{row}\omega(\delta_{q+1}\vec{h};X^{h_q}) \\
								  & =\delta\Phi_{row}\omega(\vec{h};X),
\end{align*}
and the desired result.

\end{proof}
We now turn to the $2$-van Est map restricted to the $p$-pages of the triple complex of Lie $2$-group $\G=\xymatrix{G \ar[r]^i & H}$ with values in a $2$-representation $\rho$ on the $2$-vector space $\xymatrix{W \ar[r]^\phi & V}$, and we use the above lemmas to get a first approximation of its cohomology. Since we saw that the $p$-pages arose from considering the double complex associated to the double Lie group map
\begin{eqnarray*}
\xymatrix{
\G_p\ltimes G \ar[r] & GL(W)\ltimes GL(W):(\gamma ;g) \ar@{|->}[r] & (\rho_0^1(t_p(\gamma))^{-1},\rho_0^1(i(g))) ,
}
\end{eqnarray*} 
by appropriately replacing the first column, we insert a first column replacement for the associated \LA -double complexes.
\begin{lemma}\label{VerLA r-cx}
Let $\xi_1,...,\xi_q\in\gg_p$ and $g\in G$ and $\omega\in\bigwedge^q\gg_p^*\otimes V$. Set
\begin{eqnarray*}
\partial'\omega(\xi_1,...,\xi_q;g)=\rho_1(g)\omega(\xi_1,...,\xi_q),
\end{eqnarray*}
then the first term of the $q$th row of the vertical \LA -double complex associated to a $p$-pages of the triple complex can be replaced by 
\begin{eqnarray*}
\xymatrix{
\bigwedge^q\gg_p^*\otimes V \ar[r]^{\partial'\quad} & C(G,\bigwedge^q\gg_p^*\otimes W) \ar[r]^{\delta'} & C(G^2,\bigwedge^q\gg_p^*\otimes W) \ar[r] & ...
}
\end{eqnarray*}
still yielding a complex.
\end{lemma}
\begin{proof}
We prove that, for $\omega\in \bigwedge^q\gg_p^*\otimes V$, $\delta'\partial'\omega=0$.
\begin{align*}
\delta'\partial' & \omega(\xi_1,...,\xi_q;g_0,g_1) \\
                 & =\rho_G^{q}(\xi_1,...,\xi_q;g_0)\partial'\omega(\xi_1,...,\xi_q;g_1)-\partial'\omega(\xi_1,...,\xi_q;g_0g_1)+\partial'\omega(\xi_1,...,\xi_q;g_0) \\
                 & =\rho_0^1(i(g_0))\rho_1(g_1)\omega(\xi_1,...,\xi_q)-\rho_1(g_0g_1)\omega(\xi_1,...,\xi_q)+\rho_1(g_0)\omega(\xi_1,...,\xi_q) \\
                 & =\rho_1(g_1)\omega(\xi_1,...,\xi_q)+\rho_1(g_0)\phi\rho_1(g_1)\omega(\xi_1,...,\xi_q)-\rho_1(g_0g_1)\omega(\xi_1,...,\xi_q)+\rho_1(g_0)\omega(\xi_1,...,\xi_q), 
\end{align*} 
and the result follows from the fact that $\rho_1$ is a Lie group homomorphism landing in $GL(\phi)_1$.

\end{proof}
\begin{lemma}\label{HorLA r-cx}
Let $\gamma_1,...,\gamma_q\in\G_p$ and $x\in\gg$ and $\omega\in C(\G_p^q,V)$. Set
\begin{eqnarray*}
\partial'\omega(\gamma_1,...,\gamma_q;x)=\rho_0^1(t_p(\gamma_q...\gamma_1))^{-1}\dot{\rho}_1(x)\omega(\gamma_1,...,\gamma_q),
\end{eqnarray*}
then the first term of the $q$th row of the horizontal \LA -double complex associated to a $p$-pages of the triple complex can be replaced by 
\begin{eqnarray*}
\xymatrix{
C(\G_p^q,V) \ar[r]^{\partial'\quad} & C(\G_p^q,\gg^*\otimes W) \ar[r]^{\delta'} & C(\G_p^q,\bigwedge^2\gg^*\otimes W) \ar[r] & ...
}
\end{eqnarray*}
still yielding a complex.
\end{lemma}
\begin{proof}
We prove that, for $\omega\in C(\G_p^q,V)$, $\delta'\partial'\omega=0$. Now,
\begin{align*}
\delta'\partial'\omega(\gamma_1,...,\gamma_q;x_0,x_1) & =-\partial'\omega(\gamma_1,...,\gamma_q;[x_0,x_1])+ \\
													  & \qquad\qquad+\rho_\gg^{q}(x_0)\partial'\omega(\gamma_1,...,\gamma_q;x_1)-\rho_\gg^{q}(x_1)\partial'\omega(\gamma_1,...,\gamma_q;x_0) ,
\end{align*}
but 
\begin{align*}	
\rho_\gg^{q}(x_m)\partial'\omega(\gamma_1,...,\gamma_q;x_n) & =\rho_\gg(x_m^{t_p(\gamma_q...\gamma_1)})\rho_0^1(t_p(\gamma_q...\gamma_1))^{-1}\dot{\rho}_1(x_n)\omega(\gamma_1,...,\gamma_q)
\end{align*}	
and $\rho_\gg=\dot{\rho_G}$; therefore,
\begin{align*}
\rho_\gg(x_m^{t_p(\gamma_q...\gamma_1)}) & =\frac{d}{d\tau}\rest{\tau=0}\rho_G(\exp(x_m^{t_p(\gamma_q...\gamma_1)})) \\
 & =\frac{d}{d\tau}\rest{\tau=0}\rho_0^1(i(\exp(x_m)^{t_p(\gamma_q...\gamma_1)})) \\
 & =\frac{d}{d\tau}\rest{\tau=0}\rho_0^1(t_p(\gamma_q...\gamma_1)^{-1}i(\exp(x_m))t_p(\gamma_q...\gamma_1)) \\
 & =\rho_0^1(t_p(\gamma_q...\gamma_1))^{-1}\dot{\rho}_0^1(\mu(x_m))\rho_0^1(t_p(\gamma_q...\gamma_1)).
\end{align*}
Hence,
\begin{align*}							  
\delta'\partial'\omega(\gamma_1,...,\gamma_q;x_0,x_1) & =\rho_0^1(t_p(\gamma_q...\gamma_1))^{-1}\Big{[}-\dot{\rho}_1([x_0,x_1])\omega(\gamma_1,...,\gamma_q) \\
													  & \qquad +\dot{\rho}_0^1(\mu(x_0))\dot{\rho}_1(x_1)\omega(\gamma_1,...,\gamma_q)-\dot{\rho}_0^1(\mu(x_1))\dot{\rho}_1(x_0)\omega(\gamma_1,...,\gamma_q)\Big{]} ,											
\end{align*} 
and the result follows from $\dot{\rho}_0^1\circ\mu=\dot{\rho}_1\circ\phi$ and the fact that $\dot{\rho}_1$ is a Lie algebra homomorphism landing in $\ggl(\phi)_1$.

\end{proof}
These latter lemmas will help us build the double complexes where the column and row-wise assembly of van Est maps should land when coming out of the $p$-pages of the triple complex. Thereafter, we prove that replacing the value of the $\Phi_{col}$ and $\Phi_{row}$ for the first columns still yields a map of double complexes.
\begin{prop}\label{(q,r)-VerLAdoubleCx}
For each $p$, 
\begin{eqnarray*}
\xymatrix{
\vdots                                                             & \vdots                                                         & \vdots                                                   & \\ 
\bigwedge^3\gg_p^*\otimes V \ar[r]^{\partial'\quad} \ar[u]         & C(G,\bigwedge^3\gg_p^*\otimes W) \ar[r]^{\delta'}\ar[u]          & C(G^2,\bigwedge^3\gg_p^*\otimes W) \ar[r]\ar[u]          & \dots \\
\bigwedge^2\gg_p^*\otimes V \ar[r]^{\partial'\quad}\ar[u]^{\delta} & C(G,\bigwedge^2\gg_p^*\otimes W) \ar[r]^{\delta'}\ar[u]^{\delta} & C(G^2,\bigwedge^2\gg_p^*\otimes W) \ar[r]\ar[u]^{\delta} & \dots \\
\gg_p^*\otimes V  \ar[r]^{\partial'\quad}\ar[u]^{\delta}           & C(G,\gg_p^*\otimes W) \ar[r]^{\delta'}\ar[u]^{\delta}            & C(G^2,\gg_p^*\otimes W) \ar[r]\ar[u]^{\delta}            & \dots \\
V \ar[r]^{\delta_{(1)}}\ar[u]^{\delta}                             & C(G,W) \ar[r]^{\delta'}\ar[u]^{\delta}                         & C(G^2,W)	\ar[r]\ar[u]^{\delta}                          & \dots  
}
\end{eqnarray*}
is a double complex.
\end{prop}
\begin{proof}
Due to lemmas \ref{G-coh[V]} and \ref{VerLA r-cx}, each row is a complex, and clearly so is each column. Disregarding the first column of squares, lemma \ref{Ver p-page} says that we have got a double complex. Now, in order to finish the proof, one needs to check that the generic first column square 
\begin{eqnarray*}
\xymatrix{ 
\bigwedge^{q+1}\gg_p^*\otimes V \ar[r]^{\partial'\quad}                & C(G,\bigwedge^{q+1}\gg_p^*\otimes W)          \\
\bigwedge^{q}\gg_p^*\otimes V   \ar[r]^{\partial'\quad}\ar[u]^{\delta} & C(G,\bigwedge^{q}\gg_p^*\otimes W) \ar[u]^{\delta} 
}
\end{eqnarray*}
commutes. First, for $q=0$, let $\xi\in\gg_p$, $g\in G$ and $v\in V$,
\begin{align*}
\partial'\delta v(\xi;g) & =\rho_1(g)\delta v(\xi)=\rho_1(g)\dot{\rho}_0^0(\hat{t}_p(\xi))v,
\end{align*}
whereas
\begin{align*}
\delta\delta_{(1)}v(\xi;g) & =\rho_{\gg_p}^1(\xi)\delta_{(1)}v(g)=\rho_{\gg_p}^1(\xi)\rho_1(g)v.
\end{align*}					   
Since the anchor of $\gg_p\ltimes G$ is given by
\begin{align*}
d_{(1;g)}\Lf{t}(\xi,0) & =\frac{d}{d\lambda}\rest{\lambda=0}g^{t_p(\exp(\lambda\xi))}\in T_{g}G 
\end{align*}
and for $(\xi;g)\in\gg_p\ltimes G$, $\hat{t}(\xi;g)=\xi$, then, by definition
\begin{eqnarray*}
\rho_{\gg_p}^1(\xi)\rho_1(g)v=\rho_{\gg_p}(\xi)\rho_1(g)v+\frac{d}{d\lambda}\rest{\lambda=0}\rho_1(g^{t_p(\exp(\lambda\xi))})v .
\end{eqnarray*}
Computing,
\begin{align*}
\frac{d}{d\lambda} & \rest{\lambda=0}\rho_1(g^{t_p(\exp(\lambda\xi))})v \\
                   & =\frac{d}{d\lambda}\rest{\lambda=0}\rho_0^1(t_p(\exp(\lambda\xi)))^{-1}\rho_1(g)\rho_0^0(t_p(\exp(\lambda\xi))v \\
                   & =(\frac{d}{d\lambda}\rest{\lambda=0}\rho_0^1(\exp(\lambda\hat{t}_p(\xi))^{-1}))\rho_1(g)\rho_0^0(t_p(\exp(0)))v+\rho_0^1(t_p(\exp(0)))^{-1}\rho_1(g)\frac{d}{d\lambda}\rest{\lambda=0}\rho_0^0(\exp(\lambda\hat{t}_p(\xi))v \\
                   & =-\dot{\rho}_0^1(\hat{t}_p(\xi))\rho_1(g)v+\rho_1(g)\dot{\rho}_0^0(\hat{t}_p(\xi))v,
\end{align*}
so the equality follows from $\rho_{\gg_p}(\xi)=\dot{\rho}_0^1(\hat{t}_p(\xi))$. \\
As for the other values of $q$, let $\Xi=(\xi_0,...,\xi_q)\in\gg_p^{q+1}$, then
\begin{align*}
\partial'\delta\omega(\Xi;g) & =\rho_1(g)\delta\omega(\Xi) \\
							 & =\rho_1(g)\Big{(}\sum_{m<n}(-1)^{m+n}\omega([\xi_m,\xi_n],\Xi(m,n))+\sum_{j=0}^{q}(-1)^{j+1}\dot{\rho}_0^0(\hat{t}_p(\xi_j))\omega(\Xi(j)))\Big{)}.
\end{align*}
On the other hand,
\begin{align*}
\delta\partial'\omega(\Xi;g) & =\sum_{m<n}(-1)^{m+n}\partial'\omega([\xi_m,\xi_n],\Xi(m,n);g)+\sum_{j=0}^{q}(-1)^{j+1}\rho_{\gg_p}^{1}(\xi_j)\partial'\omega(\Xi(j);g)\\
							 & =\sum_{m<n}(-1)^{m+n}\rho_1(g)\omega([\xi_m,\xi_n],\Xi(m,n))+\sum_{j=0}^{q}(-1)^{j+1}\rho_{\gg_p}^{1}(\xi_j)\rho_1(g)\omega(\Xi(j)),
\end{align*}
but we saw that $\rho_{\gg_p}^{1}(\xi_j)\rho_1(g)=\rho_1(g)\dot{\rho}_0^0(\hat{t}_p(\xi_j))$; thus, the result follows.

\end{proof}
\begin{prop}\label{(q,r)-HorLAdoubleCx}
For each $p$, 
\begin{eqnarray*}
\xymatrix{
\vdots                                             & \vdots                                                   & \vdots   & \\ 
C(\G_p^3,V) \ar[r]^{\partial'\quad} \ar[u]         & C(\G_p^3,\gg^*\otimes W) \ar[r]^{\delta'}\ar[u]          & C(\G_p^3,\bigwedge^2\gg^*\otimes W) \ar[r]\ar[u]          & \dots \\
C(\G_p^2,V) \ar[r]^{\partial'\quad}\ar[u]^{\delta} & C(\G_p^2,\gg^*\otimes W) \ar[r]^{\delta'}\ar[u]^{\delta} & C(\G_p^2,\bigwedge^2\gg^*\otimes W) \ar[r]\ar[u]^{\delta} & \dots \\
C(\G_p,V)  \ar[r]^{\partial'\quad}\ar[u]^{\delta}  & C(\G_p,\gg^*\otimes W) \ar[r]^{\delta'}\ar[u]^{\delta}   & C(\G_p,\bigwedge^2\gg^*\otimes W) \ar[r]\ar[u]^{\delta} & \dots \\
V \ar[r]^{\delta_{(1)}}\ar[u]^{\delta}             & \gg^*\otimes W \ar[r]^{\delta'}\ar[u]^{\delta}           & \bigwedge^2\gg^*\otimes W \ar[r]\ar[u]^{\delta}                  & \dots  
}
\end{eqnarray*}
is a double complex.
\end{prop}
\begin{proof}
Due to lemmas \ref{Alg r-cx} and \ref{HorLA r-cx}, each row is a complex, and clearly so is each column. Disregarding the first column of squares, lemma \ref{Hor p-page} says that we have got a double complex. Now, in order to finish the proof, one needs to check that the generic first column square 
\begin{eqnarray*}
\xymatrix{ 
C(\G_p^{q+1},V) \ar[r]^{\partial'\quad}                & C(\G_p^{q+1},\gg^*\otimes W)          \\
C(\G_p^q,V)     \ar[r]_{\partial'\quad}\ar[u]^{\delta} & C(\G_p^{q},\gg^*\otimes W) \ar[u]_{\delta} 
}
\end{eqnarray*}
commutes. First, for $q=0$, let $\gamma\in\G_p$, $x\in\gg$ and $v\in V$,
\begin{align*}
\partial'\delta v(\gamma;x) & =\rho_0^1(t_p(\gamma))^{-1}\dot{\rho}_1(x)\delta v(\gamma)=\rho_0^1(t_p(\gamma))^{-1}\dot{\rho}_1(x)(v-\rho_0^0(t_p(\gamma))v),
\end{align*}
whereas
\begin{align*}
\delta\delta_{(1)}v(\gamma;x) & =\rho_{\G_p}^1(\gamma;x)\delta_{(1)}v(x)-\delta_{(1)}v(x^{t_p(\gamma)}) \\
							  & =\rho_0^1(t_p(\gamma))^{-1}\dot{\rho}_1(x)v-\dot{\rho}_1(x^{t_p(\gamma)})v,
\end{align*}					   
so the equality follows from the equivariance $\dot{\rho}_1(x^{t_p(\gamma)})=\rho_0^1(t_p(\gamma))^{-1}\dot{\rho}_1(x)\rho_0^0(t_p(\gamma))$. \\
As for the other values of $q$, let $\vec{\gamma}=(\gamma_0,...,\gamma_q)\in\G_p^{q+1}$, then
\begin{align*}
\partial'\delta\omega(\vec{\gamma};x) & =\rho_0^1(t_p(\gamma_q...\gamma_1))^{-1}\dot{\rho}_1(x)\delta\omega(\vec{\gamma}) \\
							          & =\rho_0^1(t_p(\gamma_q...\gamma_1))^{-1}\dot{\rho}_1(x)\Big{(}\sum_{j=0}^{q}(-1)^{j}\omega(\delta_j\vec{\gamma})+(-1)^{q+1}\rho_0^0(t_p(\gamma_q))\omega(\delta_{q+1}\vec{\gamma})\Big{)} ,
\end{align*}
and
\begin{align*}
\delta\partial'\omega(\vec{\gamma};x) & =\rho_{\G_p}^1(\gamma_0;x^{\gamma_q...\gamma_1})\partial'\omega(\delta_0\vec{\gamma};x)+ \\
							          & \qquad\qquad +\sum_{j=1}^{q}(-1)^{j}\partial'\omega(\delta_j\vec{\gamma};x)+(-1)^{q+1}\partial'\omega(\delta_{q+1}\vec{\gamma};x^{t_p(\gamma_q)}) \\
							          & =\rho_0^1(t_p(\gamma_0))^{-1}\rho_0^1(t_p(\gamma_q...\gamma_1))^{-1}\dot{\rho}_1(x)\omega(\delta_0\vec{\gamma})+ \\
							          & \qquad\qquad +\sum_{j=1}^{q}(-1)^{j}\rho_0^1(t_p(\gamma_q...(\gamma_{j+1}\gamma_j)\gamma_0))^{-1}\dot{\rho}(x)\omega(\delta_j\vec{\gamma})+ \\
							          & \qquad\qquad +(-1)^{q+1}\rho_0^1(t_p(\gamma_{q-1}...\gamma_0))^{-1}\dot{\rho}_1(x^{t_p(\gamma_q)})\omega(\delta_{q+1}\vec{\gamma}) \\
							          & =\rho_0^1(t_p(\gamma_q...\gamma_0))^{-1}\dot{\rho}_1(x)\Big{(}\sum_{j=0}^{q}(-1)^{j}\omega(\delta_j\vec{\gamma})+(-1)^{q+1}\rho_0^0(t_p(\gamma_q))\omega(\delta_{q+1}\vec{\gamma})\Big{)},
\end{align*}
where the last line is again due to equivariance.

\end{proof}
We will refer to the double complexes of propositions \ref{(q,r)-VerLAdoubleCx} and \ref{(q,r)-HorLAdoubleCx}, respectively, as the vertical and horizontal \LA -double complex associated to the $p$-page, and will write them
\begin{eqnarray*}
C_{\LA}(\gg_p\ltimes G,\phi) & \textnormal{and} & C_{\LA}(\G_p\ltimes\gg,\phi).
\end{eqnarray*}
Let 
\begin{eqnarray*}
\xymatrix{
\Phi_V :C^{p,\bullet}_\bullet(\G,\phi) \ar[r] & C_{\LA}(\gg_p\ltimes G,\phi)
}
\end{eqnarray*}
be defined by assembling column-wise van Est maps, i.e. $\Phi_V$ coincides with $\Phi_{col}$ outside the first column, and in the first column it is the van Est map for the Lie group $\G_p$. Analogously, let 
\begin{eqnarray*}
\xymatrix{
\Phi_H :C^{p,\bullet}_\bullet(\G,\phi) \ar[r] & C_{\LA}(\G_p\ltimes\gg,\phi)
}
\end{eqnarray*}
be defined by assembling row-wise van Est maps, i.e. $\Phi_H$ coincides with $\Phi_{row}$ outside the first column, and restricted to the first column $\Phi_H$ is the identity of the complex of cochains of the Lie group $\G_p$.
\begin{prop}
$\Phi_V$ and $\Phi_H$ are maps of double complexes.
\end{prop}
\begin{proof}
First of all, $\Phi_V$ and $\Phi_H$ define maps of complexes when restricted to the first column by the very definition. Now, in sight of lemma \ref{partialLAVanEsts}, we just need to prove that the maps are compatible with the horizontal differential of the first column, i.e. Let $\omega\in C(\G_p^q,V)$, we prove that
\begin{eqnarray*}
\Phi_V\partial'\omega =\partial'\Phi_V\omega & \textnormal{and} & \Phi_H\partial'\omega =\partial'\Phi_H\omega.
\end{eqnarray*}
Let $\Xi=(\xi_1,...,\xi_q)\in\gg_p^q$ and $g\in G$, then
\begin{align*}
\overrightarrow{R}_{\Xi}\partial'\omega(g) & =\frac{d^I}{d\lambda_I}\rest{\lambda=0}\partial'\omega(\exp(\lambda_q\xi_q),...,\exp(\lambda_1\xi_1);g) \\
										   & =\frac{d^I}{d\lambda_I}\rest{\lambda=0}\rho_0^1(t_p(\exp(\lambda_1\xi_1)...\exp(\lambda_q\xi_q)))^{-1}\rho_1(g)\omega(\exp(\lambda_q\xi_q),...,\exp(\lambda_1\xi_1)) .
\end{align*}
As usual, we compute 
\begin{align*}
\frac{d}{d\lambda_q}\rest{\lambda_q=0} & \rho_0^1(t_p(\exp(\lambda_1\xi_1)...\exp(\lambda_q\xi_q)))^{-1}\rho_1(g)\omega(\exp(\lambda_q\xi_q),...,\exp(\lambda_1\xi_1)) \\
 & =(\frac{d}{d\lambda_q}\rest{\lambda_q=0}\rho_0^1(t_p(\exp(\lambda_1\xi_1)...\exp(\lambda_q\xi_q)))^{-1})\rho_1(g)\omega(\exp(0),\exp(\lambda_{q-1}\xi_{q-1}),...,\exp(\lambda_1\xi_1))+ \\
 & \qquad +\rho_0^1(t_p(\exp(\lambda_1\xi_1)...\exp(\lambda_{q-1}\xi_{q-1})\exp(0)))^{-1}\rho_1(g)\frac{d}{d\lambda_q}\rest{\lambda_q=0}\omega(\exp(\lambda_q\xi_q),...,\exp(\lambda_1\xi_1)) \\
 & =\rho_0^1(t_p(\exp(\lambda_1\xi_1)...\exp(\lambda_{q-1}\xi_{q-1})))^{-1}\rho_1(g)R_{\xi_q}\omega(\exp(\lambda_{q-1}\xi_{q-1}),...,\exp(\lambda_1\xi_1)),
\end{align*}
and by induction deduce
\begin{eqnarray*}
\overrightarrow{R}_{\Xi}\partial'\omega(g)=\rho_1(g)\overrightarrow{R}_{\Xi}\omega ;
\end{eqnarray*}
thus, taking the alternating sum over $S_q$ yields the first equation. \\
As for the second equation, let $\vec{\gamma}=(\gamma_1,...,\gamma_q)\in\G_p^q$ and $x\in\gg$, then
\begin{align*}
\Phi_H\partial'\omega(\gamma;x) & =R_x\partial'\omega(\vec{\gamma}) \\
 								& =\frac{d}{d\tau}\rest{\tau=0}\rho_0^1(t_p(\gamma_q...\gamma_1))^{-1}\rho_1(\exp(\tau x))\omega(\vec{\gamma}) \\
 								& =\rho_0^1(t_p(\gamma_q...\gamma_1))^{-1}\dot{\rho}_1(x)\omega(\vec{\gamma})=\partial'\Phi_H\omega(\gamma;g).
\end{align*}
\end{proof}
\begin{cor}\label{Ver p-vanEst}
If $\G_p$ is $k$-connected,
\begin{eqnarray*}
H^n_{tot}(\Phi_V)=(0), & \textnormal{for all degrees }n\leq k.
\end{eqnarray*}
\end{cor}
\begin{proof}
We compute the cohomology of the mapping cone of $\Phi_V$ using the spectral sequence of the double complex filtrated by columns, whose first page is
\begin{eqnarray*}
\xymatrix{
           & \vdots                                        & \vdots                                        & \vdots                           &                      \\ 
           & H^2_{\delta_{\Phi_V}}(C^{0,\bullet}(\Phi_V)) \ar[r] & H^2_{\delta_{\Phi_V}}(C^{1,\bullet}(\Phi_V)) \ar[r] & H^2_{\delta_{\Phi_V}}(C^{2,\bullet}(\Phi_V)) \ar[r] & \dots \\
E^{p,q}_1: & H^1_{\delta_{\Phi_V}}(C^{0,\bullet}(\Phi_V)) \ar[r] & H^1_{\delta_{\Phi_V}}(C^{1,\bullet}(\Phi_V)) \ar[r] & H^1_{\delta_{\Phi_V}}(C^{2,\bullet}(\Phi_V)) \ar[r] & \dots \\
           & H^0_{\delta_{\Phi_V}}(C^{0,\bullet}(\Phi_V)) \ar[r] & H^0_{\delta_{\Phi_V}}(C^{1,\bullet}(\Phi_V)) \ar[r] & H^0_{\delta_{\Phi_V}}(C^{2,\bullet}(\Phi_V)) \ar[r] & \dots 
}
\end{eqnarray*}
Now, $\Phi_V$ is defined column-wise by van Est maps; consequently, the columns of the mapping cone double complex coincide with the mapping cones of these. By means of the Crainic-van Est theorem (cf. Theorem \ref{CrainicVanEst}), we know that the $r$th column of the latter diagram is zero below $k$, as $\G_p$  the $s$-fibre of 
\begin{eqnarray*}
\xymatrix{
\G_p\ltimes G^r \ar@<0.5ex>[r]\ar@<-0.5ex>[r] & G^r ,
}
\end{eqnarray*}
and it is $k$-connected by hypothesis. Given that
\begin{eqnarray*}
E^{p,q}_1\Rightarrow H^{p+q}_{tot}(\Phi_V),
\end{eqnarray*}
it follows from lemma \ref{BelowDiag} that $H_{tot}^n(\Phi_V)=(0)$ for $n\leq k$ as desired.

\end{proof}
We will also have that $\Phi_H$ has trivial cohomology; however, it will follow from a tweaked van Est theorem that takes into account that the space of $0$-cochains is not the usual one. Let us point out that, though the space of $0$-cochains changes, it is not unrelated to the usual one. Indeed, for all $(p,q)$, 
\begin{eqnarray*}
\xymatrix{
C(\G_p^q,W) \ar[rd]^{\delta'}\ar[dd]_{\phi_p^q} & \\
  & C(\G_p^q\times G,W), \\
C(\G_p^q,V) \ar[ur]_{\partial'} &
}
\end{eqnarray*}
where the vertical arrow is given by
\begin{eqnarray*}
(\phi_p^q\omega)(\vec{\gamma}):=\rho_0^0(t_p(\gamma_q...\gamma_1))\phi(\omega(\vec{\gamma}))
\end{eqnarray*}
for $\omega\in C(\G_p^q,W)$ and $\vec{\gamma}=(\gamma_1,...,\gamma_q)\in\G_p^q$, commutes. Let $g\in G$ and $\vec{\gamma}$ as before, then
\begin{align*}
\delta'\omega(\vec{\gamma};g) & =\rho_0^1(i(g^{t_p(\gamma_q...\gamma_1)}))\omega(\vec{\gamma})-\omega(\vec{\gamma}) \\
					          & =\omega(\vec{\gamma})+\rho_1(g^{t_p(\gamma_q...\gamma_1)})\phi(\omega(\vec{\gamma}))-\omega(\vec{\gamma}) \\
					          & =\rho_0^1(t_p(\gamma_q...\gamma_1))^{-1}\rho_1(g)\rho_0^0(t_p(\gamma_q...\gamma_1))\phi(\omega(\vec{\gamma}))=\partial'\phi_p^q(\omega)(\gamma;g).
\end{align*}
Here, the second equality holds, as $\rho$ is a $2$-representation of Lie $2$-groups and hence $\rho_0^1(i(g))=I+\rho_1(g)\phi$. This relation immediately implies
\begin{eqnarray*}
\textnormal{If } Z^0:=\ker\delta' ,\quad\phi(Z^0)\subseteq Z^0_\phi:=\ker\partial'\textnormal{ and }B^1:=\Img(\delta)\subseteq B^1_\phi :=\Img(\partial').
\end{eqnarray*} 
Also, for any $\ker\phi$-valued cochain $\omega\in C(\G_p^q,\ker\phi )$, $\omega\in Z^0$. 
Furthermore, one can argue analogously and conclude that 
\begin{eqnarray*}
\xymatrix{
C(\G_p^q,W) \ar[rd]^{\delta'}\ar[dd]_{\phi_p^q} & \\
  & C(\G_p^q,\gg^*\otimes W), \\
C(\G_p^q,V) \ar[ur]_{\partial'} &
}
\end{eqnarray*}
also commutes. Indeed, for $x\in\gg$ and $\vec{\gamma}$ as before, 
\begin{align*}
\delta'\omega(\vec{\gamma};x) & =\dot{\rho}_0^1(\mu(x^{t_p(\gamma_q...\gamma_1)}))\omega(\vec{\gamma}) \\
		    		          & =\rho_0^1(t_p(\gamma_q...\gamma_1))^{-1}\dot{\rho}_1(x)\rho_0^0(t_p(\gamma_q...\gamma_1))\phi(\omega(\vec{\gamma}))=\partial'\phi_p^q(\omega)(\gamma;x)
\end{align*}
This time around, the second equality holds, as $\dot{\rho}$ is a $2$-representation of Lie $2$-algebras and hence $\dot{\rho}_0^1(\mu(x))=\dot{\rho}_1(x)\phi$ together with the equivariance of $\dot{\rho}_1$. Bottom line is, this relation also implies
\begin{eqnarray*}
\textnormal{If }\dot{Z}^0:=\ker\delta' ,\quad\phi( \dot{Z}^0)\subseteq\dot{Z}^0_\phi:=\ker\partial'\textnormal{ and } \dot{Z}^1:=\Img(\delta)\subseteq\dot{Z}^1_\phi :=\Img(\partial');
\end{eqnarray*}
moreover, $C(\G_p^q,\ker\phi )\subseteq\dot{Z}^0$ and appropriately joining the triangles with $\Phi_{row}$ and $\Phi_H$ makes a commuting prism. \\
To emphasize where we are taking coefficients we write 
\begin{eqnarray*}
H^n(G,\phi):=Z^n/B^n_\phi & \textnormal{ and } & H^n(\gg ,\phi)=\dot{Z}^n/\dot{B}^n_\phi 
\end{eqnarray*}  
for $n\in\lbrace 0,1\rbrace$. 
\begin{lemma}\label{tweak0}
If $G$ is connected $\Phi_H$ induces an isomorphism between $H^0(G,\phi)$ and $H^0(\gg ,\phi)$.
\end{lemma}
\begin{proof}
Since $\Phi_H$ restricted to the first column is the identity and there are no cochains of negative degree, the statement of the lemma might as well be restated as
\begin{eqnarray*}
Z^0_\phi =\dot{Z}^0_\phi .
\end{eqnarray*}
($\subseteq$) Suppose $\omega\in Z^0_\phi$, that is
\begin{eqnarray*}
\partial'\omega(\vec{\gamma};g)=0
\end{eqnarray*}
for all $(\vec{\gamma};g)\in\G_p^q\times G$. Then, for $x\in\gg$, we clearly have
\begin{eqnarray*}
\partial'\omega(\vec{\gamma};x)=\Phi_H(\partial'\omega)(\vec{\gamma};x)=\frac{d}{d\tau}\rest{\tau =0}\partial'\omega(\vec{\gamma};\exp(\tau x))=0.
\end{eqnarray*}
($\supseteq$) Conversely, suppose $\omega\in\dot{Z}^0_\phi$, that is
\begin{eqnarray*}
\partial'\omega(\vec{\gamma};x)=0
\end{eqnarray*}
for all $(\vec{\gamma};x)\in\G_p^q\times\gg$, and given that $\rho_0^1(t_p(\gamma_q...\gamma_1))^{-1}$ is an isomorphism, this equation also implies $\dot{\rho}_1(x)\omega(\vec{\gamma})$. Now, it is well-known that the exponential restricts to a diffeomorphism from some neighbourhood of $0\in\gg$, say $U$, onto a neighbourhood of $1\in G$. Since $G$ is connected, it is generated by $\exp(U)\subset G$; and therefore, for all $g\in G$, there exist $x_1,...,x_n\in\gg$ such that
\begin{eqnarray*}
g=\exp(x_1)...\exp(x_n).
\end{eqnarray*}
Now, for all ranging values of $k$,
\begin{eqnarray*}
\rho_1(\exp_G(x_k))=\exp_{GL(\phi)_1}(\dot{\rho}_1(x_k));
\end{eqnarray*}
hence, recalling from proposition \ref{TheExpOfGL(phi)} that 
\begin{eqnarray*}
\exp_{GL(\phi)_1}(\dot{\rho}_1(x_k))=\sum_{n=0}^\infty\frac{(\dot{\rho}_1(x_k)\phi)^n}{(n+1)!}\dot{\rho}_1(x_k),
\end{eqnarray*}
we conclude that 
\begin{eqnarray*}
\rho_1(\exp_G(x_k))\omega(\vec{\gamma})=0
\end{eqnarray*}
for all $1\leq k\leq n$. In so, since 
\begin{align*}
\partial'\omega(\vec{\gamma};g) & =\rho_0^1(t_p(\gamma_q...\gamma_1))^{-1}\rho_1(g)\omega(\vec{\gamma}) \\
								& =\rho_0^1(t_p(\gamma_q...\gamma_1))^{-1}\rho_1(\exp(x_1)...\exp(x_n))\omega(\vec{\gamma}).
\end{align*}
and
\begin{align*}
\rho_1(\exp(x_1)...\exp(x_n))=\rho_1(\exp(x_1)... & \exp(x_{n-1}))+\rho_1(\exp(x_n))+ \\
												  & +\rho_1(\exp(x_1)...\exp(x_{n-1}))\phi\rho_1(\exp(x_n)),						
\end{align*}
inductively we get the desired inclusion.

\end{proof}
It is worth to point out that this latter lemma is a consequence of several pieces of Lie theory, and would not hold in a general scenario in which we have got quasi-isomorphic complexes with altered $0$-cochains in the fashion we described above. In contrast, the following lemma is stated as a general result of homological algebra and will imply naturally that if $G$ is $1$-connected, $\Phi_H$ induces an isomorphism between $H^1(G,\phi)$ and $H^1(\gg ,\phi)$. 
\begin{lemma}\label{tweak1}
Let $\xymatrix{\Phi:C^\bullet \ar[r] & D^\bullet}$ be a map of complexes inducing isomorphisms in cohomology for degrees $0$ and $1$. Suppose now there is a second map of complexes $\xymatrix{\Phi':(C')^\bullet \ar[r] & (D')^\bullet}$, which coincides with $\Phi$, $C$ and $D$ in all degrees except for the $0$th degree, where there are spaces
\begin{eqnarray*}
(C')^0=C_\phi^0 & \textnormal{and} & (D')^0=D_\phi^0
\end{eqnarray*} 
that are related to the original complexes by maps that fit into
\begin{eqnarray*}
\xymatrix{
C^0 \ar[rd]^{d_C}\ar[dd]_{\phi_C}\ar[rr]^{\Phi} &                          & D^0 \ar[rd]^{d_D}\ar'[d]_{\phi_D}[dd] &      \\
                                                & C_1 \ar[rr]^{\qquad\Phi} &                                       & D_1. \\
C_\phi^0 \ar[ur]_{d_C'}\ar[rr]_{\Phi'}          &                          & D_\phi^0\ar[ur]_{d_D'}                &
}
\end{eqnarray*}
Then, $\Phi'$ also induces an isomorphism in cohomology for degree $1$.
\end{lemma}
\begin{proof}
Let $Z^1_X:=\ker(\xymatrix{d_X:X^1 \ar[r] & X^2})$ and $B^1_X:=d_X(X^0)$, for $X\in\lbrace C,D\rbrace$ and consider the map of exact sequences
\begin{eqnarray*}
\xymatrix{
0 \ar[r] & B^1_C \ar[d]_{\Phi\rest{B^1_C}}\ar[r] & Z^1_C \ar[d]_{\Phi}\ar[r] & H^1(C) \ar[d]_{[\Phi]}\ar[r] & 0 \\
0 \ar[r] & B^1_D \ar[r]                          & Z^1_D \ar[r]              & H^1(D) \ar[r] & 0 .
}
\end{eqnarray*}
By the snake lemma, there is a long exact sequence
\begin{align*}
\xymatrix{
0 \ar[r] & \ker(\Phi\rest{B^1_C}) \ar[r] & \ker\Phi } & \xymatrix{{} \ar[r] & \ker[\Phi] \ar[r] & ...} \\
 & \xymatrix{... \ar[r] & \coker(\Phi\rest{B^1_C}) \ar[r] & \coker\Phi \ar[r] & \coker[\Phi] \ar[r] & 0 ;
}
\end{align*}
however, by hypothesis $\ker[\Phi]$ and $\coker[\Phi]$ are trivial, thus implying
\begin{eqnarray*}
\ker(\Phi\rest{B^1_C})=\ker\Phi & \textnormal{and} & \coker(\Phi\rest{B^1_C})\cong\coker\Phi .
\end{eqnarray*}
Notice that the first equation can be reinterpreted as saying $\ker\Phi\subset B^1_C$. Using the analogous sequence with $B^1_{\phi_X}=d'_X(X^0_\phi)$,
\begin{eqnarray*}
\xymatrix{
0 \ar[r] & B^1_{\phi_C}\ar[d]_{\Phi\rest{B^1_{\phi_C}}}\ar[r] & Z^1_C \ar[d]_{\Phi}\ar[r] & H^1(C') \ar[d]_{[\Phi']}\ar[r] & 0 \\
0 \ar[r] & B^1_{\phi_D}\ar[r]                                 & Z^1_D \ar[r]              & H^1(D') \ar[r] & 0 
}
\end{eqnarray*}
and the snake lemma yet again, we get
\begin{align*}
\xymatrix{
0 \ar[r] & \ker(\Phi\rest{B^1_{\phi_C}}) \ar[r] & \ker\Phi } & \xymatrix{{} \ar[r] & \ker[\Phi'] \ar[r] & ...} \\
 & \xymatrix{... \ar[r] & \coker(\Phi\rest{B^1_{\phi_C}}) \ar[r] & \coker\Phi \ar[r]  & \coker[\Phi'] \ar[r] & 0 .
}
\end{align*}
This time around, nonetheless, since for every element $x\in X^0$, $d_X(x)=d'_X(\phi_X(x))$, we have got that
\begin{eqnarray*}
B^1_C\subseteq B^1_{\phi_C} & \textnormal{and} & B^1_D\subseteq B^1_{\phi_D};
\end{eqnarray*}
therefore, from the first inclusion and the one above,
\begin{eqnarray*}
\ker(\Phi\rest{B^1_\phi})=\ker\Phi\cap B^1_\phi=\ker\Phi .
\end{eqnarray*}
Thus, we conclude that $\ker[\Phi']$ is zero, or the induced map is injective, and the short exact sequence
\begin{eqnarray*}
\xymatrix{
0 \ar[r] & \coker(\Phi\rest{B^1_{\phi_C}}) \ar[r] & \coker\Phi \ar[r]  & \coker[\Phi'] \ar[r] & 0.
}
\end{eqnarray*}
Finally, using the second inclusion and the remaining sequences, we have got a map 
\begin{eqnarray*}
\xymatrix{
0 \ar[r] & \coker(\Phi\rest{B^1_C}) \ar[d]_{\alpha}\ar[r] & \coker\Phi \ar[d]_{Id}\ar[r]  & 0 \ar[d]\ar[r]       & 0 \\
0 \ar[r] & \coker(\Phi\rest{B^1_{\phi_C}}) \ar[r]         & \coker\Phi \ar[r]             & \coker[\Phi'] \ar[r] & 0 
}
\end{eqnarray*}
where $\alpha(d_D(y)+\Phi(B^1_C)):=d'_D(\phi_D(y))+\Phi(B^1_{\phi_C})$, whose associated long exact sequence by the snake lemma tells us that $\alpha$ is an isomorphism and also $\coker[\Phi']$ is trivial, or the induced map is surjective.

\end{proof}
In regards of the latter proof, the inclusions $B^1_X\subseteq B^1_{\phi_X}$ also give rise to exact sequences
\begin{eqnarray*}
\xymatrix{
0 \ar[r] & B^1_X \ar[d]\ar[r]  & Z^1_C \ar[d]_{Id}\ar[r] & H^1(X) \ar[d]\ar[r] & 0 \\
0 \ar[r] & B^1_{\phi_X} \ar[r] & Z^1_C \ar[r]            & H^1(X') \ar[r]  & 0 
}
\end{eqnarray*}
out of whose snake lemma long exact sequence one reads that 
\begin{eqnarray*}
H^1(X)\cong H^1(X')\oplus \frac{B^1_{\phi_X}}{B^1_X}.
\end{eqnarray*}
What the proof ultimately says is that the isomorphism $H^1(C)\cong H^1(D)$ is diagonal with respect to the direct sum decompositions.
\begin{prop}\label{TweakedVanEst}
If $G$ is $k$-connected, and $\Phi_H^q$ is the restriction of $\Phi_H$ to the $q$th row,
\begin{eqnarray*}
H^n(\Phi_H^q)=(0), & \textnormal{for all degrees }n\leq k.
\end{eqnarray*}
\end{prop}
\begin{proof}
Since $G$ is the $s$-fibre of 
\begin{eqnarray*}
\xymatrix{
\G_p^q\times G \ar@<0.5ex>[r]\ar@<-0.5ex>[r] & \G_p^q
}
\end{eqnarray*} 
and it is $k$-connected by hypothesis, together with $\Phi_H^q$ being defined by the van Est map from degree $1$ upwards, we can use theorem \ref{CrainicVanEst} to deduce that
\begin{eqnarray*}
H^n(\Phi_H^q)=(0), & \textnormal{for all degrees }1<n\leq k.
\end{eqnarray*}
As for degrees $0$ and $1$, that is the contents of lemmas \ref{tweak0} and \ref{tweak1}.

\end{proof}
\begin{cor}\label{Hor p-vanEst}
If $G$ is $k$-connected,
\begin{eqnarray*}
H^n_{tot}(\Phi_H)=(0), & \textnormal{for all degrees }n\leq k.
\end{eqnarray*}
\end{cor}
\begin{proof}
To make the argument run in the same spirit of the other van Est type theorems we have got, we consider the both the domain and the co-domain of $\Phi_H$ reflected by the diagonal. Then, we compute the cohomology of the mapping cone of $\Phi_H$ using the spectral sequence of the double complex filtrated by columns, whose first page is
\begin{eqnarray*}
\xymatrix{
           & \vdots                                        & \vdots                                        & \vdots                           &                      \\ 
           & H^2_{\delta_{\Phi_H}}(C^{0,\bullet}(\Phi_H)) \ar[r] & H^2_{\delta_{\Phi_H}}(C^{1,\bullet}(\Phi_H)) \ar[r] & H^2_{\delta_{\Phi_H}}(C^{2,\bullet}(\Phi_H)) \ar[r] & \dots \\
E^{p,q}_1: & H^1_{\delta_{\Phi_H}}(C^{0,\bullet}(\Phi_H)) \ar[r] & H^1_{\delta_{\Phi_H}}(C^{1,\bullet}(\Phi_H)) \ar[r] & H^1_{\delta_{\Phi_H}}(C^{2,\bullet}(\Phi_H)) \ar[r] & \dots \\
           & H^0_{\delta_{\Phi_H}}(C^{0,\bullet}(\Phi_H)) \ar[r] & H^0_{\delta_{\Phi_H}}(C^{1,\bullet}(\Phi_H)) \ar[r] & H^0_{\delta_{\Phi_H}}(C^{2,\bullet}(\Phi_H)) \ar[r] & \dots 
}
\end{eqnarray*}
and since $\Phi_H$ is defined row-wise by van Est maps, after the reflection through the diagonal, these look like column-wise van Est maps. Thus, the columns of the mapping cone double complex coincide with the mapping cones of these and proposition \ref{TweakedVanEst} tells us that the $q$th column of the latter diagram is zero below $k$. Given that
\begin{eqnarray*}
E^{p,q}_1\Rightarrow H^{p+q}_{tot}(\Phi_H),
\end{eqnarray*}
it follows from lemma \ref{BelowDiag} that $H_{tot}^n(\Phi_H)=(0)$ for $n\leq k$ as desired.

\end{proof}
Perhaps more interestingly than for this latter corollary, we can use proposition \ref{TweakedVanEst} to prove that that differentiating yet again in the other direction one lands in the $p$-page of the Lie $2$-algebra and that this process induces isomorphisms in cohomology. \\
Instead of going through the whole process of proving that indeed there is a double complex with coefficients associated to the double Lie algebroid that comes out of differentiating twice the double Lie groups that generate the $p$-pages and their particular representations again, we proceed to show that assembling row-wise van Est maps in the vertical \LA -double complex, one casually lands in the $p$-page of the triple complex of the Lie $2$-algebra with values in a $2$-representation. Before doing so, we point out that we would get the same result by assembling column-wise van Est maps in the horizontal \LA -double complex, but we refrain to choose this approach, since this would need us to develop a van Est theorem for the sub-complexes of Lie groupoid cochains that show up in the columns. This is not the case in the approach we take, as we will see shortly. \\
Recall that the generic row of the vertical \LA -double complex consists of the Lie groupoid cochains of the Lie group bundle
\begin{eqnarray*}
\xymatrix{
\gg_p^q\times G \ar@<0.5ex>[r]\ar@<-0.5ex>[r] & \gg_p^q
}
\end{eqnarray*}
which are alternating $q$-multilinear in the $\gg_p$-coordinates with values in the trivial vector bundle $\xymatrix{\gg_p^q\times W \ar[r] & \gg_p^q}$. Because of the formula for the differential, we see that these complexes coincide with the Lie group complex of $G$ with values in the representation
\begin{eqnarray*}
\xymatrix{\rho_{(q)}:G \ar[r] & GL(\bigwedge^q\gg_p^*\otimes W)} \\
\rho_{(q)}(g)\omega =\rho_0^1(i(g))\omega\qquad
\end{eqnarray*}
tweaked by $V$ at degree $0$. Differentiating this complex, one lands in the tweaked Chevalley-Eilenberg complex of $\gg$ with values in
\begin{eqnarray*}
\xymatrix{\dot{\rho}_{(q)}:\gg \ar[r] & \ggl(\bigwedge^q\gg_p^*\otimes W)} \\
\dot{\rho}_{(q)}(x)\omega =\dot{\rho}_0^1(\mu(x))\omega\qquad
\end{eqnarray*}
which is precisely the generic row in the $p$-page of the double triple complex of the Lie $2$-algebra with values in $\dot{\rho}$. Since we already know that these complexes fit together into a double complex by placing the Chevalley-Eilenberg complexes of $\gg_p$ with values in
\begin{eqnarray*}
\xymatrix{\dot{\rho}^{(r)}:\gg_p \ar[r] & \ggl(\bigwedge^r\gg^*\otimes W)} 
\end{eqnarray*}
\begin{eqnarray*}
\dot{\rho}^{(r)}(\xi)\omega(x_1,...,x_r) =\dot{\rho}_0^1(\hat{t}_p(\xi))\omega(x_1,...,x_r)-\sum_{k=1}^r\omega(x_1,...,x_{k-1},\Lie_{\hat{t}_p\xi}x_k,x_{k+1},...,x_r)
\end{eqnarray*}
for columns, we proceed to prove that the assembly of the row-wise van Est maps is indeed a map of double complexes.
\begin{prop}
\begin{eqnarray*}
\xymatrix{
\Phi_{row}:C_{\LA}(\gg_p\ltimes G,\phi) \ar[r] & C^{p,\bullet}_\bullet(\gg_1 ,\phi)
}
\end{eqnarray*}
is a map of double complexes.
\end{prop}
\begin{proof}
By definition $\Phi_{row}$ intertwine horizontal differentials. Moreover, since restricted to the first column $\Phi_{row}$ is defined to be the identity, we just need to show that the square
\begin{eqnarray*}
\xymatrix{
C(G^{r},\bigwedge^{q+1}\gg_p^*\otimes W) \ar[r]^{\Phi_{row}} & \bigwedge^{q+1}\gg_p^*\otimes\bigwedge^{r}\gg^*\otimes W \\
C(G^{r},\bigwedge^{q}\gg_p^*\otimes W) \ar[u]^\delta\ar[r]_{\Phi_{row}} & \bigwedge^{q}\gg_p^*\otimes\bigwedge^{r}\gg^*\otimes W \ar[u]_\delta
}
\end{eqnarray*}
commutes. Indeed, let $\omega\in C(G^{r},\bigwedge^{q}\gg_p^*\otimes W)$ and $(\Xi;X)=(\xi_0,...,\xi_q;x_1,...,x_r)\in\gg_p^q\oplus\gg^r$, then
\begin{align*}
\overrightarrow{R}_X\delta\omega(\Xi) & =\frac{d^J}{d\tau_J}\rest{\tau=0}\delta\omega(\Xi;\exp(\tau\cdot X)) \\
									  & =\frac{d^J}{d\tau_J}\rest{\tau=0}\sum_{m<n}\omega([\xi_m,\xi_n],\Xi(m,n);\exp(\tau\cdot X))+\sum_{j=0}^{q}(-1)^j\rho_{\gg_p}^{r}(\xi_j)\omega(\Xi(j);\exp(\tau\cdot X)) \\
									  & =\sum_{m<n}\overrightarrow{R}_X\omega([\xi_m,\xi_n],\Xi(m,n))+\sum_{j=0}^{q}(-1)^j\frac{d^J}{d\tau_J}\rest{\tau=0}\rho_{\gg_p}^{r}(\xi_j)\omega(\Xi(j);\exp(\tau\cdot X)) .
\end{align*}
Now, let $\lbrace e_a\rbrace$ be a basis for $W$, and 
\begin{align*}
\omega(\Xi(j);\exp(\tau\cdot X))=\omega^a e_a ,
\end{align*}
where $\omega^a=\omega^a(\Xi(j);\exp(\tau\cdot X))$ and we used the Einstein summation convention. Since by definition
\begin{eqnarray*}
\rho_{\gg_p}^r(\xi_j)\omega(\Xi(j);\exp(\tau\cdot X))=\omega^a\dot{\rho}_0^1(\hat{t}_p(\xi_j))e_a+(\frac{d}{d\lambda}\rest{\lambda=0}\omega^a(\Xi(j);\exp(\tau\cdot X)^{t_p(\exp(\lambda\xi_j))}))e_a ,
\end{eqnarray*}
we compute
\begin{align*}
\frac{d^J}{d\tau_J}\rest{\tau=0}\omega^a\dot{\rho}_0^1(\hat{t}_p(\xi_j))e_a & =\dot{\rho}_0^1(\hat{t}_p(\xi_j))\overrightarrow{R}_X\omega(\Xi(j)),
\end{align*}
and, using that $\exp(x)^h=\exp(x^h)$ for all $x\in\gg$ and $h\in H$,
\begin{align*}
\frac{d}{d\tau_r} & \rest{\tau_r=0}\frac{d}{d\lambda}\rest{\lambda=0}\omega^a(\Xi(j);\exp(\tau\cdot X)^{t_p(\exp(\lambda\xi_j))}) \\
 & =\frac{d}{d\lambda}\rest{\lambda=0}R_{x_r^{t_p(\exp(\lambda\xi_j))}}\omega^a(\Xi(j);\exp(\tau\cdot X(r))^{t_p(\exp(\lambda\xi_j))}) \\
 & =R_{\Lie_{\hat{t}_p(\xi_j)}x_r}\omega^a(\Xi(j);\exp(\tau\cdot X(r)))+\frac{d}{d\lambda}\rest{\lambda=0}R_{x_r}\omega^a(\Xi(j);\exp(\tau\cdot X(r))^{t_p(\exp(\lambda\xi_j))}).
\end{align*}
Notice that, if we are to continue inductively, the first term of this latter equation will immediately yield
\begin{eqnarray*}
R_{x_1}...R_{x_{r-1}}R_{\Lie_{\hat{t}_p(\xi_j)}x_r}\omega(\Xi(j));
\end{eqnarray*}
whereas, carefully using the induction hypothesis, we get 
\begin{align*}
\frac{d}{d\tau_{r-1}} & \rest{\tau_{r-1}=0}\frac{d}{d\lambda}\rest{\lambda=0}R_{x_r}\omega^a(\Xi(j);\exp(\tau\cdot X(r))^{t_p(\exp(\lambda\xi_j))}) \\
 & =R_{\Lie_{\hat{t}_p(\xi_j)}x_{r-1}}R_{x_r}\omega^a(\Xi(j);\exp(\tau\cdot X(r-1,r)))+ \\
 & \qquad +\frac{d}{d\lambda}\rest{\lambda=0}R_{x_{r-1}}R_{x_r}\omega^a(\Xi(j);\exp(\tau\cdot X(r-1,r))^{t_p(\exp(\lambda\xi_j))})
\end{align*}
and eventually, 
\begin{eqnarray*}
\sum_{k=1}^rR_{x_1}...R_{x_{k-1}}R_{\Lie_{\hat{t}_p\xi}x_k}R_{x_{k+1}}...R_{x_r}\omega(\Xi(j)).
\end{eqnarray*}
Thus, after taking the alternating sum over $S_r$,
\begin{align*}
\Phi_{row}\delta\omega(\Xi;X) & =\sum_{m<n}\Phi_{row}\omega([\xi_m,\xi_n],\Xi(m,n);X)+\sum_{j=0}^{q}(-1)^j\Big{[}\dot{\rho}_0^1(\hat{t}_p(\xi_j))\Phi_{row}\omega(\Xi(j);X)+ \\
                              & \qquad +\sum_{k=1}^r\Phi_{row}\omega(\Xi(j);x_1,...,x_{k-1},\Lie_{\hat{t}_p\xi}x_k,x_{k+1},...,x_r)\Big{]} \\
                              & =\sum_{m<n}\Phi_{row}\omega([\xi_m,\xi_n],\Xi(m,n);X)+\sum_{j=0}^{q}(-1)^j\dot{\rho}^{(r)}(\xi_j)\Phi_{row}\omega(\Xi(j);X) \\
                              & =\delta\Phi_{row}\omega(\Xi;X)
\end{align*}
as claimed.

\end{proof}
As announced, we now prove that there is an isomorphism between the total cohomologies of the vertical \LA -complex and the $p$-pages of the Lie $2$-algebra.
\begin{prop}\label{res-vanEst} % res stands for residual.
If $G$ is $k$-connected and $\Phi_{row}$ is the map of the previous proposition, then
\begin{eqnarray*}
H^n_{tot}(\Phi_{row})=(0), & \textnormal{for all degrees }n\leq k.
\end{eqnarray*}
\end{prop}
\begin{proof}
Once again, we opt to consider both the domain and the co-domain of $\Phi_{row}$ reflected by the diagonal, so that the argument runs along the same lines of the other van Est type theorems we have got. Computing the total cohomology of the mapping cone of $\Phi_{row}$ using the spectral sequence of the double complex filtrated by columns yields the familiar first page
\begin{eqnarray*}
\xymatrix{
           & \vdots                                        & \vdots                                        & \vdots                           &                      \\ 
           & H^2_{\delta_{\Phi_{row}}}(C^{0,\bullet}(\Phi_{row})) \ar[r] & H^2_{\delta_{\Phi_{row}}}(C^{1,\bullet}(\Phi_{row})) \ar[r] & H^2_{\delta_{\Phi_{row}}}(C^{2,\bullet}(\Phi_{row})) \ar[r] & \dots \\
E^{p,q}_1: & H^1_{\delta_{\Phi_{row}}}(C^{0,\bullet}(\Phi_{row})) \ar[r] & H^1_{\delta_{\Phi_{row}}}(C^{1,\bullet}(\Phi_{row})) \ar[r] & H^1_{\delta_{\Phi_{row}}}(C^{2,\bullet}(\Phi_{row})) \ar[r] & \dots \\
           & H^0_{\delta_{\Phi_{row}}}(C^{0,\bullet}(\Phi_{row})) \ar[r] & H^0_{\delta_{\Phi_{row}}}(C^{1,\bullet}(\Phi_{row})) \ar[r] & H^0_{\delta_{\Phi_{row}}}(C^{2,\bullet}(\Phi_{row})) \ar[r] & \dots 
}
\end{eqnarray*}
Since $\Phi_{row}$ is defined row-wise by van Est maps, after tilting, these can be thought of as column-wise van Est maps. Thus, the columns of the mapping cone double complex coincide with the mapping cones of these and proposition \ref{TweakedVanEst} applied to the Lie $2$-group $\xymatrix{G \ar[r] & 1}$ with values in the corresponding representations, tells us that the $q$th column of the latter diagram is zero below $k$. Given that
\begin{eqnarray*}
E^{p,q}_1\Rightarrow H^{p+q}_{tot}(\Phi_{row}),
\end{eqnarray*}
it follows from lemma \ref{BelowDiag} that $H_{tot}^n(\Phi_{row})=(0)$ for $n\leq k$ as desired.

\end{proof}
We use this latter result to prove the only missing piece of the proof of the main theorem.
\begin{prop}
If $H$ and $G$ are both $k$-connected, the restriction of the $2$-van Est map to the $p$th page induces isomorphisms
\begin{eqnarray*}
\xymatrix{
\Phi(p)^n :H_{tot}^n(C^{p,\bullet}_\bullet(\G ,\phi)) \ar[r] & H_{tot}^n(C^{p,\bullet}_\bullet(\gg_1 ,\phi)),
}
\end{eqnarray*}
for $n\leq k$, and it is injective for $n=k+1$.
\end{prop}
\begin{proof}
We start by pointing out that from the $k$-connectedness of $H$ and $G$, together with a K\"unneth type argument as the one we used in the proof of theorem \ref{2vE-vs}, we can safely conclude that $\G_p$ is also $k$-connected. We can thus invoke lemma \ref{Ver p-vanEst} and the rephrasing of proposition \ref{mapConeCoh} to conclude that 
\begin{eqnarray*}
\xymatrix{
\Phi_V^n :H_{tot}^n( C^{p,\bullet}_\bullet(\G ,\phi)) \ar[r] & H_{tot}^n( C_\LA(\gg_p\ltimes G,\phi))
}
\end{eqnarray*}
is an isomorphism for $n\leq k$, and it is injective for $n=k+1$. Analogously, using that $G$ is $k$-connected, we get from proposition \ref{res-vanEst} and the rephrasing of proposition \ref{mapConeCoh} that $\Phi_{row}$ also induces an isomorphism and an injective map in the corresponding degrees. Since from the very formula, we see that the restriction of the $2$-van Est map to each $p$-page factors through the vertical \LA -double complex
\begin{eqnarray*}
\xymatrix{
C(\G_p^q\times G^r,W) \ar[rr]^\Phi\ar[dr]_{\Phi_V} & & \bigwedge^q\gg_p^*\otimes\bigwedge^r\gg^*\otimes W , \\
                       & C(G^r,\bigwedge^q\gg_p^*\otimes W) \ar[ur]_{\Phi_{row}} &
}
\end{eqnarray*}
we get
\begin{eqnarray*}
\xymatrix{
H_{tot}^n(C^{p,\bullet}_\bullet(\G ,\phi)) \ar[rr]^{\Phi(p)^n}\ar[dr]_{\Phi_V^n} & & H_{tot}^n(C^{p,\bullet}_\bullet(\gg_1 ,\phi)) \\
                       & H_{tot}^n( C_\LA(\gg_p\ltimes G,\phi)) \ar[ur]_{\Phi_{row}^n} &
}
\end{eqnarray*}
and $\Phi(p)^n$ is a composition of isomorphisms for $n\leq k$, and a composition of injective maps for $n=k+1$.

\end{proof}
For convenience, we recast the latter result using proposition \ref{mapConeCoh}.
\begin{prop}\label{vanEst p-pages}
Under the hypothesis of the previous proposition, if $\Phi(p)$ is the restriction of the $2$-van Est map to the $p$th page, then
\begin{eqnarray*}
H^n_{tot}(\Phi(p))=(0), & \textnormal{for all degrees }n\leq k.
\end{eqnarray*}
\end{prop}

Thus, we finally conclude:
\begin{proof} (\textit{of the main Theorem})
Just as it was outlined above, let $\bar{\Phi}$ be the map of honest double complexes induced by $\Phi$. Computing the cohomology of the mapping cone of $\bar{\Phi}$ using the spectral sequence of the double complex filtrated by columns, we arrive at the first page 
\begin{eqnarray*}
\xymatrix{
           & \vdots                                        & \vdots                                        & \vdots                           &                      \\ 
           & H^2_{\delta_{\bar{\Phi}}}(C^{0,\bullet}(\bar{\Phi})) \ar[r] & H^2_{\delta_{\bar{\Phi}}}(C^{1,\bullet}(\bar{\Phi})) \ar[r] & H^2_{\delta_{\bar{\Phi}}}(C^{2,\bullet}(\bar{\Phi})) \ar[r] & \dots \\
E^{p,q}_1: & H^1_{\delta_{\bar{\Phi}}}(C^{0,\bullet}(\bar{\Phi})) \ar[r] & H^1_{\delta_{\bar{\Phi}}}(C^{1,\bullet}(\bar{\Phi})) \ar[r] & H^1_{\delta_{\bar{\Phi}}}(C^{2,\bullet}(\bar{\Phi})) \ar[r] & \dots \\
           & H^0_{\delta_{\bar{\Phi}}}(C^{0,\bullet}(\bar{\Phi})) \ar[r] & H^0_{\delta_{\bar{\Phi}}}(C^{1,\bullet}(\bar{\Phi})) \ar[r] & H^0_{\delta_{\bar{\Phi}}}(C^{2,\bullet}(\bar{\Phi})) \ar[r] & \dots 
}
\end{eqnarray*}
Since $\bar{\Phi}$ is defined for the $p$th column as the restriction of the $2$-van Est map to the corresponding $p$-page, the $p$th column of the mapping cone double complex is the total complex of the mapping cone of $\Phi(p)$. By means of proposition \ref{vanEst p-pages}, we know that the $p$th column of the latter diagram is zero below $k$, as both $H$ and $G$ are $k$-connected by hypothesis. Given that
\begin{eqnarray*}
E^{p,q}_1\Rightarrow H^{p+q}_{tot}(\bar{\Phi}),
\end{eqnarray*}
it follows from lemma \ref{BelowDiag} that $H_{tot}^n(\bar{\Phi})=(0)$ for $n\leq k$. Using proposition \ref{mapConeCoh} one last time,
the $2$-van Est map induces isomorphisms
\begin{eqnarray*}
\xymatrix{
\Phi^n :H_\nabla^n(\G ,\phi) \ar[r] & H_\nabla^n(\gg_1 ,\phi),
}
\end{eqnarray*}
for $n\leq k$, and it is injective for $n=k+1$ as desired.

\end{proof}
%*** is there need for a word about the relation between $H_\nabla$ and $H_{tot}(\bar{-})$? *** \\

\section[Integration]{Integration}\label{sec:bary}
% \ref{th:main} sup. this summons the main theorem.
\sectionmark{Integration}
%---------------------------------------------------------------------------------------------------------------------------------------------------
Having established the equivalence between the categories of crossed modules and Lie $2$-groups and Lie $2$-algebras, the problem of integrating a Lie $2$-algebras boils down to integrating the structure map defining its associated crossed module and the Lie algebra action by derivations to a Lie group action by automorphisms. The first task is easily performed due to the classical Lie theory, whereas the latter follows from the path strategy to integration (see appendix \ref{appchapter}).
\begin{theorem}\cite{ZhuInt2Alg}
Let $\xymatrix{\phi :\hh \ar[r] & \mathfrak{Der}(\gg )}$ be a Lie algebra action by derivations. Let $\xymatrix{\varphi :H \ar[r] & Aut(\gg )}$ be the unique group morphisms integrating $\phi$.
\begin{itemize}
\item Let $x\in\gg$ and $h\in H$. Then $\varphi_h(x)=\xi(1)$, where $\xi\in P(\gg)$ is the solution of the ODE:
\begin{eqnarray*}
\frac{d}{d\lambda}\xi(\lambda)=\phi_{\zeta(\lambda)}\xi(\lambda)
\end{eqnarray*}
with initial value $\xi(0)=x$. Here $\zeta\in P(\hh)$ is any representative of the $\hh$-homotopy class $h$.
\item Let $\xi\in P(\gg)$ and $h\in H$. Then $\varphi_h\xi(\lambda)=\varpi(1,\lambda)$, where $\varpi(-,\lambda)\in P(\gg)$ is the solution of the ODE:
\begin{eqnarray*}
\frac{\partial}{\partial\lambda_0}\varpi(\lambda_0,\lambda_1)=\phi_{\zeta(\lambda_0)}\varpi(\lambda_0 ,\lambda_1)
\end{eqnarray*}
with initial value $\varpi(0,\lambda)=\xi(\lambda)$. Here $\zeta\in P(\hh)$ is taken as in the previous item.
\end{itemize} 
Consequently, the corresponding group action $\xymatrix{\Phi :H \ar[r] & Aut(G)}$ is given by
\begin{eqnarray*}
\Phi _h([\xi])=[\varphi_h(\xi)]=[\varpi(1,-)].
\end{eqnarray*}
\end{theorem}
We argue that the fact that such a proof exists precludes the positivity of the integrability result that follows. In it, we will use the van Est type strategy to get a new proof of the integrability of Lie $2$-algebras.
%---------------------------------------------------------------------------------------------------------------------------------------------------
\subsection{The van Est strategy}\label{sec:VS}
%---------------------------------------------------------------------------------------------------------------------------------------------------
In \cite{VanEst}, the author used a cohomological approach to the integration of Lie algebras. Let us outline what we will call the van Est strategy. This consists of two steps:
First, suppose the adjoint representation is faithful. Then, we can regard the Lie algebra $\gg$ as a sub-algebra of $\ggl (\gg )$. As a matter of fact, all sub-algebras of the linear Lie algebra are integrable, thus yielding the result.\\
Second, suppose otherwise, that the adjoint representation is not faithful, then there is an exact sequence 
\begin{eqnarray*}
\xymatrix{
0 \ar[r] & \ker (ad) \ar[r] & \gg \ar[r] & ad(\gg ) \ar[r] & 0.
}
\end{eqnarray*}
By the usual theory of Lie algebra extensions, this exact sequence corresponds to a class 
\begin{eqnarray*}
\omega \in H^2(ad(\gg),\ker (ad)).
\end{eqnarray*}
The fact that $ad(\gg)$ is integrable, being a sub-algebra of the linear Lie algebra, still holds. Then, let $G$ be the $1$-connected integration of $ad(\gg)$. Now, there is van Est theorem relating the differential cohomology of a Lie group to the Chevalley-Eilenberg cohomology of its Lie algebra.
\begin{theorem}
Let $H$ be a Lie group with Lie algebra $\hh$. If $H$ is $k$-connected, the map
\begin{eqnarray*}
\xymatrix{
\Phi: H^n(H,V) \ar[r] & H^n(\hh ,V)
}
\end{eqnarray*}
that differentiates group cochains to get Lie algebra cochains is an isomorphism for all degrees $n\leq k$, and it is injective for $n=k+1$.
\end{theorem}
Explicitly, the map $\Phi$ is defined by the chain map
\begin{eqnarray*}
\Phi(f)(x_1,...,x_p):=\sum _{\sigma\in S_p}\abs{\sigma}R_{x_{\sigma (p)}}...R_{x_{\sigma (1)}}f,
\end{eqnarray*}
for $f\in C^p(G)$ and $x_1,...,x_p\in\gg$. Here, $\xymatrix{R_x:C^p (G) \ar[r] & C^{p-1}(G)}$ is given by $R_x(f)(g_1,...,g_{p-1}):=\vec{x}f(-,g_1,...,g_{p-1})(e)$, where $\vec{x}$ stands for the right-invariant vector field associated to $x$ and $e$ is the group identity. Both this map and the statement the latter theorem are clearly special cases of the ones in theorem \ref{CrainicVanEst}. \\
Now, since $G$ is simply connected and a Lie group, $\pi_2(G)=(0)$; hence, $G$ is at least $2$-connected. Please note that this fact is far from trivial; famous proofs use Morse theory among other tools from algebraic topology. Now, given that $G$ is $2$-connected, van Est theorem says that the class $\omega$ comes from a class $\small{\int}\omega\in H^2(G,\ker (ad))$. This class, in turn by the theory of Lie group extensions, corresponds to a sequence
\begin{eqnarray*}
\xymatrix{
1 \ar[r] & \ker (ad) \ar[r] & \G \ar[r] & G \ar[r] & 1.
}
\end{eqnarray*}
Consequently, $\G$ is a Lie group integrating $\gg$, as we wanted. Notice that because $\G$ fits as an extension of two simply connected manifolds, it is simply connected itself, so we can say that $\G$ is the $1$-connected integration of $\gg$. \\
The previous strategy was successfully applied by the author of \cite{VanEstC} to prove theorem \ref{MariusInt}. Indeed, the second cohomology of Lie groupoid and Lie algebroid cohomology also classifies extensions, and the rest of the ingredients are explicitly asked as hypothesis. Notice, though, that the main differences between the statement of the referred theorem and the proof of Lie III given by the van Est theorem are that, on the one hand, groupoids lack an adjoint representation to start with and, on the other, there is no argument to assert that a Lie groupoid (or more precisely its source fibres) is $2$-connected. One could say that precisely because of this, the van Est strategy failed to solve the integrability of Lie algebroids. Due to this failure, we find it purposeful to index the list of ingredients for the argument to run.
\subsubsection{The necessary ingredients}\label{vanEstIngr}
 Roughly, one will be needing: \\
(1) The existence of an adjoint representation, \\
(2) the integrability of sub-algebras of the linear Lie algebra, \\
(3) a cohomology theory whose second group classifies extensions, \\
(4) $\pi_2(G)=(0)$ and finally, \\
(5) a van Est theorem. \\
Of course, alternatively, any representation whose integration one can control would suffice to replace ingredients (1) and (2). Also, the cohomology theory of step (3) is actually two cohomology theories, one for the infinitesimal and one for the global counter-parts of the given geometric structure. In this, the van Est theorem of (5) is a theorem relating these cohomology theories and giving isomorphisms under the appropriate connectedness hypothesis that (4) represents. \\
One might ask whether it is possible to use this strategy to integrate \LA -groupoids directly, but pretty much as in the groupoid case, there is no natural representation and, worse yet, there is no working definition of a representation of an \LA -groupoid. \\
Now, in the body of this work we got all of the ingredients but the second one. Nevertheless, the image of any linear functor between $2$-vector spaces yields a Lie subgroupoid and in \cite{IntSubLin2}, the authors prove that Lie $2$-subalgebras of $\ggl(\phi)$ are always integrable; thus, we are ready to apply the van Est strategy.

%-------------- Should do it, now shouldn't I? At least go through the proof in the reference.

\subsection{Integrating Lie 2-algebras}
\subsubsection{Vector space coefficients}
As a corollary of Theorem \ref{2vE-vs} and its corollary, we have got the following integrability result.
\begin{theorem}
If $\xymatrix{\gg_1 \ar@<-0.5ex>[r] \ar@<0.5ex>[r] & \hh}$ is a Lie $2$-algebra such that 
\begin{eqnarray*}
\gg\cap\mathfrak{c}(\hat{u}(\hh))=(0),
\end{eqnarray*}
where $\gg$ is the core and $\mathfrak{c}(\hat{u}(\hh))$ is the centralizer of $\hat{u}(\hh)$ in $\gg_1$, then it is integrable. 
\end{theorem}
\begin{proof}
We are to use the van Est strategy. So, consider the adjoint $2$-representation of the given Lie $2$-algebra and the natural sequence
\begin{eqnarray}\label{vE2-Ext}
\xymatrix{
0 \ar[r] & \ker(ad_1) \ar[d]\ar@{^{(}->}[r] & \gg \ar[d]_\mu\ar[r]^{ad_1\quad} & ad_1(\gg) \ar[d]\ar[r] & 0  \\
0 \ar[r] & \ker(ad_0) \ar@{^{(}->}[r]       & \hh \ar[r]_{ad_0\quad}           & ad_0(\hh) \ar[r]       & 0, 
}
\end{eqnarray}
If the hypothesis holds, for each $x\in\gg$, $x\not\in\mathfrak{c}(\hat{u}(\hh))$. Now, by definition
\begin{eqnarray*}
\mathfrak{c}(\hat{u}(\hh)):=\lbrace \zeta\in\gg : [\zeta ,\hat{u}(y)]_1=0,\quad\forall y\in\hh\rbrace ;
\end{eqnarray*}
therefore, there exists a $y_x\in\hh$ such that $[x,\hat{u}(y_x)]_1=0$. However,
\begin{align*}
[x,\hat{u}(y_x)]_1 & =-[\hat{u}(y_x),x]_1 \\
				   & =-\Lie_{y_x}x=ad_1(x)(y_x)\neq 0.
\end{align*}
Hence, $ad_1(x)\neq 0\in Hom(\hh,\gg)$ for every $x\in\gg$, and this is $\ker(ad_1)=(0)$. As a consequence of this $\ker(ad)$ is an honest vector space, and the $2$-extension \ref{vE2-Ext} is classified by an $[\omega]\in H_{\nabla}^2(\gg_1,\ker(ad_0))$.\\
As we referenced, $\xymatrix{ad_1(\gg) \ar[r] & ad_0(\hh)}\leq\xymatrix{\ggl(\mu)_1 \ar[r] & \ggl(\mu)_0}$ is integrable, say by the Lie $2$-group $\G=\xymatrix{G \ar[r] & H}$ where both $G$ and $H$ are $1$-connected. Since every Lie group has trivial second cohomology, $H$ is $2$-connected and we can use theorem \ref{2vE-vs} to deduce that 
\begin{eqnarray*}
[\omega]=\Phi^2[\small{\int}\omega],\quad\textnormal{for some }[\small{\int}\omega]\in H_{\nabla}^2(\G ,\ker(ad_0)).
\end{eqnarray*}
Associated to this class, there is a $2$-extension of $\G$ by $\ker(ad_0)$,
\begin{eqnarray*}
\xymatrix{
1 \ar[r] & 0    \ar[d]\ar[r] & E_1 \ar[d]\ar[r] & G \ar[d]\ar[r] & 1  \\
1 \ar[r] & \ker(ad_0) \ar[r] & E_0 \ar[r]       & H \ar[r]       & 1, 
}
\end{eqnarray*}
and the Lie $2$-algebra of $\xymatrix{E_1\rtimes E_0 \ar@<-0.5ex>[r] \ar@<0.5ex>[r] & E_0}$ is $\xymatrix{\gg_1 \ar@<-0.5ex>[r] \ar@<0.5ex>[r] & \hh}$, thus proving the theorem.

\end{proof}
There is yet another way to cast this result in terms of the isotropy Lie algebras of the action of $\hh$ on $\gg$ which might be of interest; indeed, the hypothesis of the theorem is clearly equivalent to
\begin{eqnarray*}
\dim\hh^\Lie_x >0,\quad\textnormal{for all }x\in\gg ,
\end{eqnarray*}
though we picked the statement from the theorem as it refers to the Lie $2$-algebra structure, rather than the crossed module one. For the sake of clarity, the isotropy Lie algebra at $x\in\gg$ for the action $\Lie$ is given by 
\begin{eqnarray*}
\hh^\Lie_x :=\lbrace y\in\hh :\Lie_y x=0\rbrace .
\end{eqnarray*}

\subsubsection{$2$-vector space coefficients}
As a corollary of Theorem \ref{2-vanEstTheo}, we recover the integrability of Lie $2$-algebras.
\begin{theorem}
Every Lie $2$-algebra is integrable
\end{theorem}
\begin{proof}
We are to use the van Est strategy. So, consider the adjoint $2$-representation of the given Lie $2$-algebra and the natural sequence
\begin{eqnarray}
\xymatrix{
0 \ar[r] & \ker(ad_1) \ar[d]\ar@{^{(}->}[r] & \gg \ar[d]_\mu\ar[r]^{ad_1\quad} & ad_1(\gg) \ar[d]\ar[r] & 0  \\
0 \ar[r] & \ker(ad_0) \ar@{^{(}->}[r]       & \hh \ar[r]_{ad_0\quad}           & ad_0(\hh) \ar[r]       & 0. 
}
\end{eqnarray}
This $2$-extension is classified by an $[\omega]\in H_{\nabla}^2(\gg_1,\ker(ad))$. Let $\G=\xymatrix{G \ar[r] & H}$ be the integration of $\xymatrix{ad_1(\gg) \ar[r] & ad_0(\hh)}\leq\xymatrix{\ggl(\mu)_1 \ar[r] & \ggl(\mu)_0}$, where both $G$ and $H$ are $1$-connected. Since every Lie group has trivial second cohomology, both $H$ and $G$ are $2$-connected and we can use theorem \ref{2-vanEstTheo} to deduce that 
\begin{eqnarray*}
[\omega]=\Phi^2[\small{\int}\omega],\quad\textnormal{for some }[\small{\int}\omega]\in H_{\nabla}^2(\G ,\ker(ad)).
\end{eqnarray*}
Associated to this class, there is a $2$-extension of $\G$ by $\ker(ad)$,
\begin{eqnarray*}
\xymatrix{
1 \ar[r] & \ker(ad_1) \ar[d]\ar[r] & E_1 \ar[d]\ar[r] & G \ar[d]\ar[r] & 1  \\
1 \ar[r] & \ker(ad_0) \ar[r] & E_0 \ar[r]       & H \ar[r]       & 1, 
}
\end{eqnarray*}
and the Lie $2$-algebra of $\xymatrix{E_1\rtimes E_0 \ar@<-0.5ex>[r] \ar@<0.5ex>[r] & E_0}$ is $\xymatrix{\gg_1 \ar@<-0.5ex>[r] \ar@<0.5ex>[r] & \hh}$, thus proving the theorem.

\end{proof}

%% ------------------------------------------------------------------------- %%
\chapter{Appendices}\label{appchapter}
\chaptermark{Appendices}

%------------------------------------------------------
\section{Appendix A: the path strategy}
Historically, the only successful strategy for handling the integration of Lie algebroids is what we will refer to as the path strategy \cite{CF2}. In it, certain topological correspondence between spaces of paths on the Lie groupoid and its Lie algebroid is exploited to find out what is needed for the latter to inherit a smooth structure. We proceed to outline this strategy in more detail. \\
Starting out with a finite dimensional connected Lie group $G$, we consider the space of paths starting at the identity element $e\in G$
\begin{eqnarray*}
P(G)=\lbrace\gamma\in C^2(I,G):\gamma(0)=e\rbrace .
\end{eqnarray*}
For each element of $P(G)$, there is a corresponding path on the Lie algebra $\gg$ of $G$, given by
\begin{eqnarray}
D\gamma(\lambda):=\frac{d}{d\tau}\gamma(\tau)\gamma(\lambda)^{-1}\rest{\tau=\lambda}.
\end{eqnarray}
In plain words, this map takes the tangent vector to $\gamma$ at a given time and (right) translates it to the identity. Having started with a $C^2$ path, the resulting curve is $C^1$. Let $P(\gg)$ be the space of all $C^1$ paths on $\gg$. Remarkably, endowing $P(G)$ and $P(\gg)$ with the $C^2$ and $C^1$-topologies, $D$ turns out to be a homeomorphism. Since one builds the $1$-connected cover $\tilde{G}$ of $G$ quoting $P(G)$ by the equivalence classes defined by homotopies with fixed end-points, one can transport the equivalence relation to $P(\gg)$. Furthermore, one can also transport the group structure in $\tilde{G}$ which is defined in terms of concatenation of paths. In doing so, one is effectively building the $1$-connected Lie group integrating the Lie algebra $\gg$ without need of ever using $G$ or its structure. \\
The setting for the case of Lie groupoids and Lie algebroids does not differ much from this. One starts by defining the space of $G$-paths for an $s$-connected Lie groupoid $\xymatrix{G \ar@<-0.5ex>[r]\ar@<0.5ex>[r] & M}$ as the space of twice differentiable paths starting at a unit an lying on $s$-fibres. In symbols,
\begin{eqnarray*}
P(G)=\lbrace\gamma\in C^2(I,G):\gamma(0)=u(x)\textnormal{ for }x\in M, \gamma(\lambda)\in s^{-1}(x),\quad\forall\lambda\in I\rbrace .
\end{eqnarray*}
The thorough understanding of the Lie theory of Lie groupoids, gives an analogous correspondence to that of the Lie group case in that each $s$-connected Lie groupoid with Lie algebroid $A_G$ arises as a source-wise quotient of a unique (up to isomorphims) $s1$-connected Lie groupoid which can be described as
\begin{eqnarray}
\tilde{G}=P(G)/\sim ,
\end{eqnarray}
where this time around $\sim$ is the equivalence relation defined by source-wise homotopies. Again, the whole structure of the $\tilde{G}$ will be defined in terms of concatenations and evaluations. To transport this structure, there is a corresponding derivation that returns paths on the Lie algebroid of $G$ which is actually given by the same formula, though having a slightly different meaning, as the right multiplication is defined only between source fibres. One can see that $G$-paths correspond to Lie groupoid maps form the pair groupoid of $I$ to $G$ and accordingly, the space of $A$-paths is thus defined as the space of Lie algebroid maps from $TI$ to $A$. Having a Lie II theorem for morphisms of Lie algebroids, one easily sees that $D$ yields a homeomorphism again. Transporting homotopies is a harder task; however, once that is done, one can again spell out the full structure of the $s1$-connected groupoid of $A$-paths modulo $A$-homotopies of a given Lie algebroid $A$ in terms of algebroid data alone. Moreover, this groupoid is in fact a topological groupoid. The theory of \cite{CF2} devotes to answering the question of when it is possible to endow this groupoid with a smooth structure, culminating in the following theorem. % Should say more about where the monodromy groups come from? I already mentioned it before, in the remark right after the isotropies.
\begin{theorem}\cite{CF2}
A Lie algebroid $\xymatrix{A \ar[r] & M}$ is integrable if, and only if for each $x\in M$, the monodromy groups 
\begin{eqnarray*}
\mathcal{N}_x(A):=\lbrace v\in A_x : \rho_A(v)=0;\quad v\simeq_A \ast\rbrace 
\end{eqnarray*}
are locally uniformly discrete.
\end{theorem}
The precise meaning of the condition in the theorem is that for each $x\in M$, there must exist an open set $U\in A$ around $0\in A_x$ such that for all $y$'s close enough to $x$, $\mathcal{N}_x(A)\cap U=\lbrace 0\rbrace$. \\

% -------------------------
% -------------------------
% Appendix B
% -------------------------
% -------------------------

\section{Appendix B: cores}
The core is an ubiquitous structure in the theory of double objects, we dedicate this appendix to recall them. These were introduced by Pradines in his study of double vector bundles. Most of what follows is covered in \cite{2ndOrdGeom}. Let $\Omega$ be an \LA -groupoid, since $\hat{s}$ is a Lie algebroid map and due to the double source map hypothesis in the definition, it is fibrewise linear and its kernel is a vector bundle. 
\begin{eqnarray}\label{CoreSeq}
\xymatrix{
(0) \ar[r] & \ker(\hat{s}) \ar[d]\ar[r] & \Omega \ar[d]\ar[r]^{(\pi ,\hat{s})} & s^*A \ar[d]\ar[r] & (0) \\
           & H \ar@{=}[r]               & H \ar@{=}[r]                         & H                 &
}
\end{eqnarray} 
\begin{Def}
The \textit{core} of $\Omega$ is $\CC{\Omega}:=u^*\ker(\hat{s})$. In other words, $\CC{\Omega}$ is the restriction of the kernel of $\hat{s}$ to the space of units of $H$.
\end{Def}
We would like to point out that there is a natural isomorphism between the fibre of $\CC{\Omega}$ over a given element $y\in M$ and the fibres of $\ker(\hat{s})$ over each element in $h\in t^{-1}(y)\subset H$. Indeed, the isomorphism is given by
\begin{eqnarray*}
\xymatrix{
R_h:\CC{\Omega}_y \ar[r] & \ker(\hat{s})_h:c \ar@{|->}[r] & \hat{m}(c;0_h),
}
\end{eqnarray*}
where by $0_h$, we mean the zero of the fibre of $\Omega$ over $h$, and the multiplication is well defined, as $\hat{s}(c)=0_y=0_{t(h)}=\hat{t}(0_h)$. Further, the map lands where it should, because $\hat{s}(\hat{m}(c;0_h))=\hat{s}(0_h)=0_{s(h)}$ and the $\hat{m}$ covers the multiplication of $H$, and therefore $\hat{m}(c;0_h)\in\Omega_{u(y)h}=\Omega_h$. The inverse map is given by
\begin{eqnarray*}
\xymatrix{
R_h^{-1}:\ker(\hat{s})_h \ar[r] & \CC{\Omega}_y:\zeta \ar@{|->}[r] & \hat{m}(\zeta ;0_{h^{-1}}),
}
\end{eqnarray*}
where similar considerations are in order: First, $\hat{s}(\zeta)=0_{s(h)}=0_{t(h^{-1})}=\hat{t}(0_{h^{-1}})$, then $\hat{s}(\hat{m}(\zeta ;0_{h^{-1}}))=\hat{s}(0_{h^{-1}})=0_{s(h^{-1})}=0_{t(h)}=0_y$, and finally, $\hat{m}(\zeta ;0_{h^{-1}}))\in\Omega_{hh^{-1}}=\Omega_{u(y)}$. Put together, these maps build an isomorphism of vector bundles over $H$,
\begin{eqnarray*}
\xymatrix{
R:t^*\CC{\Omega} \ar[r] & \ker(\hat{s}):(h,c) \ar@{|->}[r] & \hat{m}(c;0_h).
}
\end{eqnarray*}
This isomorphism justifies why the sequence \ref{CoreSeq} is referred to as the core sequence of the \LA -groupoid. This construction works out for VB-groupoids in general; however, in the presence of an \LA -groupoid structure, there is an additional structure of Lie algebroid inherited by the core. We proceed to describe it. Consider the map
\begin{eqnarray*}
\xymatrix{
\Gamma(\CC{\Omega}) \ar[r] & \Gamma(\Omega):\sigma \ar@{|->}[r] & \vec{\sigma},
}
\end{eqnarray*}
where 
\begin{eqnarray*}
\vec{\sigma}_h:=R_h(\sigma_{t(h)})
\end{eqnarray*}
for $h\in H$. This map is an isomorphism onto its image, thus prompting us to define
\begin{eqnarray*}
[\sigma_1,\sigma_2]_{\CC{\Omega}}:=[\vec{\sigma}_1,\vec{\sigma}_2]_\Omega\rest{M}.
\end{eqnarray*}
As for the anchor, notice that if we regard $c\in\CC{\Omega}$ as an element of $\Omega\rest{M}$, $\rho_\Omega(c)\in\ker (ds)\rest{M}$. Indeed,
\begin{eqnarray*}
ds(\rho_\Omega(c))=\rho_A(\hat{s}(c))=0;
\end{eqnarray*}
however, using the anchor of $\Omega$ cannot be the right Ansatz as $\ker (ds)\leq TH$, and our anchor needs to have $TM$ as codomain. Now, since $\ker (ds)\rest{M}=A_H$, we use its anchor to define
\begin{eqnarray*}
\rho_{\CC{\Omega}}(c):=dt\circ\rho_\Omega(c)=\rho_A\circ\hat{t}(c)
\end{eqnarray*}
\begin{prop}
$(\CC{\Omega},[,]_{\CC{\Omega}},\rho_{\CC{\Omega}})$ is a Lie algebroid. 
\end{prop}
%\begin{rmk}As we will see (we won't, replace), the fact that the anchor of the core is borrowed from the Lie algebroid $A_H$ of the base groupoid implies that it builds a connection between the \LA -groupoid itself and $A_H$ in terms of certain relations between their isotropies and their foliations.\end{rmk} 
For the global counter-part of this construction, start out with a double Lie groupoid $D$, then
\begin{Def}
The \textit{core} of $D$ is given by $\CC{D}:=\mathbb{S}^{-1}((u\times\Rg{u})(M))$. In other words, 
\begin{eqnarray*}
\CC{D}:=\lbrace d\in D:\Tp{s}(d)\in\Rg{u}(M),\quad\Lf{s}(d)\in u(M)\rbrace.
\end{eqnarray*}
\end{Def} 
Spelling out which elements of $D$ belong to $\CC{D}$, one sees that they are precisely the squares of the form
\begin{eqnarray*}
\xymatrix{
{}^y\bullet        & \bullet^x \ar[l] \\
{}_x\bullet \ar[u] & \bullet_x \ar@{=}[l]\ar@{=}[u] .
}
\end{eqnarray*}
Thanks to the hypothesis that the double source be a submersion, one figures that the core is a manifold. We use the structure of the double Lie groupoid to endow $\CC{D}$ with a Lie groupoid structure over $M$. Clearly, the reasonable thing to do is to define
\begin{eqnarray*}
s_{\CC{D}}(d):=s^2(d), & t_{\CC{D}}(d):=t^2(d), & u_{\CC{D}}(x):=u^2(x),
\end{eqnarray*}
for $d\in\CC{D}$ and $x\in M$, where $s^2=\Rg{s}\circ\Tp{s}=s\circ\Lf{s}$, $t^2=\Rg{t}\circ\Tp{t}=t\circ\Lf{t}$ and $u^2=\Tp{u}\circ\Rg{u}=\Lf{u}\circ u$.
Now, given a pair $(d_1,d_2)\in\CC{D}_{s^2}\times_{t^2}\CC{D}$, we can think of it as the following 
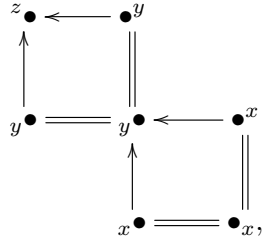
\begin{eqnarray*}
\xymatrix{
{}^z\bullet        & \bullet^y \ar[l]                 &                  \\
{}_y\bullet \ar[u] & {}_y\bullet \ar@{=}[l]\ar@{=}[u] & \bullet^x \ar[l] \\
                   & {}_x\bullet \ar[u]               & \bullet_x , \ar@{=}[l]\ar@{=}[u]                      
}
\end{eqnarray*}
we define the multiplication by filling this square with the corresponding units
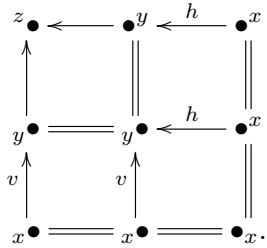
\begin{eqnarray*}
\xymatrix{
{}^z\bullet          & \bullet^y \ar[l]                 & \bullet^x \ar[l]_h             \\
{}_y\bullet \ar[u]   & {}_y\bullet \ar@{=}[l]\ar@{=}[u] & \bullet^x \ar[l]_h\ar@{=}[u]   \\
{}_x\bullet \ar[u]^v & {}_x\bullet \ar[u]^v\ar@{=}[l]   & \bullet_x . \ar@{=}[l]\ar@{=}[u]                     
}
\end{eqnarray*}
In symbols,
\begin{eqnarray*}
d_1\boxplus d_2:=(d_1\vJoin\Tp{u}\circ\Tp{t}(d_2))\Join d_2=(d_1\Join\Lf{u}\circ\Lf{t}(d_2))\vJoin d_2.
\end{eqnarray*}
Notice that these multiplications make sense. Indeed, on the one hand,
\begin{align*}
\Lf{s}(d_1)=u(s^2(d_1))=u(t^2(d_2))=u\circ\Rg{t}(\Tp{t}(d_2))=\Lf{t}(\Tp{u}\circ\Tp{t}(d_2)) & \Longrightarrow (d_1;\Tp{u}\circ\Tp{t}(d_2))\in D_{\Lf{s}}\times_{\Lf{t}}D, \\
\Tp{s}(d_1)=\Rg{u}(s^2(d_1))=\Rg{u}(t^2(d_2))=\Rg{u}\circ t(\Lf{t}(d_2))=\Tp{t}(\Lf{u}\circ\Lf{t}(d_2)) & \Longrightarrow (d_1;\Lf{u}\circ\Lf{t}(d_2))\in D_{\Tp{s}}\times_{\Tp{t}}D,
\end{align*}
and on the other,
\begin{align*}
\Tp{s}(d_1\vJoin\Tp{u}\circ\Tp{t}(d_2)) & =\Rg{m}(\Tp{s}(d_1);\Tp{s}\circ\Tp{u}(\Tp{t}(d_2))) \\
										& =\Rg{m}(\Rg{u}(s^2(d_1));\Tp{t}(d_2))=\Tp{t}(d_2)\Longrightarrow (d_1\vJoin\Tp{u}\circ\Tp{t}(d_2);d_2)\in D_{\Tp{s}}\times_{\Tp{t}}D, \\
\Lf{s}(d_1\Join\Lf{u}\circ\Lf{t}(d_2))  & =\Lf{s}(d_1)\Lf{s}\circ\Lf{u}(\Lf{t}(d_2))) \\
										& =u(s^2(d_1))\Lf{t}(d_2)=\Lf{t}(d_2)\Longrightarrow (d_1\Join\Lf{u}\circ\Lf{t}(d_2);d_2)\in D_{\Lf{s}}\times_{\Lf{t}}D.
\end{align*}
\begin{prop} % Look for a reference (I believe it is indicated in one of Raj's papers
$\xymatrix{\CC{D} \ar@<0.5ex>[r]\ar@<-0.5ex>[r] & M}$ is a Lie groupoid. 
\end{prop}
\begin{rmk}
As it was said before, the construction of the core of an \LA -groupoid works out for VB-groupoids. Interestingly enough, if one regards such a VB-groupoid as a double Lie groupoid, its cores coincide.  
\end{rmk}
As one would expect, the two constructions outlined above are related by the Lie functor.
\begin{prop}
Let $D$ be a double Lie groupoid. If we write $A_D$ for the \LA -groupoid of $D$, then
\begin{eqnarray*}
A_{\CC{D}}=\CC{A_D}.
\end{eqnarray*}
\end{prop}
\begin{rmk}
Notice that the spaces that appear as domains of the crossed module structural map in both the Lie $2$-algebra case and the Lie $2$-group case are the core of $\gg_1$ and the core of $\G$ regarded as an \LA -groupoid and as a double Lie groupoid respectively.
\end{rmk}
We close this subsection with a collection of examples. % of cores 
\begin{ex:}\label{core1,2}
When regarded as zero \LA -groupoids or as unit double Lie groupoids, Lie groupoids clearly have zero core. This is the case as well for Lie algebroids and their integrations (if exist).
\end{ex:}
\begin{ex:}\label{core3}
The tangent prolongation of $G$ has $A_G$ as core. Indeed, by definition $A_G=\ker ds\rest{M}$. According to the latter proposition any double Lie groupoid integrating it should have core an integration of $A_G$. For example, the pair double Lie groupoid $G\times G$ has core isomorphic to $G$.
% In case \Pi(G) admits a double Lie groupoid structure, what's its core? It would seem as though it is \tilde{G}. Indeed, a path \gamma in the core has source isomorphic to the constant path x. From the contractibility of s\circ\gamma, one should be able to lift the homotopy, to yield a path on s^{-1}(x) starting at a unit.
\end{ex:}
\begin{ex:}\label{core4}
Being sub-groupoids of the tangent prolongation, multiplicative foliations on $G$ have Lie sub-algebroids $F_G$ of $A_G$ as cores. We know very little of their global counterparts, but in concordance with the previous example, the core of any double Lie groupoid integrating a multiplicative foliation needs to be a Lie groupoid integrating $F_G$.
\end{ex:}
The following examples are more geometric.
\begin{ex:}\label{core7}
It is well known that a Poisson structure on a manifold $M$ gives rise to an algebroid structure on the cotangent bundle $T^*M$. In the case of a multiplicative Poisson bivector, one can prove that the dual of the linearization of the Poisson bivector at the identity element $e$ defines a Lie algebra structure on $\gg^*$. Along with this algebroid structure,
\begin{eqnarray}\label{PLG}
\xymatrix{
 T^*G \ar@<0.5ex>[r] \ar@<-0.5ex>[r] \ar[d]  & \gg^* \ar[d] \\
 G \ar@<0.5ex>[r] \ar@<-0.5ex>[r] & \ast 
}
\end{eqnarray}
turns into an \LA -groupoid, where top groupoid is the usual structure of the cotangent groupoid. The global counter-parts of these are simultaneous symplectic realizations of pairs of Poisson Lie groups in duality.\\
Poisson-Lie groups have as special defining treat that they have got no core. Indeed, this follows by counting dimensions, as the source in $\xymatrix{T^*G \ar@<0.5ex>[r] \ar@<-0.5ex>[r] & \gg^*}$ is defined by $\hat{s}(\xi)=(dL_g)^*(\xi)$ for $\xi\in T_g^*G$, which has trivial kernel. An \LA -groupoid with this property will be called vacant, and any double Lie groupoid integrating it will have trivial core as well. 
\end{ex:}
\begin{ex:}\label{core8}
The general principle that endows $\gg^*$ with a Lie algebra structure in the presence of a multiplicative Poisson bivector, carries over to the groupoid world. That is, given a Poisson-Lie groupoid, there is an induced Lie algebroid structure on the dual of its Lie algebroid; hence, the cotangent groupoid gets upgraded to an \LA -groupoid. As an extension of this, one can consider multiplicative Dirac structures. A Dirac structure on $M$ is a maximally coisotropic involutive sub-bundle of $TM\oplus T^*M$. By analogy \cite{DiracLieGps}, a multiplicative Dirac structure on a Lie group is a maximally coisotropic involutive sub-object of $TG\oplus T^*G$. The key observation is to note that the latter is no longer simply a vector bundle, but indeed a VB-groupoid $\xymatrix{TG\oplus T^*G \ar@<0.5ex>[r] \ar@<-0.5ex>[r] & \gg^*}$. The associated Lie algebroid of this Dirac structure is an \LA -groupoid. As in the case of multiplicative Poisson structures, this idea carries over to the groupoid case, where the ``background'' object is $\xymatrix{TG\oplus T^*G \ar@<0.5ex>[r] \ar@<-0.5ex>[r] & TM\oplus A_G^*}$.
As a matter of fact, morphisms of double objects induce morphisms between their cores. Since Dirac Lie groups are sub-objects of $TG\oplus T^*G$, their cores will be sub-objects of the core of the VB-groupoid $\xymatrix{TG\oplus T^*G \ar@<0.5ex>[r] \ar@<-0.5ex>[r] & \gg^*}$, which is the Lie algebra $\gg$ of $G$. Since Dirac Lie groups are known to induce \LA -groupoids, their cores will be sub-algebras of $\gg$. %In fact, by carefully noticing what the unit Lie algebra is the annihilator of the Lie algebra of the leaf through the identity, the core will be that Lie algebra... or is it not? The observation is: The structure on the base is simply given by the restriction of the Dirac structure to T_e^*G. Since the thing is isotropic, necessarily, an element $(x,\xi)\in L_e$ <(x,\xi),(x,\xi)>=2\xi(x)=0. One concludes then that the space of units is the annihilator of the core.
In the groupoid case, one will have a sub-object of the core of the VB-groupoid $\xymatrix{TG\oplus T^*G \ar@<0.5ex>[r] \ar@<-0.5ex>[r] & TM\oplus A_G^*}$, which is $A_G\otimes T^*M$. % Do I have more data about this?
\end{ex:}

% -------------------------
% -------------------------
% Appendix C
% -------------------------
% -------------------------

\section{Appendix C: VB-groupoids and representations up to homotopy}
VB-groupoids and VB-algebroids \cite{VB&Reps,2ndOrdGeom} can be thought of as vector bundles in the categories of Lie groupoids and Lie algebroids. Paradigmatic examples include the tangent and cotangent bundles of Lie groupoids and Lie algebroids and, as we will see shortly, classic representations of Lie groupoids. As we pointed out, these coincide with double Lie groupoids and \LA -groupoids whose groupoid structure is flat and abelian. In \cite{BCD}, the authors study the Lie theory for VB-algebroids successfully proving a Lie III type theorem. The strategy is to regard vector bundles as a special type of homogeneous spaces \cite{Or}, though we will not be delving deeper into this matter. \\
In this appendix, we would like to make explicit the relation between VB-groupoids and representations. Recall that, one can equivalently think of a VB-groupoid as being an abelian \LA -groupoid. Now, given a (left) representation of a Lie groupoid $\xymatrix{G \ar@<-0.5ex>[r]\ar@<0.5ex>[r] & M}$ on a vector bundle $E$ over $M$, it is easy to see that
\begin{eqnarray*}
\xymatrix{
G\ltimes E \ar@<-0.5ex>[r]\ar@<0.5ex>[r]\ar[d] & E \ar[d] \\
G \ar@<-0.5ex>[r]\ar@<0.5ex>[r]                & M
}
\end{eqnarray*}
is a VB-groupoid whose core is trivial, i.e. Seen as a double Lie groupoid, it is the unit groupoid of $M$, and if seen as an abelian \LA -groupoid, it is the zero bundle over $M$. As it turns out, classic representations of Lie groupoids are in one-to-one correspondence with these \textit{coreless} VB-groupoids. In a sense, these are too few representations. For example, with this notion a non-zero representation of the pair groupoid can only take values on trivial vector bundles. More importantly, though, with this definition, a Lie groupoid lacks an adjoint representation. Arguably motivated by these considerations, the more general representations up to homotopy were introduced in \cite{AAC, AAC2}. Representations up to homotopy take values on graded vector bundles over $M$ instead of vector bundles, say $\E=\bigoplus E_k$. There is a graded right $C(G)$-module
\begin{eqnarray*}
C(G,\E)^n=\bigoplus_{p-q=n}C^p(G,E_q),
\end{eqnarray*}
and a representation up to homotopy is defined to be a differential for this. Formally, 
\begin{Def}
A \textit{representation up to homotopy} of $G$ on $\E$ is a continuous degree $1$ operator $D$ on $C(G,\E)$ such that
\begin{eqnarray*}
D^2=0 & \textnormal{and} & D(\omega\star f)=(D\omega)\star f+(-1)^{p+q}\omega\star\partial f
\end{eqnarray*}
for all $\omega\in C^p(G,E_q)$ and $f\in C(G)$.
\end{Def}
We will chiefly be concerned about representations up to homotopy valued on graded vector bundles concentrated in degrees $0$ and $1$. These play an important r\^{o}le in the representation theory of Lie groupoids, as the adjoint representation can be thought of as being one of these. Before stating it as an example, we will spell out what the definition amounts to in this restricted case and will go over a useful correspondence. The following material can be found in \cite{VB&Reps}. Let $\E =E\oplus C[1]$, then
\begin{eqnarray*}
C(G,\E )^n=C^n(G,E)\oplus C^{n+1}(G,C)
\end{eqnarray*}
and in this direct sum decomposition a degree $1$ operator has a matrix decomposition whose four components are 
\begin{align*}
D^E    & :\xymatrix{C^n(G,E) \ar[r] & C^{n+1}(G,E) } & \partial & :\xymatrix{C^{n+1}(G,C) \ar[r] & C^{n+1}(G,E) } \\
\Omega & :\xymatrix{C^n(G,E) \ar[r] & C^{n+2}(G,C) } & D^C      & :\xymatrix{C^{n+1}(G,C) \ar[r] & C^{n+2}(G,C) }.
\end{align*}
We reinterpret the conditions in the definition for these components. First, the Leibniz rule says that $\partial$ and $\Omega$ are maps of $C(G)$-modules. In particular, for $n=-1$, we have got a map between the modules of sections which translate to a map of vector bundles over the identity
\begin{eqnarray*}
\xymatrix{
C \ar[r] & E:c \ar@{|->}[r] & \partial(c)
}
\end{eqnarray*} 
via the equivalence of categories established by the Serre-Swan theorem. As for $\Omega$, for $n=1$ it induces a section of the vector bundle $s_2^*E^*\otimes t_2^*C$ by the formula
\begin{eqnarray*}
\Omega_{(g_1,g_2)}(e)=\Omega(\epsilon)(g_1,g_2)
\end{eqnarray*}
where $(g_1,g_2)\in G^{(2)}$, $e\in E_{s(g_2)}$ and $\epsilon\in\Gamma(E)$ is any section satisfying $\epsilon(s(g_1))=e$. This formula is well-defined. Indeed, invoking Serre-Swan yet again, $\Gamma(E)$ is a finitely generated $C(M)$-module. Let $\lbrace \alpha_j\rbrace$ be a generating set, and $\epsilon_1,\epsilon_2\in\Gamma(E)$ be two sections such that $\epsilon_1(s(g))=e=\epsilon_2(s(g))$. Using the Einstein summation convention
\begin{eqnarray*}
\epsilon_1 =\epsilon_1^j\alpha_j & \textnormal{and} & \epsilon_2 =\epsilon_2^j\alpha_j.
\end{eqnarray*}
Since $\Omega$ is a $C(G)$-module homomorphism, 
\begin{align*}
\Omega(\epsilon_k^j\alpha_j)(g_1,g_2) & =\Omega(\alpha_j\star\epsilon_k^j)(g_1,g_2)=\Omega\alpha_j(g_1,g_2)\epsilon_k^j(s(g_1)),	     
\end{align*}
but the last term in this equation does not depend on $k$. \\
On the other hand, the fact that $D$ verifies the Leibniz condition of the definition, implies that there are quasi-actions
\begin{align*}
\xymatrix{
G\times_M E \ar[r] & E:(g,e) \ar@{|->}[r] & \Delta^E_g e ,
} \\
\xymatrix{
G\times_M C \ar[r] & C:(g,c) \ar@{|->}[r] & \Delta^C_g c .
}
\end{align*}
These quasi-actions are defined by the corresponding operators in the following manner. First, let $\epsilon\in\Gamma(E)$ and $\kappa\in\Gamma(C)$ be sections such that 
\begin{eqnarray*}
\epsilon(s(g))=e & \textnormal{and} & \kappa(s(g))=c .
\end{eqnarray*}
Then, define
\begin{align*}
\Delta^E_g e & := D^E(\epsilon)(g) + \epsilon(t(g)) \in E_{t(g)} \\
\Delta^C_g c & := -D^C(\kappa)(g) + \kappa(t(g))    \in C_{t(g)}.
\end{align*}
We claim that these expressions are well-defined. Indeed, using the Serre-Swan theorem one last time and the notation above, for fixed $j$, the Leibniz rule yields
\begin{align*}
D^E(\epsilon_k^j\alpha_j)(g) & =D^E(\alpha_j\star\epsilon_k^j)(g) \\
						     & =(D^E\alpha_j)\star\epsilon_k^j(g)+\alpha_j\star\partial\epsilon_k^j(g) \\
						     & =D^E\alpha_j(g)\epsilon_k^j(s(g))+\alpha_j(t(g))(\epsilon_k^j(s(g))-\epsilon_k^j(t(g))).
\end{align*}
Summing over $j$, and rearranging terms we get
\begin{eqnarray*}
D^E(\epsilon_k)(g)+\epsilon_k(t(g)) = (D^E\alpha_j(g)+\alpha_j(t(g)))\epsilon_k^j(s(g)).
\end{eqnarray*}
Since the right hand side does not depend on $k$, the claim follows. An analogous computation, where the Leibniz rule comes with a sign due to the degree of $\kappa$, or more specifically, due to the degree of the elements in the generating set of $\Gamma(C)$, proves that $\Delta^C$ is well-defined as well. \\
Turning our attention to the equation $D^2=0$, we compute 
\begin{align*}
(D^E)^2+\partial\Omega & =0 & D^E\partial+\partial D^C & =0 \\
\Omega D^E+D^C\Omega   & =0 & \Omega\partial+(D^C)^2   & =0 ;
\end{align*}
in order, this implies the following equations hold for the quasi-actions above and for every $(g_1,g_2,g_3)\in G^{(3)}$,
\begin{align*}
\Delta^E_{g_1}\Delta^E_{g_2}-\Delta^E_{g_1g_2}+\partial\Omega_{(g_1,g_2)} & =0 & \Delta^E\partial-\partial\Delta^C & =0 \\
\Delta^C_{g_1}\Omega_{(g_2,g_3)}-\Omega_{(g_1g_2,g_3)}+\Omega_{(g_1,g_2g_3)}-\Omega_{(g_1,g_2)}\Delta^E_{g_3} & =0 & \Omega_{(g_1,g_2)}\partial+\Delta^C_{g_1}\Delta^C_{g_2}-\Delta^C_{g_1g_2}   & =0. 
\end{align*}
The key observation to see how this set of equations follows from the first one is to realize that, since $\Gamma(t_p^*E)\cong C\big{(}G^{(p)}\big{)}\otimes\Gamma(E)$ regarded as a $C\big{(}G^{(p)}\big{)}$-module, $\Gamma(t_p^*E)$ is generated by $\lbrace t_p^*\alpha_j\rbrace$. %Should I add the actual computations? ***. 
As it is explained in \cite{VB&Reps}, the data $(\partial,\Delta^E,\Delta^C,\Omega)$ subject to the latter set of equations is all that a representation up to homotopy on a $2$-term vector bundle amounts to. \\
Now, given a representation up to homotopy $(\partial,\Delta^E,\Delta^C,\Omega)$, there is a notion of semi-direct product that extends the construction of a VB-groupoid out of an honest representation which we proceed to describe. Let 
\begin{eqnarray*}
\Gamma :=t^*C\oplus  s^*E\cong\lbrace(c,g,e)\in C\times G\times E : c\in C_{t(g)},e\in E_{s(g)}\rbrace
\end{eqnarray*}
and define a structure of Lie groupoid over $E$ using the maps
\begin{eqnarray*}
\hat{s}(c,g,e):=e & \hat{t}(c,g,e):=\partial c+\Delta^E_g e & \hat{u}(e)=(0,u(x),e),
\end{eqnarray*}
where in the definition of the unit, $e\in E_x$. Clearly, these maps cover the corresponding maps and are fibre-wise linear. For the multiplication, let $(c_1,g_1,e_1;c_2,g_2,e_2)\in\Gamma^{(2)}$ and define
\begin{eqnarray*}
\hat{m}(c_1,g_1,e_1;c_2,g_2,e_2):=(c_1+\Delta^C_{g_1}c_2-\Omega_{(g_1,g_2)}e_2,g_1g_2,e_2).
\end{eqnarray*}
Again, this map clearly covers the multiplication and is indeed a map of vector bundles. After a computation, one sees that this multiplication admits inverses with formula
\begin{eqnarray*}
\hat{\iota}(c,g,e)=(-\Delta^C_{g^{-1}}c+\Omega_{(g^{-1},g)}e,g^{-1},\partial c+\Delta^E_g e);
\end{eqnarray*}
which ultimately shows that this whole structure fits into a VB-groupoid over $\xymatrix{G \ar@<-0.5ex>[r]\ar@<0.5ex>[r] & M}$. \\
Conversely, starting out with a VB-groupoid $\xymatrix{\Gamma \ar@<-0.5ex>[r]\ar@<0.5ex>[r] & E}$ and regarding it as an abelian \LA -groupoid, one can consider its core sequence
\begin{eqnarray*}
\xymatrix{
(0) \ar[r] & t^*C \ar[r] & \Gamma \ar[r] & s^*E \ar[r] & (0).
}
\end{eqnarray*}
Picking a splitting $\xymatrix{h:G\times_M E \ar[r] & \Gamma:(g,e) \ar@{|->}[r] & h_g(e)}$ such that $h_{u(x)}(e)=\hat{u}(e)$, there is an isomorphism $\Gamma\cong t^*C\oplus s^*E$. Using the VB-groupoid structure, we get 
\begin{eqnarray*}
\xymatrix{
\partial :C \ar[r] & E:c \ar@{|->}[r] & \hat{t}(c).
}
\end{eqnarray*}
We also have got quasi-actions
\begin{align*}
\Delta^E_g e & := \hat{t}(h_g(e)) \\
\Delta^C_g c & := h_g(\partial(c))\cdot c\cdot 0_{g^{-1}},
\end{align*}
and finally,
\begin{eqnarray*}
\Omega_{(g_1,g_2)}e=(h_{g_1g_2}(e)-h_{g_1}(\Delta^E_{g_2}e)\cdot h_{g_2}(e))\cdot 0_{(g_1g_2)^{-1}}.
\end{eqnarray*}
In fact, the main theorem in \cite{VB&Reps} says that the above constructions are essentially reverse to one another.
\begin{theorem}
There is a one-to-one correspondence between isomorphism classes of VB-groupoids over $E$ with core $C$ and gauge-equivalence classes of representations up to homotopy on $E\oplus C[1]$.
\end{theorem}
Later \cite{CristMat}, it was proven that this correspondence can be upgraded to an equivalence of categories.

%--------------------------------------
% MULTIPLICATIVE (bi)SECTIONS 

%*** to make sense of these, there should be some definitions saved: double gpd, \LA -gpd, cores...
%--------------------------------------

\section{Appendix D: an alternative perspective}
%We would like to start this appendix by outlining the relationship between the integrability of the Lie algebroid of arrows in an \LA -groupoid and the integrability of its base and its core. (are you though?)
% If you're seriously gonna do this, you better understand the idea of Cañez. Or look at the diagrams in 04.03. relating monodromies, isotropies and fund. gps of orbits.***Next up: If the Lie algebroid of arrows of an \LA -groupoid is integrable, so are the base and the core (!). Also, the relations among the different foliations and the ones among the isotropies. *** \\
%*** Related to the integrability of the unit algebroid. Is there an easy and neglected necessary condition for the integrability in terms of the the developability of the base? cf. Ieke's work on developable foliations conclude something along the lines of ``The subgroupoid interating a given Lie subalgebroid is embedded if and only if the foliation is developable.'' *** \\ 
% If so, if it doesn't hold in general that \Gamma_mult(\Omega) is a Lie 2 algebra, how come? what changed?
%In this subsection, we propose an alternative approach to the problem of integrating \LA -groupoids. 
Recall that by definition the space of sections of a Lie algebroid has a infinite dimensional Lie algebra structure. On the other hand, if one is given a Lie groupoid, there is a group structure on the space of bisections. In general, a bisection of the Lie groupoid $\xymatrix{G \ar@<-0.5ex>[r]\ar@<0.5ex>[r] & M}$ can be defined as a map $\xymatrix{\sigma:M \ar[r] & G}$ which is a section of the source map and for which $t\circ\sigma$ is a diffeomorphism of $M$. We call the space of all such bisections $\Bis(G)$. The group structure will be given by the rule
\begin{eqnarray*}
(\sigma_1\sigma_2)(x):=\sigma_1(t(\sigma_2(x)))\sigma_2(x).
\end{eqnarray*}
This defines a bisection as
\begin{align*}
s((\sigma_1\sigma_2)(x)) & =s(\sigma_2(x))=x & t((\sigma_1\sigma_2)(x)) & = (t\circ\sigma_1)((t\circ\sigma_2)(x)).
\end{align*}
Furthermore, as the product rule is defined in terms of the multiplication of the groupoid, it will inherit its properties:
\begin{align*}
(\sigma_1(\sigma_2\sigma_3))(x) & =\sigma_1(t((\sigma_2\sigma_3)(x)))(\sigma_2\sigma_3)(x) \\
								& =\sigma_1((t\circ\sigma_2)((t\circ\sigma_3)(x)))\Big{(}\sigma_2(t(\sigma_3(x)))\sigma_3(x)\Big{)} \\
								& =\Big{(}\sigma_1(t(\sigma_2(t(\sigma_3(x)))))\sigma_2(t(\sigma_3(x)))\Big{)}\sigma_3(x) \\
								& =(\sigma_1\sigma_2)(t(\sigma_3(x)))\sigma_3(x) =((\sigma_1\sigma_2)\sigma_3)(x);
\end{align*}
therefore, the product is associative. Clearly, $u$ is a bisection and it satisfies 
\begin{eqnarray*}
(u\sigma)(x)=u(t(\sigma(x)))\sigma(x)=\sigma(x) & \textnormal{and} & (\sigma u)(x)=\sigma(t(u(x)))u(x)=\sigma(x),
\end{eqnarray*}
so it plays the r\^{o}le of the identity element. Finally, defining $\sigma^{\dagger}:=\iota\circ\sigma\circ(t\circ\sigma)^{-1}$, we get the desired
\begin{align*}
(\sigma\sigma^{\dagger})(x) & =\sigma(t(\iota(\sigma\circ(t\circ\sigma)^{-1}(x))))\iota\circ\sigma((t\circ\sigma)^{-1}(x))  &  (\sigma^{\dagger}\sigma)(x) & =\iota\circ\sigma((t\circ\sigma)^{-1}(t(\sigma(x))))\sigma(x) \\
                       & =\sigma(s(\sigma\circ(t\circ\sigma)^{-1}(x)))\iota\circ\sigma(t(\sigma)                       &   & =\iota(\sigma(x))\sigma(x) \\
                       & =\sigma((t\circ\sigma)^{-1}(x))\iota(\sigma(t\circ\sigma(x)))                                 &   & = u(s(\sigma(x))) =u(x). \\
                       & =u(t(\sigma((t\circ\sigma)^{-1}(x))))=u(x),
\end{align*}
Now, it is rather easy to see that using the compact-open % 
topology, $\Bis(G)$ has got the structure of a topological group. One needs to be extremely careful if one is to regard this object as an infinite dimensional Lie group. Just until fairly recently \cite{LieGpBis} was it proven that $\Bis(G)$ is not only a topological group, but it admits the structure of a locally convex manifold making it an infinite dimensional Lie group whose Lie algebra is $\Gamma(A_G)$. This however, is done under restrictive hypothesis. Apart from this reference (see also the references therein) and \cite{LagBis,LieAlgBis}, this subject is surprisingly missing from the literature. This hints in the direction that the problem of integrating a given Lie algebroid $A$ should be related to integrating the Lie algebra $\Gamma(A)$ with the caveat that once integrated to an infinite dimensional Lie group, this would need to be realized as the space of bisections of a Lie groupoid (direction in which there are no attempts in the literature to the best of our knowledge, although related issues are explored in \cite{FromBisToGpd}). Following an analogous reasoning, one can ask whether there is an algebraic structure on the space of compatible sections of an \LA -groupoid. This is the contents of this appendix.
\subsubsection{The category of multiplicative sections}
Let $\Omega$ be an \LA -groupoid over $G$. Define the category of multiplicative sections as the space of groupoid morphisms $\xymatrix{G \ar[r] & \Omega}$ that are also sections at both the level of objects and of morphisms. Since such maps are functors, we consider natural transformations between them as morphisms; these, however, will be asked to induce the identity after horizontal composition with the projection. We spell out this concepts carefully. \\
First, let 
\begin{eqnarray*}
\Gamma_{mult}(\Omega)_0:=\lbrace\omega\in\Gamma(\Omega):\omega\textnormal{ Lie groupoid homomorphism}\rbrace .
\end{eqnarray*}
Thanks to the linearity of the structural maps of the groupoid $\Omega$, $\Gamma_{mult}(\Omega)$ is a vector subspace of $\Gamma(\Omega)$. Given $\omega_1,\omega_2\in\Gamma_{mult}(\Omega)_0$, we define $\Gamma_{mult}(\Omega)(\omega_1,\omega_2)$ to be a subset of smooth natural transformations between $\omega_1$ and $\omega_2$. By a smooth natural transformation, we mean a smooth map $\xymatrix{\tau :M \ar[r] & \Omega}$ such that 
\begin{eqnarray*}
\hat{s}(\tau(x))=\omega_1(x) & \textnormal{and} & \hat{t}(\tau(x))=\omega_2(x),
\end{eqnarray*}
and as usual, such that for each $g\in G(x,y)$, the diagram
\begin{eqnarray*}
\xymatrix{
\omega_1(x) \ar[r]^{\tau(x)}\ar[d]_{\omega_1(g)} & \omega_2(x) \ar[d]^{\omega_2(g)} \\
\omega_1(y) \ar[r]_{\tau(y)}                     & \omega_2(y)
}
\end{eqnarray*}
commutes in $\Omega$. We write the full space of such smooth natural transformations as $Nat^\infty(\omega_1,\omega_2)$ and with this notation,
\begin{eqnarray*}
\Gamma_{mult}(\Omega)(\omega_1,\omega_2):= \lbrace\tau\in Nat^\infty(\omega_1,\omega_2):\tau\ast_{h}1_{\pi} =1_{id_G}\rbrace .
\end{eqnarray*}
The compatibility equation that appears in this definition is to make the defining equation for a section come into play. Diagrammatically, 
%\begin{eqnarray*}
%\xymatrix{
%G \ar@/^{1pc}/[rr]^{\omega_1}="a"\ar@/_{1pc}/[rr]_{\omega_2}="b" & \ar@{=>}^{\tau}"a";"b" & \Gamma \ar@/^{1pc}/[rr]^{\pi}="c"\ar@/_{1pc}/[rr]_{\pi}="d" & \ar@{=>}^{1_\pi}"c";"d" & G
%} = \xymatrix{
%G \ar@/^{1pc}/[rr]^{id_G}="e"\ar@/_{1pc}/[rr]_{id_G}="f" & \ar@{=>}^{1_{id_G}}"e";"f" & G .
%}
%\end{eqnarray*}
\begin{align*}
\xymatrix{
G \ar@/^{1pc}/[rr]^{\omega_1}\ar@/_{1pc}/[rr]_{\omega_2} & \quad\Big{\Downarrow} \small{\tau} & \Gamma \ar@/^{1pc}/[rr]^{\pi}\ar@/_{1pc}/[rr]_{\pi} & \quad\Big{\Downarrow} \small{1_{\pi}} & G
} & =\xymatrix{
G \ar@/^{1pc}/[rr]^{id_G}\ar@/_{1pc}/[rr]_{id_G} & \qquad\Big{\Downarrow} \small{1_{id_G}} & G .
}
\end{align*}
Of course, with this diagram, it is clear that the composition of two given natural transformations verifying the defining equation, will verify it as well. As for the smoothness, it will follow from the smoothness of the multiplication in $\Omega$. Finally, for any $\sigma\in\Gamma_{mult}(\Omega)$, $1_\sigma$ is defined using $\hat{u}$, which ensures its smoothness, and clearly satisfies the defining equation. Summing up, $\Gamma_{mult}(\Omega)$ is indeed a category. \\
$\Gamma_{mult}(\Omega)$ is a vector space as well. Its linear structure will be inherited from the that of the fibres of $\Omega$ in the following manner: For $\tau_1,\tau_2\tau\in\Gamma_{mult}(\Omega)$ and $\lambda\in\Rr$, set 
\begin{eqnarray*}
(\tau_1+\tau_2)(x):=\tau_1(x)+\tau_2(x) & \textnormal{and} & (\lambda\tau)(x):=\lambda\tau(x).
\end{eqnarray*}
To see that the additive structure is well-defined, notice that due to the defining equation, 
\begin{align*}
(\tau\ast_{h}1_\pi)(x) & =\pi(\tau(x))1_\pi(\omega_1(x)) \\
					   & =\pi(\tau(x))u(p(\omega_1(x))) =\pi(\tau(x))
\end{align*}
has to coincide with $1_{id_G}(x)=u(x)$ and therefore, both $\tau_1(x)$ and $\tau_2(x)$ belong to $\Omega_{u(x)}$. This not only prove that the defining formula for the addition makes sense, but incidentally proves that $\tau_1+\tau_2$ verifies the defining equation for $\Gamma_{mult}(\Omega)$. Though we did not mentioned it, it makes sense to define both source and target maps as being linear, that is, if 
\begin{eqnarray*}
\xymatrix{
\tau_1:\omega_1 \ar@{=>}[r] & \omega'_1
} & \textnormal{and} & \xymatrix{
\tau_2:\omega_2 \ar@{=>}[r] & \omega'_2 ,
}
\end{eqnarray*} 
then define
\begin{eqnarray*}
\xymatrix{
\tau_1+\tau_2:\omega_1+\omega_2 \ar@{=>}[r] & \omega'_1+\omega'_2
} & \textnormal{and} & \xymatrix{
\lambda\tau_1:\lambda\omega_1 \ar@{=>}[r] & \lambda\omega'_1 .
}
\end{eqnarray*} 
Indeed, since both $\hat{s}$ and $\hat{t}$ are linear, this prescription does not conflict with the formulas above. Moreover, since the multiplication $\hat{m}$ is linear as well, for each $g\in G(x,y)$, the diagram
\begin{eqnarray*}
\xymatrix{
\omega_1+\omega_2(x) \ar[r]^{\tau_1+\tau_2(x)}\ar[d]_{\omega_1+\omega_2(g)} & \omega'_1+\omega'_2(x) \ar[d]^{\omega'_1+\omega'_2(g)} \\
\omega_1+\omega_2(y) \ar[r]_{\tau_1+\tau_2(y)}                              & \omega'_1+\omega'_2(y)
}
\end{eqnarray*}
and its analogous for the multiplication by scalar do commute. In fact, the linearity of this maps hints in the direction that $\Gamma_{mult}(\Omega)$ is actually a $2$-vector space. The unit map is linear, because it is defined in therms of the linear $\hat{u}$, 
\begin{eqnarray*}
1_{\omega}(x):=\hat{u}(\omega(x)),
\end{eqnarray*}
as for the vertical composition, it is once again defined using the multiplication $\hat{m}$ thus proving it is linear as well. \\
We saw that there is an equivalence between the categories of Lie $2$-algebras and that of crossed modules. This correspondence carries to the infinite dimensional case, so considering the $2$-vector space $\Gamma_{mult}(\Omega)$, there is an associated $2$-term complex of vector spaces
\begin{eqnarray*}
\xymatrix{
\eth:\Gamma_{mult}(\Omega)_1 \ar[r] & \Gamma_{mult}(\Omega)_0 ,
}
\end{eqnarray*} 
where $\Gamma_{mult}(\Omega)_1$ is defined to be the kernel of the source map and the map $\eth$ is defined to be the restriction of the target. Notice that the zero element in $\Gamma_{mult}(\Omega)_0$ is the zero section, and therefore, a natural transformation $\tau\in\Gamma_{mult}(\Omega)_1$ has the property that $\hat{s}(\tau(x))=0$. Since we already saw that $\tau(x)\in\Omega_{u(x)}$, we conclude that 
\begin{eqnarray*}
\Gamma_{mult}(\Omega)_1\subseteq\Gamma(\CC{\Omega}).
\end{eqnarray*}
Notice that because of the naturallity, $\eth\tau\in\Gamma_{mult}(\Omega)_0$ is given by the formula
\begin{eqnarray*}
\eth\tau(g)=\tau(t(g))\cdot 0_g\cdot\hat{\iota}(\tau(s(g))).
\end{eqnarray*}
Conversely, given a section $\sigma\in\Gamma(\CC{\Omega})$, we can think of it as being a natural transformation from the zero section to $\omega_\sigma\in\Gamma_{mult}(\Omega)$, where
\begin{eqnarray}
\omega_\sigma(g)=\sigma(t(g))\cdot 0_g\cdot\hat{\iota}(\sigma(s(g))).
\end{eqnarray}
Thus defined, this section is obviously multiplicative and $\sigma\in\Gamma_{mult}(\Omega)(0,\omega_\sigma)$.
\begin{remark}
This construction has two parallels in the literature \cite{CristJames,EugLer}. In both references, the map that appears is the difference of the right-translation and the left-translation
\begin{eqnarray*}
\eth\sigma(g) :=\overrightarrow{\sigma}(g)-\overleftarrow{\sigma}(g)=\sigma(t(g))\cdot 0_g +0_g\cdot\iota(\sigma(s(g))).
\end{eqnarray*}
This formula coincides with the one above, as the following computation shows
\begin{align*}
\sigma(t(g))\cdot 0_g +0_g\cdot\iota(\sigma(s(g))) & =(\sigma(t(g))\cdot 0_g)\cdot\hat{u}(0_{s(g)})+(\hat{u}(0_{t(g)})\cdot 0_g)\cdot\iota(\sigma(s(g))) \\
											       & =\Big{(}(\sigma(t(g))\cdot 0_g)+(\hat{u}(0_{t(g)})0_g)\Big{)}\cdot\Big{(}0_{u(s(g))}+\iota(\sigma(s(g)))\Big{)} \\
											       & =\Big{(}(\sigma(t(g))+0_{u(t(g))})\cdot(0_g+0_g)\Big{)}\cdot\iota(\sigma(s(g))) \\
											       & =\sigma(t(g))\cdot 0_g\cdot\iota(\sigma(s(g))).
\end{align*}
\end{remark}
What is remarkable about this correspondence is that, of course the core of $\Omega$ comes with a Lie algebroid structure, and with a bit of effort, one can show that the fact that the structural maps of $\Omega$ are Lie algebroid maps, translates into the fact that $\Gamma_{mult}(\Omega)_0$ is a Lie sub-algebra of $\Gamma(\Omega)$. 
\begin{theorem}\cite{CristJames}
$\Gamma_{mult}(\Omega)$ is a Lie $2$-algebra.
\end{theorem}
\begin{proof}
The strategy is to prove that the complex 
\begin{eqnarray*}
\xymatrix{
\eth:\Gamma_{mult}(\Omega)_1 \ar[r] & \Gamma_{mult}(\Omega)_0 ,
}
\end{eqnarray*}
can be endowed with the structure of a crossed module of Lie algebras. As outlined before the statement, $\Gamma_{mult}(\Omega)_1\cong\Gamma(\CC{\Omega})$; therefore, it is indeed a Lie algebra. Let us prove that $\Gamma_{mult}(\Omega)_0$ is closed under the bracket of $\Gamma(\Omega)$. As claimed, this will follow from the fact that $\hat{u}$ and $\hat{m}$ are Lie algebroid morphisms. Consider 
\begin{eqnarray*}
\xymatrix{
\Omega \ar@<0.5ex>[r]\ar@<-0.5ex>[r]             & A  \\
G \ar@<0.5ex>[r]\ar@<-0.5ex>[r]\ar[u]^{\omega_k} & M \ar[u]_{\alpha_k}
}
\end{eqnarray*} 
for $k\in\lbrace 1,2\rbrace$, two elements in $\Gamma_{mult}(\Omega)_0$. Let $\lbrace\epsilon_j\rbrace$ be a generating set for $\Gamma(\Omega)$, then using the Einstein summation convention,
\begin{eqnarray*}
\omega_k =F_k^j\epsilon_j
\end{eqnarray*}
for some $F_k^j\in C^\infty(G)$. With this decomposition,
\begin{eqnarray*}
[\omega_1 ,\omega_2]_{\Omega}=\Lie_{\rho_\Omega(\omega_1)}(F_2^k)\epsilon_k -\Lie_{\rho_\Omega(\omega_2)}(F_1^j)\epsilon_j +F_1^jF_2^k[\epsilon_j ,\epsilon_k]_{\Omega}.
\end{eqnarray*}
Since the sections are multiplicative, we have got
\begin{eqnarray*}
\hat{u}(\alpha_k(x))=\omega_k(u(x))=F_k^j(u(x))\epsilon_j(u(x))
\end{eqnarray*}
for all elements $x\in M$. Now, $\hat{u}$ is a Lie algebroid morphism; hence,
\begin{align*}
\hat{u}([\alpha_1 ,\alpha_2]_A(x))=\Lie_{\rho_A(\alpha_1)}(F_2^k\circ u)(x)\epsilon_k(u(x))- & \Lie_{\rho_A(\alpha_2)}(F_1^j\circ u)(x)\epsilon_j(u(x)) + \\
 &  + F_1^j(u(x))F_2^k(u(x))[\epsilon_j ,\epsilon_k]_{\Omega}(u(x)).
\end{align*}
Computing, 
\begin{align*}
\Lie_{\rho_A(\alpha_1)}(F_2^k\circ u)(x) & =d_x(F_2^k\circ u)(\rho_A(\alpha_1(x))) \\
										 & =d_{u(x)}F_2^k(d_xu\circ\rho_A(\alpha_1(x)) \\
										 & =d_{u(x)}F_2^k(\rho_\Omega\circ\hat{u}(\alpha_1(x)) \\
										 & =d_{u(x)}F_2^k(\rho_\Omega(\omega_1(u(x))) = \Lie_{\rho_\Omega(\omega_1)}(F_2^k)(u(x));
\end{align*}
thus proving that $[\omega_1 ,\omega_2]_{\Omega}$ respects units. As for the multiplication, we first need to see that for every pair of composable arrows $(g_1,g_2)\in G^{(2)}$,
\begin{eqnarray*}
\hat{s}([\omega_1 ,\omega_2]_{\Omega}(g_1))=\hat{t}([\omega_1 ,\omega_2]_{\Omega}(g_2)).
\end{eqnarray*}
Indeed, proceeding analogously to the previous part, let $\lbrace\zeta_j\rbrace$ be a generating set for $\Gamma(A)$, then using the Einstein summation convention,
\begin{eqnarray*}
\alpha_k =f_k^j\zeta_j
\end{eqnarray*}
for some $f_k^j\in C^\infty(M)$. With this decomposition,
\begin{eqnarray*}
[\alpha_1 ,\alpha_2]_{A}=\Lie_{\rho_A(\alpha_1)}(f_2^k)\zeta_k -\Lie_{\rho_A(\alpha_2)}(f_1^j)\zeta_j +f_1^jf_2^k[\zeta_j ,\zeta_k]_{A}.
\end{eqnarray*}
Since the sections are multiplicative, we have got
\begin{eqnarray*}
\hat{s}(\omega_k(g))=\alpha_k(s(g))=f_k^j(s(g))\zeta_j(s(g))
\end{eqnarray*}
for all elements $g\in G$. Now, $\hat{s}$ is a Lie algebroid morphism; hence,
\begin{align*}
\hat{s}([\omega_1 ,\omega_2]_\Omega(g))=\Lie_{\rho_\Omega(\omega_1)}(f_2^k\circ s)(g)\zeta_k(s(g))- & \Lie_{\rho_\Omega(\omega_2)}(f_1^j\circ s)(g)\zeta_j(s(g)) + \\
 &  + f_1^j(s(g))f_2^k(s(g))[\zeta_j ,\zeta_k]_{A}(s(g)).
\end{align*}
Computing, 
\begin{align*}
\Lie_{\rho_\Omega(\omega_1)}(f_2^k\circ s)(g) & =d_g(f_2^k\circ s)(\rho_\Omega(\omega_1(g))) \\
										 & =d_{s(g)}f_2^k(d_gs\circ\rho_\Omega(\omega_1(g)) \\
										 & =d_{s(x)}f_2^k(\rho_A\circ\hat{s}(\omega_1(g)) \\
										 & =d_{s(x)}f_2^k(\rho_A(\alpha_1(s(g))) = \Lie_{\rho_A(\alpha_1)}(f_2^k)(s(g));
\end{align*}
thus proving that $[\omega_1 ,\omega_2]_{\Omega}$ respects the source. An identical computation shows that it does respect the target as well. As a consequence, taking the pair $(g_1,g_2)\in G^{(2)}$, we have got
\begin{eqnarray*}
\hat{s}([\omega_1 ,\omega_2]_\Omega(g_1))=[\alpha_1 ,\alpha_2]_A(s(g_1))=[\alpha_1 ,\alpha_2]_A(t(g_2))=\hat{t}([\omega_1 ,\omega_2]_\Omega(g_2))
\end{eqnarray*}
as desired. Naturally, such relation will also be verified by the diagonal sections $(\omega_k,\omega_k)$ making them take values on $\Gamma(\Omega^{(2)})$. By definition, we will have that
\begin{eqnarray*}
[(\omega_1,\omega_1),(\omega_2,\omega_2)]_{\Omega^{(2)}}=([\omega_1 ,\omega_2]_\Omega,[\omega_1 ,\omega_2]_\Omega)\in\Gamma(\Omega^{(2)});
\end{eqnarray*}
thus, using that $\hat{m}$ is a Lie algebroid map, together with 
\begin{eqnarray*}
\hat{m}(\omega_k(g_1);\omega_k(g_2))=\omega_k(g_1g_2)=F_k^j(m(g_1;g_2))\epsilon_j(g_1g_2),
\end{eqnarray*}
\begin{align*}
\hat{m}([\omega_1 ,\omega_2]_\Omega(g_1);[\omega_1 ,\omega_2]_\Omega(g_1)) & =\Lie_{\rho_\Omega^2(\omega_1,\omega_1)}(F_2^k\circ m)(g_1;g_2)\epsilon_k(g_1g_2)+ \\
 & \qquad-\Lie_{\rho_\Omega^2(\omega_2,\omega_2)}(F_1^j\circ m)(g_1;g_2)\epsilon_j(g_1g_2) + \\
 & \qquad\qquad + F_1^j(g_1g_2)F_2^k(g_1g_2)[\epsilon_j,\epsilon_k]_{\Omega}(g_1g_2).
\end{align*}
Computing one last time, 
\begin{align*}
\Lie_{\rho_\Omega^2(\omega_1,\omega_1)}(F_2^k\circ m)(g_1;g_2) & =d_{(g_1;g_2)}(F_2^k\circ m)(\rho_\Omega^2(\omega_1(g_1),\omega_1(g_2))) \\
										 & =d_{g_1g_2}F_2^k(d_{(g_1;g_2)}m\circ\rho_\Omega^2(\omega_1(g_1),\omega_1(g_2))) \\
										 & =d_{g_1g_2}F_2^k(\rho_\Omega\circ\hat{m}(\omega_1(g_1),\omega_1(g_2))) \\
										 & =d_{g_1g_2}F_2^k(\rho_\Omega(\omega_1(g_1g_2)) = \Lie_{\rho_\Omega(\omega_1)}(F_2^k)(g_1g_2);
\end{align*}
ultimately proving that $[\omega_1 ,\omega_2]_{\Omega}$ respects the multiplication, and is thus a multiplicative section. \\
Next, we prove that $\eth$ is a map of Lie algebras. This follows easily from the remark preceding the statement of the theorem; indeed, we defined the bracket in $\Gamma(\CC{\Omega})$ by embedding it into $\Gamma(\Omega)$ via $\overrightarrow{(-)}$. The rest of the proof is rather involved. One needs to define the action, prove it is by derivations and, finally, that $\eth$ is equivariant and the infinitesimal Peiffer equation is satisfied.

\end{proof}
We give some examples of these Lie $2$-algebras of multiplicative sections.\\
\begin{ex:}\label{2AlgAsSections}
The Lie $2$-algebra of multiplicative sections of a Lie $2$-algebra $\xymatrix{\gg_1 \ar@<0.5ex>[r]\ar@<-0.5ex>[r] & \hh}$ is $\gg_1$ itself. Indeed, functors from the category with one object and one arrow are points on the base with their respective identities. Therefore, $\Gamma_{mult}(\gg_1)_0\cong\hh$. On the other hand, the core of a Lie $2$-algebra is a Lie algebra. Let $\gg=\ker\hat{s}$, since the Lie algebra of sections of a Lie algebra seen as a Lie algebroid is also the Lie algebra itself, $\Gamma(\gg)=\gg$. Thus, it is not hard to see that the crossed module of the Lie $2$-algebra of multiplicative sections is $\xymatrix{\gg \ar[r]^\mu & \hh}$, which is the same crossed module associated to $\gg_1$.
\end{ex:}
\begin{ex:}\label{MultVectorFields}
The multiplicative sections of the tangent prolongation $TG$ has been studied in \cite{EugLer,hep2}. Its associated crossed module is
\begin{eqnarray*}
\xymatrix{
\Gamma(A_G) \ar[r] & \XX_{mult}(G):a \ar@{|->}[r] & \overrightarrow{a}-\overleftarrow{a}
}
\end{eqnarray*}
with action $\Lie_Xa:=[X,\overrightarrow{a}]\rest{M}$.
\end{ex:}
\begin{ex:}\label{multBisVBAlgbd}
A functor between flat abelian groupoids is but a map of vector bundles; therefore, if $\Omega$ is a VB-algebroid over the vector bundle $E$, the space of multiplicative sections is the space of linear sections of \cite{DoubleVB}, $\Gamma_{mult}(\Omega)_0=\Gamma_{l}(\Omega,E)$. It is remarked in the reference that $\Gamma_l(\Omega,E)$ is a locally free $C^\infty(M)$ module of rank $\rk(A)+\rk(E)\rk(C)$, where $C$ is the core of $\Omega$; thus, $\Gamma_l(\Omega,E)\cong\Gamma(\hat{\Omega})$ for some vector bundle $\hat{\Omega}$ over $M$ fitting in the exact sequence 
\begin{eqnarray*}
\xymatrix{
(0) \ar[r] & Hom(E,C) \ar[r] & \hat{\Omega} \ar[r] & A \ar[r] & (0) .
}
\end{eqnarray*}
On the other hand, the core of a VB-algebroid is always a vector bundle. The crossed module of the Lie $2$-algebra of multiplicative sections has zero map, and the action is the one induced by the bracket of $\Omega$. 
\end{ex:}
\begin{ex:}\label{MultSectOfPoisson-LieGp}
Using the left trivialization, one sees that the cotangent groupoid $\xymatrix{T^*G \ar@<0.5ex>[r] \ar@<-0.5ex>[r] & \gg^*}$ is isomorphic to the action Lie groupoid $\xymatrix{G\ltimes\gg^* \ar@<0.5ex>[r] \ar@<-0.5ex>[r] & \gg^*}$ associated to the co-adjoint representation
\begin{eqnarray*}
Ad^*_g(\xi):=(Ad_{g^{-1}})^*(\xi)=(dR_g)^*\circ(dL_{g^{-1}})(\xi).
\end{eqnarray*} 
Then a section is a map $\xymatrix{\omega:G \ar[r] & G\ltimes\gg^*}$ given by $\omega(g)=(g,\varphi(g))$, and, to be multiplicative, it needs to respect source and target: $\varphi(g)=\varphi(1)$ and $Ad^*_g\varphi(g)=\varphi(1)$. Thus, the multiplicative sections of the \LA -groupoid of a Poisson-Lie group consists of the Lie sub-algebra of $Ad^*$-invariant elements of $\gg^*$. Since we already saw that the core of these is trivial, the Lie $2$-algebra of multiplicative sections is $\xymatrix{0 \ar[r] & (\gg^*)^G}$. \\
Note that the space of $Ad^*$-invariant co-vectors is a subspace of the annihilator of the derived Lie algebra, $[\gg,\gg]^\circ$. Indeed, for $\xi\in(\gg^*)^G$, $x,y\in\gg$ and $\lambda\in\Rr$,
\begin{align*}
    Ad^*_{\exp(\lambda x)}\xi(y)=\xi(Ad_{\exp(-\lambda x)}y)=\xi(y).
\end{align*}
Differentiating with respect to $\lambda$ and evaluating at zero, one gets $\xi([y,x])=0$, and $x,y$ were arbitrary. The converse also holds if $G$ is connected.
\end{ex:}
\begin{ex:}\label{MSPoisson-LieGpd}
For a Poisson-Lie groupoid $(G,\pi)$ over $M$, the space of multiplicative sections is by definition the space of multiplicative $1$-forms. The core of a Poisson-Lie groupoid is the core of the cotangent groupoid, which is the cotangent space to the base $M$. Since the base inherits a unique Poisson structure making the source a Poisson map, the core inherits a Lie algebroid structure. The crossed module associated to the Lie $2$-algebra of multiplicative sections is 
\begin{eqnarray*}
\xymatrix{
\Omega^1(M) \ar[r] & \Omega_{mult}^1(G):\theta \ar@{|->}[r] & t^*\theta-s^*\theta
}
\end{eqnarray*}
with action $\Lie_\omega\theta:=[\omega,t^*\theta]_{pi}\rest{M}$.
\end{ex:}
\begin{ex:}\label{actionLASections}
The space of multiplicative sections of an action \LA -groupoid of $G$ on $A$, is $\Gamma(A)^G$. Indeed, every multiplicative section $\omega\in\Gamma_{mult}(G\ltimes A)_0$ has components $\omega_1(g,x)=(g,\omega_A(g,x))$ and has to be compatible with the source and the target, then
\begin{eqnarray*}
\hat{s}(\omega_1(g,x))=\omega_A(g,x)=\omega_0(x) & \textnormal{and} & \hat{t}(\omega_1(g,x))=g\omega_0(x)=\omega_0(g\cdot x).
\end{eqnarray*}
Therefore, $\omega$ is specified by its base component $\omega_0$ and the fact it is equivariant. Since we saw that the core of an action \LA -groupoids is trivial, the crossed module of its Lie $2$-algebra is $\xymatrix{(0) \ar[r] & \Gamma(A)^G}$.
\end{ex:}
We close this subsection with the following remarks. In sight of the discussion at the beginning of the section, one can certainly ask how the integrability of this Lie $2$-algebra of multiplicative sections of an \LA -groupoid is related to the integrability of the \LA -groupoid itself. We will see shortly that the integrability of the \LA -groupoid implies the integrability of its Lie $2$-algebra of multiplicative bisections. The converse, however, does not hold in general as shown by example e.g. \ref{MultSectOfPoisson-LieGp}. If it were the case that under some additional hypothesis the integrability of the Lie $2$-algebra of multiplicative bisections implied the integrability of the \LA -groupoid, its associated crossed module would make patent the fact that the integrability of the core plays a fundamental r\^{o}le. We will now make the relation clearer, by exploring the global counter-part of the category of multiplicative sections. %However, before doing so, we would like to point out that there is an integrability theorem for infinite dimensional Lie algebras (see e.g. \cite{IntInfDim}), in which the integration fails to be a group but is rather a $2$-group. Integrating the crossed module with this procedure, one gets a crossed module of Lie $2$-groups and a corresponding %**adjective saying that it is special** double Lie groupoid. How is this related to our integration problem? \\
%*** NOTE: It would be nice to see whether the van Est theory developed in the second chapter can be applied to this construction. For example (or specially) in the case of VB-algebroids. Maybe Alejandro and Thiago's van Est theorem is what is needed to conclude the integration of those using the techniques here. ***
\subsubsection{The category of multiplicative bisections}
We define the global counter-parts of the multiplicative sections of an \LA -groupoid; namely, the space of multiplicative bisections of a double Lie groupoid. \\
Let $D$ be a double Lie groupoid. We define the category of multiplicative bisections as the space of groupoid morphisms $\xymatrix{H \ar[r] & D}$ that are also bisections at both the level of objects and of morphisms. Since such maps are functors as well, we consider natural transformations between them as morphisms. This time around, the compatibility condition will also reflect the defining condition for a bisection. \\
Let 
\begin{eqnarray*}
\Bis_{mult}(D)_0:=\lbrace\sigma\in\Bis(D):\sigma\textnormal{ Lie groupoid homomorphism}\rbrace .
\end{eqnarray*}
The whole compatibility of the the structural maps of the double Lie groupoid will imply that the multiplicative bisections form a subgroup, as the following lemma shows.
\begin{lemma}
$\Bis_{mult}(D)_0$ is a subgroup of $\Bis(D)$.
\end{lemma}  
\begin{proof}
First, notice that $\Lf{u}$ is indeed a bisection and a Lie groupoid homomorphism by definition. Now, let 
\begin{eqnarray*}
\xymatrix{
D \ar@<0.5ex>[r]\ar@<-0.5ex>[r]                    & V \\
H \ar@<0.5ex>[r]\ar@<-0.5ex>[r]\ar[u]^{\sigma_1^k} & M \ar[u]_{\sigma_0^k}
}
\end{eqnarray*} 
for $k\in\lbrace 1,2\rbrace$, be two elements in $\Bis_{mult}(D)_0$. We will prove that $\sigma_1^1\sigma_1^2\in\Bis_{mult}(D)_0$. For $x\in M$, we've got that
\begin{eqnarray*}
\sigma_1^k(u(x))=\Tp{u}(\sigma_0^k(x));
\end{eqnarray*}
hence, using that $\Tp{u}$ is a Lie groupoid homomorphism,
\begin{align*}
(\sigma_1^1\sigma_1^2)(u(x)) & =\sigma_1^1(\Lf{t}(\sigma_1^2(u(x))))\vJoin \sigma_1^2(u(x)) \\
							 & =\sigma_1^1(\Lf{t}(\Tp{u}(\sigma_0^2(x))))\vJoin \Tp{u}(\sigma_0^2(x)) \\
							 & =\sigma_1^1(u(\Rg{t}(\sigma_0^2(x))))\vJoin \Tp{u}(\sigma_0^2(x)) \\
							 & =\Tp{u}(\sigma_0^1(\Rg{t}(\sigma_0^2(x))))\vJoin \Tp{u}(\sigma_0^2(x)) \\
							 & =\Tp{u}(\sigma_0^1(\Rg{t}(\sigma_0^2(x)))\sigma_0^2(x))=\Tp{u}((\sigma_0^1\sigma_0^2)(x)).
\end{align*}
On the other hand, let $h\in H$, then
\begin{eqnarray*}
\Tp{s}\sigma_1^k(h)=\sigma_0^k(s(h)),
\end{eqnarray*}
and using that $\Tp{s}$ is a Lie groupoid homomorphism, conclude
\begin{align*}
\Tp{s}(\sigma_1^1\sigma_1^2)(h) & =\Tp{s}(\sigma_1^1(\Lf{t}(\sigma_1^2(h)))\vJoin \sigma_1^2(h)) \\
							 & =\Tp{s}\sigma_1^1(\Lf{t}(\sigma_1^2(h)))\Tp{s}\sigma_1^2(h) \\
							 & =\sigma_0^1(s(\Lf{t}(\sigma_1^2(h))))\sigma_0^2(s(h)) \\
							 & =\sigma_0^1(\Rg{t}(\Tp{s}(\sigma_1^2(h))))\sigma_0^2(s(h)) \\
							 & =\sigma_0^1(\Rg{t}(\sigma_0^2(s(h))))\sigma_0^2(s(h))=(\sigma_0^1\sigma_0^2)(s(h)).
\end{align*}
An identical computation shows that $\Tp{t}(\sigma_1^1\sigma_1^2)(h)=(\sigma_0^1\sigma_0^2)(t(h))$. Now, for $(h_1,h_2)\in H\times_MH$,\begin{eqnarray*}
\Tp{s}(\sigma_1^1\sigma_1^2)(h_1)=(\sigma_1^1\sigma_1^2)(s(h_1))=(\sigma_1^1\sigma_1^2)(t(h_2))=\Tp{t}(\sigma_1^1\sigma_1^2)(h_2).
\end{eqnarray*}
Therefore, form the compatibility 
\begin{eqnarray*}
\sigma_1^k(h_1)\Join\sigma_1^k(h_2)=\sigma_1^k(h_1h_2)
\end{eqnarray*}
together with the interchange law and the fact that $\Lf{t}$ is a Lie groupoid homomorphism, we get
\begin{align*}
(\sigma_1^1\sigma_1^2)(h_1)\Join(\sigma_1^1\sigma_1^2)(h_2) & =(\sigma_1^1(\Lf{t}\sigma_1^2(h_1))\vJoin\sigma_1^2(h_1))\Join(\sigma_1^1(\Lf{t}\sigma_1^2(h_2))\vJoin\sigma_1^2(h_2)) \\
			 & =(\sigma_1^1(\Lf{t}\sigma_1^2(h_1))\Join\sigma_1^1(\Lf{t}\sigma_1^2(h_2)))\vJoin(\sigma_1^2(h_1)\Join\sigma_1^2(h_2)) \\
			 & =\sigma_1^1(\Lf{t}\sigma_1^2(h_1)\Lf{t}\sigma_1^2(h_2))\vJoin\sigma_1^2(h_1h_2) \\
			 & =\sigma_1^1(\Lf{t}(\sigma_1^2(h_1)\Join\sigma_1^2(h_2)))\vJoin\sigma_1^2(h_1h_2) \\
			 & =\sigma_1^1(\Lf{t}(\sigma_1^2(h_1h_2)))\vJoin\sigma_1^2(h_1h_2)=(\sigma_1^1\sigma_1^2)(h_1h_2).
\end{align*}
Finally, we prove that $\Bis_{mult}(D)_0$ is also closed under inverses. Recall that 
\begin{eqnarray*}
\sigma^{\dagger}:=\Lf{\iota}\circ\sigma\circ(\Lf{t}\circ\sigma)^{-1},
\end{eqnarray*}
This will follow from the fact that the Lie groupoid morphism $\Lf{t}\circ\sigma$ is not only a diffeomorphism, but a Lie groupoid automorphism. To prove this assertion, consider
\begin{align*}
\Lf{t}\circ\sigma_1(u(x)) & =\Lf{t}\circ\Tp{u}(\sigma_0(x))=u(\Rg{t}\circ\sigma_0(x)),
\end{align*}
evaluating at $x=(\Rg{t}\circ\sigma_0)^{-1}(y)$, and taking the inverse map 
\begin{eqnarray*}
u((\Rg{t}\circ\sigma_0)^{-1}(y))=(\Lf{t}\circ\sigma_1)^{-1}u(y).
\end{eqnarray*}
Analogously, 
\begin{align*}
s(\Lf{t}\circ\sigma_1(h)) & =\Rg{t}(\Tp{s}\circ\sigma_1(h))=\Rg{t}\circ\sigma_0(s(h))
\end{align*}
and evaluating at $h=(\Lf{t}\circ\sigma)^{-1}(g)$ and applying the inverse map,
\begin{eqnarray*}
(\Lf{t}\circ\sigma)^{-1}(s(g))=s((\Lf{t}\circ\sigma)^{-1}(g)).
\end{eqnarray*}
Lastly,
\begin{align*}
\Lf{t}\circ\sigma_1(h_1)\Lf{t}\circ\sigma_1(h_2) & =\Lf{t}(\sigma_1(h_1)\Join\sigma_1(h_2))=\Lf{t}circ\sigma_1(h_1h_2);
\end{align*}
thus, evaluating at $h_k=(\Lf{t}\circ\sigma)^{-1}(g_k)$ and taking the inverse one last time
\begin{eqnarray*}
(\Lf{t}\circ\sigma)^{-1}(g_1g_2)=(\Lf{t}\circ\sigma)^{-1}(g_1)(\Lf{t}\circ\sigma)^{-1}(g_2).
\end{eqnarray*}
\end{proof}
Given $\sigma^1,\sigma^2\in\Bis_{mult}(D)_0$, we define 
\begin{eqnarray*}
\Bis_{mult}(D)(\sigma^1,\sigma^2):= \lbrace\tau\in Nat^\infty(\sigma^1,\sigma^2):\tau\ast_{h}1_{\Lf{s}} =1_{id_H}\rbrace .
\end{eqnarray*}
The compatibility equation that appear this time around is analogous to the one of multiplicative sections. Indeed, diagrammatically, 
%\begin{eqnarray*}
%\xymatrix{
%H \ar@/^{1pc}/[rr]^{\sigma^1}="a" \ar@/_{1pc}/[rr]_{\sigma^2}="b" & \ar@{=>}^{\tau}"a";"b" & D \ar@/^{1pc}/[rr]^{\Lf{s}}="c" \ar@/_{1pc}/[rr]_{\Lf{s}}="d" & \ar@{=>}^{1_{\Lf{s}}}"c";"d" & H
%} = \xymatrix{
%H \ar@/^{1pc}/[rr]^{id_H}="e"\ar@/_{1pc}/[rr]_{id_H}="f" & \ar@{=>}^{1_{id_H}}"e";"f" & H .
%}
%\end{eqnarray*}
\begin{align*}
\xymatrix{
H \ar@/^{1pc}/[rr]^{\sigma_1}\ar@/_{1pc}/[rr]_{\sigma_2} & \quad{\Big{\Downarrow}}\small{\tau} & D \ar@/^{1pc}/[rr]^{\Lf{s}}\ar@/_{1pc}/[rr]_{\Lf{s}} & \quad\Big{\Downarrow} \small{1_{\Lf{s}}} & H
} & =\xymatrix{
H \ar@/^{1pc}/[rr]^{id_H}\ar@/_{1pc}/[rr]_{id_H} & \qquad\Big{\Downarrow} \small{1_{id_H}} & H .
}
\end{align*}
It would be reasonable to assume that such a $\tau$ should also take into account the induced functors $\Lf{t}\circ\sigma^k$; however, the only canonically induced natural transformation between them is $\Lf{t}\circ\tau$, which is already $\tau\ast_{h}1_{\Lf{t}}$. Similar considerations to those succeeding the definition of the category of multiplicative sections apply, thus settling that $\Bis_{mult}(D)$ is a category as well. \\
Notice that the condition that $\tau\ast_{h}1_{\Lf{s}} =1_{id_H}$ means that for $x\in M$, $\tau(x)$ is a square in $D$ of the form
\begin{eqnarray*}
\xymatrix{
\bullet  & \bullet \ar[l] \\
{}_x\bullet \ar[u]^{\sigma_0^2(x)} & \bullet_x \ar@{=}[l]\ar[u]_{\sigma_0^1(x)} 
}
\end{eqnarray*}
Therefore, analogous to the infinitesimal case, an element $\tau\in\Bis_{mult}(\Lf{u},\sigma)$ defines a bisection of $\CC{D}$ with induced diffeomorphism given by $\Rg{t}\circ\sigma$. Conversely, given $\tau\in\Bis(\CC{D})$, one can think ot it as being an element of $\Bis_{mult}(D)(\Lf{u},\Delta\tau)$, where
\begin{eqnarray*}
(\Delta\tau)_0(x):=\Tp{t}(\tau(x)) &\textnormal{and} & (\Delta\tau)_1(h):=\tau(t(h))\Join\Lf{u}(h)\Join\Tp{\iota}(\tau(s(h)))
\end{eqnarray*} 
for $x\in M$ and $h\in H$. Let us verify that this is well-defined. First, the multiplication in $(\Delta\tau)_1$ makes sense, in that
\begin{align*}
\Tp{s}(\tau(t(h))) & =\Rg{u}(t(h))       & \Tp{s}(\Lf{u}(h)) & =\Rg{u}(s(h))      \\
				   & =\Tp{t}(\Lf{u}(h)), &		             & =\Tp{s}(\tau(s(h)))=\Tp{t}\circ\Tp{\iota}(\tau(s(h))).
\end{align*}
Now, let us see that $\Delta\tau$ is a multiplicative bisection. For starters, it is indeed a bisection:
\begin{align*}
\Lf{s}((\Delta\tau)_1(h)) & =\Lf{s}(\tau(t(h))\Join\Lf{u}(h)\Join\Tp{\iota}(\tau(s(h)))) \\ 
				          & =\Lf{s}(\tau(t(h)))\Lf{s}(\Lf{u}(h))\Lf{s}(\Tp{\iota}(\tau(s(h)))) \\
				          & =u(t(h))h\iota(u(s(h)))=h,
\end{align*}				          
and
\begin{align*}				          
\Lf{t}((\Delta\tau)_1(h)) & =\Lf{t}(\tau(t(h))\Join\Lf{u}(h)\Join\Tp{\iota}(\tau(s(h)))) \\ 
				          & =\Lf{t}(\tau(t(h)))\Lf{t}(\Lf{u}(h))\Lf{t}(\Tp{\iota}(\tau(s(h)))) \\
				          & =\Lf{t}(\tau(t(h)))h\iota(\Lf{t}(\tau(t(h)))).
\end{align*}
This last expression defines an a diffeomorphism whose inverse is given by
\begin{eqnarray*}
(\Tp{t}\circ(\Delta\tau)_1)^{-1}(h):=\iota(\Lf{t}(\tau(t(h))))h\Lf{t}(\tau(t(h))).
\end{eqnarray*}
On the other hand, $\Delta\tau$ is multiplicative. For $x\in M$
\begin{eqnarray*}
(\Delta\tau)_1(u(x)):=\tau(x)\Join\Lf{u}(u(x))\Join\Tp{\iota}(\tau(x))=\Tp{u}(\Tp{t}(\tau(x)))=\Tp{u}((\Delta\tau)_0(x)).
\end{eqnarray*}
For $h\in H$,
\begin{eqnarray*}
\Tp{s}((\Delta\tau)_1(h))=\Tp{s}(\Tp{\iota}(\tau(s(h))))=\Tp{t}(\tau(s(h)))=(\Delta\tau)_0(s(h))
\end{eqnarray*}
and also, immediately,
\begin{eqnarray*}
\Tp{t}((\Delta\tau)_1(h))=\Tp{t}(\tau(t(h)))=(\Delta\tau)_0(t(h)).
\end{eqnarray*}
Finally, for $(h_1,h_2)\in H\times_MH$, because of the latter relations
\begin{eqnarray*}
\Tp{s}((\Delta\tau)_1(h_1))=(\Delta\tau)_0(s(h_1))=(\Delta\tau)_0(t(h_2))=\Tp{t}((\Delta\tau)_1(h_2));
\end{eqnarray*}
therefore, it makes sense to compute
\begin{align*}
(\Delta\tau)_1(h_1)\Join(\Delta\tau)_1(h_2) & =\tau(t(h_1))\Join\Lf{u}(h_1)\Join\Tp{\iota}(\tau(s(h_1)))\Join\tau(t(h_2))\Join\Lf{u}(h_2)\Join\Tp{\iota}(\tau(s(h_2))) \\
											& =\tau(t(h_1))\Join\Lf{u}(h_1)\Join\Lf{u}(h_2)\Join\Tp{\iota}(\tau(s(h_2))) \\
											& =\tau(t(h_1h_2))\Join\Lf{u}(h_1h_2)\Join\Tp{\iota}(\tau(s(h_1h_2)))=(\Delta\tau)_1(h_1h_2).
\end{align*}
It would be left to verify that $\tau$ is indeed an natural transformation verifying the defining equation, but this follows immediately from the definition. Now that we have got this identification, we know that there is a group structure on a subspace of $\Bis_{mult}(D)$. We will ultimately prove that this is no accident, as this object has the structure of a $2$-group. We will follow the same strategy as before, and prove that 
\begin{eqnarray*}
\xymatrix{
\Bis(\CC{D}) \ar[r]^\Delta & \Bis_{mult}(D)_0
}
\end{eqnarray*}
can be endowed with the structure of a crossed module of groups.
\begin{lemma}\label{DeltaHomo}
$\Delta$ is a group homomorphism
\end{lemma}
\begin{proof}
First, recall that the unit map of $\CC{D}$ defines the unit in $\Bis(\CC{D})$. Then,
\begin{eqnarray*}
(\Delta u_{\CC{D}})_0(x):=\Tp{t}(u^2(x))=\Rg{u}(x) &\textnormal{and} & (\Delta u_{\CC{D}})_1(h):=u^2(t(h))\Join\Lf{u}(h)\Join\Tp{\iota}(u^2(s(h)))=\Lf{u}(h).
\end{eqnarray*} 
Now, let $\tau_1,\tau_2\in\Bis(\CC{D})$ and $h\in H(x,y)$. We will call the isomorphisms induced by the target
\begin{eqnarray*}
\overline{\tau}_k:=t_{\CC{D}}\circ\tau_k & \textnormal{and} & \overline{\Delta\tau}_k:=\Lf{t}\circ(\Delta\tau_k)_1.
\end{eqnarray*}
Notice further that $\Rg{t}\circ(\Delta\tau_k)_0=t_{\CC{D}}\circ\tau_k$, reason because of which we need not add more notation. Computing, 
\begin{align*}
(\Delta\tau_1\tau_2)_1(h) & =(\tau_1\tau_2)(y)\Join\Lf{u}(h)\Join\Tp{\iota}((\tau_1\tau_2)(x)) \\
 & =\Big{(}\tau_1(\overline{\tau}_2(y))\boxplus\tau_2(y)\Big{)}\Join\Lf{u}(h)\Join\Tp{\iota}\Big{(}\tau_1(\overline{\tau}_2(x))\boxplus\tau_2(x)\Big{)} \\
 & =\Big{(}\tau_1(\overline{\tau}_2(y))\vJoin\Tp{u}((\Delta\tau_2)_0(y))\Big{)}\Join\tau_2(y)\Join\Lf{u}(h)\Join\Tp{\iota}\Big{[}\Big{(}\tau_1(\overline{\tau}_2(x))\vJoin\Tp{u}((\Delta\tau_2)_0(x))\Big{)}\Join\tau_2(x)\Big{]} \\
 & =\Big{(}\tau_1(\overline{\tau}_2(y))\vJoin\Tp{u}((\Delta\tau_2)_0(y))\Big{)}\Join(\Delta\tau_2)_1(h)\Join\Tp{\iota}\Big{(}\tau_1(\overline{\tau}_2(x))\vJoin\Tp{u}((\Delta\tau_2)_0(x))\Big{)} \\
 & =\Big{(}\tau_1(\overline{\tau}_2(y))\vJoin\Tp{u}((\Delta\tau_2)_0(y))\Big{)}\Join\Big{(}\Lf{u}(\overline{\Delta\tau}_2(h))\vJoin(\Delta\tau_2)_1(h)\Big{)}\Join\Big{(}\Tp{\iota}(\tau_1(\overline{\tau}_2(x)))\vJoin\Tp{\iota}(\Tp{u}((\Delta\tau_2)_0(x)))\Big{)} \\
 & =\Big{(}\big{(}\tau_1(\overline{\tau}_2(y))\Join\Lf{u}(\overline{\Delta\tau}_2(h))\big{)}\vJoin\big{(}\Tp{u}((\Delta\tau_2)_0(y))\Join(\Delta\tau_2)_1(h)\big{)}\Big{)}\Join\Big{(}\Tp{\iota}(\tau_1(\overline{\tau}_2(x)))\vJoin\Tp{u}((\Delta\tau_2)_0(x))\Big{)} \\
 & =\Big{(}\tau_1(\overline{\tau}_2(y))\Join\Lf{u}(\overline{\Delta\tau}_2(h))\Join\Tp{\iota}(\tau_1(\overline{\tau}_2(x)))\Big{)}\vJoin\Big{(}\Tp{u}((\Delta\tau_2)_0(y))\Join(\Delta\tau_2)_1(h)\Join\Tp{u}((\Delta\tau_2)_0(x))\Big{)} \\
 & =(\Delta\tau_1)_1(\overline{\Delta\tau}_2(h))\vJoin(\Delta\tau_2)_1(h) =((\Delta\tau_1)_1(\Delta\tau_2)_1)(h),
\end{align*}
thus completing the proof.

\end{proof}
We use the notation introduced in the proof of this latter lemma, to state the following.
\begin{prop}\label{acts}
There is a right action of $\Bis_{mult}(D)_0$ on $\Bis(\CC{D})$ given by
\begin{eqnarray*}
\tau^\sigma(x):=\sigma_1^\dagger(\Lf{t}(\tau(\overline{\sigma}_0(x)))\vJoin\tau(\overline{\sigma}_0(x))\vJoin\Tp{u}(\sigma_0(x)).
\end{eqnarray*} 
Moreover, thus defined, the action is by automorphisms.
\end{prop}
\begin{proof}
Maybe, we start by pointing out that the definition looks like a conjugation inside of the group of multiplicative bisections. The caveat being that $\tau$ is not an element of $\Bis_{mult}(D)_0$, prompting us to define it element by element. Now, we will start by proving that the formula above is well-defined. Since
\begin{eqnarray*}
\Lf{s}(\tau(\overline{\sigma}_0(x)))=u(\overline{\sigma}_0(x))=\Lf{t}(\Tp{u}(\sigma_0(x)))
\end{eqnarray*}
and
\begin{align*}
\Lf{s}(\sigma_1^\dagger(\Lf{t}(\tau(\overline{\sigma}_0(x)))) & =\Lf{s}\circ\Lf{\iota}(\sigma_1(\overline{\sigma}_1^{-1}\circ\Lf{t}(\tau(\overline{\sigma}_0(x)))) \\
  & =\Lf{t}\circ\sigma_1(\overline{\sigma}_1^{-1}\circ\Lf{t}(\tau(\overline{\sigma}_0(x))) =\Lf{t}(\tau(\overline{\sigma}_0(x))),
\end{align*}
the multiplication of the definition can be performed. On the other hand,
\begin{align*}
\Tp{s}(\tau^\sigma(x)) & =\Tp{s}(\sigma_1^\dagger(\Lf{t}(\tau(\overline{\sigma}_0(x)))\vJoin\tau(\overline{\sigma}_0(x))\vJoin\Tp{u}(\sigma_0(x))) \\
  & =\Tp{s}(\sigma_1^\dagger(\Lf{t}(\tau(\overline{\sigma}_0(x))))\Tp{s}(\tau(\overline{\sigma}_0(x)))\Tp{s}(\Tp{u}(\sigma_0(x))) \\
  & =\Tp{s}\circ\Lf{\iota}(\sigma_1(\overline{\sigma}_1^{-1}\circ\Lf{t}(\tau(\overline{\sigma}_0(x))))\Rg{u}(\overline{\sigma}_0(x))\sigma_0(x) \\
  & =\Rg{\iota}\circ\Tp{s}(\sigma_1(\overline{\sigma}_1^{-1}\circ\Lf{t}(\tau(\overline{\sigma}_0(x))))\sigma_0(x) \\
  & =\Rg{\iota}\circ\sigma_0(s(\overline{\sigma}_1^{-1}\circ\Lf{t}(\tau(\overline{\sigma}_0(x))))\sigma_0(x) \\
  & =\Rg{\iota}\circ\sigma_0(\overline{\sigma}_1^{-1}(s\circ\Lf{t}(\tau(\overline{\sigma}_0(x))))\sigma_0(x) \\
  & =\Rg{\iota}\circ\sigma_0(\overline{\sigma}_1^{-1}(\overline{\sigma}_0(x)))\sigma_0(x) = \Rg{u}(\Rg{s}(\sigma_0(x)))=\Rg{u}(x)  
\end{align*}
and
\begin{eqnarray*}
\Lf{s}(\tau^\sigma(x))=\Lf{s}\circ\Tp{u}(\sigma_0(x))=u\circ\Rg{s}(\sigma_0(x))=u(x);
\end{eqnarray*}
thus, $\tau^\sigma(x)\in\CC{D}$ and incidentally, due to this second relation it also defines a section of $s_{\CC{D}}$. To see it is a bisection, we compute
\begin{align*}
t_{\CC{D}}(\tau^\sigma(x)) & =t\circ\Lf{t}(\tau^\sigma(x)) \\
  						   & =t\circ\Lf{t}(\sigma_1^\dagger(\Lf{t}(\tau(\overline{\sigma}_0(x)))) \\
						   & =t(\Lf{t}\circ\Lf{\iota}(\sigma_1(\overline{\sigma}_1^{-1}\circ\Lf{t}(\tau(\overline{\sigma}_0(x))))) \\
						   & =t(\Lf{s}(\sigma_1(\overline{\sigma}_1^{-1}\circ\Lf{t}(\tau(\overline{\sigma}_0(x))))) \\
						   & =t(\overline{\sigma}_1^{-1}\circ\Lf{t}(\tau(\overline{\sigma}_0(x)))) \\
						   & =\overline{\sigma}_1^{-1}(t\circ\Lf{t}(\tau(\overline{\sigma}_0(x)))) \\
						   & =\overline{\sigma}_1^{-1}(t_{\CC{D}}\circ\tau(\overline{\sigma}_0(x))))=\overline{\sigma}_1^{-1}\circ\overline{\tau}\circ\overline{\sigma}_0(x)
\end{align*}
and conclude it is a diffeomorphism being a composition of them. \\
We move on to show that the formula does define an action. First,
\begin{eqnarray*}
\tau^{\Lf{u}}(x)=\Lf{u}^\dagger(\Lf{t}(\tau(\overline{\Rg{u}}(x)))\vJoin\tau(\overline{\Rg{u}}(x))\vJoin\Tp{u}(\Rg{u}(x))=\Lf{u}(\Lf{t}(\tau(x))\vJoin\tau(x)=\tau(x).
\end{eqnarray*}
For $\sigma^1,\sigma^2\in\Bis_{mult}(D)$, 
\begin{eqnarray*}
\tau^{\sigma^1\sigma^2}(x)=(\sigma^1\sigma^2)_1^\dagger(\Lf{t}(\tau((\overline{\sigma^1\sigma^2})_0(x))))\vJoin\tau((\overline{\sigma^1\sigma^2})_0(x))\vJoin\Tp{u}((\sigma^1\sigma^2)_0(x)),
\end{eqnarray*}
whereas 
\begin{align*}
(\tau^{\sigma^1})^{\sigma^2}(x)=(\sigma^2)_1^\dagger((\overline{\sigma}_1^1)^{-1}(\tau(\overline{\sigma}^1_0(\overline{\sigma}^2_0( & x)))))\vJoin(\sigma^1)_1^\dagger(\Lf{t}(\tau(\overline{\sigma}^1_0(\overline{\sigma}^2_0(x)))))\vJoin \\
 & \vJoin\tau(\overline{\sigma}^1_0(\overline{\sigma}^2_0(x)))\vJoin\Tp{u}(\sigma^1_0(\overline{\sigma}^2_0(x)))\vJoin\Tp{u}(\sigma^2_0(x)).
\end{align*}
Hence, since 
\begin{align*}
(\overline{\sigma^1\sigma^2})_0(x) & =\Rg{t}((\sigma^1\sigma^2)_0(x)) \\
								   & =\Rg{t}(\sigma^1(\overline{\sigma}^2_0(x))\sigma^2_0(x)) \\
								   & =\Rg{t}(\sigma^1(\overline{\sigma}^2_0(x))=\overline{\sigma}^1_0(\overline{\sigma}^2_0(x))
\end{align*}
and $\Tp{u}(v_1v_2)=\Tp{u}(v_1)\vJoin\Tp{u}(v_2)$ for all $v_1,v_2\in V$, it is left to see that 
\begin{eqnarray*}
(\sigma^1\sigma^2)_1^\dagger(\Lf{t}(\tau(z)))=(\sigma^2)_1^\dagger((\overline{\sigma}_1^1)^{-1}(\Lf{t}(\tau(z))))\vJoin(\sigma^1)_1^\dagger(\Lf{t}(\tau(z))),
\end{eqnarray*}
but $(-)^\dagger$ is a group anti-homomorphism and for all $h\in H$,
\begin{eqnarray*}
\Lf{t}((\sigma^1)_1^\dagger(h))=\Lf{t}(\Lf{\iota}\circ\sigma_1((\overline{\sigma}^1_1)^{-1}(h))),
\end{eqnarray*}
so the desired equality holds. \\
Finally, we give reasons for the second part of the statement. We start again by showing that the identity element is taken to itself; indeed,
\begin{align*}
(u_{\CC{D}})^\sigma(x) & =\sigma_1^\dagger(\Lf{t}(u^2(\overline{\sigma}_0(x)))\vJoin u^2(\overline{\sigma}_0(x))\vJoin\Tp{u}(\sigma_0(x)) \\
				     & =\Lf{\iota}\circ\sigma_1((\overline{\sigma}_1)^{-1}(u(\overline{\sigma}_0(x))))\vJoin\sigma_1(u(x)) \\
				     & =\Lf{\iota}\circ\sigma_1(u((\overline{\sigma}_0)^{-1}(\overline{\sigma}_0(x))))\vJoin\sigma_1(u(x))=\Lf{u}\circ\Lf{s}(\sigma_1(u(x)))=u^2(x).
\end{align*}
Here, the pass to the last line holds because, as we saw, because $(\overline{\sigma}^1_1)^{-1}$ is a Lie groupoid automorphism. Now, let $\tau_1,\tau_2\in\Bis(\CC{D})$ and, to ease notation a bit, write $z:=\overline{\sigma}_0(x)$. Then, 
\begin{align*}
(\tau_1\tau_2)^\sigma(x) & =\sigma_1^\dagger(\Lf{t}((\tau_1\tau_2)(z)))\vJoin(\tau_1\tau_2)(z)\vJoin\Tp{u}(\sigma_0(x)) \\
					     & =\sigma_1^\dagger(\Lf{t}((\tau_1\tau_2)(z)))\vJoin\Big{(}\tau_1(\overline{\tau}_2(z))\boxplus\tau_2(z)\Big{)}\vJoin\Tp{u}(\sigma_0(x)) \\
					     & =\sigma_1^\dagger(\Lf{t}((\tau_1\tau_2)(z)))\vJoin\Big{(}\tau_1(\overline{\tau}_2(z))\Join\Lf{u}\circ\Lf{t}(\tau_2(z))\Big{)}\vJoin\tau_2(z)\vJoin\Tp{u}(\sigma_0(x)) \\
					     & =\sigma_1^\dagger(\Lf{t}((\tau_1\tau_2)(z)))\vJoin\Big{(}\tau_1(\overline{\tau}_2(z))\Join\Lf{u}\circ\Lf{t}(\tau_2(z))\Big{)}\vJoin\Lf{\iota}\big{(}\sigma_1^\dagger(\Lf{t}(\tau_2(z)))\big{)}\vJoin\tau_2^\sigma(x) .
\end{align*}
We consider the first terms separately, 
\begin{align*}
\sigma_1^\dagger & (\Lf{t}((\tau_1\tau_2)(z)))\vJoin\Big{(}\tau_1(\overline{\tau}_2(z))\Join\Lf{u}\circ\Lf{t}(\tau_2(z))\Big{)}\vJoin\Lf{\iota}\big{(}\sigma_1^\dagger(\Lf{t}(\tau_2(z)))\big{)} \\
  & =\sigma_1^\dagger(\Lf{t}(\tau_1(\overline{\tau}_2(z))\boxplus\tau_2(z)))\vJoin\Big{(}\tau_1(\overline{\tau}_2(z))\Join\Lf{u}\circ\Lf{t}(\tau_2(z))\Big{)}\vJoin\sigma_1(\overline{\sigma}^{-1}_1(\Lf{t}(\tau_2(z))) \\
  & =\sigma_1^\dagger\big{(}\Lf{t}(\tau_1(\overline{\tau}_2(z)))\Lf{t}(\tau_2(z))\big{)}\vJoin\Big{(}\tau_1(\overline{\tau}_2(z))\Join\Lf{u}\circ\Lf{t}(\tau_2(z))\Big{)}\vJoin\sigma_1(\overline{\sigma}^{-1}_1(\Lf{t}(\tau_2(z))) \\
  & =\Big{(}\sigma_1^\dagger(\Lf{t}(\tau_1(\overline{\tau}_2(z))))\Join\sigma_1^\dagger(\Lf{t}(\tau_2(z))\Big{)}\vJoin\Big{(}\tau_1(\overline{\tau}_2(z))\Join\Lf{u}\circ\Lf{t}(\tau_2(z))\Big{)}\vJoin\sigma_1(\overline{\sigma}^{-1}_1(\Lf{t}(\tau_2(z))) \\
  & =\Big{(}\big{(}\sigma_1^\dagger(\Lf{t}(\tau_1(\overline{\tau}_2(z))))\vJoin\tau_1(\overline{\tau}_2(z))\big{)}\Join\big{(}\sigma_1^\dagger(\Lf{t}(\tau_2(z))\vJoin\Lf{u}\circ\Lf{t}(\tau_2(z))\big{)}\Big{)}\vJoin\sigma_1(\overline{\sigma}^{-1}_1(\Lf{t}(\tau_2(z))) \\
  & =\Big{(}\big{(}\sigma_1^\dagger(\Lf{t}(\tau_1(\overline{\tau}_2(z))))\vJoin\tau_1(\overline{\tau}_2(z))\big{)}\Join\sigma_1^\dagger(\Lf{t}(\tau_2(z))\Big{)}\vJoin\Big{(}\Tp{u}\circ\Tp{t}(\sigma_1(\overline{\sigma}_1^{-1}(\Lf{t}(\tau_2(z)))\Join\sigma_1(\overline{\sigma}^{-1}_1(\Lf{t}(\tau_2(z)))\Big{)} \\
  & =\Big{(}\sigma_1^\dagger(\Lf{t}(\tau_1(\overline{\tau}_2(z))))\vJoin\tau_1(\overline{\tau}_2(z))\vJoin\Tp{u}(\sigma_0\circ t(\overline{\sigma}_1^{-1}(\Lf{t}(\tau_2(z)))\Big{)}\Join\Big{(}\sigma_1^\dagger(\Lf{t}(\tau_2(z))\vJoin\sigma_1(\overline{\sigma}^{-1}_1(\Lf{t}(\tau_2(z)))\Big{)} \\
  & =\Big{(}\sigma_1^\dagger(\Lf{t}(\tau_1(\overline{\tau}_2(z))))\vJoin\tau_1(\overline{\tau}_2(z))\vJoin\Tp{u}(\sigma_0(\overline{\sigma}_0^{-1}(t^2(\tau_2(z))))\Big{)}\Join\Lf{u}\circ\Lf{s}\Big{(}\sigma_1(\overline{\sigma}^{-1}_1(\Lf{t}(\tau_2(z)))\Big{)} \\
  & =\tau_1^\sigma(\overline{\sigma}_0^{-1}(\overline{\tau}_2(z)))\Join\Lf{u}(\overline{\sigma}^{-1}_1(\Lf{t}(\tau_2(z))) .
\end{align*}
We finish the proof by explicitly showing that
\begin{align*}
\Lf{t}(\tau_2^\sigma(x)) & =\Lf{t}(\sigma_1^\dagger(\Lf{t}(\tau_2(z)))) \\
						 & =\Lf{t}\circ\Lf{\iota}\circ\sigma_1(\overline{\sigma}_1^{-1}(\Lf{t}(\tau_2(z)))) \\
						 & =\overline{\sigma}_1^{-1}(\Lf{t}(\tau_2(z))))
\end{align*}
and consequently
\begin{align*}
t_{\CC{D}}(\tau_2^\sigma(x)) & =t\circ\Lf{t}(\tau_2^\sigma(x)) \\
							 & =t(\overline{\sigma}_1^{-1}(\Lf{t}(\tau_2(z)))) \\
							 & =\overline{\sigma}_0^{-1}(t(\Lf{t}(\tau_2(z))))=\overline{\sigma}_0^{-1}(\overline{\tau}_2(z));
\end{align*}
thus,
\begin{align*}
(\tau_1\tau_2)^\sigma(x) & =\Big{(}\tau_1^\sigma(\overline{\sigma}_0^{-1}(\overline{\tau}_2(z)))\Join\Lf{u}(\overline{\sigma}^{-1}_1(\Lf{t}(\tau_2(z)))\Big{)}\vJoin\tau_2^\sigma(x) \\
						 & =\Big{(}\tau_1^\sigma(\overline{\tau}_2^\sigma(x))\Join\Lf{u}\circ\Lf{t}(\tau_2^\sigma(x))\Big{)}\vJoin\tau_2^\sigma(x) \\
						 & =\tau_1^\sigma(\overline{\tau}_2^\sigma(x))\boxplus\tau_2^\sigma(x)=(\tau_1^\sigma\tau_2^\sigma)(x).
\end{align*}
\end{proof}
To easily prove the formulas involved, we prove the following algebraic lemma. Its contents will let us multiply $3\times 3$ arranges of squares in the preferred direction.
\begin{lemma}\label{shuffle}
Let $\begin{pmatrix}d_{11} & d_{12} & d_{13} \\
					d_{21} & d_{22} & d_{23} \\
					d_{31} & d_{32} & d_{33}\end{pmatrix}\in D^3_3$, then
\begin{align*}
(d_{11}\Join d_{12}\Join d_{13})\vJoin( & d_{21}\Join d_{22}\Join d_{23})\vJoin(d_{31}\Join d_{32}\Join d_{33}) \\
 & = (d_{11}\vJoin d_{21}\vJoin d_{31})\Join(d_{12}\vJoin d_{22}\vJoin d_{32})\Join(d_{13}\vJoin d_{23}\vJoin d_{33}).
\end{align*}					
\end{lemma}
\begin{proof}
This follows easily from the interchange law, 
\begin{align*}
(d_{11}\Join d_{12}\Join d_{13})\vJoin(d_{21}\Join d_{22}\Join d_{23}) & = ((d_{11}\Join d_{12})\Join d_{13})\vJoin((d_{21}\Join d_{22})\Join d_{23}) \\
 & = ((d_{11}\Join d_{12})\vJoin(d_{21}\Join d_{22}))\Join(d_{13}\vJoin d_{23}) \\
 & = (d_{11}\vJoin d_{21})\Join(d_{12}\vJoin d_{22})\Join(d_{13}\vJoin d_{23}) ,
\end{align*}	
so multiplying by the remaining term
\begin{align*}
((d_{11}\vJoin d_{21})\Join(d_{12} & \vJoin d_{22})\Join(d_{13}\vJoin d_{23}))\vJoin(d_{31}\Join d_{32}\Join d_{33}) \\
 & = \big{(}((d_{11}\vJoin d_{21})\Join(d_{12}\vJoin d_{22}))\Join(d_{13}\vJoin d_{23})\big{)}\vJoin\big{(}(d_{31}\Join d_{32})\Join d_{33}\big{)} \\
 & = \big{(}((d_{11}\vJoin d_{21})\Join(d_{12}\vJoin d_{22}))\vJoin(d_{31}\Join d_{32})\big{)}\Join\big{(}(d_{13}\vJoin d_{23})\vJoin d_{33}\big{)} \\
 & = \big{(}((d_{11}\vJoin d_{21})\vJoin d_{31})\Join((d_{12}\vJoin d_{22})\vJoin d_{32})\big{)}\Join(d_{13}\vJoin d_{23}\vJoin d_{33})
\end{align*}
which is what we wanted.

\end{proof}
We are ready to state the the theorem.
\begin{theorem}
$\Bis_{mult}(D)$ is a $2$-group.
\end{theorem}
\begin{proof}
We prove that the space of multiplicative sections coincide with the (semi-direct) product of the terms of the crossed module given by $\Delta$.
Due to lemma \ref{DeltaHomo} and proposition \ref{acts}, we are left to prove that $\Delta$ is equivariant and that the Peiffer equation is verified. To do so, let $\tau\in\Bis(\CC{D})$ and $\sigma\in\Bis_{mult}(D)_0$. On the one hand, for $x\in M$, 
\begin{align*}
(\Delta\tau^\sigma)_0(x) & =\Tp{t}(\tau^\sigma(x)) \\
						 & =\Tp{t}(\sigma_1^\dagger(\Lf{t}(\tau(\overline{\sigma}_0(x)))\vJoin\tau(\overline{\sigma}_0(x))\vJoin\Tp{u}(\sigma_0(x))) \\
						 & =\Tp{t}(\sigma_1^\dagger(\Lf{t}(\tau(\overline{\sigma}_0(x))))\cdot\Tp{t}(\tau(\overline{\sigma}_0(x)))\cdot\Tp{t}(\Tp{u}(\sigma_0(x))) \\
						 & =\sigma_0^\dagger(t(\Lf{t}(\tau(\overline{\sigma}_0(x))))\cdot(\Delta\tau)_0(\overline{\sigma}_0(x)))\cdot\sigma_0(x) \\
						 & =\sigma_0^\dagger(\Rg{t}((\Delta\tau)_0(\overline{\sigma}_0(x))))\cdot(\Delta\tau)_0(\overline{\sigma}_0(x)))\cdot\sigma_0(x)=(\sigma_0^\dagger(\Delta\tau)_0\sigma_0)(x).
\end{align*}
On the other, for $h\in H(x,y)$,
\begin{eqnarray}\label{EquivLHS}
(\Delta\tau^\sigma)_1(h)=\tau^\sigma(y)\Join\Lf{u}(h)\Join\Tp{\iota}(\tau^\sigma(x)).
\end{eqnarray}
We can factor the middle term as 
\begin{align*}
\Lf{u}(h) & = \Lf{\iota}(\sigma_1(h))\vJoin\sigma_1(h) \\
          & = \Lf{\iota}(\sigma_1(h))\vJoin\Lf{u}\circ\Lf{t}(\overline{\sigma}_1(h))\vJoin\sigma_1(h);
\end{align*} 
therefore, since clearly
\begin{eqnarray*}
\big{(}\Tp{u}(\sigma_0(y));\sigma_1(h);\Tp{u}(\sigma_0(x))\big{)},\big{(}\tau(\overline{\sigma}_0(y));\Lf{u}\circ\Lf{t}(\overline{\sigma}_1(h));\Tp{\iota}(\tau(\overline{\sigma}_0(x)))\big{)}\in D\times_V D\times_V D
\end{eqnarray*}
and
\begin{align*}
\Tp{s}(\sigma_1^\dagger(\Lf{t}(\tau(\overline{\sigma}_0(y)))) & =\sigma_0^\dagger(s(\Lf{t}(\tau(\overline{\sigma}_0(y)))) \\
 & =\sigma_0^\dagger(\Rg{t}(\Tp{s}(\tau(\overline{\sigma}_0(y)))) \\
 & =\sigma_0^\dagger(\Rg{t}(\Rg{u}(\overline{\sigma}_0(y)))) \\
 & =\sigma_0^\dagger(\overline{\sigma}_0(y)))) \\
 & =\Rg{\iota}(\sigma_0(y)) \\
 & =\Rg{\iota}(\Tp{t}(\sigma_1(h)))=\Tp{t}(\Lf{\iota}(\sigma_1(h)))
\end{align*}
as well as 
\begin{align*}
\Tp{t}(\Tp{\iota}(\sigma_1^\dagger(\Lf{t}(\tau(\overline{\sigma}_0(x)))))) & =\Tp{s}(\sigma_1^\dagger(\Lf{t}(\tau(\overline{\sigma}_0(x)))) \\
 & =\Rg{\iota}(\sigma_0(x)) \\
 & =\Rg{\iota}(\Tp{s}(\sigma_1(h)))=\Tp{s}(\Lf{\iota}(\sigma_1(h))),
\end{align*}
we have got that the right hand side of equation \ref{EquivLHS} is a $3\times 3$ arrange and we can use lemma \ref{shuffle}. In doing so, we get the vertical multiplication  
\begin{align*}
\Big{(}\sigma_1^\dagger(\Lf{t} & (\tau(\overline{\sigma}_0(y)))\Join\Lf{\iota}(\sigma_1(h))\Join\Tp{\iota}(\sigma_1^\dagger(\Lf{t}(\tau(\overline{\sigma}_0(x)))))\Big{)}\vJoin \\
 & \vJoin\Big{(}\tau(\overline{\sigma}_0(y))\Join\Lf{u}\circ\Lf{t}(\overline{\sigma}_1(h))\Join\Tp{\iota}(\tau(\overline{\sigma}_0(x)))\Big{)}\vJoin \\
 & \quad\qquad\qquad\qquad\vJoin\Big{(}\Tp{u}(\sigma_0(y))\Join\sigma_1(h)\Join\Tp{\iota}(\Tp{u}(\sigma_0(x)))\Big{)}.
\end{align*}
The bottom term is
\begin{eqnarray*}
\Tp{u}(\sigma_0(y))\Join\sigma_1(h)\Join\Tp{u}(\sigma_0(x))=\sigma_1(h),
\end{eqnarray*}
the middle term is 
\begin{eqnarray*}
\tau(\overline{\sigma}_0(y))\Join\Lf{u}\circ\Lf{t}(\overline{\sigma}_1(h))\Join\Tp{\iota}(\tau(\overline{\sigma}_0(x)))=(\Delta\tau)_1(\overline{\sigma}_1(h)),
\end{eqnarray*}
and the top term is
\begin{align*}
\sigma_1^\dagger(\Lf{t}( & \tau(\overline{\sigma}_0(y)))\Join\Lf{\iota}(\sigma_1(h))\Join\Tp{\iota}(\sigma_1^\dagger(\Lf{t}(\tau(\overline{\sigma}_0(x)))) \\
 & =\Lf{\iota}\Big{(}\sigma_1(\overline{\sigma}_1^{-1}\circ\Lf{t}(\tau(\overline{\sigma}_0(y)))\Join\sigma_1(h)\Join\Tp{\iota}(\sigma_1(\overline{\sigma}_1^{-1}\circ\Lf{t}(\tau(\overline{\sigma}_0(x))))\Big{)} \\
 & =\Lf{\iota}\circ\sigma_1\Big{(}\overline{\sigma}_1^{-1}\circ\Lf{t}(\tau(\overline{\sigma}_0(y)))\cdot\overline{\sigma}_1^{-1}(\overline{\sigma}_1(h))\cdot\iota(\overline{\sigma}_1^{-1}\circ\Lf{t}(\tau(\overline{\sigma}_0(x)))\Big{)} \\
 & =\sigma_1^\dagger\big{(}\Lf{t}(\tau(\overline{\sigma}_0(y)))\cdot\overline{\sigma}_1(h)\cdot\iota(\Lf{t}(\tau(\overline{\sigma}_0(x))))\big{)} \\
 & =\sigma_1^\dagger\big{(}\Lf{t}(\tau(\overline{\sigma}_0(y))\Join\Lf{u}(\sigma_1(h))\Join\Tp{\iota}(\tau(\overline{\sigma}_0(x))))\big{)}=\sigma_1^\dagger(\Lf{t}((\Delta\tau)_1)(\overline{\sigma}_1(h)));
\end{align*}
yielding 
\begin{eqnarray*}
(\Delta\tau^\sigma)_1(h)=(\sigma_1^\dagger(\Delta\tau)_1\sigma_1)(h)
\end{eqnarray*} 
and the the equivariance of $\Delta$. \\
We now turn to prove that the Peiffer equation holds, so let $\tau_1,\tau_2\in\Bis(\CC{D})$ and write for convenience 
\begin{eqnarray*}
\Lf{t}(\tau_k(x))=h_k(x), & \Tp{t}(\tau_k(x))=(\Delta\tau_k)_0(x)=v_k(x)\textnormal{ and} & t_{\CC{D}}(\tau_k(x))=x_k,
\end{eqnarray*}
that is $\tau_k(x)$ is a square of the form
\begin{eqnarray*}
\xymatrix{
{}^{x_k}\bullet            & \bullet^x \ar[l]_{h_k(x)} \\
{}_x\bullet\ar[u]^{v_k(x)} & \bullet_x \ar@{=}[l]\ar@{=}[u]
}
\end{eqnarray*}
Now, consider
\begin{align*}
(\tau_2^\dagger\tau_1\tau_2)(x) & =\tau_2^\dagger(x_{21})\boxplus\tau_1(x_2)\boxplus\tau_2(x) \\
 & =\Big{(}\tau_2^\dagger(x_{21})\Join\Lf{u}\circ\Lf{t}(\tau_1(x_2)\boxplus\tau_2(x))\Big{)}\vJoin(\tau_1(x_2)\boxplus\tau_2(x)) \\
 & =\Big{(}\tau_2^\dagger(x_{21})\Join\Lf{u}(h_1(x_2)h_2(x))\Big{)}\vJoin\Big{(}\big{(}\tau_1(x_2)\vJoin\Tp{u}(v_2(x))\big{)}\Join\tau_2(x)\Big{)} \\
\end{align*}
We now operate to the right with
\begin{eqnarray*}
\Join u^2(x)=\Join(\Lf{\iota}(\Tp{\iota}(\tau_2(x)))\vJoin\Tp{\iota}(\tau_2(x))).
\end{eqnarray*}
Since, 
\begin{align*}
\Tp{s}(\tau_2^\dagger(x_{21})\Join\Lf{u}(h_1(x_2)h_2(x))) & =\Tp{s}(\Lf{u}(h_2(x))) \\
	& =\Rg{u}(x) \\
	& =\Tp{s}(\tau_2(x)) \\
	& =\Rg{\iota}(\Tp{s}(\tau_2(x))) \\	
	& =\Rg{\iota}(\Tp{t}(\Tp{\iota}(\tau_2(x)))=\Tp{t}(\Lf{\iota}(\Tp{\iota}(\tau_2(x))) 
\end{align*}
and
\begin{eqnarray*}
\Tp{s}(\big{(}\tau_1(x_2)\vJoin\Tp{u}(v_2(x))\big{)}\Join\tau_2(x))=\Tp{s}(\tau_2(x))=\Tp{t}(\Tp{\iota}(\tau_2(x))),
\end{eqnarray*}
we can use the interchange law, thus getting
\begin{align*}
(\tau_2^\dagger\tau_1\tau_2)(x)=\Big{(}\tau_2^\dagger(x_{21})\Join & \Lf{u}(h_1(x_2)h_2(x))\Join\iota^2(\tau_2(x))\Big{)}\vJoin\\
 & \vJoin\Big{(}\big{(}\tau_1(x_2)\vJoin\Tp{u}(v_2(x))\big{)}\Join\tau_2(x)\Join\Tp{\iota}(\tau_2(x))\Big{)} \\
\end{align*}
The second term is clearly
\begin{eqnarray*}
\tau_1(x_2)\vJoin\Tp{u}(v_2(x))=\tau_1((\overline{\Delta\tau_2})_0(x))\vJoin\Tp{u}((\Delta\tau_2)_0(x));
\end{eqnarray*}
whereas, calling $y=\overline{\tau}_2^{-1}(x_{21})$ and expanding the first one, we get
\begin{align*}
\tau_2^\dagger(x_{21}) & \Join\Lf{u}(h_1(x_2)h_2(x))\Join\iota^2(\tau_2(x)) \\
 & =\iota_{\CC{D}}(\tau_2(y))\Join\Lf{u}(h_1(x_2)h_2(x))\Join\iota^2(\tau_2(x)) \\
 & =\Lf{\iota}(\tau_2(y))\Join\Lf{u}(\iota(\Lf{t}(\tau_2(y))))\Join\Lf{u}(h_1(x_2)h_2(x))\Join\Lf{\iota}(\Tp{\iota}(\tau_2(x))) \\
 & =\Lf{\iota}\Big{(}\tau_2(y))\Join\Lf{u}(\iota(h_2(y))h_1(x_2)h_2(x))\Join\Tp{\iota}(\tau_2(x))\Big{)} \\
 & =\Lf{\iota}((\Delta\tau_2)_1(\iota(h_2(y))h_1(x_2)h_2(x))).
\end{align*}
However, we claim that 
\begin{eqnarray*}
\iota(h_2(y))h_1(x_2)h_2(x) & =(\Lf{t}\circ(\Delta\tau_2)_1)^{-1}(h_1(x_2));
\end{eqnarray*}
thus, the result follows. Indeed, since $t(\iota(h_2(y)))=s(h_2(y))=y$, 
\begin{align*}
\Lf{t}\circ(\Delta\tau_2)_1(h_2(y)h_1(x_2)h_2(x)) & =\Lf{t}(\tau_2(y)\Join\Lf{u}(\iota(h_2(y))h_1(x_2)h_2(x))\Join\Tp{\iota}(\tau_2(x))) \\
												  & =h_2(y)\iota(h_2(y))h_1(x_2)h_2(x)\iota(h_2(x))=h_1(x_2),
\end{align*}
ultimately finishing the proof, as
\begin{eqnarray*}
(\tau_2^\dagger\tau_1\tau_2)(x)=(\Delta\tau_2)_1^\dagger(h_1(x_2))\vJoin\tau_1((\overline{\Delta\tau_2})_0(x))\vJoin\Tp{u}((\Delta\tau_2)_0(x))=(\tau_1)^{\Delta\tau_2}(x).
\end{eqnarray*}
\end{proof}
Let us see some examples of these $2$-groups of multiplicative bisections.
%It'd be kool to see what the corresponding thing is for the more geometrical objects such as Lu's double Lie group, Dirac double Lie groups,...; also, VB-groupoids. \\
\begin{ex:}\label{2GpAsBisections}
The $2$-group of multiplicative bisections of a Lie $2$-group $\xymatrix{\G \ar@<0.5ex>[r]\ar@<-0.5ex>[r] & H}$ is $\G$ itself. Indeed, functors from the category with one object and one arrow are points on the base with their respective identities. Therefore, $\Bis_{mult}(\G)_0\cong H$. On the other hand, the core of a Lie $2$-group is a Lie group. Let $G=\ker s$, since the group of bisections of a Lie group regarded as a Lie groupoid is the Lie group itself, $\Bis(G)=G$. Thus, it is not hard to see that the crossed module of the $2$-group of multiplicative bisections is $\xymatrix{:G \ar[r] & H}$, which is the same crossed module associated to $\G$.
\end{ex:}
\begin{ex:}\label{MultBiPair}
Interestingly, the multiplicative bisections of the pair double Lie groupoid $G\times G$ coincides with the group of Lie groupoid automorphisms. Indeed, a bisection $\sigma_G(g)=(\sigma_1(g),g)$ with base map $\sigma_M(x)=(\sigma_0(x),x)$ is multiplicative if 
\begin{eqnarray*}
s(\sigma_1(g))=\sigma_0(s(g)), & t(\sigma_1(g))=\sigma_0(t(g)), & u(\sigma_0(x))=\sigma_1(u(x))
\end{eqnarray*}
and if additionally, $\sigma_1(g_1g_2)=\sigma_1(g_1)\sigma_1(g_2)$ for all pairs of composable arrows $(g_1,g_2)\in G^{(2)}$. This extends the case of the group of bisections of the pair groupoid of a manifold $M$, which coincides with $Diff(M)$. The associated crossed module of the group of multiplicative bisections of the pair double Lie groupoid is 
\begin{eqnarray*}
\xymatrix{
\Bis(G) \ar[r] & Aut(G):\tau \ar@{|->}[r] & C_\tau ,
}
\end{eqnarray*}
where $C_\tau(g):=\tau(t(g))\cdot g\cdot \iota(\tau(s(g))$ is the conjugation induced by the bisection $\tau$. In this sense, the image of the crossed module map returns inner automorphisms. The action is a conjugation itself, this time around in the group of diffeomorphisms. It is given by $\tau^\sigma(x):=(\sigma_1^{-1}\circ\tau\circ\sigma_0)(x)$.
% The example of the fundamental double groupoid is harder already computed the one for S^1\times S^1=>S^1 (see 05.03), but the techniques involved 
\end{ex:}
\begin{ex:}\label{multBisVBGpd}
The simplest possible VB-groupoid is a $2$-vector space
\begin{eqnarray*}
\xymatrix{
\mathbb{V} \ar@<0.5ex>[r]\ar@<-0.5ex>[r]\ar@<0.5ex>[d]\ar@<-0.5ex>[d] & \ast \ar@<0.5ex>[d]\ar@<-0.5ex>[d] \\
V \ar@<0.5ex>[r]\ar@<-0.5ex>[r] & \ast
}
\end{eqnarray*}
Call $W:=\ker\Lf{s}$ and $\phi:=\Lf{t}\rest{W}$. Then $\mathbb{V}\cong W\oplus V$ and one can write the structural maps as 
\begin{eqnarray*}
s(w,v)=v, & \textnormal{and} & t(w,v)=v+\phi(w). 
\end{eqnarray*}
As remarked before, a functor between flat abelian Lie groupoids is but a map of vector bundles. In this particular case, it will be simply a linear map $\xymatrix{\sigma:V \ar[r] & W\oplus V}$. If $\sigma(v)=(\sigma_W(v),\sigma_V(v))$, the fact it is a section of the source map says that $\sigma_V=v$; on the other hand, the defining condition of a bisections implies $\Bis_{mult}(\mathbb{V})_0=\lbrace A\in Hom(V,W):I+\phi A\in GL(V)\rbrace$. The core of $\mathbb{V}$ will be the vector space $W$ regarded as a Lie group; hence, its space of bisections is $W$ again. The crossed module of the group of multiplicative bisections will have zero structural map, and action $w^A=w-A(I+\phi A)^{-1}\phi(w)$. We would like to remark that the formula for the group structure in the group of bisections is given by
\begin{eqnarray*}
A_1\odot A_2=A_1+A_2+A_1\phi A_2.
\end{eqnarray*}
In general, for a VB-groupoid over an arbitrary base, after having picked an splitting (or equivalently, determined an associated representation up to homotopy), the group of multiplicative bisections coincides with 
\begin{eqnarray*}
\Gamma_{mult}(s^*E\oplus t^*C)=\lbrace (\sigma_0,\sigma_1)\in\Bis(G)\times Hom(E,C):\partial\circ\sigma_1 +\Delta^E_{\sigma_0}\in GL(E)\rbrace
\end{eqnarray*}
To make sense of the condition, it should be noted that $Hom$ stands for the space of all vector bundle morphisms, as opposed to only those covering the identity. Indeed, in order for the sum to be possible, $\sigma_1$ needs cover $t\circ\sigma_0$. The core of such a VB-groupoid is the vector bundle $C$ regarded as a Lie groupoid; thus, its space of bisections is the vector space of sections $\Gamma(C)$. The crossed module associated to the $2$-group of multiplicative bisections will have trivial structural map, and the action will be a rather involved formula, though analogous to that of the $2$-vector space.
% Find the formula, and verify that the structural map is indeed zero.
\end{ex:}
\begin{ex:}\label{MultBiSectOfPoisson-LieGp}
A functor from a Lie group to a Lie groupoid is a group homomorphism to the isotropy Lie group of the image of its identity element. Say, $\Gamma$ is the double Lie groupoid integrating a Poisson Lie group $G$, and let $\sigma\in\Bis_{mult}(\Gamma)_0$. Call $v:=\sigma(1)$, then $\sigma(g)\in\Gamma_{v}$. By definition,
\begin{eqnarray*}
\Gamma_v\cong\lbrace(g,h)\in G\times G:\dinc{v}\inc{g}=\inc{h}\dinc{v}\rbrace .
\end{eqnarray*}
Using this identification, $\sigma(g)=(\varphi(g),g)$, where $\varphi\in Diff(G)$. The defining property says that 
\begin{eqnarray*}
\inc{\varphi(g)}=\dinc{v}\inc{g}(\dinc{v})^{-1}
\end{eqnarray*}
in $D$. In so, $\varphi$ can only be defined in terms of inner automorphisms of $D$ by the elements of $G^*$ regarded as embedded in $D$. These inner automorphisms do not necessarily leave (the image of) $G$ fixed; thus, we can write
\begin{eqnarray*}
\Bis_{mult}(\Gamma)_0:=\lbrace v\in G^*:C_{\dinc{v}}^D(G)=(G)\rbrace,
\end{eqnarray*} 
where $C^D$ stands for the conjugation in $D$. Poisson-Lie groups and their integrations are vacant; hence, the $2$-group of multiplicative bisections is the unit group $\Bis_{mult}(\Gamma)_0$.
\end{ex:}
\begin{ex:}\label{actionLABiSections}
The space of multiplicative bisections of a Lie groupoid integrating an action \LA -groupoid of $G$ on an integrable Lie algebroid $A$, is $\Bis(\check{\G}_A)^G$. Indeed, every multiplicative bisection $\sigma\in\Bis_{mult}(G\ltimes\check{\G}_A)$ has components $\sigma_1(g,x)=(g,\sigma_{\check{\G}_A}(g,x))$ and has to be compatible with the source and the target, then
\begin{eqnarray*}
\Tp{s}(\sigma_1(g,x))=\sigma_{\check{\G}_A}(g,x))=\sigma_0(x) & \textnormal{and} & \Tp{t}(\sigma_1(g,x))=g\bullet\sigma_0(x)=\sigma_0(g\cdot x).
\end{eqnarray*}
Therefore, $\sigma$ is specified by its base component $\sigma_0$ and the fact it is equivariant. Since we saw that the core of an action \LA -groupoids is trivial, the crossed module of its $2$-group is $\xymatrix{1 \ar[r] & \Bis(\check{\G}_A)^G}$.
\end{ex:}
%In sight of these, it appears that there is a relation. Is the space of multiplicative bisections a Lie $2$-group of some kind?

%----------------------------------------------------
%---------------------------------------------------
%------------------------------------------------------
%------------------------------------------------------
%\section*{Sugestões para Pesquisas Futuras} 

%Oh! The amount of questions left unanswered.\\
%With no particular order of preference,\\
%1) Is the deformation theory of (strict) Lie $2$-algebras controlled y $H^2(\mu ,ad)$? \\
%2) Are this cohomology classes somehow related to some geometric data? e.g. X. Tang's non-commutative Poisson structures. (I am thinking of the Poisson structure as a Lie algebroid cohomology class on $T^*M$... not any class of course... is it even a class?) \\
%3) Is there a way of relating X. Tang's work with some sort of section on an \LA -groupoid over the orbifold? Furthermore, if one starts out with a compatible action of a finite Lie group on a Poisson manifold, so that the action groupoid has a compatible structure, does the associated \LA -groupoid has the prescribed sections? Does this help when quantizing? \\
%4) What the hell is Santiago Canez reduction? \\
%5) Is it worth understanding the enveloping algebras of more general classes of $2$-term  $L_\infty$-algebras? \\
%6) About the product in cohomology: Check the first notebook, right before Aline's Li\c{c}\~ao 2. 
       % associado ao arquivo: 'cap-conclusoes.tex'
% ---------------------------------------------------------------------------- %
% Bibliografia
\backmatter 
%---------------%
\fancyhead[RE,LO]{} % ESTE FANCYHEAD QUITA LA NUMERACION DE LA BIBLIOGRAFIA
%---------------%
\singlespacing                               % espaçamento simples
\bibliographystyle{plain}
\bibliography{biblos.bib} 		             % associado ao arquivo: 'bibliografia.bib' I hope biblos.bib

% ---------------------------------------------------------------------------- %
\end{document}